\documentclass[11pt]{article}

\usepackage{amssymb,amsfonts,amsmath,mathrsfs,mathtools}
\usepackage{graphicx,epsfig,subfigure}
\usepackage{bm}
\usepackage[margin=1in,papersize={8.5in,11in}]{geometry}
\usepackage{times}

\newcommand{\vs}{\vspace{0.7cm}}
\makeatother

\DeclareMathOperator*{\argmin}{arg\,min}
\DeclareMathOperator*{\argmax}{arg\,max}

\title{The Numerical Approximation of Nonlinear Functionals and Functional Differential Equations}

\author{Daniele Venturi\\
\em  Department of Applied Mathematics and Statistics \\
University of California, Santa Cruz}

\usepackage[usenames,dvipsnames]{xcolor}
\definecolor{r}{rgb}{0.0,0.0,0.0} 
\definecolor{b}{rgb}{0.0,0.0,0.0} 

\begin{document}
 
\maketitle

\begin{abstract}
The fundamental importance of functional differential equations has 
been recognized in many areas of mathematical physics, such as 
fluid dynamics (Hopf characteristic functional equation), 
quantum field theory (Schwinger-Dyson equations) and
statistical physics (equations for generating functionals and 
effective Fokker-Planck equations).
However, no effective numerical method has yet
been developed to compute their solution.
The purpose of this report is to fill this gap, and 
provide a new perspective on the problem of numerical approximation 
of nonlinear functionals and functional differential equations. 
\end{abstract}

\vs
\noindent
\centerline{\line(1,0){200}}\vspace{0.0cm}
\tableofcontents

\vspace{0.6cm}
\noindent
\centerline{\line(1,0){200}}\vspace{0.0cm}

\section{Introduction}
\label{sec:motivation}

In this report we address a rather neglected but very important 
research area in computational mathematics, namely the numerical 
approximation of nonlinear functionals and functional differential 
equations (FDEs). 
FDEs are arise naturally in many different areas of mathematical physics.
For example, in the context of fluid dynamics,
the Hopf equation \cite{Hopf}  
{\color{r}
\begin{equation}
\frac{\partial \Phi([\bm \theta],t)}{\partial t}=
\sum_{k=1}^3\int_V\theta_k(\bm x)\left(i \sum_{j=1}^3\frac{\partial }{\partial x_j}
\frac{\delta^2 \Phi([\bm \theta],t)}{\delta \theta_k(\bm x)\delta\theta_j(\bm x)}
+\nu \nabla^2\frac{\delta \Phi([\bm \theta],t)}{\delta \theta_k(\bm x)}\right)d\bm x,
\label{hopfns}
\end{equation}
was deemed by Monin and Yaglom (\cite{Monin2}, Ch. 10) 
to be ``the most compact formulation 
of the general turbulence problem'', which is the problem of 
determining the statistical properties of the velocity 
and the pressure fields of Navier-Stokes equations given 
statistical information on the initial state\footnote{Stani\v{s}i\'c \cite{Stanisic} refers to the Hopf equation \eqref{hopfns} 
as the ``only exact formulation in the entire field of 
turbulence''  (Ch. 12, p. 233).}.
In equation \eqref{hopfns} $V\subset \mathbb{R}^3$ is 
a periodic box, $\bm 
\theta(\bm x)=(\theta_1(\bm x), \theta_2(\bm x), 
\theta_3(\bm x))$ 
is a vector-valued test function in a suitable divergence-free 
space, and $\Phi$ is a nonlinear complex-valued functional 
known as Hopf functional \cite{Hopf}. Remarkably, 
with such functional available it is possible to 
compute any statistical property of 
the velocity field that solves the Navier-Stokes 
equations (see \cite{Monin2}). 
This is of great conceptual importance: the solution to 
one single linear functional differential equation can 
describe all statistical features of turbulence and there 
is no need to refer back to the Navier-Stokes 
equations.
From a mathematical viewpoint the Hopf
functional is basically a time-dependent nonlinear operator 
in the space of test functions $D(\Phi)$ (domain of the operator 
$\Phi$) with range in the complex plane (see Figure \ref{fig:1}).
\begin{figure}
\centerline{
\includegraphics[height=6.5cm]{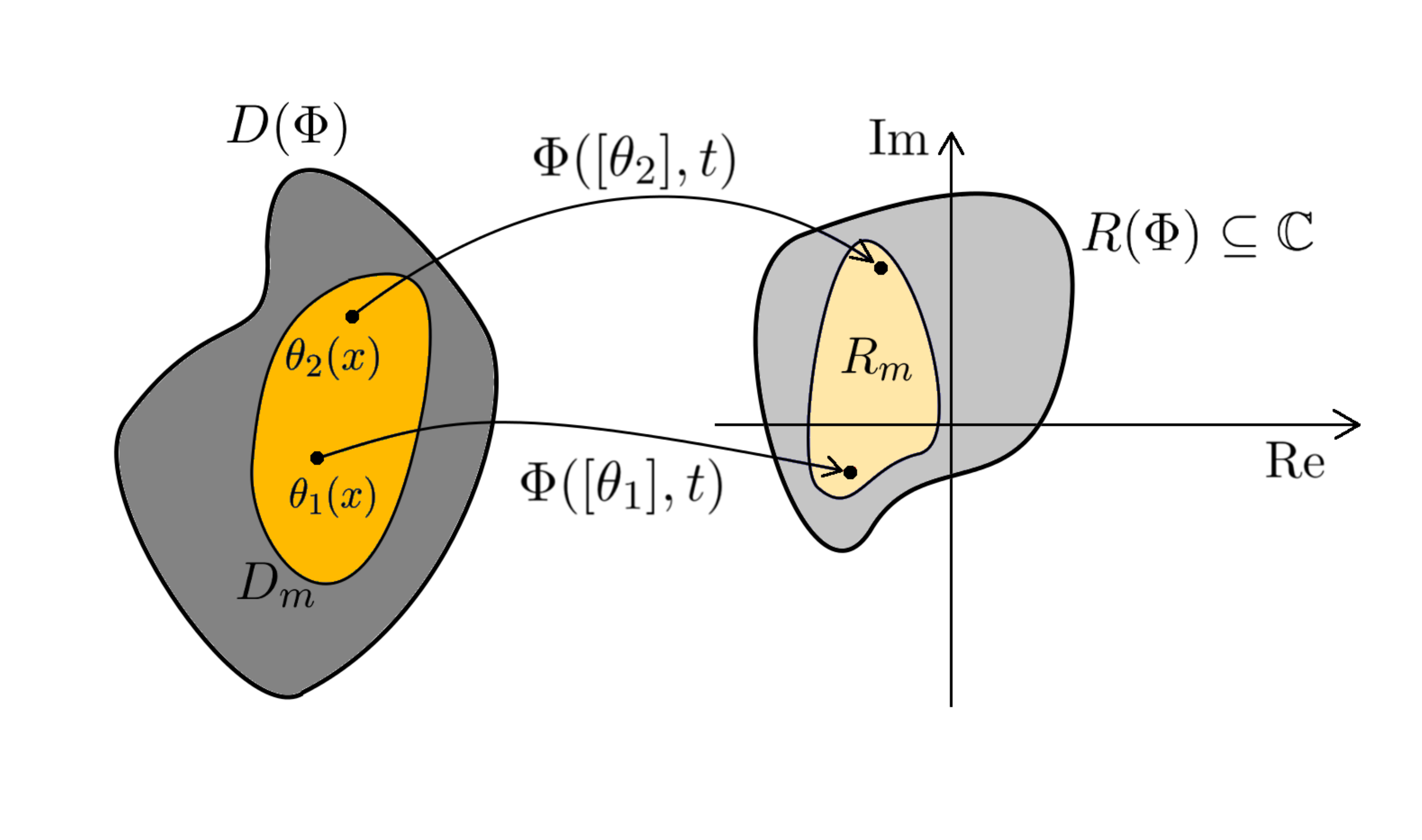}}
\caption{ 
Sketch of the mapping at the basis of the 
Hopf functional. The domain of the functional $\Phi$, 
denoted as $D(\Phi)$ is a suitable 
of space of functions while the range of $\Phi$, denoted as 
$R(\Phi)$, is a subset of the complex plane $\mathbb{C}$. The 
approximation space $D_m$
is a subset of $D(\Phi)$, which is mapped into $R_m$.}
\label{fig:1}
\end{figure}
The operator can be formally defined as a functional integral
\begin{equation}
\Phi([\bm \theta],t)=\int_{\Omega} \exp\left[i\int_{V}
\bm u (\bm x,t;\omega)\cdot \bm \theta(\bm x)d\bm x\right]P([\bm u_0]) \mathcal{D}[\bm u_0]
\label{HcF}
\end{equation}
where $\bm u(\bm x,t; \omega)$ is a stochastic solution to the 
Navier-Stokes equations and $P([\bm u_0])$ is the 
probability functional of the random initial state 
(assuming it exists). 
Thus, computing the solution to the Hopf equation \eqref{hopfns} 
is equivalent to compute a (complex-valued) 
time-dependent nonlinear operator $\Phi$ from an 
equation that involves classical 
partial derivatives with respect to space and time variables 
as well as derivatives with respect to functions, i.e., 
functional derivatives $\delta /\delta 
\theta_j(\bm x)$ \cite{Vainberg,Nashed}.}

Another example of functional differential equation 
is the Schwinger-Dyson equation of quantum field 
theory \cite{Easther,Justin}. Such equation describes the 
exact dynamics of the Green functions of a general field theory, 
and it allows us to propagate field interactions, either in a 
perturbation setting \cite{Okopinska} (weak coupling regime) or 
in a strong coupling regime \cite{Swanson}. The Schwinger-Dyson 
formalism is also useful in computing the statistical properties of 
stochastic dynamical systems. For example, consider Langevin 
equation 
\begin{equation}
\frac{d \bm \psi(t)}{dt}=\bm G(\bm \psi(t),t)+\bm f(t;\omega),
\label{ODE}
\end{equation}
where $\bm f(t;\omega)$ is random noise. 
Define the generating functional \cite{Phythian,Jensen}
\begin{equation}
Z([\bm \xi,\bm \eta])=Z_0 \int \mathcal{D}[\bm \psi] \mathcal{D}[\bm \beta]A([\bm \psi,\bm \beta])\exp\left[\int_0^t d\tau \left(\bm \xi(\tau)\cdot \bm \psi(\tau)+\bm \eta(\tau)\cdot \bm \beta(\tau)
\right)\right],
\label{GF}
\end{equation}
 where $Z_0$ is a normalization constant and
 \begin{equation}
A([\bm \psi,\bm \beta])= C([\bm \beta])\exp\left[-\frac{1}{2}
\int_0^t d\tau \nabla \cdot \bm G(\bm \psi(\tau),\tau) -i 
\int_0^t d\tau \bm \beta(\tau)\cdot\left(\frac{d\bm \psi(\tau)}{dt}-\bm G(\bm \psi(\tau),\tau) \right)\right].
\label{AU} 
 \end{equation}  
The functional $C([\bm \beta])$ in \eqref{AU} denotes the (known) characteristic functional of the external random noise $\bm f(t;\omega)$. Clearly, if we have available the stochastic solution to \eqref{ODE}, then we can construct the functional $Z([\bm \xi, \bm \eta])$ and compute all statistical properties we are interested in. On the other hand, it is straightforward to show that $Z([\bm \xi, \bm \eta])$ satisfies the following system of linear FDEs\footnote{
The expression
\begin{equation}
G_k\left(\frac{1}{i}\frac{\delta}{\delta \bm \xi (\tau)},\tau\right)
\end{equation}
in equations \eqref{SD1} and \eqref{SD2} has
to be interpreted in the sense of symbolic operators. 
For example, in one dimension, if $G\left(z,t\right)=z+z^2$ then 
\begin{equation}
G\left(\frac{1}{i}\frac{\delta}{\delta \xi (\tau)},t\right)Z = \frac{1}{i}\frac{\delta  Z}{\delta 
\xi(\tau)}-\frac{\delta^2 Z}{\delta \xi(\tau)^2}.
\end{equation} }
(Schwinger-Dyson equations)
\begin{align}
\frac{1}{i}\frac{\partial }{\partial \tau}\frac{\delta Z}{\delta \xi_k(\tau)}
&=\eta_k(\tau)Z + G_k\left(\frac{1}{i}\frac{\delta}{\delta \bm \xi(\tau)},
\tau\right)Z - iD_k\left(\left[\frac{1}{i}\frac{\delta}{\delta \bm 
\eta(\tau)}\right],\tau\right)Z,\label{SD1}\\
\frac{1}{i}\frac{\partial }{\partial \tau}\frac{\delta Z}{\delta \eta_k(\tau)}&=- \xi_k(\tau)Z + i\sum_{j=1}^n\frac{\delta}{\delta\eta_j(\tau^+)}\frac{\partial G_k}{\partial \psi_j}\left(\frac{1}{i}\frac{\delta}{\delta \bm \xi(\tau)},\tau\right)Z,\label{SD2}
\end{align}
where
\begin{equation}
D_i([\bm \beta],\tau)=\frac{\delta C([\bm \beta])}{\delta \beta_i(\tau)}.
\end{equation}
The solution to the Schwinger-Dyson equations 
\eqref{SD1}-\eqref{SD2} is a nonlinear functional 
(i.e., a nonlinear operator) $Z([\bm \xi, \bm \eta])$ which allows 
us to compute all statistical properties of the system without any 
knowledge of the stochastic process $\bm \psi(t;\omega)$ defined 
implicitly by the stochastic ODE \eqref{ODE}. By generalizing 
\eqref{GF}, it is possible to derive a functional formalism for any 
classical field theory or stochastic system. This yields, in particular, 
Schwinger-Dyson-type equations for generating functionals associated 
with the solution to stochastic partial differential equations (SPDEs). 
If the SPDE admits an action functional, then the construction of the 
generating functional as well as the derivation of the corresponding 
Schwinger-Dyson equation are rather straightforward (see 
\cite{Jensen,Amit,Kleinert}). 

The usage of functional differential equations
grew very rapidly during the sixties, when it became
clear that techniques developed for quantum field theory
by Dyson, Feynman, and Schwinger could be applied, at least
formally, to other branches of mathematical physics.
The seminal work of Martin, Siggia, and Rose \cite{Martin} became
a landmark on this subject, since it revealed the possibility of
applying (at least formally) quantum field theoretic methods, such as
functional integrals and diagrammatic
expansions \cite{Phythian,Jensen,Phythian1,Jouvet},
to classical physics. Relevant applications of
these techniques can be found in non-equilibrium statistical mechanics
\cite{Jensen,Phythian,Phythian1,Jouvet,Eyink_1996,Langouche,Ueda,Ueda1},
stochastic dynamics \cite{Hochberg,Venturi_PRS,Klyatskin1}, 
and turbulence theory
\cite{Frisch,McComb,Eyink,Chen_KraichnanPRL,Giles,Dopazo,Alankus,Lewis, Monin1,Monin2,Rosen_1960,Rosen_1967,Rosen_1969,Hopf_book,Hosokawa}.

An open question that has persisted over the years is:  How do we 
compute the solution to a functional differential equation? 
From the fifties to the eighties, 
researchers were of course investigating analytical methods, e.g., 
based on functional power series 
\cite{Rosen_1960,Rosen_1967,Volterra,Nelles}, functional integrals 
\cite{Kleinert,Lobanov,Popov,Egorov}, transforms with respect to 
appropriate measures (\cite{Monin2}, p. 802), and diagrammatic 
expansions. 
More recently, Waclawczyc and Oberlack \cite{Oberlack,Waclawczyk} 
proposed a Lie group analysis and applied it to the Hopf-Burgers equation, 
which represents a step forward toward developing new analytical solution 
methods. Specifically, invariant solutions of the Hopf-Burgers equation were 
found based on the analysis of the infinitesimal generator of 
suitable symmetry transformations. 
From a numerical viewpoint, 
recent advances in computational mathematics -- in particular in 
numerical tensor methods \cite{Hackbusch_book} -- open the possibility 
to solve functional differential equations on a computer. 
In this report, we will present state-of-the-art 
mathematical techniques and numerical algorithms to 
represent nonlinear functionals and compute the 
numerical solution to functional differential equations. 

If FDEs are so important, why do they not have a prominent role in 
computational mathematics? There are several  possible answers to this 
question. First of all, FDEs are infinite-dimensional equations, in the sense that 
they are, in principle, equivalent to an infinite-dimensional system of PDEs, or 
PDEs in an infinite number of variables. This 
may have understandably discouraged researchers in numerical analysis to 
even attempt a numerical discretization. Most schemes proposed so far are 
based on truncations of infinite hierarchies of PDEs obtained, e.g,
from functional power series expansions 
\cite{Rosen_1960,Rosen_1967,Rosen_1969,Rosen_1971,Monin2,Ahmadi,Frisch}, 
or Lundgren-Monin-Novikov hierarchies 
\cite{Wacawczyk1,Friedrich,Lundgren,Hosokawa,Rosteck}. 
Other approaches are based on a direct discretization of the 
functional integral \cite{Kleinert,Egorov,Popov,Langouche} 
that defines the field theory (e.g. $Z$ in 
equation \eqref{GF}), and its evaluation using Monte Carlo 
methods, or source Galerkin 
methods \cite{Guralnik1,Guralnik2}.
Dealing with systems of infinitely many PDEs or 
very high-dimensional PDEs can indeed be discouraging,  
but nowadays it is quite common, 
for example when discretizing stochastic systems driven by colored 
random noise or stochastic partial differential 
equations (SPDEs) \cite{HandyGK,Db_book,Xiao2,db2,DanieleWick}. 
Another reason why FDEs have not yet been numerically studied 
extensively may be due to a lack of awareness of their existence 
within the computational mathematics community. Also, there is no 
universal agreement across scientific disciplines as to even the basic 
definition of an FDE. For example, most applied mathematics literature 
refers to FDEs as ordinary differential equations with memory or delay 
terms \cite{Wu,Azbelev,Hale}. In the pure mathematics community, 
functional equations have been studied in the context of approximate 
homeomorphisms (the Ulam stability problem) \cite{Rassias1,Rassias}, or 
more generally within problems where the unknown is a function, e.g., 
Cauchy or d’Alambert functional equations \cite{Sahoo}. The physics 
literature, on the other hand, clearly identifies FDEs as those equations 
whose unknown is a functional (i.e., a nonlinear operator) and that involve 
partial derivatives with respect to independent variables (e.g., space and 
time), as well as derivatives with respect to functions (functional 
derivatives). These kinds of equations are usually far more challenging 
than the functional equations studied by the pure mathematics community, 
and indeed there are very few general theorems on the existence and the uniqueness of their solution \cite{Hale,Foias}.

In this report we take the physicist viewpoint
and consider linear functional
differential equations in the form
\begin{equation}
\frac{\partial F([\bm \theta],t)}{\partial t} =L([\bm \theta ],t)F([\bm \theta ],t)+H([\bm \theta θ],t), \qquad F([\bm \theta],0)=F_0([\bm \theta]),
\label{functDE}
\end{equation}
where $F ([\bm \theta], t)$ is a real or complex-valued functional (time-dependent nonlinear operator in a space of functions), $F_0([\bm \theta])$ is a given initial condition, $L([\bm \theta],t)$ is a linear operator in the space of nonlinear functionals, and $H([\bm \theta],t)$ is a known forcing functional. 
The linear operator $L([\bm \theta], t)$ usually involves functional derivatives with respect to $\bm \theta(\bm x)$ as well as partial derivatives with respect to independent variables, e.g,  
space and time coordinates. 
For example, $L([\bm \theta],t)$ 
could be the linear operator defining 
the right hand side of equation \eqref{hopfns}.
We emphasize that the class of equations in the 
form \eqref{functDE} is very broad as it encompasses 
FDEs describing many physical systems, including statistical properties of 
nonlinear SODEs and SPDEs (e.g., Hopf characteristic functional equations
\cite{Hopf,Hopf_book,Monin2,Klyatskin1} or equations for probability 
density functionals \cite{Giles,Beran,Dopazo}), functional equations arising 
in control theory \cite{Bensoussan_DaPrato}, 
generalized principles of least actions \cite{Daniele_JMathPhys}, 
and functional equations of quantum field theory
\cite{Justin,Itzykson,Phythian,Jouvet,Jensen,Langouche}.

To the best of our knowledge, no effective numerical methods have yet 
been developed to compute the solution to linear functional differential 
equations in the form \eqref{functDE}, and little has been done 
for functional differential equations in general, despite their 
fundamental importance in many areas of 
mathematical physics.
The purpose of this report is to fill this gap and 
present state-of-the-art mathematical techniques, including 
new classes of numerical algorithms, to approximate
nonlinear functionals and the numerical 
solution to functional differential equations in 
the form \eqref{functDE}
This report is organized in two parts: 
\begin{enumerate}

\item{\bf Approximation of Nonlinear Functionals}
A nonlinear functional is a particular type 
of nonlinear operator from a space of functions into 
a vector space, e.g., $\mathbb{R}$ or $\mathbb{C}$.
Therefore, the process of approximating a nonlinear functional
is basically the same as approximating a nonlinear operator
\cite{Howlett,Torokhti,Bensoussan_DaPrato}.
In this report we will present various techniques for nonlinear functional 
approximation, ranging from polynomial functional series expansions, 
to expansions based on stochastic processes, and 
functional tensor methods. Within the context of polynomial 
functional series expansions we will discuss in particular Lagrange 
interpolation in Hilbert and Banach spaces 
\cite{Makarov,Prenter,Prenter1,Bertuzzi,Allasia,Chaika}, 
where the interpolation ``nodes'' are functions in a suitable 
function space. We will also discuss series expansions based 
on functional tensor methods \cite{Hackbusch_book}. 
This class of methods relies on recent developments on 
multivariate function approximation such as canonical polyadic (CP) 
\cite{Reynolds,Beylkin,Kolda} and hierarchical Tucker (HT) 
\cite{Bachmayr,Grasedyck2017} series expansions.

\item {\bf Approximation of Functional Differential Equations}
\label{item:2}
A functional differential equation is an equation whose 
solution is a nonlinear functional, i.e., a nonlinear operator. 
The equation usually involves functional derivatives of 
such functional and partial derivatives 
and integrals with respect to independent variables. 
The goal of approximation theory for functional differential 
equations is therefore to determine an approximation of 
such nonlinear functional, e.g., its time evolution given an 
initial state.
In this report  we will discuss new classes of methods that
extend classical Galerkin, least-squares, and collocation
techniques to functional differential equations.
These methods are based on suitable 
representations of the solution functional, e.g., 
in terms of polynomial functionals or functional 
tensor networks.
\end{enumerate}

This paper is organized as follows:
In Section \ref{sec:NonlinearFunctionals} we briefly review 
what nonlinear functionals, functional derivatives and 
provide useful examples of nonlinear functionals in physics.
In Section \ref{sec:rep-Hopf} we address the approximation 
of nonlinear functionals and functional derivatives.
In particular, we discuss Lagrange interpolation in
spaces of infinite dimensions (function spaces),
series expansions in terms of polynomial functionals, and
functional tensor methods such as canonical polyadic 
and hierarchical Tucker expansions. 
In Section \ref{sec:FDEs approximation} we discuss functional 
differential equations. We begin by presenting several examples 
of FDEs and show how they arise in the context of 
well-known physical theories. 
In Section \ref{sec:FDEapprox} we address the problem 
of computing the numerical solution to an FDE. Specifically, we introduce infinite-dimensional extensions of least squares, Galerkin, and collocation methods. 
Finally, in Section \ref{sec:numerical results functionals} and 
Section \ref{sec:numerical results functional 
equations} we present numerical results on nonlinear functional 
approximation and also compute the numerical solution 
to a prototype functional advection-reaction problem.

\section{Nonlinear Functionals}
\label{sec:NonlinearFunctionals}

Let $X$ be a Banach space of functions. A nonlinear functional 
on $X$ is a nonlinear operator $F$ that takes in an element
$\theta$ of $X$ (i.e., a function), and returns a real or a 
complex number. The functional $F$ usually does not operate 
on the entire linear space $X$ but rather on a subset set of$X$, 
which we denote as $D(F)\subseteq X$ (domain of the functional).   
Let us first provide simple examples of nonlinear functionals.

\vs
\noindent
{\em Example 1:} 
Consider 
\begin{equation}
F([\theta])=\int_0^1 x^3 e^{-\theta(x)}dx, \qquad 
\theta\in D(F)=C^{(0)}([0,1]),
\label{functional1}
\end{equation}
where $C^{(0)}([0,1])$ is the space of continuous 
functions in $[0,1]$. 

\vs
\noindent
{\em Example 2 (homogeneous polynomial functional of order $n$):} Consider
\begin{equation}
P_n([\theta])=\underbrace{\int_0^1\cdots \int_0^1}_{\textrm{$n$ times}}
K_n(x_1,...,x_n)\theta(x_1)\cdots \theta(x_n)dx_1\cdots dx_n,
\label{functional2} 
\end{equation}
where $K_n$ are given kernel functions.

{\color{r}
\vs
\noindent
{\em Example 3 (Ginzburg-Landau energy functional):} 
The Ginzburg-Landau theory describes  
describes phase transitions and critical phenomena in a great 
variety of statistical systems ranging from magnetic systems, to diluted polymers and superconductors \cite{KleinertPhi4,Amit}. At the basis of the theory is the energy functional 
\begin{equation}
E([\phi]) = \int_{\mathbb{R}^3} 
\left[\frac{1}{2}\left(\left|\nabla \phi(\bm x)\right|^2+m^2\phi^2(\bm x)\right)+\frac{\lambda}{4!}\phi^4(\bm x)\right]d\bm x
\end{equation}
where $m$ is the ``mass'' of the field $\phi$ and $\lambda$ 
is a coupling constant. In this case, the domain of the 
functional $E$ can be chosen as $D(E)=C^{(2)}(\mathbb{R}^3)$.

\vs
\noindent
{\em Example 4 (Hopf characteristic functional of a
Gaussian random field):} 
Consider a scalar Gaussian  random field $u(\bm x;\omega)$ defined on a domain $V\subseteq \mathbb{R}^3$ it can be shown (see, e.g., \cite{Klyatskin1}) that the Hopf characteristic functional of such 
random field is
\begin{equation}
\Phi([\theta])=\exp\left[i\int_V\mu(\bm x)\theta(\bm x)d\bm x -
\frac{1}{2}\int_V \int_V C(\bm x,\bm y)\theta(\bm x)\theta(\bm y) d\bm xd\bm y\right],
\label{white_0}
\end{equation}
where $\mu(\bm x)$ and $C(\bm x,\bm y)$ denote, respectively, 
the mean and the covariance function of the field $u(\bm x;\omega)$.

\vs
\noindent
Analysis of nonlinear functionals in Banach spaces 
is a well-developed subject \cite{Vainberg,Nashed,Schwartz1964}.
In particular, the classical definition of continuity and 
differentiability at a point that holds for real-valued functions 
can extended in a more or less straightforward way to nonlinear 
functionals. For instance, we say that a functional 
$F([\theta])$ is {\em continuous} at a point $\theta(x)\in D(F)$ 
if for any sequence of functions 
$\{\theta_1(x), \theta_2(x), ... \}$ in $D(F)$) converging 
to $\theta$ we have that the sequence $\{F([\theta_n])\}$ 
converges to $\{F([\theta])\}$, i.e., 
\begin{equation}
\lim_{n\rightarrow \infty}
\left|\theta_n(x)-\theta(x)\right|\rightarrow 0 \quad \Rightarrow \quad \lim_{n\rightarrow \infty}
\left|F([\theta_n])-F([\theta])\right|\rightarrow 0.
\label{continuity}
\end{equation}
The functionals we discussed in Examples 1-4 are all continuous.
From the continuity definition \eqref{continuity} it follows, 
in particular, that if $F([\theta])$ is continuous on a compact 
function space $D(F)$ then $F$ is {\em bounded}. 

\vs
\noindent
{\em Example 4:}
Another example of a continuous functional 
in $D(F)=C^{(\infty)}([-1,1])$ is
\begin{equation}
 F([\theta]) = \theta(0)^2+\int_{-1}^1 \sin(\theta(x))dx.
 \label{ex4fun}
\end{equation}
In fact, consider any sequence of functions 
$\{\theta_n\}$ in $D(F)$. Also choose 
an integrable function $g(x)$ such that 
$\left|\theta_n(x)\right|\leq g(x)$ for all $n\in \mathbb{N}$. 
In these conditions, the Lebesgue dominated convergence 
theorem applies and we have
\begin{equation}
\lim_{n\rightarrow\infty} \int_{-1}^1 \sin(\theta_n(x))dx=\int_{-1}^1 
\sin(\theta(x))dx, 
\end{equation}
i.e.,  
\begin{equation}
\theta_n\substack\rightarrow \theta\,  
\Rightarrow \, F([\theta_n])\rightarrow F([\theta]).
\end{equation}

}

\subsection{Functional Derivatives}
\label{app:functional derivatives}
Consider a real or a complex valued functional $F$ defined on 
the function space $D(F)$ (domain of the functional). For 
simplicity, let us assume that $D(F)$ is a space of real valued 
functions $\theta(x)$ on the real line. 
We say that the functional $F$ is differentiable 
at $\theta(x)$ if the limit  
\begin{equation}
\lim_{\epsilon\rightarrow 0}\frac{F([\theta(x)+\epsilon \eta(x)])
-F([\theta(x)])}{\epsilon}
\label{gateaux}
\end{equation}
exists and it is finite. The quantity \eqref{gateaux} 
is known as {\em G\^ateaux differential} of $F$ in the 
direction of $\eta(x)$ (see, e.g., \cite{Vainberg,Schwartz1964}). 
Under rather general assumptions such derivative can 
be represented as a linear operator \cite{Nashed,Vainberg,Daniele_JMathPhys} 
acting on $\eta(x)$. For small $\epsilon$ we have
\begin{equation}
F([\theta(x)+\epsilon \eta(x)])=F([\theta(x)])+\epsilon
L([\theta(x)])\eta(x)+R_1([\epsilon\eta(x);\theta(x)]),
\end{equation}
In this series $L([\theta(x)])$ is a linear operator that sends 
the function $\eta(x)$ to a real 
or a complex number, while $R_1$ represents a reminder term. 
It is clear that $L([\theta(x)])$ involves 
integration with respect to $x$, since $L([\theta])\eta(x)$ is a real or complex number. Thus, we look for a representation of $L([\theta])\eta(x)$ in the form
\begin{equation}
 L([\theta(x)])\eta(x)=\int \frac{\delta F([\theta])}{\delta \theta(x)}\eta(x) dx,
\label{first_order_L}
 \end{equation}
where the kernel function $\delta F([\theta])/\delta \theta(x)$ 
is a functional of $\theta$ and a function of $x$. 
At this point it is convenient to establish 
a parallel between functionals and functions in $m$ variables. 
Recall that the differential of a scalar field $f(a_1,..,a_m)$ 
in the direction $\widehat{n}=(n_1,...,n_m)$  is the scalar product of the 
gradient $\nabla f$ and $\widehat{n}$, i.e., 
\begin{equation}
df_{\widehat{n}}=\nabla f\cdot \widehat{n}.
\end{equation}
By analogy, $\delta F([\theta])/\delta \theta(x)$ in equation \eqref{first_order_L} can be considered as an infinite-dimensional gradient, known as {\em first-order functional derivative of $F$ with respect 
to $\theta(x)$}. 

\begin{figure}[t]
\centerline{\includegraphics[height=4cm]{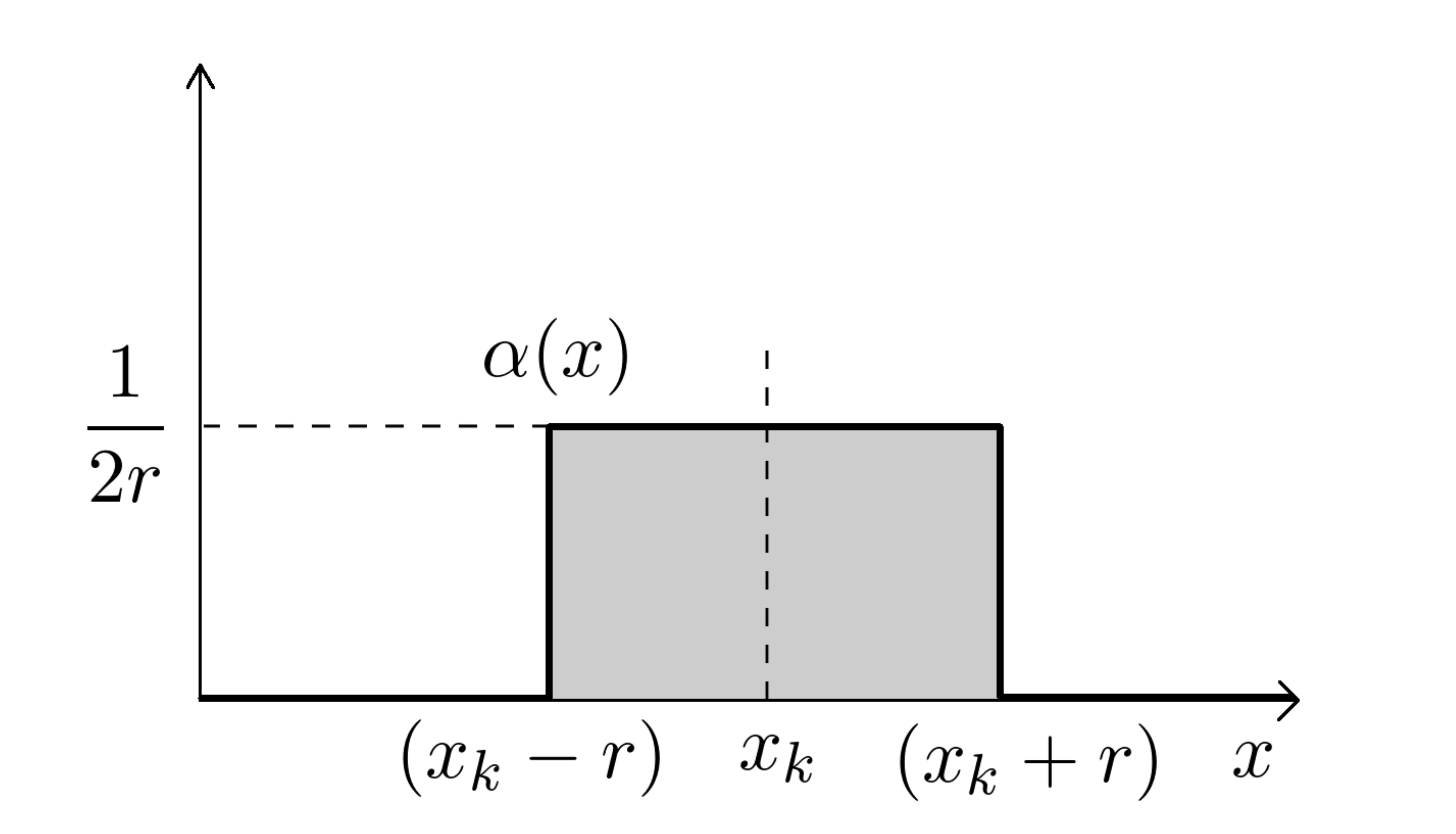}}
\caption{A possible function for the calculation of the kernel
\eqref{A} at $x_k$. If we set
$ \mathcal{I}_{x_k}(r)=\{x\in\mathbb{R}: |x-x_k|\leq r\}$, then 
we have $\displaystyle\int_{\mathcal{I}_{x_k}(r)}\alpha(x)dx=1$, 
independently on $r$.}
\label{fig:3}
\end{figure}

The next question is: how do we compute such functional derivative? 
A possible way is to use the definition \eqref{gateaux} and 
a compactly supported class of test functions, e.g., functions that 
are nonzero only in a small neighbor $\mathcal{I}_x(r)$ 
of radius $r$ centered at $x$.
Such functions could be compactly supported 
elements of a Dirac delta sequence (see Figure \ref{fig:3}), 
or even a delta function itself. This allows us to write
\begin{equation}
\frac{\delta F([\theta(x)])}{\delta \theta(x)}=\lim_{\epsilon,r\rightarrow 0}
\frac{F([\theta(x)+\epsilon \alpha(x)])-F([\theta(x)])}
{\displaystyle\epsilon\int_{\mathcal{I}_x(r)}\alpha(y)dy}
\label{A}
\end{equation}
In particular, if we set $\alpha(y)=\delta(x-y)$ then the denominator 
in \eqref{A} simply reduces to $\epsilon$, yielding the formula 
\begin{align}
\frac{\delta F([\theta])}{\delta \theta(x)}=&\lim_{\epsilon\rightarrow 0}
\frac{F([\theta(y)+\epsilon \delta(x-y)])-F([\theta(y)])}
{\epsilon}.
\label{A2}
\end{align}
Functional derivatives of higher order can be defined 
in a similar manner. For example, the second order functional 
derivative of $F$ is 
\begin{equation}
\frac{\delta^2 F([\theta])}{\delta \theta(x)\delta \theta(y)}=
\frac{d}{d\epsilon}\left[\frac{\delta F([\theta(z)+\epsilon \delta(z-y)])}
{\delta\theta(x)}\right]_{\epsilon=0}.
\label{2ND}
\end{equation}
Note that \eqref{2ND} is a function and $x$ and $y$ and a 
functional of $\theta(x)$.

\vs 
\noindent
{\em Example 1:} The first-order functional 
derivative of the nonlinear functional \eqref{functional1} can 
be obtained as follows. We first compute the G\^ateaux differential 
\begin{align}
\frac{d}{d\epsilon}\left. F([\theta+\epsilon \eta])\right|_{\epsilon=0}=&
-\int_{0}^1 x^{3}e^{-\theta(x)}\eta(x)dx\nonumber\\
=&\int_{0}^1 \frac{\delta F([\theta])}{\delta \theta(x)}\eta(x)dx.
\end{align} 
Therefore,
\begin{equation}
\frac{\delta F([\theta])}{\delta \theta(x)} = -x^{3}e^{-\theta(x)}.
\end{equation} 
Note that the functional derivative is a function of $x$ and a (local) 
functional of $\theta(x)$.

\vs 
\noindent
{\em Example  2  (Functional Derivatives of the Hopf Functional):} 
Consider the Hopf characteristic functional of a random function 
$u(x;\omega)$ defined in $[0,1]$
\begin{equation}
\Phi([\theta]) = \left<\exp\left[i\int_{0}^1 u(x;\omega)\theta(x)dx\right]\right>.
\label{Hopf76}
\end{equation}
The average operator $\left<\cdot \right>$ here 
denotes a functional integral over the probability 
functional of $u(x,\omega)$. The G\^ateaux differential 
of $\Phi([\theta])$ along $\eta$ is 
\begin{align}
\left.\frac{d} {d\epsilon}\Phi([\theta(x)+\epsilon\eta(x)])\right|_{\epsilon=0} 
=&i\left<\exp\left[i\int_0^1 u(x;\omega)\theta(x)dx\right]
\int_0^1u(x;\omega)\eta(x)dx\right>\nonumber \\
=&\int_0^1 i\left<\exp\left[i\int_0^1 u(x;\omega)\theta(x)dx\right]
u(x;\omega)\right> \eta(x)dx\nonumber\\
=& \int_0^1 \frac{\delta \Phi([\theta])}{\delta \theta(x)} \eta(x)dx.
\label{first-order_FD}
\end{align}
This implies that  
\begin{equation}
\frac{\delta \Phi([\theta])}{\delta \theta(x)}=i\left<u(x;\omega)
\exp\left[i\int_a^b u(x;\omega)\theta(x)dx\right]\right>
\label{kerneldef}
\end{equation}
is the first-order functional derivative of $\Phi([\theta])$ at $\theta(x)$. 
Note that \eqref{kerneldef} is itself a functional of $\theta(x)$, which depends also on $x$. As a result, 
$\delta\Phi([\theta])/\delta \theta(x)$ has two types of derivatives: an ordinary one 
with respect to $x$, and a functional one with respect to $\theta(x)$. 
The latter is the the second-order functional derivative 
of $\Phi([\theta(x)])$. A simple calculation 
shows that 
\begin{align}
\frac{\delta^2 \Phi([\theta])}
{\delta \theta(x)\delta \theta(y)}=
i^2\left<u(y;\omega)u(x;\omega)\exp\left[i\int_a^b u(x;\omega)\theta(x)dx\right]\right>.
\label{B}
\end{align}
Proceeding similarly, we can obtain the expression 
of higher-order functional derivatives. For instance, the 
third-order one is explicitly given as
\begin{align}
\frac{\delta^3 \Phi([\theta])}
{\delta \theta(x)\delta \theta(y)\delta\theta(z)}=
i^3\left<\exp\left[i\int_a^b u(x;\omega)\theta(x)dx\right]
u(y;\omega)u(x;\omega)u(z;\omega)\right>.
\end{align}
Now, suppose we have available $\delta \Phi([\theta])/\delta\theta(x)$. 
Based on the definition \eqref{kerneldef}, we see that  
\begin{equation}
\left<u(x;\omega)\right>=\frac{1}{i}\frac{\delta \Phi([0])}{\delta\theta(x)}.
\end{equation}
Similarly, higher order moments and cumulants of the 
random function $u(x;\omega)$ can be obtained 
by computing higher order functional derivatives of $\Phi([\theta(x)])$ and 
$\ln \Phi([\theta(x)])$, respectively, and evaluating them at $\theta(x)=0$.
In particular, the second- and third-order correlation functions are, 
respectively 
\begin{equation}
\left<u(x;\omega)u(y;\omega)\right>=\frac{1}{i^2}
\frac{\delta^2 \Phi([0])}{\delta\theta(x)\delta\theta(y)},
\quad 
\left<u(x;\omega)u(y;\omega)u(z;\omega)\right>=\frac{1}{i^3}
\frac{\delta^3 \Phi([0])}{\delta\theta(x)\delta\theta(y)\delta\theta(z)}.  
\end{equation}

{\color{r}
\vs
\noindent
{\em Example 3:} The Gateaux differential of the 
nonlinear functional \eqref{ex4fun} in the direction $\eta(x)$ is
\begin{align}
\left.\frac{d}{d\epsilon} 
F([\theta+\epsilon \eta])\right|_{\epsilon=0} = &
\frac{d}{d\epsilon}\left( 
\theta(0)+\epsilon \eta(0) +\int_{-1}^1 
\sin(\theta(x)+\epsilon \eta(x))dx
\right)_{\epsilon=0}\nonumber\\
&= \int_{-1}^1 \left(\delta(x)+\cos(\theta(x))\right)\eta(x)dx
\end{align}
Therefore the first-order functional derivative is 
\begin{equation}
\frac{\delta F([\theta])}{\delta \theta(x)}=\delta(x)+\cos(\theta(x)),
\end{equation}
where $\delta(x)$ at the right hand side is the Dirac delta 
function \cite{Kanwal}. 
}

\paragraph{Regularity of Functional Derivatives} 
The last example clearly shows that functional derivatives 
of nonlinear functionals can easily be {\em distributions}
\cite{Kanwal}, e.g., Dirac delta functions. 
For example, let $h\in C^{(\infty)}(\mathbb{R})$. 
Then any functional in the form $F([\theta])=h((\theta,\theta))$, 
where $(,)$ is an inner product in $C^{(\infty)}(\mathbb{R})$, 
has a singular second-order functional derivative. In fact, 
\begin{equation}
\frac{\delta F([\theta])}{\delta \theta(x)}=2\left.
\frac{\partial h}{\partial a}\right|_{a={(\theta,\theta)}}\theta(x)
\end{equation}
\begin{equation}
\frac{\delta^2 F([\theta])}{\delta \theta(x)\delta \theta(y)}=
2\delta(x-y)\left.\frac{\partial h}{\partial a}\right|_{a={(\theta,\theta)}}+
4\left.\frac{\partial^2 h}{\partial a^2}\right|_{a={(\theta,\theta)}}\theta(x)\theta(y).
\end{equation}
The characteristic functional of zero-mean Gaussian white noise 
\begin{equation}
 \Phi([\theta])=e^{-(\theta,\theta)/2}
\end{equation}
belongs to this class, i.e., it has a ``singular'' second-order functional 
derivative. 
On the other hand, the functional 
\begin{equation}
 G([\theta])=h\left( (K,\theta)\right),
\end{equation}
where $K(x)$ is a given smooth kernel, has smooth functional 
derivatives 
\begin{equation}
\frac{\delta G([\theta])}{\delta \theta(x)}=\left.
\frac{\partial h}{\partial a}\right|_{a={(K,\theta)}}K(x),\qquad 
\frac{\delta^2 G([\theta])}{\delta \theta(x)\delta \theta(y)}=
\left.\frac{\partial^2 h}{\partial a^2}\right|_{a={(K,\theta)}}K(x)K(y).
\end{equation}

\section{Approximation of Nonlinear Functionals}
\label{sec:rep-Hopf}

Approximation theory for nonlinear functionals is 
strongly related to approximation theory of nonlinear 
operators \cite{Howlett,Torokhti,Bensoussan_DaPrato}. 
A nonlinear functional $F$ is in fact 
a particular type of nonlinear operator from a space 
of functions $D(F)$ (the domain of the functional $F$) 
into a vector space, e.g., $\mathbb{R}$ or $\mathbb{C}$. 
Thus, the problem of approsimating nonlinear functionals 
is basically the same as approximating nonlinear operators.
This topic has been studied extensively by different 
scientific communities (see, e.g., \cite{Torokhti,Rugh,Nelles,Shetzen,Galman,
Prenter,Bertuzzi,Makarov,Khlobystov1}) for obvious reasons. 
What does it mean to approximate a nonlinear functional? 
Consider, as an example,  the functional \eqref{functional1}, hereafter rewritten for convenience
\begin{equation}
F([\theta])=\int_0^1 x^3 e^{-\theta(x)}dx, \qquad  D(F)=C^{(0)}([0,1]).
\label{functional11}
\end{equation}
Approximating $F([\theta])$ in this case means that 
we are aiming at constructing a nonlinear 
operator $\widehat{F}([\theta])$ that allows us 
to compute an approximation of {\em all} possible integrals in the 
form \eqref{functional11}, for arbitrary continuous functions 
$\theta\in C^{(0)}([0,1])$. This challenging problem 
includes cases in which $F([\theta])$ admits 
an analytical solution, e.g., $F([x])$ or $F([\sin(x)])$, as well 
as cases where no analytical solution is available, e.g., $F([x^2])$.

Perhaps, the most classical and widely used approach to 
represent nonlinear functionals relies on functional power 
series\footnote{Functional power series have been 
widely used in the turbulence theory to obtain 
moment and cumulant expansions (see \cite{Frisch,Monin2,Rosen_1967}).}. The method was originally developed by Volterra \cite{Volterra}, and it represents  
the counterpart of power series expansions in the theory of 
functions. In practice, the functional of interest is represented in terms 
of a series of integral operators involving increasing powers of the test 
function and kernels that need to be determined. The canonical form 
of the power series expansion is
\begin{equation}
F([\theta])=\sum_{k=0}^\infty P_k([\theta]),\quad \textrm{where}\quad 
P_k([\theta])=\int_{-\infty}^\infty\cdots\int_{-\infty}^\infty 
K_k(x_1,...,x_k)\theta(x_1)\cdots\theta(x_k)dx_1\cdots dx_k.
\label{power_series}
\end{equation}

\noindent
{\em Remark:}
Functional power series are known to have bad approximation 
properties and other issues. For example, they often do not preserve 
important properties of the functional, e.g., positive definiteness 
or normalization in the case of Hopf functionals.

{\color{r}
\subsection{Functional Approximation in Finite-Dimensional Function Spaces}
\label{sec:Finite_Dim_Approx}
The simplest ways to establish a closed functional 
representation is to restrict the domain of the functional 
to a finite-dimensional function space spanned by the 
basis $\{\varphi_1(x),...,\varphi_m(x)\}$, i.e.,  
\begin{equation}
D_m = \textrm{span}\{\varphi_1(x),...,\varphi_m(x)\}.
\label{Dspan}
\end{equation}
In this way, any element in $D_m$ can be 
represented as\footnote{
\color{r}
If $D(F)$ is a space of multivariate 
functions defined on some subset of $V\subseteq \mathbb{R}^d$ 
then \eqref{testf} takes the form 
\begin{equation}
\theta(\bm x)=\sum_{k=1}^m a_k \varphi_k(\bm x), \qquad \bm x\in V.
\end{equation}
More generally, $\theta(\bm x)$ can be represented by 
series expansions based on tensor products, 
or more advanced expansions that rely on HDMR \cite{Rabitz,Li1}
or tensor methods \cite{Hackbusch_book,Kolda} (see also Section 
\ref{sec:HDMR} and Section \eqref{sec:tensor}). 
The latter techniques are recommended when operating 
on test function spaces defined on high-dimensional 
domains $V$.
}
\begin{equation}
\theta_m(x)=\sum_{k=1}^m a_k \varphi_k(x).
\label{testf}
\end{equation}
Possible choices of $\varphi_k(x)$ are:
\begin{enumerate}
\item {\em Lagrange Characteristic Polynomials}.
Given a set of $m$ distinct interpolation nodes $\{x_j\}$ in the interval 
$[a,b]$, for example Gauss-Chebyshev-Lobatto nodes, we set
\begin{equation}
 \varphi_j(x)=
 \prod_{\substack{i=1\\i\neq j}}^{m}\frac{(x-x_i)}{(x_i-x_j)}.
\end{equation}
\item {\em Jacobi Polynomials.}  
The function space $D_m\subseteq D(F)$ can be also represented by 
a finite set of Jacobi polynomials $J^{(\alpha,\beta)}_j(x)$ (see 
\cite{Gautshi,Hesthaven}), i.e.
\begin{equation}
 \varphi_j(x)=J^{(\alpha,\beta)}_j(x).
\end{equation}
As is well known, Jacobi polynomials include 
many other families of widely used polynomials such as Gegenbauer, 
Legendre and Chebyshev.
\item {\em Trigonometric Polynomials}. If $D(F)$ is the space of
periodic functions in $[0,2\pi]$, then a convenient choice 
for $\varphi_k(x)$ may be the set of (nodal) trigonometric 
polynomials \cite{Hesthaven}
\begin{equation}
\varphi_j(x)=\frac{1}{m}\sin\left(m\frac{x-x_j}{2}\right)
\cot\left(\frac{x-x_j}{2}\right) \qquad x_j=\frac{2\pi}{m}j\qquad j=0,...,m
\label{nodalpoly}
\end{equation}
or, equivalently, classical Fourier modes
\begin{equation}
\varphi_0(x) = 1,\qquad \varphi_k(x)=\sin(kx)\qquad \varphi_{m+k}(x)=\cos(kx) 
\qquad k=1,...,m.
\end{equation} 

\end{enumerate}

\noindent
For each specific choice of the basis 
set $\{\varphi_1(x),...,\varphi_m(x)\}$, the test function 
\eqref{testf} lies on a {\em parametric manifold} of 
dimension $m$, i.e., a hyperplane. 
Any discretization of the function space $D(F)$ in terms 
of a finite-dimensional basis, reduces the 
functional $F$ into a multivariate function with 
domain $D_m$ and range $F([D_m])$. 
Such function depends on as many 
variables as the number of degrees of freedom 
we consider in the finite-dimensional 
approximation of $D(F)$.

\vs
\noindent
{\em Example 1:} A substitution of \eqref{testf} 
into the Hopf functional \eqref{Hopf76} yields the 
complex-valued multivariate function
\begin{equation}
\phi(a_1,...,a_m) = \left<
\exp\left(i\sum_{k=1}^m a_j U_j(\omega)\right)
\right>, \qquad
U_j(\omega) = \int_{0}^1 u(x;\omega)\varphi_j(x)dx,
\label{missdef}
\end{equation}
i.e., the joint characteristic function of the Fourier coefficients 
$U_j(\omega)$.
Note that $\phi$ depends on $m$ real variables $(a_1,...,a_m)$. 
Such multivariate function can be seen as a $m$-dimensional 
parametrization of the mapping $\Phi$ shown in Figure \ref{fig:1}, i.e., 
a parametrization of the nonlinear transformation 
$D_m \rightarrow \Phi(D_m)$.
In this setting, approximation of nonlinear functionals is 
equivalent to approximation of a real- or complex-valued 
multivariate functions. The question of whether a functional 
can be approximated by evaluating it in a space spanned by a 
finite-dimensional basis is different, and will be addressed 
in Section \ref{sec:approximability_of_functionals}.

\vs
\noindent
{\em Example 2:} Consider the nonlinear functional \eqref{functional1} 
and let $D_m\subseteq D(F)$ be the function space 
spanned by a suitable set of orthogonal 
polynomials\footnote{Given any positive measure in a one-dimensional 
interval, it is always possible to construct a set of polynomials that 
is orthogonal with respect to such measure \cite{Gautshi1,Gautshi}.} 
in $[0,1]$. Evaluating $F([\theta])$ in $D_m$, i.e., considering test 
functions $\theta_m$ in the form \eqref{testf}, yields 
the multivariate function 
\begin{equation}
f(a_1,....,a_m) = \int_0^1 x^3 \prod_{k=1}^m e^{-a_k\varphi_k(x)}dx.
\end{equation}
This is the exact form of the functional $F$, evaluate in $D_m$, i.e., 
\begin{equation}
F([\theta_m])=f(a_1,....,a_m).
\end{equation}

\subsubsection{Functional Derivatives} 
Evaluating a nonlinear functional $F([\theta])$ 
in a finite-dimensional function space 
$D_m$ allows for a simple and effective representation of 
functional derivatives. In particular, it can be shown that
\begin{equation}
\left.\frac{\delta F([\theta])}{\delta \theta(x)}
\right|_{\theta\in D_m}=\sum_{k=1}^m \varphi_k(x)\frac{\partial f}{\partial a_k},
\label{ha1}
\end{equation}
\begin{equation}
 \left.\frac{\delta^2 F([\theta])}{\delta \theta(x)\delta\theta(y)}
\right|_{\theta\in D_m}=\sum_{j,k=1}^m \varphi_k(x)\varphi_j(y)
\frac{\partial f}{\partial a_j\partial a_k}.
\label{ha2}
\end{equation}
Here, $f(a_1,...,a_m)=F([\theta_m])$ is the function we obtain 
by evaluating the functional $F$ in the finite dimensional function 
space $D_m$. The meaning of \eqref{ha1} and \eqref{ha2} is the 
following: if we evaluate the functional derivatives of 
$F$ in the finite dimensional space $D_m$ (recall that the functional derivatives are themselves nonlinear functionals) 
then we can represent them in terms of classical partial derivatives 
of $f(a_1,...,a_m)$. Note that the basis function 
$\varphi_j$ spanning $D_m$ also appear in \eqref{ha1}-\eqref{ha2}, 
suggesting that the accuracy of the functional derivatives depend on 
the choice of such basis functions.
Rather than proving \eqref{ha1} and \eqref{ha2} in a general setting, 
let we provide two constructive examples that yield 
expressions in the form \eqref{ha1} and \eqref{ha2}.

\vs
\noindent
{\em Example 1:} Consider the Hopf functional \eqref{Hopf76}. 
By evaluating the analytical expression of the first- and second-order  
functional derivatives \eqref{kerneldef}-\eqref{B} in 
the finite-dimensional function space \eqref{Dspan} we obtain
\begin{align}
\left.\frac{\delta \Phi(\theta)}{\delta \theta(x)}
\right|_{\theta\in D_m}=& i\left<u(x)e^{i(a_1U_1(\omega)+\cdots 
+ a_m U_m(\omega))}\right>,\qquad \\
\left.\frac{\delta^2 \Phi(\theta)}
{\delta \theta(x)\delta \theta(y)}
\right|_{\theta\in D_m}=&-\left<u(x)u(y)
e^{i(a_1U_1(\omega)+\cdots+ a_m U_m(\omega))}\right>,
\end{align}
where $U_j(\omega)$ are random variables 
defined in \eqref{missdef}.
By using the definition of the characteristic 
function \eqref{missdef} we have that
\begin{align}
\frac{\partial \phi}{\partial a_k}=&\int_{0}^1\varphi_k(x)
\left.\frac{\delta \Phi(\theta)}{\delta \theta(x)}
\right|_{\theta\in D_m}dx,\\ 
\frac{\partial^2 \phi}{\partial a_k\partial a_j}=&
\int_0^1\int_{0}^1\varphi_k(x)\varphi_j(y)\left.\frac{\delta^2 \Phi(\theta)}
{\delta \theta(x)\delta \theta(y)}
\right|_{\theta\in D_m}dxdy.
\end{align}
This means that the partial derivatives of the characteristic function 
are nothing but the  projection of the Hopf functional derivatives 
onto the space $D_m$. Clearly, if the random 
function $u(x;\omega)$ is $D_m$ then the following 
inverse formulas hold 
\begin{equation}
 \left.\frac{\delta \Phi([\theta])}{\delta \theta(x)}
\right|_{\theta,u\in D_m}=\sum_{k=1}^m \varphi_k(x)\frac{\partial \phi}{\partial a_k},
\label{ha111}
\end{equation}
\begin{equation}
 \left.\frac{\delta^2 \Phi([\theta])}{\delta \theta(x)\delta\theta(y)}
\right|_{\theta,u\in D_m}=\sum_{j,k=1}^m \varphi_k(x)\varphi_j(y)
\frac{\partial \phi}{\partial a_j\partial a_k}.
\label{ha222}
\end{equation}

 \vs
 \noindent
{\em Example 2:} Consider the sine functional
\begin{equation}
 F([\theta])=\sin\left(\int_a^b K(x)\theta(x)dx\right)
 \label{sinFunctional}
\end{equation}
where $K(x)$ is a given kernel function.
Evaluating $F$ in $D_m$ yields the multivariate function 
\begin{equation}
 f(a_1,...,a_m)=\sin\left(\sum_{i=1}^m a_i \int_a^b K(x)\varphi_i(x)dx\right).
 \label{fdr}
\end{equation}
Similarly, evaluating the functional derivative of $F$ in $D_m$ yields 
\begin{equation}
 \left.\frac{\delta F([\theta])}{\delta \theta(x)}\right|_{\theta\in D_m}=
 \cos\left(\sum_{i=1}^m a_i \int_a^b K(x)\varphi_i(x)dx\right)K(x).
 \label{fde}
\end{equation}
A comparison between \eqref{fdr} and \eqref{fde} immediately yields
\begin{equation}
\frac{\partial f}{\partial a_k}=\int_a^b\left.\frac{\delta F([\theta])}{\delta \theta(x)}\right|_{\theta\in D_N}
\varphi_j(x)dx,
\label{gg1}
\end{equation}
i.e., the gradient of $f$ is the projection of the 
functional derivative of $F([\theta])$ (evaluated in $D_m$) 
onto $D_m$. On the other hand, if $K(x)$ is a function in $D_m$ 
then 
\begin{equation}
\left.\frac{\delta F([\theta]}{\delta \theta(x)}\right|_{\theta\in D_m}= \sum_{k=1}^m
\frac{\partial f}{\partial a_k} \varphi_k(x).
\label{gg2}
\end{equation}
This clarifies the meaning of the functional derivative 
in both finite- and infinite-dimensional ($m\rightarrow \infty$) 
cases.

\subsubsection{Distances between Function Spaces and Approximability of Functionals}
\label{sec:approximability_of_functionals}
A key concept when approximating a nonlinear functional $F([\theta])$ by restricting its domain $D(F)$ to a finite-dimensional space functions 
$D_m$ is the distance between $D_m$ and $D(F)$. Such distance can be quantified in different ways (see, e.g., \cite{Pinkus}). For example we can define the {\em deviation} of $D_m$ from $D(F)$ as
\begin{equation}
E(D_m,D(F)) = \adjustlimits \sup_{\theta\in D(F)}
\inf_{\theta_m\in D_m} \left\|\theta-\theta_m\right\|
\label{functionaldeviation}
\end{equation}
The number $E$ measure the extent to which the worst 
element of $D(F)$ can be approximated from $D_m$. One 
may also ask how well we can approximate $D(F)$ with  
$m$-dimensional subspaces of $D(F)$ which are allowed to vary within $D(F)$. A measure of such approximation is 
given by the Kolmogorov $m$-width 
\begin{equation}
d_m(D_m,D(F))= \adjustlimits \inf_{D_m} \sup_{\theta\in D(F)}
\inf_{\theta_m\in D_m} \left\|\theta-\theta_m\right\|_{D(F)}
\label{knw}
\end{equation}
which quantifies the error of the {\em best approximation} 
to the elements of $D(F)$ by elements in a vector 
subspace $D_m$ of  dimension at most $m$.
The Kolmogorov $m$-width can be rigorously defined, e.g., for 
nonlinear functionals in Hilbert spaces (\cite{Pinkus}, Ch. 4). 
In simpler terms we can define the notion of approximability of 
a nonlinear functional as follows. Let $F([\theta])$ be a 
continuous nonlinear functional with domain $D(F)$, and consider 
a finite-dimensional subspace $D_m\subseteq D(F)$, for example 
$D_m=\textrm{span}\{\varphi_1,...,\varphi_m\}$. 
We say that $F([\theta])$ is approximable in $D_m$ if for all $
\theta\in D(F)$ and $\epsilon>0$, there exists $m$ (depending on 
$\epsilon$) and an element $\theta_m\in D_m$ such that
\begin{equation}
\left\| F([\theta])-F([\theta_m])\right\| \leq \epsilon.
\label{Ferr}
\end{equation}
Clearly if $F$ is continuous and $\theta_m$ is close to $\theta$, 
i.e., the deviation \eqref{functionaldeviation} between $D(F)$ 
and $D_m$  is small, then we expect $\epsilon$ to be small. 
It is important to emphasize that the approximation error and the computational complexity of approximating a nonlinear functional 
depends on the choice of $D_m$. In particular, a functional may be low-dimensional in one function space and high-dimensional in another. The following example clarifies this question.

\vs
\noindent
{\em Example 1:} Consider the sine functional  
\begin{equation}
F([\theta])=\sin\left(\int_{0}^{2\pi}\theta(x)dx\right)
\label{Fsin}
\end{equation}
in the space $D(F)$ of periodic functions in $[0,2\pi]$. 
If we represent $\theta$ in terms of orthonormal Fourier modes, 
i.e., we consider 
\begin{equation}
D_{2m+1}= \textrm{span}\left\{\frac{1}{\sqrt{2\pi}},\frac{\sin(x)}{\sqrt{\pi}},..., \frac{\sin(mx)}{\sqrt{\pi}},\frac{\cos(x)}{\sqrt{\pi}},...,\frac{\cos(mx)}{\sqrt{\pi}}\right\}
\label{D2m}
\end{equation}
and 
\begin{equation}
\theta_{2m+1}(x)=\frac{a_0}{\sqrt{2\pi}} + \frac{1}{\sqrt{\pi}}
\sum_{j=1}^m a_j \sin(jx)+  \frac{1}{\sqrt{\pi}}
\sum_{j=1}^m b_j \cos(jx)
\end{equation}
In this setting, we obtain
\begin{equation}
F([\theta_{2m+1}])=\sin(\sqrt{2\pi} a_0).
\end{equation}
This means that \eqref{Fsin} is approximable in the function space 
\eqref{D2m}. Moreover, the approximation is 
{\em exact} and just one-dimensional.  
On the other hand, if we set the space 
$D_{2m+1}$ to be the span of a normalized nodal Fourier basis 
$\{\varphi_0, ...,\varphi_{2m}\}$, e.g., the normalized odd 
expansion discussed in \cite{Hesthaven}, then the functional 
\eqref{Fsin} technically requires an infinite number of variables 
variables. In fact, in this case we have
\begin{equation}
F([\theta_{2m+1}])=\sin\left(\eta \sum_{k=0}^{2m} a_k\right), 
\quad \textrm{where}\quad \eta = \int_{0}^{2\pi}\varphi_k(x)dx =
\left(\frac{2\pi}{2m+1}\right)^{\frac{1}{2}}.
\end{equation}

\vs
\noindent
{\em Example 2:}
Consider the characteristic functional of zero-mean 
Gaussian white noise (see equation \eqref{white_0}),
\begin{equation}
\Phi([\theta])=\exp\left[-\frac{1}{2}\int_0^{2\pi} \theta(x)^2dx\right],
\label{white}
\end{equation}
where $D(\Phi)$ is the space of periodic functions 
in $[0,2\pi]$. Let $D_m$ be the space spanned by 
any finite orthonormal set of periodic functions. 
The deviation between $D(\Phi)$ and $D_m$ 
this case yields a functional approximation 
error of order 1. 
To show this in a simple way, evaluate the functional 
\eqref{white} in both $D(\Phi)$ and $D_m$. This yields 
\begin{equation}
 \Phi([\theta])=\exp\left[-\frac{1}{2}\sum_{k=1}^\infty a_k^2\right],\qquad
 \Phi([\theta_m])=\exp\left[-\frac{1}{2}\sum_{k=1}^m a_k^2\right].
\end{equation}
If we measure the error between $\Phi([\theta])$ 
and $\Phi([\theta_m])$ in the uniform operator 
norm then we have 
\begin{equation}
\left\|\Phi([\theta_m])-\Phi([\theta])\right\|_{\infty} = 1,
\end{equation}
independently on $m$. In other words, \eqref{white} 
is not approximable in any finite-dimensional subset of $D(\Phi)$. 
This result is consistent with white-noise 
theory \cite{Stratonovich}. Recall, in fact, that a 
delta-correlated Gaussian process has a flat Fourier 
power spectrum. This implies that any finite
truncation of the Fourier series of such process 
yields a systematic error that is not small.

\vs
\noindent
{\em Remark:} In some cases the effects of 
the distance between $D(F)$ and $D_m$ can be 
mitigated by the presence of smooth functions 
appearing the in functional $F$. For example, consider 
the sine functional
\begin{equation}
F([\theta])=\exp\left(\int_{-1}^{1} \sin(x) \theta(x)dx\right)
\end{equation}
and let $D(F)$ be the space of infinitely differentiable functions
 in $[-1,1]$. If we expand $\theta$ in terms of Legendre 
 polynomials $\{\varphi_k\}$, i.e., 
 \begin{equation}
 \theta_m(x)=\sum_{k=0}^m a_k \varphi_k(x)
 \end{equation}
 then,
 \begin{equation}
 \int_{-1}^{1} \sin(x) \theta(x) dx = \sum_{k=0}^m a_k
 \int_{-1}^1\sin(x) \varphi_k(x)dx.
 \end{equation}
As is well known, the coefficients  $\int_{-1}^1\sin(x) \varphi_k(x)dx$ 
decay to zero exponentially fast with $m$ \cite{Hesthaven}. 
This implies that convergence of $F([\theta_m])$ to $F([\theta])$ 
is exponentially fast in the number of dimensions $m$, that is the 
error \eqref{Ferr} goes to zero exponentially fast with $m$. 

\vs
\noindent
{\em Example 3:}
Consider the Hopf characteristic functional of a zero-mean 
correlated Gaussian process in $[0,2\pi]$ 
\begin{equation}
\Phi([\theta])=\exp\left[-\frac{1}{2}\int_0^{2\pi}\int_{0}^{2\pi} C(x,y)\theta(x)\theta(y)dxdy\right],
\label{hopfc}
\end{equation}
where $C(x,y)$ is a smooth covariance function, and 
$D(\Phi)$ is the space of periodic functions in $[0,2\pi]$. 
If the projection of $C(x,y)$ onto the span of 
an orthonormal set  $D_m\subseteq D(F)$ 
decays with $m$, then the functional $\Phi$ is 
approximable in $D_m$. Note that this is indeed 
the case if the covariance is smooth (and periodic) 
and $\varphi_j(x)$ are the Fourier modes in \eqref{D2m}. 
The smoother the covariance the smaller the number 
of Fourier modes we need to achieve a certain
accuracy \cite{Hesthaven}, i.e., 
the smaller the number of dimensions. 
If $C(x,y)=1$ then the Hopf functional \eqref{hopfc} 
is effectively {\em  one-dimensional}.

}

\subsection{Functional Interpolation Methods}
\label{sec:Functional Collocation Methods}
 
In this Section we discuss how to construct an approximation 
of a nonlinear functional $F([\theta])$ in terms of 
a {\em functional interpolant} $\Pi([\theta])$, i.e., a functional that interpolates $F([\theta])$ at at a given set 
of nodes $\{\theta_1(x),...,\theta_m(x)\}\in D(F)$ 
\begin{equation}
F([\theta_j])=\Pi([\theta_j]) \qquad j=1,...,m.
\end{equation}
Differently from interpolation methods in spaces of 
finite dimension (e.g., $d$-dimensional Euclidean spaces), 
interpolation here is in a space of functions, i.e., 
the {interpolation nodes} $\theta_k(x)$ are 
functions in a Hilbert or a Banach space. 
Over the years, the problem of constructing 
a functional interpolant through suitable nodes 
in Hilbert or Banach spaces has been studied by several
authors and convergence results were established in rather 
general cases
\cite{Makarov,Khlobystov0,Prenter,Porter,Kaplitskii,Allasia,
Khlobystov,Khlobystov1,PorterSIAM,Torokhti}. 

Before discussing functional interpolation in detail, 
let us provide some geometric intuition on what functional 
interpolation is and what kind of representations 
we should expect.
To this end, let us first recall recall that 
a hyperplane in a $d$-dimensional space is a linear 
manifold defined uniquely by $d$ interpolation nodes, 
each node being a vector of $\mathbb{R}^{d}$. 
If we send $d$ to infinity then the hyperplane intuitively 
becomes a linear functional, which is therefore defined uniquely 
by an infinite number of $\infty$-dimensional nodes, 
i.e., an {\em infinite number of functions}.
This suggests that if we consider any finite number of  
nodes in a function 
space, say $\{\theta_1(x),...,\theta_m(x)\}$, then 
we cannot even represent {\em linear functionals} in 
an exact way\footnote{We recall that 
the variational form of nonlinear PDEs is defined by 
linear functionals on test function spaces. 
In this setting, classical Galerkin methods 
to solve PDEs (see Section \ref{sec:variational}) are basically 
identification problems for linear functionals.}, 
i.e., functionals in the form 
\begin{equation}
P_1([\theta])=\int_a^bK_1(x_1)\theta(x_1)dx_1,
\end{equation}
where $K_1(x)$ is a given kernel. 
The same conclusion obviously holds 
for nonlinear functionals, with the aggravating factor 
that the number of test functions theoretically required 
for the exact representation grows significantly. 
For example, quadratic and cubic forms in $d$-dimensions are identified by  
by $d^2$ and $d^3$ interpolation nodes, respectively, where each node is 
vector of $\mathbb{R}^{d}$. When we send $d$ to infinity, we intuitively obtain
homogeneous polynomial functionals of second- and third-order, respectively.
These functionals are in the form 
\begin{align}
 P_2([\theta])&= \int_a^b\int_a^b K_2(x_1,x_2)\theta(x_1)\theta(x_2)dx_1dx_2,
 \label{sof}\\
 P_3([\theta])&= \int_a^b\int_a^b\int_a^b K_3(x_1,x_2,x_3)
              \theta(x_1)\theta(x_2)\theta(x_3)dx_1dx_2dx_3.
\end{align}
Thus, to represent $P_2$ and $P_3$ exactly we need $\infty^2$ 
and $\infty^3$ nodes in a function space. Why $\infty^2$ and $\infty^3$? 
Consider $P_2([\theta])$ and assume that 
the kernel function $K_2(x_1,x_2)$ is in a separable 
Hilbert space. Represent $K_2$ relative to any complete 
orthonormal basis $\{\varphi_k\}_{k=1,...,\infty}$ 
\begin{equation}
K_2(x_1,x_2)=\sum_{i,j=1}^\infty a_{ij}\varphi_i(x_1)\varphi_j(x_2). 
\label{K2e}
\end{equation}
A substitution of \eqref{K2e} into \eqref{sof} yields
\begin{equation}
P_2([\theta])= \sum_{i,j=1}^\infty a_{ij}\int_{a}^b\varphi_i(x)\theta(x)dx 
\int_a^b \varphi_j(x)\theta(x)dx.
\label{ssof}
\end{equation}
Without loss of generality we can assume that $K_2$ is symmetric, i.e., that $a_{ij}$
is a symmetric matrix. To represent $P_2([\theta])$ exactly by means of a functional 
interpolant we need enough nodes $\{\theta_p(x)\}$ to determine each 
$a_{ij}$ in \eqref{K2e} uniquely. Clearly, the choice
$\theta_p(x)=\varphi_p(x)$ ($p=1$,..., $\infty$) is not sufficient for 
this purpose, since it allows us to determine only the diagonal 
entries $a_{pp}$. Therefore we need to construct a larger set of collocation 
nodes, e.g., the set $\theta_p(x)=\varphi_{p_i}(x)+\varphi_{p_j}(x)$, where
$p_i \geq p_j$ and $p_i\in \{1,...,\infty\}$. 
This is the $\infty^2$ number of functions we have 
mentioned above. Similarly, to identify 
$P_3([\theta])$ through functional interpolation we need
$\infty^3$ nodes, e.g., in the form $\theta_p(x)=\varphi_{p_i}(x)
+\varphi_{p_j}(x)+\varphi_{p_k}(x)$, where $p_i\geq p_j\geq p_k$ 
and $p_k\in \{1,...,\infty\}$.
As we shall see later in this Section, determining a polynomial interpolant 
of an unknown functional $F([\theta])$, i.e., determining the kernels 
$K_j(x_1,...,x_j)$ in \eqref{power_series} from input-output relations 
is a {\em linear problem} that involves high-dimensional systems and 
big data.

\subsubsection{Interpolation Nodes in Function Spaces}
\label{sec:interpolationnodes}
Let $F([\theta])$ be a continuous functional  
with domain $D(F)$. Within $D(F)$ we define the spaces 
of functions
\begin{equation}
 S^{(m)}_q=\left\{\theta(x)\in D(F)\,\,|\,\, \theta(x)=a_{i_1}\varphi_{i_1}(x)+\cdots + 
 a_{i_q}\varphi_{i_q}(x)\right\}.
 \label{SNq}
\end{equation}
where $i_j\in\{1,...,m\}$, $a_{i_j}\in \mathbb{R}$ 
and $\{\varphi_{1}(x),...,\varphi_{m}(x)\}\in D(F)$. The elements of $ S^{(m)}_1$, 
$S^{(m)}_2$ and $S^{(m)}_3$ are in the form 
\begin{equation}
 \begin{array}{cl}
 S^{(m)}_1:& \theta(x) = a_{i_1}\varphi_{i_1}(x),\nonumber\\
 S^{(m)}_2:& \theta(x) = a_{i_1}\varphi_{i_1}(x)+a_{i_2}\varphi_{i_2}(x),\nonumber\\
 S^{(m)}_3:& \theta(x) = a_{i_1}\varphi_{i_1}(x)+a_{i_2}\varphi_{i_2}(x)+
 a_{i_3}\varphi_{i_3}(x).\nonumber
\end{array}
\label{snq}
\end{equation}
Clearly, if $\theta_1\in S^{(m)}_j$ and $\theta_2\in S^{(m)}_q$, then  
$(\theta_1+\theta_2)\in S^{(m)}_{j+q}$ (if $q+j< m$).
Also, note that the sequence of spaces $S_j^{(m)}$, $j=1, 2, ...$  is 
{\em hierarchical} in the sense that the following
chain of embeddings hold
\begin{equation}
S^{(m)}_1\subset  S^{(m)}_2\subset \cdots \subset S^{(m)}_m\subset D(F).
\label{embedding}
\end{equation}

\begin{figure}[!t]
\centerline{
\includegraphics[height=6cm]{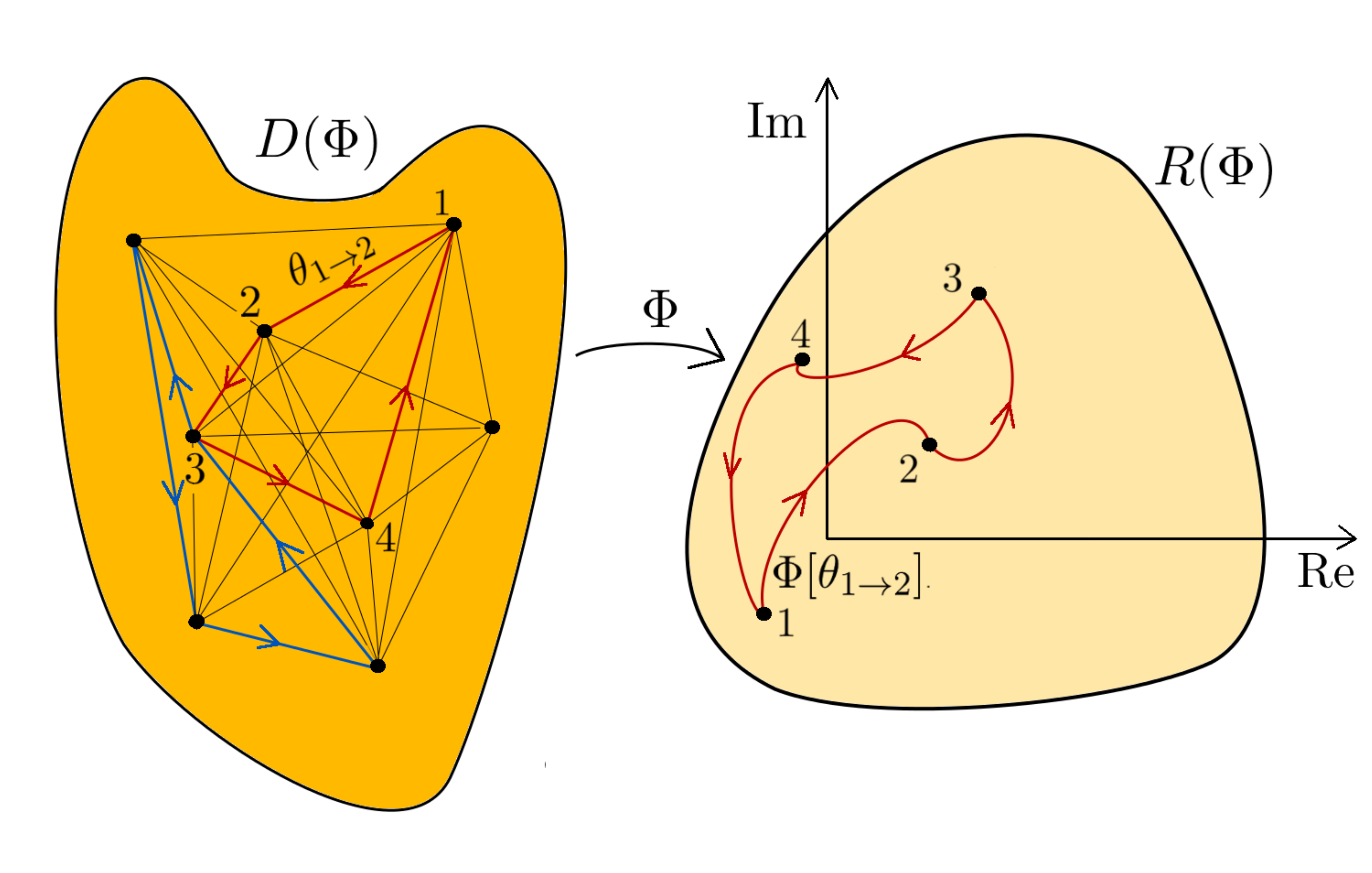}
\includegraphics[height=5cm]{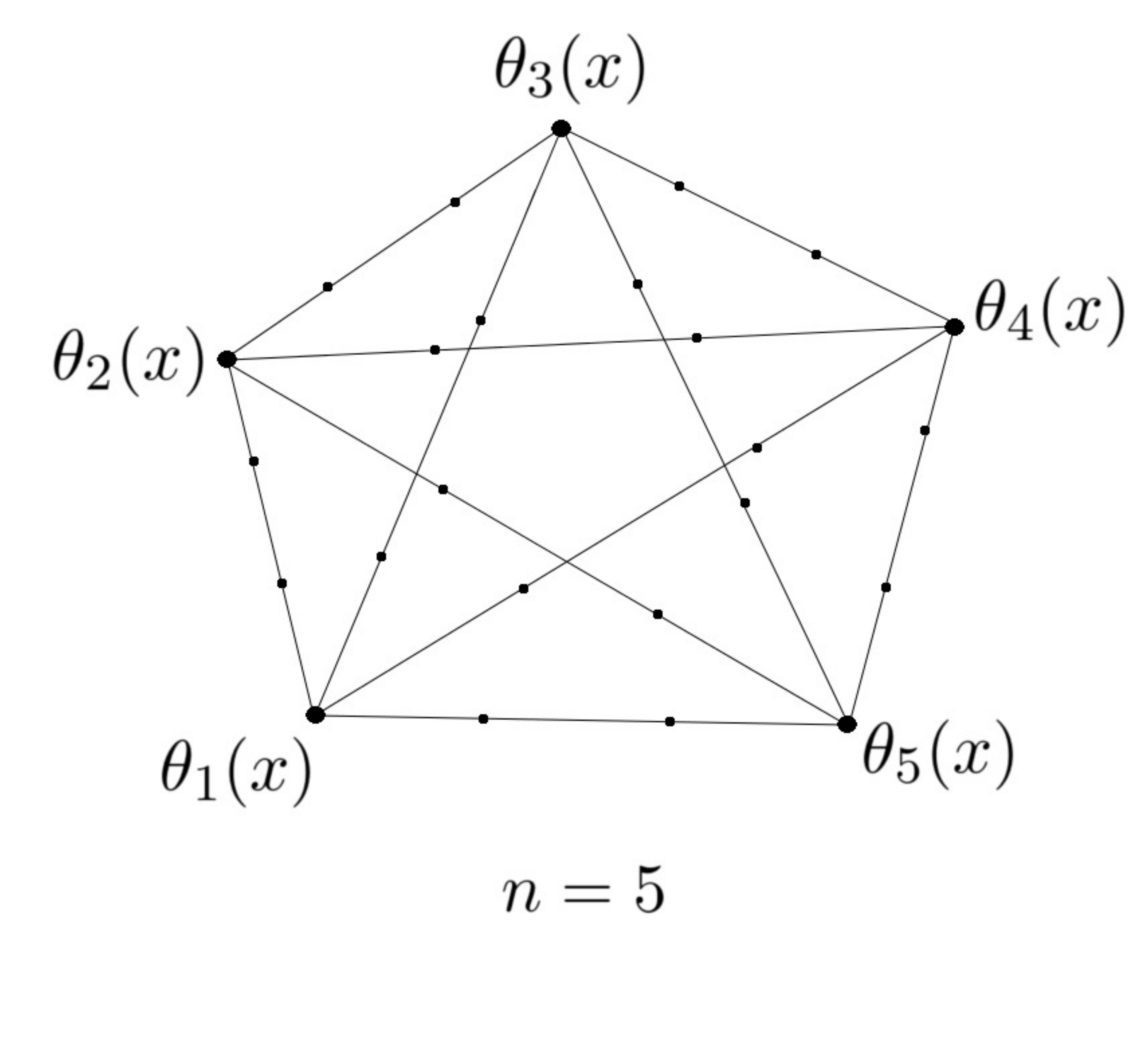}
}
\caption{Network of test functions: 
A closed path in the function space $D(\Phi)$ made of four lines 
(e.g., the red polygon) is mapped through the functional $\Phi$ 
into a closed curve in the complex plane. Vertex and edges of the 
network are elements of $S_2^{(m)}$. We also show a complete 
graph with $n=5$ nodes, where each edge is discretized 
with $p=4$ points.}
\label{fig:8}
\end{figure}
\begin{figure}[!t]
\centerline{
\includegraphics[height=6cm]{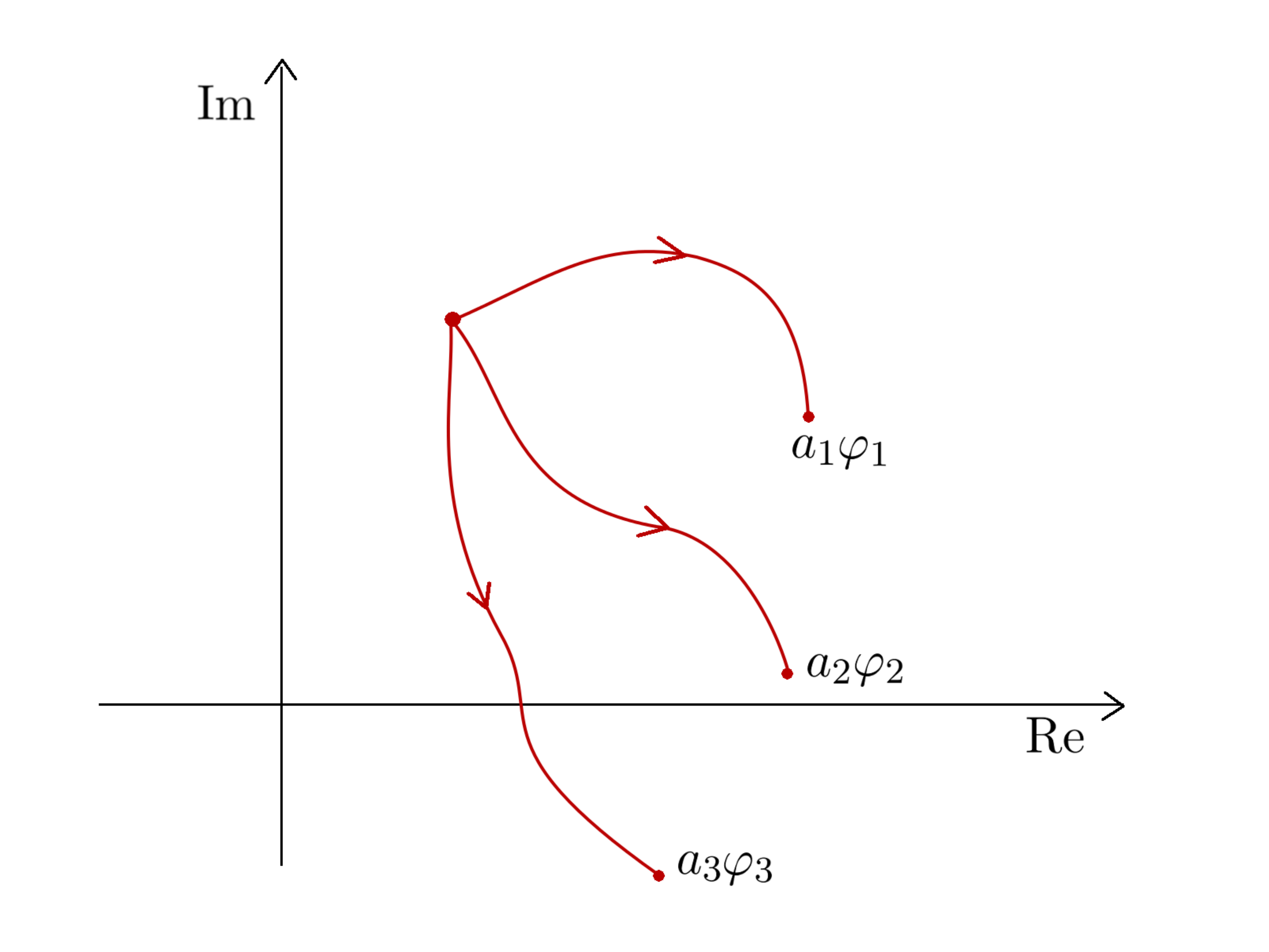}
\includegraphics[height=5.5cm]{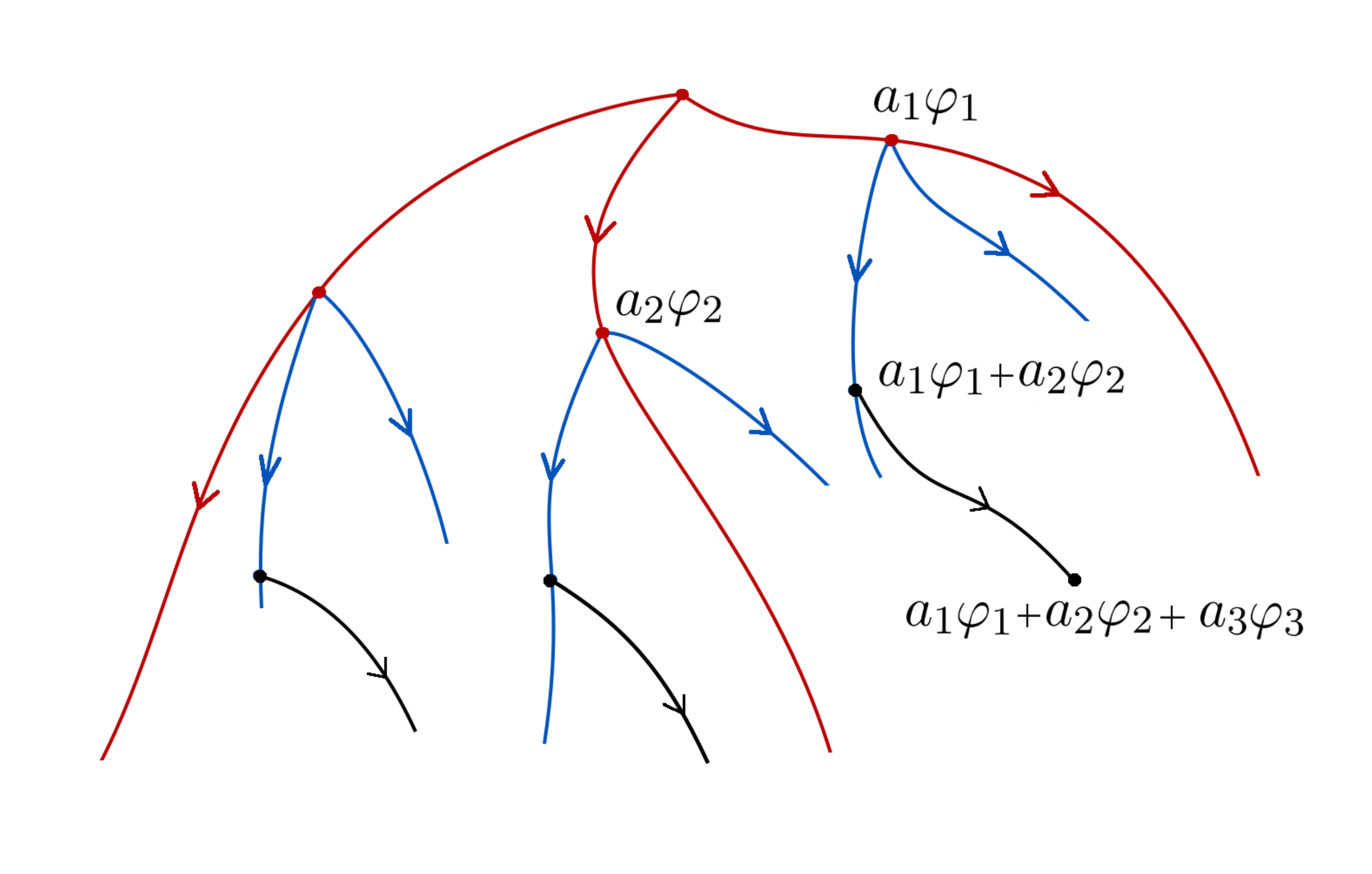}
}
\caption{Trajectories in the range of  a complex valued 
functional $\Phi$ corresponding 
to functions in the sets $S^{(m)}_1$ (left) and $S^{(m)}_3$ (right).
In the case of $S^{(m)}_1$ each curve is parametrized 
by only one parameter $a_j$ and it cannot branch into two distinct curves. 
On the contrary, if we consider $S_3^{(m)}$ we are adding 
two more degrees of freedom and each curve can branch at most twice.
}
\label{fig:branching}
\end{figure}

The function space $S_q^{(n)}$ admits a simple yet powerful graphical 
representation in terms of {\em trajectories of functions} 
\cite{Magri,Tonti,Tonti1,Tonti2,Daniele_JMathPhys}. To illustrate 
such representation,  consider a complex-valued functional $\Phi$. 
A trajectory of functions in the space $D(\Phi)$ is a curve 
in $D(\Phi)$, which is mapped to a curve 
in $\mathbb{C}$ (see Figure \ref{fig:8}). 
Furthermore, if the functional is continuous and differentiable, 
a smooth curve in $D(F)$ is mapped onto a smooth curve in $\mathbb{C}$.
The set of trajectories in the complex plane associated with $S^{(m)}_1$ 
is shown in Figure  \ref{fig:branching}. 
Each curve is parametrized by only one parameter $a_j$ and it 
cannot branch into two distinct curves . 
On the contrary, if we consider $S_2^{(m)}$ we are adding 
one more degree of freedom and each curve departing from $1$ can branch, 
but only once. Similarly, we can have three branches in $S_3^{(m)}$, etc.
A remarkable distribution of nodes in $S_2^{(m)}$ is 
associated with {\em networks of test functions}, i.e., {\em graphs} 
in the function space $D(F)$. Vertex and edges are elements 
of $S_2^{(m)}$. A simple example is shown in Figure \ref{fig:8}. 
In mathematics such network is called {\em complete graph}, i.e., 
an undirected graph in which every pair of distinct nodes is 
connected by a unique edge -- the edge being the trajectory 
(straight line) of functions connecting $\theta_i$ to $\theta_j$. 
If we discretize each edge with $p$ 
collocation points (including the endpoints) and we have $m$ nodes then 
the number of degrees of freedom is $m(m-1)(p-2)/2+m$.

Another interesting set of interpolation nodes in $D(F)$ 
is the one obtained by setting all 
coefficients $a_{i_p}$ in \eqref{SNq} 
equal to $1$, i.e.,
\begin{equation}
 \widehat{S}^{(m)}_q=\left\{\theta_i(x)\in D(F)\,\,|\,\, \theta_i(x)=\varphi_{i_1}(x)+\cdots + 
 \varphi_{i_q}(x)\right\},\qquad i_s=0,...,m,
 \label{SNq1}
\end{equation}
If the set $\{\varphi_j\}$ includes the null element $\{0\}$, 
then by symmetry \eqref{SNq1} is equivalent to
\begin{equation}
 \widehat{S}^{(m)}_q=\left\{\{0\},\{\varphi_{1},..,\varphi_{m}\},
 \{2\varphi_{1},..,2\varphi_{m}\}, 
 \{(\varphi_{1}+\varphi_{2}),...,(\varphi_{1}+\varphi_{m})\},
 \{(\varphi_{2}+\varphi_{3}),...,(\varphi_{2}+\varphi_{m})\},...\right\}.
 \label{SNq2}
\end{equation}
In this case, the number of elements (cardinality) of $\widehat{S}^{(m)}_q$ is
\begin{equation}
\# \widehat{S}^{(m)}_q=\sum_{j=0}^q\binom{j+m-1}{j},\qquad \textrm{where} 
\,\,\, \binom{i}{j}\,\,\,
\textrm{is the binomial coefficient.}
\label{SCardinality}
\end{equation} 
For example, 
\begin{equation}
\# \widehat{S}^{(10)}_1=11,\qquad \#\widehat{S}^{(10)}_2=66\qquad
\#\widehat{S}^{(10)}_3=286,\qquad \#\widehat{S}^{(10)}_5=3003,\qquad
\#\widehat{S}^{(10)}_{10}=184756.
\label{cardS}
\end{equation}
In addition, the cardinality of $\widehat{S}^{(m)}_q$ satisfies the recursion relation
\begin{equation}
 \# \widehat{S}^{(m)}_q= \# \widehat{S}^{(m)}_{q-1}+\binom{q+m-1}{q}.
\end{equation}
The set of functions $\widehat{S}^{(m)}_q$ is sufficient 
to uniquely identify a polynomial functional of order $q$ in which each 
kernel function is represented relative to tensor product basis with $m$ 
elements in each variable.
However, $\widehat{S}^{(m)}_q$ is, in general, {\em not} sufficient 
to accurately interpolate nonlinear functionals. The main problem is that the 
set of nodes \eqref{SNq1} may not be large enough or may not 
cover the function space $D(F)$ appropriately. 
Another open question is related to the selection of 
optimal interpolation nodes in $D(F)$ yielding highly accurate 
representations. This question is addressed in Section 
\ref{sec:optimal interpolation nodes}.

\begin{figure}[t]
\vspace{-0.5cm}
\centerline{\hspace{0.4cm}level 5\hspace{4.8cm} level 7\hspace{4.8cm} level 10}
\centerline{\includegraphics[height=4.5cm]{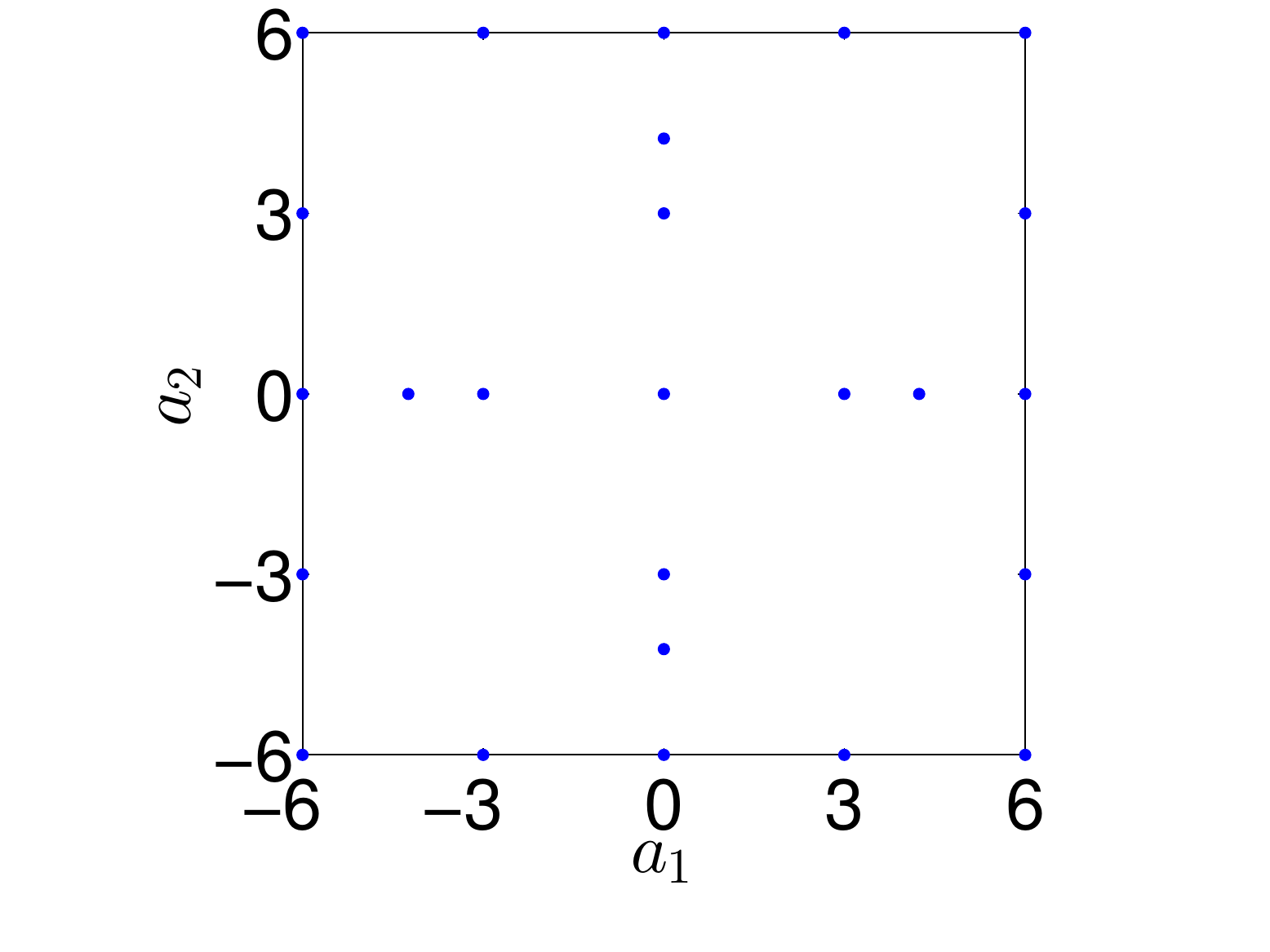}
            \includegraphics[height=4.5cm]{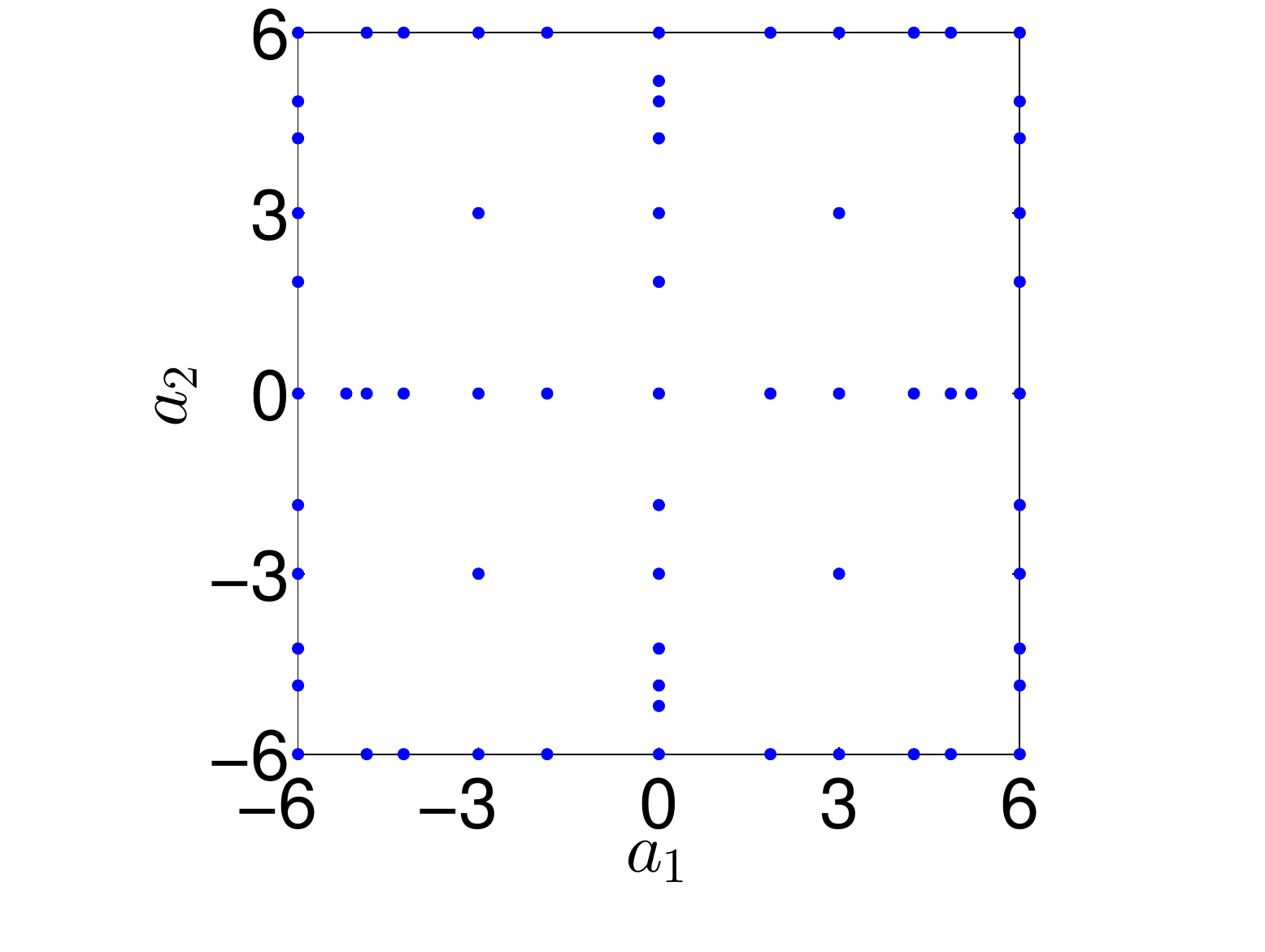}
            \includegraphics[height=4.5cm]{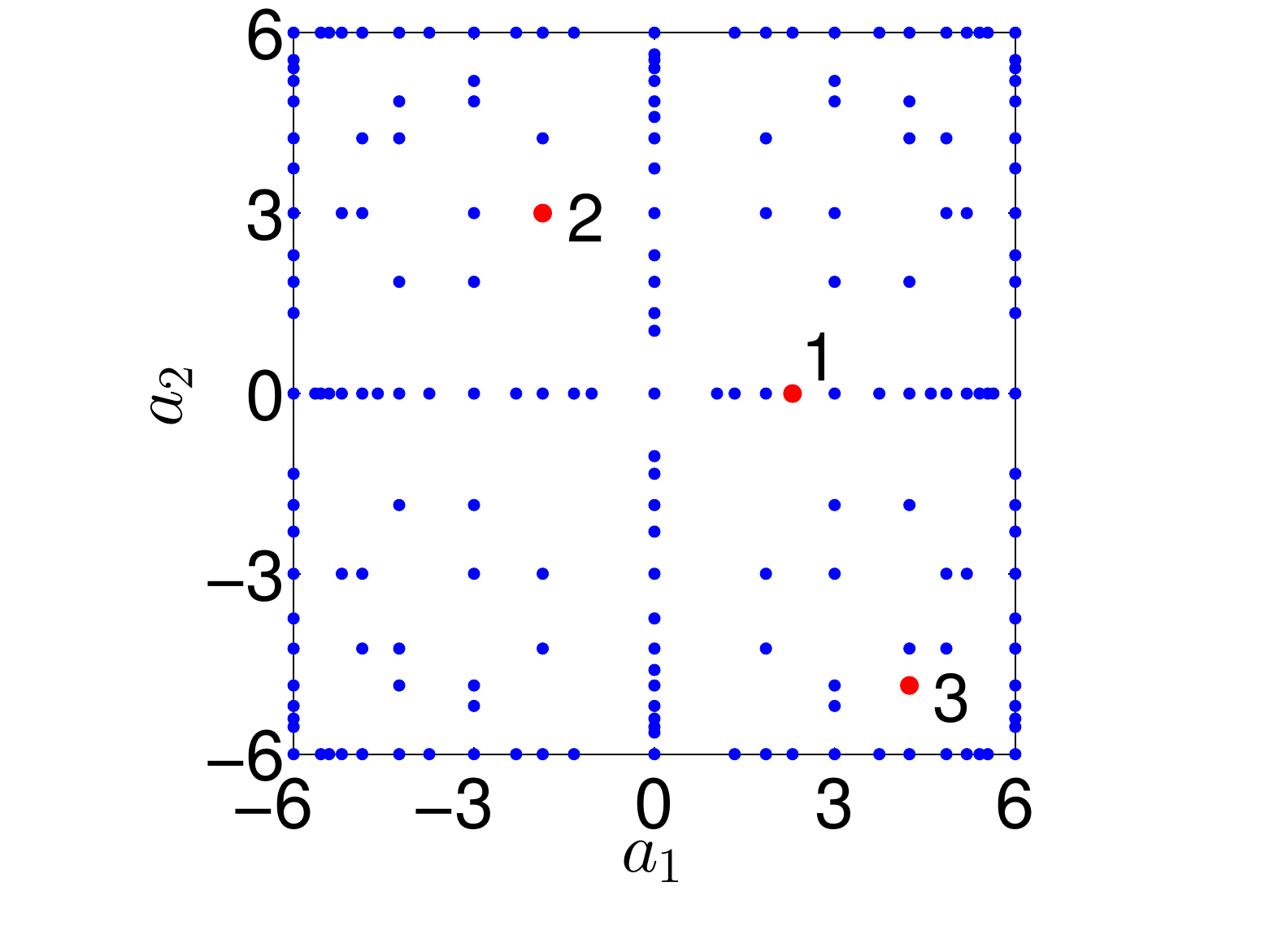}
} 
\centerline{\includegraphics[height=4cm]{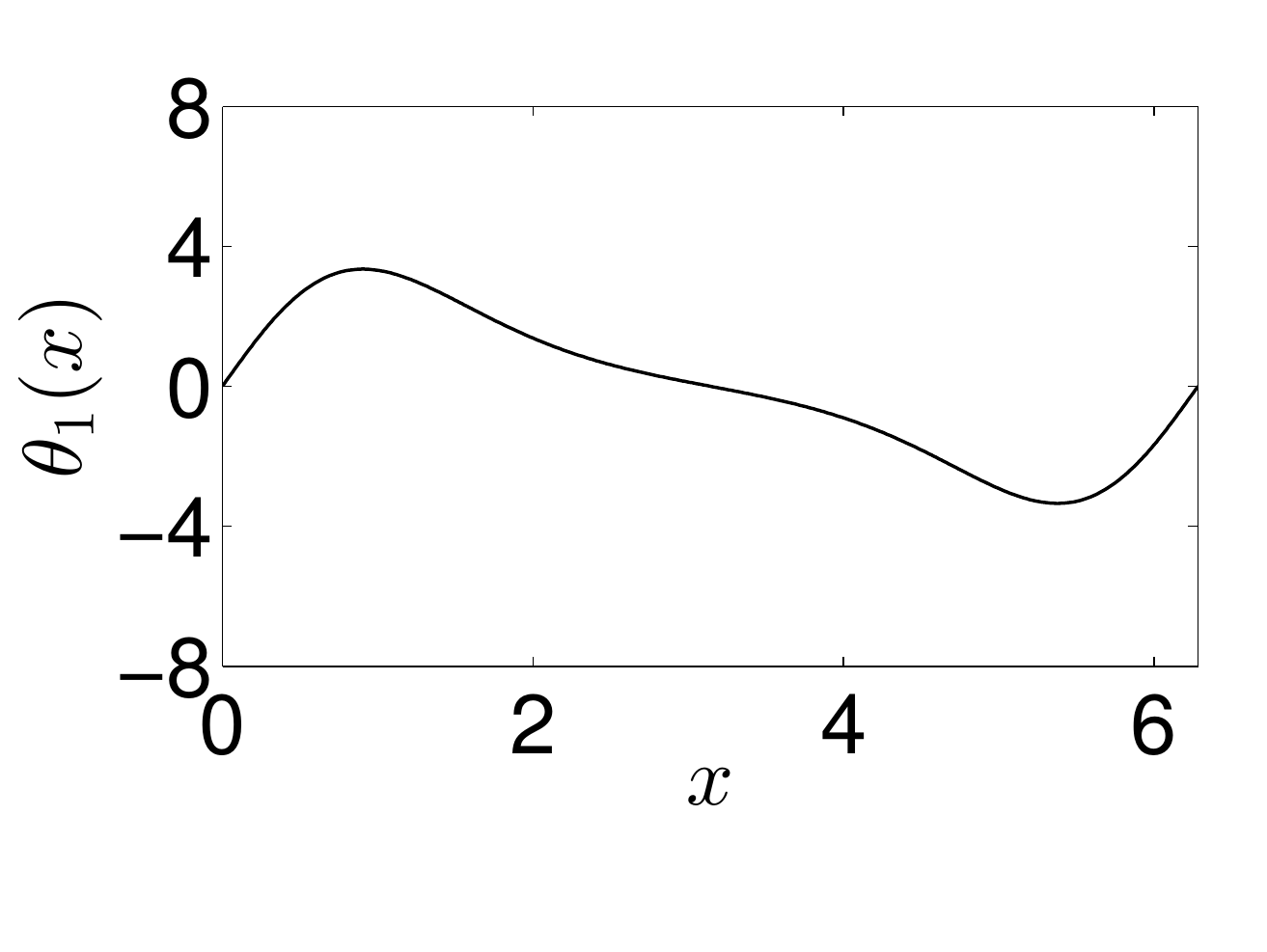} 
            \includegraphics[height=4cm]{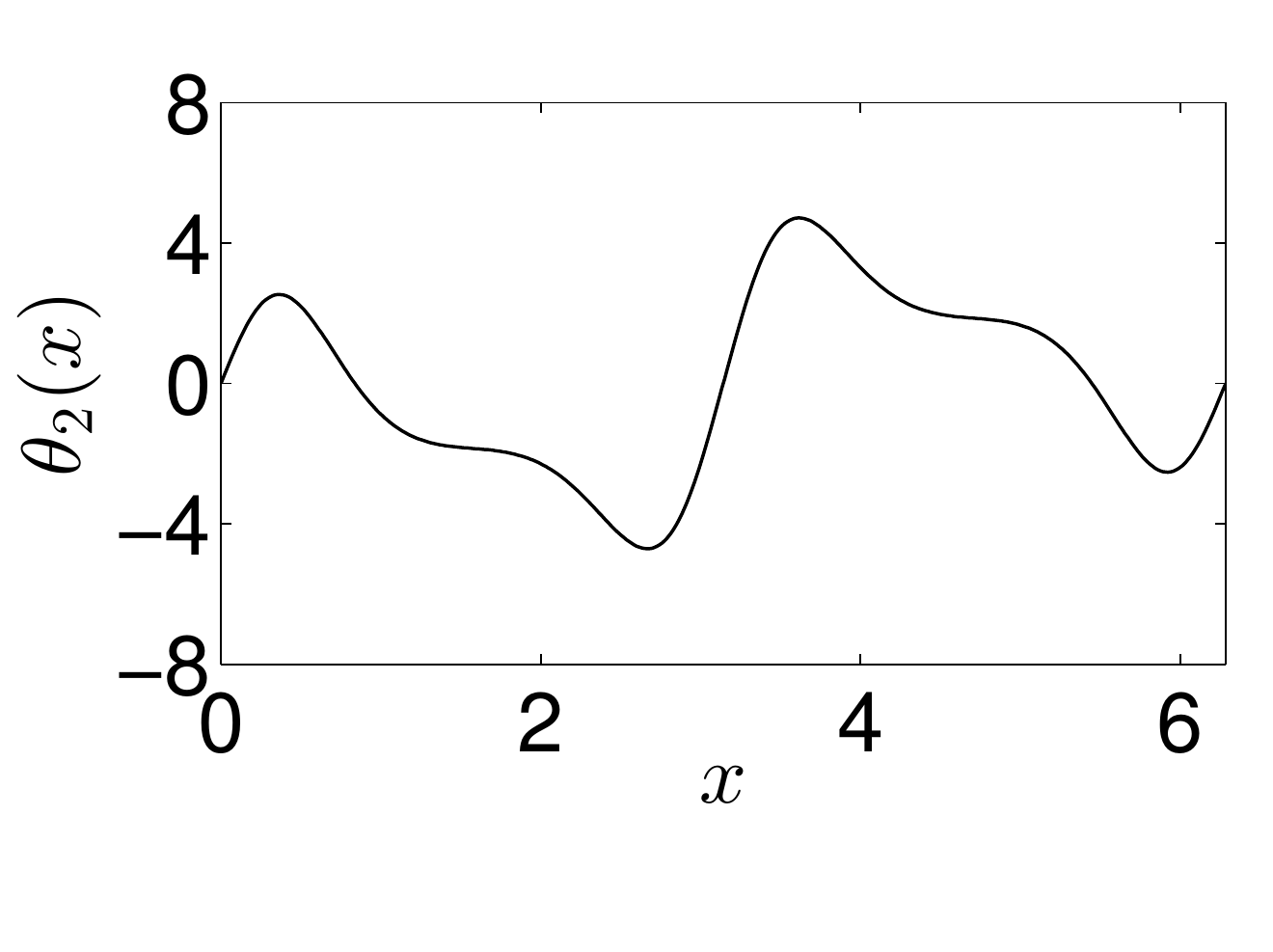}
            \includegraphics[height=4cm]{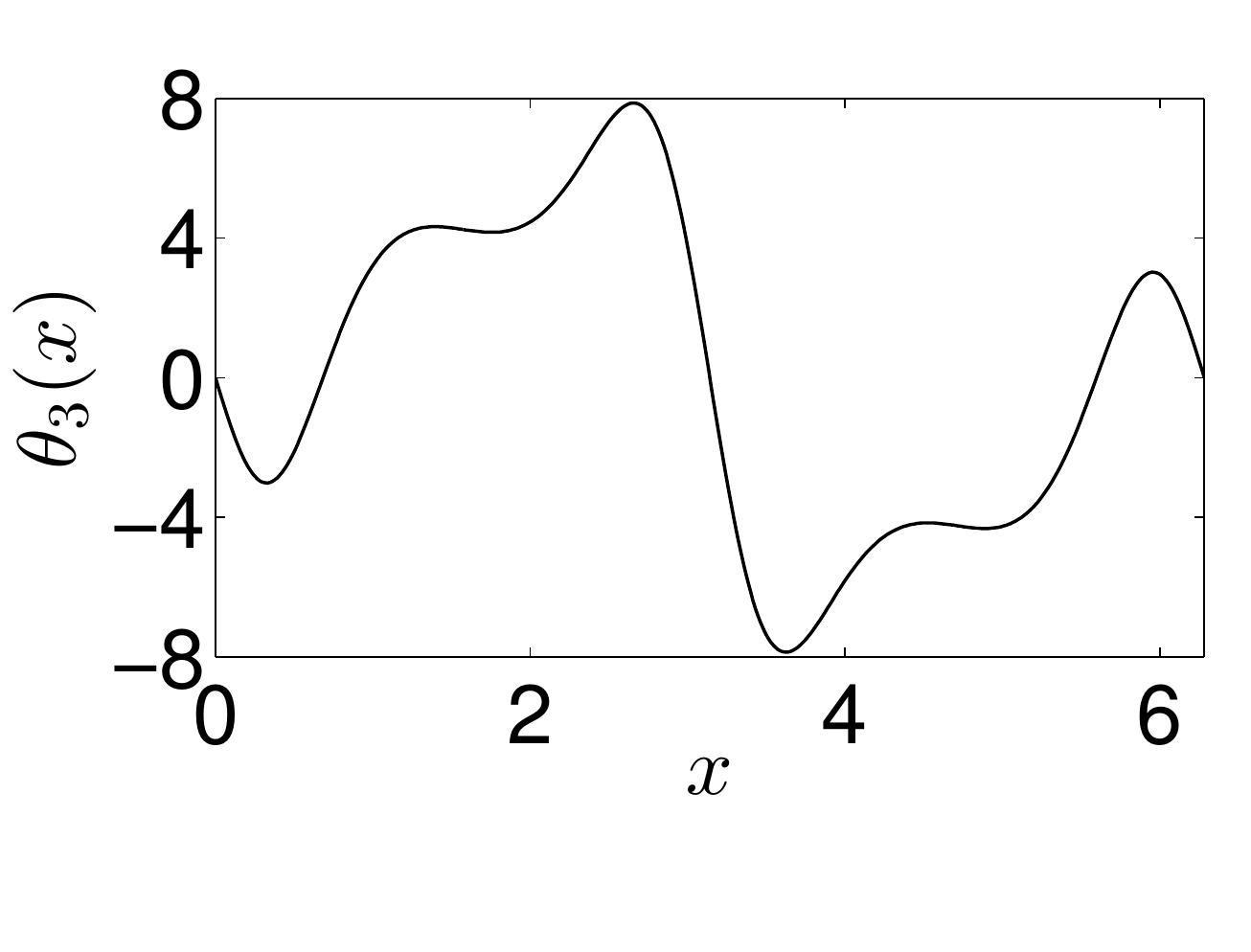}
}
\caption{Clenshaw-Curtis point sets and corresponding nodes in  
$\theta_i(x)=a_{1_i}\varphi_1(x)+a_{2_i}\varphi_2(x)$ 
($i=1, 2, 3$) in the test function space $\widetilde{S}^2_2$. Specifically, 
we plot $\theta_i(x)$ corresponding to the nodes $1$, $2$ and $3$ marked 
in red in the top right Figure. The basis elements $\varphi_k(x)$ here are 
chosen as $\sin(kx)\exp(\cos(kx))$  ($k=1,2$).}
\label{fig:optimal_nodes_2D}
\end{figure}

In a finite-dimensional setting, we can sample 
the coefficients $a_{i_p}$ in \eqref{SNq}, e.g.,  at sparse 
grids locations \cite{Barthelmann,Bungartz}, thai is at 
unions of appropriate tensorizations of one-dimensional point 
sets such as Gauss-Hermite, Clenshaw-Curtis, 
Chebyshev or Leja \cite{Akil}. This yields the set 
\begin{equation}
\widetilde{S}^{(m)}_q=\left\{\theta_i(x)\in D(F)\,\,|\,\, \theta_i(x)=
a_{i_1}\varphi_{i_1}(x)+\cdots + 
a_{i_q}\varphi_{i_q}(x)\right\},\qquad i_s=1,...,m,
\label{SNq3}
\end{equation}
where the vector $(a_{i_1},...,a_{i_q})$ takes discrete values 
at sparse grid nodes.
As an example, in Figure \ref{fig:optimal_nodes_2D} we plot three 
Clenshaw-Curtis grids and few samples of the corresponding 
interpolation nodes in $\widetilde{S}^{(2)}_2$.
For illustration purposes, the basis 
elements $\varphi_1(x)$ and $\varphi_2(x)$ here are chosen as
\begin{equation}
 \varphi_1(x)=\sin(x)e^{\cos(x)},\qquad \varphi_2(x)=\sin(2x)e^{\cos(2x)}.
\end{equation}
The construction of sparse grids usually follows 
the Smolyak algorithm. Other dimension-adaptive schemes and greedy Leja 
rules were recently proposed in \cite{Chkifa,Akil,Akil1}.

\subsubsection{Polynomial Interpolation of Nonlinear Functionals}
\label{sec:functional_polynomial_interpolation}
Let $F:D(F)\rightarrow Y$ ($Y=\mathbb{R}$ or $\mathbb{C}$) 
be a nonlinear real- or complex-valued functional on $D(F)$. 
Consider the set of polynomial functionals 
of degree $n$ 
\begin{equation}
\Pi_n([\theta])=L_0+L_1([\theta])+L_2([\theta],[\theta])+\cdots +
L_n([\theta],...,[\theta]),\label{pfun}
\end{equation}
where $L_0\in Y$, and $L_k:D(F)^k\rightarrow Y$ are $k$-linear 
{\em symmetric} functionals
\begin{equation}
L_k([\theta_1],...,[\theta_p])=\int\cdots\int K_k(x_1,..,x_k)
\theta_1(x_1)\cdots 
\theta_k(x_k)dx_1\cdots dx_k,\qquad k=1,...,n.
\label{lP}
\end{equation}
A comparison between equations \eqref{pfun} and \eqref{power_series}
yields
\begin{equation}
P_k([\theta])= L_k([\theta],...,[\theta]), 
\end{equation}
and therefore we can equivalently write \eqref{pfun} as 
\begin{equation}
\Pi_n([\theta])=\sum_{k=1}^n P_k([\theta])
\label{pfun1}.
\end{equation}
The symmetry assumption on $L_k$ implies that 
$K_k(x_1,...,x_k)$ are symmetric kernels, i.e., 
any permutation of $x_1$, ..., $x_n$ leaves $K_k$ unchanged. 
It is obviously possible to define polynomial functionals 
with non-symmetric kernels. However, such functionals can be {\em always} 
written in a symmetric form by rearranging the kernel functions appropriately. 
For example, let $H_2(x_1,x_2)$ be non-symmetric. It is easy to verify that 
\begin{equation}
\int_a^b \int_a^b H_2(x_1,x_2)\theta(x_1)\theta(x_2)dx_1dx_2=
\int_a^b \int_a^b K_2(x_1,x_2)\theta(x_1)\theta(x_2)dx_1dx_2,\qquad\forall \theta\in D(F),
\end{equation}
where
\begin{equation}
 K_2(x_1,x_2)=\frac{1}{2}\left(H_2(x_1,x_2)+H_2(x_2,x_1)\right).
\end{equation}
In other words, the value of the integral does not change if we 
replace $H_2(x_1,x_2)$ with its symmetrized version $K_2(x_1,x_2)$. 
More generally, we can symmetrize any kernel $H_p(x_1,...,x_p)$ 
by summing up all terms corresponding to all possible permutations of $(x_1,...,x_p)$ 
and then dividing up by the factorial of $p$. For example, 
\begin{align}
 K_3(x_1,x_2,x_3)=\frac{1}{3!}& \left(H_{3}(x_1,x_2,x_3)+H_{3}(x_1,x_3,x_2)+
 H_{3}(x_2,x_1,x_3)+H_{3}(x_2,x_3,x_1) + \right.\nonumber\\ 
 &\left.H_{3}(x_3,x_1,x_2)+ H_{3}(x_3,x_2,x_1)\right).
\end{align}
The symmetry of the operators $L_k$ significantly 
reduces the number of collocation nodes in the function space $D(F)$ needed 
to identify kernels $K_k$, provided these are of finite-rank. 

The polynomial functional interpolation problem can be stated as follows:
Given a set of $m$ nodes $\{\theta_1(x),...,\theta_m(x)\}$  in $D(F)$, 
find a polynomial functional in the form \eqref{pfun} 
satisfying the interpolation conditions
\begin{equation}
\Pi_n([\theta_i])=F([\theta_i]),\qquad i=1,...,m.
\label{interpolation2}
\end{equation}

\paragraph{The Stone-Weierstrass Approximation Theorem} 
The possibility of approximating an arbitrary continuous functional 
in Hilbert or Banach spaces in terms polynomial functionals
is justified by theorems analogous to the classical 
Weierstrass theorem for continuous functions. 
We recall that such theorem states that 
if $f (x)$ is a continuous, real-valued function on
the closed interval $[a ,b ]$, then given any $\epsilon > 0$ 
there exists a real polynomial $p(x)$ such
that $| f (x)- p (x) | < \epsilon $ for all $x \in [a ,b ]$. 
A remarkable generalization of this result has its roots in the
{\em Stone-Weierstrass theorem} \cite{Stone}, 
which can be stated as follows: suppose that
$X$ is a compact metric space and $K$ is an algebra of continuous, 
real-valued functions on $X$ that separates points\footnote{The 
algebra $K$ separates points if for any two distinct 
elements $u_1$ ,$u_2 \in X$ there exists $\Pi_n \in K$
such that $\Pi_n ([u_1]) -\Pi_n([x_2]) \neq 0$ of $X$ and that
contains the constant function.}.
Then for any continuous, real-valued functional $F$ on $X$ and
for any $\epsilon > 0$ there exists a polynomial functional $\Pi \in K$ 
such that $\left\| F ([u]) - \Pi([u]) \right\| < \epsilon$ 
for all $u \in X$. 
The first paper dealing with this subject 
is due Frech\'et \cite{Frechet}. He showed that any continuous 
functional can be represented by a series of polynomial functionals 
whose convergence is uniform in all compact sets 
of {\em  continuous functions}. 
This result was generalized to compact sets of functions in 
Hilbert and Banach spaces by Prenter \cite{Prenter1970} and 
Istratescu \cite{Istratescu}, respectively. Other relevant work in this 
area  is \cite{Galman,PorterClark,Bertuzzi,Chaika,Poggio,Porter}.

\subsubsection{Porter Interpolants} 
\label{sec:porter}
An effective way to construct 
finite-order polynomial functionals with minimal norm 
interpolating arbitrary continuous functionals in Hilbert spaces 
was proposed by W. Porter in \cite{PorterSIAM}. The key idea relies on 
minimizing the norm of \eqref{pfun} subject to the interpolation 
conditions \eqref{interpolation2}. A natural way to impose such 
conditions is through Lagrange multipliers. This yields the variational 
principle
\begin{equation}
 \min_{K_1,..,K_n}\left\|\Pi_n\right\|^2+\sum_{j=1}^m
 \lambda_i \left(\Pi_n([\theta_i])-F([\theta_i])\right),
 \label{minpr}
\end{equation}
The minimum is relative to arbitrary variations of the kernel 
functions $K_j(x_1,..., x_j)$. Also, $\left\|\Pi_n\right\|^2$
is the norm of the polynomial functional \eqref{pfun}, which is 
defined as 
\begin{equation}
\left\|\Pi_n\right\|^2=\sum_{p=0}^n\int\cdots 
\int \left|K_p(x_1,...,x_p)\right|^2 dx_1\cdots dx_p< \infty 
\end{equation}
The solution to the variational principle \eqref{minpr} 
allows us to identify the kernel functions $K_j(x_1,..., x_j)$ and, 
correspondingly, the polynomial functional
with {\em minimal norm} interpolating $F([\theta])$ at the $m$ 
nodes $\{\theta_1,...,\theta_m\}$. 
Specifically, we obtain
\begin{align}
 K_p(x_1,...,x_p)=\sum_{j,k=1}^m \pi^{(p)}_j(x_1,..,x_p)H^{-1}_{jk} F([\theta_k]), 
 \qquad i=0,...,n.
 \label{Kp}
\end{align}
In this equation, 
\begin{align}
\pi^{(0)}_j=1,\qquad 
 \pi^{(p)}_j(x_1,..,x_p)=\theta_j(x_1)\cdots\theta_j(x_p)\qquad p=1,2, ...
\end{align}
while the (symmetric) matrix $H_{ij}$ is defined as 
\begin{equation}
 H_{ij}=1+\left(\theta_i,\theta_j\right)+\left(\theta_i,\theta_j\right)^2+
 \cdots+\left(\theta_i,\theta_j\right)^n,
\label{tt1}
 \end{equation}
where $(,)$ denotes the $L_2(V)$ inner product, $V$ being 
the domain of the interpolation nodes $\theta_j(x)$.
The polynomial functional $\Pi_n$ constructed in this way exists 
if $F([\theta_k])$ is in the range of the matrix $\bm H$ 
(see \cite{PorterSIAM} for further details).
If one wants to approximate $F([\theta])$ in terms of 
a superimposition of monomials with orders defined by 
an index set $\mathcal{I}$ then
\begin{equation}
H_{ij}=\sum_{p\in\mathcal{I}} \left(\theta_i,\theta_j\right)^p
\label{matH}
\end{equation}
The {\em total degree} of the polynomial functional is the largest 
number in the index set $\mathcal{I}$.
It is convenient to write Porter's interpolant in terms of basis 
functionals $g_i([\theta])$ as
\begin{equation}
\Pi_n \left([\theta]\right)=\sum_{k=1}^m F\left([\theta_k]\right) g_{k}([\theta]),
\label{Porter_interpolant}
\end{equation}
where 
\begin{equation}
 g_{k}([\theta])=\sum_{j=1}^m H_{jk}^{-1}\sum_{p\in\mathcal{I}}\left(\theta_j,\theta\right)^p.
 \label{gi_0}
\end{equation}
The functional interpolant \eqref{Porter_interpolant}-\eqref{gi_0} 
has the following properties: 
\begin{enumerate}
\item $\{g_k([\theta])\}$ is a set of cardinal basis functionals, i.e., $g_{k}([\theta_q])=\delta_{kq}$. This implies that Porter's 
interpolant is a {\em cardinal Lagrangian interpolant}.

{\color{r}
\item If the interpolation nodes $\{\theta_1,...,\theta_m\}$ are 
orthonormal with respect to the inner product $(,)$ then 
$H_{ii}=\# \mathcal{I}$ (cardinality of the index set $\mathcal{I}$) 
and $H_{ij}$ ($i\neq j$) either equal to one or zero, 
depending on whether we have $\{0\}$ in the 
set $\mathcal{I}$ or not.  In every case, $\bm H$ is a matrix with 
diagonal entries equal to $\# \mathcal{I}$ and off-diagonal entries 
equal to wither zero or one. Such matrix is {\em always invertible} 
provided $\mathcal{I}$ does not reduce to the single element $\{0\}$.
}

\item Porter's interpolant is degenerate for $\mathcal{I}=\{0\}$ as the matrix 
$\bm H$ is rank one and therefore it is not invertible. The Moore-Penrose 
pseudoinverse $\bm H^+$, however, exists and it provides the 
correct form of the interpolant. To show this in a simple case, consider 
the constant functional $F([\theta])=c\in \mathbb{R}$ and the 
zero-order polynomial interpolant at $\{\theta_1,...,\theta_m\}$ 
\begin{equation}
\Pi_0([\theta])=K_0=\sum_{k,j=1}^m H^{-1}_{jk} F([\theta_j]).
\qquad H_{ij}=1.
\end{equation}
Clearly, $H^{-1}_{jk}$ does not exist since $H_{ij}=1$ and therefore 
$\textrm{rank}(H)=1$. However, the Moore-Penrose pseudoinverse 
of $\bm H$  has components $H^+_{ij}=1/m^2$, and therefore $\Pi_0([\theta])=c=F([\theta])$.

\item The polynomial functional \eqref{Porter_interpolant}-\eqref{gi_0} 
is an interpolant if and only if the matrix \eqref{matH} is invertible, i.e., 
full rank. This depends on both the choice of interpolation 
nodes and on the index set $\mathcal{I}$. 
The Moore-Penrose pseudoinverse $\bm H^+$, in general, {\em does 
not allow} to satisfy the interpolation condition. Indeed, by evaluating 
\eqref{Porter_interpolant} at the interpolation nodes $\theta_k$ we obtain
\begin{equation}
\Pi_n \left([\theta_i]\right)=\sum_{k=1}^m F\left([\theta_k]\right) 
\sum_{j=1}^m H_{jk}^{-1}H_{ji} = F([\theta_i]).
\end{equation}
However, if we replace $H_{jk}^{-1}$ with $H_{jk}^{+}$ 
then we get $\Pi_n \left([\theta_i]\right)\neq F \left([\theta_i]\right)$, 
since $\bm H^+\bm H\neq \bm I$.

\item Consider the approximation of constant functionals $F([\theta])=c$ by 
polynomial functionals of order one, i.e.,  $\Pi_{1}([\theta])=L_0+L_1([\theta])$.
Assume that the interpolation nodes $\{\theta_k\}_{k=1,...,m}$ are orthogonal functions 
with norm that decreases with $m$ as $1/m$ (see Section 
\ref{sec:results linear functionals}). Let $H_{ij}=1+h\delta_{ij}$ 
where $h=\left\|\theta_j\right\|^2$ be the matrix \eqref{matH}
corresponding to the index set $\mathcal{I}=\{0,1\}$.
By using the identity 
\begin{equation}
 \lim_{h\rightarrow 0}\frac{1}{h}\left(1+h\delta_{ij}\right)^{-1}=\delta_{ij}-\frac{1}{m}
\end{equation}
we obtain
\begin{equation}
\Pi_{1}([\theta])\simeq c\sum_{k,j=1}^m 
\left(h\delta_{jk}-\frac{h}{m}\right)
\left(1+(\theta,\theta_j)\right) \qquad \textrm{as}\qquad h\rightarrow 0\quad \textrm{and} 
\quad m\rightarrow\infty.
\end{equation}
This means that $\displaystyle\lim_{m\rightarrow \infty} \Pi_{1}([\theta])=c$, i.e., the 
polynomial interpolant is {\em consistent} with the 
functional $F$ in the sense that 
the linear term becomes smaller and smaller as we increase the number 
of test functions $\theta_k$. In the limit $m\rightarrow \infty$ (infinite number of 
test functions) we see the linear term is absent, and we correctly 
recover the constant functional. 
\end{enumerate}
By extending these arguments to higher-order polynomial functionals 
in Hilbert spaces, one can show that Porter's interpolants of order $n$ 
{\em converge pointwise to entire functionals or any polynomial 
functional of order $n$ or less} as the number of interpolation 
nodes $\theta_k$ goes to infinity (see \cite{Khlobystov2} and 
Theorem 1 in \cite{Khlobystov0}).

\paragraph{Functional Derivatives} The functional derivatives of Porter 
interpolants can be easily determined by computing the functional 
derivatives of the basis functionals $g_i([\theta])$ defined in \eqref{gi_0}. 
To this end, we first notice that
\begin{equation}
\frac{\delta(\theta,\theta_k)^p}{\delta\theta(x)}=(\theta,\theta_k)^{p-1}p\theta_k(x), 
\qquad p\geq1.
\end{equation}
A substitution of this formula into \eqref{gi_0} yields 
\begin{equation}
 \frac{\delta g_{k}([\theta])}{\delta\theta(x)}=
 \sum_{j=1}^m \theta_j(x)H_{jk}^{-1}\sum_{p\in\mathcal{I}}\left(\theta_j,\theta\right)^{p-1}p
 \label{dgi_0}
\end{equation}
By evaluating $g_{k}([\theta])$ at the nodes $\theta_j(x)$ we obtain a 
functional generalization of the classical differentiation matrix \cite{Hesthaven} 
\begin{equation}
D^{(1)}_{ji}(x)=\frac{\delta g_i([\theta_j])}
{\delta \theta(x)}. 
\label{Hint}
\end{equation}
Similarly, the second-order functional derivative of $g_k([\theta])$ is 
\begin{equation}
 \frac{\delta^2 g_{k}([\theta])}{\delta\theta(x)\delta\theta(y)}=
 \sum_{j=1}^m \theta_j(x) \theta_j(y)H_{jk}^{-1}
 \sum_{p\in\mathcal{I}}p(p-1)\left(\theta_j,\theta\right)^{p-2},
 \label{dgi_00}
\end{equation}
and it yields the following second-order functional differentiation matrix
\begin{equation}
D^{(2)}_{ji}(x,y)= \frac{\delta^2 g_{i}([\theta_j])}{\delta\theta(x)\delta\theta(y)}. 
\label{H2int}
\end{equation}
At this point it is useful to provide simple examples 
of functional interpolation in Hilbert spaces.
\vs
\noindent
{\em Example 1:} 
Consider the first-order polynomial functional
\begin{equation}
P_1([\theta]) =\int_a^b K_1(x)\theta(x)dx.
\end{equation}
To represent $P_1([\theta])$ in terms of a functional interpolant it is 
sufficient to consider the set of orthonormal functions 
$\widehat{S}^{(m)}_1=\{\varphi_1,...,\varphi_m\}$ (see Eq. \eqref{SNq1}).
In this case, Porter's cardinal basis functionals \eqref{gi_0} reduce to 
\begin{equation}
g_i([\theta])=\left(\varphi_i,\theta\right),
\end{equation}
and the functional interpolant can be written as 
\begin{equation}
\Pi_1([\theta])=\sum_{i=1}^m P_1([\varphi_i])g_i([\theta]).
\end{equation}
Clearly, we have that $\Pi_1([\theta])\rightarrow P_1$ as 
$m\rightarrow \infty$.

\begin{figure}
\centerline{\footnotesize\hspace{1cm}$\hat{\hat{S}}^{(10)}_1$ \hspace{5.cm}$
\hat{\hat{S}}^{(10)}_2$
\hspace{5.cm} $\hat{\hat{S}}^{(10)}_3$}
\centerline{\includegraphics[height=5.5cm]{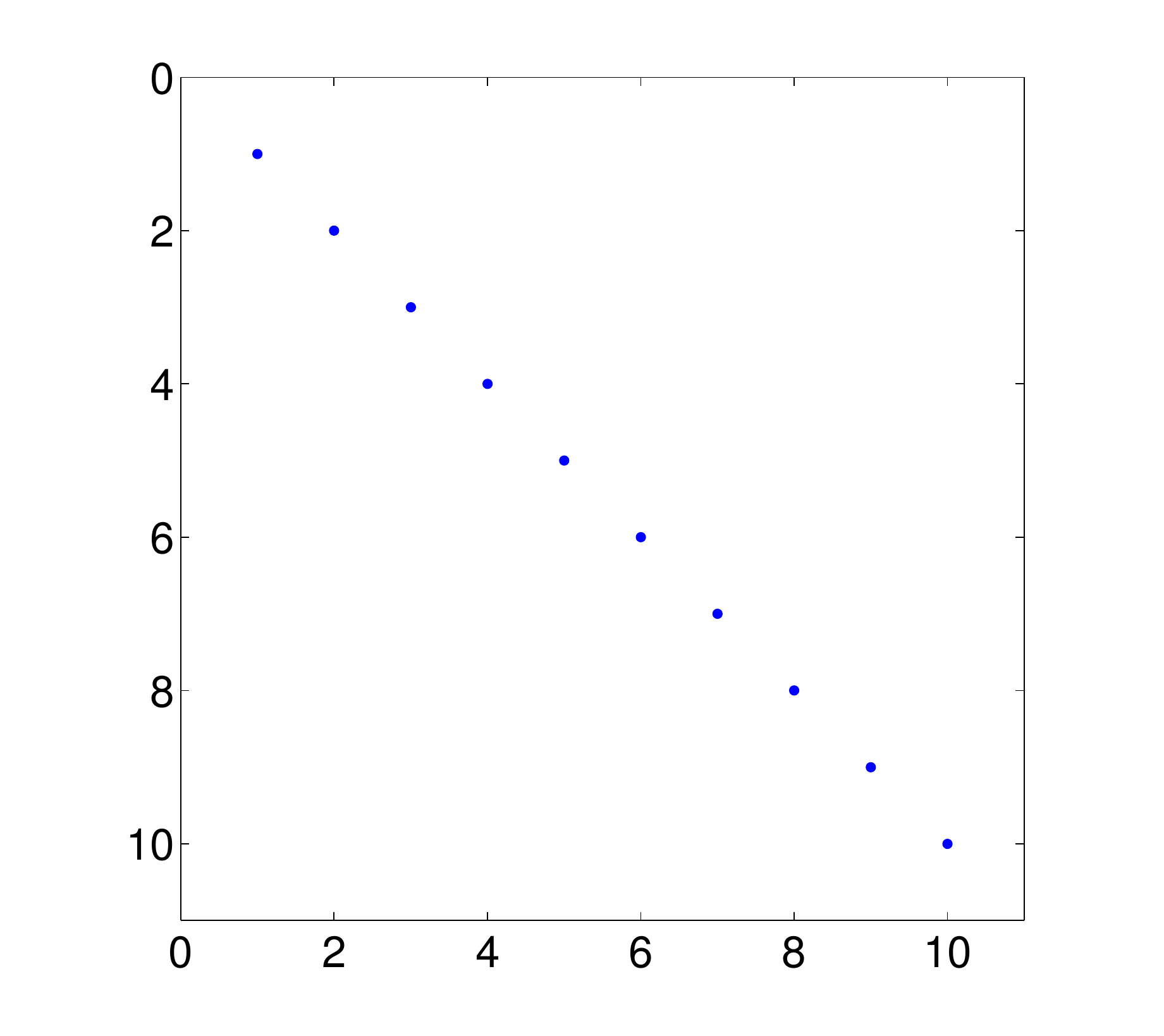}
            \includegraphics[height=5.5cm]{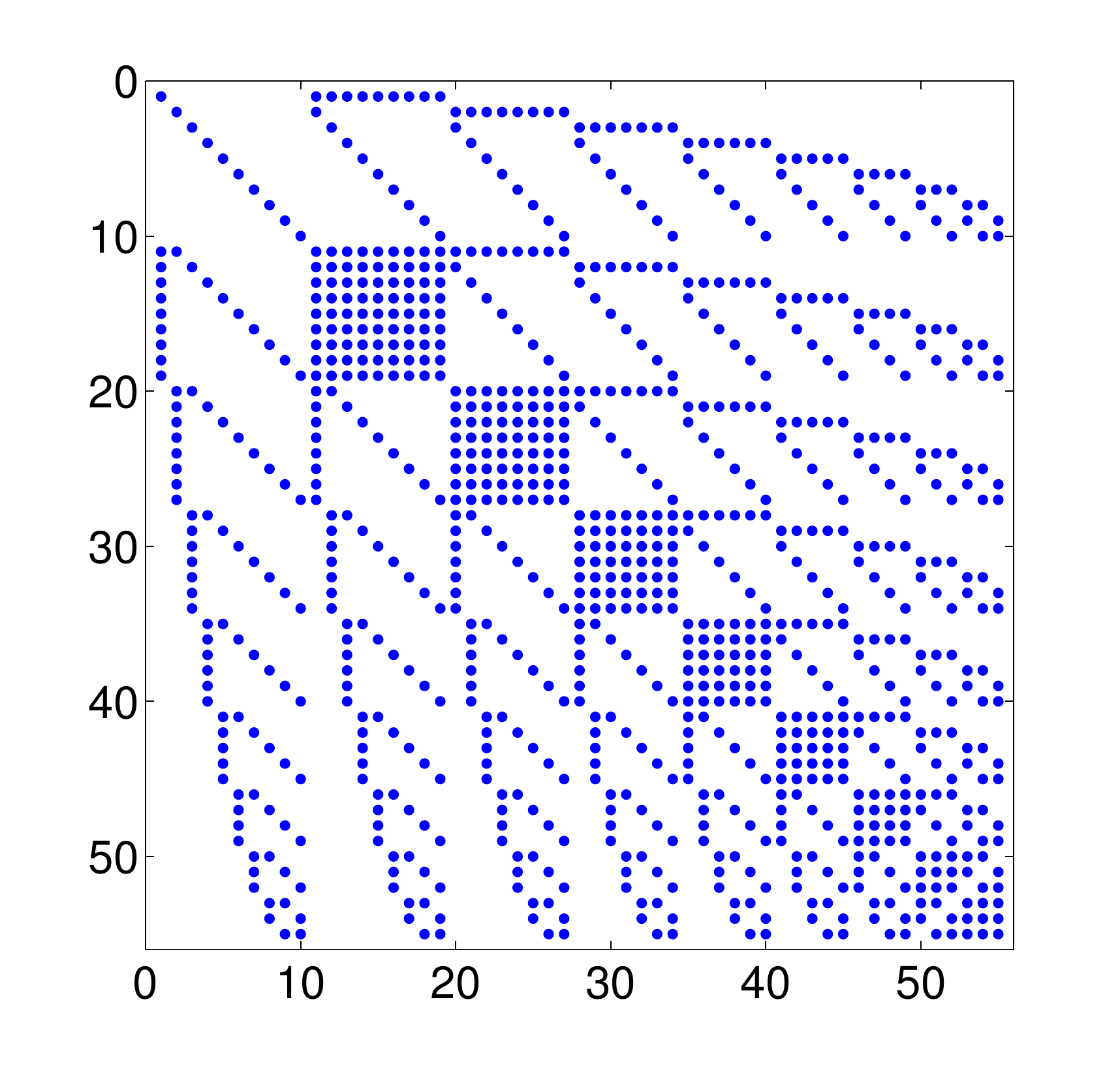}
            \includegraphics[height=5.5cm]{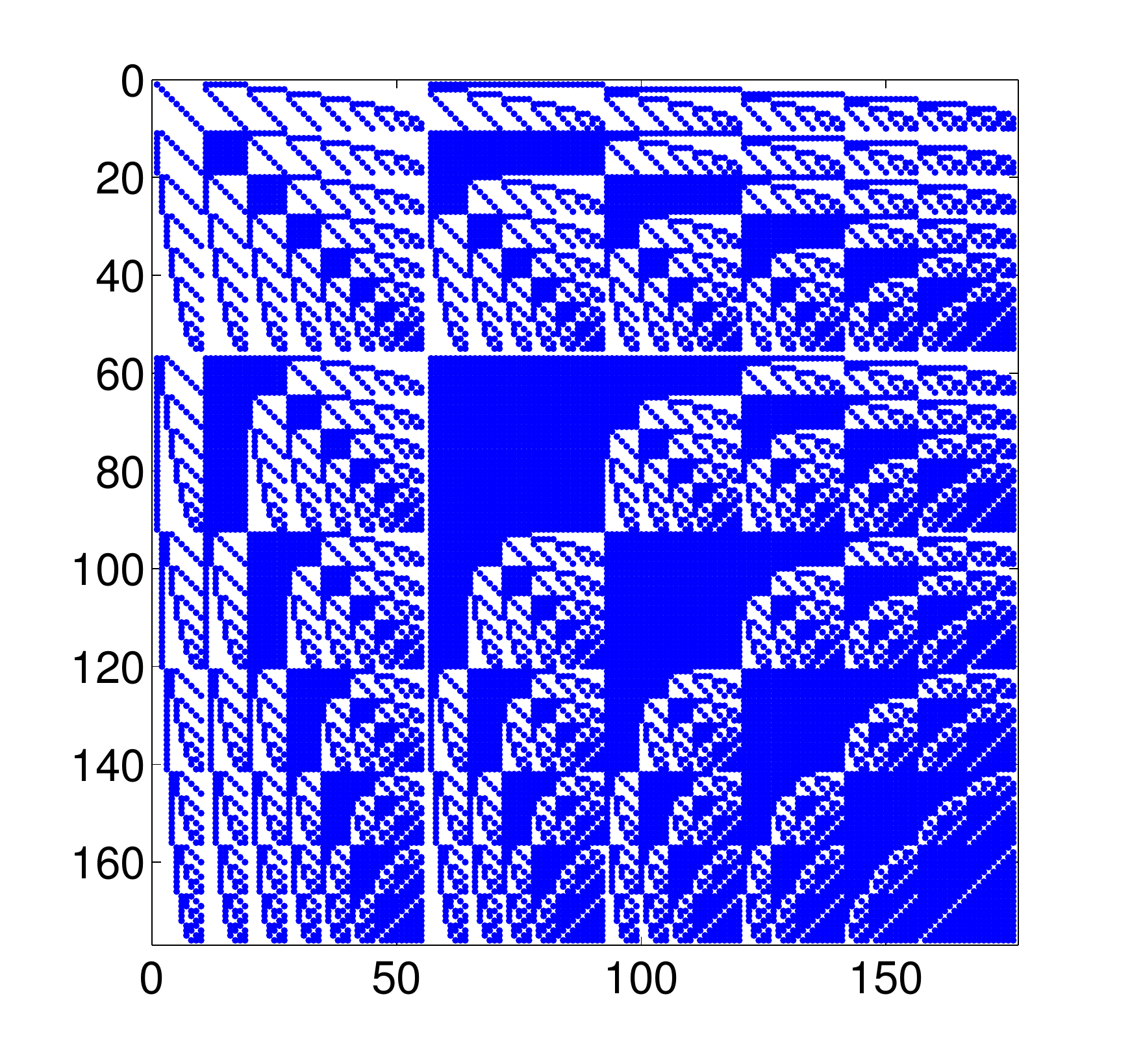}}
\caption{Structure of the $H$-matrix \eqref{matH} associated with 
the function set $\hat{\hat{S}}^{(10)}_1$, $\hat{\hat{S}}^{(10)}_2$ and 
$\hat{\hat{S}}^{(10)}_3$ defined in Eq. \eqref{sset}. 
We consider $\mathcal{I}=1$, $\mathcal{I}=2$ 
and $\mathcal{I}=3$, respectively.}
\label{fig:porter_matrix_structure}
\end{figure}

\vs
\noindent
{\em Example 2:}
Consider the second-order polynomial functional
\begin{equation}
P_2([\theta]) =\int_a^b\int_a^b K_2(x,y)\theta(x)\theta(y)dxdy.
\end{equation}
To represent $P_2([\theta])$ in terms of a functional interpolant it is 
sufficient to consider the set of orthonormal functions 
\begin{equation}
\hat{\hat{S}}^{(m)}_2=\{\varphi_1,...,\varphi_m,(\varphi_1+\varphi_2),...,
(\varphi_1+\varphi_m),(\varphi_2+\varphi_3),  ..., (\varphi_2+\varphi_m), ...\}
\label{sset}
\end{equation}
which is similar (but not equal) to \eqref{SNq1}.
The matrix \eqref{matH} associated with this set has the structure shown in 
Figure \ref{fig:porter_matrix_structure}.
The cardinal basis functionals \eqref{gi_0} reduce to 
\begin{equation}
 g_{k}([\theta])=\sum_{j=1}^{\#\hat{\hat{S}}^{(m)}_2} H_{jk}^{-1}\left(\theta_j,\theta\right)^2,
\end{equation}
where $\#\hat{\hat{S}}^{(m)}_2$ is the number of elements 
of $\hat{\hat{S}}^{(m)}_2$, 
$\theta_k=\varphi_k$ ($k=1,...,m$), $\theta_{m+1}=\varphi_1+\varphi_2$, 
etc. The functional interpolant can be written as 
\begin{equation}
\Pi_2([\theta])=\sum_{i=1}^{\#\hat{\hat{S}}^{(m)}_2} P_2([\theta_i])g_i([\theta]),
\end{equation}
and it converges $P_2$ as $m\rightarrow \infty$ (see 
Section \ref{sec:numerical results functionals}).

\vs
\noindent
{\em Example 3:} Consider the third-order polynomial functional
\begin{equation}
P_3([\theta]) =\int_a^b\int_a^b\int_a^b K_3(x_1,x_2,x_3)
\theta(x_1)\theta(x_2)\theta(x_3)dx_1dx_2dx_3.
\label{p3}
\end{equation}
To represent $P_3([\theta])$ in terms of a functional interpolant it is 
sufficient to consider the set of functions $\hat{\hat{S}}^{(m)}_3\subset\hat{S}^{(m)}_3$
defined as 
\begin{equation}
 \hat{\hat{S}}^{(m)}_3=\hat{\hat{S}}^{(m)}_2\cup \{(\varphi_i+\varphi_j+\varphi_k), 
 \quad k>j>i\}
\end{equation}
where $\hat{\hat{S}}^{(m)}_2$ is as in \eqref{sset}.
The matrix \eqref{matH} associated with this set has the structure shown in 
Figure \ref{fig:porter_matrix_structure}
The cardinal basis functionals \eqref{gi_0} reduce to 
\begin{equation}
 g_{k}([\theta])=\sum_{j=1}^{\#\hat{\hat{S}}^{(m)}_3}
 H_{jk}^{-1}\left(\theta_j,\theta\right)^3.
\end{equation}
The functional interpolant can be written as 
\begin{equation}
\Pi_3([\theta])=\sum_{i=1}^{\#\widehat{S}^{(m)}_3} P_3([\theta_i])g_i([\theta]).
\end{equation}

\paragraph{More-Penrose Pseudoinverse and Non-Cardinal Basis Functionals}
We emphasized that the matrix $H_{ij}$ 
defined in \eqref{matH} may be not invertible in some 
cases. This happens, for example, if the interpolation 
nodes $\theta_k(x)$ are linearly dependent or if there exist 
a symmetry such the inner product of $\theta_k$ and $\theta_j$ 
yields linearly dependent rows/columns
in \eqref{matH}. In this cases we can still construct a 
polynomial functional with {\em minimal norm} which, however, 
{\em does not} interpolate $F$ at $\theta_k$. 
To this end, we simply use the Moore-Penrose 
pseudoinverse of $H_{ij}$ to obtain a representation 
of Porter's polynomial functionals in terms of 
a {\em non-cardinal basis} $g^+_{k}([\theta])$ as
\begin{equation}
\Pi_n \left([\theta]\right)=\sum_{k=1}^m 
F\left([\theta_k]\right) g^+_{k}([\theta]).
\label{Porter_interpolant_noncardinal}
\end{equation}
where 
\begin{equation}
 g^+_{k}([\theta])=\sum_{j=1}^m H_{jk}^{+}\sum_{p\in\mathcal{I}}
 \left(\theta_j,\theta\right)^p.
 \label{gi_0_noncardinal}
\end{equation}
In the last equation $H_{ij}^+$ denotes the Moore-Penrose 
pseudoinverse of $H_{ij}$. 
The approximation properties of polynomial functionals 
in the form \eqref{Porter_interpolant_noncardinal} will be 
studied in Section \ref{sec:Hopf polynomial} and Section 
\ref{sec:sine functional}.

\paragraph{Recursive Porter Interpolation}

The number of interpolation nodes 
$\theta_k$ required to represent exactly a  
polynomial functional of order $n$ is given in 
\eqref{SCardinality}.
For example, if we set $m=10$ elementary functions 
and polynomial order $n=12$ such formula 
yields $646646$ interpolation nodes!
Such large number of nodes may be an issue 
when computing Porter's interpolants. In fact, 
computing the inverse of \eqref{matH} 
rapidly becomes intractable as we increase 
either $m$ or $n$. To overcome this problem 
we can split the process of inverting the matrix \eqref{matH} 
into a recursive algorithm, e.g., by using Schur 
complements and blockwise inversion. 
To this end, consider the set of nodes
\begin{equation}
\left\{\{\theta_1,...,\theta_m\},\{\theta_{m+1},...,\theta_{2m}\},
\{\theta_{2m+1},...,\theta_{3m}\}, ...,\{\theta_{(N-1)m+1},...,\theta_{Nm}\}\right\},
\label{partitioned_nodes}
\end{equation}
and define the matrices
\begin{align}
H^{(I,J)}_{ij}=\sum_{p\in\mathcal{I}}(\theta_i,\theta_j)^p\qquad i=(I-1)m+1,..,Im\quad j=(J-1)m+1,..,Jm,
\end{align}
where $\mathcal{I}$ denotes the index set of Porter's monomials 
while $I$ and $J$ run from $1$ to $N$. The   
$H$-matrix \eqref{matH} corresponding to the 
set \eqref{partitioned_nodes} can be represented in a block-wise 
form  as
\begin{equation}
 H_N=\left[
 \begin{array}{ccc}
H^{(1,1)} & \cdots & H^{(1,N)}\\
\vdots & & \vdots\\
H^{(N,1)} & \cdots & H^{(N,N)}\\
 \end{array}
 \right].
\end{equation}
The computational cost of inverting such matrix is prohibitive if $Nm$ 
is large. However, we can use the following recursive procedure. We first build and invert 
$H_1=H^{(1,1)}$, corresponding to the first set of $m$. This allows us to 
determine Porter's interpolant on the first set of $m$ nodes in \eqref{partitioned_nodes}.
Next we add the second set, i.e., the nodes $\{\theta_{m+1},...,\theta_{2m}\}$.
The matrix \eqref{matH} corresponding to the nodal set $\{\theta_1,...,\theta_{2m}\}$ is 
\begin{equation}
H_2=\left[
\begin{array}{cc}
H_1 &  C_1\\
C_1^T &  H^{(2,2)}\\
\end{array}
\right]\qquad C_1=H^{(1,2)}
\end{equation}
and it can be inverted by using the block-wise formula \cite{Miao}
\begin{equation}
H^{-1}_2=\left[
\begin{array}{cc}
H_1^{-1}+ B_1A_1B_1^T &  -B_1A_1\\
-A_1B_1 &  A_1 \\
\end{array}
\right]
\end{equation}
where 
\begin{align}
B_1&= H_1^{-1}  C_1,\\
A_1&=\left(H^{(2,2)}- C_1^TB_1\right)^{-1} \qquad 
\textrm{(inverse of the Schur complement)}.
\end{align} 
Now we bring in the third set of nodes $\{\theta_{2m+1},...,\theta_{3m}\}$.
The matrix \eqref{matH} corresponding to the nodal set $\{\theta_1,...,\theta_{3m}\}$ is 
\begin{equation}
H_3=\left[
\begin{array}{cc}
H_2                    &  C_2  \\
C_2^T &   H^{(3,3)}\\
\end{array}
\right]\qquad 
C_2=\left[
\begin{array}{c}
H^{(1,3)}\\H^{(2,3)}
\end{array}
\right],
\end{equation}
and its inverse is, as before,
\begin{equation}
H^{-1}_3=\left[
\begin{array}{cc}
H_2^{-1}+ B_2A_2B_2^T & -B_2A_2\\
-A_2B_2 &  A_2 \\
\end{array}
\right]
\end{equation}
where 
\begin{align}
B_2= H_2^{-1}  C_2,\qquad A_2=\left(H^{(3,3)}- C_2^TB_2\right)^{-1}.
\end{align} 
At this point it is clear that the procedure can be iterated as many times 
as needed. This generates a sequence of basis functionals \eqref{gi_0}, and  
an interpolating polynomial functional with minimal norm that
passes through an increasing number of nodes. 
In this way, we have reduced  
the problem of computing Porter's interpolant through  a very large number 
of nodes into a sequence of matrix inversions of dimension at most $m\times m$. 
The storage requirements of the algorithm just described, however, 
is not small as because Porter's basis functionals are 
ultimately defined by $H_N^{-1}$. An open question is 
the identification of interpolation nodes $\theta_j$ 
leading to minimal complexity/storage requirements 
for $H_N^{-1}$. 

\begin{figure}
\centerline{\includegraphics[height=6cm]{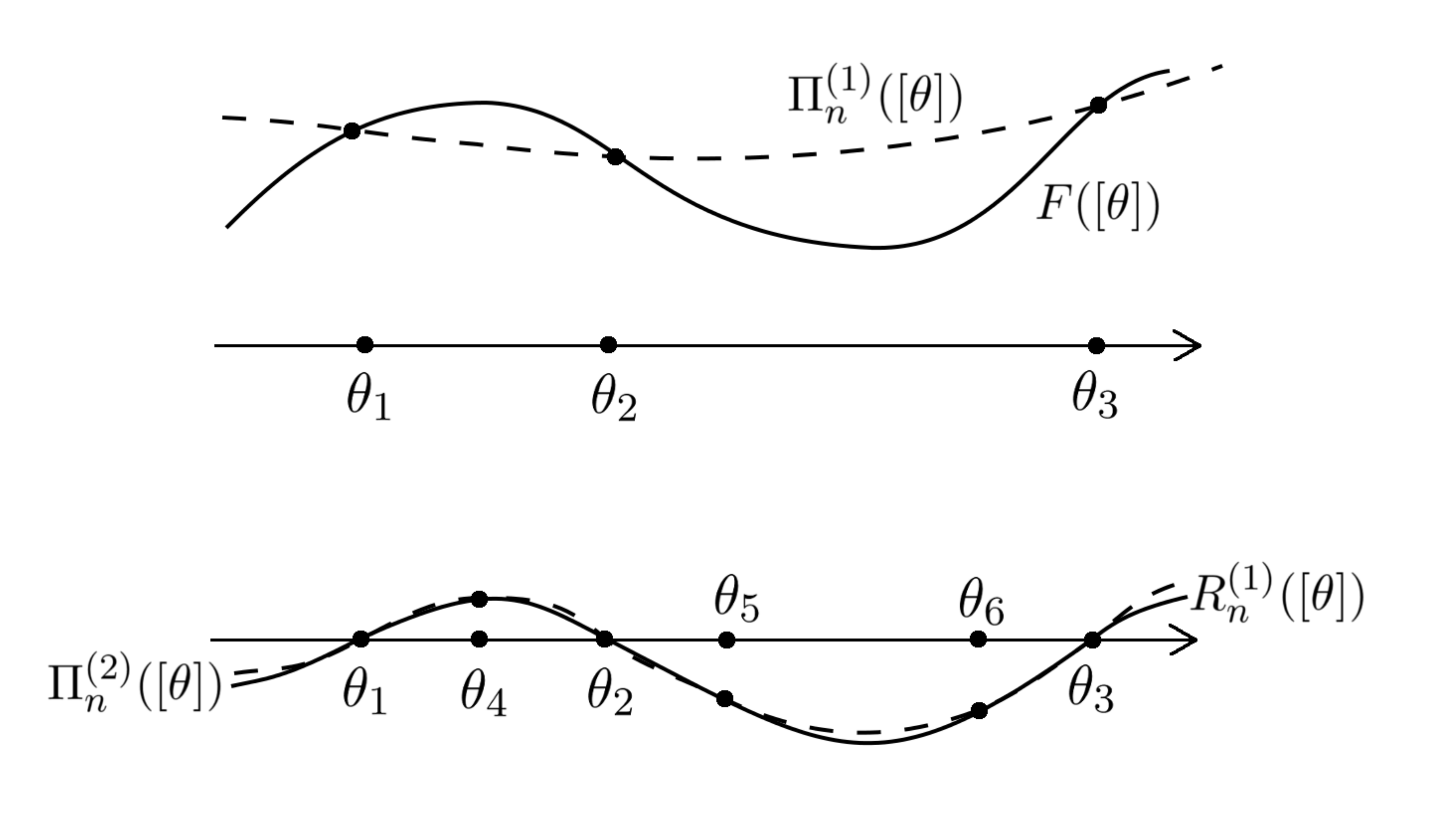}}
\caption{Recursive polynomial interpolation: Minimal residual approach.}
\label{fig:recursive_sketch}
\end{figure}

An alternative interpolation method makes use of residuals. 
The main idea is sketched in Figure \ref{fig:recursive_sketch}. 
The functional $F([\theta])$ is interpolated by a polynomial   
of order $n$, denoted as $\Pi^{(1)}_n([\theta])$, at just three 
nodes $\{\theta_1,\theta_2,\theta_3\}$. 
Subtracting $\Pi^{(1)}_n([\theta])$ from $F([\theta])$ yields 
the functional residual 
\begin{equation}
 R^{(1)}_n([\theta]) = F([\theta])-\Pi^{(1)}_n([\theta])
\end{equation}
which is zero at 
$\theta_1$, $\theta_2$ and $\theta_3$ because of the interpolation 
condition. Now we add three more nodes $\{\theta_4,\theta_5,\theta_6\}$
and construct a Porter's interpolant of $R^{(1)}_n([\theta])$ at 
$\{\theta_1,...,\theta_6\}$. We denote such interpolant 
by $\Pi^{(2)}_n([\theta])$. The computation of 
$\Pi^{(2)}_n([\theta])$ can be carried out as above by using Schur complements 
and block-wise inversion of \eqref{matH}. The polynomial functional 
$\Pi^{(1)}_n([\theta])+\Pi^{(2)}_n([\theta])$ 
interpolates $F([\theta])$ at $\{\theta_1,...,\theta_6\}$. The recursive construction 
proceeds with the definition of the new residual
\begin{equation}
R^{(2)}_n([\theta])= F([\theta])- \Pi^{(1)}_n([\theta])- \Pi^{(2)}_n([\theta]),
\end{equation}
three more nodes $\{\theta_7,\theta_8,\theta_9\}$, 
and a Porter's polynomial $\Pi^{(3)}_n([\theta])$, interpolating 
$R^{(2)}_n([\theta])$ at $\{\theta_1,...,\theta_9\}$. 
Proceeding recursively with higher-order 
residuals up to order $r$, we obtain the polynomial functional
\begin{equation}
Q_n([\theta])= \sum_{k=1}^r \Pi^{(k)}_n([\theta]).
\end{equation}
Clearly $Q_n([\theta])$ interpolates $F([\theta])$ at all 
nodes $\theta_k$ and therefore it is completely 
equivalent to $\Pi_n([\theta])$, i.e., 
$Q_n([\theta])=\Pi_n([\theta])$.

\paragraph{Hierarchical Matrices}
The algorithm we just described aim at reducing 
the computational cost of computing polynomial 
functional interpolants by inverting the matrix $\bm H$ 
defined in \eqref{matH} in a block-wise fashion or recursively.  
The structure of such matrix obviously depends on how we select 
the interpolation nodes in the function space $D(F)$. 
An interesting open question is whether we can determine 
sets of nodes for which the interpolation problem 
can be solved at a minimal cost. Stated in matrix terms, 
can we identity sets of nodes yielding structured matrices that 
can be easily inverted, e.g., {\em hierarchical matrices}? 
The matrices shown in Figure \ref{fig:porter_matrix_structure}
have indeed a self similar structure which can be used to speed up 
their inversion. We leave this question open for future research.
For uniquely solvable interpolation problems we could equivalently 
construct Khlobystov polynomial functionals (see 
Section \ref{sec:Khlobystov}), and determine the coefficients 
of the expansion by using the method of moments.

\subsubsection{Prenter Interpolants}
\label{sec:prenter}
Another method to determine polynomial 
interpolants in Hilbert and Banach spaces was introduced by 
Prenter in \cite{Prenter}. She proved that 
if $F([\theta])$ is a functional in a Hilbert space $D(F)$, 
and $\theta_j\in D(F)$ are interpolation nodes, 
then there exists a $n$th-order functional 
interpolant in the form
\begin{equation}
\Pi_n ([\theta])=\sum_{i=1}^n F([\theta_i])g_i([\theta]), 
\label{prenter}
\end{equation}
where
\begin{equation}
g_i([\theta]) = \prod_{\substack{j=1\\j\neq i}}^n \frac{
\left(\theta -\theta_j, \theta_i -\theta_j\right)}
{\left(\theta_i -\theta_j, \theta_i -\theta_j\right)} \qquad \textrm{(cardinal basis)}.
\label{gi_prenterer}
\end{equation}
As before, $(,)$ denotes the inner product in 
$L_2(V)$, where $V$ is the domain of $\theta_k(x)$.
Note that each basis element $g_i([\theta])$ is a 
polynomial functional of order $n$. 
On the other hand, Porter's method yields polynomial basis 
functionals of total degree $\max(\mathcal{I})$, where the index set 
$\mathcal{I}$ does not depend on the number of collocation points.
Porter \cite{Porter} applied Prenter's theorems to causal systems, 
while Bertuzzi, Gandolfi and Germani \cite{Bertuzzi1,Bertuzzi} 
extended Prenter's results to causal approximation of 
input-output maps in Hilbert spaces. Generalizations to Banach spaces
can be found in Chapter 3 of \cite{Torokhti} (see also \cite{Makarov}
and the references therein).
The functional derivatives of Prenter's polynomial functionals can 
be obtained by computing the functional derivatives of \eqref{gi_prenterer}. 
This yields\footnote{These expressions can be easily proved by noting that
\begin{equation}
 \left.\frac{d}{d\epsilon}g_i([\theta+\epsilon \eta])\right|_{\epsilon=0}=\sum_{\substack {k=1\\k\neq i}}^n
 \frac{(\eta,\theta_i-\theta_k)}{\left\|\theta_i-\theta_k\right\|^2}
 \prod_{\substack{j=1\\j\neq k,i}}^n
 \frac{(\theta-\theta_j,\theta_i-\theta_j)}{\left\|\theta_i-\theta_j\right\|^2}.
\end{equation}
}
\begin{equation}
\frac{\delta g_i([\theta])}{\delta \theta(x)}=\sum_{\substack {k=1\\k\neq i}}^n
\frac{\theta_i(x)-\theta_k(x)}
{\left\|\theta_i-\theta_k\right\|^2}g_{ik}([\theta]),\label{fd1}
\end{equation}
\begin{equation}
\frac{\delta^2 g_i([\theta])}{\delta \theta(x)\delta \theta(y)}=
\sum_{\substack {k=1\\k\neq i}}^n
\frac{\theta_i(x)-\theta_k(x)}
{\left\|\theta_i(x)-\theta_k(x)\right\|^2}
\sum_{\substack {s=1\\s\neq k}}^n
\frac{\theta_i(y)-\theta_s(y)}
{\left\|\theta_i(y)-\theta_s(y)\right\|^2}
g_{iks}([\theta]),\label{fd2}
\end{equation}
where
\begin{equation}
 g_{ik}([\theta])=
\prod_{\substack{j=1\\j\neq i,k}}^n
 \frac{(\theta-\theta_j,\theta_i-\theta_j)}{\left\|\theta_i-\theta_j\right\|^2},\qquad 
g_{iks}([\theta])=
\prod_{\substack{j=1\\j\neq i,k,s}}^n
 \frac{(\theta-\theta_j,\theta_i-\theta_j)}{\left\|\theta_i-\theta_j\right\|^2}.
\end{equation}
Note that the functional derivatives \eqref{fd1}-\eqref{fd2} are more 
complicated than the ones we obtained in the case Porter's basis \eqref{dgi_0}. 
As noted by Allasia and Bracco in \cite{Allasia}, Prenter's interpolants 
are badly conditioned as the number of interpolation nodes in 
function space increases. This unfortunate feature is common 
to many Lagrange interpolants.

\subsubsection{Khlobystov Interpolants}
\label{sec:Khlobystov}
We have seen that any continuous 
functional in a Hilbert or a Banach space 
can be approximated uniformly in terms of polynomial 
functionals \cite{Prenter1,Bertuzzi,Istratescu,Makarov,Howlett,Torokhti}. 
Such polynomial functionals may be built based 
on an interpolation process (see Section 
\ref{sec:porter} and Section \ref{sec:prenter}). 
What are the convergence properties of 
such expansions? 
To address this question, let us consider a polynomial interpolation 
problem of a given polynomial functional $\Pi_n([\theta])$ in 
the form \eqref{pfun1}. Assume that $D(\Pi_n)$ is a Hilbert 
space of test functions, introduce the orthonormal 
basis $\{\varphi_1,...,\varphi_m\}$ in $D(\Pi_n)$, and 
consider the following expansion
\begin{equation}
\Pi^I_n([\theta])=\sum_{k=1}^n \sum_{i_1=1}^m\cdots \sum_{i_k=1}^m 
 (\theta,\varphi_{i_1})\cdots(\theta,\varphi_{i_k})a_k(\varphi_{i_1},...,\varphi_{i_k}),
 \label{interp4}
\end{equation}
where $(,)$ is an inner product in $D(\Pi_n)$ and 
$a_k(\varphi_{i_1},...,\varphi_{i_k})$ are real- or complex-valued coefficients.
Clearly $\Pi^I_n([\theta])$ is an interpolant 
of $\Pi_n([\theta])$ at $\varphi_k$, i.e., 
\begin{equation}
\Pi^I_n([\varphi_k])=\Pi_n([\varphi_k]).
\end{equation}
The next step is to write the 
coefficients $a_k(\varphi_{i_1},...,\varphi_{i_k})$ in 
terms of $\Pi_n([\varphi_j])$. To this end, we can use the method 
for finding orthonormal moments of regular polynomial 
functionals (Theorem 1 in \cite{Kashpur}). This yields
\begin{align}
n!a_n(\varphi_{i_1},...,\varphi_{i_n})=&
\Pi_n([\varphi_{i_1}+\cdots +\varphi_{i_n}])-
\left\{\Pi_n([\varphi_{i_1}+\cdots +\varphi_{i_{n-1}}])+\nonumber\right.\\
& \left.\Pi_n([\varphi_{i_1}+\cdots +\varphi_{i_{n-2}}+\varphi_{i_{n}}])+\cdots
+\Pi_n([\varphi_{i_2}+\cdots +\varphi_{i_{n}}])\right\}+\nonumber\\
&\left\{ \Pi_n([\varphi_{i_1}+\cdots +\varphi_{i_{n-2}}])+
\Pi_n([\varphi_{i_1}+\cdots +\varphi_{i_{n-3}}+\varphi_{i_{n-1}}])+\cdots+\right.\nonumber\\
&\left.\Pi_n([\varphi_{i_3}+\cdots +\varphi_{i_{n}}])\right\} \cdots
+(-1)^{n-1}\left\{\Pi_n([\varphi_{i_1}])+\cdots +\Pi_n([\varphi_{i_n}])\right\}+\nonumber\\
& (-1)^{n}\Pi_n([0]).
\label{coeffKL}
\end{align}
Once $a_n(\varphi_{i_1},...,\varphi_{i_n})$ 
are available, we can construct the polynomial 
$$\Pi_n-a_n(\varphi_{i_1},...,\varphi_{i_n})
\prod_{z=1}^n(\theta,\varphi_{i_z})$$
and apply \eqref{coeffKL} again to determine 
$a_{n-1}(\varphi_{i_1},...,\varphi_{i_{n-1}})$.
After $n$ iterations, we have available all 
coefficients to construct the 
polynomial functional \eqref{interp4}.
It was shown in  \cite{Kashpur} (Theorem 2) that 
the following error estimate holds
\begin{equation}
\left\|\Pi^I_n([\theta])-\Pi_n([\theta])\right\|\leq \sum_{k=1}^n\left\|
P_k\right\|\left[\left(\left\|\theta\right\|+\epsilon_m([\theta])\right)^k-
\left\|\theta\right\|^k\right],\label{conver}
\end{equation}
where 
\begin{equation}
 \left\|P_k\right\|=\int\cdots\int \left|K_k(x_1,...,x_k)\right|dx_1\cdots dx_k,
\end{equation}
\begin{equation}
 \epsilon_m([\theta])=\left\|\theta-\sum_{k=1}^m(\theta,\varphi_k)\varphi_k\right\|.
\end{equation}
This means that the interpolant \eqref{interp4} converges 
{\em pointwise} to the polynomial functional 
$\Pi_n([\theta])$ in \eqref{pfun} as the number of nodes $\{\varphi_k\}$ 
in the function space $D(\Pi_n)$ goes to infinity\footnote{Note that 
$\epsilon_m([\theta])\rightarrow 0$ as $m\rightarrow \infty$. 
Thus, the right hand side of the error estimate \eqref{conver} 
goes to zero as $m\rightarrow \infty$, i.e., 
$\Pi^I_n([\theta])\rightarrow \Pi_n([\theta])$ 
pointwise as $m\rightarrow \infty$.}.
As noted by Khlobystov in \cite{Khlobystov0}, 
Porter's interpolant \eqref{Porter_interpolant} can be written 
exactly in the form \eqref{interp4} if we consider test 
functions in the form $\theta_i=\varphi_{i_1}+...+\varphi_{i_n}$ 
(see Lemma 1 in \cite{Khlobystov0}). In this sense,  Porter's 
interpolants represent the Lagrangian form of 
Khlobystov's interpolants. 
An interesting and very important question is the approximation 
of polynomial functionals of order $n$ in terms interpolants 
over nodal sets $\widehat{S}^{(m)}_q$, with $q<n$  (see Eq. \eqref{SNq1}). 
A specific example would be the approximation of $P_3([\theta])$ 
in \eqref{p3} by using an interpolant built on the  
set of nodes $\widehat{S}^{(m)}_1$. 
In some very special cases we may get uniform convergence  
as $m\rightarrow \infty$, for example when $K_3$ is diagonal
\begin{equation}
 K_3(x_1,x_2,x_3)=\sum_{j=1}^m a_j\varphi_j(x_1)\varphi_j(x_2)\varphi_j(x_3).
 \label{almostCP}
\end{equation}
However, this won't happen in general, and therefore
we will have to accept a systematic truncation error in representing 
continuous nonlinear functionals. Such error is similar the 
error committed when we approximate infinite-variate functions, e.g.,  
by second-order HDMR \cite{Wasilkowski1} (see Section 
\ref{sec:HDMR}).

\paragraph{Khlobystov interpolants with Hilbert-Schmidt Kernels} 
In this Section we present a procedure to construct interpolants of 
polynomial operators in Hilbert spaces with separable kernels. 
To this end, let us consider an orthonormal 
basis $\{\varphi_1,\varphi_2,...,\varphi_m\}$ 
and represent each kernel in \eqref{lP} as 
\begin{align}
K_1(x_1)=&\sum_{i=1}^m a_{i}\varphi_{i}(x),\label{k1}\\
K_2(x_1,x_2)=&\sum_{i,j=1}^m a_{ij}\varphi_{i}(x_1)\varphi_{j}(x_2),
\label{k2}\\
K_3(x_1,x_2,x_3)=&\sum_{i,j,k=1}^m a_{ijk}\varphi_{i}(x_1)
\varphi_{j}(x_2)\varphi_{k}(x_3),\label{k3}\\
\cdots .\nonumber&
\end{align}
Without loss of generality we can assume that $K_2$, $K_3$, ...
are symmetric, i.e., that the coefficients $a_{ijk\cdots}$ are invariant under 
any permutation of the indices $i$, $j$, $k$, etc.
A substitution of \eqref{k1}, \eqref{k2}, etc.,  into 
\eqref{pfun} yields the polynomial functional 
\begin{align}
\Pi_n([\theta]) = K_0 
+\sum_{i=1}^m a_i(\varphi_i,\theta)
+\sum_{i,j=1}^m a_{ij}(\varphi_i,\theta)(\varphi_j,\theta)
+\sum_{i,j,k=1}^m a_{ijk}(\varphi_i,\theta)(\varphi_j,\theta)
(\varphi_k,\theta)+\cdots.
\end{align}
At this point we pose the following question: How many interpolation 
nodes do we need to identify the 
coefficients $a_i$, $a_{ij}$, $a_{ijk}$, etc., uniquely? 
To clarify the question, consider the following second-order 
polynomial functional
\begin{align}
\Pi_2([\theta]) = K_0 
+\sum_{i=1}^m a_i(\varphi_i,\theta)
+\sum_{i,j=1}^m a_{ij}(\varphi_i,\theta)(\varphi_j,\theta).
\label{p2interp}
\end{align}
The total number of degrees of freedom (number of independent 
coefficients $K_0$, $a_i$ and $a_{ij}$) is 
\begin{equation}
 1+m+m+\binom{m}{2}=1+\frac{3m+m^2}{2}.
\end{equation}
To determine 
such coefficients, we need to test 
$\Pi_2([\theta])$ at $(2+3m+m^2)/2$ distinct nodes 
\footnote{Note, that testing a second-order polynomial functional 
in $\varphi_i(x)$ or $2\varphi_i(x)$ yields different results.}, e.g.,
\begin{equation}
\left\{0,\{\varphi_i\}_{i=1}^m,\{\varphi_i+
\varphi_j\}_{\substack{i,j=1\, (j\geq i)}}^m\right\}.
\label{test_f_s}
\end{equation}
This yields the linear system
\begin{align}
\Pi_2([0])         &= K_0\nonumber\\
\Pi_2([\varphi_p]) &= K_0+a_p+a_{pp}\nonumber\\
\Pi_2([\varphi_p+\varphi_q]) &= K_0+a_p+a_q+ a_{pp}+
a_{qq}+2a_{pq}\nonumber
\end{align}
which can be immediately solved for $a_{pq}$, $a_p$ and $K_0$
\begin{align}
a_{pq}&=\frac{1}{2}\left(\Pi_2([\varphi_p+\varphi_q])-
\Pi_2([\varphi_p])-\Pi_2([\varphi_q])+\Pi_2([0])\right),\label{apq}\\
a_{p}&=-\frac{1}{2}\Pi_2([2\varphi_p])+2\Pi_2([\varphi_p])-\frac{3}{2}
\Pi_2([0]),\label{ap}\\
K_0  &= \Pi_2([0])\label{k0}.  
\end{align}
In this way, we can identify the kernels 
\eqref{k1}-\eqref{k3} and therefore 
any polynomial functional in the form \eqref{p2interp}. 
Note that if the basis elements in \eqref{test_f_s} are not 
normalized (but still orthogonal) then we simply need to 
replace $a_{pq}$ and $a_p$ in \eqref{apq}-\eqref{k0} 
with $\left\|\varphi_p\right\|^2\left\|\varphi_q\right\|^2a_{pq}$ and 
$\left\|\varphi_p\right\|^2a_{p}$, respectively. 
Higher-order polynomial functionals can be constructed  
in a similar way. However, the number of degrees of freedom 
may increase significantly with the order of the polynomial 
(see equation \eqref{SCardinality}). 
For example, a twelve-order polynomial functional in which 
each kernel is represented relative to $m=10$ basis functions 
(tensor product) yields $646646$ degrees of freedom!

\vs
\noindent
{\em Remark:} 
The fact that we can get analytical expressions for the 
coefficients of the polynomial interpolant means that 
Porter's matrix \eqref{Hint} is {\em invertible analytically} 
for uniquely solvable interpolation problems and orthonormal 
bases.

\vs 
\noindent
The procedure we just described to identify the coefficients 
of the symmetric kernel functions $K_n$ relies on 
tensor product representations, i.e., series expansions in the form 
\begin{align}
K_n(x_1,x_2,...,x_n)=\sum_{i_1,...,i_n=1}^{m}a_{i_1\cdots i_n} 
\varphi_{i_1}(x_1)\cdots \varphi_{i_n}(x_n).
\label{eq:tp}
\end{align}
The number of independent coefficients $a_{i_1\cdots i_n}$  is 
\begin{equation}
\binom{n+m-1}{n}=\frac{(n+m-1)!}{n!(m-1)!}.
\label{dof_sym}
\end{equation}
Such number can pose serious computational challenges even for moderate 
values of $m$ and $n$.
To alleviate this problem one could use 
HDMR expansions \cite{Rabitz,Li,Li1}, i.e., represent  
$K_n(x_1,x_2,\cdots,x_n)$ in terms of a superimposition of 
functions involving a lower number of variables 
(interaction terms). This yields, for example
\begin{align}
K_n(x_1,x_2,...,x_n)=f_{0}+\sum_{i=1}^{n}f_{i}(x_{i})
+\sum_{i<j}^{n}f_{ij}(x_{i},x_{j})+\sum_{i<j<k}^{n}f_{ijk}(x_{i},x_{j},x_k)
+ \cdots\,.
\label{eq:anova}
\end{align}
The function $f_{0}$ is a constant. The functions $f_i(x_i)$ 
(first-order interactions) give us the overall effects of the 
variables $x_i$ in $f$ as if they were acting 
independently of the other input variables. The functions $f_{ij}(x_i, x_j)$ 
(second-order interactions) describe the cooperative effect 
of the variables $x_i$ and $x_j$. Similarly, higher-order terms 
reflect the cooperative effects of an increasing number of variables. 
Representing $f_{i}$, $f_{ij}$, $f_{ijk}$ 
relative to the orthonormal basis $\{\varphi_j(x)\}$ 
yields the series
\begin{align}
K_n(x_1,x_2,...,x_n)=& f_{0}+\sum_{i=1}^{n}\sum_{s=1}^m a^{i}_s
\varphi_s(x_{i})
+\sum_{i<j}^{n}\sum_{s,q=1}^ma^{ij}_{sq}\varphi_s(x_{i})
\varphi_q(x_{j}) + \nonumber \\
&\sum_{i<j<k}^{n} \sum_{s,q,r=1}^m a^{ijk}_{sqr}\varphi_s(x_{i})
\varphi_q(x_{j})
\varphi_r(x_{j})+\cdots.
\label{eq:anova1}
\end{align}
Given the symmetry of each function $f_{i}$, $f_{ij}$, $f_{ijk}$, the total 
number of degrees of freedom of an HDMR expansion of order $Z$ is
\begin{equation}
\sum_{k=0}^Z \binom{n}{k} \binom{m+k-1}{k} = 1+\binom{n}{1}
\binom{m}{1} +\binom{n}{2}\binom{m+1}{2} + 
\binom{n}{3}\binom{m+2}{3} +\cdots. 
\end{equation}
For example, the second- and third-order 
truncations of a $n=10$ dimensional 
kernel relative to $m=10$ yield $2576$ and $28976$ 
degrees of freedom, respectively. On the other hand, 
the tensor product representation yields $92378$ 
degrees of freedom. 
{\color{r}
Alternatively, one can use tensor expansions 
(see Section \eqref{sec:tensor}), e.g., canonical polyadic or 
hierarchical-Tucker series, of each kernel to 
reduce the number of degrees of freedom. 
The tensor expansion can be fit to data 
by interpolation, least-squares or 
projection \cite{Nouy2017,NouyHUQ,Nouy2,Nouy}.
}

{\color{r}

\subsection{Functional Approximation by Tensor Methods}
\label{sec:tensor}

Computing high-order polynomial functional expansions 
requires representing kernel functions in many 
independent variables. To get an feeling 
of how serious this problem could be, simply consider that 
representing a polynomial functional of order $6$ is as computationally 
expensive as representing the solution to the 
steady-state Boltzmann equation \cite{Dimarco}, a well-known 
challenging problem in computational physics.
Expanding each kernel of the polynomial functional 
in terms of HDMR \cite{Rabitz} or canonical polyadic 
decompositions -- i.e. separated series expansions  \cite{Beylkin} -- 
can mitigate the dimensionality problem, but 
it may not be the most efficient way to proceed. 
In this Section we discuss nonlinear functional approximation 
based on tensor methods.  To introduce these methods, consider 
the Hilbert space of functions spanned by the finite-dimensional 
basis $\{\varphi_1,...,\varphi_m\}$ (e.g., an orthonormal basis) 
and look for an approximant of $F([\theta])$, say $f$, 
in the form
\begin{equation}
F([\theta])= f(a_1([\theta]), ..., a_m([\theta]))
+R([\theta]),\qquad a_k([\theta])=(\theta,\varphi_k).
\label{functional-approx}
\end{equation}
In this equation, $f$ is a multivariate function of 
$a_k([\theta])$ (linear functionals of $\theta$), $(,)$ denotes an 
inner product in the Hilbert Space $D(F)$, e.g., 
\begin{equation}
(\theta,\varphi_k)=\int_a^b \theta(x)\varphi_k(x)dx,
\qquad k=1,...,m 
\end{equation}
and $R$ is a (functional) reminder term.
The functionals $a_k([\theta])=(\theta,\varphi_k)$ can be either real or 
complex-valued. In the theory of stochastic processes 
the set 
\begin{equation}
 \left\{\theta\in D(F) \,|\, ((\theta,\varphi_1),...,(\theta,\varphi_m))\subseteq\mathbb{R}^m\right\}
\end{equation}
is known as {\em cylindrical set} (see, e.g., \cite{Neerven} p. 56 
or \cite{Skorokhod} p. 45). Therefore, according to 
Eq. \eqref{functional-approx}, we are looking for an approximant 
of $F([\theta])$ in the space of {\em cylindrical functionals}, i.e., 
functionals defined on cylindrical sets.  
Thanks to the Stone-Weierstrass theorem, 
cylindrical functionals in the form \eqref{cf} can approximate 
any continuous functional in a Hilbert or a Banach space.
The representation \eqref{functional-approx} 
is very general. For example, it includes the 
case where $f$ is a polynomial functional, 
e.g.,  \eqref{interp4} or \eqref{Porter_interpolant}, 
or the case where the functional is defined 
on a finite-dimensional function space (see Section \ref{sec:Finite_Dim_Approx}).
For notational convenience we will often drop the 
the functional dependence of $a_k([\theta]$ 
and write \eqref{functional-approx} as
\begin{equation}
F([\theta]) \simeq f(a_1,...,a_m),\qquad a_j=(\theta,\varphi_j).
\label{cf}
\end{equation} 
Hereafter, we discuss effective numerical 
algorithms to compute an approximation 
of the multivariate function $f(a_1,....,a_m)$ based on 
{\em high-dimensional model representations (HDMR)}, and {\em tensor 
methods}\footnote{If we ask the question ``what is a tensor?'' 
to an engineer, a physicist or a mathematician we usually get different 
answers. The engineer usually  says  ``a tensor is a matrix''.
On the other hand, the physicist would say that a tensor is a 
mathematical object that has the fundamental property 
of transforming in a very specific way when the coordinate 
system is changed. He or she would point out that 
the laws of physics are built upon the {\em principle of general 
covariance} \cite{Weinberg,VenturiJPA2009} 
(invariance of physical laws 
relative to coordinate transformations)  which is formulated in 
a natural way in terms of tensors. The word ``tensor'' 
has recently spilled in the multi-linear algebra community 
to represent multi-dimensional arrays. Hereafter we will adopt 
such terminology, and refer to tensors as multi-dimensional arrays.}
\cite{Kolda,Hackbusch_book}, including 
canonical tensor decompositions, hierarchical Tucker formats,
and tensor networks.
Such algorithms rely on optimization (e.g., the alternating 
least squares methods \cite{Acar,Reynolds}), or multilinear 
algebra techniques such as high-order singular value decomposition \cite{Grasedyck2017}, randomized block sampling, 
or generalized Schur decompositions.

\paragraph{Functional Derivatives}
Let us compute the functional derivatives of the cylindrical 
functional approximation \eqref{cf}. To this end, we evaluate 
the G\^ateaux differential of $f(a_1([\theta]), ..., a_m([\theta])$ 
in the direction of an arbitrary function $\eta(x)$ to obtain 
\begin{align}
\int_a^b \frac{\delta F([\theta])}{\delta \theta(x)}\eta(x)dx&\simeq 
\frac{d}{d\epsilon}\left[f((\theta+\epsilon\eta,\varphi_1),...,
(\theta+\epsilon\eta,\varphi_m))\right]_{\epsilon=0}\\
&=\sum_{k=1}^m\frac{\partial f}{\partial a_k}(\eta,\varphi_k).
\end{align}
Setting $\eta(x)=\delta(x-y)$ yields 
the following approximation for the first-order functional derivative
\begin{align}
\frac{\delta F([\theta])}{\delta \theta(y)}\simeq \sum_{k=1}^m
\frac{\partial f}{\partial a_k}\varphi_k(y),
\label{ff1}
\end{align}
where $a_k=(\theta,\varphi_k)$. On the other hand, 
by projecting \eqref{ff1} onto $\varphi_j$
we obtain
\begin{equation}
 \int_a^b \frac{\delta F([\theta])}{\delta \theta(x)}\varphi_j(x)dx\simeq
 \frac{\partial f}{\partial a_j}.
 \label{ff2}
\end{equation}
This means that partial derivative of $f$ relative to 
$a_k=(\theta,\varphi_k)$ approximates 
the projection of the functional derivative
of $F$ along $\varphi_k$. Equations \eqref{ff1} and \eqref{ff2} 
are consistent with previous results on functional derivatives 
in finite-dimensional spaces (see Eqs. \eqref{gg1} and \eqref{gg2}).
By following the same procedure, we can construct functional 
derivatives of $f$ of higher-order, e.g., 
\begin{align}
\frac{\delta^2 F([\theta])}{\delta \theta(x)\delta \theta(y)}\simeq&\sum_{k,j=1}^m
\frac{\partial^2 f}{\partial a_k\partial a_j}\varphi_k(x)\varphi_j(y),
\label{ff5}\\
\frac{\delta^3 F([\theta])}{\delta \theta(x)\delta \theta(y)\delta \theta(z)}
\simeq
&\sum_{k,j,i=1}^m
\frac{\partial^3 f}{\partial a_k\partial a_j\partial a_i}
\varphi_k(x)\varphi_j(y)\varphi_i(z).
\label{ff6}
\end{align}

\paragraph{Choice of the Number of Active Dimensions}
The choice of the basis functions $\varphi_k$ and 
the number of active dimensions, i.e., the 
parameter $m$ in \eqref{functional-approx}, 
is {\em critical} for the accuracy of the cylindrical representation. 
For a fixed basis $\{\varphi_j\}$, smaller values of $m$ 
yield faster computations but at the same time can 
lead to functional approximation problems with poor 
approximation errors.

\subsubsection{Canonical Tensor Decomposition}
\label{sec:CP}
The canonical tensor decomposition of the multivariate 
function $f(a_1,...,a_m)$ in \eqref{functional-approx} 
is a series expansion in the form 
\begin{equation}
f(a_1,...,a_m)=\sum_{l=1}^r
\prod_{i=1}^mG^l_i(a_i)\qquad a_j=(\theta,\varphi_j),
\label{functional-SSE}
\end{equation}
where $G^l_i(a_i)$ are one-dimensional 
functions usually represented relative to a 
known basis $\{\phi_1,...,\phi_Q\}$, i.e., 
\begin{equation}
G^l_i(a_i)=\sum_{k=1}^Q\beta_{ik}^l \phi_k(a_i).
\label{gfun}
\end{equation}
The quantity $r$ in \eqref{functional-SSE} is 
called {\em separation rank}.
In the statistics literature, representations of 
the form \eqref{functional-SSE}
are known as {\em parallel factorizations} (see 
\cite{Kronenberg,Leugarans}). 
They are also known as proper generalized 
decompositions \cite{ChinestaBook}, 
canonical polyadic expansions (CP) \cite{Karlsson}, 
separated series \cite{HeyrimJCP_2014}, 
and Kruskal tensor formats \cite{Kolda}.
Although there are at present no useful theorems on the size of 
the separation rank $r$ needed to represent with accuracy 
general classes of functionals $f$, there are cases where the 
expansion \eqref{functional-SSE} is {\em exponentially 
more efficient} than one would expect a priori. 
The basis 
functions $\phi_k$ appearing in \eqref{gfun} represent the 
variability of the functional $f$ along different directions 
$a_k=(\theta,\varphi_k)$ in the test function space $D_m$.
As such, they have to satisfy appropriate boundary conditions. 
For example, if $f$ is periodic in the hypercube $[-b,b]^m$ 
then we could use rescaled  trigonometric polynomials
\begin{equation}
\phi_k(a)=l_k\left(\pi\left(\frac{a}{b}+1\right)\right)\qquad a\in[-b,b],
\label{Fourier}
\end{equation}
where 
\begin{equation}
l_{k+1}(x)=\frac{1}{Q}\displaystyle
\frac{\sin\left(\displaystyle\frac{Q}{2}\left(x-x_k\right)\right)}
{\sin\left(\displaystyle\frac{x-x_k }{2} \right)},
\qquad x_k=\frac{2\pi}{Q}k\qquad k=0,...,(Q-1).
\end{equation}
For more general boundary conditions we can use a polynomial basis, 
e.g., rescaled Legendre orthogonal polynomials
\begin{equation}
\phi_k(a)=L_k\left(\frac{a}{b}\right)\qquad a\in[-b,b],
\label{Legendre}
\end{equation}
where 
\begin{align}
L_0(x)=1, \quad
L_1(x)=x,\quad
\cdots, \quad 
L_{n+1}(x)=\frac{2n+1}{n+1}xL_{n}(x)-\frac{n}{n+1}L_{n-1}(x), \quad  
(n=1,...,Q-1).
\end{align}
The $L_2$ norm of \eqref{Fourier} and \eqref{Legendre} is easily 
obtained as
\begin{equation}
\eta_k= \frac{b}{\pi}\frac{2\pi}{Q}\quad \textrm{(Fourier)},\qquad 
\eta_k= \frac{2b}{2h+1} \quad \textrm{(Legendre)}, \qquad 
k=1,...,Q. 
\end{equation}
}

\vs
\noindent
{\em Example 1:} Consider the sine functional 
\eqref{sinFunctional}, hereafter rewritten for convenience
\begin{equation}
 F([\theta])=\sin\left(\int_a^b K(x)\theta(x)dx\right).
 \label{sinFunctional1}
\end{equation}
Assuming that the kernel $K(x)$ admits the 
finite-dimensional expansion
\begin{equation}
 K(x)=\sum_{j=1}^m k_j\varphi_j(x),
 \label{kkk}
\end{equation}
and substituting it into \eqref{sinFunctional1} we obtain 
\begin{align}
f(a_1,...,a_m)=&\sin\left(\sum_{j=1}^m k_j a_j\right)\nonumber \\
=&\sum_{j=1}^m\sin(k_ja_j)
\prod_{\substack{i=1\\i\neq j}}^m \frac{\sin(k_ia_i+\chi_i-\chi_j)}
{\sin(\chi_i-\chi_j)}, 
\end{align}
for all $\chi_i$ such that $\chi_i\neq \chi_j$. The last equality was 
derived in \cite{Mohlenkamp}, and it shows that the separation rank 
of the canonical tensor decomposition of 
\eqref{sinFunctional1}-\eqref{kkk} is exactly $r=m$. In other words, 
we can represent the nonlinear functional \eqref{sinFunctional1} 
exactly in terms of superimposition of $m$ terms. 
Furthermore, if we allow the $G^l_i(a_i)$ 
to be complex-valued\footnote{Constraints on the 
functions $G^l_i$ such as positivity can be also imposed, e.g.,  
if one is interested in probability functionals.} 
then
\begin{equation}
\sin\left(\sum_{j=1}^m k_ja_j\right)=
\prod_{j=1}^m\frac{e^{ik_ja_j}}{2i}-
\prod_{j=1}^m\frac{e^{-ik_ja_j}}{2i},
\end{equation}
i.e., we can reduce the separation rank to $r=2$. 
In general, the separation rank depends on the 
complexity of the nonlinear function $f(a_1,...,a_m)$.

\paragraph{Functional Derivatives}
With the canonical tensor decomposition \eqref{functional-SSE}
available, it is straightforward to compute an 
approximation of the functional derivatives 
of $F([\theta])$. Recalling that the canonical tensor 
decomposition is a cylindrical representation of the functional 
$F([\theta])$, we have the expressions 
\eqref{ff1}, \eqref{ff5} and \eqref{ff6}. 
The unknowns are the partial derivatives 
of $f$ with respect to $a_j([\theta])$, which can be 
computed based on \eqref{functional-SSE} as
\begin{equation}
\frac{\partial f}{\partial a_i}=
\sum_{l=1}^r\frac{\partial G^l_i}{\partial a_i}
\prod_{\substack{k=1\\k\neq i}}^m G_k^l(a_k),
\label{pd1}
\end{equation}
and
\begin{equation}
\frac{\partial f}{\partial a_i\partial a_j}=
\begin{cases}
\displaystyle\sum_{l=1}^r\frac{\partial^2 G^l_i}{\partial a_i^2}
\prod_{\substack{k=1\\k\neq i}}^m G_k^l(a_k) & i=j,\\
\displaystyle\sum_{l=1}^r\frac{\partial G^l_i}{\partial a_i}
\frac{\partial G^l_j}{\partial a_j}\prod_{\substack{k=1\\k\neq i,j}}^m G_k^l(a_k) 
& i\neq j.
\end{cases}
\label{pd2}
\end{equation}
These derivatives are evaluated at $a_j([\theta])=(\theta,\varphi_j)$.

\paragraph{Alternating Least Squares (ALS) Formulation}

The development of robust and efficient algorithms to 
compute \eqref{functional-SSE} to any desirable accuracy is 
still a relatively open question (see
\cite{Acar,Espig1,Karlsson,DoostanIaccarino2013,HeyrimJCP_2014} 
for recent progresses). Computing the tensor 
components $G_k^l(a_k)$ usually relies on (greedy) 
optimization techniques such as alternating least 
squares (ALS) \cite{Reynolds,Battaglino,Acar,Beylkin} 
or regularized Newton methods \cite{Espig1}, which are only locally 
convergent \cite{Uschmajew2012} (i.e., the final result may depend 
on the initial condition of the algorithm). Hereafter we describe 
the simplest form of the ALS algorithm. 
To this end, consider the functional residual
\eqref{functional-approx}, with $f$ 
given in \eqref{functional-SSE}
\begin{equation}
R([\theta])=F([\theta])-\sum_{l=1}^r 
\prod_{k=1}^m G^l_k((\theta,\varphi_k)).
\label{Fresidual}
\end{equation}  
We aim at computing the tensor components $G_k^l(a_k)$
by minimizing the norm of $R([\theta])$ relative to 
independent variations of $G^l_k$. 
In particular, if we assume that $G^l_k$ are 
in the form \eqref{gfun}, then variations 
of $G^l_k$ are generated by variations of $\beta_{hj}^n$.
Therefore, the canonical tensor decomposition of $F([\theta])$
can be computed by the variational principle\footnote{The Euler-Lagrange equations associated with \eqref{als} are nonlinear $\beta_{hj}^n$.}
\begin{equation}
 \min_{\beta_{hj}^n}\left\|R([\theta])\right\|_W^2.
 \label{als}
\end{equation}
The norm in $\|\cdot \|_W$ may be defined by a weighted 
functional integral (see Appendix \ref{sec:functional integrals}) 
in the form 
\begin{equation}
\left\| R([\theta])\right\|^2_W=\int 
R([\theta])^2 W([\theta])\mathcal{D}[\theta],
\label{continuous_norm_W}
\end{equation}
where $W([\theta])$ is the functional integration measure, or by a 
discrete sum (functional collocation setting) 
\begin{equation}
\left\|R([\theta])\right\|^2_W=
\sum_{j=1}^N R([\theta_j])^2 w_j,
\label{discrete_norm_W}
\end{equation}
where $\{\theta_1, ..., \theta_N\}$ are collocation nodes 
in the function space $D(F)$, and $w_j$ are 
integration weights.
{\color{r}
In the alternating least squares paradigm, 
we compute the minimizer of the residual \eqref{Fresidual}
by splitting the non-convex optimization problem \eqref{als} into 
a {\em sequence} of convex low-dimensional optimization 
problems (see Figure \ref{fig:ALSfigure}).
To illustrate the method, let us define the vectors 
\begin{align}
\bm \beta_1&=[\beta^1_{11}, ..., \beta^1_{1Q},...,\beta^r_{11}, ..., \beta^r_{1Q}]^T,\nonumber\\
\bm \beta_2&=[\beta^1_{21}, ..., \beta^1_{2Q},...,\beta^r_{21}, ..., \beta^r_{2Q}]^T,\nonumber\\
\vdots &\nonumber\\
\bm \beta_m&=[\beta^1_{m1}, ..., \beta^1_{mQ},...,\beta^r_{m1}, ..., \beta^r_{mQ}]^T. 
\end{align}
Note that $\bm \beta_i$ collects the degrees of freedom 
of all functions $\{G^1_{i}, ..., G^r_i\}$ depending 
on $a_i=(\theta,\varphi_i)$. Next, we split the optimization 
problem \eqref{als} into the following sequence of 
convex optimization problems 
\begin{equation}
\min_{\bm\beta_1}\left\| R([\theta])\right\|^2_W, \qquad 
\min_{\bm\beta_2}\left\| R([\theta])\right\|^2_W, \qquad\cdots, \qquad 
\min_{\bm\beta_m}\left\| R([\theta])\right\|^2_W.
\label{als11}
\end{equation}
We emphasize that the system of equations 
\eqref{als11} is not equivalent to the full problem \eqref{als}. 
In other words, the sequence of low-dimensional 
optimization problems \eqref{als11} in general 
does not allow us to compute the minimizer of \eqref{als} \cite{Espig,Bezdek,Rohwedder,Uschmajew2012}. 
\begin{figure}[t]
\centerline{\includegraphics[height=6cm]{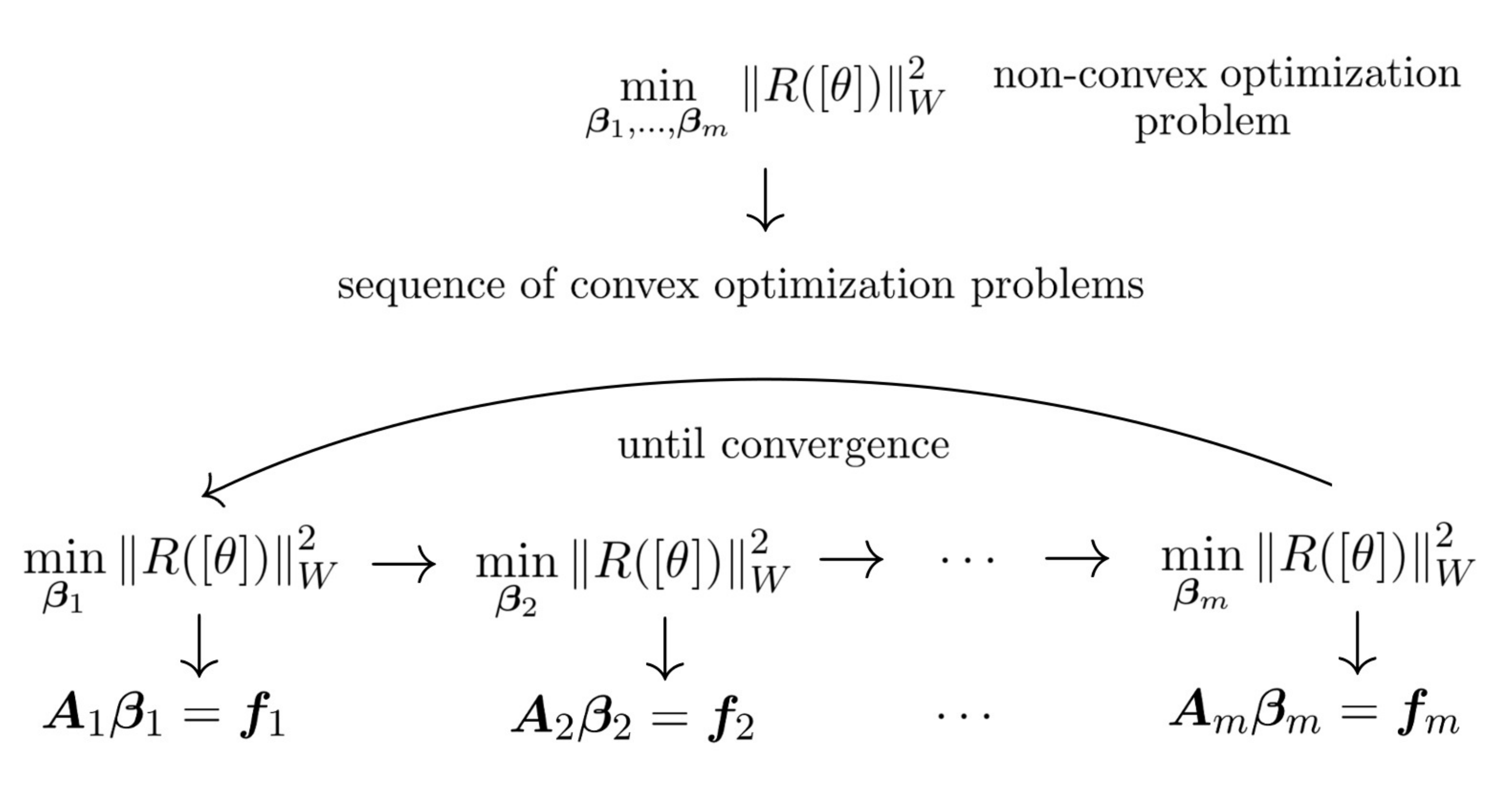}}
\caption{\color{r}Sketch of the alternating least squares (ALS) algorithm for 
the minimization of the functional residual $R([\theta])$. 
Convergence analysis of ALS  be found in 
\cite{Espig,Bezdek,Rohwedder,Uschmajew2012}.}
\label{fig:ALSfigure}
\end{figure}
The Euler-Lagrange equations associated 
with \eqref{als11} are in the form
\begin{equation}
\bm A_j \bm \beta_j =\bm f_j\qquad j=1,...,m, 
\label{SSYSTEM}
\end{equation}
where
\begin{equation}
\bm A_j=\left[
\begin{array}{cccc}
A^{11}_{j11}\cdots A^{11}_{j1Q} & A^{12}_{j11}\cdots A^{12}_{j1Q} 
& \cdots &A^{1r}_{j11}\cdots A^{1r}_{j1Q}\\
\vdots & \vdots & \vdots & \vdots \\
A^{11}_{jQ1}\cdots A^{11}_{jQQ} & A^{12}_{jQ1}\cdots A^{12}_{jQQ} 
&\cdots & A^{1r}_{jQ1}\cdots A^{1r}_{jQQ}\\
\vdots & \vdots & \vdots & \vdots \\
A^{r1}_{j11}\cdots A^{r1}_{j1Q} & A^{r2}_{j11}\cdots A^{r2}_{j1Q} 
&\cdots & A^{rr}_{j11}\cdots A^{rr}_{j1Q}\\
\vdots & \vdots & \vdots & \vdots \\
A^{r1}_{jQ1}\cdots A^{r1}_{jQQ} & A^{r2}_{jQ1}\cdots A^{r2}_{jQQ} 
&\cdots & A^{rr}_{jQ1}\cdots A^{rr}_{jQQ}
\end{array}
\right],\qquad 
\bm f_j=\left[
\begin{array}{c}
 f_{j1}^1\\
 \vdots\\
f_{jQ}^1\\
 \vdots\\
 f_{j1}^r\\
 \vdots\\
f_{jQ}^r\\ 
\end{array}
\right],
\label{ORDERING}
\end{equation}
and 
\begin{equation}
A^{ln}_{jhs}=\int\mathcal{D}[\theta]W([\theta])
\phi_h((\theta,\varphi_j))\phi_s((\theta,\varphi_j))\prod_{\substack{k=1\\k\neq j}}^m
G_k^l((\theta,\varphi_k))G_k^n((\theta,\varphi_k)),
\label{A^{ln}_{jhp}}
\end{equation}
\begin{equation}
f_{jh}^n = \int\mathcal{D}[\theta]W([\theta])F([\theta])
\phi_h((\theta,\varphi_j))\prod_{\substack{k=1\\k\neq j}}^m
G_k^n((\theta,\varphi_k)).
\label{f_{jh}^n}
\end{equation}

\noindent
The matrices $\bm A_j$ are symmetric, positive 
definite and of size $rQ\times rQ$. Also,  the functional integrals 
defining the matrix entries can be simplified and eventually 
computed by using techniques for high-dimensional integration, 
such as the quasi-Monte Carlo method \cite{Dick}.
Indeed, if we restrict the residual \eqref{Fresidual} 
to the finite-dimensional function space
$D_m=\textrm{span}\{\varphi_1,...,\varphi_m\}$ (Section 
\ref{sec:Finite_Dim_Approx}), and assume that $F([\theta])$ 
restricted to $D_m$ is compactly supported within 
the hypercube $[-b,b]^m$, then $A^{ln}_{jhp}$ and $f_{jh}^n$ 
can be written in the form
\begin{equation}
A^{ln}_{jhs}=\int_{-b}^b \phi_h(a_j) \phi_s(a_j)da_j
\prod_{\substack{k=1\\k\neq j}}^m
\int_{-b}^b G_k^l(a_k)G_k^n(a_k)da_k,
\label{Alnjhp1}
\end{equation}
\begin{equation}
f_{jh}^n = \int_{-b}^b \cdots \int_{-b}^b f(a_1,...,a_m)
\phi_h(a_j) \prod_{\substack{k=1\\k\neq j}}^m
G_k^n(a_k) da_1\cdots da_m,
\label{fh1}
\end{equation}
provided we select the integration measure $W([\theta])$ 
appropriately (see Appendix \ref{sec:functional integrals}).

}

\vs
\noindent 
{\em Remark:}
Minimizing the residual \eqref{Fresidual}
with respect to $\beta_{hj}^n$ 
is equivalent to imposing orthogonality 
relative to the space spanned by the functionals
\begin{equation}
\phi_h((\theta,\phi_j))\prod^m_{\substack{k=1\\k\neq j}}G^{n}_k((\theta,\varphi_k))
\end{equation}
To show this in a simple and intuitive way, consider the following  
example in just two dimensions. Let $f(a_1,a_2)$ be a regular function defined on the unit square $[0,1]^2$. We look for a canonical tensor 
decomposition of $f(a_1,a_2)$ in the form 
\begin{equation}
f(a_1,a_2)=\sum_{j=1}^r G^{j}_1(a_1)G^{j}_2(a_2), 
 \quad \textrm{where}\quad G^{j}_k(a_k)=\sum_{p=1}^P \beta^j_{kp}\phi_p(a_k).
\end{equation}
To determine $\beta^j_{kp}$ we first define the residual, 
\begin{equation}
 R(a_1,a_2)= f(a_1,a_2)-\sum_{j=1}^r G^{j}_1(a_1)G^{j}_2(a_2)
\end{equation}
and then minimize its weighed norm 
\begin{equation}
\left\| R\right\|^2_w=\int_0^1\int_0^1R(a_1,a_2)^2 w(a_1,a_2)da_1da_2 
\end{equation}
relative to independent variations of $\beta^j_{1p}$ and $\beta^j_{2p}$. This yields 
\begin{equation}
\begin{array}{ll} 
\displaystyle \int_0^1\int_0^1 w(a_1,a_2) R(a_1,a_2)\phi_p(a_1)G_2^{j}(a_2)da_1da_2=0 &\quad \textrm{linear system for $\beta^j_{1p}$},\vspace{0.2cm}\\
\displaystyle  \int_0^1\int_0^1 w(a_1,a_2) R(a_1,a_2)G_1^{j}(a_1)\phi_p(a_2)da_1da_2=0 &\quad \textrm{linear system for $\beta^j_{2p}$.}
\end{array}
\label{orthogR}
\end{equation}
Thus, minimizing the residual with respect to independent 
variations of $\beta^j_{1p}$ and $\beta^j_{2p}$ is equivalent impose 
Galerkin orthogonality relative to a space spanned by the basis 
functions $\phi_p(a_1)G_2^{j}(a_2)$ and $G_1^{j}(a_1)\phi_p(a_2)$, 
respectively.

{
\color{r}
\paragraph{Convergence of the ALS Algorithm} 
The ALS algorithm we just described is  
an alternating optimization scheme, i.e.,  
a nonlinear block Gauss--Seidel 
method (\cite{Ortega}, \S 7.4). 
There is a well--developed local 
convergence theory for this type of method 
(see \cite{Ortega,Bezdek}). 
In particular, it can be shown that ALS is locally 
equivalent to the linear block Gauss--Seidel 
iteration applied to the Hessian matrix. 
This implies that ALS is linearly convergent in the iteration 
number \cite{Uschmajew2012}, 
provided that the Hessian of the residual 
is positive definite (except on a trivial null 
space associated with the scaling non-uniqueness 
of the canonical tensor decomposition). The last assumption 
may not be always satisfied\footnote{It was shown in 
\cite{Uschmajew2012} that the classical alternating least squares
algorithm does not converge in the iteration number for functionals 
in the form \eqref{sinFunctional1}.}. Therefore, convergence 
of the ALS algorithm cannot be granted in general.
Another potential issue of the ALS algorithm is 
the poor conditioning of the matrices $\bm A_j$ 
in \eqref{SSYSTEM}, which can addressed by 
regularization \cite{Reynolds,Battaglino}.
The canonical tensor decomposition 
\eqref{functional-SSE} in $m$ dimensions has relatively 
small memory requirements. In fact, the number of degrees of freedom 
that we need to store is $r\times m \times Q$, where $r$ is 
the separation  rank, and $Q$ is the number of degrees of 
freedom employed in each tensor component $G_k^l(a_k)$ \eqref{gfun}.
Despite the relatively low-memory requirements, it is often desirable 
to employ scalable parallel versions the ALS algorithm 
\cite{Karlsson} to compute the canonical 
tensor expansion \eqref{functional-SSE}.
}

{
\color{r}

\subsubsection{Tucker Decomposition}
The Tucker decomposition of the cylindrical functional 
\eqref{cf} is a series expansion in the form 
\begin{equation}
f(a_1,...,a_m)=\sum_{l_1=1}^{r_1}\cdots 
\sum_{l_m=1}^{r_m} C_{l_1\cdots l_m}
T^{l_1}_{1}(a_1)\cdots T^{l_m}_{m}(a_m), 
\label{Tucker}
\end{equation}
where $C_{l_1...l_m}$ is a $r_1\times \cdots \times r_m$ 
real- or complex-valued tensor -- often referred to 
as {\em core tensor} \cite{Kolda} --  
and $T^{l_m}_{m}(a_m)$ are unknown functions.  
Tucker decomposition is known 
as {\em high-order Schmidt decomposition} in the 
context of quantum mechanics \cite{Carteret}.
It is important to emphasize such decomposition is, in general,  
are {\em not unique}\footnote{The classical Schmidt decomposition,
i.e., the bi-orthogonal decomposition of bi-variate functions is 
not unique either, and defined up to two rotations in Hilbert spaces 
\cite{Venturi2,Peres}}. 
As pointed out by Kolda and Bader in \cite{Kolda}, 
this freedom opens the possibility to choosing transformations 
that simplify the core structure in some way so that most 
of the elements of the core tensor $C_{l_1\cdots  l_m}$ are 
zero, thereby eliminating interactions between corresponding 
components. Diagonalization of the core is, in general, 
impossible\footnote{The canonical tensor decompositon 
\eqref{functional-SSE} is in the form of a fully diagonal 
high-order Schmidt decomposition, i.e, 
\begin{equation}
f(a_1,...,a_m)=\sum_{l=1}^{r} C_{l\cdots l}
G^{l}_{1}(a_1)\cdots G^{l}_{m}(a_m).
\end{equation}
The fact that diagonalization of $C_{l_1\cdots  l_m}$ is, in general, 
impossible in dimension larger than 2 implies that it is impossible 
to compute canonical tensor decompositions by standard 
linear algebra techniques. Indeed, the best 
low-rank approximation problem is {\em ill-posed} for real 
tensors with dimension $m>2$ (see, e.g., \cite{Silva,Kolda,Hillar}), 
and for complex tensors \cite{Vannieuwenhoven}.} 
\cite{Peres}, but it is possible to 
try to make as many elements either zero or as small 
as possible  (see, e.g., \cite{Carteret} or \cite{Moravitz}). 
For general tensors we have that the multilinear rank 
$r=r_j$ ($j=1,...,m$) is upper semi-continuous, i.e., 
the Tucker expansion is closed. 
Several efficient algorithms are currently available to compute 
the series expansion \eqref{Tucker}. 
For instance, Lathauwer {\em et al.} 
proposed in \cite{Lathauwer} a high-order singular value 
decomposition method to determine the 
components $T_k^{l_k}$ and the core tensor 
$C_{l_1\cdots  l_m}$ in a discrete setting. 
Such algorithm is simple, robust and it yields 
{\em quasi-optimal} low-rank 
approximations. 

To illustrate the procedure to compute the Tucker decomposition
let us first assume that the basis functions $T_k^{l_k}(a_k)$
in \eqref{Tucker} are {\em orthonormal} and {\em known}. 
In this case, the expansion \eqref{Tucker} 
is simply a  tensor product representation of a 
multivariate function relative to the basis $T_k^{l_k}$. 
The core tensor $C_{l_1\cdots  l_m}$ can be immediately 
obtained by projecting $f(a_1,...,a_m)$ onto the 
basis $T_k^{l_k}(a_k)$, i.e., 
\begin{equation}
C_{l_1\cdots  l_m}=\int  f(a_1([\theta]),...,a_m([\theta])) 
T_1^{l_1}(a_1([\theta]))
\cdots T_m^{l_m}(a_m([\theta]))W([\theta]) \mathcal{D}([\theta]).
\label{FI9}
\end{equation}
Evaluating the functional integral 
in $D_m=\textrm{span}\{\varphi_1,...,\varphi_m\}$
and rescaling the integration measure $W([\theta])$ 
properly (see Appendix \ref{sec:functional integrals})
allows us to rewrite \eqref{FI9}
\begin{equation}
C_{l_1\cdots  l_m}=\int_{-b}^b \cdots \int_{-b}^b f(a_1,...,a_m) 
T_1^{l_1}(a_1)\cdots T_m^{l_m}(x_m)da_1\cdots da_m, 
\end{equation}
where we assumed $f$ to be compactly supported in $[-b,b]^m$.
Next, suppose that each function $T_k^{l_k}(a_k)$ is a 
linear combination of $Q$ orthonormal basis functions 
$\phi_j$, i.e., 
\begin{equation}
T_k^{l_k}(a_k)=\sum_{s=1}^Q \alpha_{sk}^{l_k} \phi_{s}(a_k)
\qquad \textrm{(fixed $k$).}
\end{equation} 
If $\phi_{s}$ is a {\em cardinal basis} associated with 
a set of interpolation nodes along $a_k$, then the matrix 
$\alpha_{sk}^{l_k}$ represents the set of 
functions $T_k^{l_k}(a_k)$ for fixed $k$. We can sort arrange 
the matrix $\alpha_{sk}^{l_k}$ in a way where 
the $l_k$-column represents the value 
of $T_k^{l}(a_k)$ at $Q$ collocation nodes 
along $a_k$. This yields the matrix with entries $[T_k]_{il}$ 
($i=1,...,Q$). If we evaluate the multivariate field $f(a_1,..,a_m)$ 
at the same collocation nodes we employed to 
construct the interpolants of $T_k^{l_k}$, then we can 
rewrite \eqref{Tucker} in a full tensor notation as 
\begin{equation}
f_{i_1\cdots i_m}=\sum_{l_1=1}^{r_1}\cdots 
\sum_{l_m=1}^{r_m} C_{l_1,...,l_m}
[T_1]_{i_1l_1 }\cdots [T_m]_{i_ml_m }, 
\label{Tucker_matrix}
\end{equation}
where $i_j=1,...,Q$ are indices identifying the interpolation node 
along the axis $a_j$. 

\vs
\noindent
{\em Remark:} 
The expansion \eqref{Tucker_matrix}
is a ``matricization'' of continuum series \eqref{Tucker} obtained by 
representing each basis element in terms of collocation nodes.
Clearly, we can also set up a matricization of \eqref{Tucker} 
based on Galerkin projection. To this end, it is sufficient 
to project both the left and the right hand sides 
of \eqref{Tucker} onto 
the orthonormal basis elements
$\phi_{i_1}(a_1)\phi_{i_2}(a_2)\cdots \phi_{i_m}(a_m)$ 
to obtain an expression in the form \eqref{Tucker_matrix}. 
The meaning of $f_{i_1 \cdots i_n}$ in the case is the 
Fourier coefficients of the projection, i.e., 
\begin{equation}
f_{i_1 \cdots i_n}=\int_{-b}^b\cdots\int_{-b}^b f(a_1,...,a_m)
\phi_{i_1}(a_1)\cdots \phi_{i_m}(a_m)dx_1\cdots dx_m.
\end{equation}
Thus, \eqref{Tucker_matrix}  represents a multivariate
expansion of Fourier coefficients in a Tucker tensor format. 
In such finite-dimensional setting, we basically transformed 
the problem of decomposing a multivariate function  to 
a {\em multi-linear algebra problem}. It is immediate to 
see that the discrete Tucker format 
a two-dimensional function $f(a_1,a_2)$ is 
\begin{equation}
f_{i_1i_2}=\sum_{l_1=1}^{r_1}\sum_{l_2=1}^{r_2}  C_{l_1l_2}
[T_1]_{i_1l_1 }[T_2]_{i_2l_2 }, 
\end{equation}
We can obviously diagonalize the core tensor by using singular value 
decomposition.

\vs
\noindent
A drawback of the Tucker decomposition is the storage 
requirement of the core tensor $C_{l_1\cdots l_m}$, which is $O(r^m)$.  
Such problem can be 
mitigated by attempting to make zero as many entries of 
$C_{l_1\cdots l_m}$ as possible through suitable linear 
transformations. Another option is to introduce 
further separability properties of the core tensor. 
This yields a multitude of possible 
expansions, including hierarchical Tucker \cite{Grasedyck,Bachmayr} and  
Tucker tensor train \cite{Nouy2017,Rohwedder}. 
All these series expansions can be conveniently visualized 
by suitable graphs, and as such they fall within the setting 
of {\em tensor-networks}.

\subsubsection{Hierarchical Tucker Decomposition}
\label{sec:HT}
We have seen that the canonical tensor 
decomposition of a nonlinear functional 
has a relatively small number of degrees of freedom,
but its computation is not straightforward. 
The most reliable algorithms are based 
on regularized optimization and randomized methods 
(see \cite{Reynolds,Battaglino} for recent results), 
which often yield different series expansions for 
different initial conditions.
On the other hand, we have seen the that high-order 
Schmidt decomposition  suffers from the curse of dimensionality 
(dimension of the core tensor), but it is straightforward 
to compute by {\em nearly optimal} and {\em robust} 
techniques such as high-order singular value 
decomposition \cite{Lathauwer}. 
The hierarchical Tucker format was introduced 
by Hackbush and K\"{u}hn in \cite{Hackbusch2009} (see also \cite{Hackbusch_book,Grasedyck}) to mitigate 
the dimensionality problem in the core-tensor of 
the classical high-order Schmidt decomposition. 
The key idea is to perform a sequence 
of classical bi-variate Schmidt decompositions until the 
approximation problem is reduced to a product 
of one-dimensional functions. To illustrate the method 
is a simple way, consider the following cylindrical functional  
\begin{equation}
f(a_1,a_2,a_3)=\sum_{i,j=1}^r A^{\{1\}}_{ij}T^{\{1\}}_i(a_1) T^{\{2,3\}}_j(a_2,a_3), \qquad a_j=(\theta,\varphi_j).
\end{equation}
The matrix representation of $f(a_1,a_2,a_3)$ (either Galerkin or collocation form) relative to the given 
bases $T^{\{1\}}_i$ and $T^{\{2,3\}}_j$ 
is $A_{ij}^{\{1\}}$, i.e.,  it has two indices\footnote{The 
second index labels points on the plane $(a_2,a_3)$ 
(collocation setting),  or the projection of $f$ onto the two-dimensional basis function $T_j^{2,3}(a_2,a_3)$ (Galerkin setting).}. We decompose $T^{\{2,3\}}_j$ further by another Schmidt expansion
\begin{equation}
T^{\{2,3\}}_j(a_2,a_3)=\sum_{l,n=1}^r A^{\{2\}}_{jln} 
T^{\{2\}}_l(a_2)T^{\{3\}}_n(a_3),
\end{equation}
 to obtain
 \begin{equation}
f(a_1,a_2,a_3)=\sum_{i,j,l,n=1}^r A^{\{2\}}_{jln} 
A^{\{1\}}_{ij}T^{\{1\}}_i(x_1)T^{\{2\}}_l(a_2)T^{\{3\}}_n(a_3).
\end{equation}
The procedure just described is at the basis of the hierarchical 
Tucker decomposition and it yields an expansion in the 
form \eqref{Tucker}, where the core tensor is factored 
as a {\em product of at most three-dimensional tensors}. 
To show this in a higher-dimensional case, consider the 
six-dimensional cylindrical functional $f(a_1,...,a_6)$. 
By performing a sequence of Hilbert-Schmidt 
factorizations, while keeping the separation 
rank $r$ constant in each factorization, we obtain 
\begin{align}
f(a_1,...,a_6)= &\sum_{i_7,i_8=1}^r A^{\{1\}}_{i_7i_8} T^{\{1,2,3\}}_{i_7}
(a_1,a_2,a_3)T^{\{4,5,6\}}_{i_8}(a_4,a_5,a_6),\nonumber\\
=&\sum_{i_7,i_8=1}^r A^{\{1\}}_{i_7i_8} 
\sum_{i_1,i_9=1}^r  A^{\{2\}}_{i_7i_1i_9}
 T^{\{1\}}_{i_1}(a_1)
 T^{\{2,3\}}_{i_9}(a_2,a_3)
\sum_{i_4,i_{10}=1}^r A^{\{3\}}_{i_8 i_4i_{10}}
T^{\{4\}}_{i_4}(a_4)T^{\{5,6\}}_{i_{10}}(a_5,a_6),\nonumber\\
&=\sum_{i_1,\cdots, i_{6}=1}^r
C_{i_1\cdots i_6}
T^{\{1\}}_{i_1}(a_1)T^{\{2\}}_{i_2}(a_2)T^{\{3\}}_{i_3}(a_3)
T^{\{4\}}_{i_4}(a_4)T^{\{5\}}_{i_5}(a_5)
T^{\{6\}}_{i_{6}}(a_6),
\label{HTD}
\end{align}
where the $6$-dimensional core tensor is explicitly given 
\begin{equation}
C_{i_1\cdots i_{6}} =
\sum_{i_7,i_8=1}^r A^{\{1\}}_{i_7i_8} 
\sum_{i_9=1}^r  A^{\{2\}}_{i_7i_1i_9} A^{\{4\}}_{i_9i_2i_3}
\sum_{i_{10}=1}^r A^{\{3\}}_{i_8 i_4i_{10}}A^{\{5\}}_{i_{10}i_5i_6}.
\end{equation}
Note that such kernel has a recursive structure and it is factored 
as a product of at most three-dimensional matrices. 
Parallel algorithms to perform basic operations on  
hierarchical Tucker expansions such as addition, multiplication, 
and rank reduction  were recently proposed in 
\cite{Grasedyck2017,Etter}. Also, an optimization 
framework that leverages on the on recursive subspace 
factorizations of Hierarchical Tucker expansions was 
developed in \cite{DaSilva2015}.

\subsubsection{Tensor Train Decomposition} 
If we single-out one variable at a time and perform a 
sequential Schmidt decomposition of the remaining 
variables we obtain the so-called {\em tensor-train} 
decomposition \cite{OseledetsTT,Rohwedder}. 
Tensor train decomposition is a subcase case of hierarchical 
Tucker decomposition. For example, the tensor train 
decomposition of the four-dimensional 
cylindrical functional $f(a_1,...,a_4)$ reads 
\begin{align}
f(a_1,...,a_4)= &\sum_{i1,i_2=1}^r A^{\{1\}}_{i_1i_2} 
T^{\{1\}}_{i_1}(a_1)T^{\{2,3,4\}}_{i_2}(a_2,a_3,a_4),\nonumber\\
=&\sum_{i_1,i_2=1}^r A^{\{1\}}_{i_1i_2}
T^{\{1\}}_{i_1}(a_1) \sum_{i_3,i_4=1}^r  A^{\{2\}}_{i_2i_3i_4}
 T^{\{2\}}_{i_3}(a_2) T^{\{3,4\}}_{i_4}(a_3,a_4),\nonumber \\
=&\sum_{i_1,...,i_4=1}^r A^{\{1\}}_{i_1i_2}A^{\{2\}}_{i_2i_3i_4}
T^{\{1\}}_{i_1}(a_1)    T^{\{2\}}_{i_3}(a_2) 
\sum_{i_5,i_6}^rA^{\{3\}}_{i_4i_5i_6} T^{\{3\}}_{i_5}(a_3) 
T^{\{4\}}_{i_6}(a_4),\nonumber \\
=&\sum_{i_1,...,i_6=1}^r A^{\{1\}}_{i_1i_2}A^{\{2\}}_{i_2i_3i_4}
A^{\{3\}}_{i_4i_5i_6}
T^{\{1\}}_{i_1}(a_1)    T^{\{2\}}_{i_3}(a_2)  T^{\{3\}}_{i_5}(a_3) T^{\{4\}}_{i_6}(a_4).
\label{TTD}
\end{align}
Tensor train and hierarchical Tucker expansions
can be conveniently visualized by {\em graphs}.

\subsubsection{Tensor Networks}
During the last decades, the field of Tensor Networks has 
undergone an explosion of results in several different directions, 
e.g., in the study of quantum many-body systems, and more generally in multivariate functional approximation. Tensor networks
can be conveniently represented by {\em undirected graphs}. 
To this end, we adopt the following standard convention: 
\begin{itemize}

\item a node in a graph represents a tensor 
in as many variables as the number of the 
edges connected to it; 

\item connecting two tensors by an edge represents 
a tensor contraction over the index associated to certain variable.
\end{itemize}

\noindent
In this way, a three-dimensional Tucker format is represented by one 
node with three edges emanating from it. If we connect two of such 
tensors with one edge, then we obtain a tensor in four variables. In particular, the connection operation here results in a 
Tucker format with the core tensor 
represented by a product of two three-dimensional 
tensors in which we contract one 
index (see Figure \ref{Fig:TN}).
\begin{figure}[t]
\centerline{\hspace{-1.2cm}Tucker Format\hspace{6.cm} Tensor Network}
\centerline{\hspace{-1cm}
\includegraphics[width=3.5cm]{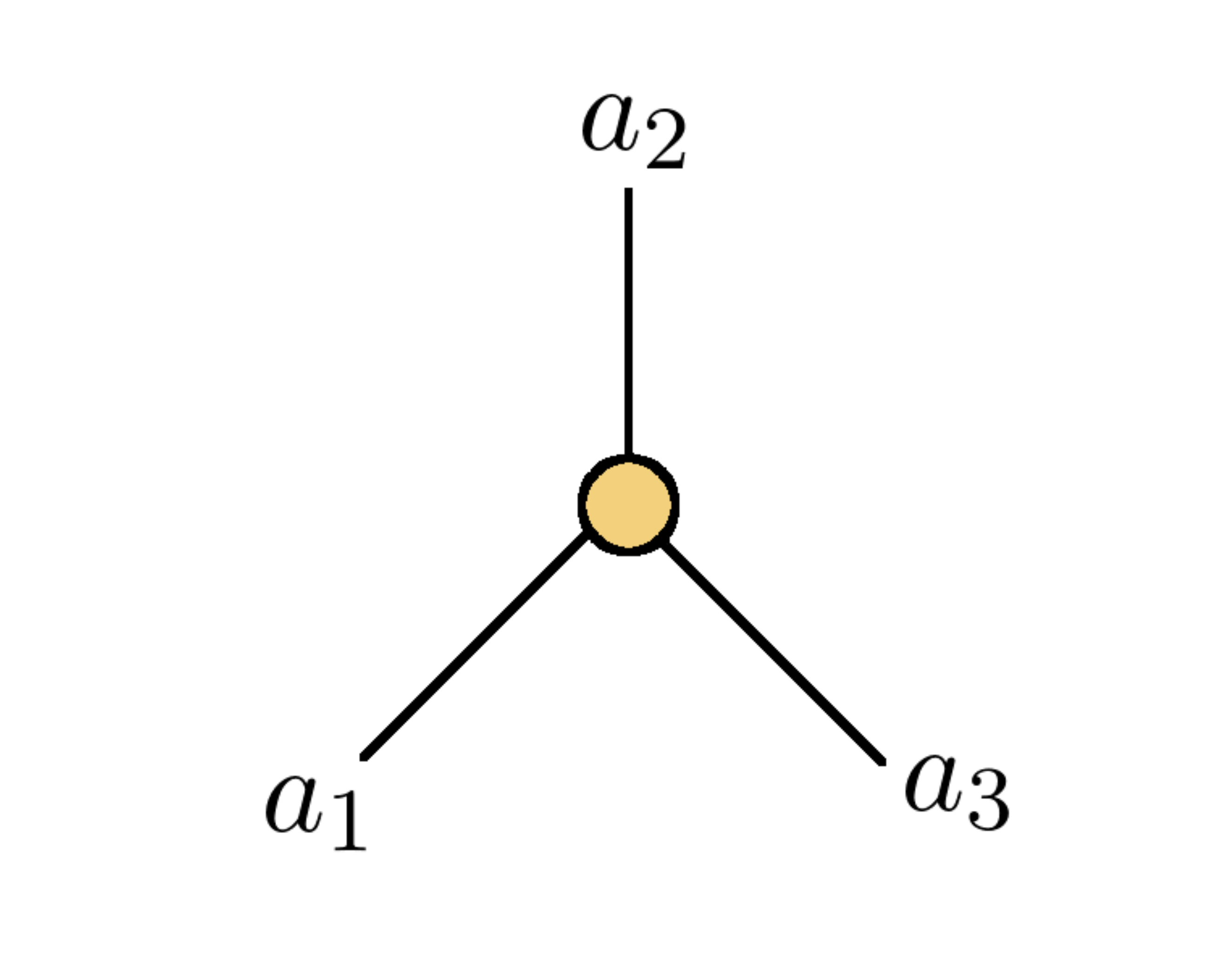}
\hspace{4.2cm}
\includegraphics[width=4.5cm]{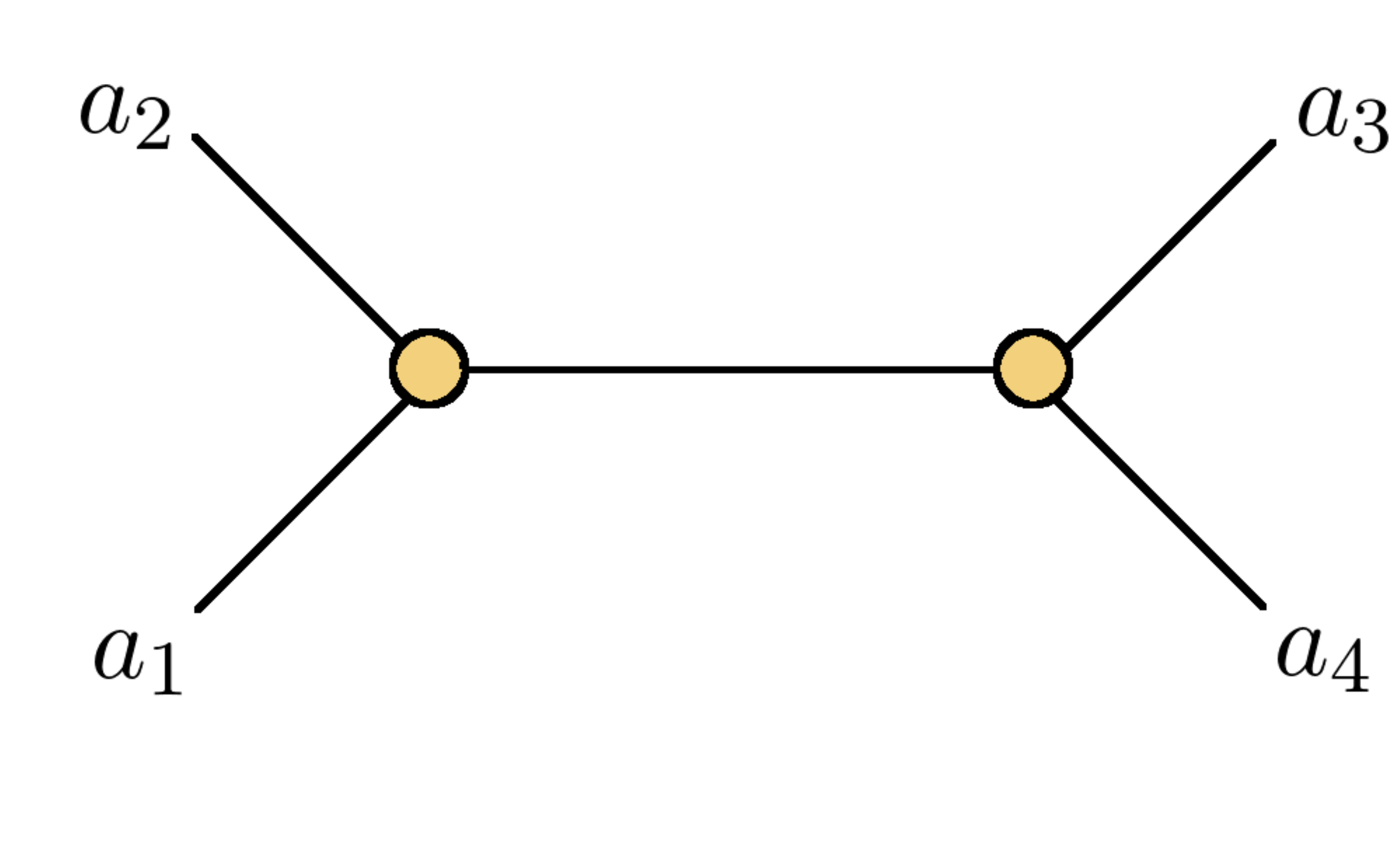}
}
{\footnotesize $\displaystyle \sum_{i_1,i_2,i_3=1}^r
A_{i_1i_2i_3}T^{\{1\}}_{i_1}(a_1)    T^{\{2\}}_{i_2}(a_2)  T^{\{3\}}_{i_3}(a_3)$}\hspace{1.5cm}
{\footnotesize $\displaystyle \sum_{i_1, ... ,i_4=1}^r
\left(\sum_{j=1}^rA^{\{1\}}_{i_1i_2j}A^{\{2\}}_{ji_3i_4}\right)T^{\{1\}}_{i_1}(a_1)    T^{\{2\}}_{i_2}(a_2)  T^{\{3\}}_{i_3}(a_3)T^{\{4\}}_{i_4}(a_4)$
}
\caption{\color{r}Examples of tensor networks representing a 
cylindrical functional in three and four variables. }
\label{Fig:TN}
\end{figure}
The graph representation of the hierarchical Tucker 
and tensor train decomposition is shown in 
Figure \ref{fig:HT_TT}.
\begin{figure}[ht]
\centerline{\hspace{-1.2cm}Hierarchical Tucker\hspace{4.2cm} Tensor Train}
\centerline{
\includegraphics[width=6.5cm]{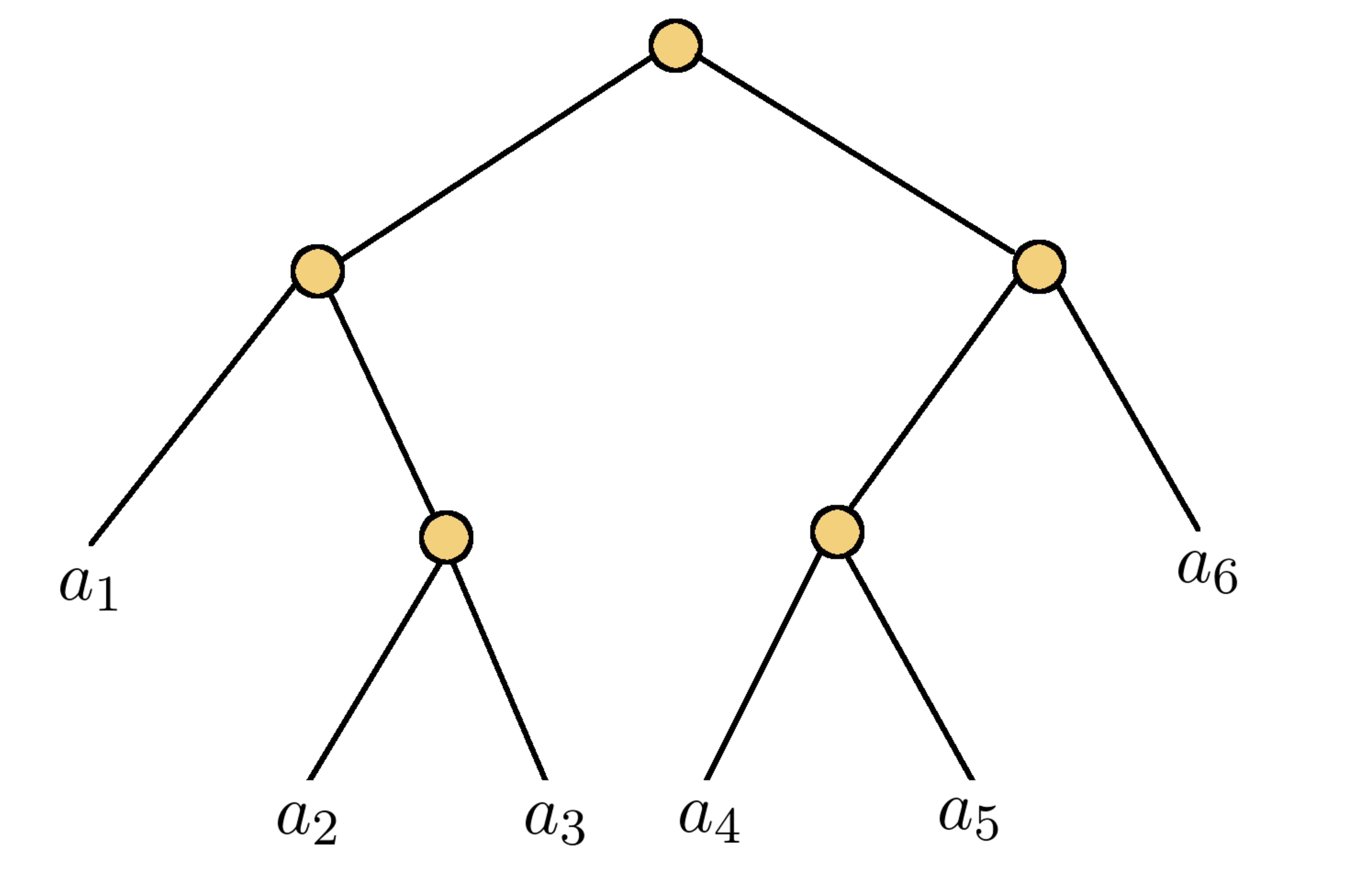}
\hspace{1cm}
\includegraphics[width=6cm]{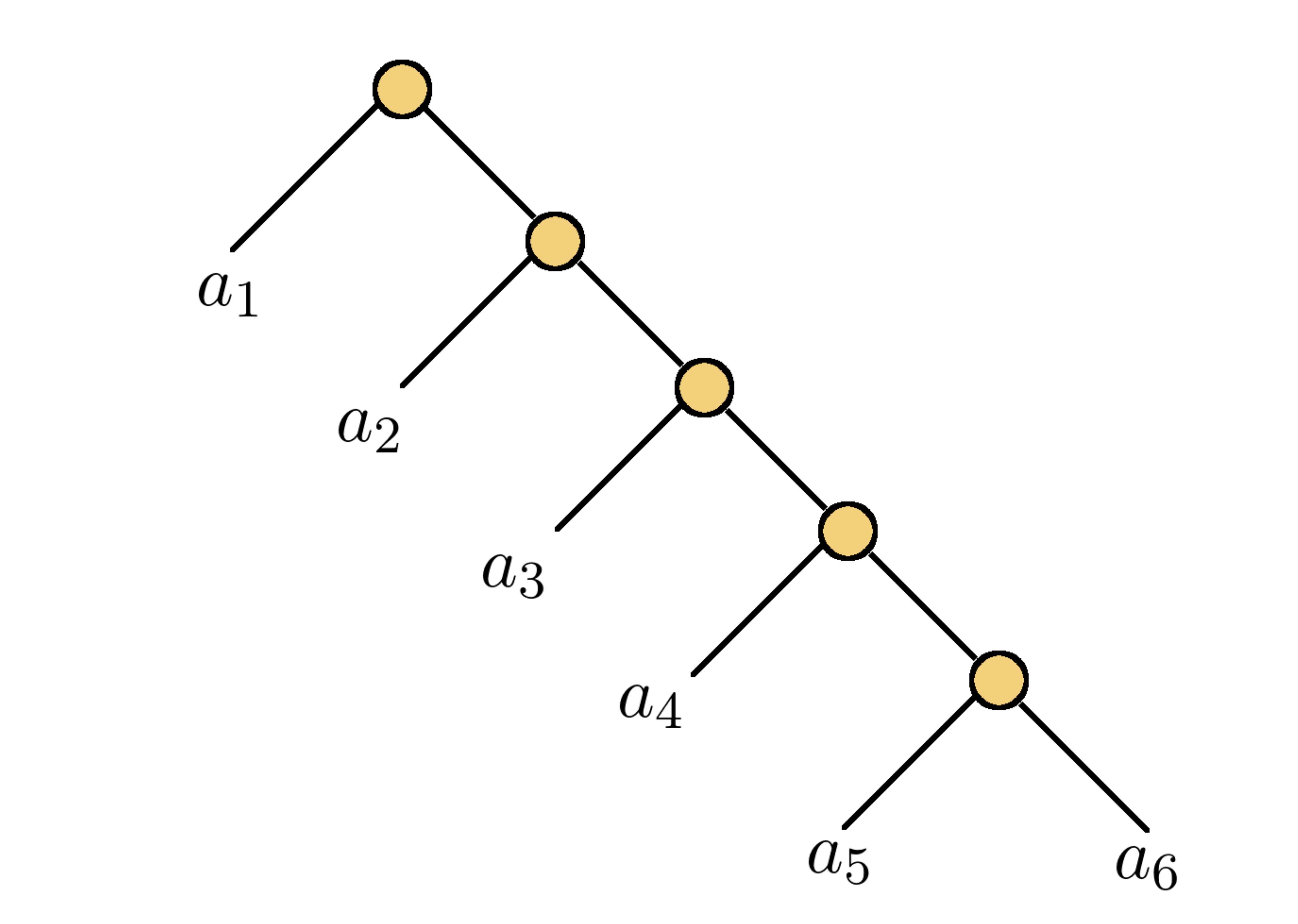}
}
\caption{\color{r}Graph representation of the Hierarchical Tucker (HT) and tensor train (TT) decomposition of a six-dimensional cylindrical functional.}
\label{fig:HT_TT}
\end{figure}
A more general tensor network representation of  
a cylindrical functional in five variables is shown in Figure 
\ref{Fig:TN1}.
\begin{figure}[t]
\centerline{\hspace{0.5cm}Tensor Network\hspace{3.cm} 
Equivalent Tucker Format}
\centerline{\hspace{-2cm}
\includegraphics[width=6.0cm]{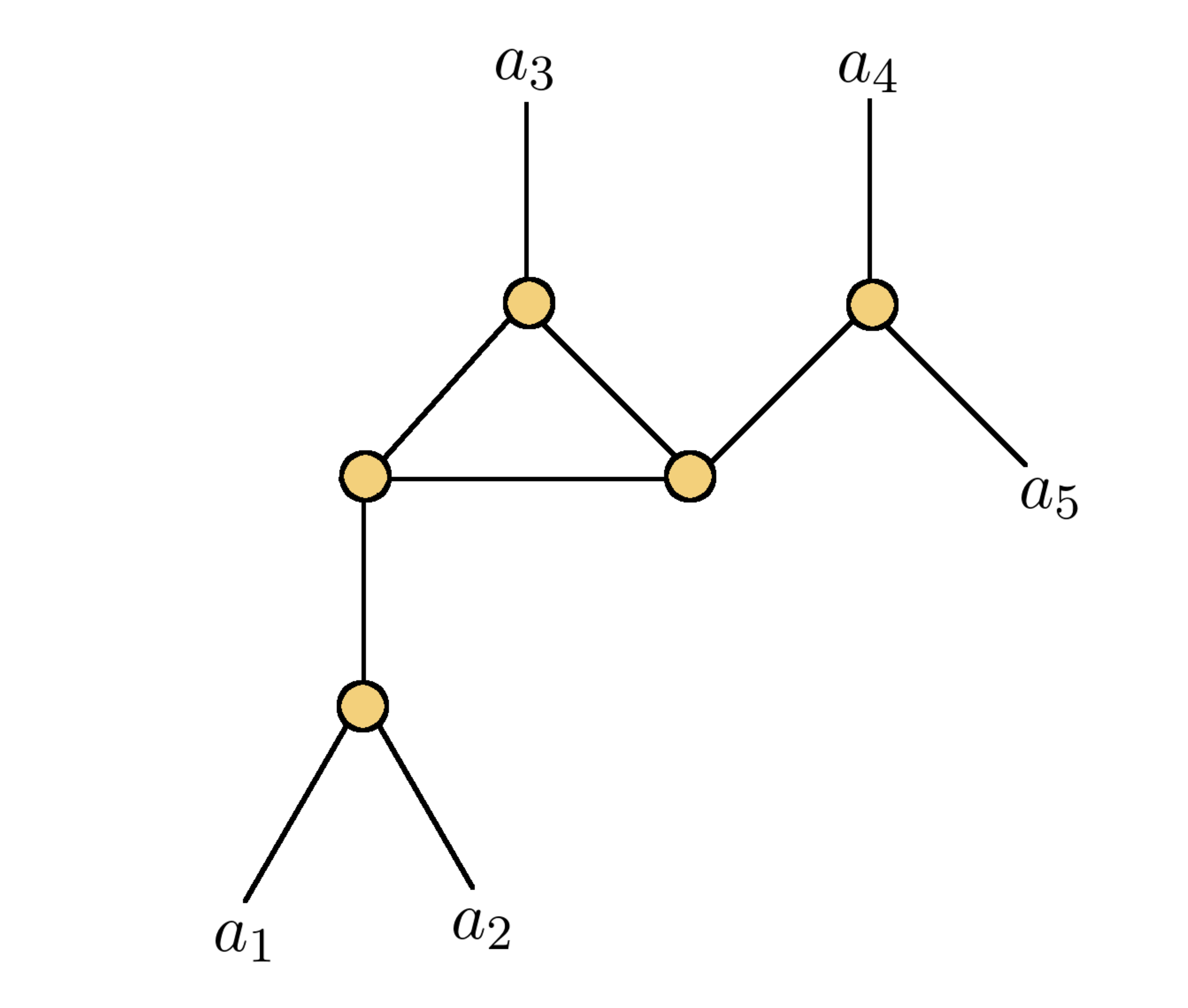}
\hspace{2cm}
\includegraphics[width=3.0cm]{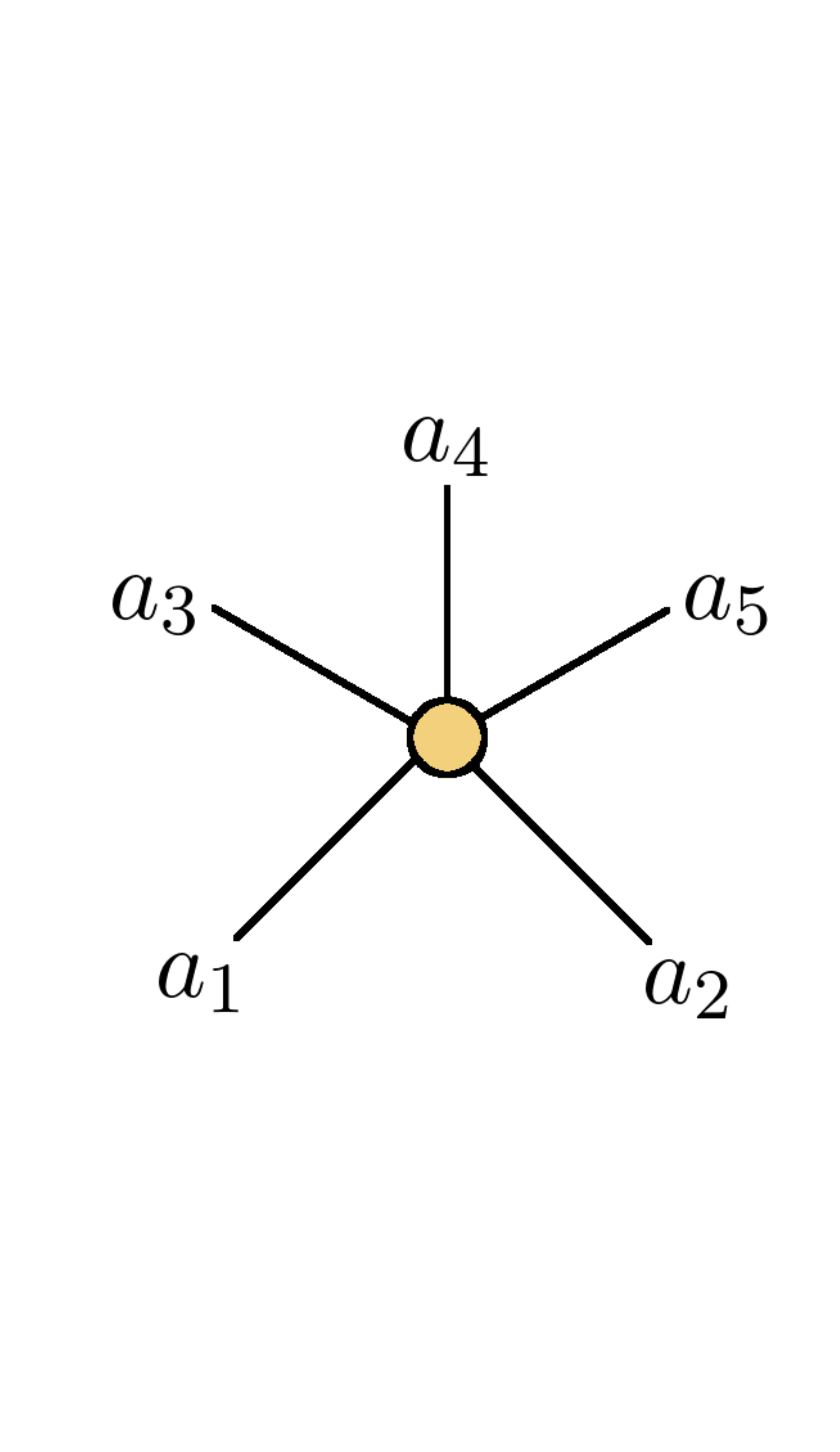}}
\caption{\color{r}Tensor network representation of  a five-dimensional cylindrical functional. After reduction (elimination of the cycle and connection of nodes), it is easy to see the graph is  equivalent to 
the graph corresponding to a Tucker format of dimension five (right).
The core tensor of such Tucker format is expressed as products and contractions of at most three-dimensional tensors.
}
\label{Fig:TN1}
\end{figure}
Note that each node is connected with at most 
three edges (coordination number equal to three) 
as we employed a Schmidt decomposition to represent
each function with more than one variable. Clearly, 
the graph can be reduced to one 
node  with five leaves (i.e. a Tucker series) 
by eliminating the cycle and clustering all nodes that are connected 
by edges \cite{Kressner2014} . 
The graph can be therefore reduced to a node 
with five leaves -- i.e. a Tucker series, with a particular 
structure of the core tensor. At this point the 
graph is irreducible (one node with five leaves).
Efficient algorithms that implement basic operations 
between tensors, such as addition, orthogonalization, rank 
reduction, scalar products, multiplication, and linear 
transformations are discussed in 
\cite{Kressner2014,Etter,Grasedyck2015,Nouy2017,Kolda}.

}

\subsection{Generalized Lagrangian Interpolation}
\label{sec:gen_Lagr}
We have seen that the classical Lagrange interpolation 
problem for multivariate functions can be generalized 
to functionals in Hilbert or Banach spaces, i.e., 
objects depending on an infinite number of variables. 
Given a real-valued continuous functional 
\begin{equation}
 F:D(F)\rightarrow \mathbb{R},
\end{equation}
and a set of elements $\{\theta_1(x),...,\theta_n(x)\}$ in $D(F)$, 
the Lagrangian interpolant of $F$ is a continuous functional  
\begin{equation}
\Pi: D(F) \rightarrow \mathbb{R},
\end{equation}
such that
\begin{equation}
\Pi ([\theta_i])=F([\theta_i]). 
\label{interpolation problem}
\end{equation}
A simple representation of $\Pi([\theta])$ 
can be given in terms of {\em cardinal basis functionals} 
$g_i:D(F)\rightarrow \mathbb{R}$ as
\begin{equation}
\Pi ([\theta])=\sum_{i=1}^n F([\theta_i])g_i([\theta]), \quad
\textrm{where}\quad g_i([\theta_j])=\delta_{ij}.
\label{interpolant}
\end{equation} 
Similarly to classical Lagrange interpolation for 
functions from $\mathbb{R}^n$ to $\mathbb{R}^m$, we 
have a great freedom in selecting $g_i([\theta])$. We have 
already seen, for example, Porter's and Prenter's cardinal 
bases defined in equation \eqref{gi_0} and equation \eqref{gi_prenterer}, respectively.
Now, let 
\begin{equation}
\kappa:D(F)\times D(F)\rightarrow\mathbb{R}
\end{equation}
be any continuous (not necessarily symmetric) functional subject to the sole 
requirement that 
\begin{equation}
\kappa([\theta],[\eta])=0\quad \Leftrightarrow \quad \theta(x)=\eta(x).
\label{cond1}
\end{equation}
Examples of such functional are 
\begin{align}
\kappa_1([\theta],[\eta])&=\left\|\theta(x)-\eta(x)\right\|^{p},
\qquad  p>0\label{positiveker}\\
\kappa_2([\theta],[\eta])&= \left\|\theta(x)-\eta(x)\right\|+
\left\|\theta(x)\right\|-\left\|\eta(x)\right\|.\\
\kappa_3([\theta],[\eta])&= 1- \exp\left(-\left\|\theta(x)-\eta(x)\right\|^2\right).
\end{align}
By using $\kappa([\theta],[\eta])$ we  define the following 
cardinal basis associated with the interpolation nodes
$\{\theta_1(x),...,\theta_n(x)\}$ in $D(F)$ 
(see \cite{Cheney}, Ch. 10)
\begin{equation}
g_i([\theta]) = \prod_{\substack{j=1\\j\neq i}}^n
\frac{\kappa([\theta],[\theta_j])}{\kappa([\theta_i],[\theta_j])}
\qquad i=1,...,n.\label{cardinal_cheeney}
\end{equation}
A substitution of \eqref{cardinal_cheeney} into 
\eqref{interpolant} yields a {\em Lagrangian functional interpolant}  that depends on the choice of $\kappa([\theta],[\eta])$.
An advantage of such interpolant over Porter's interpolant 
is that the basis functionals \eqref{cardinal_cheeney} 
are given analytically and do not require 
any computation such as the inversion (and storage) of 
the $H$-matrix \eqref{matH}. On the other hand, we have no 
guarantee that \eqref{cardinal_cheeney} have good 
approximation properties away from the interpolation nodes. 
The Lagrange interpolation formula \eqref{interpolant}-\eqref{cardinal_cheeney}
can be written in a {\em barycentric form}. 
To this end, we first rewrite the numerator 
of \eqref{cardinal_cheeney} as 
$\Omega([\theta])/\kappa([\theta],[\theta_i])$, where 
\begin{equation}
 \Omega([\theta])=\prod_{i=k}^n\kappa([\theta],[\theta_k]).
\end{equation}
Then define 
the {\em barycentric weights} as
\begin{equation}
 w_i = \prod_{\substack{j=1\\j\neq i}}^n
\frac{1}{\kappa([\theta_i],[\theta_j])}.\label{weights}
\end{equation}
Clearly, 
\begin{equation}
 g_i([\theta])=\Omega([\theta])\frac{w_i}{\kappa([\theta],[\theta_i])}. 
\end{equation}
Substituting this into \eqref{interpolant} and applying it 
to the constant functional $F([\theta])=1$ yields 
\begin{equation}
\Omega([\theta])=\frac{1}{\displaystyle\sum_{i=1}^n\frac{w_i}{\kappa([\theta],[\theta_i])}}.
\end{equation}
Therefore, we obtain the interpolant
\begin{equation}
\Pi([\theta])=
\left(\sum_{i=1}^n\frac{w_i}{\kappa([\theta],[\theta_i])}\right)^{-1}
\sum_{i=1}^n \frac{w_iF([\theta_i])}{\kappa([\theta],[\theta_i]).}
\end{equation}
This expression can be evaluated in $\mathcal{O}(n)$ operations, provided we have available 
$w_i$, which needs $\mathcal{O}(n^2)$ operations.
What about updating? Adding a new interpolation node $\theta_{n+1}$ entails two 
calculations: 1) divide each $w_j$ by $\kappa([\theta_j],[\theta_{n+1}])$ and 
2) compute $w_{n+1}$ using the formula \eqref{weights}.

\vs
\noindent 
{\em Remark:} The interpolation processes defined by \eqref{interpolant} 
and \eqref{cardinal_cheeney} has a variation known as
Shepard's method \cite{Farwig,Lazzaro}. Such method was studied by 
Allasia and Bracco \cite{Allasia} for  functional interpolation 
in Banach spaces and it relies on the cardinal functionals
\begin{equation}
g_i([\theta]) = \frac{1}{\displaystyle1+\sum_{\substack{ k=1\\k\neq i}}^n 
\frac{{\kappa}([\theta],[\theta_i])}{{\kappa}
([\theta],[\theta_k])}}, 
\qquad i=1,...,n.\label{Shepardb}
\end{equation}
The kernel ${\kappa}([\theta],[\eta])$ here satisfies \eqref{cond1} and 
the additional positivity requirement  
\begin{equation}
{\kappa}([\theta],[\eta])>0\qquad \textrm{if} \qquad 
\theta(x)\neq\eta(x).
\label{cond2}
\end{equation} 
This implies that the basis functionals \eqref{Shepardb} 
have the following properties 
\begin{equation}
 0\leq g_i([\theta])\leq 1\quad \forall \theta\qquad 
 g_i([\theta_j])=\delta_{ij}.
\end{equation}
One unfortunate consequence of these conditions is that 
the functional derivative of $g_i([\theta])$ (if it exists) 
is identically zero at all nodes $\theta_j$. 
In fact, $g_i([\theta])$ is always positive, it has 
a maximum (equal to 1) at $\theta_i$ and 
many minima (equal to 0) at all other nodes. 
This yields functional differentiation matrices \eqref{Hint} 
that are identically zero.
Functional Lagrangian interpolation can also be  defined in terms of 
{\em radial basis functionals} 
$\kappa(\left\|\theta-\theta_j\right\|)$, 
for suitable choices of $\kappa$.
In this cases, the cardinal basis $g_i([\theta])$ can be obtained 
by standard techniques such as the ratio of two Vandermonde 
determinants. Interpolation using radial 
basis functionals can converge towards polynomial functional 
interpolants if we take increasingly flat radial bases 
(see, e.g., \cite{Schaback,Driscoll,Song}). 

\vs
\noindent
The Lagrangian interpolation process defined 
by \eqref{interpolant}, \eqref{cond1} and
\eqref{cardinal_cheeney} can be generalized further 
by considering a family of continuous functionals
\begin{equation}
\kappa_i:D(F)\times D(F)\rightarrow\mathbb{R},\qquad i=1,...,n
\end{equation}
as many as the number of interpolation nodes, satisfying 
$\kappa_i([\theta],[\theta_j])\neq 0$ for $i\neq j$, 
$\kappa_i([\theta],[\theta_i])=0$ and 
$\kappa_i([\theta],[\theta])=0$.
The cardinal basis corresponding 
to $\{\kappa_1,...,\kappa_n\}$ is 
\begin{equation}
g_i([\theta])=\prod_{\substack{j=1\\j\neq i}}^n \frac{\kappa_i([\theta],[\theta_j])}
{\kappa_i([\theta_i],[\theta_j])}.
\end{equation}
In particular, the choice
\begin{equation}
 \kappa_i([\theta],[\eta])=(\theta-\eta,\theta_i-\eta)
 \label{kprenter}
\end{equation}
yields Prenter's cardinal basis \eqref{gi_prenterer}. It fact, 
\begin{align}
\kappa_i([\theta],[\theta_j])  &=(\theta-\theta_j,\theta_i-\theta_j),\label{kk1}\\
\kappa_i([\theta_i],[\theta_j])&=\left\|\theta_i-\theta_j\right\|^2,\label{kk2}\\
\kappa_i([\theta],[\theta_i])  &=0.\label{kk3}
\end{align}
\noindent
Another generalization of the Lagrange interpolation process 
can be obtained along the lines of Porter's functional interpolants. To this end, 
let $\kappa([\theta],[\eta])$ be any symmetric functional 
such that the matrix 
\begin{equation}
H_{ij}=\kappa([\theta_i],[\theta_j])
\label{kH}
\end{equation}
is invertible for every set of distinct nodes $\{\theta_j(x)\}$. 
Then the functional 
\begin{equation}
\Phi_n([\theta])=\sum_{k=1}^m F(\theta_{k})\underbrace{\sum_{j=1}^m
\kappa([\theta_k],[\theta_j])^{-1}\kappa([\theta],[\theta_j])}_{g_k([\theta])},
\end{equation}
where $\kappa([\theta_k],[\theta_j])^{-1}$ denotes the inverse matrix of 
\eqref{kH}, interpolates $F([\theta])$ at $\theta_k(x)$ ($k=1,...,m$). 
The question of how to select the symmetric functional $\kappa$
depends on the requirements we impose on the behavior of the functional interpolant away from the interpolation nodes. 
In particular, 
the choice 
\begin{equation}
\kappa([\theta],[\eta])=\sum_{p=0}^n (\theta,\eta)^p
\end{equation}
yields Porter's interpolants, i.e., polynomial functionals of 
order $n$ with minimal norm.

\subsubsection{Optimal Interpolation Nodes}
\label{sec:optimal interpolation nodes}
In this Section we briefly address the question 
of how to select the interpolation nodes in the function space
$D(F)$ in such a way that the corresponding 
functional interpolant, e.g., the Porter's one 
\eqref{Porter_interpolant}-\eqref{gi_0}, 
exhibits good approximation properties.

\paragraph{Adaptive Leja Sequences} 
An effective approach to select sub-optimal interpolation nodes 
relies on adaptive Leja sequences \cite{Akil,Bos}.
In this setting, given the 
set of nodes $\{\theta_1(x),...,\theta_n(x)\}$ in $D(F)$, 
the problem is to find a new node $\theta^*(x)$
satisfying the maximization principle  
\begin{equation}
\argmax_{\theta^*(x)\in  D(F)}\left|\det ( V)\right|
\label{determinant_VDM}
\end{equation}
where $V$ is a Vandermonde-like operator 
(infinite-dimensional matrix) (see \cite{Bos,Akil,Akil1,Marchi} 
for further details).
This results in  the so-called {\em Leja construction} of the 
test function $\theta^*(x)$, a greedy version of the well-known 
Fekete construction involving nonlinear optimization. 
There are also non-determinant versions 
of multivariate Leja interpolation nodes, which are related 
to {\em potential theory} \cite{Garcia,Ishizaka}. 
In this setting, given a symmetric functional 
$\kappa([\theta(x)],[\eta(x)])$ one can construct a 
{\em greedy $\kappa$-energy sequence} in which 
$\theta^*(x)$ is chosen to satisfy
\begin{equation}
\sum_{i=1}^{n}\kappa([\theta^*(x)],[\theta_i(x)])=
\argmin_{\theta^*\in D(F)}\sum_{i=1}^{n} \kappa([\theta^*(x)],[\theta_i(x)]), 
\end{equation}
assuming $\{\theta_1(x),...,\theta_n(x)\}$ are available.
The element $\theta^*(x)$ obtained in this way is not unique.
Also, depending on the choice of the kernel 
function $\kappa([\theta(x)],[\eta(x)])$, we can 
have different sets of nodes.  A relevant kernel for 
the $\kappa$-energy sequence is the so-called {\em Riesz $p$-kernel} 
\begin{equation}
 \kappa_p([\theta(x)],[\eta(x)])=
 \begin{cases}
\displaystyle \frac{1}{\left\|\theta(x)-\eta(x)\right\|^p} & \textrm{if $p>0$}\\
\displaystyle -\ln(\left\|\theta(x)-\eta(x)\right\|) & \textrm{if $p=0$}\\
\end{cases}
\label{Rieszp}
\end{equation}
Such kernel avoids placing points too close 
to each other (the potential $\kappa_p$ is rather large in such regions), 
while decreasing monotonically at larger distances.
The asymptotic behavior of greedy $\kappa$-energy sequences has 
been studied theoretically by L\'opez-Garc\'ia and Saff \cite{Garcia}, 
in a finite-dimensional setting.

\paragraph{Minimization of Lebesgue Functionals} Another approach 
to determine sub-optimal nodes in the function space $D(F)$ 
relies on greedy algorithms minimizing {\em Lebesgue-like functionals}. 
This problem has been recently addressed in a 
finite-dimensional setting by Narayan and Xiu \cite{Akil1,Akil2}, 
Maday {\em et al.} \cite{Maday,Maday1} and 
Van Barel {\em et. al} \cite{Barel}. 
The basic idea is simple and can be generalized to 
the infinite-dimensional case. Suppose we have available 
$\{\theta_1(x)$, ..., $\theta_{n}(x)\}$. We look for a new 
node $\theta^{*}(x)$ maximizing a suitable objective 
functional, e.g., related to the well-known Lebesgue function 
in finite-dimensional interpolation problems. 
In particular, following the weighed approach of Narayan 
and Xiu \cite{Akil1}, we can look for a new node $\theta^{*}(x)$
satisfying the following (greedy) optimization problem 
\begin{equation}
\argmax_{\theta^{*}\in D(F)} \chi([\theta^*]),\qquad 
\chi([\theta])=W([\theta])^2\sum_{k=1}^n 
\frac{g_k([\theta])^2}{W([\theta_k])^2},
\label{lebesgue function}
\end{equation}
where $W([\theta])$ is a positive 
functional and $g_i([\theta])$ is a cardinal basis. 
Other choices of $\chi([\theta])$ yield different 
sets of nodes. For example,  
\begin{equation}
\chi([\theta])=W([\theta])\sum_{k=1}^n 
\frac{g_i([\theta])^2}{W([\theta_k])} 
\label{fejerOPT}
\end{equation}
yields an infinite-dimensional extension of the Fej\'er points.
The sequence of nodes $\theta_i(x)$ we obtain in this 
way {\em strongly depends on the initial set of nodes}, 
on the weight $W([\theta])$ and 
on the cardinal basis $\{g_i(\left[\theta\right])\}$. 
Finding the maximum of \eqref{lebesgue function} involves 
computing the solution to an optimization 
problem in infinite dimensions \cite{Bensoussan_DaPrato}. 
A necessary condition for a stationary point 
of \eqref{lebesgue function} is 
\begin{equation}
\frac{\delta \chi([\theta]) }{\delta \theta(x)}=0.
\label{Euler-Lagrange}
\end{equation}

\vs
\noindent
{\em Example 1:} The Euler-Lagrange equation \eqref{Euler-Lagrange} 
can be written down explicitly, e.g., for Porter's 
interpolants. To this end, simply substitute \eqref{gi_0}  in \eqref{fejerOPT} and set $W=1$ to obtain
\begin{equation}
\chi([\theta]) = \sum_{k=1}^n
\sum_{j,z=1}^m H_{jk}^{-1}H_{zk}^{-1} \sum_{p,q=0}^h\left(\theta_j,\theta\right)^p
\left(\theta_z,\theta\right)^q.
\end{equation}
This allows us to write explicitly the conditions identifying 
stationary points of $\chi([\theta])$, e.g., in the 
function space 
\begin{equation}
D_N=\left\{\theta(x)\in D(F)\,\left|  
\, \theta(x)=\sum_{k=1}^{N}a_k\varphi_k(x)\right.\right\}. 
 \label{Gqsig}
\end{equation}

\subsection{High-Dimensional Model Representation} 
\label{sec:HDMR}
The high-dimensional model representation (HDMR) \cite{Rabitz,Li,Li1,CaoCG09} of the cylindrical functional \eqref{functional-approx} is a series 
expansion in the form
\begin{equation}
f(a_1, ..., a_m)=f_0+
\sum_{k=1}^m f_k(a_k) +
\sum_{\substack{k,j=1\\k< j}}^m 
f_{kj}(a_k,a_j)+\cdots.
\label{functional-HDMR}
\end{equation}
The functional $f_{0}$ is simply a constant. 
The functionals $f_i(a_i)$, which we shall call 
first-order interaction terms, give us the 
overall effects of the variables $a_i=(\theta,\varphi_i)$ in $f$ 
as if they were acting independently of the other variables. 
The functions $f_{ij}(a_i, a_j)$ represent the interaction effects 
of $a_i$ and $a_j$, and therefore they are usually 
called second-order interactions. 
Similarly, higher-order terms reflect the cooperative 
effects of an increasing number of variables.
The interaction terms in the HDMR expansion can  be 
easily computed if we assume that the domain of 
$f$ is compactly supported and included in the hypercube 
$[-b,b]^m$ (see Appendix \ref{sec:functional integrals}). 
By introducing the finite-dimensional integration 
measure $w(\bm a)=w(a_1,...,a_m)$ and the vector 
notation $d\bm a_{-i}=da_{i-1}da_{i+1}\cdots da_m$,
we have 
\begin{align}
f_0 &= \int_{[-b,b]^m}  w(\bm a) f(\bm a) d\bm a,
\label{anova0}\\
f_i(a_i) &=  \int_{[-b,b]^{m-1}}  w(\bm a) f(\bm a) d\bm a_{-i} -f_0,
\label{anova1}\\
f_{ij}(a_i,a_j) &= \int_{[-b,b]^{m-2}}  w(\bm a) f(\bm a) d\bm a_{-ij}-f_i(a_i)-f_j(a_j) -f_0,
\label{anova2}\\
&\cdots\nonumber . 
 \end{align}
This procedure generates, by construction, 
terms that are orthogonal in the weighted $L_2$ 
sense. The HDMR series \eqref{functional-HDMR} with 
components \eqref{anova0}-\eqref{anova2}
is often called ANOVA-HDMR expansion.
If we consider a Dirac delta measure 
$w(\bm a)=\delta(a_1-c_1)\cdots  \delta(a_m-c_m)$
with ``anchor point'' $\bm c=(c_1,...,c_m)$, then the 
HDMR series is called CUT-HDMR \cite{Dick1,Handy_ANOVA_2010}. 
There is also a random-sampling version 
of HDMR -- namely the RS-HDMR \cite{Li1} expansion -- in which 
the high-dimensional integrals in \eqref{anova0}-\eqref{anova1} 
are computed by Monte-Carlo or quasi-Monte Carlo methods. 
The HDMR expansion \eqref{functional-HDMR} is usually 
truncated at some interaction order, and the interactions
$f_i$, $f_{ij}$ are expanded relative to a certain basis. 
For example, the first-order ANOVA-HDMR expansion of 
$f$ reads 
\begin{equation}
f(a_1,...,a_m)\simeq f_0+
\sum_{k=1}^m f_k(a_k),
\label{functional-HDMR-1}
\end{equation}
where  
\begin{equation}
 f_k(a_k)=\sum_{i=1}^Q 
 \alpha_{ik}\phi_i(a_k),
\end{equation}
and  $\varphi_i$ are suitable basis 
functions in the variable $a_k-(\theta,\varphi_k)$.

\subsection{Cluster Expansion}
\label{sec:cluster_expansion}

The cluster expansion is a functional approximation 
method for Hopf characteristic functionals 
that leverages on the structure of the functional 
itself to provide an approximation that preserves 
important properties such as normalization and 
marginalization rules. 
To illustrate the method, let us 
consider the restriction of the Hopf functional to 
the finite-dimensional function space $D_m\subseteq D(\Phi)$
(see Example 1 in Section \ref{sec:Finite_Dim_Approx}, and Figure \ref{fig:1}). 

\begin{equation}
\phi(a_1,...,a_m)=\left<e^{i(a_1U_1(\omega)+\cdots+ a_m U_m(\omega))}\right>,
\qquad U_k(\omega)=\int_a^b u(x;\omega)\varphi_k(x)dx.
\label{jointCHF}
\end{equation}
To represent the high-dimensional (complex-valued) 
function $\phi(a_1,...,a_m)$ one can use techniques such as 
HDMR (Section \ref{sec:HDMR}), functional tensor methods 
(Section \ref{sec:tensor}), or sparse collocation 
\cite{Bungartz,DoostanOwhadi_2011,NovakR_96}. 
One of the problems with such techniques is that they 
do not preserve important properties of the 
characteristic function, normalization $\phi(0,...,0)=1$, 
and marginalization rules. 
However, in the case of \eqref{jointCHF}, we do know 
the {\em structure} of cylindrical functional we are 
approximating. Therefore, we can leverage on such 
structure to build an expansion that preserves 
marginalization rules and other properties.
To this end, let us first define the following 
reduced-order characteristic functions 
\begin{align}
\phi_n(a_n)&=\left<e^{ia_nU_n(\omega)}\right>,\label{short1}\\
\phi_{nm}(a_n,a_m)&=\left<e^{ia_nU_n(\omega)+ia_mU_m(\omega)}\right>,\label{short2}\\ 
\phi_{nmk}(a_n,a_m,a_k)&=\left<e^{ia_nU_n(\omega)+ia_mU_m(\omega)+ia_kU_k(\omega)}\right>,\label{short3}\\
&\cdots\nonumber.
\end{align}
Clearly,  we have
\begin{equation}
 \phi_{nmk}(a_n,0,0)=\phi_{nm}(a_n,0)= \phi_{n}(a_n), \qquad\textrm{(marginalization rule).}
\end{equation}
By using well-known cumulant series representations \cite{Kubo}, we expand \eqref{jointCHF} as  
\begin{equation}
 \left<\exp\left[i\sum_{j=1}^N a_jU_j(\omega) \right]\right>=
 \exp\left[\sum_{\nu_1,...,\nu_N=0}^\infty\left<U_1^{\nu_1}\cdots U_N^{\nu_N}\right>_c
 \prod_{k=1}^N\frac{(ia_k)^{\nu_k}}{\nu_k!}\right],
 \label{Kubo_cumulant}
\end{equation}
where the last summation is over $\nu_1$, ..., $\nu_N$ 
excluding the case in which all indices are zero, i.e., 
excluding $\nu_1=\cdots=\nu_N=0$. The cumulant averages 
in \eqref{Kubo_cumulant} can be written in terms of 
moments of $U_j$. For example, we have 
\begin{align}
\left<U_j\right>_c=&\left<U_j\right>,\nonumber\\
\left<U_jU_k\right>_c=&\left<U_jU_k\right>-\left<U_j\right> \left<U_k\right>,\nonumber\\
\left<U_jU_kU_m\right>_c=&\left<U_jU_kU_m\right>-\left<U_j\right> \left<U_kU_m\right>-
\left<U_k\right> \left<U_jU_m\right>-\left<U_m\right> \left<U_kU_j\right>+
2\left<U_j\right>\left<U_k\right>\left<U_m\right>.\nonumber
\end{align}
At this point, it is useful to write the \eqref{Kubo_cumulant} explicitly 
for simple cases, e.g., for \eqref{short2} and \eqref{short3}. We have,
\begin{equation}
\phi_{nm}(a_n,a_m)= \phi_n(a_n)\phi_m(a_m)
\exp\left[\sum_{j,k=1}^\infty\left<U_m^{k}U_n^{j}\right>_c
\frac{(ia_m)^{k}(ia_n)^{j}}{k!j!}\right]\quad \textrm{and}
\end{equation}
\begin{align}
\phi_{nmk}(a_n,a_m,a_k)= & \phi_n(a_n)\phi_m(a_m)\phi_k(a_k)
{\color{r}\frac{\phi_{nm}(a_n,a_m)}{\phi_n(a_n)\phi_m(a_m)}
\frac{\phi_{mk}(a_m,a_k)}{\phi_m(a_m)\phi_k(a_k)}
\frac{\phi_{nk}(a_n,a_k)}{\phi_m(a_n)\phi_k(a_k)}}\times\nonumber\\
& \exp\left[\sum_{j,z,q=1}^\infty\left<U_m^{j}U_n^{z}U_k^{q}\right>_c
\frac{(ia_m)^{j}(ia_n)^{z}(ia_k)^{q}}{j!z!q!}\right].
\label{ioio}
\end{align}
By generalizing these results to $N$ variables, we find the following 
{\em cluster expansion}
\begin{align}
\phi(a_1,...,a_N)=&\prod_{n=1}^N\phi_n(a_n)
\prod_{n<m}^N\frac{\phi_{nm}(a_n,a_m)}{\phi_n(a_n)\phi_m(a_m)}\times\nonumber\\
&\prod_{n<m<k}^N\frac{\phi_{nmk}(a_n,a_m,a_k)}{\phi_n(a_n)\phi_m(a_m)
\phi_k(a_k)\phi_{nm}(a_n,a_m)\phi_{mk}(a_m,a_k)\phi_{nk}(a_n,a_k)}\times\cdots .
\label{cluster_expansion}
\end{align}
Any truncation of \eqref{cluster_expansion} to a certain order 
in the multi-point characteristic functions $\phi_{lmn\cdots}$ 
yields approximations\footnote{The rationale behind this 
approximation relies on the fact that the joint cumulants 
of $U_k$ often decay with the order, and therefore
the exponential function (see \eqref{ioio}) 
tends to $1$ quickly.}. 
In particular, the second-order truncation
\begin{align}
\phi(a_1,...,a_N)\simeq&\prod_{n=1}^N\phi_n(a_n)
\prod_{n<m}^N\frac{\phi_{nm}(a_n,a_m)}{\phi_n(a_n)\phi_m(a_m)}
\label{cluster_expansion_3}
\end{align}
involves $N(N-1)/2$ functions $\phi_{ln}$ ($\phi_l$ can be 
determined from $\phi_{lm}$),  and it can be defined as {\em  Gaussian approximation}.  
The reason for such definition is that 
if $U_1$,..., $U_N$ are jointly Gaussian,
then \eqref{cluster_expansion_3} is exact 
(the $n$th-order cumulants of a multivariate Gaussian are zero 
if $n\geq 3$). 
We can establish an interesting connection between the networks 
of test functions we discussed in Section \ref{sec:interpolationnodes} 
and the truncation order in \eqref{cluster_expansion}. 

\paragraph{Representation in $S^{(m)}_1$} 
The representation of the Hopf functional in $S^{(m)}_1$ 
(see Eq. \eqref{SNq}) is defined in terms of one-dimensional functions 
\begin{equation}
\Phi([a_j\varphi_j(x)]),\qquad j=1,...,m.
\end{equation}
The corresponding approximation takes the form
\begin{align}
\Phi([a_1\varphi_1(x)+\cdots +a_m\varphi_m(x)])\simeq&
\prod_{j=1}^m\Phi([a_j\varphi_j(x)]).
\label{cluster_expansion_1}
\end{align}
In statistical physics this is known as {\em mean-field approximation}, 
and it relies on a statistical independence hypothesis between the
random variables $U_k(\omega)$ in \eqref{jointCHF}.

\paragraph{Representation in $S^{(m)}_2$} 
\label{sec:S2N}
The representation of the Hopf functional in $S^{(m)}_2$ 
involves two-dimensional functions in the form
\begin{equation}
\Phi([a_j\varphi_j(x)+a_k\varphi_k(x)]), \qquad k,j=1,...,m.
\label{twopoint}
\end{equation}
With this set we  consider the second-order truncation 
of \eqref{cluster_expansion}, i.e., 
\begin{align}
\Phi([a_1\varphi_1(x)+\cdots +a_m\varphi_m(x)])\simeq&\prod_{j=1}^m\Phi([a_j\varphi_j(x)])
\prod_{j<k}^m\frac{\Phi([a_j\varphi_j(x)+a_k\varphi_k(x)])}
{\Phi([a_j\varphi_j(x)])\Phi([a_k\varphi_k(x)])}\nonumber\\
=&\prod_{j<k}^m \Phi([a_j\varphi_j(x)+a_k\varphi_k(x)])
\prod_{j=1}^m\Phi([a_j\varphi_j(x)])^{2-m}.
\label{cluster_expansion_2}
\end{align}
Clearly, the one-dimensional functions $\Phi([a_j\varphi_j(x)])$
can be obtained from \eqref{twopoint} by setting $a_k=0$.

\subsection{Functional Approximation Based on Random Processes }
\label{sec:Stochastic Functional Methods}

So far we discussed representations of nonlinear functionals 
based on power series, Lagrangian interpolants, and tensor methods. 
In this Section we briefly discuss another way of constructing 
polynomial functional expansions by using stochastic processes. 
The method was pioneered by N. Wiener in \cite{Wiener66}, 
Suppose we are given a random function $u(x;\omega)$, 
with know statistical properties and a real-valued nonlinear 
functional $F([\theta])$. If $u(x;\omega)$ is in the domain of 
$F$ then we can evaluate $F([u(x;\omega)])$, which is 
a real-valued random variable. The set 
\begin{equation}
\{u(x;\omega),F([u(x;\omega)])\}
\end{equation}
can be considered as an {\em infinite collection} of 
input-output signals from which we would like to determine 
determine a polynomial approximation of $F$.
The key point is that the stochastic signal $u(x;\omega)$ is 
equivalent to an infinite collection of functions 
that span the domain of $F$, hopefully in a way that is 
sufficient to identify $F([\theta])$. 
To this end, we look for an expansion in the form  
\begin{equation}
F([\theta])=\sum_{k=0}^\infty G_{k}([\theta]),
\label{stochastic_exp}
\end{equation}
where $G_{k}([\theta])$ $k=0, 1, ...$ is a complete 
set of orthogonal polynomial functionals \cite{Wiener66,Segall,Ogura,Ernst}. 
Orthogonality here is relative to the probability measure $P[u]$ of the 
random process $u(x;\omega)$, i.e., relative to the inner product 
\begin{equation}
\left(G_{k}([u]),G_{j}([u])\right)_{dP[u]}=\int_{\Omega} G_{k}([u])G_{j}([u])dP[u].
\label{IP2}
\end{equation}
It was shown in \cite{Wiener66} that if $u(x;\omega)$ is Brownian motion, then $G_k$ are Hermite polynomial functionals, 
and \eqref{stochastic_exp} becomes the celebrated 
{\em Wiener-Hermite expansion} \cite{Wiener66,Cameron,Db_book}.
For completeness, we recall that the first- and 
the second-order Hermite polynomial 
functionals are defined as (see \cite{Wiener66}, p. 32)
\begin{align}
G_{1}([u])=&\int K_1(x_1)du(x_1;\omega),
\label{D10}\\
G_{2}([u])=&\int \int K_2(x_1,x_2)du(x_1;\omega)du(x_2;\omega)-
\int K_2(x_1,x_1) dx_1,
\label{D2}
\end{align}
where $ K_2(x_1,x_2)$ is subject to the normalization condition
\begin{equation}
\int K_2(x_1,x_2)^2 dx_1dx_2=\frac{1}{2}.
\end{equation}
Clearly, if $u(x;\omega)$ is Brownian motion then 
the integrals in \eqref{D10}-\eqref{D2} do not exist 
in the ordinary Stieltjes sense because $u(x;\omega)$ 
is nowhere differentiable. However, we can interpret 
the derivative of $u(x;\omega)$ in a distributional 
sense to obtain a generalized process, i.e., the white noise. 
Wiener has shown that in such generalized setting, the 
integral \eqref{D10}-\eqref{D2} are perfectly well defined 
for kernel functions $K_1(x_1)$ and $K_2(x_1,x_2)$ 
in $L_2$. A detailed mathematical analysis of the 
Wiener-Hermite expansion can be found in \cite{Ernst,Cameron,Wiener66,Poggio}.  

The process of determining the kernels of $G_{1}([\theta])$, 
$G_{2}([\theta])$, etc., is known as {\em identification process} 
in the theory of nonlinear systems, and it has been 
studied extensively for obvious reasons (see \cite{Nelles,Rugh,Shetzen}).
By leveraging on the orthogonality of $G_k$
relative to the inner product \eqref{IP2} one can 
show that (see, e.g., Eq. (4.4) in \cite{Wiener66}) 
\begin{equation}
K_p(x_1,...,x_p)=\frac{1}{p!\epsilon^p}\left(F([u]),G_{p}(H_p,[u])\right)_{dP[u]},\qquad
H_p(z_1,...,z_p)=
\begin{cases}
1 & x_i\leq z_i\leq x_i+\epsilon \\
0 & \textrm{otherwise}
\end{cases}.
\label{Kpp}
\end{equation}
This result relies on the fact that orthogonality of 
$G_k$ and $G_j$ in the sense of \eqref{IP2} is not dependent on 
the actual kernel functions appearing in $G_k$ 
and $G_j$. In fact, if we consider the inner product of \eqref{stochastic_exp} 
with $G_p(H_p,[u])$ -- arbitrary kernel $H_p(x_1,...,x_n)$ -- 
then the only term that survives at the right hand side 
is $\left(G_p(K_p,[u])G_p(H_p,[u])\right)_{dP[u]}$. 
Such inner product can be written as
\begin{equation}
\left(G_p(K_p,[u])G_p(H_p,[u])\right)_{dP[u]}=p!
\int_a^b\cdots \int_a^bK_p(x_1,...,x_p)H_p(x_1,...,x_p)dx_1\cdots dx_p, 
\end{equation}
i.e., all contributions of lower-order terms are identically zero. 
Thus, if we select $H_p$, properly (e.g., $H_p(x_1,...,x_p)=1$ in 
a small hypercube centered at $(x_1,...,x_p)$ with side length 
$\epsilon$ and zero otherwise) then we can 
extract $K_p(x_1,...,x_p)$ as in equation \eqref{Kpp}.

This means that if we know the response $F([u])$ corresponding 
to the stochastic process $u(x;\omega)$, then we can 
identify the kernels \eqref{Kpp} and therefore 
construct the corresponding polynomial functional 
expansion \eqref{stochastic_exp}.
The path integrals in \eqref{Kpp} can be evaluated 
numerically, e.g., by using Monte-Carlo or quasi-Monte Carlo 
methods (see Appendix \ref{sec:functional integrals}). 
To this end, suppose we have available a 
collection of $N_s$ of response-excitation signals 
$\{u(x;\omega_i),F([u(x;\omega_i)])\}$ 
($i=1,..., N$), where $u(x;\omega_i)$ is a
realization of the process $u(x;\omega)$ obtained 
by sampling the probability functional $P[u]$. 
The Monte-Carlo estimate of 
the kernels \eqref{Kpp} is simply 
\begin{equation}
K_p(x_1,...,x_p)\simeq \frac{1}{p! \epsilon^p N}\sum_{i=1}^N 
F([\xi(x,\omega_i)])G_{k}(H_p,[\xi(x,\omega_i)])\qquad p=0,1,...
\label{KppMC}
\end{equation}
where $\xi(x,\omega_i)$ is a realization of the 
distributional derivative of $u(x;\omega_i)$. 
In particular, if $u(x;\omega)$ is
Brownian motion then $\xi(x;\omega)$ is 
{\em spatial white noise}\footnote{A simple numerical 
approximation of white 
noise $\xi(x;\omega)=du(x;\omega)/dx$ can be obtained as 
\begin{equation}
\xi_{\Delta x}(x;\omega) = \frac{u(x+\Delta x;\omega)-u(x;\omega)}{\Delta x},
\qquad \textrm{($u$ Brownian motion).}
\end{equation}
This is a Gaussian process with zero mean, variance $1/\Delta x$ and covariance
\begin{equation}
\left<\xi_{\Delta x}(x;\omega)\xi_{\Delta x}(y;\omega)\right>=
\begin{cases}
\displaystyle\frac{1}{\Delta x}\left(1-\frac{|x-y|}{\Delta x}\right)& \textrm{if $|x-y|\leq \Delta x$}\\
0 & \textrm{otherwise}
\end{cases}
\end{equation}
The Fourier transform of $\left<\xi_{\Delta x}(x;\omega)\xi_{\Delta x}(y;\omega)\right>$ 
yields a power spectrum which is not flat as it is supposed to be for white noise, 
but decays at sufficiently large frequencies.
}. 
In this case, we have that the first two 
functionals $G_1(H_1,[\xi])$ and $G_2(H_2,\xi)$ 
in \eqref{KppMC} are given by 
\begin{align}
G_{1}(H_1,[\xi])&=\epsilon \xi(x_1;\omega), \nonumber\\
G_{2}(H_2,[\xi])&=\epsilon^2 \left(\xi(x_1;\omega)
\xi(x_2;\omega)-1\right), \nonumber
\end{align}
and therefore
\begin{align}
K_1(x_1)\simeq &\frac{1}{ N}\sum_{i=1}^N F([\xi^{(i)}])
h_1(\xi^{(i)}(x_1;\omega)),\label{KppMC1}\\
K_2(x_1,x_2)\simeq & \frac{1}{ N}\sum_{i=1}^N F([\xi^{(i)}])
h_2(\xi^{(i)}(x_1;\omega),\xi^{(i)}(x_2;\omega)).
\label{KppMC2}
\end{align}
Here $h_1(x_1)=x_1$ and $h_2(x_1,x_2)=(x_1x_2-1)/2$ are the 
classical first- and second-order Hermite polynomials. 
A very insightful discussion on how to compute the kernels 
$K_j$ by using white noise is given by Rugh in \cite{Rugh}, \S 7.4.
The expansion \eqref{stochastic_exp} obtained in this way 
can, in principle, be used to compute the value of the functional 
$F([\theta])$ corresponding to any deterministic 
function $\theta(x)$. However, since the 
kernels \eqref{Kpp} are built upon stochastic processes 
and their averages, it is not clear in what 
sense \eqref{stochastic_exp} will converge for 
deterministic input functions $\theta(x)$. 
This question was addressed by 
Palm and Poggio in \cite{Poggio} (Theorems 4 and 5),
where necessary and sufficient conditions for pointwise 
convergence of Wiener-Hermite expansions are provided
(see also \cite{Ernst}).
The functional expansion in terms of orthogonal 
polynomial functionals can be generalized
to random processes other than Brownian motion and 
white noise\footnote{The construction 
of such generalized expansion proceeds as follows. 
Starting from the constant functional $F_0$ we look for 
\begin{equation}
F_1([u])=C_1 \int K_1(x) du(x;\omega)+ F_0,
\end{equation}
where $u(x;\omega)$ is a generalized random process, and 
we make it orthogonal to $F_0$ in the sense 
of $dP[u(x;\omega)]$. This yields $C_1$. 
Then we construct $F_2([u])$ and we make it 
orthogonal to both $F_1([u])$ and $F_0$. 
All these conditions are ultimately expressed 
analytically in terms of multi-point averages of 
the random process $u(x;\omega)$.}.
However, one has to be very careful when 
performing such generalizations. There are indeed 
random processes that do not allow for a complete 
representation of the functional $F([\theta])$ 
(see \cite{Ernst,Segall,Ogura} for details).

\section{Functional Differential Equations}
\label{sec:FDEs approximation}

A functional differential equation (FDE) is an equation involving  
a nonlinear functional (i.e., a nonlinear operator), 
derivatives with respect to functions (functional derivatives) 
and derivatives with respect to independent variables, e.g., 
space and time variables. 
In this report we will study linear FDEs in the form 
\begin{equation}
\frac{\partial F([\theta],t)}{\partial t}=L([\theta],t) F([\theta],t)+H([\theta],t)
\qquad 
F([\theta],0)=F_0([\theta]),
\label{linfde00}
\end{equation}
where $F_0([\theta])$ is an initial condition,  
$H([\theta],t)$ is a known forcing term 
and $L([\theta],t)$ is a linear operator 
in the space of functionals. This class of FDEs 
is very broad, and it includes many well-known 
functional equations of theoretical and quantum physics.
The solution to initial value problem \eqref{linfde00} 
(assuming it exists) depends on the 
{\em initial condition} $F_0$ as well as on the 
choice of the {\em function space} $D(F)$ (domain 
of the functional). The following 
examples clarifies this point. 

\vs
\noindent
{\em Example 1 (Functional Advection Equation):}
Consider the following FDE
\begin{equation}
 \frac{\partial F([\theta],t)}{\partial t}+\int_0^{2\pi}\theta(x)
 \frac{\partial }{\partial x}\frac{\delta F([\theta],t)}{\delta \theta(x)}dx=0, \qquad 
 F([\theta],0)=F_0([\theta]).\label{FDEln}
\end{equation}
Clearly, this equation is in the form \eqref{linfde00}, with 
$H([\theta],t)=0$ and 
\begin{equation}
 L([\theta],t)=-\int_0^{2\pi}\theta(x)
 \frac{\partial }{\partial x}\frac{\delta }{\delta \theta(x)} dx.
 \label{FDOP}
\end{equation}
The solution to \eqref{FDEln} depends 
on the {\em initial condition} $F_0([\theta])$ 
as well as on the choice of 
the {\em function space} $D(F)$. To show this, 
let us set 
\begin{equation}
F_0([\theta])=\int_0^{2\pi}K(x)\theta(x)dx,
\label{ic0}
\end{equation}
where $K(x)$ is a given kernel function\footnote{
The solution functional $F([\theta],t)$ corresponding to 
the initial condition \eqref{ic0} is linear and 
homogeneous for all $t\geq0$, i.e., it is in the form
\begin{equation}
F([\theta],t)=\int_0^{2\pi}R(x,t)\theta(x)dx,\qquad R(x,0)=K(x).
\label{K1lin}
\end{equation}
To see this, it is sufficient 
to write the first-order Euler scheme 
\begin{equation}
 F([\theta],\Delta t)=F([\theta],0) -\Delta t \int_0^{2\pi}\theta(x)
 \frac{\partial K(x,0)}{\partial x}dx.
\end{equation}
and note that the term within the integral at the right hand side 
is a linear functional of $\theta$. This implies that $F([\theta],\Delta t)$ is 
a linear functional of $\theta$. By applying this argument over 
and over we see that  $F([\theta],t)$ is a linear functional of 
$\theta$ for all $t\geq0$, i.e., it is in the form \eqref{K1lin}.
A substitution of \eqref{K1lin} into \eqref{FDEln} yields
\begin{equation}
\int_0^{2\pi}\theta(x)\left(\frac{\partial R(x,t)}{\partial t}+
\frac{\partial R(x,t)}{\partial x}\right)dx=0,\qquad \textrm{i.e.} \qquad
\frac{\partial R(x,t)}{\partial t}+
\frac{\partial R(x,t)}{\partial x}=0.\label{adv4}
\end{equation}
Therefore, \eqref{K1lin} is a solution to \eqref{FDEln}-\eqref{ic0} 
if $R(x,t)$ solves a simple linear advection equation on the real line, 
with initial condition $R(x,0)=K(x)$. 
On the other hand, if the initial condition $F_0([\theta])$ is 
a nonlinear functional of $\theta$, then $F([\theta],t)$ 
is nonlinear functional of $\theta(x)$. We will discuss this 
case extensively in Section \ref{sec:numerical results functional equations}).
}. 
As we will see in Section \ref{sec:ADVR}, 
if we solve \eqref{FDOP} in the space of periodic functions 
\begin{equation}
D(F) =\{\theta(x)\in C^{(\infty)}([0,2\pi])\,|,\ \theta(0)=\theta(2\pi)\},  
\end{equation}
then we obtain the solution 
\begin{equation}
F([\theta],t)=F_0(\theta),
\end{equation}
i.e., the constant functional. On the other hand, 
if we consider the function space 
\begin{equation}
D(F) =\{\theta(x)\in {C}^{(\infty)}([0,2\pi])\,|,\ \theta(x)=0\},  
\end{equation}
then the solution to \eqref{FDOP} is 
\begin{equation}
F([\theta],t)=
\begin{cases}
\displaystyle \int_t^{2\pi}K(x) \theta(x)dx & t\leq 2\pi\\
0  & t > 2\pi
\end{cases}.
\end{equation}

\paragraph{Existence and Uniqueness of the Solution} 
A fundamental question that has lasted over the years 
is whether the solution to linear FDEs are unique 
or not, given the initial state $F_0([\theta])$ 
and the function space $D(F)$. This is an unsolved mathematical 
problem which we hope will be addressed systematically 
in near future. At today, there are few general theorems and 
results on the existence and the uniqueness of solutions to
functional differential equations (see, e.g., 
\cite{Foias,Hale,Wu}).

Before addressing the question of numerical approximation 
functional differential equations, it is useful 
to show how such equations look like and, more importantly, 
how they arise in the context of well-known mathematical theories.

\subsection{Variational Form of PDEs}
\label{sec:variational}
Perhaps, the simplest example
of a (algebraic) functional equation is the 
variational form of a PDE. 
To show how such equation 
looks like, consider the scalar PDE 
\begin{equation}
N(u)=0,\qquad
\label{PDE}
\end{equation}
where $N$ is a nonlinear differential operator subject 
to appropriate initial/boundary conditions. For example, 
\begin{equation}
N(u)=\frac{\partial u}{\partial t} +u\frac{\partial u}{\partial x}-
\frac{\partial^2 u}{\partial x^2},\qquad u(x,0)=u_0(x), \qquad u(0,t)=u(L,t).
\label{firsteq}
\end{equation}
We multiply \eqref{PDE} by the test function $\theta(x)$ 
in the space 
\begin{equation}
D(F)= \left\{
\theta \in C^{(2)} ([0,2\pi])\,|\, \theta(0)=\theta(2\pi)\right\}
\end{equation}
and integrate over $[0,2\pi]$ we obtain 

\begin{equation}
F([\theta])=\left(N(u),\theta\right) =0, \qquad \forall\theta\in D(F).
\end{equation}
This is the starting point of the well-known method 
of {\em weighed residuals} \cite{Finlayson}, 
from which classical Galerkin, 
collocation, least-squares and finite volumes 
schemes can be derived (see \cite{GKSS_2005}, p.18).
As we will see in Section \ref{sec:rep-Hopf}, to identify 
the functional $F$ in this case it is sufficient  
to test it relative to a set of linearly independent 
functions $\theta_i$ ($i=1, 2,...$), e.g., an 
orthonormal basis of $D(F)$.

\subsection{Schwinger-Dyson Equations}
The Schwinger-Dyson equations are functional differential 
equations for the generating functional of a field theory. 
They arise in both classical statistical physics as 
well as in quantum field theory.  Hereafter we review the main 
aspects of such equations.

\paragraph{Statistical Physics}
The functional integral approach to classical statistical dynamics 
\cite{Jensen,Phythian,Jouvet,Langouche} allows us derive formally exact evolution equations for phase space functions 
such as the mean and the correlation function of 
the solution to SODEs and SPDEs.  
The standard approach relies on a {\em generating functional} $Z$.
For stochastic dynamical systems in the form
\begin{equation}
 \frac{d{\bm \psi}(t)}{dt}=\bm \Lambda(\bm \psi(t),t)+\bm f(t;\omega),\label{eqofm}
\end{equation}
$Z$ can be written as a functional integral 
(see \cite{Phythian,Jensen,Justin,Amit,Stratonovich}) 
\begin{equation}
Z([\bm \xi,\bm \eta])=Z_0\int \int\mathcal{D}[\bm \psi]
\mathcal{D}[\bm \chi]A([\bm \psi,\bm \chi])
\exp\left(\int_0^t d\tau (\bm \xi(\tau)\cdot\bm \psi(\tau)+
\bm \eta(\tau)\cdot\bm \chi(\tau))\right),
\label{gf}
\end{equation}
where {\color{r} 
\begin{equation}
\frac{1}{Z_0} = \int \int\mathcal{D}[\bm \psi]\mathcal{D}[\bm \chi] 
A([\bm \psi,\bm \chi]),
\end{equation}}
and
\begin{equation}
A([\bm \psi,\bm \chi]) = C([\bm \chi])\exp\left(-\frac{1}{2}\int_0^t d\tau 
\nabla\cdot \bm \Lambda(\bm \psi(\tau),\tau)-i\int_0^t d\tau \bm \chi(\tau)\cdot 
\left[ \frac{d{\bm \psi}(\tau)}{dt}-\bm \Lambda(\bm \psi(\tau),\tau)\right]\right).
\end{equation}
Here $C([\bm \chi])$ denotes the (known) characteristic 
functional of the external random noise $\bm f(t;\omega)$. 
Clearly, if we have available the solution to the stochastic dynamical 
system \eqref{eqofm} then we can construct 
the functional $Z([\bm \xi,\bm \eta])$,  and from it 
compute any statistical property we may be interested in. 
On the other hand, it is straightforward to show 
that $Z([\bm \xi,\bm \eta])$ satisfies a 
coupled system of {\em linear} functional differential equations, 
known as {\em Schwinger-Dyson equations} in quantum field theory \cite{Itzykson}.
The equations are in the form
\begin{align}
\frac{\partial}{\partial \tau}\left(\frac{1}{i}\frac{\delta Z}{\delta \xi_k(\tau)}\right)&=
\eta_k(\tau)Z+\Lambda_k\left(\frac{1}{i}\frac{\delta}{\delta \bm \xi(\tau)},\tau\right)Z-
iD_k\left(\left[\frac{1}{i}\frac{\delta}{\delta\bm \eta(\tau)}\right],\tau\right)Z,
\label{f1}\\
\frac{\partial}{\partial \tau}\left(\frac{1}{i}\frac{\delta Z}{\delta \eta_k(\tau)}\right)&=
-\xi_k(\tau)Z+i\frac{\delta}{\delta \eta_j(\tau^+)}\frac{\partial \Lambda_k}{\partial \psi_j}
\left(\frac{1}{i}\frac{\delta}{\delta \bm \xi(\tau)},\tau\right)Z
\label{f2},
\end{align}
where
\begin{equation}
D_i([\bm \chi],\tau)=\frac{\delta}{\delta \chi_i(\tau)}\ln C[\bm \chi].
\end{equation}
The quantities $\delta/\delta \xi_k$ and $\delta/\delta \eta_k$ are 
first-order functional derivative operators, defined in Appendix
\ref{app:functional derivatives}. 
By solving \eqref{f1}-\eqref{f2} we can identify 
$Z([\bm \xi,\bm \eta])$ without any knowledge of the 
stochastic process $\bm \psi(t;\omega)$.    
By generalizing \eqref{gf}, it is possible to derive a functional integral 
formalism to classical statistical dynamics yielding Schwinger-Dyson equations 
for generating functionals associated with SPDEs (see, e.g., \cite{Jensen,Martin}).
In particular, if the SPDE admits an action functional $A[\phi]$ 
(see \cite{Daniele_JMathPhys,Eyink,Gomes,Funaki,Yasue,Eyink_1996,Amit}), 
then the construction of the generating functional as well as the derivation 
of the corresponding Schwinger-Dyson equations are standard.
In this setting, the Schwinger-Dyson functional 
differential equations provide a non-perturbative 
formulation of the problem of computing the statistical 
properties of nonlinear random systems, including 
stochastic ODEs and stochastic PDEs.

{\color{r}
\paragraph{Quantum Field Theory}
In quantum field theory, the 
Schwinger-Dyson equations govern the dynamics 
of the Green functions and they characterize
the propagation field interactions \cite{Easther,Justin,Weinberg}. 
Such functional equations can be employed in 
a perturbation setting \cite{Okopinska} 
(weak coupling regime), but they 
show their true strength in the strong coupling 
regime \cite{Swanson,Guralnik1}. 
The starting point to derive the Schwinger-Dyson 
equations is the generating functional of 
the correlation functions (Green functions), which 
can be often expressed as a functional integral
\begin{equation}
\displaystyle
Z([j(\bm x)]) = \int \mathcal{D}[\phi] e^{i A([\phi]) + i 
\int j(\bm x)\phi(\bm x)d\bm x}. 
\label{GF0}
\end{equation}
where $A([\phi])$ is the {\em action functional}. 
In the context of quantum $\phi^4$-theory \cite{KleinertPhi4} 
the action $A$ for a field with mass $m$ is given by
\begin{equation}
A([\phi])  =  \int\left[  \frac{1}{2}\left(\nabla \phi(\bm x)\right)^2- 
\frac{1}{2} m^2 \phi^2(\bm x) - 
\frac{\lambda}{4!}\phi^4(\bm x)\right] d\bm x.
\label{phi4}
\end{equation}
By employing \eqref{GF0}, 
we can express the Green functions of the 
quantum field theory at any order by functional 
differentiation\footnote{The so-called {\em connected} Green 
functions of the field theory are obtained by functional differentiation 
of  $\log Z([j])$. The main reason for such definition is that if 
we represent the functional derivatives of $\log Z$ 
in terms of Feynman diagrams then only the 
connected diagrams contribute to the expansion.}, i.e., 
\begin{align}
 Z([0]) G(\bm x_1,...,\bm x_n)= &\displaystyle\int\mathcal{D}[\phi]\phi(\bm x_1)\cdots \phi(\bm x_n) 
e^{i A([\phi])},\nonumber\\
=&\frac{1}{i^n}\left.\frac{\delta^n Z([j(\bm x)])}{\delta j(\bm x_1)\cdots \delta j(\bm x_n)}\right|_{j(\bm x)=0}.
\end{align}
Roughly speaking, the correlation functions are 
averages of products of fields $\phi(x_1)\cdots \phi(x_n)$ 
with respect to the functional measure $\exp{iA([\phi])}$. 
From a physical viewpoint they represent  the transition 
amplitude for the propagation of a particle or 
excitation between different points 
in space-time. 
The generating functional \eqref{GF0} satisfies the 
Schwinger-Dyson  equation
\begin{equation}
\frac{\delta A([\phi])}{\delta \phi(\bm x)}\left(\left[-i\frac{\delta}{\delta j
(\bm x)}\right]\right)Z([j])+j(\bm x)Z([j])=0.
\end{equation}
For example, if the action $A$ is the one given 
in \eqref{phi4}, then the Schwinger-Dyson equation 
takes the form  
\begin{equation}
 \square
\frac{\delta Z([j])}{\delta j(\bm x)}+ 
 m^2 \frac{\delta Z([j])}{\delta j(\bm x)}-
\frac{\lambda}{3!} 
\frac{\delta^3 Z([j])}{\delta j(\bm x)^3 }-i j(\bm x) Z([j])=0.
\end{equation}
A substitution of the functional Taylor expansion 
of $Z([j])$ into this equation yields an 
infinite-dimensional coupled PDE system for the 
correlation functions $G(\bm x_1,...,\bm x_n)$. 
Here $\bm x_i$ is a quadruple of coordinates. 
Solution methods for the Schwinger-Dyson equations 
rely on truncated series expansions \cite{Okopinska,Bender}, 
or renormalized expansions of the generating functional 
$Z([\eta])$ (see \cite{Kleinert1} and \cite{Amit}, p. 385), 
or numerical algorithms \cite{Easther,Guralnik1,Guralnik2}.
}

\subsection{Hopf Characteristic Functional Equations}
\label{sec:hopf_equations}
The Hopf characteristic functional of a random field 
is the functional Fourier transform of the probability 
density functional (see \cite{Rosen1,Venturi_PRS} and 
Appendix \ref{app:functional fourier transform}). 
To introduce this mathematical object in a simple way,  
let us consider an integrable random function $u(x;\omega)$ 
on an interval $x\in[a,b]$.
The Hopf functional associated with $u(x;\omega)$ is defined as
\begin{equation}
 \Phi([\theta(x)])=\left<\exp \left[ i\int_a^b u(x;\omega)\theta(x)dx\right]\right>,
\label{hopf_functional}
 \end{equation}
where $\theta(x)$ is a deterministic function (test function), 
$i$ is the imaginary unit and the average is 
defined as a functional integral over the probability 
functional of $u(x;\omega)$.
Equation \eqref{hopf_functional} assigns to each 
function $\theta(x)$ a complex number $\Phi([\theta(x)])$ 
(see Figure \ref{fig:1}).
Similarly to the probability density functional, the 
Hopf characteristic functional \eqref{hopf_functional} encodes 
the {\em full statistical information} of the random function 
$u(x;\omega)$, including multi-point moments, joint 
characteristic functions and probability density 
functions (see \cite{Monin1,Monin2}). 
If we consider instead of a random function  $u(x;\omega)$ 
a random vector field $\bm u(\bm x,t;\omega)$, e.g., a 
stochastic solution to the Navier-Stokes equations, then 
we define\footnote{Lewis and Kraichnan \cite{Lewis,Rosen1} introduced a space-time generalization 
of \eqref{hopf_functional2}, namely
\begin{equation}
\Phi([\bm \theta(\bm x,t)])=\left<\exp \left[ i\int_V\int_T \bm u(\bm x,t;\omega)\cdot \bm \theta(\bm x,t)
d\bm x dt\right]\right>.
\label{hopf_functional1}
\end{equation}
This functional allows us to determine joint multi-point statistics 
of the random field $\bm u(\bm x,t;\omega)$ at different times.} 
\begin{equation}
\Phi([\bm \theta(\bm x)],t)=\left<\exp \left[ i\int_V 
\bm u(\bm x,t;\omega)\cdot \bm\theta(\bm x)
d\bm x\right]\right>,
\label{hopf_functional2}
\end{equation}
where $V$ is a spatial domain in $\mathbb{R}^d$ ($d=2,3$).
By the Riemann-Lebesgue lemma we have that 
\begin{equation}
 \Phi([\bm \theta(\bm x)],t)\rightarrow 0\quad \textrm{as}\quad
 \left\|\bm \theta(\bm x)\right\|\rightarrow \infty,
\end{equation}
with a rate that depends on the regularity of 
the underlying probability density functional.
The derivation of the Hopf characteristic functional equation 
is relatively straightforward if the random field 
$\bm u(\bm x,t;\omega)$ satisfies a nonlinear PDE with 
polynomial nonlinearities. Hereafter we provide some 
examples.

\paragraph{Burgers-Hopf Equation}
Consider the Burgers equation
\begin{equation}
\frac{\partial u}{\partial t}+u\frac{\partial u}{\partial x}=
\nu\frac{\partial^2u}{\partial x^2},
\label{BP}
\end{equation}
in a periodic spatial domain $[0,2\pi]$, 
and let the initial condition $u_0(x,\omega)$ be 
random. 
By differentiating the Hopf functional 
\begin{equation}
\Phi([\theta(x)],t)=\left<\exp\left[i\int_{0}^{2\pi}
\theta(x)u(x,t;\omega)dx\right]\right>
% \label{Hopf1}
\end{equation}
with respect to time 
{\color{r}
and using \eqref{BP} we obtain 
\begin{align}
\frac{\partial \Phi([\theta],t)}{\partial t} = &
i \int_{0}^{2\pi} \theta(x) 
\left<\frac{\partial u(x,t;\omega) }{\partial t} 
\exp\left[i\int_{0}^{2\pi}u(x,t;\omega) \theta(x) dx  \right] 
\right>dx\nonumber\\
=& i\int_{0}^{2\pi} \theta(x) 
\left<\left(-\frac{1}{2}\frac{\partial (u(x,t;\omega)^2)}{\partial x}
+\nu\frac{\partial ^2u(x,t;\omega)}{\partial x^2}\right)
\exp\left[i\int_{0}^{2\pi}u(x,t;\omega) \theta(x) dx  \right] 
\right>dx
\label{PdFeq1}
\end{align}
i.e., 
\begin{align}
 \frac{\partial \Phi([\theta],t)}{\partial t}=
 \int_a^{b} \theta(x)\left(\frac{i}{2}\frac{\partial}{\partial x}
 \frac{\delta^2 \Phi([\theta],t)}{\delta\theta(x)^2}+
\nu\frac{\partial^2}{\partial x^2}\frac{\delta \Phi([\theta],t)}{\delta\theta(x)}\right)dx. 
\label{burgersfunctional}
\end{align}}

\noindent
This equation is known as Burgers-Hopf equation and it 
has been the subject of numerous investigations 
(see, e.g., \cite{Ahmadi,Hosokawa2,Monin2}).

\paragraph{Navier-Stokes-Hopf Equation} 
The problem of determining the evolution of the Hopf 
characteristic functional for the 
Navier-Stokes equations 
\begin{equation}
 \frac{\partial \bm u}{\partial t}+(\bm u\cdot \nabla)\bm u=-\nabla p +\nu \nabla^2 \bm u, 
 \qquad \nabla\cdot \bm u=0
\end{equation}
was deemed by Monin and Yaglom as the most compact formulation of the {\em turbulence 
problem}, which is the problem of determining the statistical 
properties of the velocity and pressure field given statistical 
information on the initial condition. The Hopf functional 
differential equation corresponding to the Navier-Stokes 
equations is derived in \cite{Monin2} Ch. 10 (see also \cite{Beran}, \S 3.1.4), and it is hereafter summarized for convenience
\begin{equation}
\frac{\partial \Phi([\bm \theta],t)}{\partial t}=
\sum_{k=1}^3\int_V\theta_k(\bm x)\left(i \sum_{j=1}^3\frac{\partial }{\partial x_j}
\frac{\delta^2 \Phi([\bm \theta],t)}{\delta \theta_k(\bm x)\delta\theta_j(\bm x)}
+\nu \nabla^2\frac{\delta \Phi([\bm \theta],t)}{\delta \theta_k(\bm x)}\right)d\bm x.
\label{hopfns1}
\end{equation}
Here $V$ is a periodic three-dimensional box and 
$\bm \theta(\bm x)$ is chosen in a divergence free 
space of test functions 
(see Section \eqref{sec:Navier-Stokes-Hopf}).
Functional formulations of non-isothermal turbulent 
reactive flow have been also considered, leading to more complicated Hopf equations \cite{Dopazo}. 
As we will see in Section \ref{app:equivalenceHopf}, 
Hopf equations are equivalent to PDEs in an infinite 
number of variables or to infinite-dimensional 
systems of coupled PDEs, e.g., the 
Monin-Lundgren-Novikov hierarchy 
\cite{Hosokawa,Montgomery,Wacawczyk1,Lundgren}. 
It is interesting to note that the structure 
of the Hopf equation somehow resembles the weak form 
of a PDE. However, there is a remarkable 
difference: in the Hopf equation, both the 
solution and its functional derivatives depend on 
the test function. In other words, the 
test function appears in the functional 
equation a nonlinear way. In addition, the 
equation involves derivatives with respect to 
functions (functional derivatives) which are 
not present in classical PDEs. 

\vs 
\noindent
It is straightforward to derive 
Hopf functional differential equations corresponding 
to linear or nonlinear evolution PDEs with polynomial 
nonlinearities, such as the Kuramoto-Sivashinsky equation, the nonlinear wave equation, and Maxwell's equations subject 
to random boundary or random initial conditions. 

{\color{r}
\paragraph{Hopf Equations Defining Random Processes}
Hopf functional equations arise naturally also in the context of 
random processes. For example, the Hopf equation defining 
the characteristic functional of a zero-mean Gaussian 
process is (see \cite{Klyatskin1}, p. 61)
\begin{equation}
\frac{d}{dt}\Phi([\theta(t)],t)=-\frac{1}{2}\Phi([\theta(t)],t)\theta(t)
\int_0^t C(t,\tau)\theta(\tau)d\tau,\qquad \Phi([\theta(t)],0)=1,
\label{FDEGauss}
\end{equation}
where $C(t,\tau)$ is the covariance function of the process. 
In fact, the solution to \eqref{FDEGauss} is the well-known Hopf functional
\begin{equation}
\Phi([\theta(t)],t)=\exp\left[-\frac{1}{2}\int_0^t\int_0^\tau C(\tau_1,\tau_2)
\theta(\tau_1)\theta(\tau_2)\right].
\end{equation}
Similar equations can be derived for orther processes, such as 
the telegrapher's random process and general Markov processes.
Note that \eqref{FDEGauss} does not involve any functional 
derivative. 
}

\subsection{Probability Density Functional Equations}
\label{sec:PDENFUNEQ}
In statistical mechanics, a system of $n$ particles can 
be described by the joint probability density function
$p( \bm x_1, ...,\bm  x_n,  \bm v_1, ..., \bm v_n,t)$ 
where $( \bm x_i,\bm v_i)$ denotes the position 
and the velocity of the $i$-th particle. 
Similarly, the phase space associated 
with any finite-dimensional approximation of the solution to a 
stochastic PDE can be described by the joint probability 
density function of the corresponding phase variables, e.g., 
the Fourier coefficients of the series expansion of the 
solution \cite{Edwards, Herring,Montgomery}.
Consider a scalar random field $u(x,t;\omega)$.
The full statistical information of the random field 
$u(x,t;\omega)$ at time $t$ is encoded in the 
probability density functional\footnote{
The probability density functional $P([a(x)],t)$ 
allows us to compute all moments, cumulants   
and joint PDFs by using functional integration. For example, 
\begin{equation}
p(b_1,b_2,t)=\int \delta(b_1-a(x_1))\delta(b_2-a(x_2))
P([a(x)],t)\mathcal{D}[a(x)].
\end{equation}
Computing functional integrals often requires a careful definition 
of the integration measure as it may be possible to run into convergence issues (see Appendix \ref{sec:functional integrals}, 
and \cite{Beran} \S 2.2.4). 
}  
\begin{align}
P([a(x)],t)=&\left<
\delta\left[a(x)-u(x,t;\omega)\right]\right>, 
\label{probability_functional}
\end{align}
where $\delta[\cdot]$ denotes a {\em Dirac delta functional} 
(see Appendix \ref{app:functional fourier transform}), 
and the average $\langle \cdot\rangle$ is a 
functional integral over the probability 
measure of $u(x,t;\omega)$.
Any well-defined nonlinear stochastic PDE for 
the scalar random field $u(x,t;\omega)$
can be rewritten as a {\em linear functional differential 
equation} for $P([a(x)],t)$. 
Such equation can be obtained by taking an appropriate 
continuum limit of a finite-dimensional joint probability density 
equation \cite{Venturi_PRS,Venturi_JCP}. 
In practice, we can replace the discrete set of 
variables $\{a_1=a(x_1), a_2=a(x_2) ,... , a_n=a(x_n)\}$ 
in the joint PDF equation by a continuous set 
denoted by a continuously indexed set $a(x)$, i.e., a 
function. Partial derivatives with respect to 
$a_i=a(x_i)$ can then be replaced by functional 
derivatives with respect to $a(x)$ (see Section 
\ref{sec:Finite_Dim_Approx}), etc. 
An alternative method to derive the probability 
functional equation was proposed by Beran in \cite{Beran}. 
The derivation parallels classical the Liouville theory, 
where a nonlinear dynamical system is converted into 
a linear PDE for the joint probability density function of 
of the state vector \cite{Venturi_PRS}. To convert 
a nonlinear PDE into a linear equation for the probability 
density functional, we need to move one level up and 
look for a functional differential equation.
{\color{r}
Perhaps, the simplest way to derive a probability density 
functional equation is to inverse Fourier transform the 
corresponding Hopf equation, and use functional integration by 
parts. To illustrate the procedure, we first differentiate the 
Hopf functional $\Phi([\theta],t)$ with respect to 
time to obtain\footnote{
\color{r}
In equation \eqref{43R}, we employed 
the identity
\begin{equation}
\int \exp\left[i\int_{0}^{2\pi}\theta(x)u(x,t;\omega)dx  \right]P([u_0])\mathcal{D}[u_0] = \int \exp\left[i\int_{0}^{2\pi}\theta(x)a(x) dx \right]P([a],t)\mathcal{D}[a],
\end{equation}
where $P([u_0]$ and $P([a],t)$ are, respectively, 
the probability functionals of $u_0(x;\omega)$ 
(initial condition) and $u(x,t;\omega)$.
}
\begin{align}
\frac{\partial \Phi([\theta],t)}{\partial t} = &
\int\exp\left[i \int \theta(x)a(x)dx\right] 
\frac{\partial P([a(x)],t)}{\partial t}\mathcal{D}[a].
\label{43R}
\end{align}
Then we set the equality with the functional equation 
that defines the evolution of $\Phi([\theta],t)$ 
for a particular nonlinear PDE. 
For instance, if we consider the Burgers equation, then we have 
that $\partial \Phi/\partial t$ is given by 
equation \eqref{burgersfunctional}, i.e., 
\begin{align}
&\int\exp\left[i \int \theta(x)a(x)dx\right] 
\frac{\partial P([a(x)],t)}{\partial t}\mathcal{D}[a] =\nonumber\\
 &\int_{0}^{2\pi} 
\int \left(-\frac{1}{2}\frac{\partial (a(x)^2)}{\partial x}
+\nu\frac{\partial ^2 a(x)}{\partial x^2}\right)
\frac{\delta}{\delta a(x)}\left(\exp\left[i\int_{0}^{2\pi}a(x) \theta(x) dx  \right]\right)dx P([a],t)
\mathcal{D}[a]
\end{align}
Performing a functional integration by parts and assuming that 
the boundary terms are zero yields 
\begin{align}
\frac{\partial P([a(x)],t)}{\partial t} = -\int_{0}^{2\pi} 
\frac{\delta}{\delta a(x)}\left(\left[-a(x)\frac{\partial a(x)}{\partial x}
+\nu\frac{\partial ^2 a(x)}{\partial x^2}\right] P([a],t)\right)dx.
\label{PrDFeq2}
\end{align}

\paragraph{The Method of Continuum Limits}
The formal procedure to derive probability density functional 
equations we just described can be justified 
in a finite dimensional setting. To this end, let us 
consider the one-dimensional diffusion problem   
\begin{equation}
 \frac{\partial u(x,t;\omega)}{\partial t}=\frac{\partial^2 u (x,t;\omega)}{\partial x^2}, 
 \qquad u(x,0;\omega)=u_0(x;\omega)\quad \textrm{(random)} 
 \label{1dheat}
\end{equation}
in the real line $x\in \mathbb{R}$. The probability 
density functional equation of the solution is a subcase of 
equation \eqref{PrDFeq2}, namely
\begin{align}
\frac{\partial P([a(x)],t)}{\partial t} = -\int_{-\infty}^{\infty} 
\frac{\delta}{\delta a(x)}\left(
\frac{\partial^2 a(x)}{\partial x^2} P([a],t)
\right)dx.
\label{PrDFeq3}
\end{align}
To derive \eqref{PrDFeq3}, we first 
discretize \eqref{1dheat} in space, e.g., 
on an spatial grid with evenly spaced 
nodes $x_j$ ($j=1,...,n$), the spacing between 
the nodes being $\Delta x$. 
If we use second-order finite differences 
we obtain  
\begin{equation}
 \frac{d u(x_k,t,\omega)}{dt}= \frac{u(x_{k+1},t;\omega)-2u(x_k,t;\omega)+u(x_{k-1},t;\omega)}{\Delta x^2}.
\end{equation}
Now, let $p(a_1,...,a_n,t)$ be the joint PDF of 
$\{u(x_1,t;\omega),...,u(x_n,t;\omega)\}$, i.e.,
\begin{equation}
p(a_1,...,a_n,t) = \left<\prod_{k=1}^n \delta(a_k-u(x_k,t;\omega)\right>. 
\end{equation}
We think of $a_k$ as the value of some 
function $a(x)$ at $x_k$, that is $a_k=a(x_k)$ 
(see Figure \ref{fig:ProbDF}).
\begin{figure}[t]
\centerline{
\includegraphics[height=3.7cm]{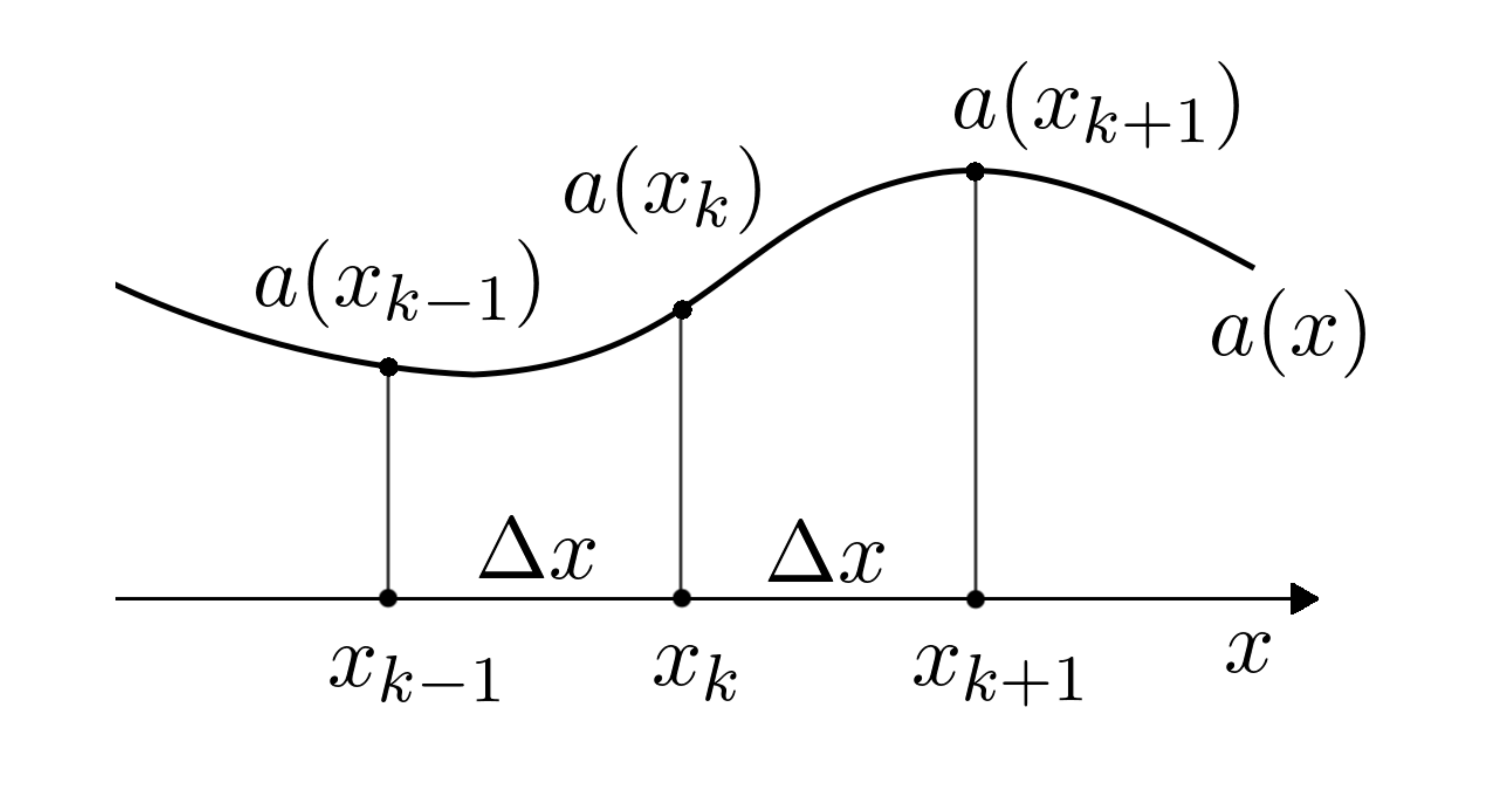}
}
\caption{Sketch of the variables 
$a(x_k)$ representing the random field $u(x,t;\omega)$ at 
$x_k$. When we send $\Delta x$ to zero, 
the number of variables $a(x_k)$ goes to infinity and 
the joint probability density function becomes a probability density 
functional (see equation \eqref{probability_functional}).} 
\label{fig:ProbDF}
\end{figure}
By using well-known identities involving the Dirac delta function \cite{Khuri} it can be shown that 
\begin{align}
&\frac{a(x_{k+1})-2a(x_k)+a(x_{k-1)}}{\Delta x^2}
\left<\prod_{k=1}^n \delta\left(a(x_k)-u(x_k,t;\omega)\right)\right>
=\nonumber\\
&\left<\frac{u(x_{k+1},t;\omega)-2u(x_k,t;\omega)+u(x_{k-1},t;\omega)}{\Delta x^2}
\prod_{j=1}^n \delta\left(a(x_j)-u(x_j,t;\omega)\right)\right>.
\end{align}
This yields,
\begin{align}
\frac{\tilde{\partial}^2 a(x_k)} {\tilde{\partial} x^2}
\left<\prod_{k=1}^n \delta\left(a(x_k)-u(x_k,t;\omega)\right)\right>
=
\left<\frac{\tilde{\partial}^2 u(x_k,t;\omega)}{\tilde{\partial} x^2}
\prod_{j=1}^n \delta\left(a(x_j)-u(x_j,t;\omega)\right)\right>,
\label{secondderivative}
\end{align}
where $\tilde{\partial}^2/\tilde{\partial} x^2$ is 
the numerical differentiation operator, i.e., the 
approximation of the second-order derivative 
operator by using finite differences, or other 
differentiation schemes such as pseudospectral 
collocation  \cite{Hesthaven}. By extending these arguments 
to higher-order derivatives, we obtain} 
\begin{equation}
 \frac{\tilde{\partial}^s a(x_k)}{\tilde{\partial} x^s}
\left<\prod_{k=1}^n \delta\left(a(x_k)-u(x_k,t;\omega)\right)\right>
=
\left<\frac{\tilde{\partial}^s u(x_k,t;\omega)}{\tilde{\partial} x^s}
\prod_{j=1}^n \delta\left(a(x_j)-u(x_j,t;\omega)\right)\right>, \quad 
s=1, 2, ...
\end{equation}
\noindent
The joint probability density function of the 
state vector $u(x_i,t;\omega)$ ($i=1,...,n$)
satisfies the equation \cite{Venturi_PRS,Venturi_JCP,Venturi_JCP2013}
\begin{align}
\frac{\partial p}{\partial t}=&-\sum_{k=1}^n 
\frac{\partial }{\partial a_k}\left<\frac{\tilde{\partial}^2 u(x_i,t;\omega)}
{\tilde{\partial} x^2}
\prod_{j=1}^n \delta\left(a_j-u(x_j,t;\omega)\right) \right>, 
\nonumber\\
=&-\sum_{k=1}^n 
\frac{\partial }{\partial a_k}
\left(\frac{\tilde{\partial}^2 a(x_i)}{\tilde{\partial} x^2}p\right),
\label{jointPDF}
\end{align}
where $p=p(a_1,...,a_n,t)$. The last equality follows 
from  \eqref{secondderivative}. By taking the continuum limit, 
i.e., by sending $\Delta x$ to zero (and correspondingly $n$ 
to infinity), we obtain the following functional differential equation
for the probability density functional of the solution 
to equation \eqref{1dheat}
\begin{align}
\frac{\partial P([a(x)],t)}{\partial t}=-\int_{-\infty}^\infty
\frac{\delta}{\delta a(x)}\left(\frac{\partial^2 a(x)}{\partial x^2}
P([a(x)],t)\right)dx.
\label{PFDE0}
\end{align}
This equation is in agreement with \eqref{PrDFeq3}.

\vs
\noindent
{\em Remark:} The functional equation \eqref{PFDE0} 
is linear in $P([a(x)],t)$, but it involves a singular term.
Such term is generated by the derivative 
\begin{align}
\frac{\delta}{\delta a(x)}
\left(\frac{\partial^2 a(x)}{\partial x^2}\right)=&
\frac{\delta}{\delta a(x)}\int_{-\infty}^\infty a(y)
\delta''(x-y)dy,\nonumber\\
=&\int_{-\infty}^\infty \delta(x-y)\delta''(x-y)dy.
\label{singular}
\end{align}
The last integral is equivalent to the second derivative of
the Dirac delta function evaluated at zero $x=0$. 
Such singularity can also be seen from a purely discrete 
viewpoint. To this end, substitute the (second-order) 
finite-difference approximation to the second-order derivative into
\eqref{jointPDF}. This yields the equation
\begin{equation}
 \frac{\partial p}{\partial t}=-\sum_{k=1}^n 
\frac{a(x_{k+1})-2a(x_k)+a(x_{k-1)}}{\Delta x^2}\frac{\partial p}{\partial a_k}-
\frac{2}{\Delta x^2}p.
\end{equation}
It is clear that as $\Delta x$ goes to zero (continuum limit), 
the term $2p/\Delta x^2$ generates a singularity. 

\paragraph{Regularity of the Probability Functional} 
The solution to a probability functional equation 
may be an irregular functional. To understand 
why, consider Figure \ref{fig:ProbDF}. When we send $\Delta x$ to 
zero we have that $x_k$ approaches $x_{k+1}$.
Correspondingly the random variables $u(x_k;t,\omega)$ 
and $u(x_{k+1};t,\omega)$ tend to be the same random variable.  
In this situation, the joint PDF of $u(x_k;t,\omega)$ and 
$u(x_{k+1};t,\omega)$ involves a Dirac delta function 
as $x_k\rightarrow x_{k+1}$ (see \cite{Lundgren}, p. 970). 
In a continuum setting, the phenomenon we just 
described happens at each point $x$.
Therefore the probability density functional can 
be an irregular mathematical object.

\vs
\noindent
{\em Example 1:} 
The probability density functional
of a zero mean Gaussian random function $u(x;\omega)$ 
($x\in \mathbb{R}$) with covariance $C(x,y)$ is 
{proportional} to
\begin{equation}
P([a(x)])\sim  \exp\left[-\frac{1}{2}\int_{-\infty}^\infty 
\int_{-\infty}^\infty 
C^{-1}(x,y)a(x)a(y)dxdy\right),
 \label{PDFunct_for_Gaussian}
\end{equation}
where $C^{-1}(x,y)$ is the inverse covariance 
function. 
Such inverse covariance may be obtained by  
solving the Fredholm integral equation 
of the first kind 
\begin{equation}
 \int_{-\infty}^{\infty } C(x,y) C^{-1}(y,z)dy=\delta(x-z).
 \label{diffinv}
\end{equation} 
If $C(x,y)$ is smooth then its differential 
inverse $C^{-1}(y,z)$ must have serious 
singularities in order for the integral in \eqref{diffinv} 
to yield a Dirac delta function (see Table \ref{tab:4}). 
If $C(x,y)$ is homogeneous, i.e., if $C(x,y)=C(x-y)$, 
then $C^{-1}(y,z)$ is called {\em convolution inverse} 
\cite{Hohlfeld,Murthy}. 
This method was pioneered by Hirschman and 
Widder \cite{Hirschman} in the late forties.
Relevant cases of convolution inverses are summarized 
in Table \ref{tab:4}. Note that the convolution 
inverses of smooth coavariance functions -- 
such as the Mat\'e rn covariance in 2D -- 
are  {\em rough} functions involving  
Laplacians and bi-harmonic operators 
applied to Dirac delta functions. 
\begin{table}
\centering 
\begin{tabular}{lcccc}
& Covariance  & & Inverse Covariance  \vspace{0.2cm}\\
\hline\\
\begin{minipage}{3cm}
 \centering Exponential (1D)\\\vspace{0.1cm }
\end{minipage}\hspace{1cm}
 &  $\displaystyle \sigma^2 e^{-|x|/h}$  & &
$\displaystyle \frac{1}{2h \sigma^2}\left(\delta(x)-h\delta''(x)\right)$ \\\\
\begin{minipage}{3cm}
 \centering Mat\'ern (2D)\\\vspace{0.1cm }\footnotesize
 (polar coordinates)
\end{minipage}\hspace{1cm}
 &  $\displaystyle \sigma^2 r B_1\left(\frac{r}{h}\right)$  & &
$\displaystyle \frac{1}{4\pi \sigma^2h^2}\left(\frac{\delta(r)}{\pi r}-
 2h^2\nabla^2\frac{\delta(r)}{\pi r}+h^4 \nabla^4\frac{\delta(r)}{\pi r}\right)$ \\\\
\begin{minipage}{3cm}
\centering Exponential (3D)\\ \vspace{0.1cm }\footnotesize
 (polar coordinates)
\end{minipage}
& $\displaystyle \sigma^2 e^{-r/h}$ &\hspace{0.5cm} & 
$\displaystyle \displaystyle \frac{1}{8\pi \sigma^2h^3}\left(\frac{\delta(r)}{\pi r}-
 2h^2\nabla^2\frac{\delta(r)}{\pi r}+h^4 \nabla^4\frac{\delta(r)}{\pi r}\right)$\\
 \end{tabular}
\caption{Convolution inverses of well-known covariance functions \cite{Oliver}. Here 
$B_1$ denotes the modified Bessel function of the first kind. 
It is seen that the convolution inverse of smooth functions are rough 
functions involving Laplacians and bi-harmonic operators 
applied to Dirac delta functions.}
\label{tab:4}
\end{table}
However, such rough functions appear 
within integrals in \eqref{PDFunct_for_Gaussian}, and therefore  
we expect some regularization. For instance, if 
we assume that the covariance function $C(x-y)$ is
exponential (see Table \ref{tab:4}), then from 
\eqref{PDFunct_for_Gaussian} we obtain  
\begin{equation}
P([a(x)])\sim  \exp\left[
-\frac{1}{4h\sigma^2}
\int_{-\infty}^\infty \left(
a(x)^2 + h^2\left[\frac{da(x)}{dx}\right]^2\right)dx
\right].
\label{PDFforGexp}
\end{equation}

Probability density functional equations represent  
an excellent starting point to obtain effective 
approximations. To this end, one needs 
to follow (by analogy) the route taken in 
classical statistical mechanics in which we start with 
the Liouville equation, and make approximations 
in order to derive an computable equation for a quantity of 
interest. Such {\em coarse-graining} process for 
functional differential equations is discussed 
by McComb \cite{McComb} in the context of fluid 
turbulence. Probability density functional 
equations were derived and studied in the context turbulent 
flows by Dopazo and O'Brien \cite{Dopazo},  and
Rosen \cite{Rosen_1960,Rosen_1969,Rosen_1967}.

\vs
\noindent
{\em Hopf Functionals and Probability Density Functionals:}
We have seen that Hopf equations and probability density 
functional equations are related by a functional Fourier transform. 
Therefore, from a purely mathematical viewpoint they are 
completely equivalent. However, from the viewpoint 
of approximation theory they are not equivalent at all.
Hopf functionals may be hard to resolve due to 
high-frequencies related to the complex exponential.
On the other hand, probability density functional equations 
may have non-smooth solutions.
The statistical properties of a random field can be equivalently 
computed by using the Hopf functional or the probability 
density functional. In the first case, we simply need to 
take functional derivatives and evaluate them at $\theta(x)=0$ 
(see Section \ref{app:functional derivatives}). 
In the second case, we need to compute functional 
integrals, i.e., integrals in an infinite number of variables. 
This requires requires a careful definition of the integration 
measure (\cite{Beran}, \S 2.2.4).

{\color{r}
\subsection{Effective Fokker-Planck Systems}
Consider the stochastic dynamical system \eqref{eqofm}. Suppose 
we are interested in determining an evolution equation for the joint 
probability density function of the state vector $\bm \psi(t)$.
To this end, we think of $\bm \psi(t)$ as a nonlinear 
functional of the random noise $\bm f(t)$, i.e., we can consider 
the map $\bm \psi(t)=\bm \Psi(t;[\bm f(t)])$.
The specific form of $\bm \Psi$ depends on the system, 
in particular on the nonlinear map $\bm \Lambda(\bm \psi,t)$ in 
\eqref{eqofm}.
The probability density function of the $\bm \Psi(t;[\bm f(t)])$ 
can be expressed as a functional integral over the probability 
density functional of the noise (assuming it exists)
\begin{equation}
p(\bm \psi,t) = \int \delta(\bm \psi -\bm \Psi(t;[\bm f])) P([\bm f])\mathcal{D}[\bm f], 
\label{PRS}
\end{equation}
where $\delta$ here is a multivariate Dirac delta function. 
From this expression, it is clear that the random noise $\bm f(t)$ 
determines $p(\bm \psi,t)$, and therefore the structure of 
probability density function equation that evolves 
$p(\bm \psi,t)$ in time. 
Such equation can be derived by using functional calculus.
\cite{Venturi_PRS,Fox,Hanggi4,Moss1}, and it takes the form
\begin{equation}
\frac{\partial p(\bm \psi,t)}{\partial t} + 
\nabla_{\bm \psi}\cdot\left(\bm \Lambda (\bm \psi ,t) 
p(\bm \psi,t) + \left<\delta\left(\bm \psi-
\bm \Psi(t;[\bm f])\right)\bm f(t)\right>_{\bm f}\right)=0,
\label{PDF}
\end{equation}
where
\begin{equation}
 \left<\delta\left(\bm \psi-
\bm \Psi(t;[\bm f])\right)\bm f(t)\right>_{\bm f}=\int 
\delta(\bm \psi -\bm \Psi(t;[\bm f])) \bm f(t)P([\bm f])\mathcal{D}[\bm f].
\label{corr}
\end{equation}
The quantity \eqref{corr} represents the {\em correlation} 
between two functionals of the random noise, namely  
$\delta(\bm \psi -\bm \Psi(t;[\bm f]))$ and the noise itself $\bm f(t)$. 
Such correlation can be disentangled 
by using functional integral techniques 
(see, e.g., \cite{Bochkov,Venturi_PRS,Klyatskin}). 
In particular, if $\bm f(t)$ is Gaussian then \eqref{corr}  
can be expressed by the well-known Furutsu-Novikov-Donsker 
formula \cite{Venturi_PRS,Furutsu,Novikov,Donsker}. 
Similarly, if $\bm f(t)$ is Gaussian white noise 
then \eqref{eqofm} is a Markovian system and 
the correlation \eqref{corr} reduces to a simple 
diffusion term \cite{Venturi_PRS}. In this 
case \eqref{PDF} coincides with the classical 
Fokker-Plank equation \cite{Risken}. We remark that 
computing the solution to \eqref{PDF} in the general (non-Markovian) setting 
is very challenging. Possible techniques rely on {\em data-driven} models 
that employ random paths of the SODE \eqref{eqofm}, 
e.g., path integral methods \cite{Hanggi3,Pesquera,McKane}.

\paragraph{Non-Markovian Random Processes on Random Graphs}

Consider a {\em non-Markovian} stochastic process 
$\bm \psi(t;\omega)\in\mathbb{R}^n$ evolving 
on a {\em random graph} 
with $n$ nodes (see Figure \ref{fig:graph}). 
Such process could model, e.g., the propagation of 
epidemics in interacing individuals \cite{Andersson,Daley}. 
Introducing uncertainty in the graph allows us to take into 
account uncertainty in the interconnections between different 
nodes, which is fundamentally important when modeling 
dynamics of {\em social networks} and disease propagation. 
\begin{figure}
\centerline{
\includegraphics[height=4cm]{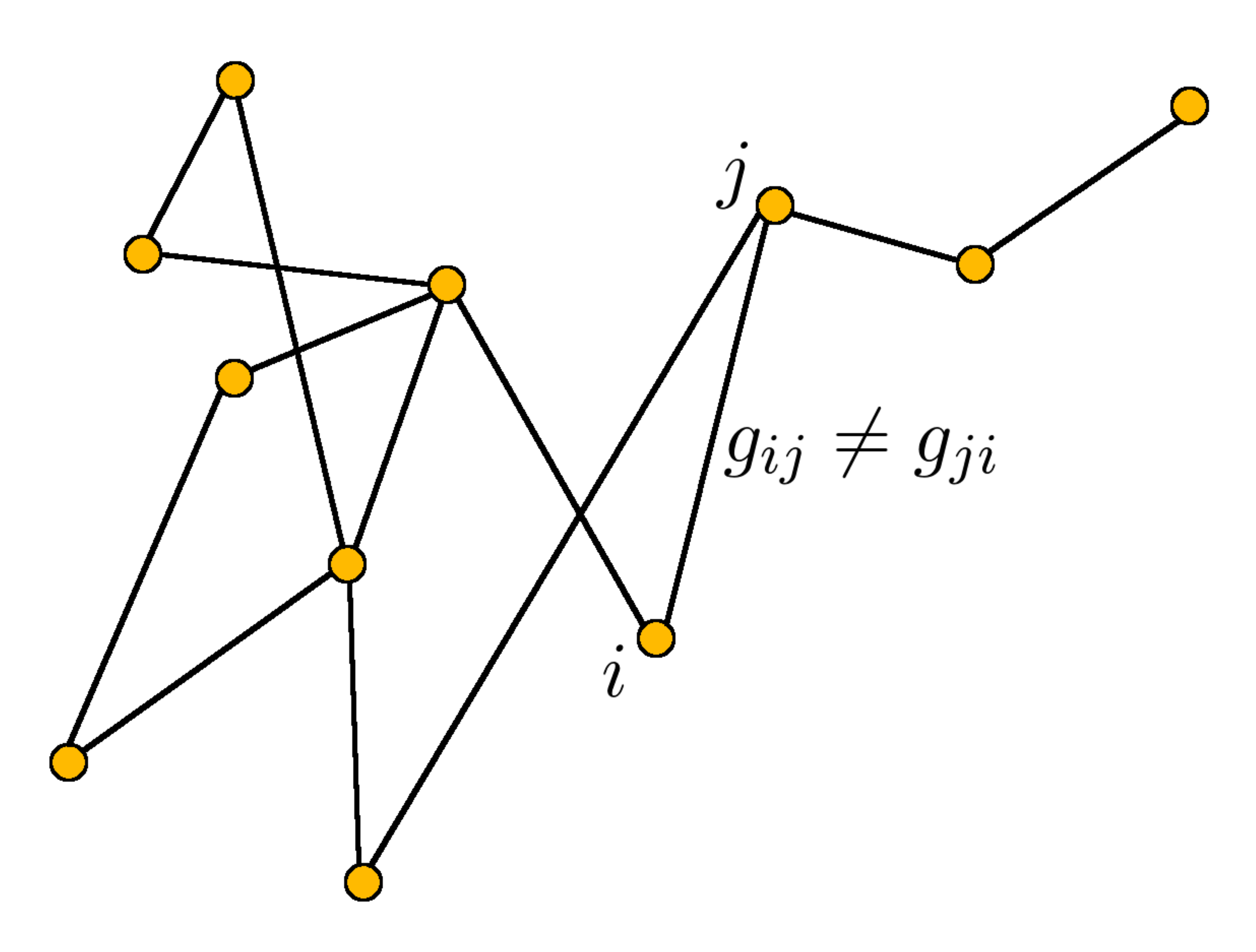}
\includegraphics[height=4cm]{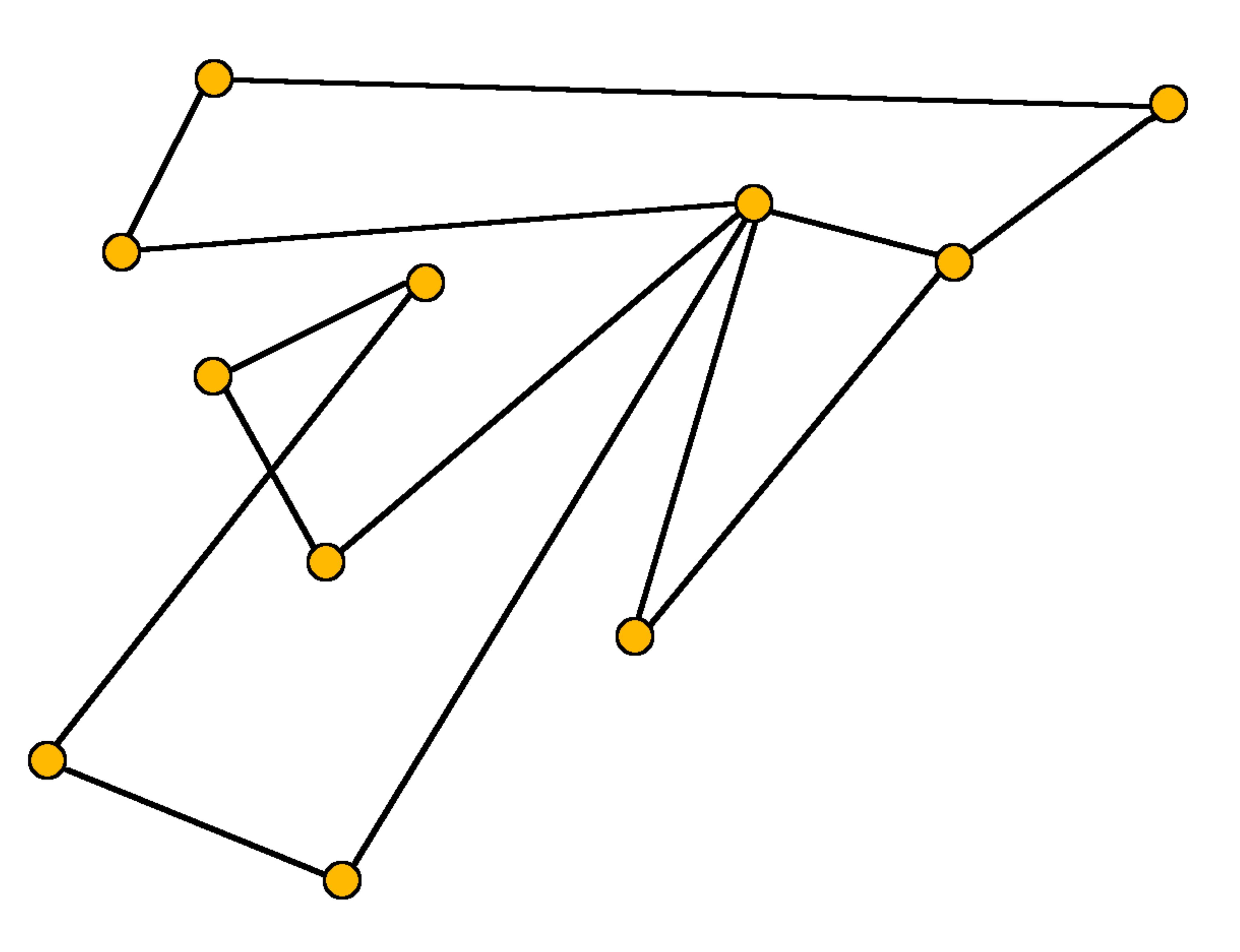}
\includegraphics[height=4cm]{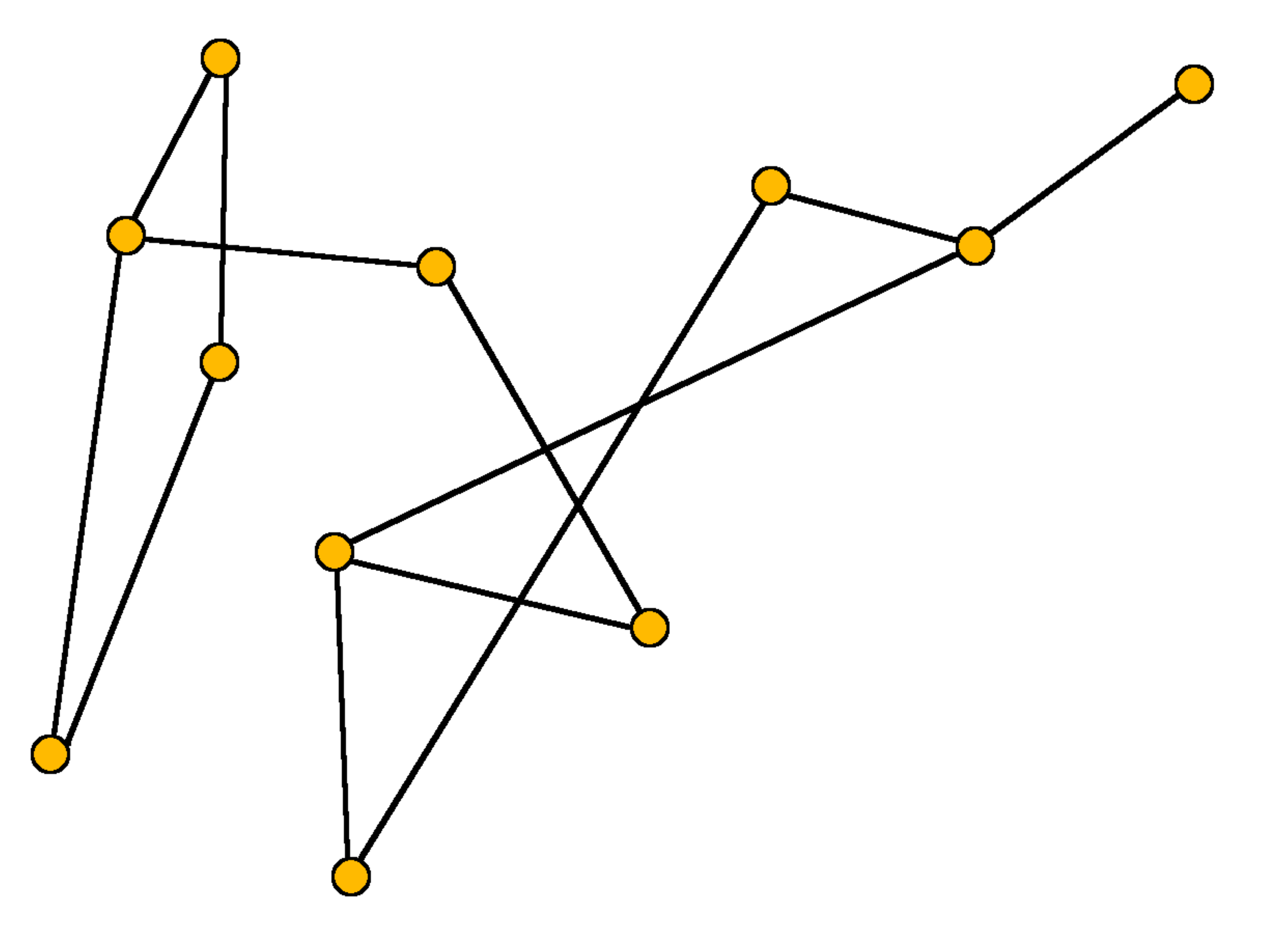}
}
\caption{Samples of a random graph with fixed 
number of nodes. Introducing uncertainty in the edges that connect 
different nodes allows us to take into account uncertainty 
in the interactions between different components of a
stochastic model evolving on the graph. This is important 
when modeling dynamics of social networks and disease propagation.}
\label{fig:graph}
\end{figure}
We can characterize a random graph mathematically in terms 
of random edges defined by {\em mixed} random variables, 
i.e., random variables with continuous/discrete probability distribution. 
To this end, let $g_{ij}(\omega)$ represent the interaction between 
the node $i$ and the node $j$, i.e., the flow of information from $i$ 
to $j$\footnote{For example, the probability density of $g_{ij}$ 
could be in the form $p(g_{ij})=\delta(0)/2 + U_{[1,2]}(g_{ij})/2$, in which case 
we have $50\%$ of chances that $i$ does not influence $j$ at all 
(term $\delta(0)/2$), and $50\%$ of chance that $g_{ij}$ is uniformly 
distributed in $[1,2]$. If the graph has a time-evolving random 
structure, i.e., some nodes in the network have a time-varying set 
of neighbors, then we can introduce time-dependence in the random 
variables $g_{ij}$. This makes $g_{ij}$ a set of stochastic processes.}. 
Such interaction can be of different types, but roughly speaking 
it just characterizes how the random 
process $\psi_i(t)$ (at node $i$) influences 
the random process $\psi_j(t)$ (at node $j$). 
The type of influence is defined by a stochastic model 
on the random graph, e.g., a system of 
stochastic differential equations in the form
\begin{equation}
 \frac{d\bm \psi }{dt}= \bm \Lambda(\bm \psi,t;{\bm G}(\omega)) + \bm f(t),
\label{Eq1}
 \end{equation}
where $\Lambda$ is a nonlinear map defining the model, 
${\bm G}(\omega)=\{g_{12}(\omega),g_{21}(\omega), ...\}\in\mathbb{R}^{n(n-1)/2}$, 
is the random vector representing all interactions among the nodes 
in the random graph, and $\bm f(t)$ is colored random noise. 
The stochastic system \eqref{Eq1} can be obviously generalized to cases 
where we have multiplicative noises, 
such as in tumororal cell growth models \cite{Zeng,Wang,Fiasconaro}.
The solution to \eqref{Eq1} (assuming it exists) is a nonlinear function  
of the random vector ${\bm G}$ defining the graph, 
and a nonlinear functional of the stochastic 
process $\bm f(t)$. We write such functional 
as $\widehat{\bm \psi}(t;{\bm G},[\bm f(t)])$.
The probability density function of $\bm \psi(t)$ 
then can be obtained by integrating out the graph 
and the noise over the corresponding probability distribution, 
i.e.,
\begin{equation}
 p(\bm \psi,t)= \int \delta\left(\bm \psi-\widehat{\bm \psi}(t;{\bm G},[\bm f])\right)
 p({\bm G})P([\bm f])\mathcal{D}[\bm f] d{\bm G}.
\end{equation}
Note that here we assumed that the noise and the graph are 
statistically independent. Integration over noise is 
(in general) a functional integral. The exact evolution 
equation for $p(\bm X,t)$ can be obtained by using a 
functional calculus approach \cite{Venturi_PRS,Fox,Hanggi4,Moss1}. This yields,
\begin{equation}
\frac{\partial p(\bm \psi,t)}{\partial t} + 
\int \left[\nabla_{\bm \psi}\cdot\left(\bm \Lambda (\bm \psi ,t;\bm G) 
p(\bm \psi,\bm {G},t) \right)+ \nabla_{\bm \psi}\cdot \left<\delta\left(\bm \psi-
\widehat{\bm \psi}(t;\bm {G},[\bm f])\right)\bm f(t)\right>_{\bm f}\right] 
d\bm G=0,
\label{PDF_ha}
\end{equation}
where $p(\bm \psi,\bm G,t)$ is the joint probability 
density of $\bm \psi(t,\omega)$ and $\bm G(\omega)$, 
while $\left<\cdot\right>_{\bm f}$ is defined in 
\eqref{corr}. As before, the correlation  $\left<\delta\left(\bm \psi-
\widehat{\bm \psi}(t;\bm {G},[\bm f])\right)\bm f(t)\right>_{\bm f}$
can be disentangled by using
functional calculus \cite{Bochkov,Venturi_PRS,Klyatskin}, or 
computed by using data-driven methods.
}

\begin{table}
{\small
\begin{center}\begin{tabular}{lcl}
% \hline
Field Equation & & Functional Differential Equation\\
\hline  \\
\begin{minipage}{4cm}
\begin{equation}
\frac{\partial \bm u}{\partial t}+(\bm u\cdot \nabla)\bm u=-\nabla p +\nu \Delta \bm u\nonumber
\end{equation}
\end{minipage}
 & &
\begin{minipage}{4cm}
 \begin{align}
\frac{\partial \Phi([\theta],t)}{\partial t}=
\sum_{k=1}^3\int_V\theta_k(\bm x)\left(i \sum_{j=1}^3\frac{\partial }{\partial x_j}
\frac{\delta^2 \Phi([\bm \theta],t)}{\delta \theta_k(\bm x)\delta\theta_j(\bm x)}
+\nu \nabla^2\frac{\delta \Phi([\bm \theta],t)}{\delta \theta_k(\bm x)}\right)d\bm x
\nonumber
 \end{align}
\end{minipage}
\\\\\\
$\displaystyle \frac{\partial u}{\partial t}+
u\frac{\partial u }{\partial x}=\nu \frac{\partial^2u}{\partial x^2}$ 
&  & 
$\displaystyle\frac{\partial P([a],t)}{\partial t} = -\int_{-\infty}^{\infty} 
\frac{\delta}{\delta a(x)}\left(\left[-a(x)\frac{\partial a(x)}{\partial x}
+\nu\frac{\partial ^2 a(x)}{\partial x^2}\right] P([a],t)\right)dx$
\\\\
$\displaystyle \square \phi+m^2 \phi = 
\frac{\lambda}{6}\phi^3$ &\hspace{1.0cm} & 
$\displaystyle
 \square
\frac{\delta Z([a])}{\delta a(\bm x)}+ 
 m^2 \frac{\delta Z([a])}{\delta a(\bm x)}-
\frac{\lambda}{3!} 
\frac{\delta^3 Z([a])}{\delta a(\bm x)^3 }-i a(\bm x) Z([a])=0$

\end{tabular}
\end{center}
}
\caption{Examples of functional differential equations. In this table,  
$\Phi$ denotes the Hopf characteristic functional, 
$P$ is the probability density functional, and $Z$ is the generating functional of the quantum $\phi^4$-theory.}
\label{tab:2}
\end{table}

\subsection{Conjugate Flow Action Functionals}
In a recent paper \cite{Daniele_JMathPhys}, we have shown how to 
construct an action functional for a non-potential 
field theory by using methods of differential geometry and nonlinear 
functional analysis \cite{Vainberg,Nashed}. The key idea is to 
represent the governing equations of the field theory 
relative to a diffeomorphic flow of curvilinear coordinates 
which is assumed to be functionally dependent on the 
field equations, i.e., on their solution. 
Such flow evolves in space and time similarly to a 
physical fluid flow of classical mechanics 
and it can be chosen to symmetrize the G\^ateaux derivative 
of the field equations relative to suitable local bilinear forms. 
This is equivalent to require that the governing equations 
of the field theory can be derived from a principle of 
stationary action on a flow, which we called 
the {\em conjugate flow of the theory}. The determining 
equations of the conjugate flow are functional 
differential equations. 
In particular, for a second-order nonlinear scalar field theory
\begin{equation}
 f\left(u; u_{,\mu}; u,_{\mu\nu}; \widehat{x}^\mu_{,\nu}; \widehat{x}^\mu_{,\nu\lambda} \right)=0,
 \label{fieldeq}
\end{equation}
we obtain
\begin{equation}
 R^{\mu\nu}_{,\nu}+R^{\mu\nu}\Gamma^\lambda_{\lambda\nu}=Z^\nu,\label{cfaf}
\end{equation}
where the comma denotes differentiation with respect to the 
independent variable $\sigma^\nu$, 
$\Gamma^\lambda_{\lambda\nu}$ is the Christoffel symbol of the 
second kind and 
\begin{align}
Z^\nu&=\frac{\partial f}{\partial u_{,\nu}}+
\frac{\partial f}{\partial \widehat{x}^\mu_{,\nu}}\frac{\delta \widehat{x}^\mu}{\delta u}+
\left(\frac{\partial f}{\partial \widehat{x}^\mu_{,\nu\rho}}+
\frac{\partial f}{\partial \widehat{x}^\mu_{,\rho\nu}}\right)
\left(\frac{\partial}{ \partial \sigma^\rho}\frac{\delta \widehat{x}^\mu}{\delta u}+
\frac{\delta^2\widehat{x}^\mu}{\delta u^2}\frac{\partial u}{\partial \sigma^\rho}\right),\label{Bs}\\
R^{\rho\nu}&=\frac{\partial f}{\partial u_{,\nu\rho}}+
\frac{\partial f}{\partial \widehat{x}^\mu_{,\rho\nu}}\frac{\delta \widehat{x}^\mu}{\delta u}. \label{Fs}
\end{align}
Given a solution to the field equation \eqref{fieldeq}, the 
system of functional differential equations \eqref{cfaf} 
allows us to identify the {\em conjugate flow} $\widehat{x}^\mu(\sigma^\nu;[u])$, i.e., the functional 
relation between the flow $\widehat{x}^\mu$ and 
the solution  $u$ for which the PDE \eqref{fieldeq} 
can be derived from a {\em principle of least action} 
(see \cite{Daniele_JMathPhys} for further details). 
The identification of transformation groups leaving the 
conjugate flow action functional invariant could lead to 
the discovery of new conservation laws.

\subsection{Large Deviation Theory and Minimum Action Methods}
Large deviations theory deals with the probabilities of rare 
events that are exponentially small as a function of some 
parameter. To illustrate the theory,  
consider a nonlinear PDE perturbed by space-time additive 
random noise of small amplitude $\epsilon$
\begin{equation}
 \frac{\partial \bm u}{\partial t}+\bm G(\bm u)=\sqrt{\epsilon}\bm f(\bm x,t;\omega).
\end{equation}
We define the set trajectories connecting two arbitrary states $\{\bm u_1,\bm u_2\}$
\begin{equation}
B=\{\bm u(\bm x,t;\omega) | \quad \bm u(\bm x,0;\omega)=\bm u_1, 
\quad \bm u(\bm x,T;\omega)=\bm u_2\}
\end{equation}
If $\bm f$ is white noise, then the Freidlin-Wentzell 
theory gives us the following large deviation principle
\begin{equation}
\lim_{\epsilon \rightarrow 0^+} \epsilon \textrm{Pr}(\bm u\in B)=\inf_{\bm u\in B}
\frac{1}{2}\int_0^T \left\|\frac{\partial \bm u}{\partial t}+
\bm G(\bm u)\right\|^2_{L_2}dt
\end{equation}
where $\textrm{Pr}(A)$ denotes the probability of the event $A$, 
$\left\|\cdot\right\|^2_{L_2}$ indicates the $L_2$ norm in space. 
The large deviation principle is equivalent to minimum action 
principle
\begin{equation}
\min_{\substack{\bm u(0,\bm x;\omega)=\bm u_1\\\bm u(T,\bm x;\omega)=\bm u_2}} 
S_T([\bm u]),\qquad S_T([\bm u])=\frac{1}{2}\int_0^T 
\left\|\frac{\partial \bm u}{\partial t}+\bm G(\bm u)\right\|^2_{L_2}dt.\label{DV}
\end{equation}
The minimizer of \eqref{DV} is called {\em minimum action path}, 
and it satisfies the functional differential equation
\begin{equation}
\frac{\delta S_T([\bm u])}{\delta \bm u(\bm x,t)}=0 \quad \bm u\in B.
\label{DV1}
\end{equation}
The minimum action path is the most probable transition 
path from $\bm u_1$ to $\bm u_2$. A method to 
solve \eqref{DV1} is based on a direct discretization 
of $S_T(\bm u)$ \cite{E1}.

\section{Approximation of Functional Differential Equations}
\label{sec:FDEapprox}
In this Section we address the numerical 
approximation of linear FDEs in the 
form \eqref{linfde00}. To this end, 
we develop a method of weighted 
residuals \cite{Finlayson,GKSS_2005} 
in the space of functionals that allows us to 
derive functional least squares, functional 
Galerkin and functional collocation methods 
to FDEs is a unified and straightforward way.

\subsection{The Method of Weighted Functional Residuals}
\label{sec:method of weighted residuals}
The method of weighted residuals illustrates how 
the choice of different weight (or test) functionals 
can be used to construct different classes of methods 
extending Galerkin, least-squares and collocation 
methods for PDEs to functional differential equations.
The general framework resembles the classical one for PDEs,
in which one minimizes a residual (least-squares method) 
or imposes its orthogonality relative to a suitable space of 
test functions (Galerkin or collocation methods). 
To describe the weighted residuals technique, let us consider 
the linear functional differential equation \eqref{linfde00}.
In approximating its solution numerically we are typically 
replacing $F([\theta],t)$ with an approximation 
\begin{equation}
 \widehat F ([\theta],t) \simeq F([\theta],t),
 \label{FDE_exp}
\end{equation}
e.g., a tensor canonical tensor decomposition 
(Section \ref{sec:CP}) or a Lagrangian 
interpolant (Section \ref{sec:Functional Collocation Methods}) 
with $N$ degrees of freedom. Substitution of 
the approximation \eqref{FDE_exp} 
into equation \eqref{linfde00} 
yields the (functional) residual 
\begin{equation}
R([\theta],t)= \frac{\partial  \widehat F ([\theta],t)}{\partial t}-L([\theta],t) 
\widehat F([\theta],t)-H([\theta],t).
\label{rres}
\end{equation}
At this point, we introduce the following inner 
product in the space of functionals\footnote{
The inner product \eqref{rres} is a functional integral, 
which is usually defined in terms of a limiting procedure \cite{Popov,Egorov}. From a mathematical viewpoint, 
the limiting procedure defining the functional integral 
measure in terms of an infinite products of elementary 
measures should be handled with care. 
In fact,  the classical Lebesgue measure does not exist 
in spaces ofinfinite dimension \cite{Marzucchi}. 
On the other hand, Gaussian measures 
are still well defined in such setting. This is why we included
$W([\theta])$ in \eqref{continuous_ip_W0}. The argument leading 
to the result on non-existence of an analogue to the Lebesgue 
measure in infinite dimension is related to the argument showing 
that the Heine-Borel theorem does not hold in infinite-dimensional 
normed linear spaces.} (see Appendix \ref{sec:functional integrals}) 
\begin{equation}
\left(F,G\right)_W=\int 
F([\theta])G([\theta])W([\theta])\mathcal{D}[\theta],\qquad \textrm{(functional integral)}
\label{continuous_ip_W0}
\end{equation}
where $W([\theta])$ is a known weight functional, and consider the set of equations
\begin{equation}
\left(R([\theta],t),h_k([\theta])\right)_W=0\qquad k=1,...,N
\label{WR}
\end{equation}
where $h_k([\theta])$ are test functionals. 
There is no particular restriction on $h_k([\theta])$. 
For example, they can be 
cardinal basis functionals, orthogonal polynomial functionals 
or other basis functionals.
The system \eqref{WR} allows us to determine the $N$ 
degrees of freedom in the functional approximation 
$\widehat{F}([\theta],t)$. Specifically, 
we are imposing that the residual of the FDE is 
orthogonal to the span of the functionals 
$\{h_1,...,h_N\}$. The nature of the numerical 
scheme is determined by the choice of the test 
functionals $h_j([\theta])$ in \eqref{WR}. 

Evaluating the functional integrals in \eqref{WR} 
is challenging, but there are approximation methods 
that allow us to compute them. 
For example, several algorithms have been recently 
proposed for high-dimensional (possibly infinite-dimensional) 
integration \cite{Baldeaux,Wasilkowski,Wasilkowski1,Dick,Dick1}
 (see also Apendix \ref{sec:functional integrals} and 
 Chapter 4 in \cite{Egorov}).
The system of equations \eqref{WR} defines 
a {\em functional Galerkin method}\footnote{
Stochastic Galerkin methods \cite{Db_book} 
are functional Galerkin methods. 
Essentially, these approaches are based 
on stochastic representations of the solution functional (Section 
\ref{sec:Stochastic Functional Methods}), and functional 
inner products involving probability measures. 
Stochastic Galerkin methods have been studied extensively 
in the theory of turbulence  \cite{Meecham,Meecham1,Lee,Bodner}, 
and in uncertainty quantification \cite{db1,Xiao1}.
}.

\subsubsection{Functional Collocation Methods}
\label{sec:collocation}
In this class of methods the test 
functionals $h_j([\theta])$  are 
chosen to be Dirac delta functionals \cite{Jouvet} centered 
at $\theta_j(x)$, i.e., $h_j([\theta])=\delta[\theta_j(x)-\theta(x)]$. 
In this setting, the orthogonality 
condition \eqref{WR} can be 
written as 
\begin{equation}
 \int R([\theta],t) \delta[\theta_j(x)-\theta(x)]W([\theta]) \mathcal{D}[\theta]= 0\quad \Rightarrow \quad R([\theta_j],t)=0
\end{equation}
In other words, in the functional collocation method we 
impose that the residual $R([\theta],t)$ vanishes 
at $N$ collocation nodes in $D(F)$, 
i.e., $N$ functions $\{\theta_1(x),...,\theta_N(x)\}$. 
This yields a system of $N$ equations for the 
unknowns $\{\alpha_1(t), ..., \alpha_N(t)\}$. 
The solution we obtain from the functional collocation 
method obviously interpolates the exact solution at the 
nodes $\theta_i(x)$.

\subsubsection{Functional Least Squares}
\label{sec:functional_least_squares}
In this class of methods we look for an
approximate solution functional that minimizes 
the norm of the residual $R([\theta],t)$. 
Such norm may be defined in terms of the functional inner 
product \eqref{continuous_ip_W0}, i.e., 
$\left\|R\right\|_W^2=(R,R)_W$. 
In this case, we obtain the following variational 
principle involving a functional integral 
\begin{equation}
\min_{\widehat{F}\in D_N(F)}\left\|R([\theta],t)\right\|_W^2=
\min_{\widehat{F}\in D_N(F)} \int 
R([\theta],t)^2 W([\theta])\mathcal{D}[\theta].
\label{FI}
\end{equation}
The stationary points of \eqref{FI} corresponding to 
variations of the degrees of freedom $\alpha_k(t)$ 
in \eqref{FDE_exp} satisfy the the Euler-Lagrange 
equations 
\begin{equation}
\left(R([\theta],t),\frac{\partial R([\theta],t)}{\partial \alpha_k(t)}\right)_W=0,\qquad k=1,...,N.
\label{MinR}
\end{equation}
A comparison between \eqref{WR} and \eqref{MinR} suggests that the test 
functionals $h_j([\theta])$ in this case are equation-dependent, i.e., 
they depend on the residual $R([\theta],t)$ through the formula
\begin{equation}
 h_k([\theta],t) = \frac{\partial R([\theta],t)}{\partial \alpha_k}.
\end{equation}

\vs
\noindent
{\em Remark:} 
Error analysis, stability and consistency of 
functional Galerkin, functional collocation and functional 
least squares methods is an {open question}.

\vs
\noindent
{\em Remark:}  If we restrict $D(F)$ to a finite-dimensional function space, e.g, the span of a finite-dimensional basis, then weighted residual formulation we just discussed reduces to the weighted residual formulation for multivariate linear PDEs.

\subsection{Temporal Discretization}
\label{sec:ADI_SSE}

The FDE \eqref{linfde00} can be discretized in time with 
different numerical schemes such as Adams-Bashforth,  
Adams-Multon or BDF methods \cite{Quarteroni}. 
Such discretization is quite classical in numerical 
analysis, and it represents an important building 
block in the development of efficient algorithms 
to compute the numerical solution to FDEs. 
Hereafter we discuss functional tensor 
methods built upon explicit and implicit linear multistep  
schemes (see \cite{Quarteroni}, p. 497). 
As we will see, such algorithms have significant 
advantages over other approaches in terms of accuracy 
and computational cost.  
Given an evenly-spaced sequence of time 
instants $t_k= k\Delta t$ ($k=0,1,...$)  we write the formal 
solution to to the FDE \eqref{linfde00} as 
\begin{equation}
 F([\theta],t_{n})=F([\theta],t_{n-1})+\int_{t_{n-1}}^{t_n} 
 \left(L([\theta],\tau)F([\theta],\tau)+H([\theta],\tau)\right)d\tau. 
\label{FSint}
\end{equation}
By approximating the temporal integral with a quadrature rule 
we obtain a fully discrete time-integration scheme. 
For example, if we replace
\begin{equation} 
S([\theta],\tau)=L([\theta],\tau)F([\theta],\tau)+H([\theta],\tau)
\label{SintS}
\end{equation}
by the interpolating polynomial at $t_{n-1}$, ..., $t_{n-q}$, extrapolate in $[t_{n-1},t_n]$, and integrate in time 
we obtain the $q$-th order Adams-Bashforth scheme. 
Hereafter we provide some examples.

{
\color{r}
\subsubsection{Second-order Adams-Bashforth (AB2) method} 
We replace \eqref{SintS} with the polynomial interpolating 
$S([\theta],t)$ at $\{t_{n-1},t_{n-2}\}$, extrapolate such polynomial to to $[t_{n-1},t_n]$ and compute the integral in \eqref{FSint}. 
This yields the second-order explicit explicit scheme
\begin{equation}
 F([\theta],t_{n})=F([\theta],t_{n-1})+\frac{\Delta t}{2}
\left[3S([\theta],t_{n-1})-S([\theta],t_{n-2})\right] +
\Delta t \tau_n([\theta]),
\label{FSint1}
\end{equation}
where 
\begin{equation}
\tau_n([\theta])=\frac{5\Delta t^3}{12}\frac{\partial^3 F([\theta],t_{n-2})}{\partial t^3}+\mathcal{O}(\Delta t^4).
\end{equation} 
The quantity $\tau_n$ is the {\em local truncation error} 
at time $t_n$ (\cite{Quarteroni}, p. 499). 
Clearly, if the operator $L([\theta],t)$ is time-independent and $H=0$, then \eqref{FSint} has the simpler form 
\begin{equation}
 F([\theta],t_{n})=F([\theta],t_{n-1})+\frac{\Delta t}{2}
L([\theta])\left[3F([\theta],t_{n-1})-F([\theta],t_{n-2})\right] +
\Delta t \tau_n([\theta]).
\label{FSintti}
\end{equation}

\subsubsection{Crank-Nicolson Method} 
The Crank Nicolson method is an implicit 
second-order method of the Adams-Multon family. 
The scheme can be easily derived by discretizing  
the time integral in \eqref{FSint} with the trapezoidal 
rule. This yields
\begin{align}
F([\theta],t_n)=F([\theta],t_{n-1}) + \frac{\Delta t}{2}\left[S([\theta],t_{n})+S([\theta],t_{n-1})\right] + \Delta t\tau_n([\theta]),
\label{CN1}
\end{align} 
where 
\begin{equation}
\tau_n([\theta])=
-\frac{\Delta t^3}{12}\frac{\partial^3 F([\theta],t_{n-1})}{\partial t^3}+\mathcal{O}(\Delta t^4)
\end{equation} 
is the local truncation error (\cite{Quarteroni}, p. 499). 
The scheme \eqref{CN1} can be rewritten as 
\begin{align}
\left[I-\frac{\Delta t}{2} L([\theta],t_n)\right]F([\theta],t_n)=&
\left[I+\frac{\Delta t}{2} L([\theta],t_{n-1})\right] F([\theta],t_{n-1})+
\nonumber\\
&\frac{\Delta t }{2}\left[H([\theta],t_{n})+H([\theta],t_{n-1})\right]+
\Delta t\tau_{n}([\theta]),
\label{CN2}
\end{align}
The integration process proceeds as follows: 
Given $F([\theta],t_{n-1})$ we build the right hand side 
of \eqref{CN2} and then solve for $F([\theta],t_n)$. 
This involves inverting the following (functional differential) 
linear operator  
\begin{equation}
A([\theta],t_n)=\left[I-\frac{\Delta t}{2} L([\theta],t_n)\right].
\label{An}
\end{equation}
It is convenient to rewrite \eqref{CN2} as 
\begin{equation}
 A([\theta],t_n)F([\theta],t_n)=E([\theta],t_n)+\Delta t \tau_n([\theta]),
 \label{AFE}
\end{equation}
where 
\begin{equation}
 E([\theta],t_n)=\left[I+\frac{\Delta t}{2} L([\theta],t_{n-1})\right] 
 F([\theta],t_{n-1})+ \frac{\Delta t }{2}\left[H([\theta],t_{n})+
 H([\theta],t_{n-1})\right] .
\end{equation}
is a known functional, provided $F([\theta],t_{n-1})$ is known (solution at time $t_{n-1}=(n-1)\Delta t$). 
}

\subsection{Functional Approximation}
The solution to the FDE \eqref{linfde00} can be approximated at each time step by using the functional approximation methods we discussed 
in Section \ref{sec:rep-Hopf}. For example, if we restrict the 
domain the solution functional $F$ to the finite-dimensional 
Hilbert space spanned by the  orthonormal 
basis $\{\varphi_1,...,\varphi_m\}$, i.e., 
\begin{equation}
D_m=\textrm{span}\{\varphi_1,...,\varphi_m\}\subseteq D(F), 
\label{SoF}
\end{equation}
then the functional becomes a multivariate function. 
Alternatively, we can look for an approximant of $F$ in the space 
of cylindrical functionals. This yields a functional in the form 
(see Section \ref{sec:tensor}) 
\begin{equation}
f(a_1, ..., a_m, t_n)\simeq F([\theta],t_n),\qquad a_k=(\theta,\varphi_k).
\label{functional-approx_cyl}
\end{equation} 
For example, in a canonical tensor decomposition setting we have 
\begin{equation}
f(a_1, ..., a_m, t_n)\simeq \sum_{l=1}^r \prod_{j=1}^m
G_j^{l}(a_j,t_n),
\qquad G^{l}_{j}(a_j,t_n)=
\sum_{p=1}^P\beta_{jp}^{l}(t_n)\phi_p(a_j),\qquad 
\label{ff9}
\end{equation}
where $r$ is the separation rank. 
Replacing with $F([\theta],t)$ with $f(a_1,...,a_m,t)$ 
in \eqref{FSint1} or \eqref{AFE} yields 
a functional residual $R([\theta],t_n)$. 
For example, a substitution of \eqref{ff9} into 
the Crank-Nicolson scheme \eqref{AFE} yields
\begin{equation}
 R([\theta],t_{n})= \sum_{l=1}^r A([\theta],t_n)
 G^l_1((\theta,\varphi_1),t_n)\cdots G^l_m((\theta,\varphi_m),t_n) -E([\theta],t_n).
 \label{RAFE}
\end{equation}
Note that we have incorporated the local truncation 
error $\Delta t \tau_n$ within the residual $R([\theta],t_{n})$.

\subsection{CP-ALS Algorithm for FDEs with Implict Time Stepping} 
\label{sec:ALS_formulation}
{\color{r}
We have seen in Section \ref{sec:functional_least_squares} 
that the functional least squares fomulation of the FDE 
\eqref{linfde00} relies on minimizing the norm of the 
residual. Such residual can have different forms. In particular, if we discretize the FDE \eqref{linfde00} in time with 
the Crank-Nicolson method and represent its solution 
by a canonical polyadic (CP) tensor expansion 
then the residual takes the form \eqref{RAFE}. 
Its norm can be defined as 
\begin{equation}
 \left\|R([\theta],t_n)\right\|_W^2=
 \left\|\sum_{l=1}^r A([\theta],t_n)
 G^l_1((\theta,\varphi_1),t_n)\cdots G^l_m((\theta,\varphi_m),t_n) -
 E([\theta],t_n)\right\|_W^2,
\label{resnorm}
 \end{equation}
where $\|\cdot \|^2_W$ is induced by the functional 
inner product \eqref{continuous_ip_W0}. 
Recall that the functions $G_{k}^l((\theta,\varphi_k),t_n)$ 
are in the form 
\begin{equation}
G_{k}^l((\theta,\varphi_k),t_n)=\sum_{s=1}^Q 
\beta^l_{ks}(t_n)\phi_s((\theta,\varphi_k)),
\end{equation}
$\beta^l_{ks}(t_n)$ ($l=1,...,r$, $k=1,...,m$, $s=1,...,Q$) 
being the degrees of freedom. 
We look for a minimizer of \eqref{resnorm} 
computed in a {\em parsimonious way}. 
The key idea is to split the large scale optimization problem
\begin{equation}
\min_{\beta_{ks}^l(t_n)} \left\|R([\theta],t_n)\right\|_W^2
\end{equation}
into a sequence of optimization problems problems of smaller
dimension, which are solved sequentially and eventually in 
parallel \cite{Karlsson}. To this end, we define
\begin{equation}
\bm \beta_k(t_n)=[\beta^1_{k1}(t_n), ..., \beta^1_{kQ}(t_n),...,
\beta^r_{k1}(t_n), ..., \beta^r_{kQ}(t_n)]^T \qquad k=1,...,m.
\end{equation}
Note that the vector $\bm \beta_k(t_n)$ collects 
the degrees of freedom representing the solution 
functional along $(\theta,\varphi_k)$ at time $t_n$, i.e., 
the set of functions $\{G_k^1((\theta,\varphi_k),t_n),...,G^r_k((\theta,\varphi_k),t_n)\}$. 
Minimization of \eqref{resnorm} with respect to 
{\em  independent} variations of $\bm \beta_k(t_n)$
yields the sequence of {\em convex} optimization problems 
\begin{equation}
\min_{\bm \beta_1(t_n)} \left\|R([\theta],t_n)\right\|_W^2, \qquad 
\min_{\bm \beta_2(t_n)} \left\|R([\theta],t_n)\right\|_W^2, \qquad 
\cdots,\qquad 
\min_{\bm \beta_m(t_n)} \left\|R([\theta],t_n)\right\|_W^2.
\label{sequential_min}
\end{equation}
This is the set of equations defining the alternating 
least-squares (ALS) method. 
The Euler-Lagrange equations identifying the stationary points 
of \eqref{sequential_min} are 
\begin{equation}
\bm M_h(t_n) \bm \beta_h(t_n) = \bm f_h(t_n), \qquad h=1,...,m
\label{ADIS}
\end{equation}
where,
\begin{equation}
[\bm M_h(t_n)]_{qs}^{zl}=\int 
Q_{qh}^z([\theta],t_n)
Q_{sq}^l([\theta],t_n)
W([\theta])\mathcal{D}[\theta],
\label{Mn1}
\end{equation}
\begin{equation}
[\bm f_h(t_n)]^{z}_{q}=\int E([\theta],t_n)
Q_{qh}^z([\theta],t_n)W([\theta])D[\theta],
\label{kn1}
\end{equation}
\begin{equation}
Q_{qh}^z([\theta],t_n)= A([\theta],t_n)\phi_q((\theta,\varphi_h))
\prod_{\substack{j=1\\j\neq h}}^m G_j^z((\theta,\varphi_j),t_n).
\label{QQQ}
\end{equation}
The ordering of the matrix elements $[\bm M_h(t_n)]_{qs}^{zl}$ 
and the vector $[\bm f_h(t_n)]^{z}_{q}$ is the same as in
 \eqref{ORDERING}.
Note that minimizing the norm of the residual \eqref{resnorm} with 
respect $\bm \beta_h(t_n)$ is equivalent 
to impose orthogonality of \eqref{RAFE} with respect 
to the space spanned by the basis functionals 
$Q_{qh}^z([\theta],t_n)$. 
Indeed, the system \eqref{ADIS} is equivalent to 
\begin{equation}
\left(R([\theta],t_n),Q_{qh}^z([\theta],t_n)\right)_W=0 
\qquad \textrm{(fixed $h=1,...,m$),} 
\label{ALES}
\end{equation}
where $(,)_W$ is the functional inner 
product \eqref{continuous_ip_W0}.
The system of equations \eqref{ADIS} is 
symmetric and positive definite
\footnote{\color{r}The alternating least-squares formulation 
\eqref{ADIS} is very similar to alternating least squares 
formulation for nonlinear functionals we studied 
in Section \ref{sec:CP}. The main difference is that 
here we are trying to determine the canonical tensor decomposition 
of a linearly mapped functional, i.e., we solving the 
linear system $A([\theta],t_n)F([\theta],t)=E([\theta],t_n)$ 
with ALS, $F$ being represented as a 
canonical tensor decomposition. On the other hand, 
in Section \ref{sec:CP} we addressed the problem of 
computing the canonical tensor decomposition of $F([\theta])$ 
given $E([\theta])$ (compare the residual \eqref{RAFE} with \eqref{Fresidual}). These two problems are equivalent 
if the linear operator $A$ is invertible.}. 
This allows us to use well-known high-performance algorithms 
to compute the solution, e.g., the conjugate gradient 
method (\cite{Quarteroni}, p.152). 
It is important to emphasize that the minimization of the residual 
\eqref{resnorm} is basically a {\em fixed point problem} 
involving the vector 
\begin{equation}
\bm \beta(t_n)=[\bm \beta_1(t_n)\quad \cdots \quad \bm\beta_m(t_n)]
\end{equation}
The ALS method aims at solving such fixed point 
problem by splitting it into a sequence of 
linear problems \eqref{ADIS} which are solved iteratively 
for each $h$. The criterion is to freeze 
all $\bm \beta_j(t_n)$ ($j=1,...,m$, 
$j\neq h$) when solving for $\bm \beta_h(t_n)$.  
We recall that convergence of CP-ALS iterations is, in general,  
not granted (see Section \ref{sec:CP}). To overcome 
this problem, additional regularization terms may be necessary 
\cite{Acar,Battaglino}.

\paragraph{Evaluation of the Functional Integrals}
The ALS coefficients \eqref{Mn1} and \eqref{kn1}
are defined by functional integrals involving cylindrical functionals. 
The computation of such integrals is briefly addressed in 
Appendix \ref{app:change of variables in SSE}. Hereafter 
we summarize the main results. We first restrict 
$\theta(x)$ to the space of functions 
\eqref{SoF},  i.e., we assume that $\theta(x)$ can be written 
as
\begin{equation}
{\theta}(x)=\sum_{k=1}^m a_k \varphi_k(x),\qquad a_k= (\theta,\varphi_k).
\label{KOE}
\end{equation}
In this hypothesis, the basis functionals \eqref{QQQ} 
become multivariate functions of $(a_1,...,a_m)$, i.e., 
\begin{equation}
Q^z_{qh}(a_1,....,a_m,t_n)=A(a_1,...,a_m,t_n)\phi_q(a_h)\prod_{\substack{j=1\\j\neq h}}^m 
G_j^z(a_j,t_n).
\label{QQ}
\end{equation}
Similarly, the coefficients \eqref{Mn1} and \eqref{kn1} 
can be explicitly written as multivariate 
integrals\footnote{\color{r} In equations 
\eqref{Mn1f} and \eqref{kn1f} we have we set the 
weight function $W(a_1,...,a_n)$ equal to one. This is always 
possible provided the support of the integrands is compact 
(see Appendix \ref{app:change of variables in SSE}).}
\begin{align}
[\bm M_h(t_n)]_{qs}^{zl}=\int_{-b}^b\cdots\int_{-b}^b 
Q^z_{qh}(a_1,....,a_m,t_n) Q^l_{sq}(a_1,...,a_m,t_n)da_1\cdots da_m,
\label{Mn1f}
\end{align}
\begin{equation}
 [\bm f_h(t_n)]^{z}_{q}=\int_{-b}^b\cdots\int_{-b}^b E(a_1,...,a_m,t_n)
Q^z_{qh}(a_1,..,a_m,t_n) da_1\cdots da_m,
\label{kn1f}
\end{equation}
The quantity $A(a_1,...,a_m,t_n)$ appearing in \eqref{QQ} 
is the discrete version of the functional differential operator 
$A([\theta],t_n)$, i.e., it is a linear operator in the form
 \begin{equation}
 A(a_1,...,a_m,t_n)=I-\frac{\Delta t}{2} L(a_1,...,a_m,t_n)
 \end{equation}
where $L(a_1,...,a_m,t_n)$ is the finite-dimensional 
version of the operator $L([\theta],t_n)$ in \eqref{linfde00}.
In general, the integrals \eqref{Mn1f} and \eqref{kn1f} can 
be computed numerically only for a relatively small number of 
variables $a_k$ ($k=1,...,m$). However, if we 
assume that the operator $A(a_1,...,a_m,t_n)$ is {\em separable}, i.e., 
\begin{equation}
A(a_1,...,a_m,t_n)=\sum_{k=1}^{r_A} A_1^k(a_1,t_n)\cdots A_m^k(a_m,t_n),
\label{A_}
\end{equation}
then the cost of computing such integrals scales {\em linearly} with the 
dimension $m$ of the space, since 
\eqref{Mn1f} and \eqref{kn1f} can be factored as 
products of one-dimensional integrals.  
In equation \eqref{A_}, $A_i^k(a_i,t_n)$ are 
one-dimensional linear operators, while $r_A$ is the separation 
rank of the operator $A(a_1,...,a_m,t_n)$.
With the operator decomposition \eqref{A_} available, 
we can represent the multivariate fields \eqref{QQQ} as
\begin{equation}
Q_{qs}^l(a_1,...,a_m,t_n)=\sum_{k=1}^{r_A} A_q^k(a_q,t_n)\phi_s(a_q)\prod_{\substack{j=1\\j\neq q}}^m
A_j^k(a_j,t_n)G_j^{l}(a_j,t_n).\label{qqqs}
\end{equation}
This yields the following representation of the matrix coefficients \eqref{Mn1f} 
\begin{align}
\left[\bm M_q(t_n)\right]_{sh}^{lz}=&
\int_{-b}^b\cdots \int_{-b}^b Q_{qs}^l(a_1,...,a_m,t_n)
Q_{qh}^z(a_1,...,a_m,t_n)da_1\cdots da_m\nonumber\\
=& \sum_{k,e=1}^{r_A}\int_{-b}^b
 A_q^k(a_q,t_n)\phi_s(a_q) A_q^e(a_q,t_n)\phi_h(a_q) da_q \times\nonumber \\ 
& \prod_{\substack{j=1\\j\neq q}}^m\int_{-b}^b
A_j^k(a_j,t_n)G_j^{l}(a_j,t_n)A_j^e(a_j,t_n)G_j^{z}(a_j,t_n)da_j.
\label{integrals_9}
\end{align}
i.e., sums of products of {one-dimensional 
integrals}\footnote{Note that we can express all integrals 
at the right hand side of \eqref{integrals_9} in terms of the integrals 
\begin{equation}
\int_{-b}^b A_q^k(a_q,t_n)\phi_s(a_q)A_q^e(a_q,t_n)\phi_h(a_q)da_q
\end{equation}
where $q=1,...,m$, $k,e=1,...,r_A$, $s,h=1,...,Q$. In fact,
\begin{equation}
\int_{-b}^b A_j^k(a_j,t_n)G_j(a_j,t_n) A_j^e(a_j,t_n)G^z_j(a_j,t_n) da_j = \sum_{s,h=1}^Q \beta_{js}^l(t_n) \beta_{jh}^z(t_n) 
\int_{-b}^b A_j^k(a_j,t_n)\phi_s(a_j)A_j^e(a_j,t_n)\phi_h(a_j)da_j.
\end{equation} 
}. 
Let us provide a simple example.

\vs
\noindent 
{\em Example 1:} Consider the time-independent 
functional differential operator \eqref{FDOP} and the 
associated operator defined in \eqref{An}
\begin{equation}
A([\theta])=I+\frac{\Delta t}{2}\int_{0}^{2\pi} \theta(x) 
\frac{\partial }{\partial x}\frac{\delta }{\delta \theta(x)} dx.
\end{equation}
Evaluating $A([\theta])$ in the finite-dimensional function space 
\eqref{SoF} yields 
\begin{equation}
 A(a_1,...,a_m)=I+\frac{\Delta t}{2}\sum_{k=1}^m \left(\sum_{j=1}^m a_j 
 \int_{0}^{2\pi}\varphi_k(x)\frac{d\varphi_j(x) }{dx}dx\right)
\frac{\partial }{\partial a_k}.
\end{equation}
Therefore, $A(a_1,...,a_m)$ is a separable operator with separation rank $r_A=m^2+1$.

\subsubsection{Collocation Setting}
Consider the sequence of linear systems \eqref{ADIS}, and let 
$\{\phi_s(a_k)\}$ be a cardinal basis associated with the
set of collocation nodes $\{a_{k1},...,a_{kQ}\}$. 
For simplicity, we consider the same set of nodes in each dimension. 
In this assumption, the integrals defining the the matrix entries 
\eqref{integrals_9} can be significantly simplified. For example,
\begin{equation}
\int_{-b}^b A_q^k(a)\phi_s(a) A_q^e(a)\phi_h(a)da \simeq 
\underbrace{\bm A_q^k \bm W \bm A_q^e}_{\bm K_q^{ke}}
\label{eq:50}
\end{equation}
where $\bm A_q^k$ is the matrix representation of the 
operator $A_q^k(a)$ (collocation version), and $\bm W$ 
is a diagonal matrix of integration weights. 
The matrix $\bm K_q^{ke}$ is $Q\times Q$ for all 
$k,e=1,...,r_A$ and all $q=1,...,m$. With $\bm K_q^{ke}$ 
available, it is easy to determine the matrix representation of the 
integral 
\begin{equation}
\int_{-b}^b A_q^k(a)G^l_j(a) A_q^e(a)G_j^z(a)da=
\sum_{s,h=1}^Q \beta_{js}^l \beta_{jh}^z \int_{-b}^b A_q^k(a)\phi_s(a) A_q^e(a)\phi_h(a)da \simeq  
\bm \beta_j^T \bm K_q^{ke} \bm \beta_j
\end{equation}
where $\bm \beta_j$ here is a matrix that has  $G_j^z(x_{jp})$ 
($p=1,..,Q$) as $z$-th column.
This yields the following matrix 
\begin{equation}
\bm M_q = \sum_{k,e=1}^{r_A}  [\bm R^{ke}_q]^T\otimes \bm K_q^{ke}.
\end{equation}
where 
\begin{equation}
\bm R^{ke}_q=\prod_{\substack{j=1\\ 
j\neq q}}^m \bm \beta_j^T \bm K_q^{ke}\bm \beta_j.
\end{equation}

\subsection{Tensor Formats for FDEs with Explicit Time Stepping}
Consider the following finite-dimensional form of the linear 
FDE \eqref{linfde00} 
\begin{equation}
\frac{\partial f}{\partial t}=L f +h,
\label{PDEL}
\end{equation}
where $f(a_1,..,a_m,t)$ is a {\em tensor format} that approximates 
a functional $F([\theta],t)$ and $L(a_1,...,a_m,t)$ is the 
linear operator arising from the discretization of the functional 
linear operator $L([\theta],t)$, and $h(a_1,...,a_m,t)$ is the 
tensor format that appoximates $H([\theta],t)$. 
In particular, consider the case where $h=0$, $L$ is time-independent 
and time-integration follows the Adams-Bashforth scheme 
\eqref{FSintti}, i.e., 
 \begin{equation}
f_{n} = f_{n-1}+\frac{\Delta t}{2} L\left(3f_{n-1}-f_{n-2}\right).
\end{equation}
In the last equation we employed the shorthand notation 
$f_n=f(a_1,...,a_m,t_n)$. Assuming, that the 
operator $L$ is separable, e.g., 
\begin{equation}
L = -\sum_{i,j=1}^m a_iC_{ij}\frac{\partial }{\partial a_j},
\end{equation}
then a greedy computation of the 
tensor format $f_{n}$ involves the following steps:
\begin{enumerate}
\item compute a low rank representation of $\hat{f}_{n-1}=3f_{n-1}/2 - f_{n-2}/2$,
\item compute a low rank representation of 
$f_{n}=p_{n-1} +\Delta t L \hat{f}_{n-1}$.
\end{enumerate}
The need for a low rank representation is clear: 
any algebraic operation between tensors, including the 
application of a linear operator, increases the separation 
rank. Therefore, efficient {\em rank reduction methods} are needed to 
avoid an explosion of the number of terms when solving 
the PDE \eqref{PDEL} with tensor methods. 
Among them we recall methods based on alternating least 
squares \cite{Karlsson,Battaglino,Etter}, 
hierarchical Tucker formats \cite{Grasedyck2015,Kolda} 
or block coordinate descent methods \cite{Xu}. Disregarding the 
particular tensor format employed to represent the solution 
functional, we emphasize that the development of rubust and 
efficient rank-reduction algorithms is an active area 
of research \cite{Bachmayr,NouyHUQ}.
}

\section{Numerical Results: Functionals}
\label{sec:numerical results functionals}
In this Section we provide numerical results and examples on  
functional approximation. In particular, we discuss  
polynomial functional interpolants and functional 
tensor methods.

\subsection{Linear functionals}
\label{sec:results linear functionals}
Consider the linear functional 
\begin{align}
F([\theta])=\int^{2\pi}_0 K_1(x) \theta(x)dx
\label{linF}
\end{align}
on the Hilbert space of square integrable periodic functions 
in $[0,2\pi]$ 
\begin{equation}
 D(F)=\{\theta\in L_2([0,2\pi])\,|\, \theta(0)=\theta(2\pi)\}.
 \label{thespace}
\end{equation}
Our aim is to represent $F([\theta])$ in terms of a functional 
interpolant in $D(F)$, i.e.,  
\begin{equation}
F([\theta])= \sum_{k=1}^\infty  F([\theta_i])g_i([\theta]).
\end{equation}
where $g_i([\theta])$ are cardinal basis functionals 
and $\theta_i(x)$ are interpolation nodes in $D(F)$. 
In particular,  we choose $\theta_i(x)=\varphi_i(x)$ 
where $\{\varphi_1(x),\varphi_2(x),...\}$ is an orthonormal  
basis in $D(F)$.
Assuming that $K_1(x)$ is in $D(F)$, i.e.,
\begin{equation}
 K_1(x)=\sum_{k=1}^\infty (K_1,\varphi_k)\varphi_k(x), 
\end{equation}
it follows from \eqref{linF} that  
\begin{align}
F([\theta])=&\sum_{k=1}^\infty (K_1,\varphi_k)(\varphi_k,\theta),\nonumber\\
=&\sum_{k=1}^\infty F([\varphi_k])(\varphi_k,\theta).
\label{exactF}
\end{align}
This can be written as 
\begin{align}
F([\theta])=\sum_{k=1}^\infty F([\varphi_k])g_k([\theta]),  \qquad\textrm{where}\qquad
g_k([\theta])=(\varphi_k,\theta).
\label{seriesLF}
\end{align}
Note that this representation coincides with Porter's series 
expansion \eqref{Porter_interpolant}-\eqref{gi_0} on the index 
set $\mathcal{I}=\{1\}$.

\subsubsection{Polynomial Functional Interpolation} 
\begin{figure}
\centerline{\hspace{0.5cm}(a)\hspace{7cm}(b)}
\centerline{\includegraphics[height=5.5cm]{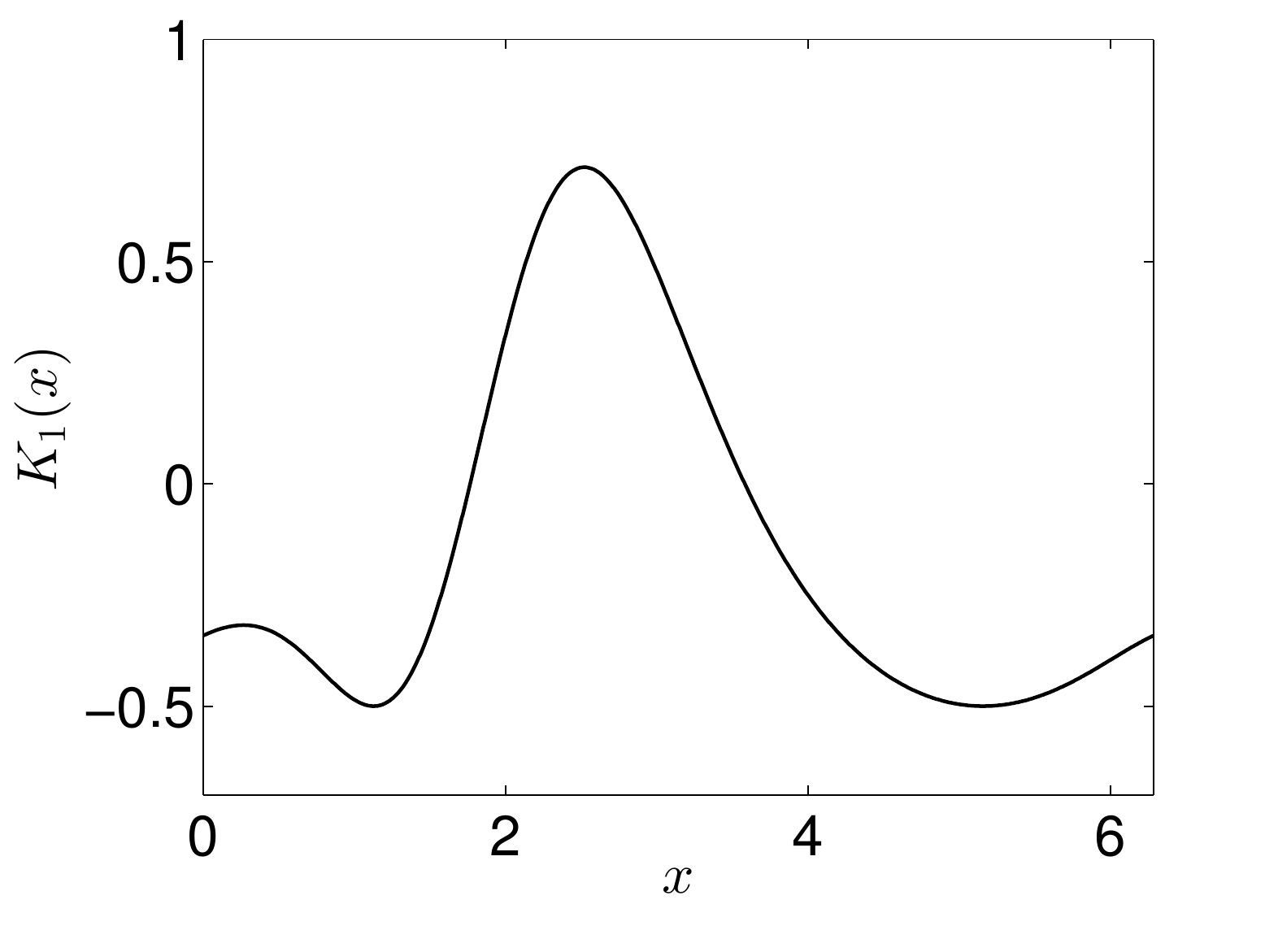}
            \includegraphics[height=5.5cm]{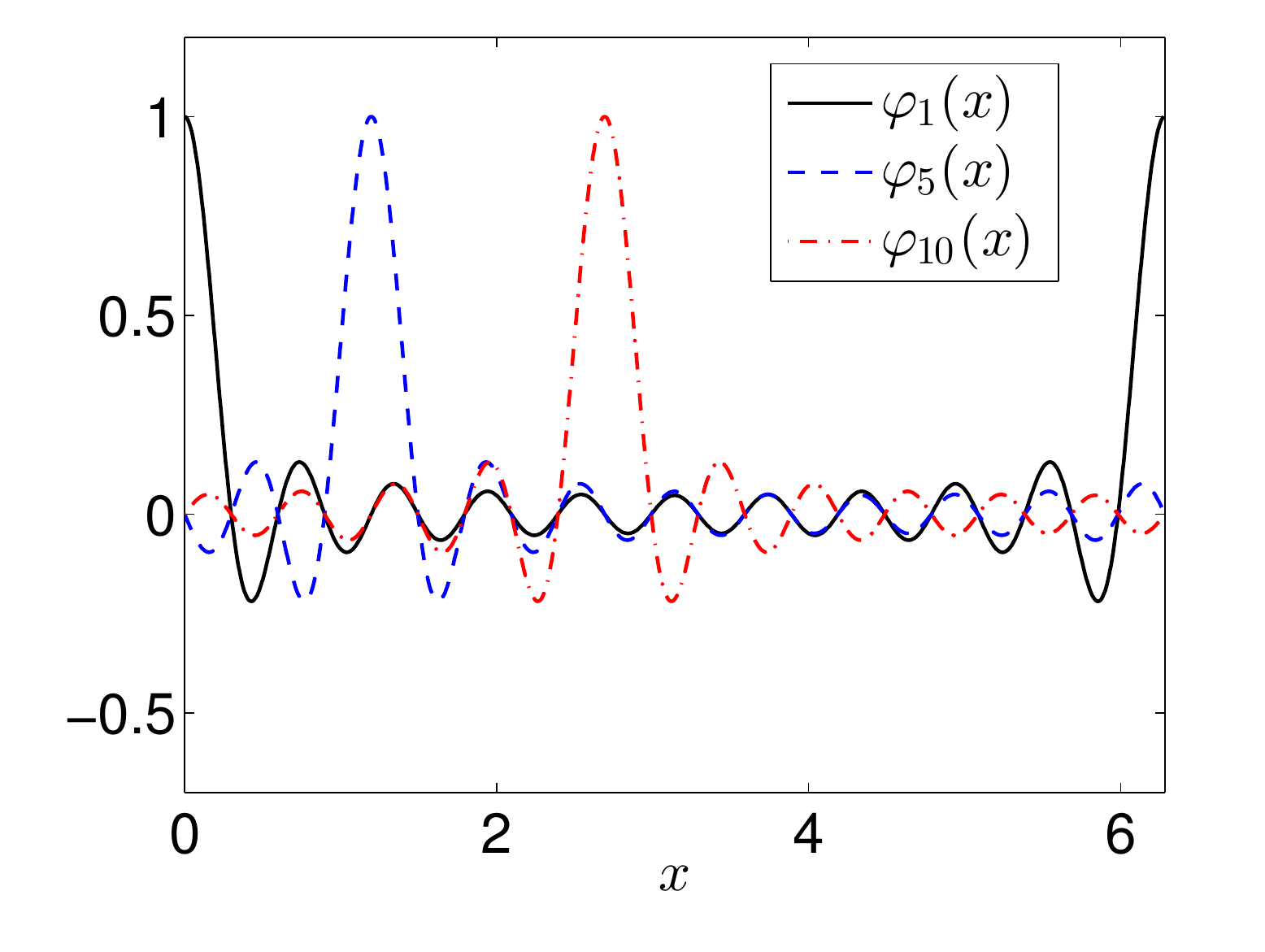}}
\caption{(a) Kernel function defining the linear functional \eqref{linF}. (b) 
Three elements of the orthogonal basis function set \eqref{setbf} ($m=20$).}
\label{fig:kernel1}
\end{figure}

\begin{figure}
\centerline{\hspace{0.5cm}(a)\hspace{8cm}(b)}
\centerline{\includegraphics[height=6.5cm]{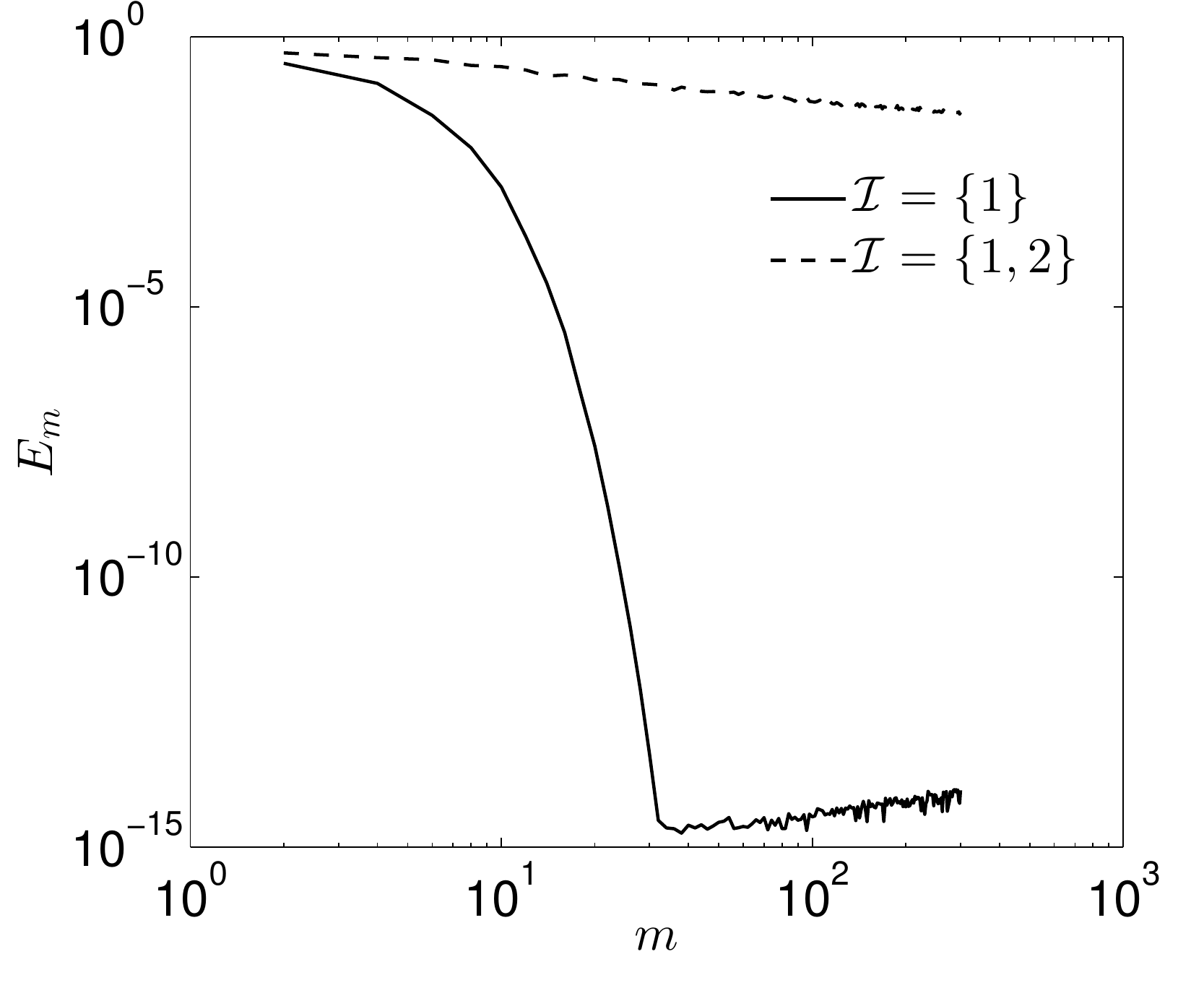}
            \includegraphics[height=6.5cm]{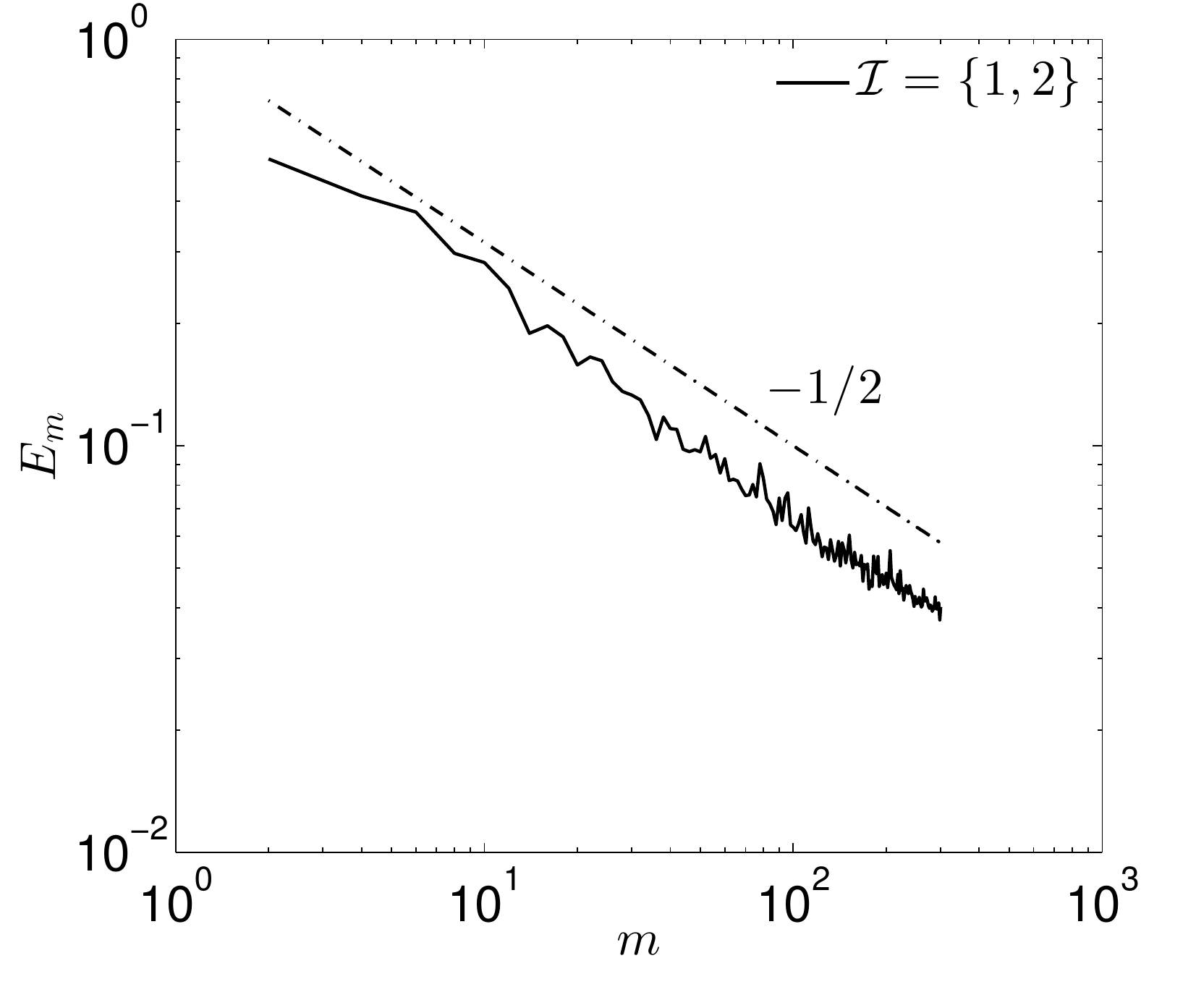}}
\caption{Linear functionals. Functional interpolation errors obtained by 
using Porter's method. Shown are the pointwise 
errors \eqref{Em} versus $m$ for the Gaussian ensemble 
$\mathcal{G}_{300}(1)$ and Porter's basis functionals
\eqref{linPorter} ($\mathcal{I}=\{1\}$), and 
\eqref{quadPorter} ($\mathcal{I}=\{1,2\}$). 
The expansion corresponding to $\mathcal{I}=\{1\}$ converges 
exponentially for obvious reasons, while the expansion corresponding 
to $\mathcal{I}=\{1,2\}$ has a $O(m^{-1/2})$ convergence rate 
due to insufficient interpolation nodes (Fig. (b)).}
\label{fig:results_linear}
\end{figure}

Let us study numerically an interpolation problem involving
a specific kernel. To this end, we set 
\begin{equation}
K_1(x)=e^{\sin(x)}(1+\sin(\cos(x)-2)-\frac{1}{2}\qquad \textrm{(Fig. \ref{fig:kernel1})}, 
\label{K1}
\end{equation}
and define the following interpolation nodes in $D(F)$ 
\begin{equation}
\varphi_{k+1}(x)=\frac{1}{m+1}\displaystyle\frac{\sin\left(\displaystyle\frac{m+1}{2}\left(x-x_k\right)\right)}
{\sin\left(\displaystyle\frac{x-x_k }{2} \right)},\qquad x_k=\frac{2\pi}{m+1}k\qquad k=0,1,...,m
\label{setbf}
\end{equation}
$m$ being an even natural number. These are well-known 
nodal trigonometric polynomials (\cite{Hesthaven}, p. 28) 
satisfying the orthogonality conditions
\begin{equation}
(\varphi_j,\varphi_k)=\frac{2\pi}{m+1}\delta_{kj}, \qquad k,j=1,..., m+1.
\end{equation} 
Clearly, if $K_1(x)$ is in the span of 
$\{\varphi_1,...,\varphi_{m+1}\}$ then the exact 
representation of the linear functional $F([\theta])$ 
involves no more than $m+1$ terms, i.e., a truncation 
of the series \eqref{seriesLF} to $m+1$ terms is exact. 
The kernel \eqref{K1} is not in such span. 
Next, we construct a polynomial functional 
that interpolates $F([\theta])$ at the nodes \eqref{setbf}. 
Specifically we consider Porter's construction 
(Section \ref{sec:porter}), which is very easy to implement and 
equivalent to Khlobystov polynomials for uniquely solvable  
interpolation problems. Such interpolants are dense 
in the space of linear functionals if we consider the 
set of nodes \eqref{SNq1}, i.e., 
\begin{equation}
 \widehat{S}^{(m)}_1=\left\{\theta_i(x)\in D(F)\,\,|\,\, 
 \theta_i(x)=\varphi_{i}(x) \quad i=0,...,m\right\}.
\end{equation}
This means that to identify linear functionals  
it is sufficient to represent them relative to orthogonal bases 
(see Section \ref{sec:variational}). As we have noticed 
in Section \ref{sec:functional_polynomial_interpolation}, 
this is not the case for higher-order polynomial functionals
or general nonlinear functionals. 
Now, let us consider the general expression of Porter's interpolants 
\begin{equation}
F\left([\theta]\right)\simeq\sum_{k=1}^{m+1} F\left([\varphi_k]\right) g_{k}([\theta]),\qquad 
 g_{k}([\theta])=\sum_{j=1}^{m+1} H_{jk}^{-1}
 \sum_{p\in\mathcal{I}}\left(\varphi_j,\theta\right)^p,
\end{equation}
where the matrix $H_{ij}$ is given in \eqref{matH}.
Depending on how we choose the index set $\mathcal{I}$ we 
have different expressions for the basis functionals $g_i([\theta])$. 
Specifically,

\begin{enumerate}

\item Constant polynomial functionals ($\mathcal{I}=\{0\}$).
This case is degenerate and it requires Moore-Penrose pseudo-inversion of 
the matrix \eqref{matH}. This yields the basis functionals
\begin{equation}
g_k([\theta])=\frac{1}{m+1}.
\end{equation}

\item Homogeneous polynomial functionals of first order ($\mathcal{I}=\{1\}$): 
\begin{equation}
g_k([\theta])= \left[\frac{2\pi}{m+1}\right]^{-1} \left(\theta,\varphi_k\right).
\label{linPorter}
\end{equation}
\item Quadratic polynomial functionals ($\mathcal{I}=\{1,2\}$) on the basis 
set $\left\{\varphi_1,...,\varphi_{m+1}\right\}$. This yields basis 
functionals
\begin{equation}
g_k([\theta]) = \left[\left(\frac{2\pi}{m+1}\right)+
\left(\frac{2\pi}{m+1}\right)^2\right]^{-1} 
\left(\left(\theta,\varphi_k\right)+\left(\theta,\varphi_k\right)^2\right).
\label{quadPorter}
\end{equation}

\end{enumerate}
In Figure \ref{fig:results_linear} we plot the pointwise error 
\begin{equation}
E_m= \sup_{\theta\in \mathcal{G}_q}\left|F([\theta])-
\sum_{k=1}^{m+1}F([\varphi_k])g_k([\theta])\right| 
\label{Em}
\end{equation}
versus $m$ for test functions $\theta$ in the Gaussian ensemble
\begin{equation}
 \mathcal{G}_q(\sigma)=\left\{\theta(x)\in D(F)\,\left|  
 \, \theta(x)=\sigma\sum_{k=1}^{q+1}b_k\varphi_k(x),\quad 
 \{b_1,...,b_{m+1}\}\quad \textrm{i.i.d. Gaussian}\right.\right\}, 
 \label{Gq}
\end{equation}
and interpolants corresponding to the index sets $\mathcal{I}=\{1\}$ and 
$\mathcal{I}=\{1,2\}$ (bases  
\eqref{linPorter} and \eqref{quadPorter}).
The reason why we obtain exponential convergence with the index set 
$\mathcal{I}=\{1\}$ is obvious: convergence of the 
polynomial interpolant is basically defined by the 
convergence of the trigonometric series of $K_1(x)$.  
On the other hand, Porter's interpolant corresponding
to $\mathcal{I}={1,2}$ shows an algebraic convergence at rate 
$O(m^{-1/2})$. The reason is that the basis 
function set $\{\varphi_1,...,\varphi_{m+1}\}$ 
does not have enough elements to correctly identify 
the quadratic part of the interpolant, which is zero in this case. 
In fact, the off-diagonal terms in \eqref{sof} or \eqref{ssof} 
cannot be identified by using an orthogonal basis, 
unless we consider a set of nodes in the form 
$\left\{\varphi_i, (\varphi_{i}+\varphi_j)\right\}$ 
where $i,j=1,...,m+1$, and  $j\geq i$.

\subsubsection{Canonical Tensor Decomposition} 
We look for a representation of \eqref{linF}  
in the form \eqref{functional-SSE}. 
The basis functions $G_i^l$ can be equivalently determined 
in an alternating Galerkin or least squares setting 
by solving the system of equations \eqref{SSYSTEM} 
with forcing given by
\begin{equation}
f^n_{jh}=\int_{-b}^b\cdots\int_{-b}^b 
\sum_{p=1}^{m+1} k_p a_p 
\phi_h(a_j)\prod_{\substack{k=1\\k\neq j}}^{m+1} G^n_k(a_k)da_1\cdots da_{m+1},
\label{fn3}
\end{equation}
\begin{equation}
k_p=\int_{0}^{2\pi} K_1(x)\varphi_p(x)dx.
\end{equation}
All integrals in \eqref{fn3} can be reduced to 
products of one-dimensional integrals (linear functionals are 
fully separable). In addition, if we use a polynomial basis 
for $G^l_i$, then we can represent \eqref{linF} exactly 
as a product of constants and 
linear polynomials. In fact,
\begin{align}
\sum_{l=1}^rG^l_1(a_1)\cdots G^l_{m+1}(a_{m+1})
=\sum_{l=1}^{m+1} k_l a_l, 
\end{align}
which means that the exact separation rank is $r={m+1}$, 
\begin{equation}
G^l_l(a_l)=a_l \qquad \textrm{and}\qquad  
\prod_{\substack{k=1\\k\neq l}}^{m+1} G_k^l(a_k)=k_l.
\end{equation}

\subsection{Quadratic Functionals}

\begin{figure}[t]
\centerline{$H_2(x,y)$\hspace{5.5cm}$K_2(x,y)$}
\centerline{\includegraphics[height=5cm]{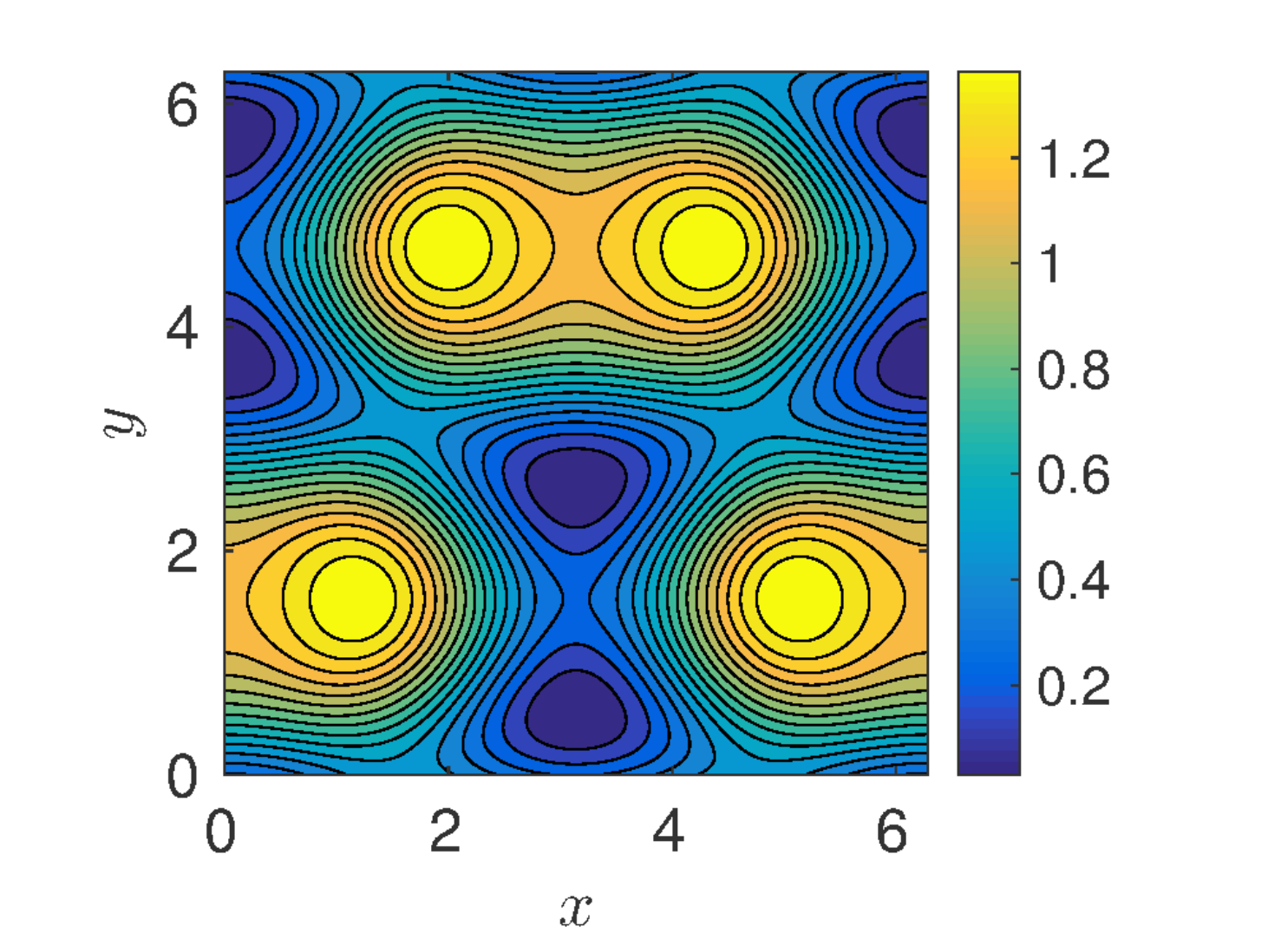} 
\includegraphics[height=5cm]{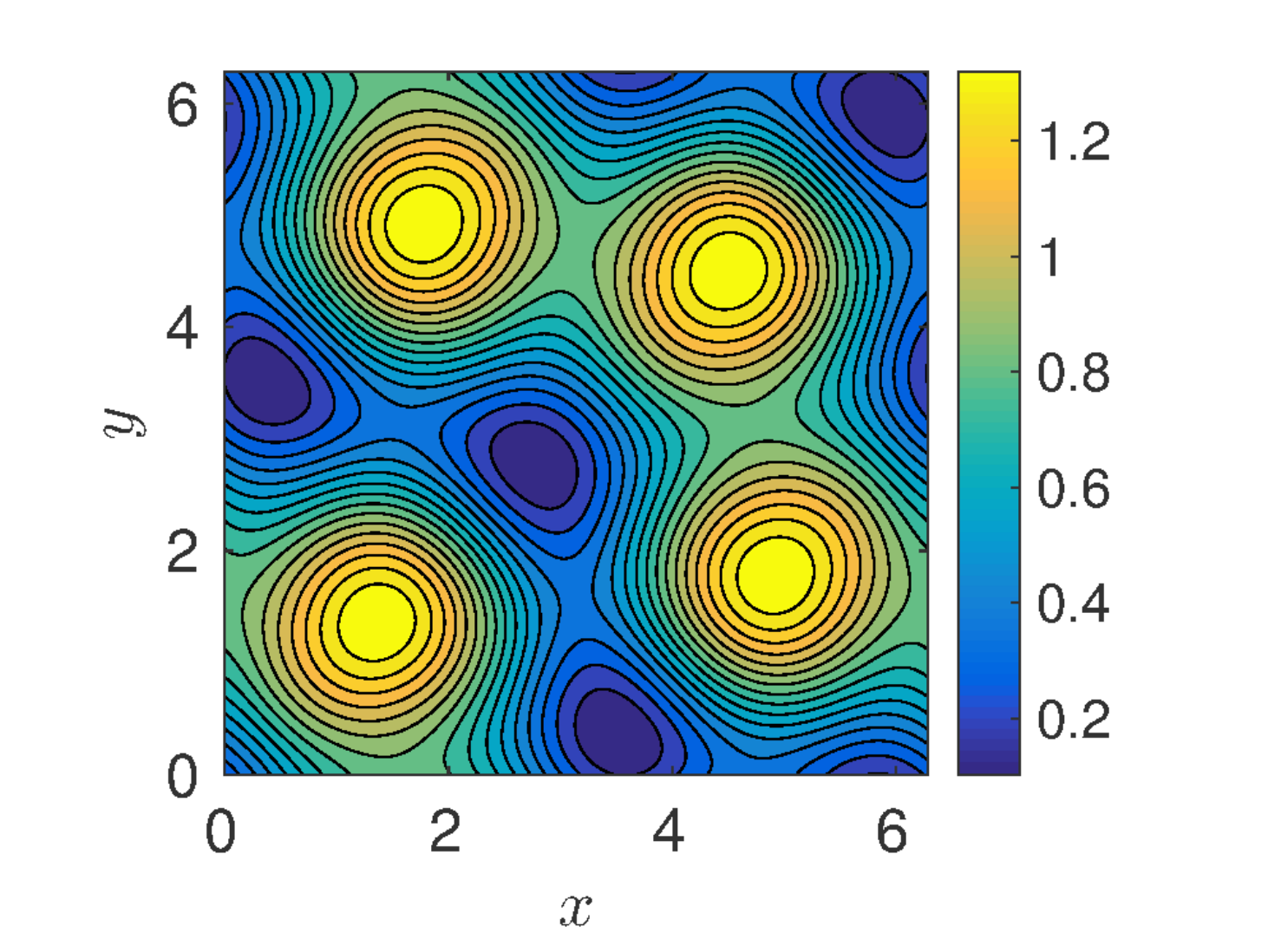}}
\caption{Quadratic functional \eqref{nonlinF}. Kernel function \eqref{K2} and its 
symmetrized version \eqref{K2sym}. Integrating $H_2(x,y)\theta(y)\theta(x)$ 
or $K_2(x,y)\theta(y)\theta(x)$ over $x$ and $y$ produces 
exactly the same result, independently on $\theta$.}
\label{fig:kernel2}
\end{figure}

Consider the quadratic functional
\begin{align}
F([\theta]) =  5 + \int^{2\pi}_0 K_1(x) \theta(x)dx+
\int^{2\pi}_0\int^{2\pi}_0 H_2(x,y)\theta(x)\theta(y)dxdy
\label{nonlinF}
\end{align}
defined on the space of square integrable periodic functions \eqref{thespace}. 
We set $K_1(x)$ as in \eqref{K1} and $H_2(x,y)$ as
\begin{equation}
 H_2(x,y)=\sin(\cos(x)+\sin(y))\sin(y)+\frac{1}{2}\cos(\cos(x)).
 \label{K2}
\end{equation}
Replacing $H_2(x,y)$ with the symmetrized kernel 
(see Figure \ref{fig:kernel2})
\begin{equation}
 K_2(x,y)=\frac{1}{2}\left(H_2(x,y)+H_2(y,x)\right)
 \label{K2sym}
\end{equation}
does not change the functional \eqref{nonlinF}. 
\subsubsection{Polynomial Functional Interpolation}
\begin{figure}[t]
\centerline{\includegraphics[height=5cm]{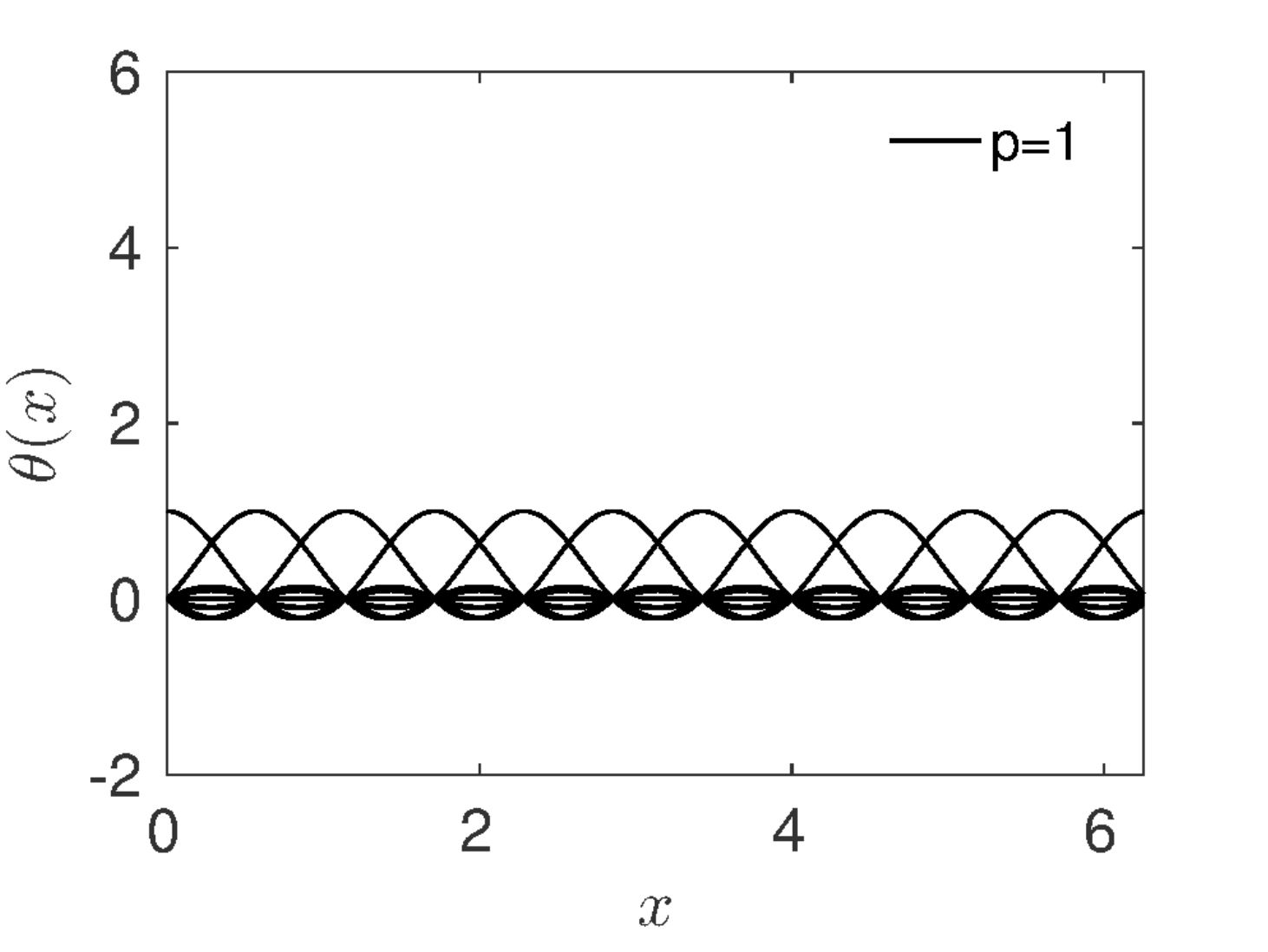}
            \includegraphics[height=5cm]{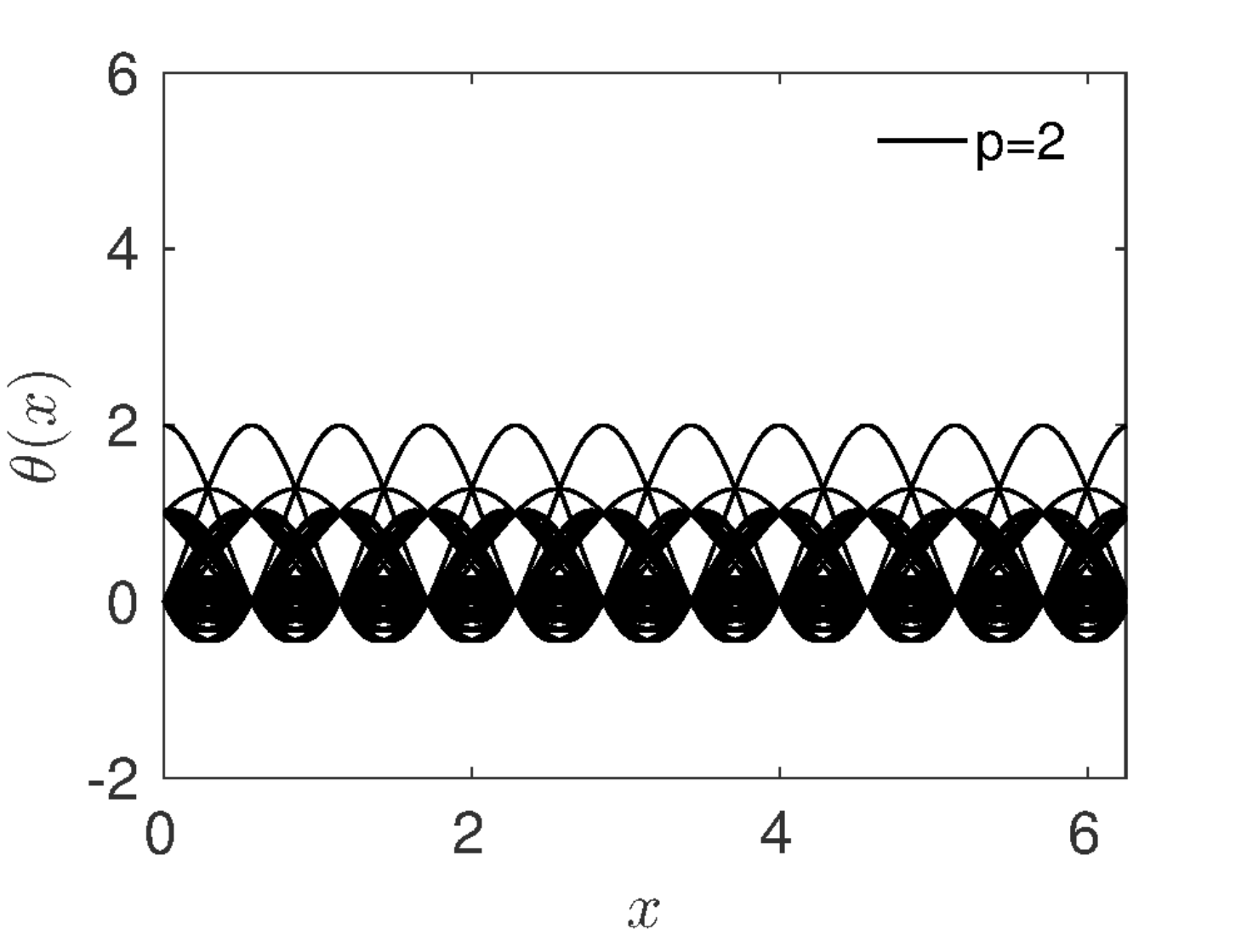}}
\centerline{\includegraphics[height=5cm]{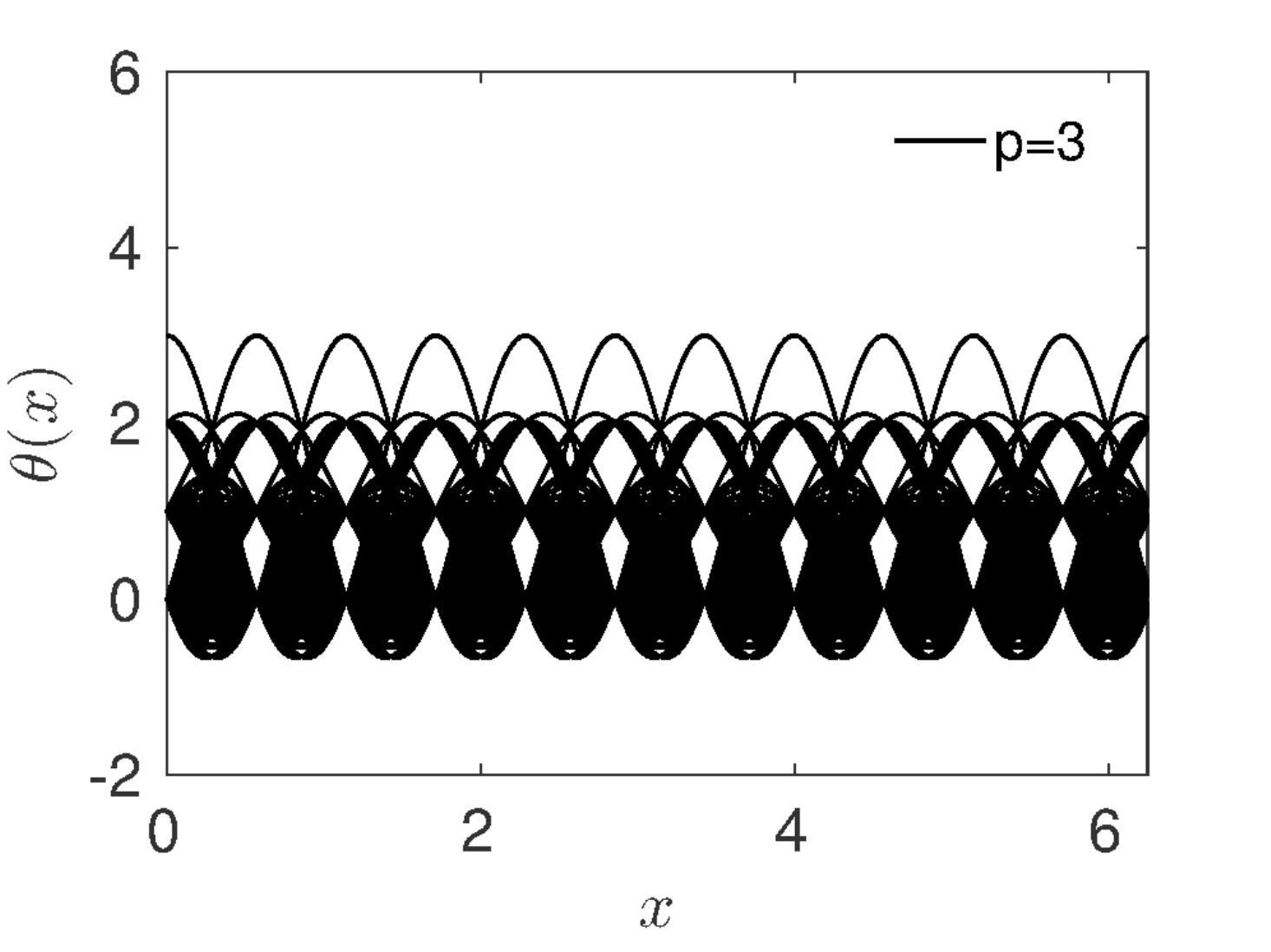}
            \includegraphics[height=5cm]{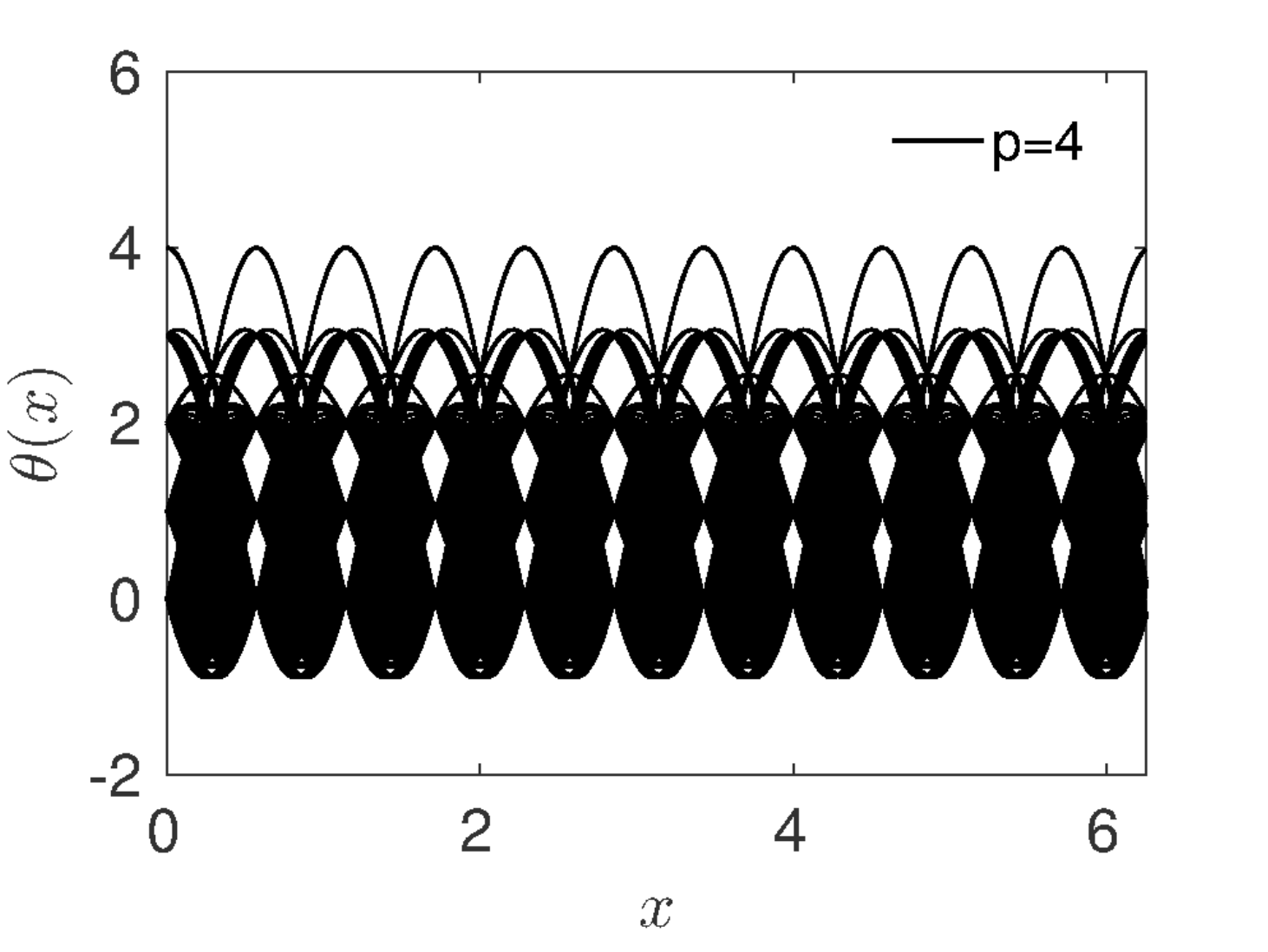}}            
\caption{Interpolation nodes in the space of periodic functions (each function is a node). 
Here we plot the elements in the sets $\widehat{S}^{(m+1)}_p$ (see equation \eqref{SNq2}) 
for $m=10$ and $p=1,2,3,4$. The total number of elements within each 
set is $\#\widehat{S}^{(11)}_1=12$, $\#\widehat{S}^{(11)}_2=78$, 
$\#\widehat{S}^{(11)}_3=364$ and $\#\widehat{S}^{(11)}_4=1365$.}
\label{fig:Khlobystov_basis}
\end{figure}

\begin{figure}[t]
\centerline{\includegraphics[height=6.5cm]{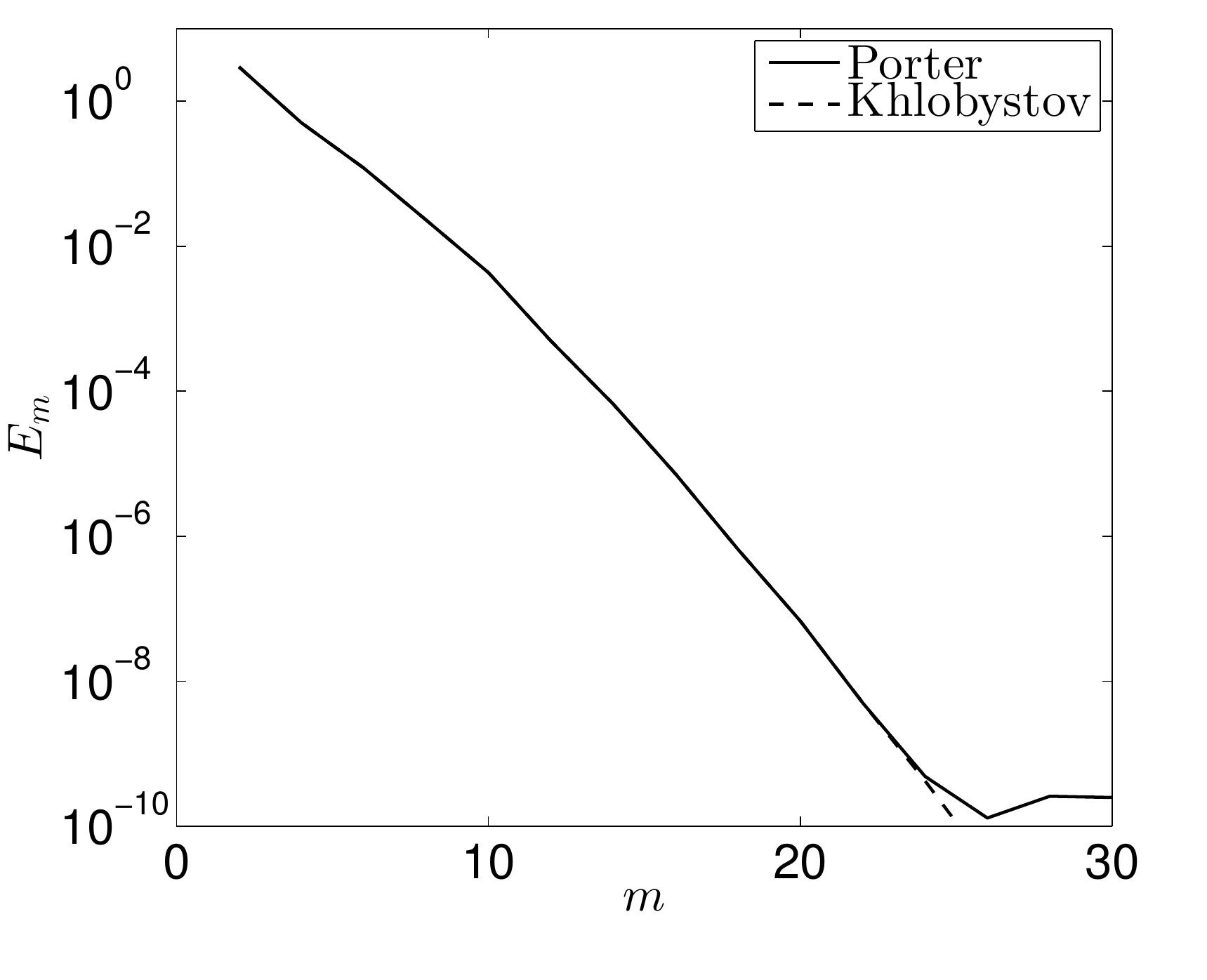}}
\caption{Interpolation of the quadratic functional \eqref{nonlinF} 
using Porter and Khlobystov approaches. 
Shown are the pointwise errors \eqref{Em2} versus $m$ for the Gaussian ensemble 
$\mathcal{G}_{50}(1)$. The small deviation between Porter's and Khlobystov's error plots
observed at $m=24$ is due to inaccuracies in the computation of the inverse of \eqref{matH}.
The number of interpolation nodes required to achieve 
accuracy of about $10^{-7}$ ($m=20$) is 
$\#\widehat{S}^{(21)}_2=253$ (see equation \ref{SCardinality}).}
\label{fig:results_quadratic}
\end{figure}

We first show that Porter's interpolants corresponding 
to the index set $\mathcal{I}=\{0,1,2\}$ and Khlobystov 
interpolants \eqref{p2interp} coincide when the interpolation 
problem is uniquely solvable, e.g., when we consider the set 
of nodes  \eqref{SNq2}, hereafter rewritten for convenience
\begin{equation}
\widehat{S}^{(m+1)}_2 = 
\left\{0,\{\varphi_i\}_{i=1}^{m+1},\{\varphi_i+\varphi_j\}_{j\geq i=1}^{m+1}\right\}.
\label{test_f_s1}
\end{equation}
In Figure \ref{fig:Khlobystov_basis} we plot the elements of such set 
for $m=10$. 
In Figure \ref{fig:results_quadratic} we plot the pointwise error  
\begin{equation}
E_m= \sup_{\theta\in \mathcal{G}_{50}(1)}\left|F([\theta])-\Pi_2([\theta])\right| 
\label{Em2}
\end{equation}
versus $m$ for functions $\theta(x)$ in the Gaussian ensemble
$\mathcal{G}_{50}(1)$ (see equation \eqref{Gq}). In \eqref{Em2}, $\Pi_2$ represents either 
Porter's or Khlobystov interpolant at the set of nodes \eqref{test_f_s1}.
It is seen that both methods achieve exponential convergence rate 
when interpolating the polynomial functional \eqref{nonlinF}. 
This is due to the fact that such functional expansions are 
basically approximating the kernel functions 
\eqref{K1} and \eqref{K2} in terms of a Fourier spectral basis. 
This also establishes the full equivalence of Porter's and Khlobystov's 
interpolants for uniquely solvable interpolation problems.

{\color{r}
\subsubsection{Canonical Tensor Decomposition}
Evaluation of \eqref{nonlinF} in the finite-dimensional 
function space $D_m=\textrm{span}\{\varphi_1,...,\varphi_{m+1}\}$ 
yields the multivariate function
\begin{equation}
f(a_1,...,a_m) = 5+\sum_{p=1}^{m+1} k_p a_p+ \sum_{q,p=1}^{m+1} 
q_{pq} a_pa_q, 
\label{fdimQ}
\end{equation}
where $a_j=(\varphi_j,\theta)$ and
\begin{equation}
k_p= \int_{0}^{2\pi} K_1(x)\varphi_p(x)dx,\qquad
q_{pq}= \int_{0}^{2\pi}\int_{0}^{2\pi} K_2(x,y)\varphi_p(x)\varphi_q(x)dxdy.
\label{cfoQ}
\end{equation}
The forcing term at the right hand side of 
\eqref{SSYSTEM} can be written as
\begin{equation}
f^n_{jh}=\int_{-b}^b\cdots\int_{-b}^b 
\left(5+\sum_{p=1}^{m+1} k_p a_p+ \sum_{q,p=1}^{m+1} 
q_{pq} a_pa_q\right)
\phi_h(a_j)\prod_{\substack{k=1\\k\neq j}}^{m+1} G^n_k(a_k)da_1\cdots da_m,
\label{fn4}
\end{equation}
As before, the integrals \eqref{fn4} can be reduced to 
products of one-dimensional integrals. If we use a Legendre basis 
$\phi_h$, then we can represent \eqref{nonlinF} 
exactly as a product of constants,  linear and quadratic 
polynomials. To show this, let us consider the 
quadratic part of the functional \eqref{nonlinF}. 
We have,
\begin{align}
\sum_{l=1}^r G^l_1(a_1)\cdots G^l_{m+1}(a_{m+1})
=\sum_{q,p=1}^{m+1} q_{qp} a_qa_p. 
\end{align}
Given the symmetry of $q_{pq}$, the separation 
rank for such quadratic part is $(m+1)(m+2)/2$. 
More generally, the function \eqref{fdimQ} is separable, 
with separation rank $r=1+(m+1)+(m+1)(m+2)/2$. In fact,
it can be written in the form 
\begin{equation}
f(a_1,...,a_m)=\sum_{l=1}^r \alpha_l G^l_1(a_1)\cdots G^l_m(a_m).  
\label{CTDQ}
\end{equation}
A possible ordering of the series could be the one in Table 
\ref{tab:ordering_functional}.
\begin{table}
\color{r}
\centering
\begin{tabular}{c|ccccc}
$q$ & $\alpha_l$ & $G_1^l$ & $G_2^l$ & $\cdots$ & $G_m^l$\\
\hline\\
$1$ &  $5$   &  $1$ & $1$ & $\cdots$& $1$\\
$2$ &  $k_1$ &  $a_1$ & $1$ & $\cdots$& $1$\\
$3$ &  $k_2$ &  $1$ & $a_2$ & $\cdots$& $1$\\
$\vdots $&  $\vdots $ & $\vdots$ & $\vdots$ &$\vdots$ &$\vdots$ \\
$m+2 $&  $k_{m+1}$ & $1$ & $1$ &$\cdots$ &$a_{m+1}$ \\
$m+3 $&  $q_{11}$ & $a_1^2$ & $1$ &$\cdots$ &$1$ \\
$m+4 $&  $q_{12}+q_{21}$ & $a_1$ & $a_2$ &$\cdots$ &$1$ \\
$\vdots $&  $\vdots $ & $\vdots$ & $\vdots$ &$\vdots$ &$\vdots$ \\
$2m+4 $&  $q_{1(m+1)}+q_{(m+1)1}$ & $a_1$ & $1$ &$\cdots$ &$a_{m+1}$\\
$2m+5 $&  $q_{22}$ & $1$ & $a_2^2$ &$\cdots$ &$1$\\
$\vdots $&  $\vdots $ & $\vdots$ & $\vdots$ &$\vdots$ &$\vdots$ \\
$m+2+(m+1)(m+2)/2 $&  $q_{(m+1)(m+1)}$ & $1$ & $1$ &$\cdots$ &$a_{m+1}^2$\\
\end{tabular}
\caption{Ordering of the terms in the canonical tensor 
expansion \eqref{CTDQ} of the quadratic 
functional \eqref{nonlinF}.}
\label{tab:ordering_functional}
\end{table}
In Figure \ref{fig:HT quadratic} we study the accuracy 
of the canonical tensor decomposition  
in representing the quadratic functional \eqref{nonlinF}. Specifically,
we plot the relative pointwise error 
\begin{equation}
e_m=\left|\frac{F([\theta])-f(a_1,...,a_m)}{F([\theta])}\right|\qquad 
a_j=(\theta,\varphi_j)
\label{ptwer}
\end{equation}
at  $\theta(x)=\sin(\cos(2x))+\sin(4x)$
versus the number of basis functions $m$ for 
different separation ranks.
}

{\color{r}
\subsubsection{Hierarchical Tucker Expansion}
The hierarchical Tucker expansion aims at mitigating 
the dimensionality of  the core-tensor of multivariate 
Schmidt decompositions (Section \ref{sec:HT}). It has 
advantages over the canonical tensor decomposition in terms 
of robustness and computational efficiency. In Figure \ref{fig:HT quadratic} 
we study convergence of the hierarchical Tucker 
expansion (see Section \ref{sec:HT}) in representing 
the quadratic functional \eqref{nonlinF}. Specifically,
we plot the relative pointwise error \eqref{ptwer} 
at $\theta(x)=\sin(\cos(2x))+\sin(4x)$  
versus the number of dimensions $m$ for different 
separation ranks.

\begin{figure}
\centerline{\hspace{-0.5cm}Canonical Tensor Decomposition\hspace{3.9cm} Hierarchical Tucker}
\centerline{
\includegraphics[height=6cm]{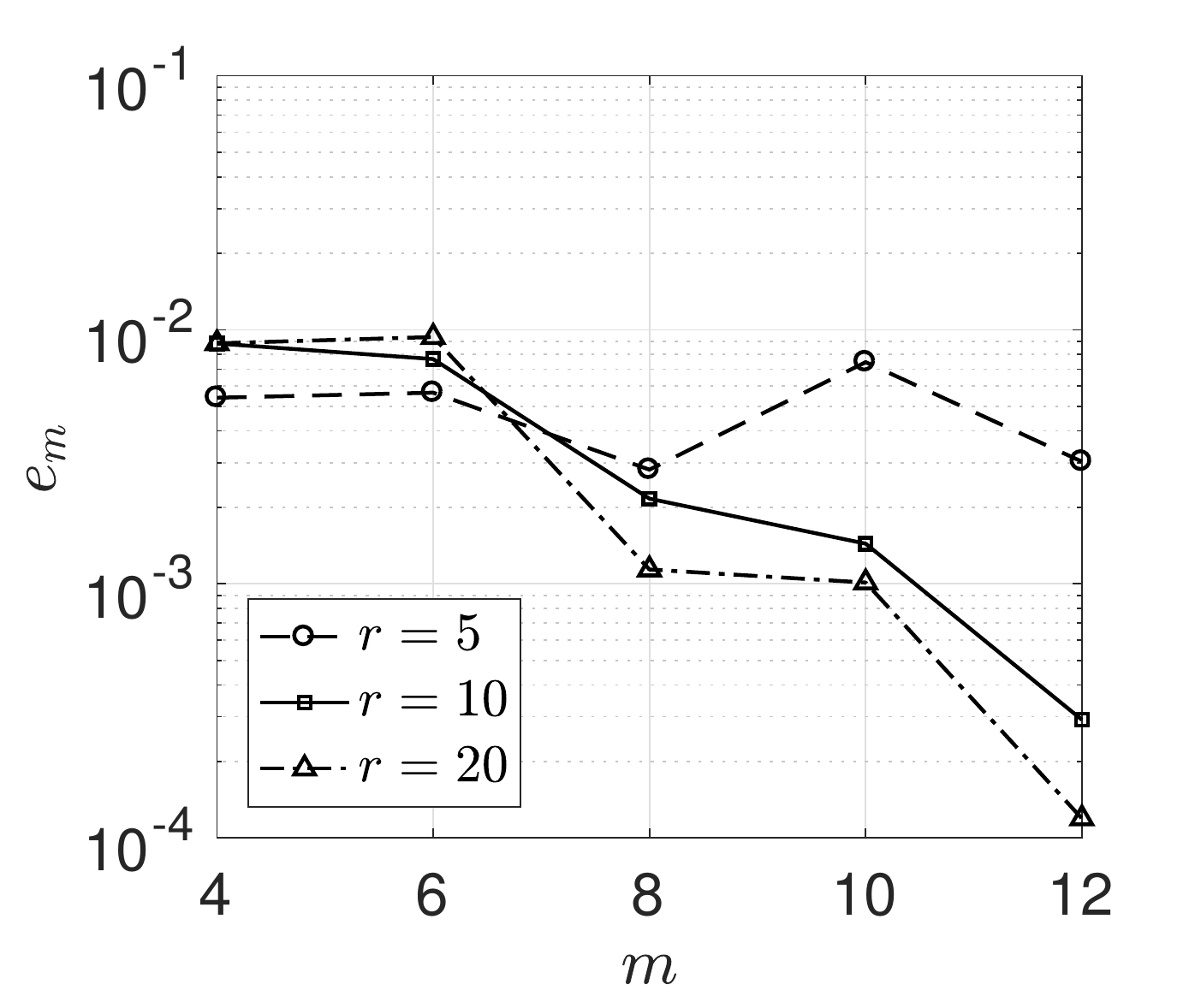}\hspace{0.8cm}
\includegraphics[height=6cm]{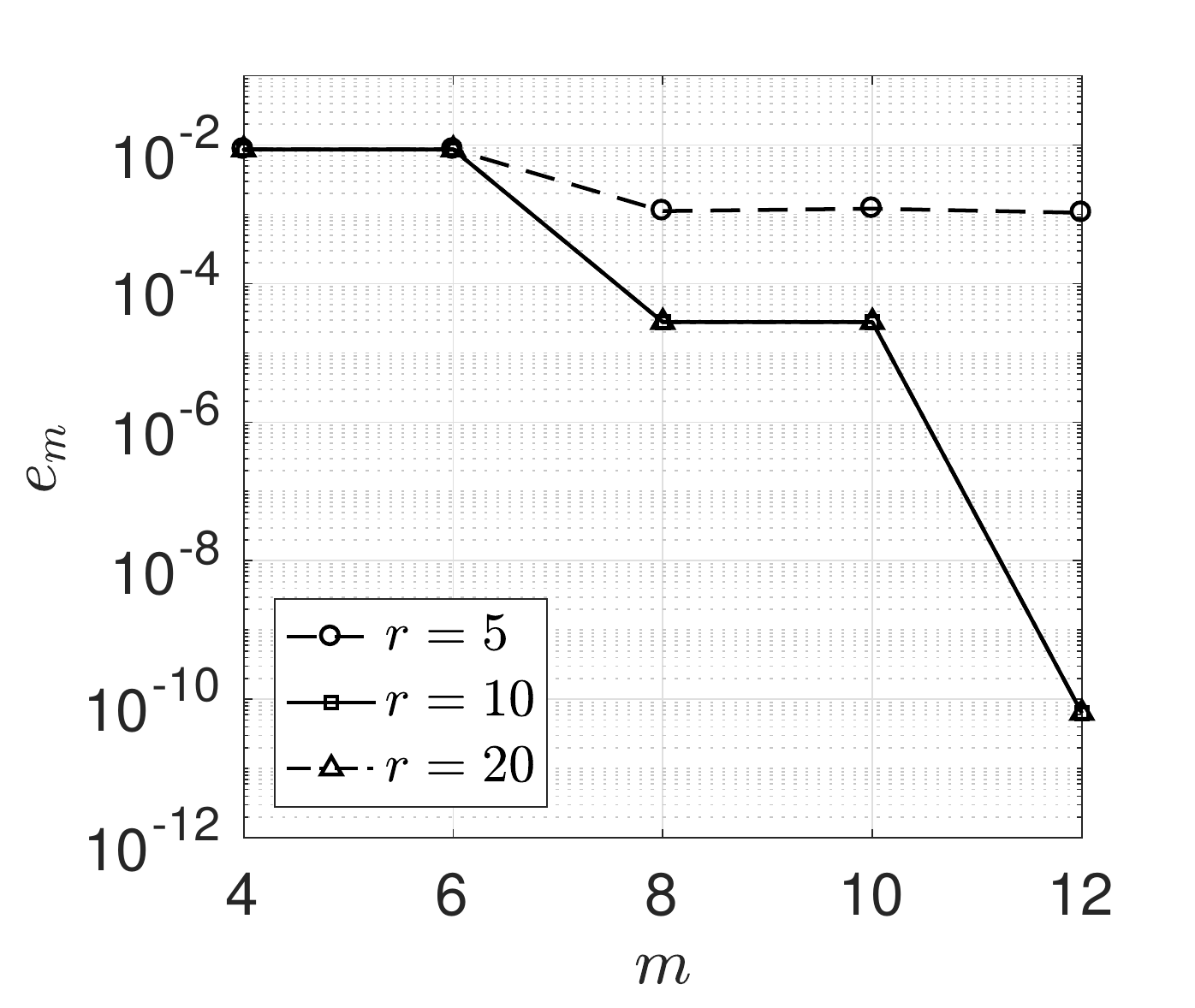}
}
\caption{\color{r}Accuracy of canonical tensor decomposition and 
hierarchical Tucker expansions in representing the quadratic 
functional \eqref{nonlinF}. Specifically, we plot the relative  
pointwise error at $\theta(x)=\sin(\cos(2x))+\sin(4x)$ versus the 
number of dimensions $m$, and for different separation ranks.}
\label{fig:HT quadratic}
\end{figure}
}

\subsection{Hopf Characteristic Functionals}
\begin{figure}[t]
\centerline{\footnotesize\hspace{0.5cm}$q=5$\hspace{4.5cm}$q=20$\hspace{4.5cm}$q=60$}
\centerline{\includegraphics[height=4.5cm]{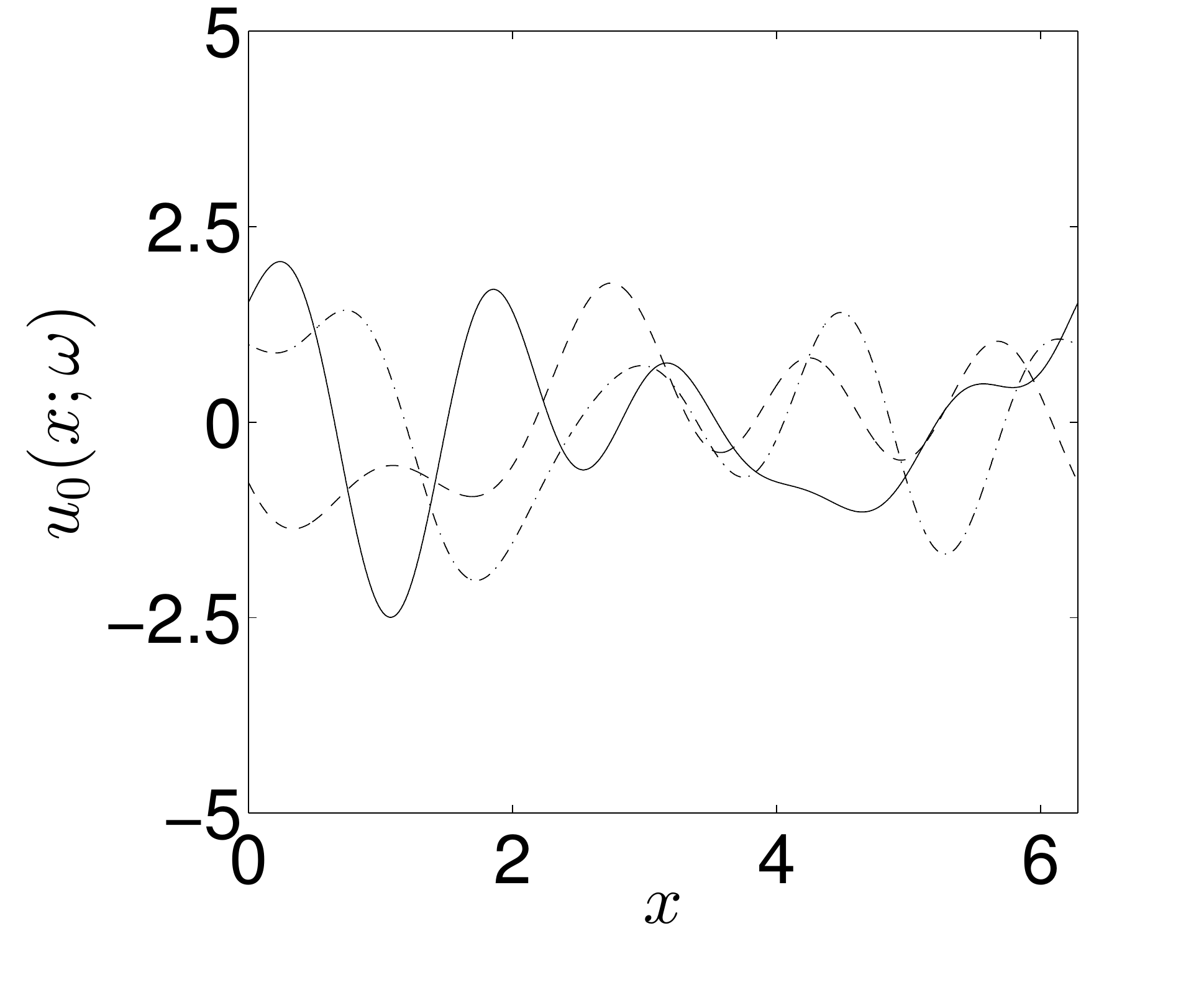}
            \includegraphics[height=4.5cm]{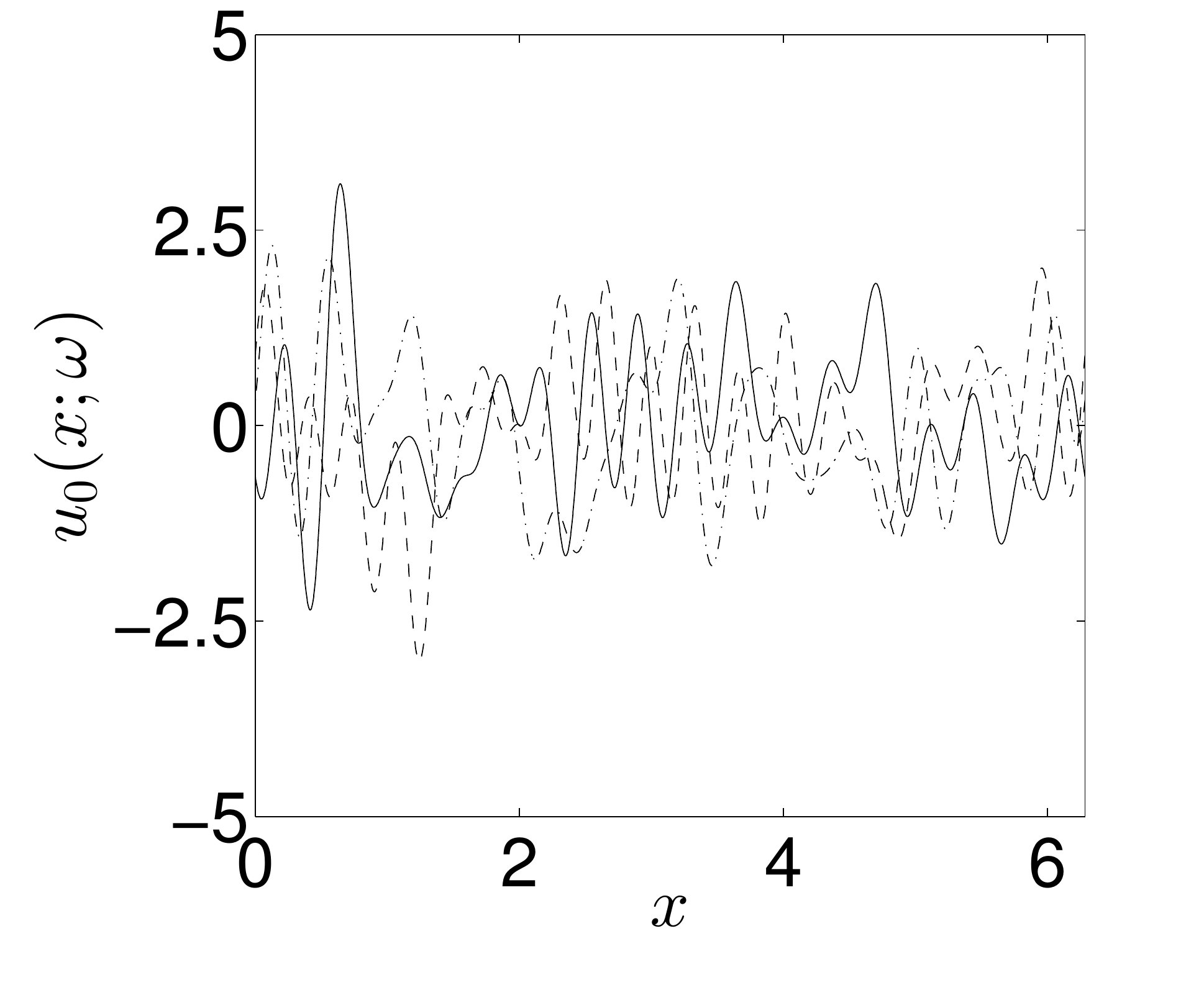}
            \includegraphics[height=4.5cm]{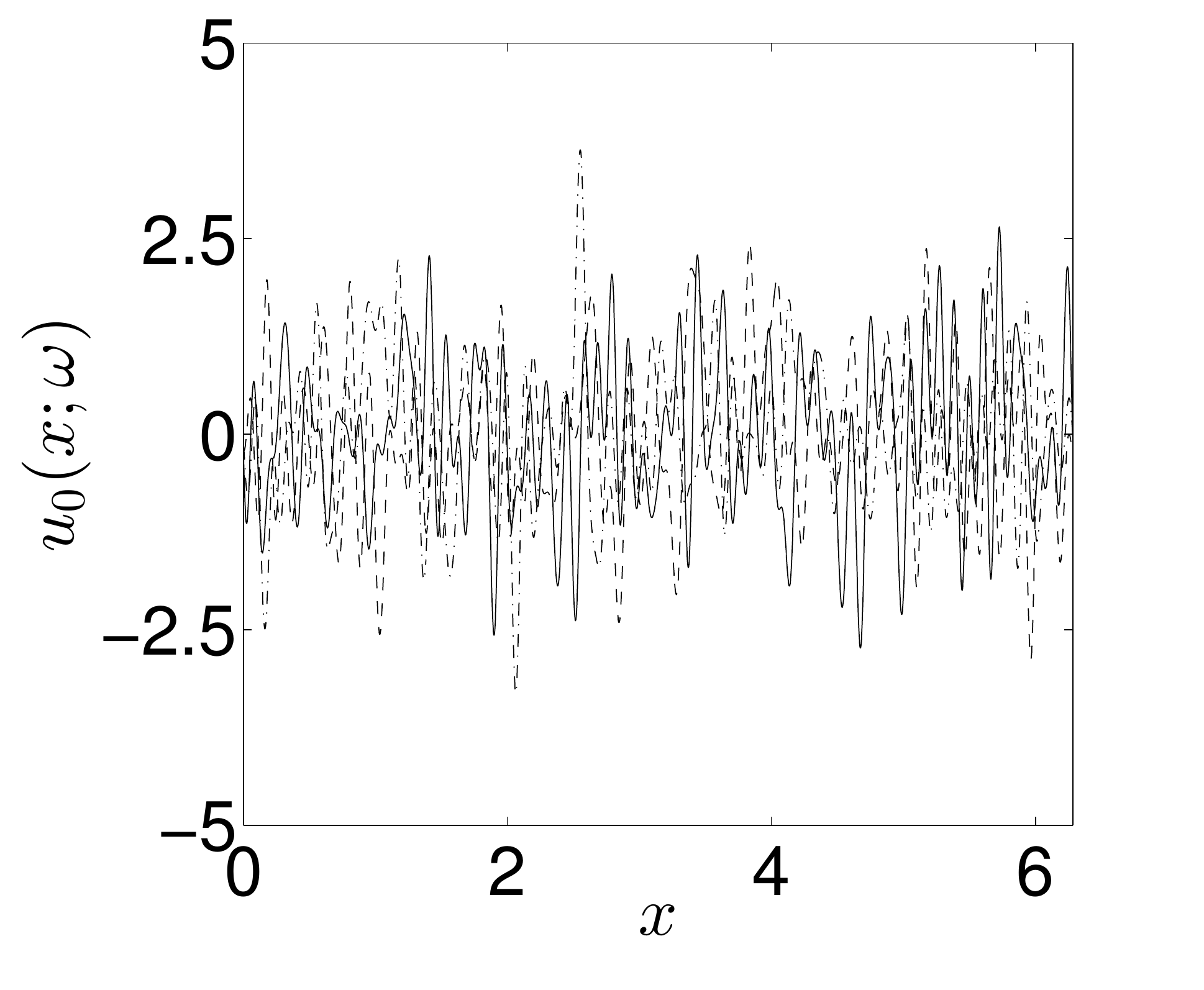}
}
\caption{Samples of the random function \eqref{Burgers_initial} for 
different values of $q$ (Gaussian $\eta_k$ and $\xi_k$).}
\label{fig:samples_u0}
\end{figure}

\begin{figure}
\centerline{\footnotesize$q=1$\hspace{5cm}$q=5$\hspace{5cm}$q=20$}
\centerline{\includegraphics[height=4.5cm]{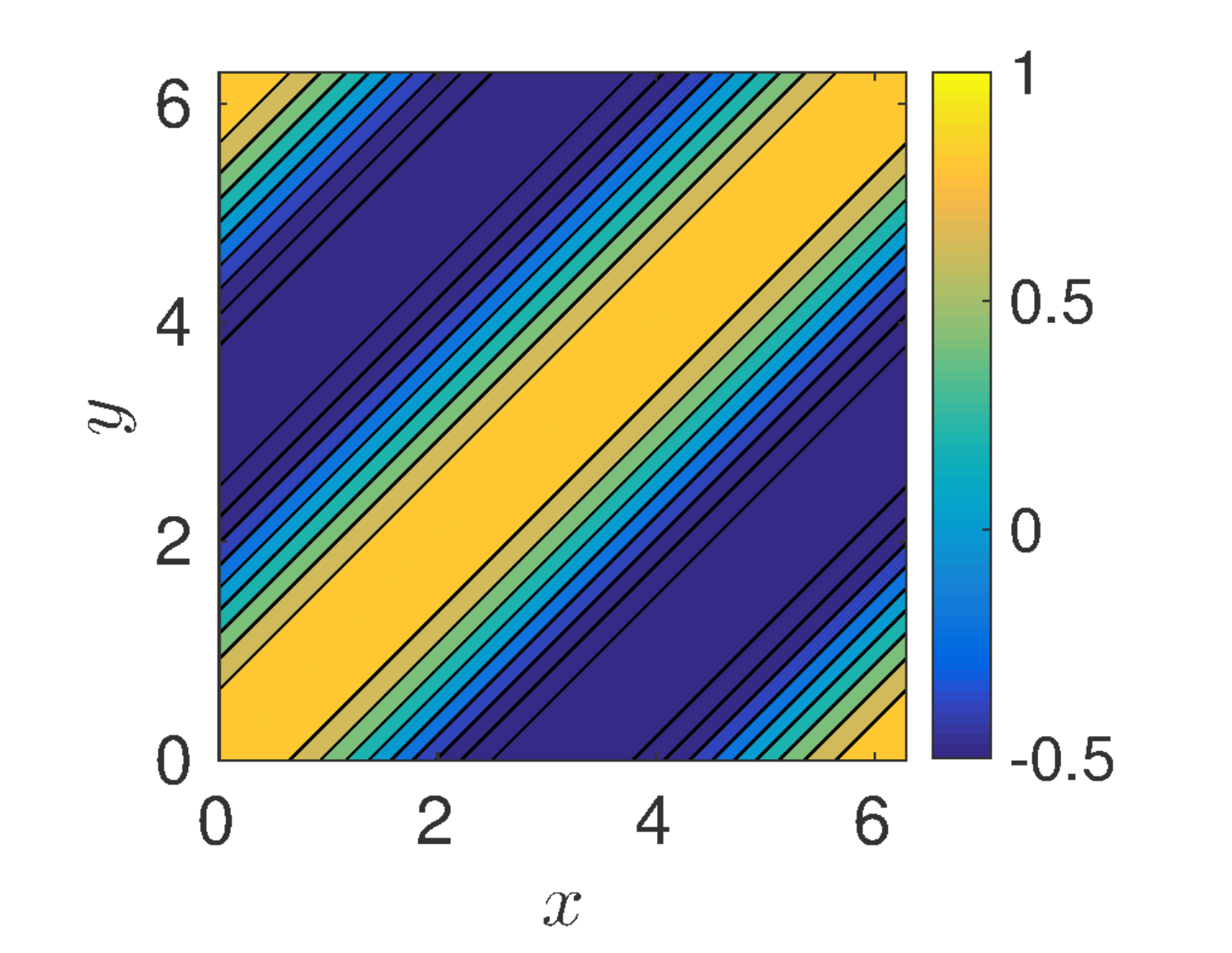}
            \includegraphics[height=4.5cm]{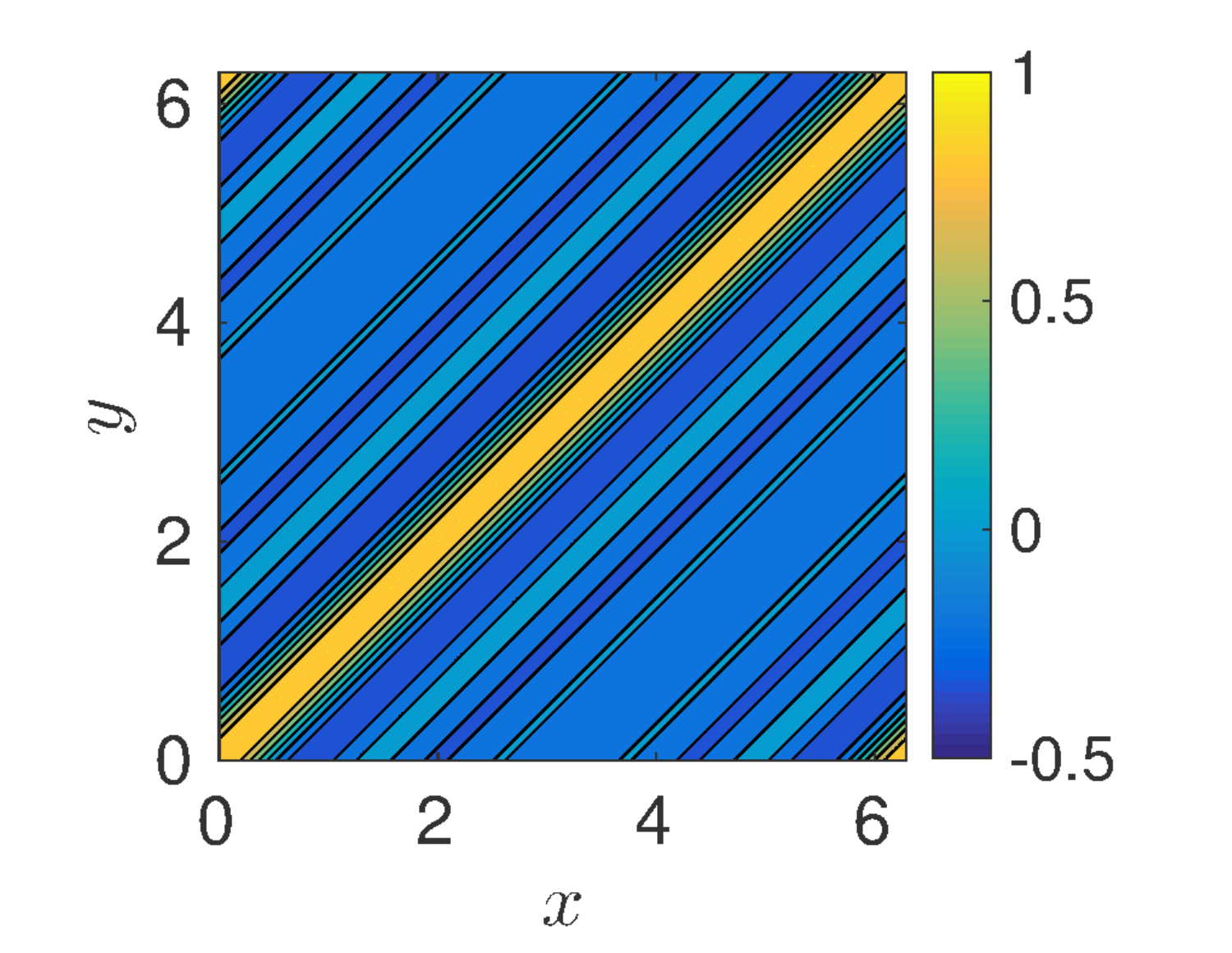}
            \includegraphics[height=4.5cm]{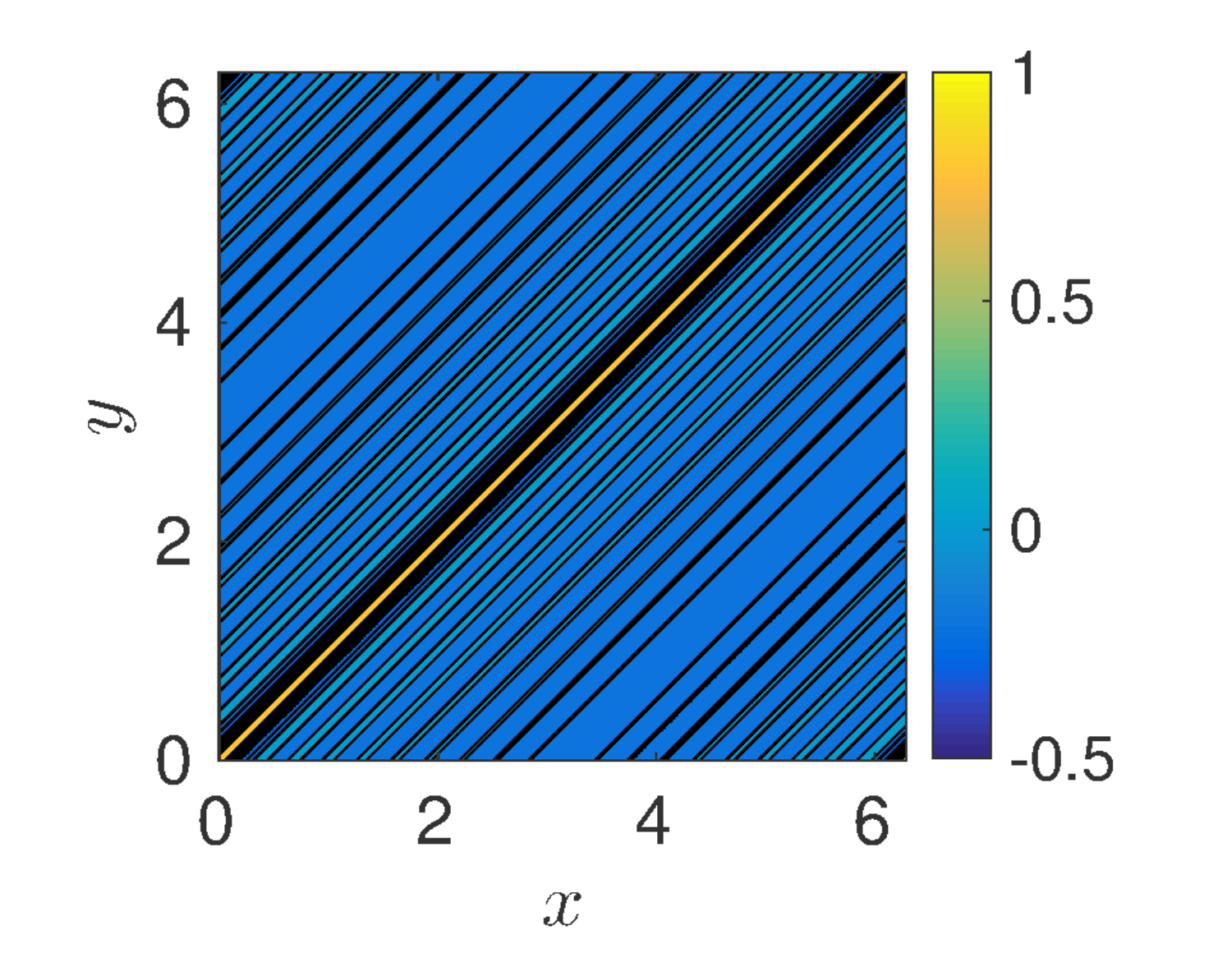}
}
\caption{Covariance function \eqref{testcorrelation} for 
different values of $q$.}
\label{fig:covariance}
\end{figure}

In this Section we study the accuracy of polynomial 
functional interpolation and tensor methods 
in representing Hopf functionals. To this end, consider 
the random function
\begin{equation}
u_0(x;\omega)=\frac{1}{\sqrt{q}}\sum_{k=1}^q\left[\eta_{k}\sin(kx)+\xi_{k}\cos(kx)\right],
\label{Burgers_initial}
\end{equation}
where $\{\eta_k,\xi_k\}_{k=1,...,q}$ are  i.i.d. random variables satisfying
\begin{equation}
 \langle\xi_{k}\rangle=0,\quad \langle\eta_{k}\rangle=0,\quad \langle\xi_{k}^2\rangle=1,
 \quad \langle\eta_{k}^2\rangle=1.
 \label{assumptions}
\end{equation}
In Figure \ref{fig:samples_u0} we plot few samples of \eqref{Burgers_initial}, 
for different values of $q$ in the hypothesis that $\{\xi_k\}$ and 
$\{\eta_k\}$ are i.i.d. Gaussian random variables.
From \eqref{assumptions} it follows that the first two statistical 
moments of \eqref{Burgers_initial} are independent of $q$, i.e., we have
\begin{equation}
\left<u_0(x;\omega)\right>=0, \qquad  
\left<u_0(x;\omega)^2\right>=1,\quad 
\textrm{for all $q$.}
\end{equation}
On the other hand, the covariance function
\begin{equation}
 C_0(x,y)=\frac{1}{q}\sum_{k=1}^q[\sin(kx)\sin(ky)
 +\cos(kx)\cos(ky)]
 \label{testcorrelation}
\end{equation}
does depend on $q$, as shown in Figure \ref{fig:covariance}.
Higher-order moments and cumulants can be computed analytically in a similar way. 
The Hopf functional of $u_0(x;\omega)$ is defined as
\begin{equation}
 \Phi([\theta(x)])=\left<\exp\left[i\int_{0}^{2\pi}
 \theta(x)u_0(x;\omega)dx\right]\right>,
 \label{H_initial}
\end{equation}
where the average $\left<\cdot\right>$ is a multi-dimensional integral 
with respect to the joint probability density function of $\xi_k$ and $\eta_k$.
By substituting \eqref{Burgers_initial} 
into \eqref{H_initial} we obtain  
\begin{equation}
 \Phi([\theta(x)])=
 \prod_{k=1}^q\int e^{ias_k[\theta]}p_{\eta_k}(a)da 
 \prod_{k=1}^q\int e^{iac_k[\theta]}p_{\xi_k}(a)da, \qquad 
 \label{H_initial1}
\end{equation}
where 
\begin{equation}
s_k[\theta]=\frac{1}{\sqrt{q}}\int_{0}^{2\pi}\theta(x)\sin(kx)dx,\qquad 
c_k[\theta]=\frac{1}{\sqrt{q}}\int_{0}^{2\pi}\theta(x)\cos(kx)dx.
\label{coeffskck}
\end{equation}
Depending on the probability density functions 
$p_{\eta_k}(a)$ and $p_{\xi_k}(a)$ appearing in \eqref{H_initial1}, 
we have different expressions of $F([\theta])$. 
\paragraph{Gaussian Random Fields} Let us assume 
 \begin{equation}
  p_{\eta_k}(a)=\frac{1}{\sqrt{2\pi}}e^{-a^2/2}\qquad 
  p_{\xi_k}(a)= p_{\eta_k}(a).
 \end{equation}
The integrals in  \eqref{H_initial1} are in the form
 \begin{equation}
\frac{1}{\sqrt{2\pi}}\int_{-\infty}^{\infty} e^{ias_k[\theta]-a^2/2}da=
e^{-s^2_k[\theta]/2}.
 \end{equation}
Therefore, we obtain 
\begin{align}
\Phi([\theta(x)])=&\exp\left[-\frac{1}{2}\sum_{k=1}^q
\left(s^2_k[\theta]+c^2_k[\theta]\right)\right]\label{H_initial3}\\
=&\exp\left[-\frac{1}{2} \int_{0}^{2\pi}\int_{0}^{2\pi}\theta(x)\theta(y)
C_0(x,y)dxdy\right],\label{FDG}
\end{align} 
in agreement with well-known results for Hopf characteristic functionals 
of Gaussian random fields.

\paragraph{Uniform Random Fields} Let us assume 
 \begin{equation}
  p_{\eta_k}(a)=
  \begin{cases}
1/(2\sqrt{3})& a\in[-\sqrt{3},\sqrt{3}]\\
0            & \textrm{otherwise}   
  \end{cases}\qquad   p_{\xi_k}(a)= p_{\eta_k}(a)
 \end{equation}
In this way the assumptions \eqref{assumptions} are satisfied.
The integrals in \eqref{H_initial1} are easily obtained as
 \begin{equation}
\frac{1}{2\sqrt{3}}\int_{-\sqrt{3}}^{\sqrt{3}} e^{ias_k[\theta]}da=
\frac{\sinh(i\sqrt{3}s_k[\theta])}{i\sqrt{3}s_k[\theta]}=
\frac{\sin(\sqrt{3}s_k[\theta])}{\sqrt{3}s_k[\theta]}.
 \end{equation}
A substitution of this formula into \eqref{H_initial1} yields
\begin{align}
 \Phi([\theta(x)])
% =& -\prod_{k=1}^m\frac{\sinh(i\sqrt{3}s_k[\theta])\sinh(i\sqrt{3}c_k[\theta])}
%  {3s_k[\theta]c_k[\theta]}
%  \label{H_initial2}\\
 =& \prod_{k=1}^q\frac{\sin(\sqrt{3}s_k[\theta])\sin(\sqrt{3}c_k[\theta])}
 {3s_k[\theta]c_k[\theta]}.
 \label{H_initial2}
\end{align}

\vspace{0.5cm}
\noindent
Note that in both cases we just discussed
the Hopf functional turns out to be 
real-valued. Moreover, $\Phi\rightarrow 0$ as 
$\left\|\theta\right\|\rightarrow\infty$ (Riemann-Lebesgue lemma), 
at a rate that depends on the regularity of the underlying 
probability density functional. In particular, in 
the Gaussian case $\Phi$ goes to zero faster than in 
the uniform case.  
Note also, that both \eqref{FDG} are \eqref{H_initial2} are 
{\em entire functionals} (i.e., analytic on the complex plane). 
This implies that the polynomial interpolation 
process converges {\em pointwise} \cite{Khlobystov2}.

\subsubsection{Effective Dimension} 
\label{sec:Hopf_effective_dimensionality}
Let us represent 
$\theta(x)$ in the space of periodic functions in $[0,2\pi]$.
Possible bases are the discrete trigonometric polynomials
\eqref{setbf} or the more classical Fourier modes
\begin{equation}
1,\qquad \sin(kx),\qquad \cos(kx), \qquad k=1,2, ... . \label{FB}
\end{equation}
Let us now ask the following question: 
what is the effective dimension of the Hopf 
functional \eqref{H_initial} in the space of periodic 
functions? 
Such dimension is clearly determined by the dimension of 
the linear functionals $s_k$ and $c_k$ in \eqref{coeffskck}.
If we expand the test function $\theta$ in a 
classical Fourier series 
\begin{equation}
\theta(x)= a_0 + \sum_{k=1}^N a_k\sin(kx)+ \sum_{k=1}^N b_k\cos(kx)
\end{equation}
and we substitute it into \eqref{coeffskck} then we obtain
\begin{equation}
s_k[\theta]=\frac{\pi}{\sqrt{q}} a_k, \qquad\textrm{and}\qquad  c_k[\theta]=
\frac{\pi}{\sqrt{q}} b_k.   
\end{equation}
This means that the {\em effective dimension} of the Hopf functional 
\eqref{H_initial} is {\em exactly} $2q$. 
To show this numerically, we consider the Hopf functional 
\eqref{FDG} (Gaussian case) and plot the spectrum of the covariance matrix 
\begin{equation}
Z_{ij} = \int_0^{2\pi}\int_0^{2\pi}C_0(x,y)\varphi_i(x)\varphi_j(y)dxdy,
\label{covariance}
\end{equation}
obtained by projecting the covariance function 
\eqref{testcorrelation} onto the 
the Fourier modes \eqref{FB} or, equivalently, 
onto the discrete trigonometric polynomials \eqref{setbf}. 
\begin{figure}[t]
\vspace{0.5cm}
\centerline{
\includegraphics[height=6cm]{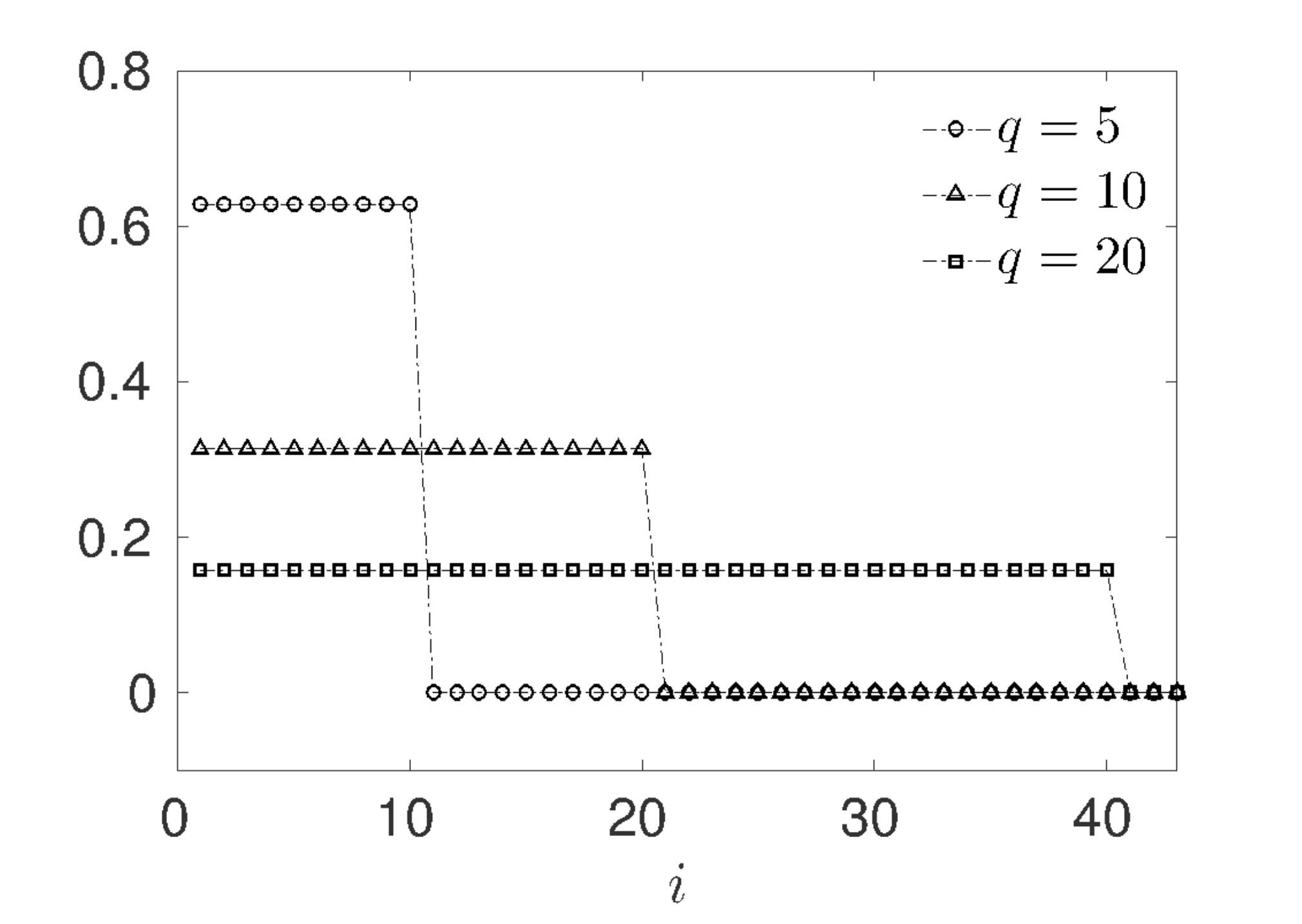}
}

\caption{ 
Spectra of the covariance matrix \eqref{covariance}
obtained by projecting the covariance function 
\eqref{testcorrelation} onto the Fourier modes \eqref{FB}, 
or the discrete trigonometric polynomials \eqref{setbf}. 
It is seen that the number of active components is exactly 
$2q$ and all variables are {\em equally important}. This has 
important consequences when we approximate the Hopf 
functional using tensor methods or polynomial functionals. 
In particular, any approximation in a function space with less 
than $2q$ dimensions yields a systematic error.}
\label{fig:covariancespectra}
\end{figure}
As it is clearly seen from 
Figure \ref{fig:covariancespectra} the number of 
active components is exactly $2q$. Note also that 
the spectrum is {flat} in all cases, which means 
that all active variables are {equally important}. 
This has important consequences when
approximating the functional \eqref{FDG} by 
polynomial functionals or SSE. In particular, if we use 
functionals involving less than $2q$ components, e.g., 
if $m<2q$ in equation \eqref{functional-SSE}, then we cannot 
approximate \eqref{FDG} accurately, no matter how we push the 
expansion order.

\subsubsection{Polynomial Functional Interpolation}
\label{sec:Hopf polynomial}
\begin{figure}[t]
\centerline{\hspace{0.5cm} Gaussian Functional \hspace{4.2cm} Uniform Functional}
\centerline{
            \includegraphics[height=6cm]{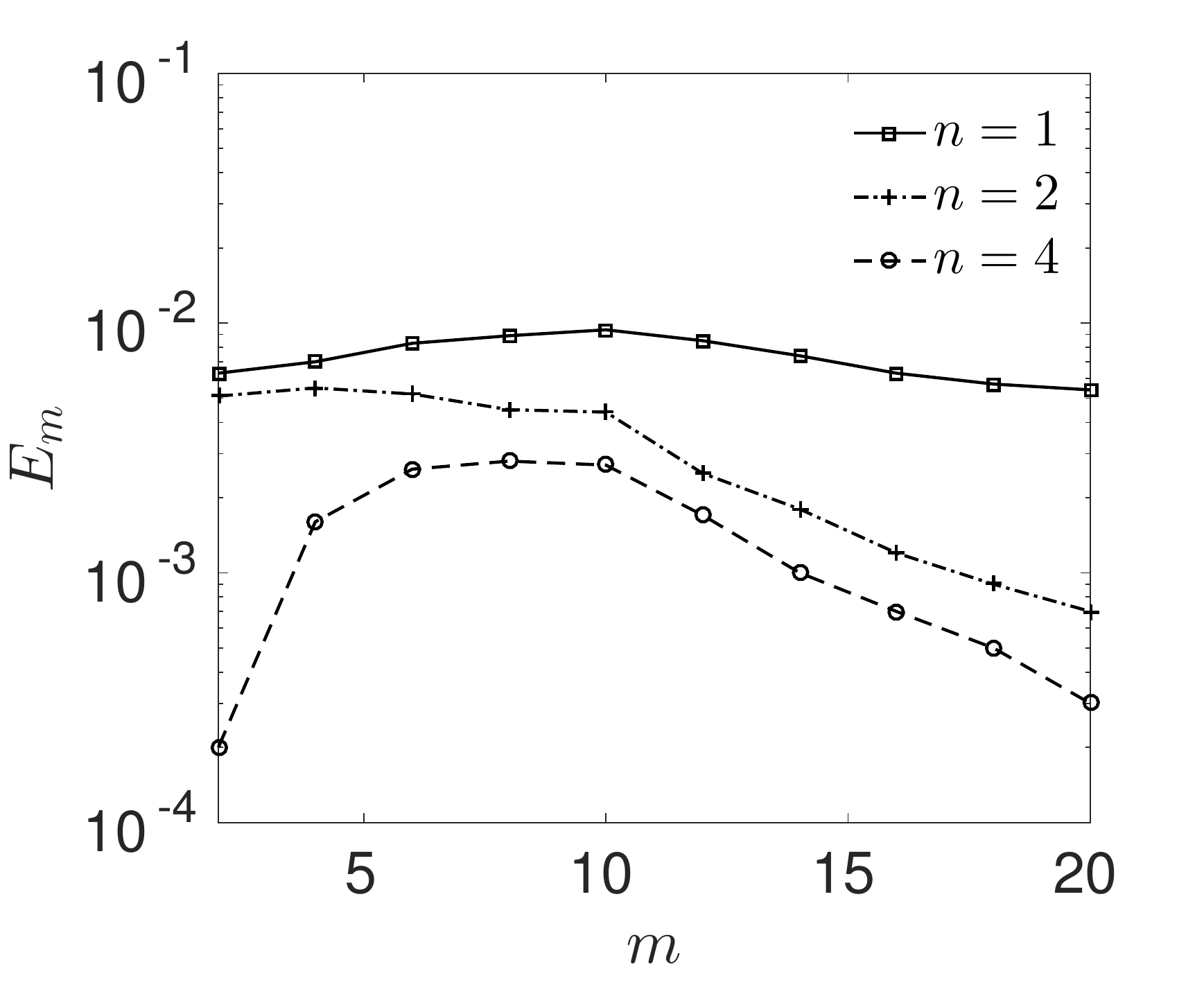}
            \includegraphics[height=6cm]{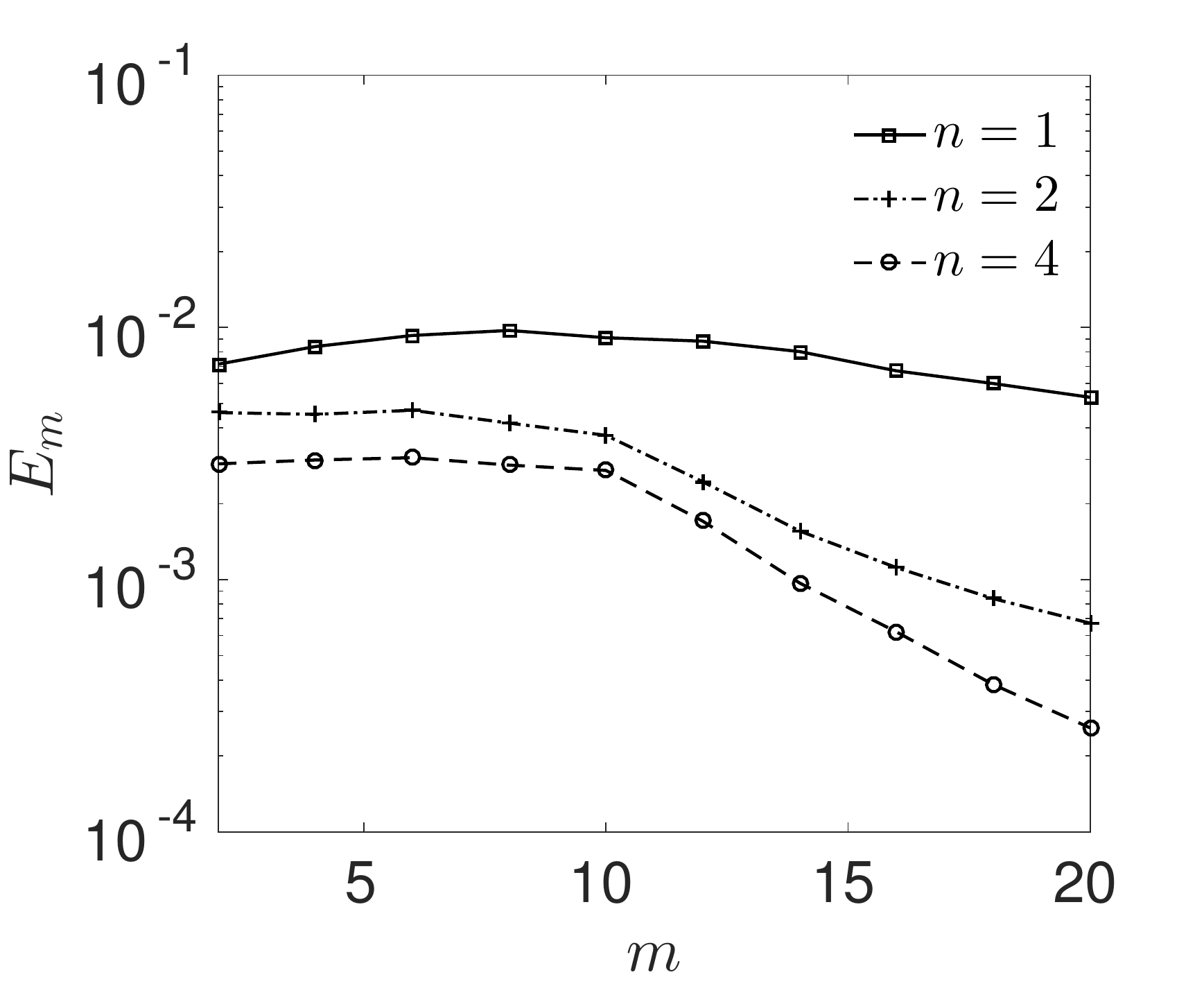}
}

\caption{Accuracy of Porter's polynomial functional interpolants in 
representing  the Hopf functionals \eqref{H_initial3} and \eqref{H_initial2}. 
Shown are the pointwise errors \eqref{Em3} versus $m$ 
for different orders of the polynomial interpolant. 
The errors are computed by sampling $20000$ functions from 
the set $\mathcal{G}_{50}(\sigma)$ with $\sigma=0.01$
and then computing the supremum \eqref{Em3}.}
\label{fig:Hopf_porter}
\end{figure}

\begin{figure}[!ht]

\centerline{\hspace{0.12cm}Nodes in Fourier space \hspace{2.4cm} Nodes in function space}
\centerline{\includegraphics[height=4.9cm]{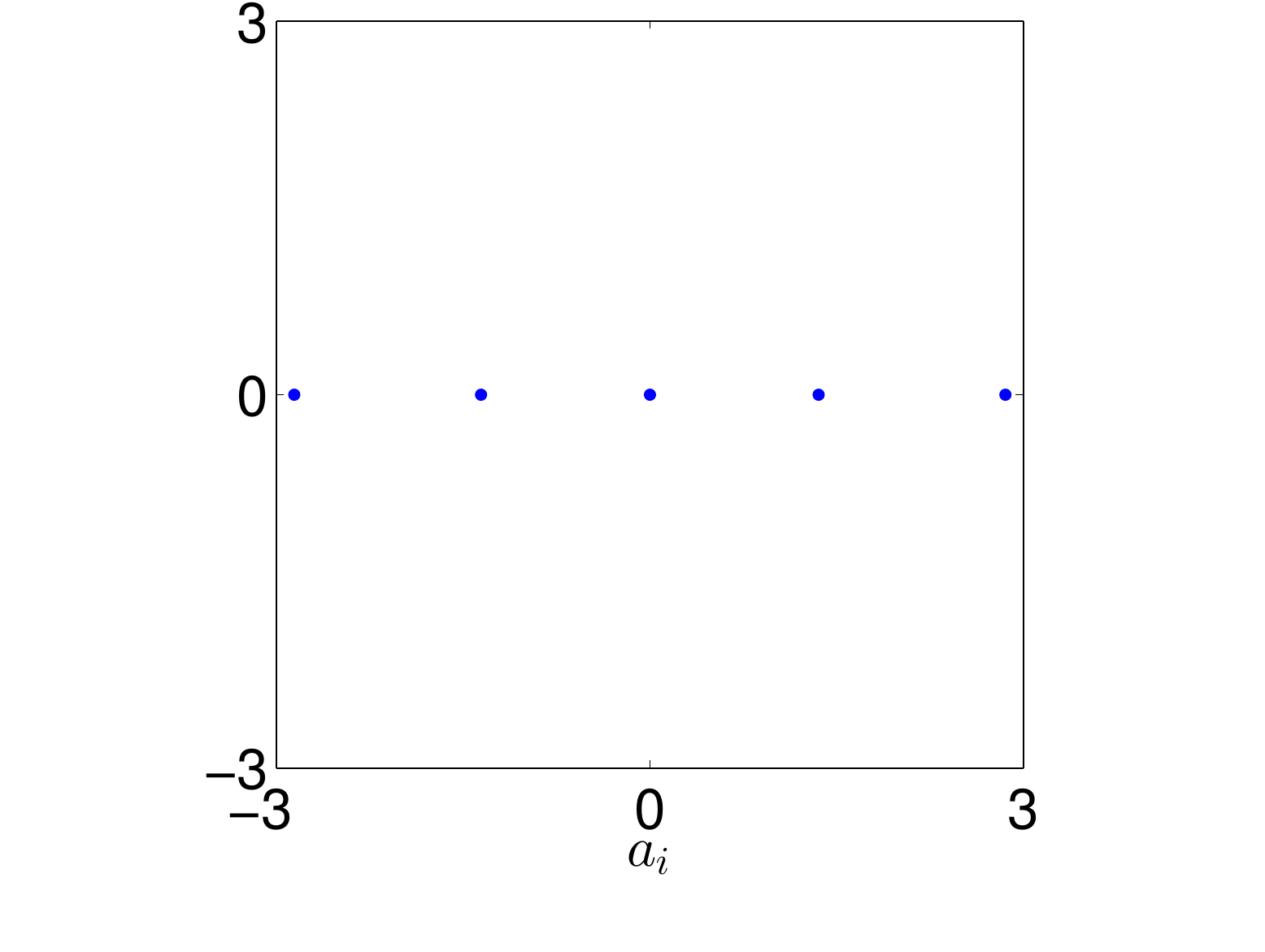}
            \includegraphics[height=4.9cm]{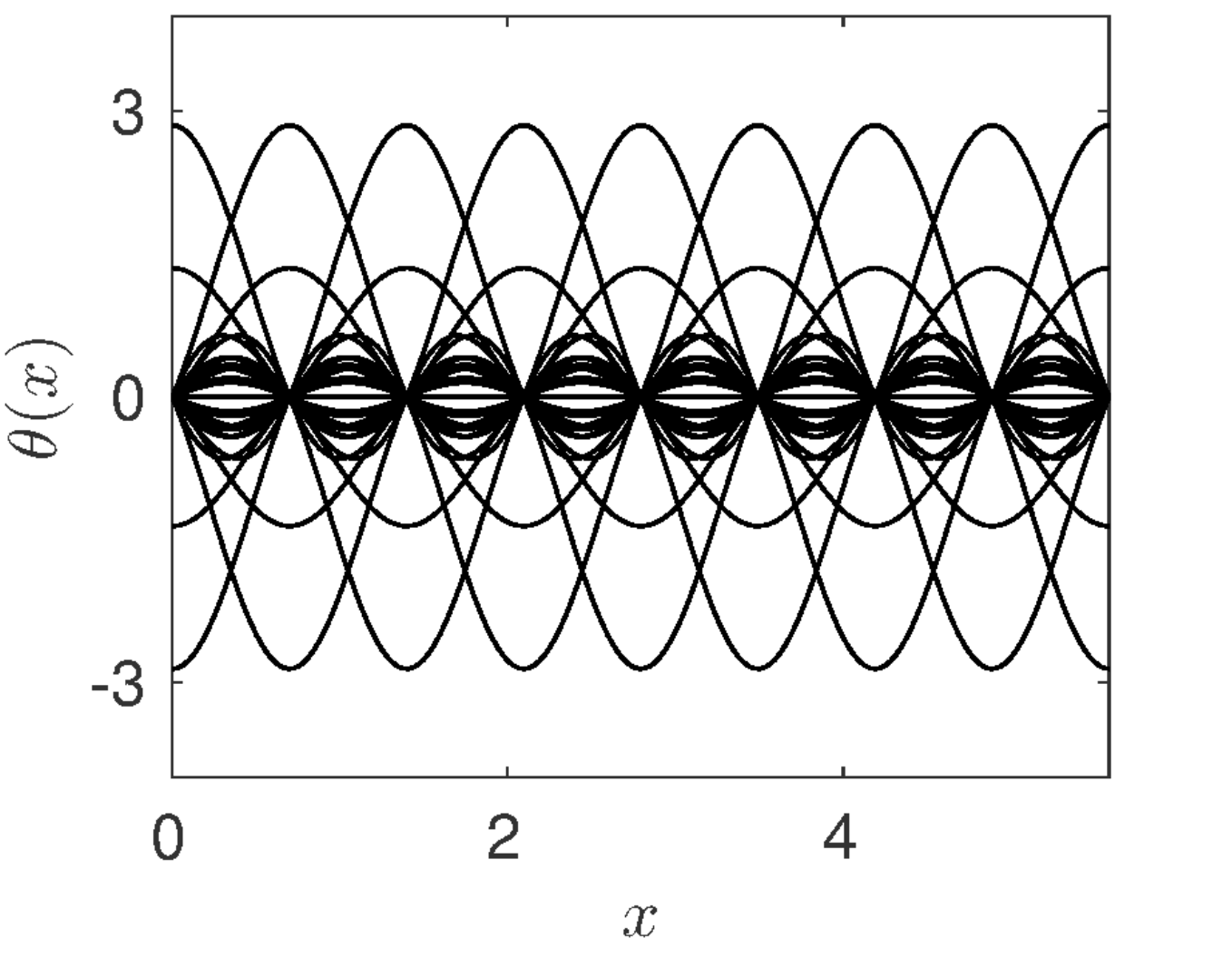}}
\centerline{\includegraphics[height=4.9cm]{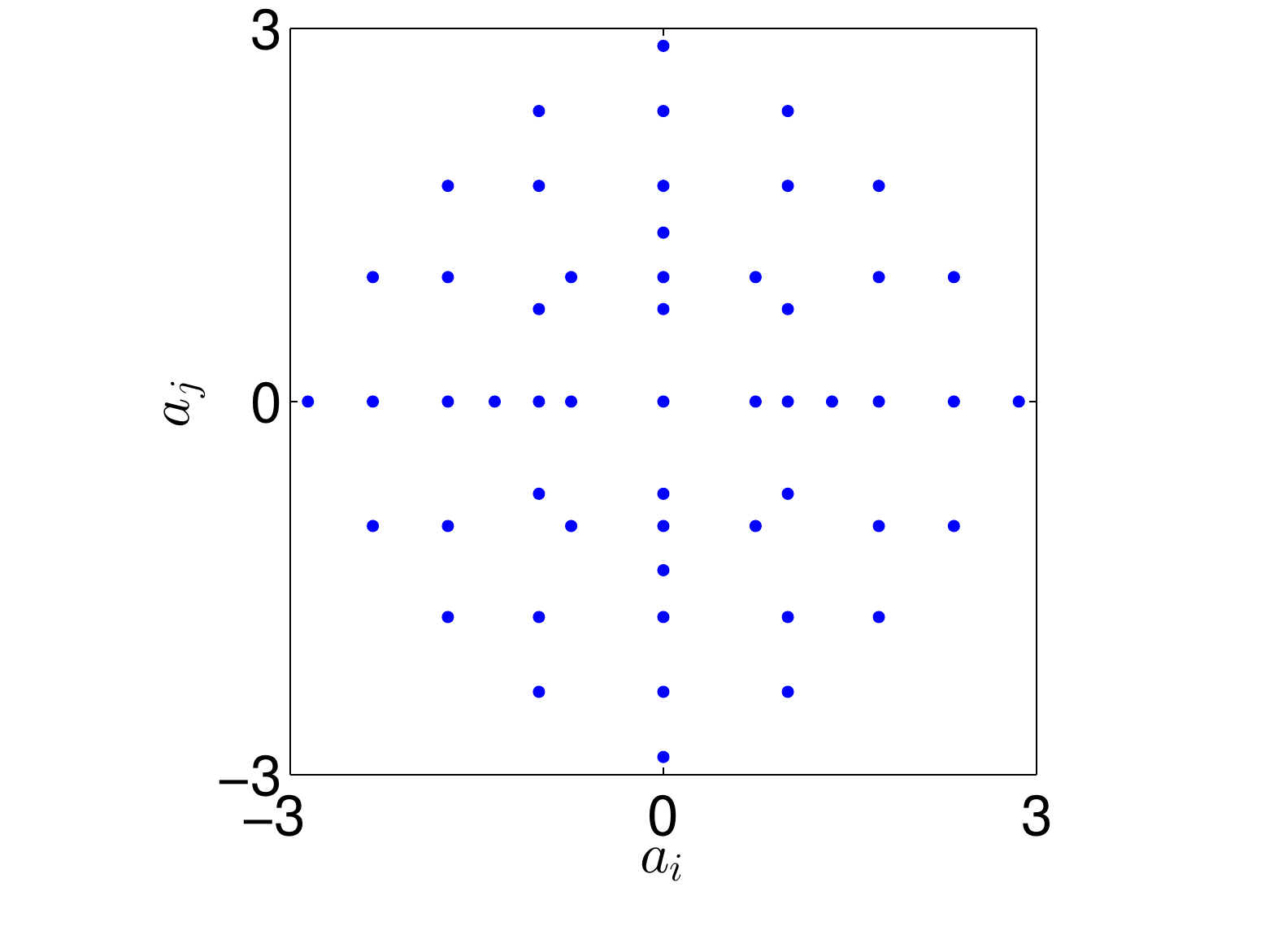}
            \includegraphics[height=4.9cm]{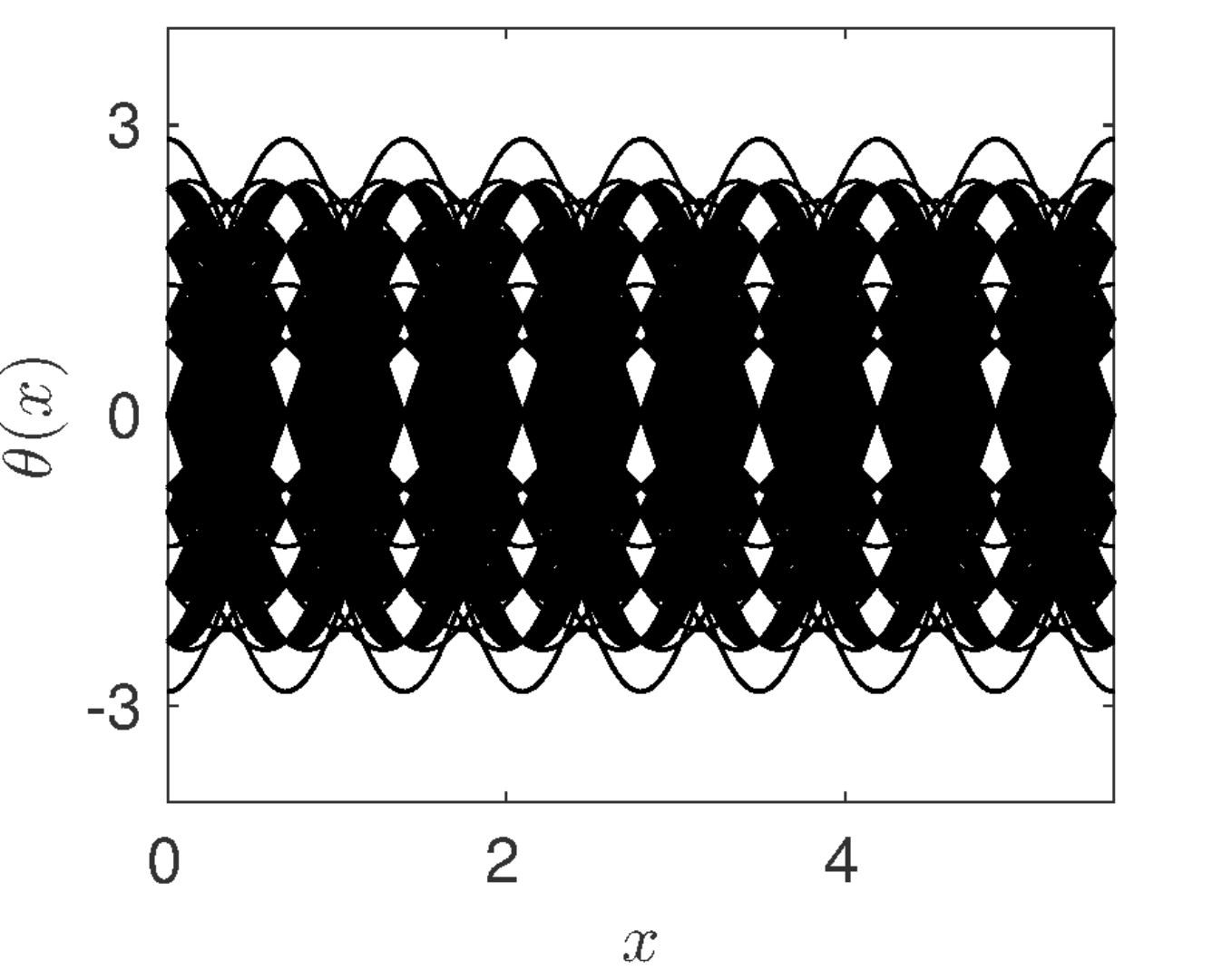}}
\centerline{\includegraphics[height=4.9cm]{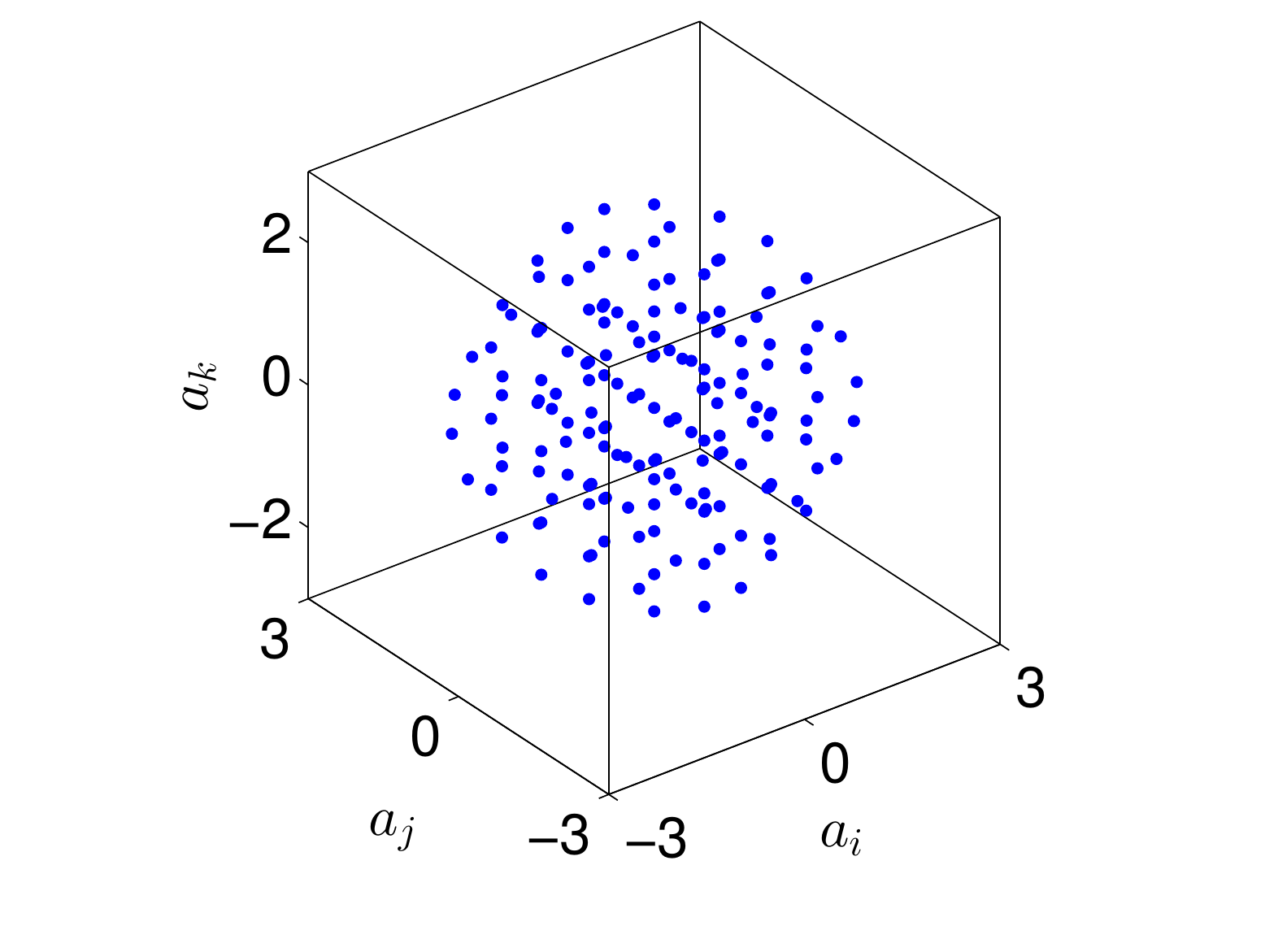}
            \includegraphics[height=4.9cm]{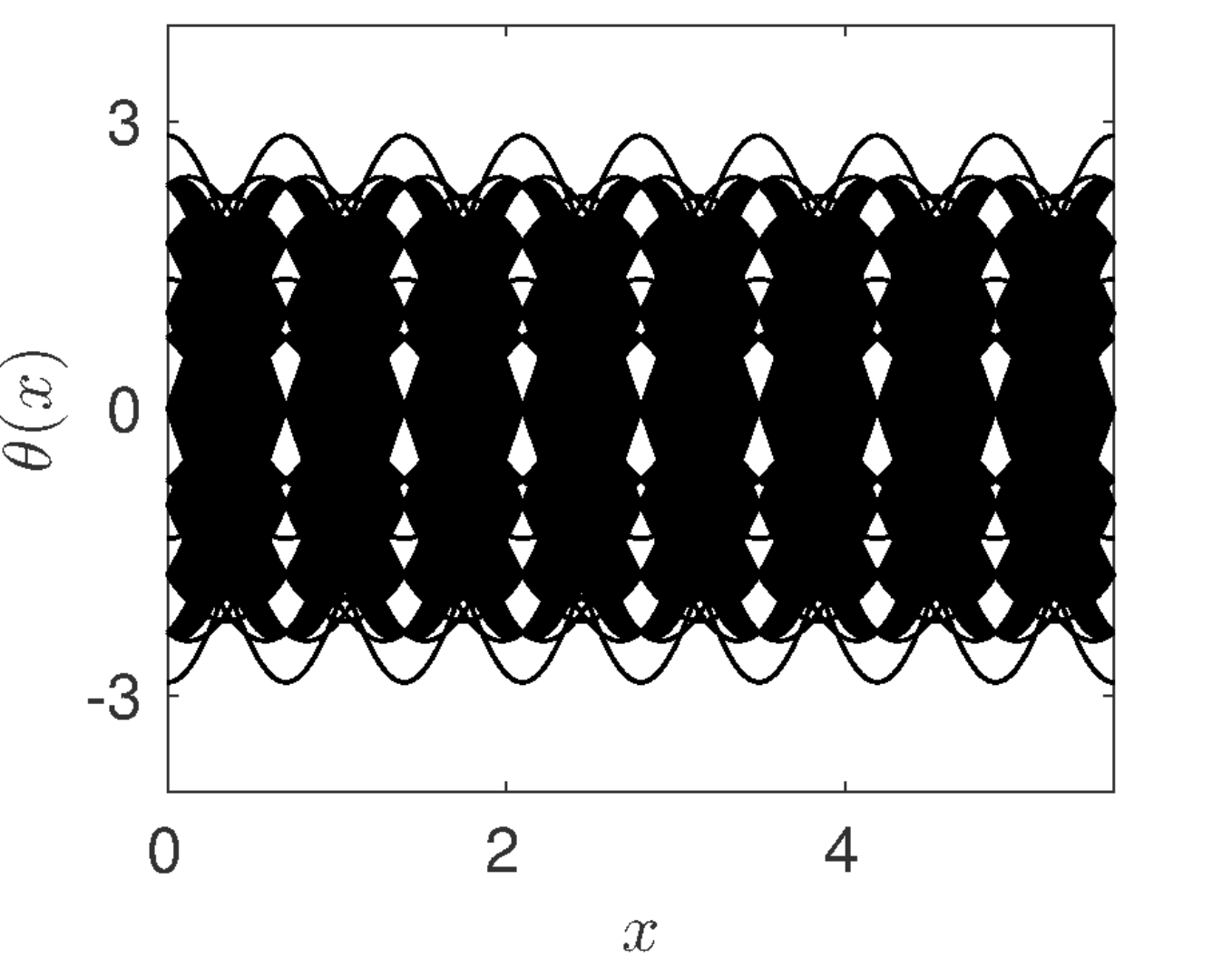}}
 
\caption{Gauss-Hermite sparse grids (left) and corresponding interpolation 
nodes in the function space $\widetilde{S}^{(m)}_n$ (right). 
Specifically, we plot all elements in 
$\widetilde{S}^{(9)}_1$ (first row, $45$ elements), 
$\widetilde{S}^{(9)}_2$ (second row, $1117$ elements) and
$\widetilde{S}^{(9)}_3$ (third row, $3949$ elements).}
\label{fig:sparse_grid_sets}
\end{figure}

\begin{figure}[!t]
\centerline{\hspace{0.5cm} Gaussian Functional \hspace{5cm} Uniform Functional}
\centerline{
\includegraphics[height=6.5cm]{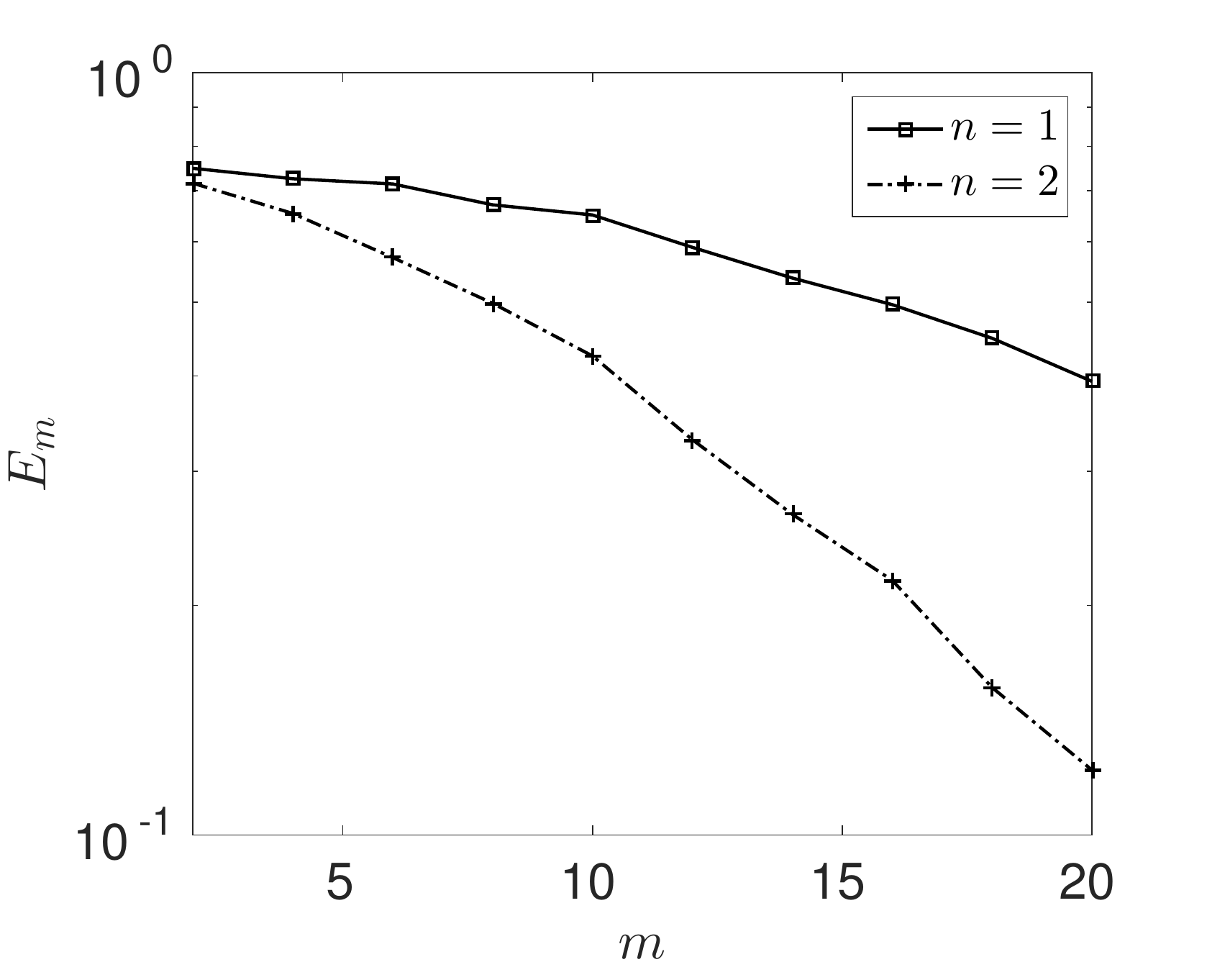}
\includegraphics[height=6.5cm]{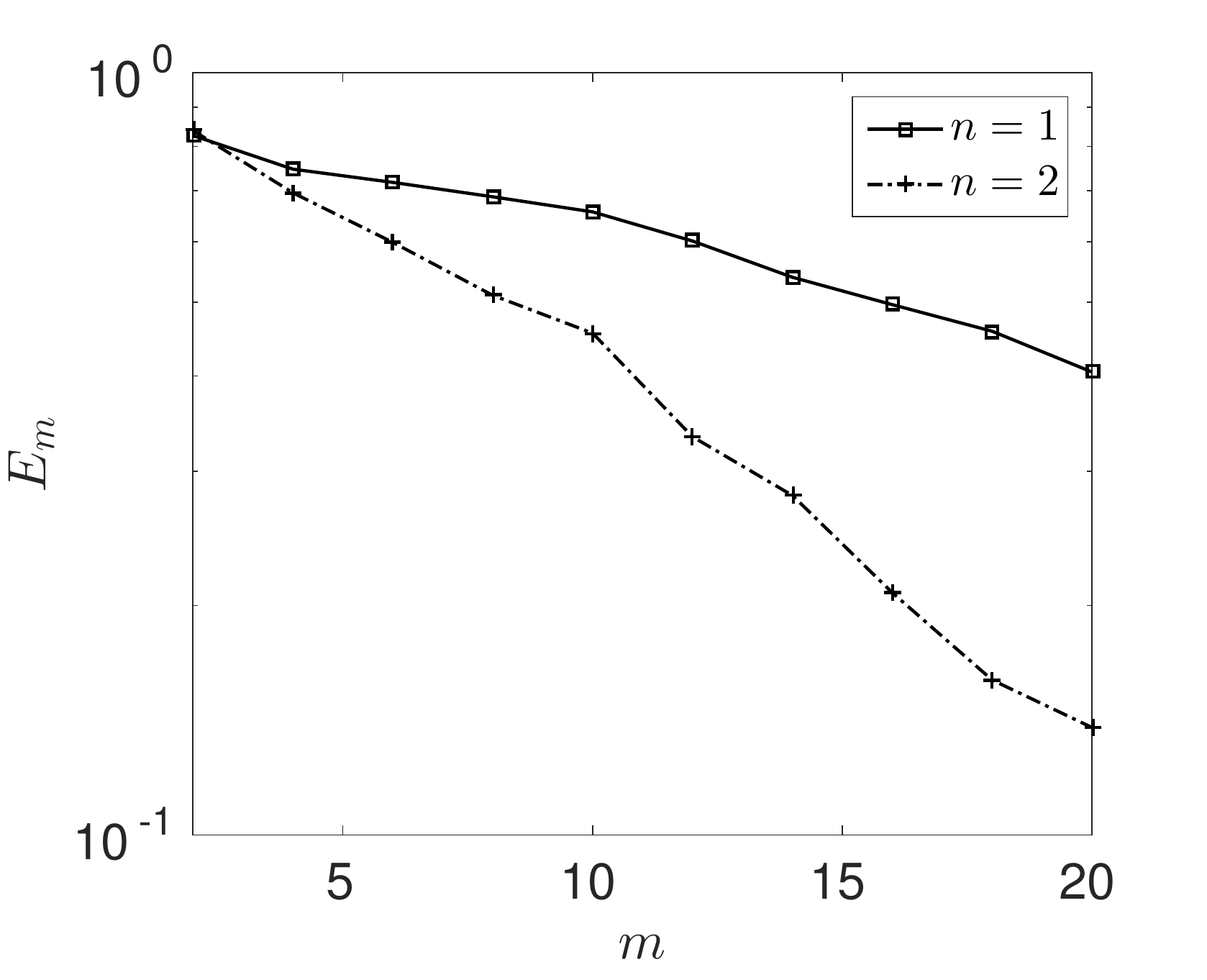}
}

\caption{Approximation of the Hopf functionals 
\eqref{H_initial3} and \eqref{H_initial2} by using 
polynomial functionals of total degree $n=1$ and $n=2$. 
The polynomial are in the form \eqref{Porter_interpolant_noncardinal} 
(i.e., non-interpolatory) and they are constructed by using the set of nodes 
$\widetilde{S}^{(m)}_1$ and $\widetilde{S}^{(m)}_2$ 
defined in equation \eqref{SNq3}, with $a_{i_j}$ sampled at Gauss-Hermite 
sparse grids of level $5$ \cite{Bungartz}. 
Shown are the pointwise errors \eqref{Em3} versus $m$ for functions
$\theta(x)$ in the ensemble $\mathcal{G}_{50}(1)$.}
\label{fig:Hopf_porter_SG}
\end{figure}

Let us expand \eqref{FDG} in a power series
\begin{align}
\Phi([\theta(x)])=&1-\frac{1}{2} \int_{0}^{2\pi}\int_{0}^{2\pi}
C_0(x_1,x_2)\theta(x_1)\theta(x_2)dx_1dx_2+\nonumber\\
&\frac{1}{4}\int_{0}^{2\pi}\int_{0}^{2\pi}\int_{0}^{2\pi}\int_{0}^{2\pi}C_0(x_1,x_2)C_0(x_3,x_4)
\theta(x_1)\theta(x_2)\theta(x_3)\theta(x_4)dx_1dx_2dx_3dx_4+\cdots.
\end{align}
For small $\theta(x)$ we can truncate the series, 
and represent $\Phi$ in terms of an interpolating polynomial functional 
of relatively small order, provided we have enough interpolation 
nodes nearby $\theta(x)=0$. 
This is demonstrated in Figure \ref{fig:Hopf_porter}, where we 
study the accuracy of Porter's interpolant through the set of nodes  
$\widehat{S}^{(m+1)}_n$ (see Eq. \eqref{SNq2}). Specifically, 
we plot the error 
\begin{equation}
E_m= \sup_{\theta\in \mathcal{G}_{50}(\sigma)}\left|\Phi([\theta])-\Pi_n([\theta])\right| 
\label{Em3}
\end{equation}
versus $m$ for $q=5$, $n=1,2,3,4$ and $\sigma=0.01$. 
It is seen that the interpolants converge in both $n$ (polynomial order) 
and $m$ (number of basis functions). Convergence in $m$ 
becomes monotonic for $m\geq10$. 
This is related to the fact that the effective dimension of 
the Hopf functionals \eqref{H_initial3} and \eqref{H_initial2} 
is $10$ for $q=5$ (see Section \ref{sec:Hopf_effective_dimensionality}). 
Therefore for $m\geq 10$ we have enough basis functions 
to fully resolve them. The accuracy of the polynomial 
interpolants then depends only on $n$ (polynomial order), 
the number of interpolation nodes and their location. 
The number of interpolation nodes in $\widehat{S}^{(m+1)}_n$
is given in \eqref{SCardinality}. For example, the case 
$m=20$ and $n=4$ yields $12650$ nodes and an 
interpolation matrix \eqref{matH} of size $12650\times12650$\footnote{Recall 
that the matrix \eqref{matH} has  to be inverted to compute the 
cardinal basis \eqref{gi_0}.}.
When evaluating the error \eqref{Em3} it is important to
select $\mathcal{G}_{50}(\sigma)$ within the convex hull 
of the interpolation nodes, e.g.,  
by choosing $\sigma$ small enough. In this way we avoid 
inaccuracies due to polynomial extrapolation. 
To mitigate this phenomenon, we can also consider different sets 
of interpolation nodes, e.g., the set $\widetilde{S}^{(m+1)}_n$ 
defined in \eqref{SNq3} with $a_{i_q}$ sampled at 
Gauss-Hermite sparse-grid points. 
In Figure \ref{fig:sparse_grid_sets} we plot the elements 
in $\widetilde{S}^{(m+1)}_1$,  $\widetilde{S}^{(m+1)}_2$ 
and $\widetilde{S}^{(m+1)}_3$ 
for Gauss-Hermite sparse grids of level 5 (see also Table \ref{tab:thenumber}). 
\begin{table}
\centering
\begin{tabular}{c|cccccccc}
 $m$   & $2$   & $4$ & $6$ & $8$ & $10$ & $16$ & $20$\\
 \hline\vspace{-0.4cm}\\
 $\#\widetilde{S}^{(m+1)}_2$  & $121$ & $341$ & $673$ & $1117$ & $1673$ & $4013$ & $6133$
\end{tabular}
\caption{Number of interpolation nodes in $\widetilde{S}^{(m+1)}_2$ with  
$a_{i_j}$ sampled at Gauss-Hermite sparse grids (level 5).}
\label{tab:thenumber}
\end{table}
The symmetry of the nodes $a_{i_j}$ in the Fourier space yields 
linearly dependent nodes in $\widetilde{S}^{(m+1)}_n$. 
Correspondingly, the interpolation matrix \eqref{matH} is 
rank-deficient, i.e., it cannot be inverted in a classical sense. 
We can overcome this issue by either removing 
some nodes from the set  $\widetilde{S}^{(m+1)}_n$, or by 
taking the More-Penrose pseudo-inverse of \eqref{matH}. In the latter 
case, we obtain a non-cardinal basis \eqref{gi_0_noncardinal} and a  
polynomial functional in the form \eqref{Porter_interpolant_noncardinal}.
In Figure \ref{fig:Hopf_porter_SG} we demonstrate convergence 
of such polynomial functional (order 1 and 2) to the
Hopf functionals \eqref{FDG} and \eqref{H_initial2}. 
Specifically, we plot the pointwise errors \eqref{Em3} 
versus $m$ for a set of $20000$ randomly generated evaluation 
nodes in $\mathcal{G}_{50}(\sigma)$, with $\sigma=1$.

\subsubsection{Canonical Tensor Decomposition}

We evaluate the Hopf functionals \eqref{FDG} and \eqref{H_initial2} 
in the function space spanned by the finite-dimensional 
orthonormal basis
\begin{equation}
\varphi_0(x)=\frac{1}{2\pi}, \quad 
\varphi_{k-1}(x)=\frac{\sin(kx/2)}{\sqrt{\pi}},\quad
\varphi_{k}(x)=\frac{\cos((1+k/2)x)}{\sqrt{\pi}}\qquad 
\textrm{($k$ even)}. 
\label{basisE}
\end{equation}
Hereafter we prove that this yields exact rank-one representations 
of both Gaussian and uniform Hopf functionals. 

\paragraph{Gaussian Functionals} 
By evaluating \eqref{FDG} in the space 
spanned by \eqref{basisE} we obtain 
\begin{equation}
f(a_1,...,a_m) = \exp\left[-\frac{1}{2}\sum_{p,q=1}^mZ_{pq} a_{p}a_q
\right],
\end{equation}
where 
\begin{equation}
Z_{ij}= \int_0^{2\pi}\int_0^{2\pi} C_0(x,y)\varphi_i(x)\varphi_j(y)dxdy.
\end{equation}
It can be verified that 
\begin{equation}
Z_{ij}=
\begin{cases}
 \displaystyle\frac{\pi}{q}\delta_{ij} & \quad i,j\leq 2q\\
 0 & \quad \textrm {otherwise}
\end{cases}
\nonumber
\end{equation}
and therefore the Gaussian functional \eqref{FDG} 
is {\em rank one} relative to the basis \eqref{basisE}, i.e.,
\begin{equation}
\Phi([\theta])=\exp\left[-\frac{\pi}{2q}\sum_{p=1}^{2q}
 a_{p}^2\right], \qquad a_p=(\theta,\varphi_p).
\end{equation}
Thus, a rank one canonical decomposition with $m=2q$ variables 
is exact in the case. It's important to remark that the 
Gaussian functional \eqref{FDG} is {\em not} rank one 
relative to other bases, for example \eqref{setbf}.

\paragraph{Uniform Functionals} 
By evaluating \eqref{H_initial2} in the space 
spanned by \eqref{basisE} we obtain 
\begin{equation}
f(a_1,...,a_m)=\prod_{k=1}^q\frac{\sin(\sqrt{3}s_k(a_1,..,a_m))\sin(\sqrt{3}c_k(a_1,..,a_m))} {3s_k(a_1,..,a_m)c_k(a_1,..,a_m)} ,
\end{equation}
where $s_k(a_1,..,a_m)$ and $c_k(a_1,..,a_m)$ are 
defined as  (see \eqref{coeffskck})
\begin{align}
s_k(a_1,..,a_m)= &\frac{1}{\sqrt{q}}\sum_{p=0}^m 
(\sin(kx),\varphi_p(x)) a_p,\label{re1}\\
c_k(a_1,..,a_m)= &\frac{1}{\sqrt{q}}\sum_{p=0}^m  
(\cos(kx),\varphi_p(x)) a_p.\label{re2} 
\end{align}
A substitution of \eqref{basisE} into \eqref{re1}-\eqref{re2} yields
\begin{align}
s_k(a_1,..,a_m)=& \sqrt{\frac{\pi}{q}} a_k,\quad k=1,...,q,
\quad (\textrm{zero otherwise}),\label{ss1}\\
c_k(a_1,..,a_m)=& \sqrt{\frac{\pi}{q}} a_{k+q},\quad k=1,...,q,
\quad (\textrm{zero otherwise}).\label{ss2}
\end{align}
This means that, relative to the  basis \eqref{basisE}
the functional \eqref{H_initial2} is {\em rank one}, i.e., 
\begin{equation}
\Phi([\theta])=\prod_{k=1}^{2q}
\frac{\sin(\sqrt{3\pi} a_k/\sqrt{q})}{3\sqrt{\pi}a_k/\sqrt{q}}, \qquad 
a_k=(\theta,\varphi_k).
\end{equation}
Therefore, a rank one canonical tensor decomposition with 
$m=2q$ variables is exact\footnote{Equations  
\eqref{ss1} and \eqref{ss2} suggest that 
all the $2q$ variables $a_1,...,a_{2q}$ are 
equally important.
}. Similarly to the Gaussian case, 
the uniform functional is not rank-one  relative to other 
bases, for example \eqref{setbf}.

\subsubsection{Functional Derivatives} 
\begin{figure}
\small
% \centerline{\line(1,0){430}\vspace{0.1cm}}
\centerline{Gaussian \hspace{4.3cm} Uniform}\vspace{-0.2cm}
\centerline{\line(1,0){300}}
\vspace{0.5cm }
\centerline{$\theta(x)=0$}
\centerline{\includegraphics[height=5cm]{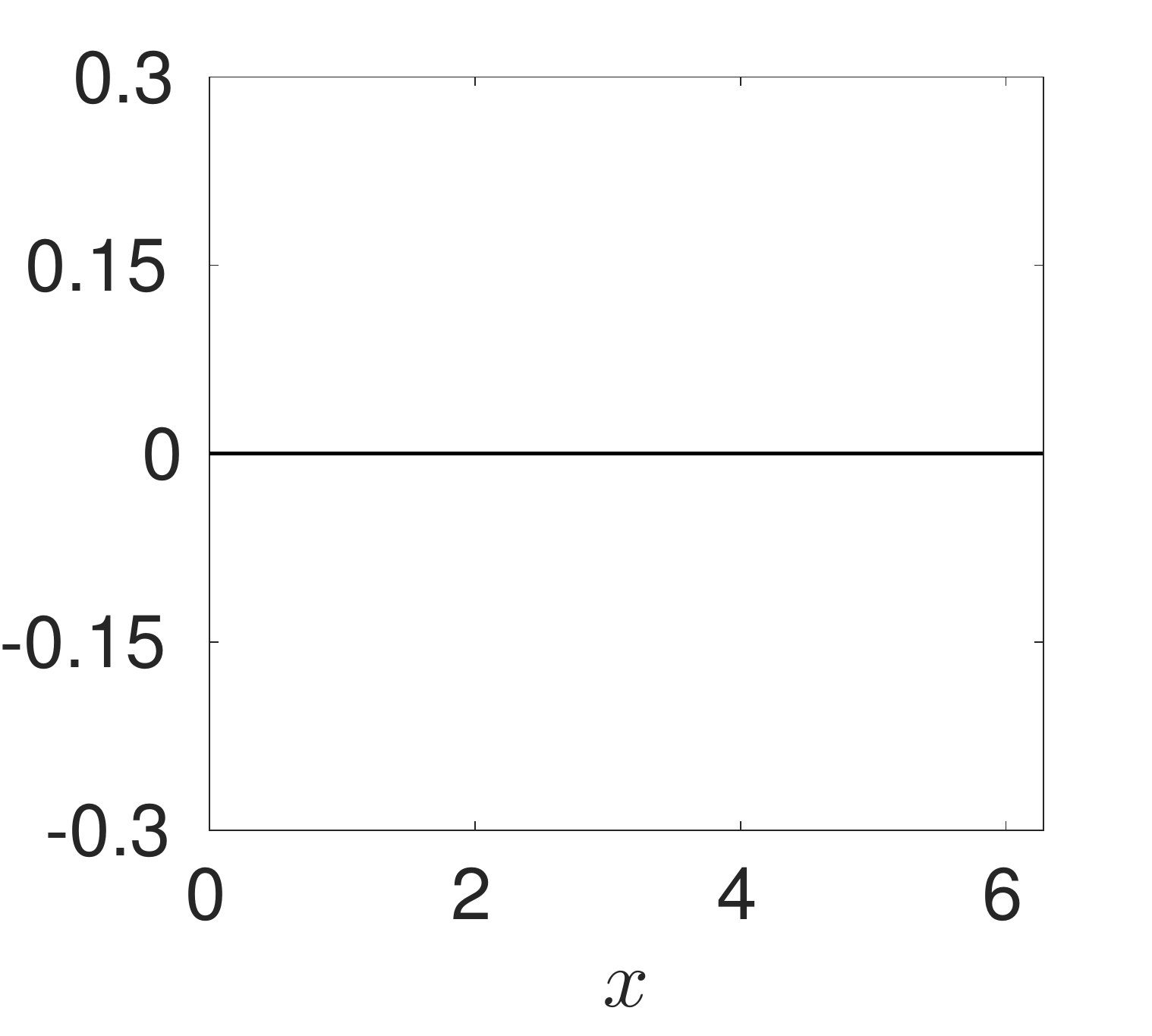}
            \includegraphics[height=5cm]{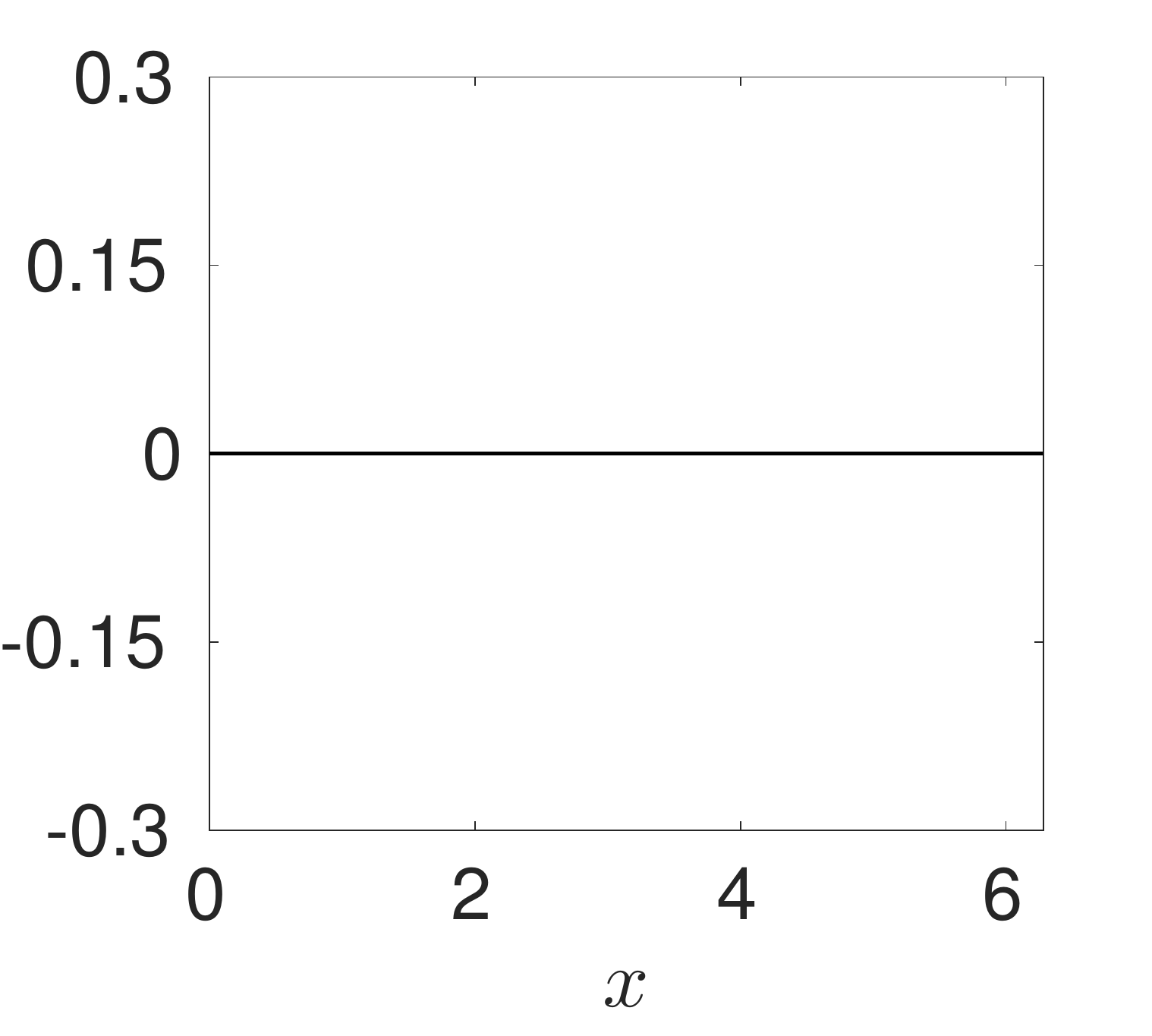}} 
\vspace{0.5cm}
\centerline{$\theta(x)=40\sin(6x)^2+5\cos(x)^4$}
\centerline{\includegraphics[height=5cm]{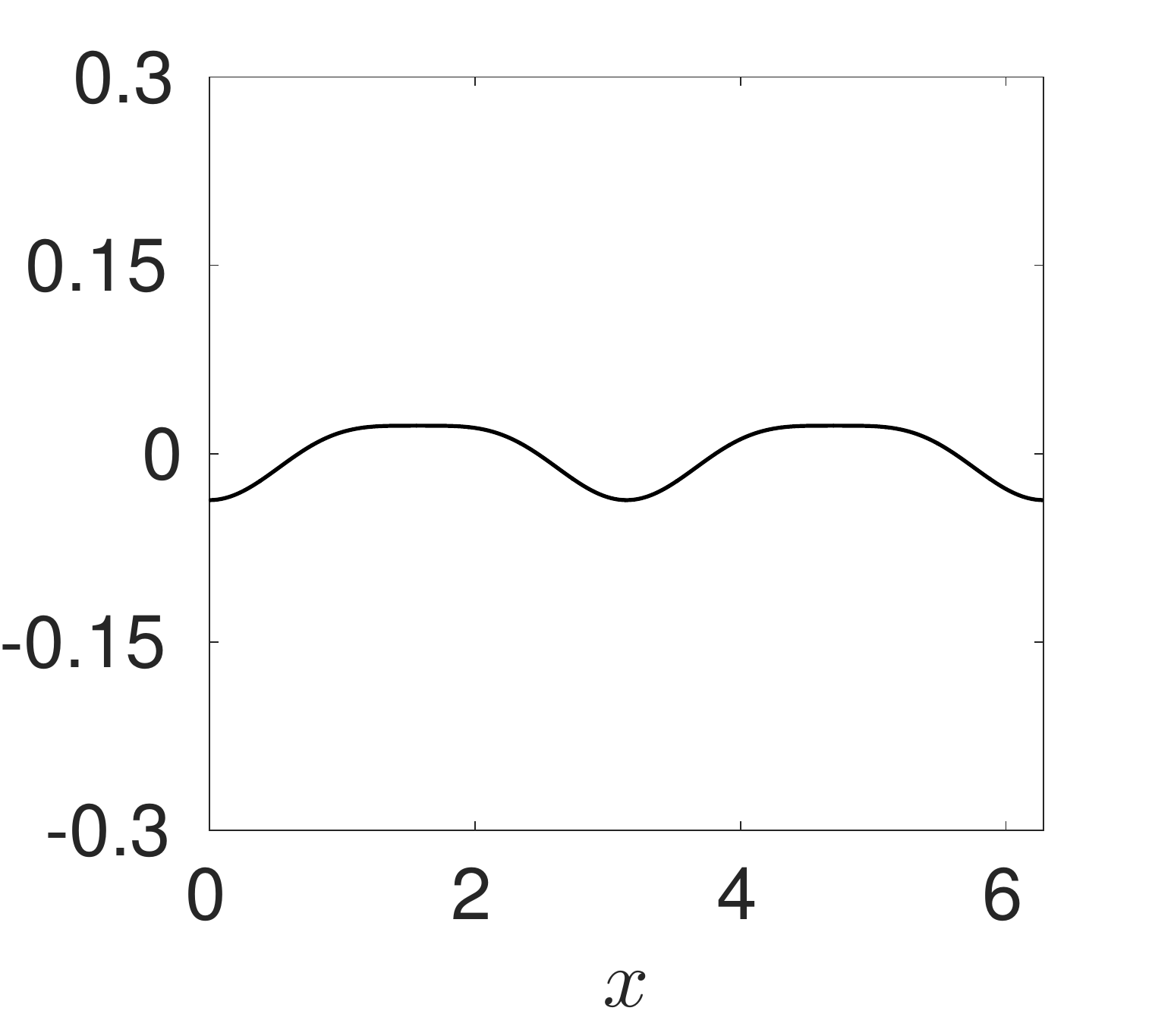}
            \includegraphics[height=5cm]{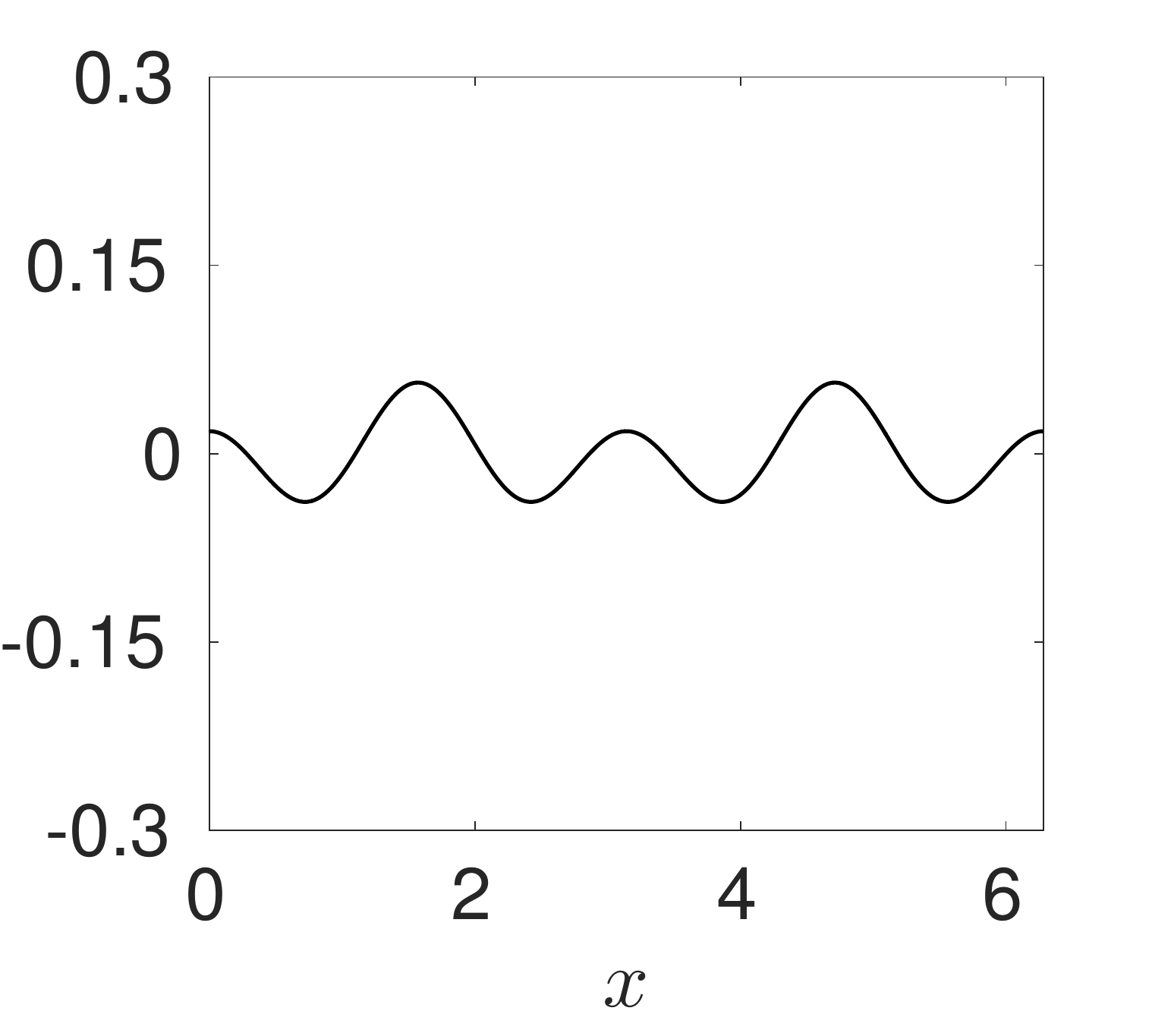}}
\vspace{0.5cm}
\centerline{$\theta(x)=2\exp(-\sin(4x))+\cos(\cos(4x))$}
\centerline{\includegraphics[height=5cm]{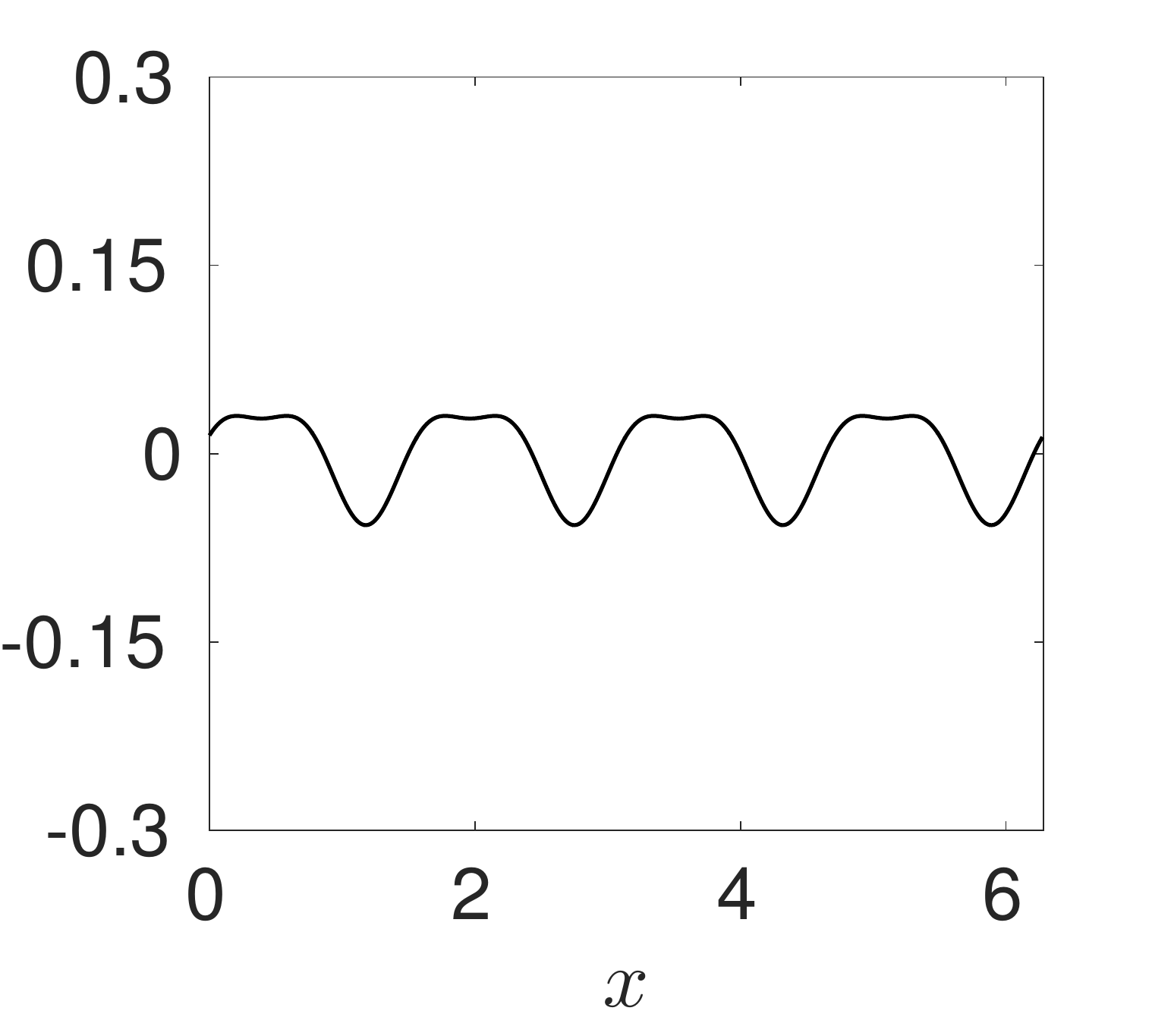}
            \includegraphics[height=5cm]{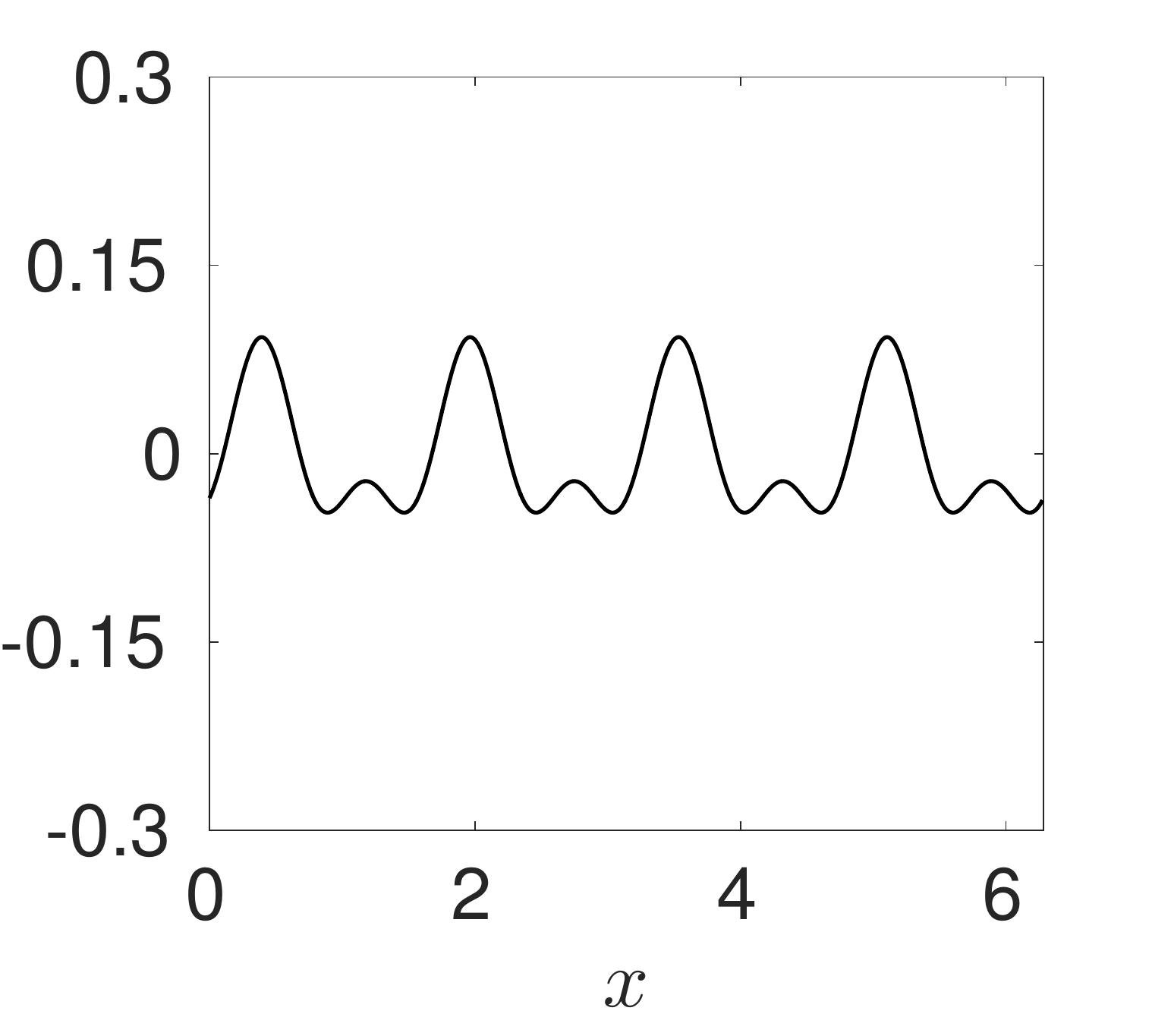}}
% \vspace{0.5cm}
\caption{First-order functional derivatives of the 
Hopf functionals \eqref{H_initial3} (Gaussian) and 
\eqref{H_initial2} (Uniform)  for $q=10$ evaluated at different test 
functions. The Hopf functionals we are considering here are real-valued 
and therefore the functional derivatives are real. Note that evaluating 
the first-order functional derivative at $\theta(x)=0$ (first row) 
yields the mean field $\left<u_0(x)\right>=0$.}
\label{fig:functional_derivatives1}
\end{figure}

\begin{figure}
\small
\centerline{Gaussian \hspace{4.8cm} Uniform}\vspace{-0.2cm}
\centerline{\line(1,0){300}}
\vspace{0.5cm }
\centerline{$\theta(x)=0$}
\centerline{\includegraphics[height=5cm]{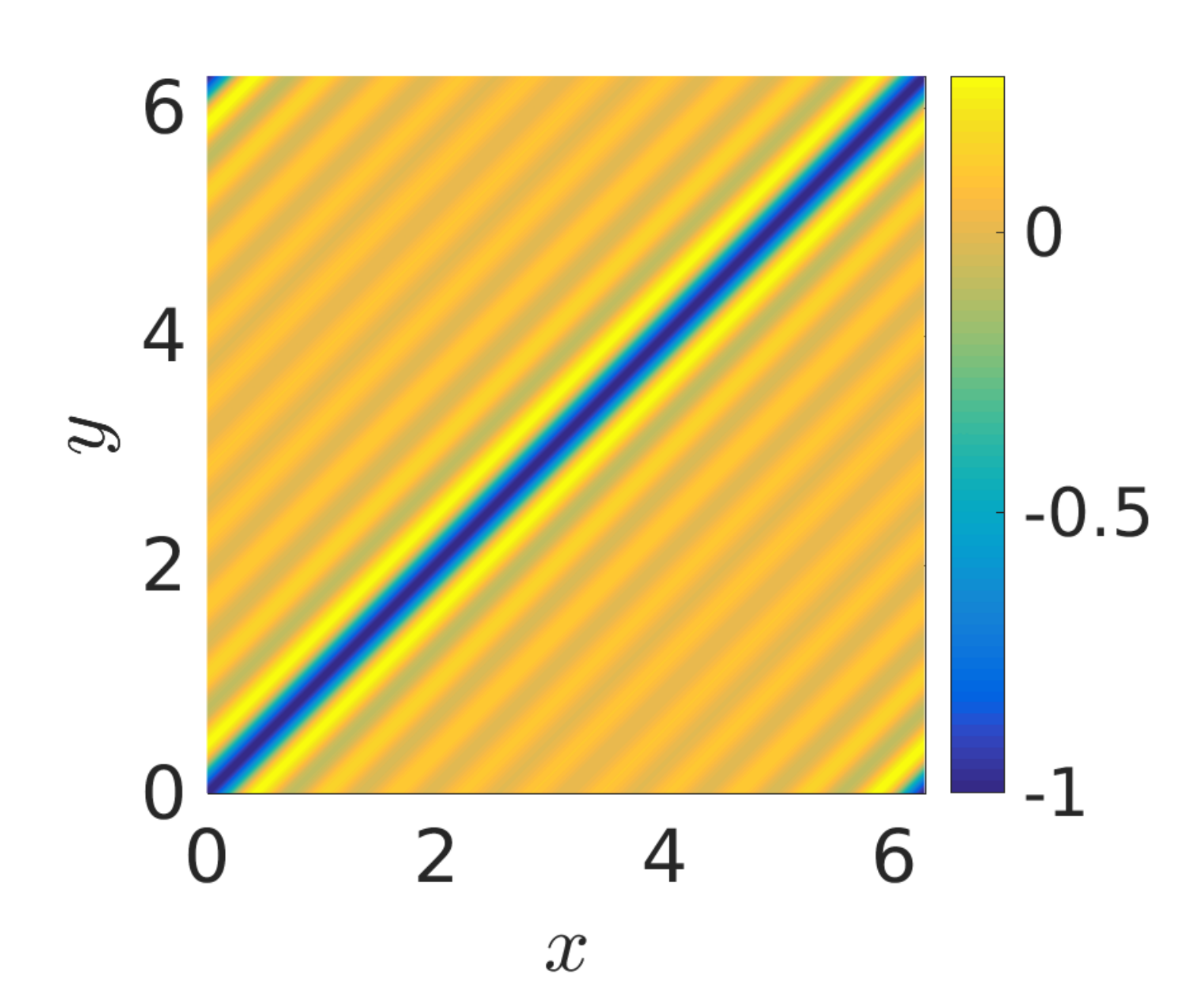}
            \includegraphics[height=5cm]{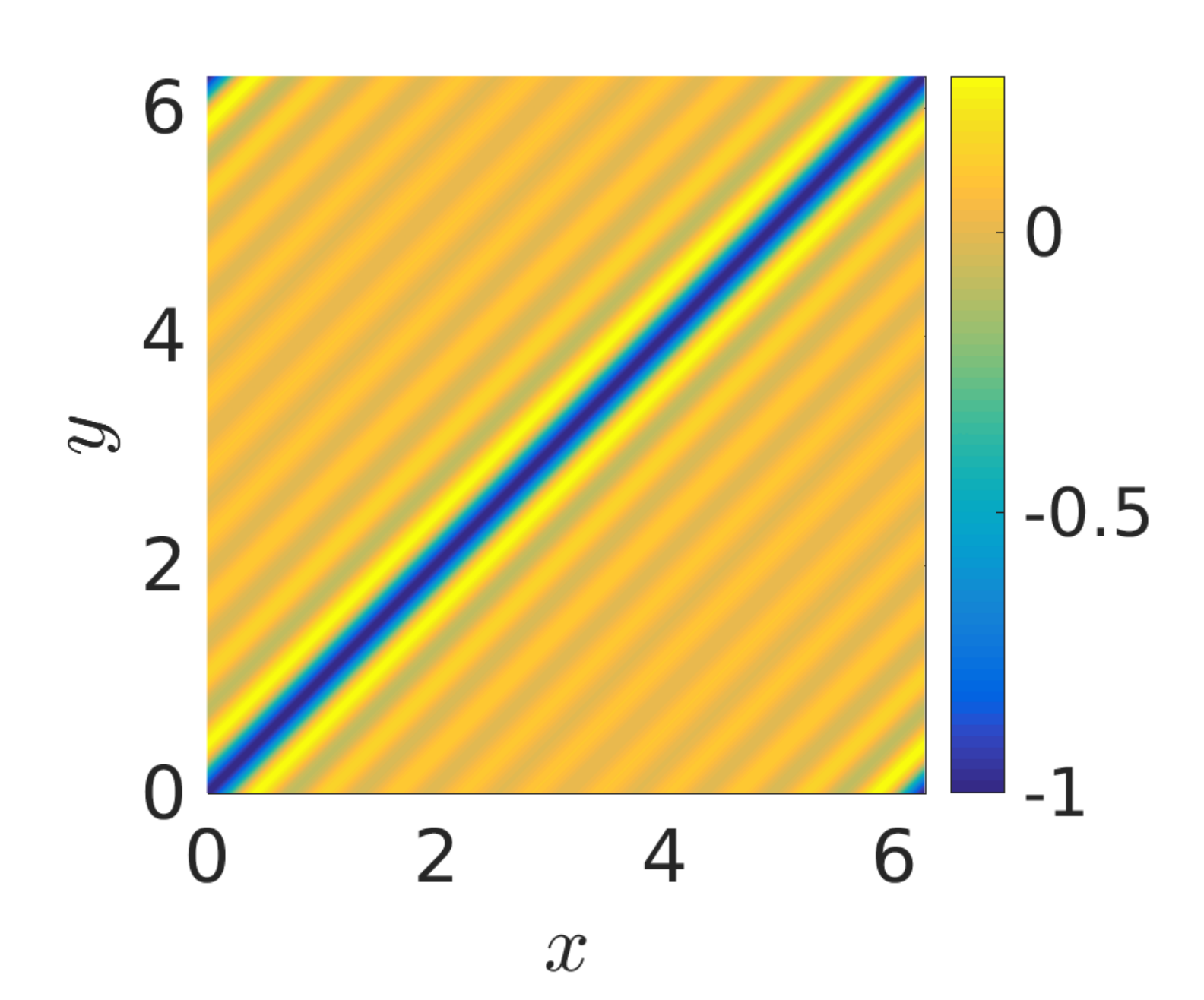}} 
\vspace{0.5cm}
\centerline{$\theta(x)=40\sin(6x)^2+5\cos(x)^4$}
\centerline{\includegraphics[height=5cm]{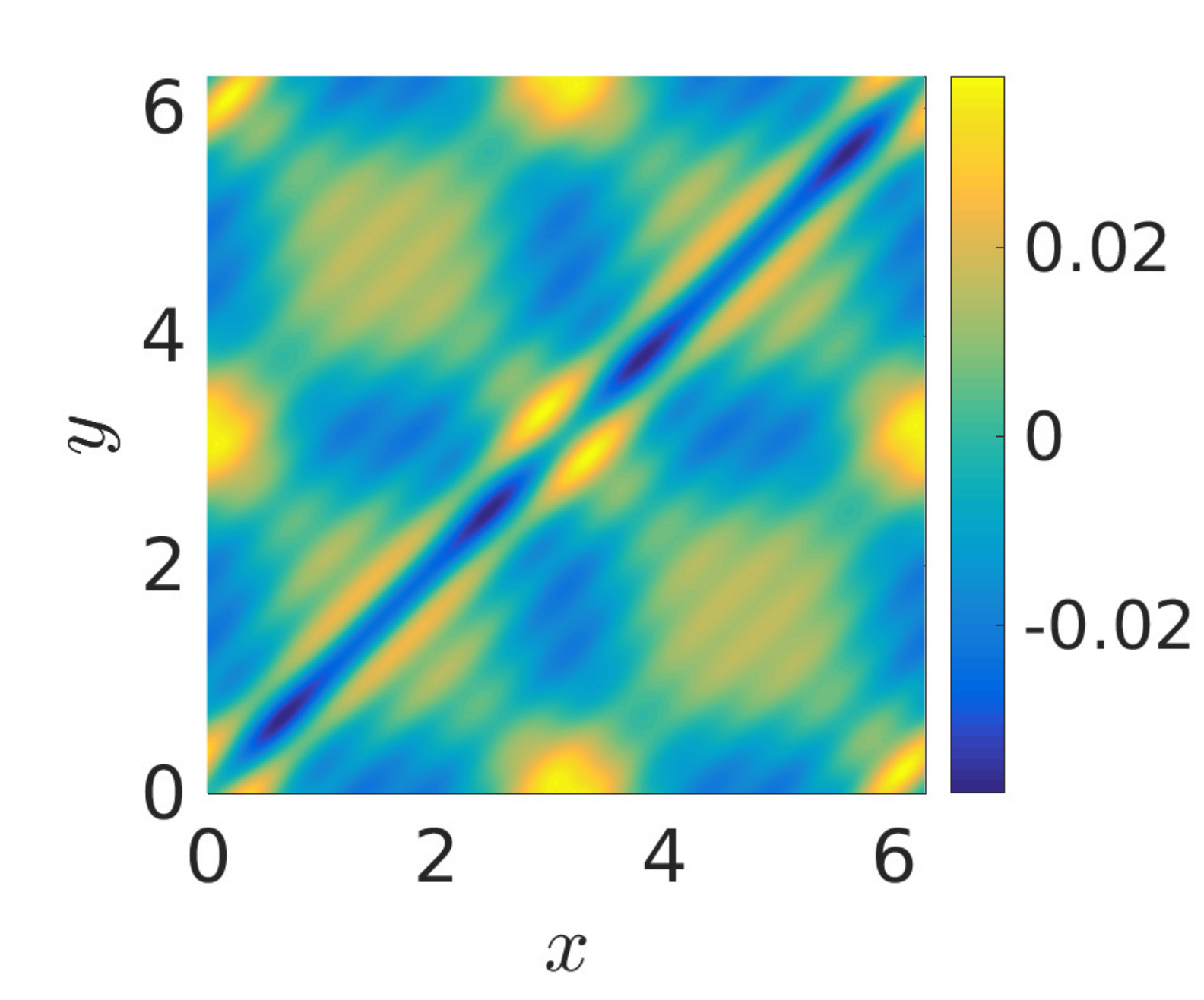}
            \includegraphics[height=5cm]{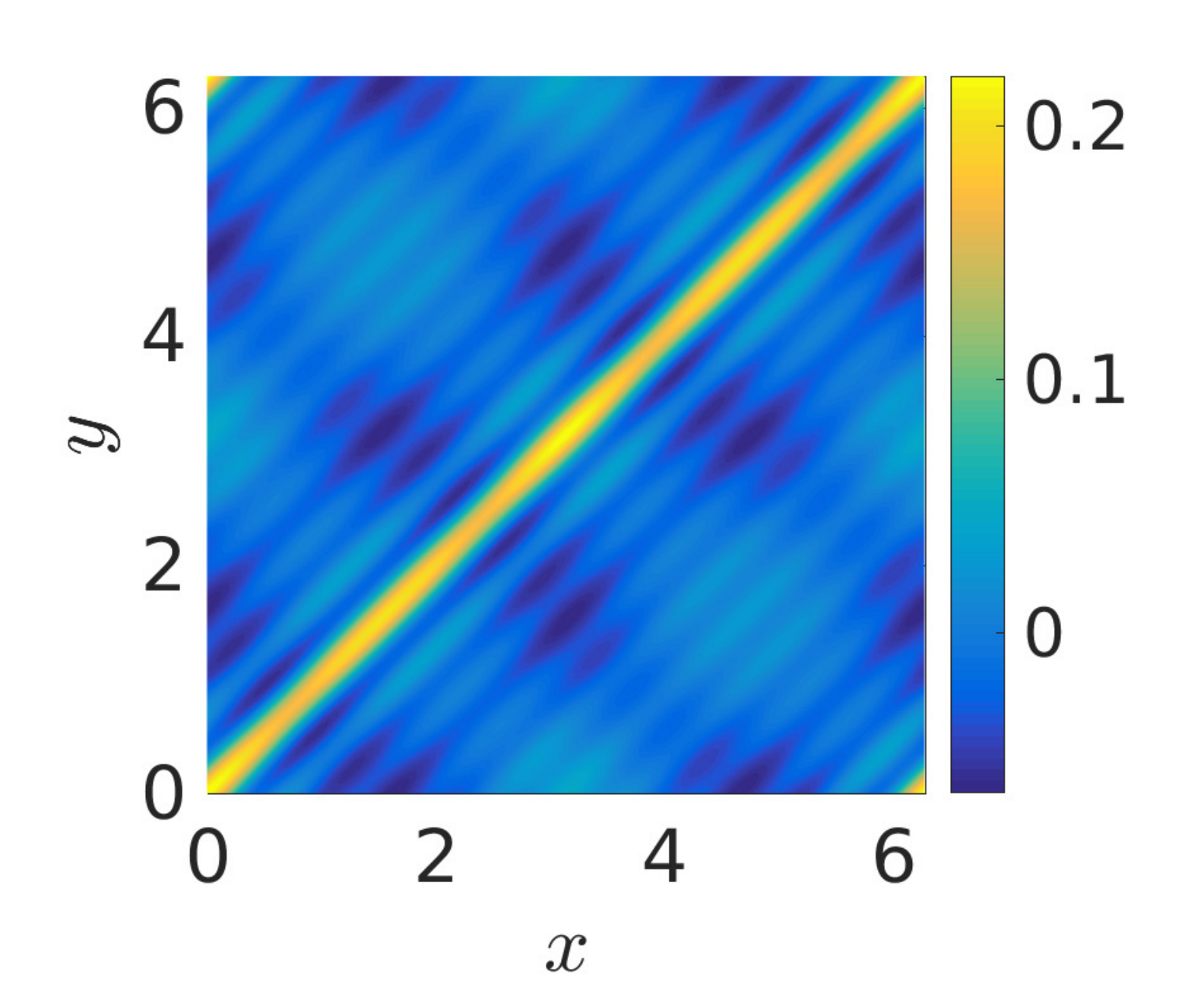}}
\vspace{0.5cm}
\centerline{$\theta(x)=2\exp(-\sin(4x))+\cos(\cos(4x))$}
\centerline{\includegraphics[height=5cm]{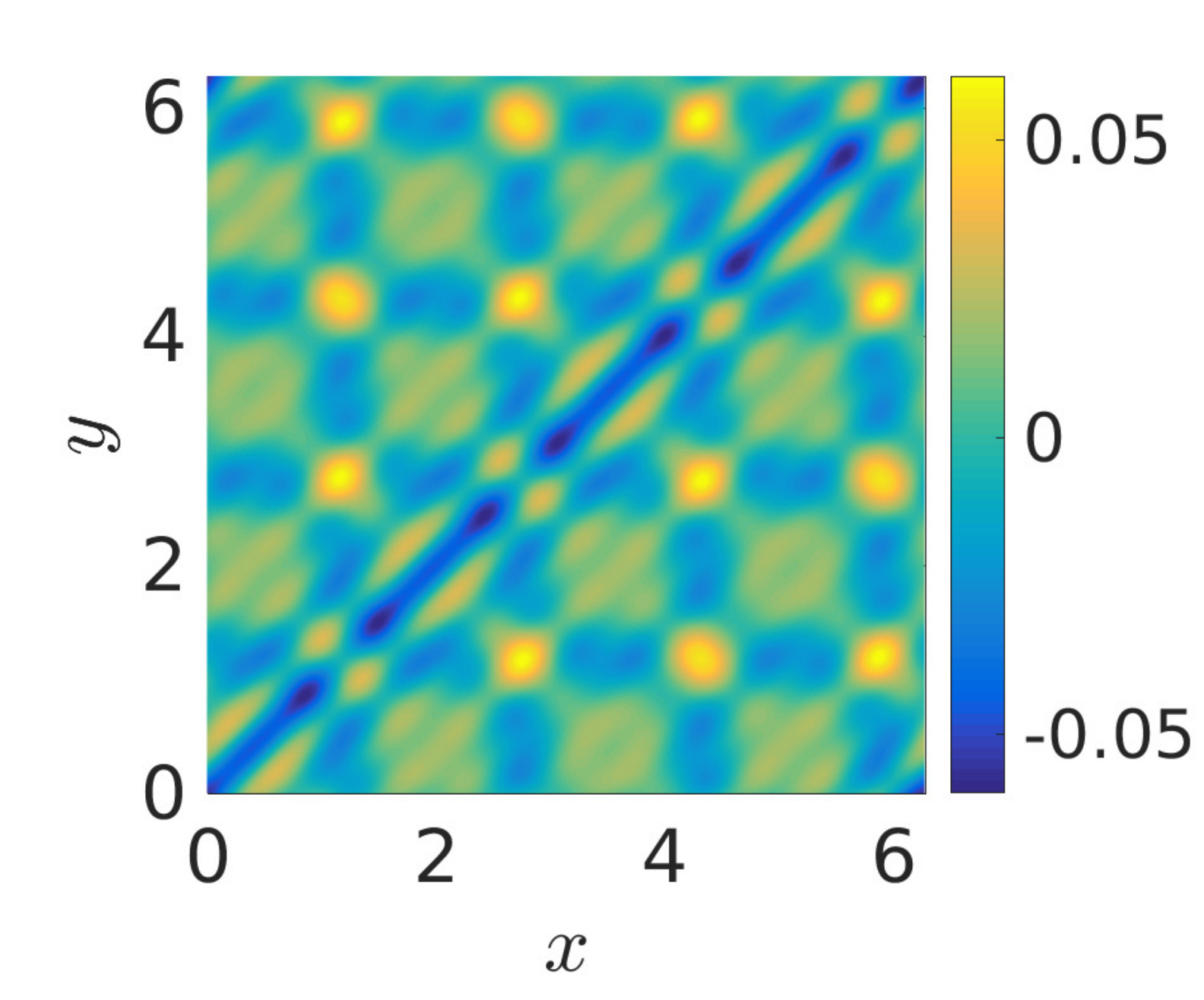}
            \includegraphics[height=5cm]{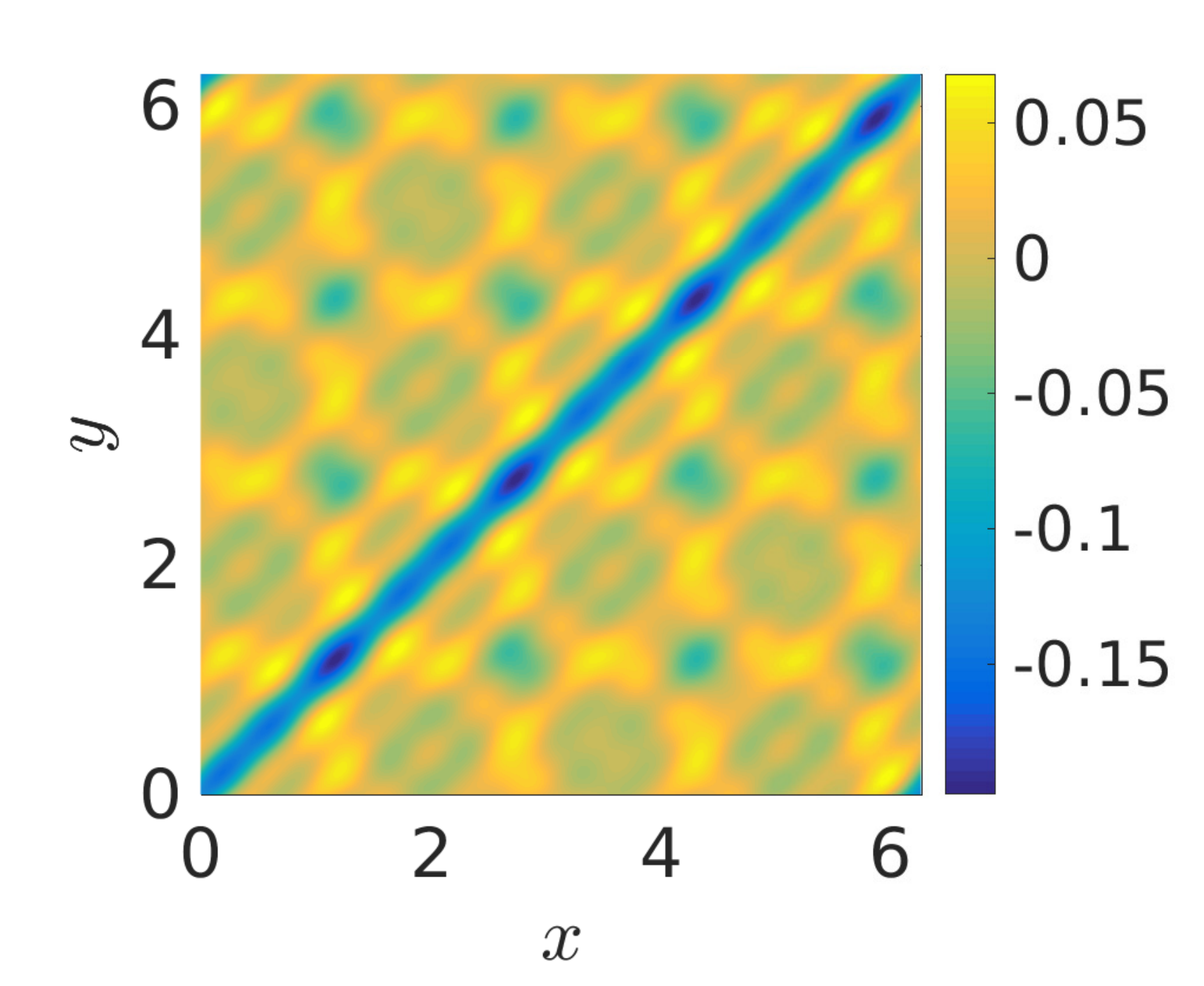}}
% \vspace{0.5cm}
\caption{Second-order functional derivatives of the 
Hopf functionals \eqref{H_initial3} (Gaussian) and 
\eqref{H_initial2} (Uniform) for $q=10$ evaluated at different test 
functions. Both functionals are real-valued 
and therefore the functional derivatives are real. Note that evaluating 
the second-order functional derivative at $\theta(x)=0$ (first row) 
yields the opposite of the correlation function $C_0(x,y)$ in \eqref{testcorrelation}.}
\label{fig:functional_derivatives2}
\end{figure}

The first- and second-order functional derivatives 
of the Hopf functional \eqref{H_initial} are defined as 
\begin{align}
\frac{\delta \Phi([\theta])}{\delta \theta(x)}&=i\left<u_0(x;\omega)\exp\left[i\int_{0}^{2\pi}
 \theta(x)u_0(x;\omega)dx\right]\right>,\label{Ai}\\
\frac{\delta^2 \Phi([\theta])}{\delta \theta(x)\delta\theta(y)}&=
-\left<u_0(x;\omega)u_0(y;\omega)\exp\left[i\int_{0}^{2\pi}
 \theta(x)u_0(x;\omega)dx\right]\right>.\label{Bi}
\end{align}
Remarkably, these derivatives can 
be expressed analytically in terms of simple functions.
To this end, we need a formula to compute 
the averages in \eqref{Ai} and \eqref{Bi}. 
A lengthy calculation shows that 
\begin{equation}
\frac{\delta \Phi([\theta])}{\delta \theta(x)}=
\frac{i}{\sqrt{q}}\left(I^{(c)}[\theta]\sum_{k=1}^q\sin(kx)I^{(s)}_k[\theta]+
 I^{(s)}[\theta]\sum_{k=1}^q\cos(kx)I^{(c)}_k[\theta]\right),\label{anA}
\end{equation}
\begin{align}
\frac{\delta^2 \Phi([\theta])}{\delta \theta(x)\delta\theta(y)}=&
-\frac{1}{q}\sum_{k,h=1}^q\left(\sin(kx)\cos(hy)I^{(s)}_{k}[\theta]I^{(c)}_{h}[\theta]+
\cos(kx)\sin(hy)I^{(c)}_{k}[\theta]I^{(s)}_{h}[\theta]\right)\nonumber\\
&-\frac{1}{q}\sum_{\substack{k,h=1\\k\neq h}}^q\left(
I^{(c)}[\theta]\sin(kx)\sin(hy)I^{(s)}_{kh}[\theta]+
I^{(s)}[\theta]\cos(kx)\cos(hy)I^{(c)}_{kh}[\theta]\right)\nonumber \\
&-\frac{1}{q}\sum_{k=1}^q\left(
I^{(c)}[\theta]\sin(kx)\sin(ky)J^{(s)}_{k}[\theta]+
I^{(s)}[\theta]\cos(kx)\cos(ky)J^{(c)}_{k}[\theta]\right),
\label{anB}
\end{align}
where the functionals $I^{(c)}[\theta]$, $I^{(s)}[\theta]$, $I_k^{(c)}[\theta]$, etc., 
are defined  in Table \ref{tab:coefficients}.
These expressions are general and they hold for any random function in the 
form \eqref{Burgers_initial}, with i.i.d. random variables $\eta_k$ and $\xi_j$. 
\begin{table}[!t]
\small
\centerline{\line(1,0){430}\vspace{0.2cm}}
\centerline{Coefficients Appearing in the Functional Derivatives \eqref{anA} and \eqref{anB}}
\centerline{\line(1,0){430}}
% $\displaystyle G_n(z)=\int_{-\infty}^\infty a^n e^{z a}p(a)da$
\vspace{0.5cm}
\centerline{\small$\displaystyle G_n(z)=\int_{-\infty}^\infty a^n e^{iz a}p(a)da$}
\vspace{1cm}
\centerline{
\begin{tabular}{ll}
\small
Gaussian PDF
& $G_0(z)= e^{-z^2/2}$  \\ 
            & $G_1(z)= ize^{-z^2/2}$ \\
            & $G_2(z)=(1-z^2)e^{-z^2/2}$\\\\
Uniform PDF & $\displaystyle G_0(z)= \frac{\sin(\sqrt{3}z)}{\sqrt{3}z}$ \\
            & $\displaystyle G_1(z)= i\frac{\sin(\sqrt{3}z)}{\sqrt{3}z^2}
             -i\frac{\cos(\sqrt{3}z)}{z}$ \\ 
            & $\displaystyle G_2(z)= \frac{2\cos(\sqrt{3}z)}{z^2}
             +\frac{\sqrt{3}\sin(\sqrt{3}z)}{z}
             -\frac{2\sin(\sqrt{3}z)}{\sqrt{3}z^3}$            
\end{tabular}}
\vspace{1cm}
\centerline{
\begin{tabular}{lcl}
\small
 $\displaystyle I^{(s)}[\theta]=\prod_{r=1}^q G_0(s_r[\theta])$ & &
$\displaystyle I^{(c)}[\theta]=\prod_{r=1}^qG_0(c_r[\theta])$\\
$\displaystyle I_k^{(s)}[\theta]=G_1(s_k[\theta]) 
  \prod_{\substack{r=1\\r\neq k}}^qG_0(s_r[\theta])$ & &
$\displaystyle I_k^{(c)}[\theta]=G_1(c_k[\theta]) 
  \prod_{\substack{r=1\\r\neq k}}^qG_0(c_r[\theta])$\\
$\displaystyle J_k^{(s)}[\theta]=G_2(s_k[\theta]) 
  \prod_{\substack{r=1\\r\neq k}}^qG_0(s_r[\theta])$ &  &
$\displaystyle J_k^{(c)}[\theta]=G_2(c_k[\theta]) 
  \prod_{\substack{r=1\\r\neq k}}^qG_0(c_r[\theta])$\\
$\displaystyle I_{kh}^{(s)}[\theta]=G_1(s_k[\theta]) G_1(s_h[\theta]) 
  \prod_{\substack{r=1\\r\neq k,h}}^qG_0(s_r[\theta])$ & &
$\displaystyle I_{kh}^{(c)}[\theta]=G_1(c_k[\theta]) G_1(c_h[\theta]) 
  \prod_{\substack{r=1\\r\neq k,h}}^qG_0(c_r[\theta])$ \\\\\hline\\
\end{tabular}}
\caption{Coefficients appearing in the functional derivatives 
\eqref{anA} and \eqref{anB}.}
\label{tab:coefficients}
\end{table}
Note that in the case of uniform PDF the coefficients 
in Table \ref{tab:coefficients} satisfy
\begin{equation}
\lim_{z\rightarrow 0}G_0(z)=1,\qquad
\lim_{z\rightarrow 0}G_1(z)=0,\qquad
\lim_{z\rightarrow 0}G_2(z)=1.
\end{equation}
and therefore there are no singularities at $z=0$. Also, all coefficients 
decay to $0$ when $z$ goes to infinity.
The expressions \eqref{anA} and \eqref{anB} can be simplified 
significantly for Gaussian Hopf functionals \eqref{FDG}. In particular, we obtain\footnote{Note that the first-order derivative reduces 
to $0$ (mean field) at $\theta=0$ while the second-order 
derivative at $\theta=0$ reduces to the opposite of the 
correlation function $C_0(x,y)$.}
\begin{align}
\frac{\delta \Phi([\theta])}{\delta \theta(x)}=-
\left[\int_{0}^{2\pi}C_0(x,y)\theta(y)dy\right]
\exp\left[-\frac{1}{2}\int_{0}^{2\pi} \int_{0}^{2\pi}C_0(x,y)\theta(x)\theta(y)dxdy\right]
\label{fD1},
\end{align}
\begin{align}
\frac{\delta^2 \Phi([\theta])}{\delta \theta(x)\delta\theta(y)}=
&\left[-C_0(x,y)+\int_{0}^{2\pi}C_0(x,y)\theta(y)dy\int_{0}^{2\pi}C_0(x,y)\theta(x)dx\right]
\times\nonumber\\
& \exp\left[-\frac{1}{2}\int_{0}^{2\pi} \int_{0}^{2\pi}C_0(x,y)\theta(x)\theta(y)dxdy\right].
\label{fD2}
\end{align}
In Figure \ref{fig:functional_derivatives1} and 
Figure \ref{fig:functional_derivatives2} we plot 
the first- and second-order functional derivatives of 
the Hopf functionals \eqref{H_initial3} and 
\eqref{H_initial2},  evaluated at different test functions. 

\begin{figure}[!ht]
\centerline{\hspace{1.8cm}Analytical \hspace{3.7cm} Canonical Tensor Decomposition}
\vspace{0.cm}
\centerline{\hspace{-.6cm}
 \includegraphics[height=4.3cm]{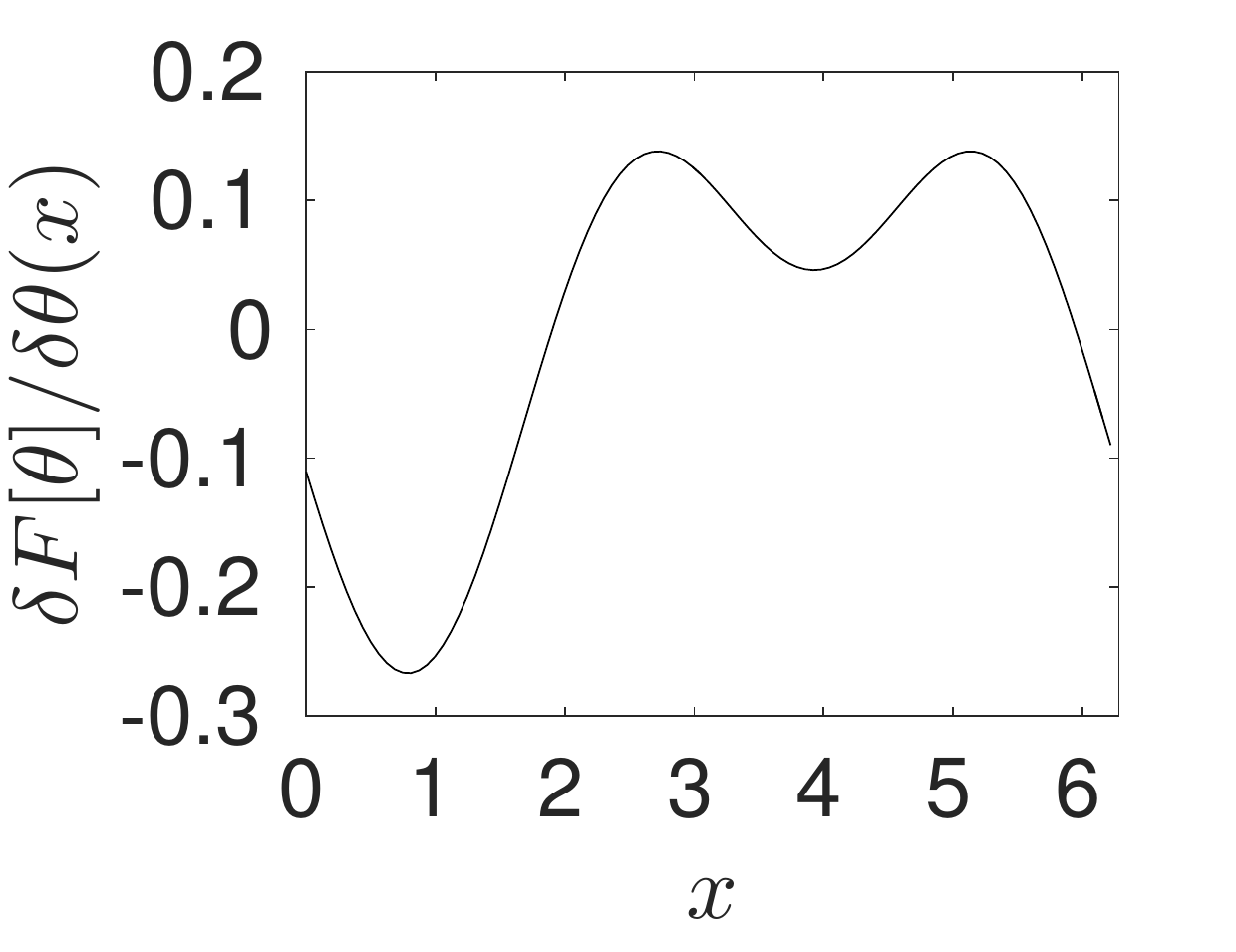}\hspace{1.8cm}
\includegraphics[height=4.3cm]{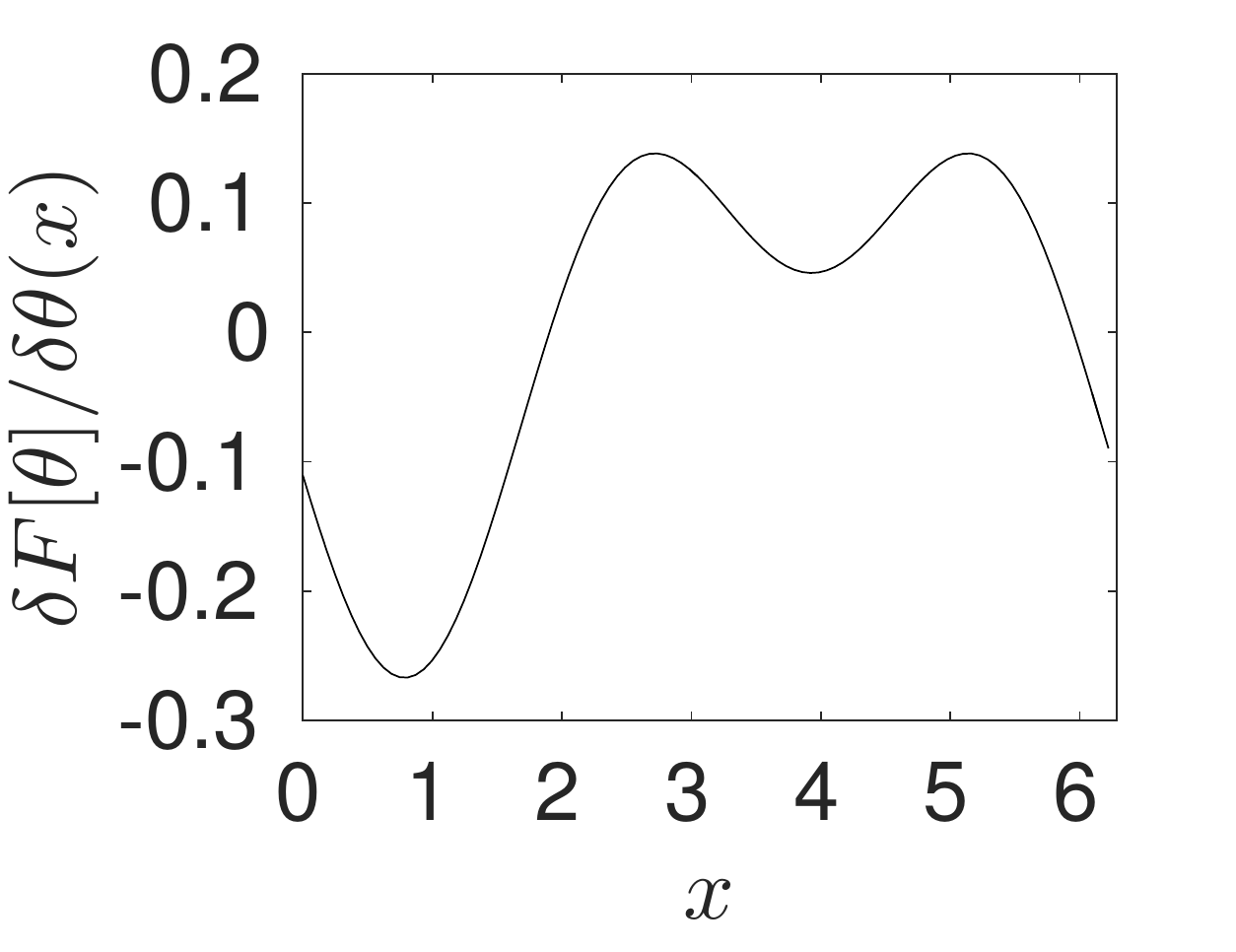}
 }

\centerline{
 \includegraphics[height=5.cm]{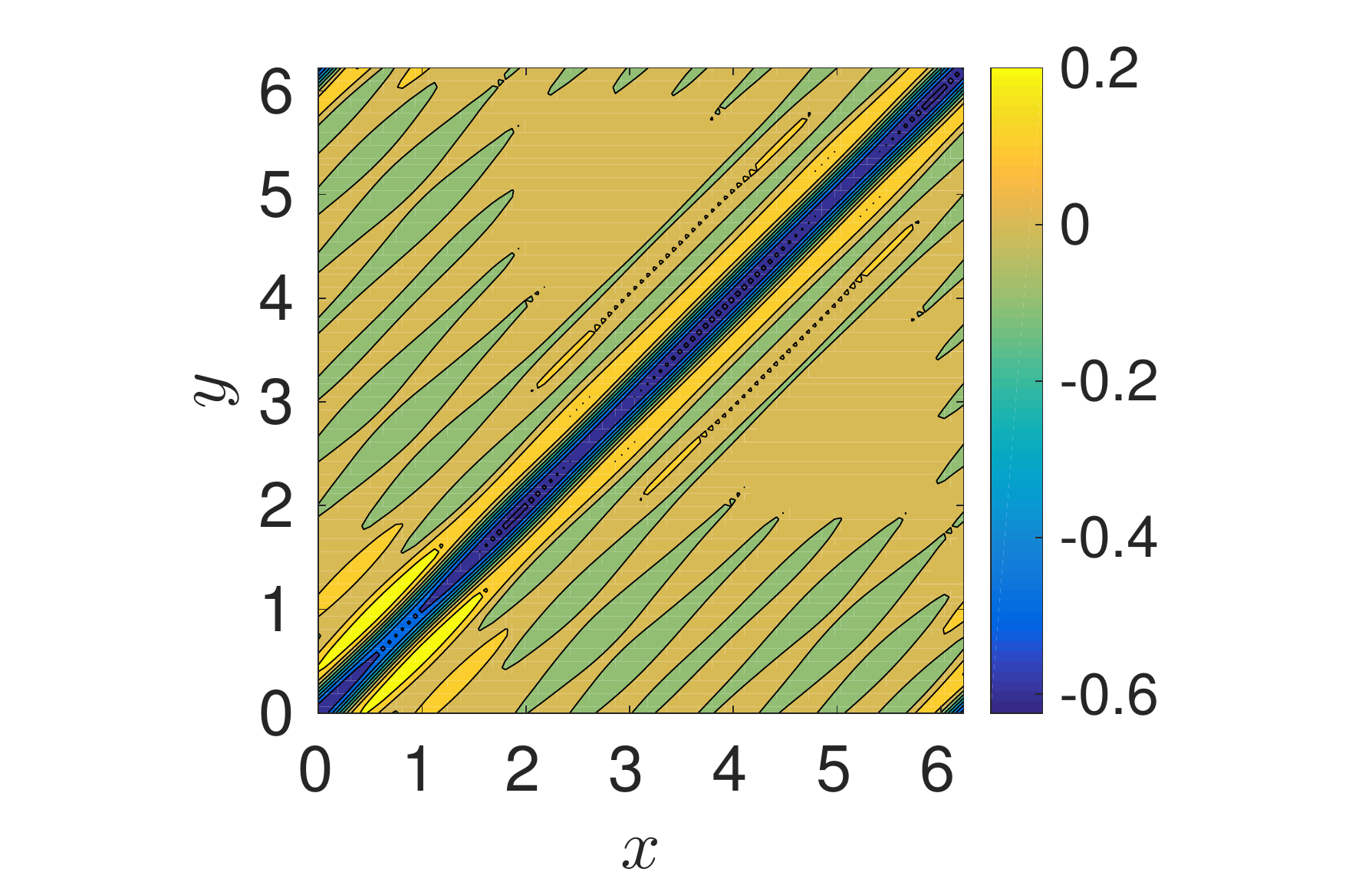}
 \includegraphics[height=5.cm]{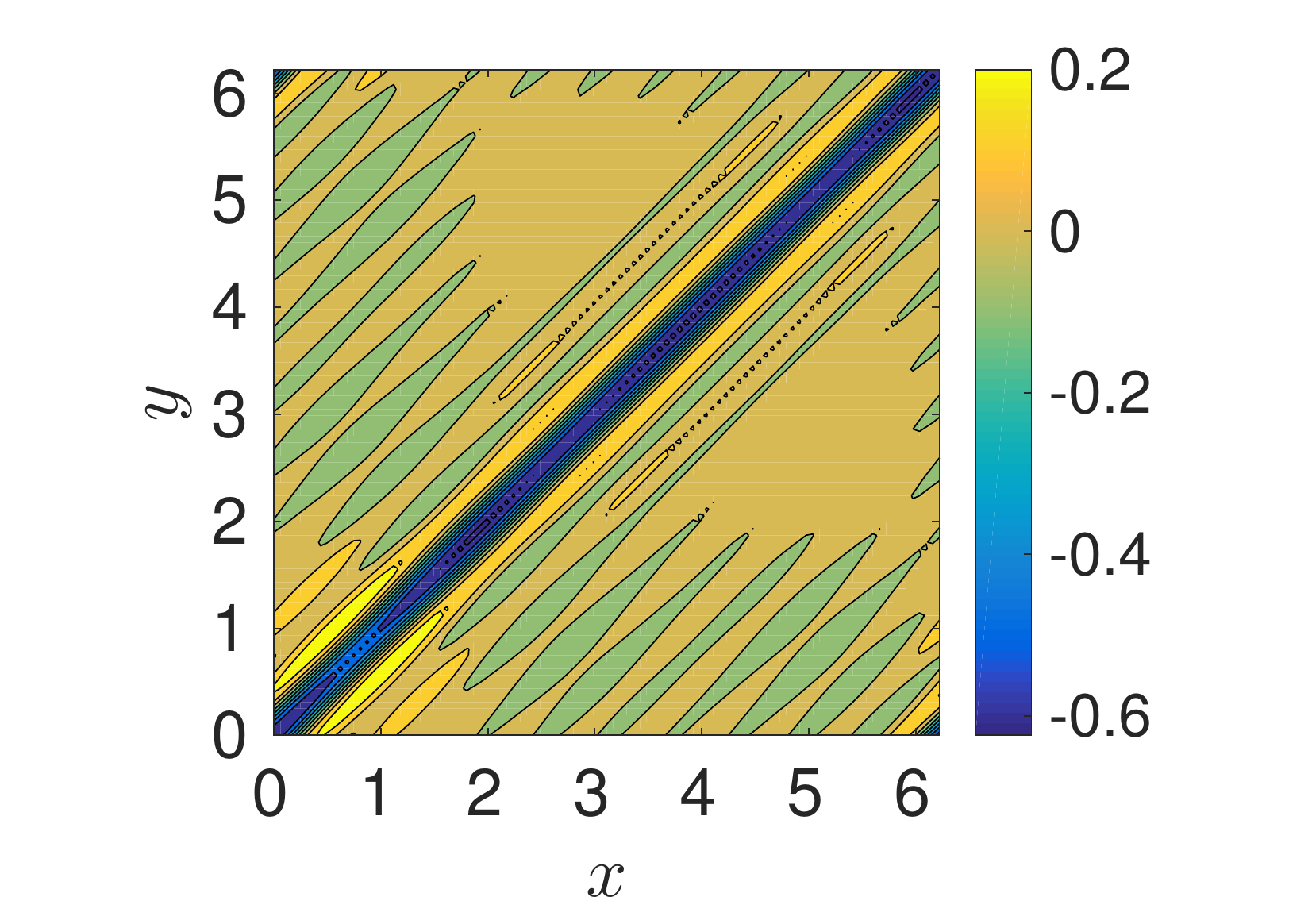}
 }
 
\caption{First- and second-order functional 
derivatives of the Gaussian Hopf functional \eqref{FDG} evaluated 
at $\theta(x)=\sin(x)+\sin(2x)+\sin(3x)$ . 
Shown are analytical results \eqref{fD1} 
and \eqref{fD2} versus numerical result obtained by 
the canonical tensor decomposition 
(equations \eqref{SSEfd1} and \eqref{SSEfd2}).}
\label{fig:accuracyHopfFDs}
\end{figure}

\paragraph{Canonical Tensor Expansion of Functional Derivatives} 
We have seen in Section \ref{sec:tensor} that the functional 
derivatives of any cylindrical approximant of $\Phi$ can be expressed as 
\begin{equation}
\frac{\delta {\Phi}([\theta])}{\delta\theta(x)}\simeq \sum_{k=1}^m
\frac{\partial f}{\partial a_k}\varphi_k(x),\label{SSEfd1}
\end{equation}
\begin{equation}
\frac{\delta^2 {\Phi}([\theta])}{\delta\theta(x)\delta\theta(y)}\simeq \sum_{k,j=1}^m
\frac{\partial^2 f}{\partial a_k\partial a_j}\varphi_k(x)\varphi_j(y),
\label{SSEfd2}
\end{equation}
where $a_k =(\theta,\varphi_k)$. In particular, if we consider 
a canonical tensor expansion of $\Phi$ then $f$ has the form 
\eqref{functional-SSE}, and the partial derivatives in 
\eqref{SSEfd1} and \eqref{SSEfd2} can be easily computed. 
In Figure \ref{fig:accuracyHopfFDs} we 
show that \eqref{SSEfd1} and \eqref{SSEfd2} 
provide a very accurate approximation of the functional 
derivatives \eqref{fD1} and \eqref{fD2}. 

\subsection{Sine Functional}
\label{sec:sine functional}
\begin{figure}[!ht]
\centerline{\hspace{1cm}Polynomial functionals \hspace{3.5cm} Canonical Tensor Decomposition}
\vspace{0.cm}
\centerline{
 \includegraphics[height=6.cm]{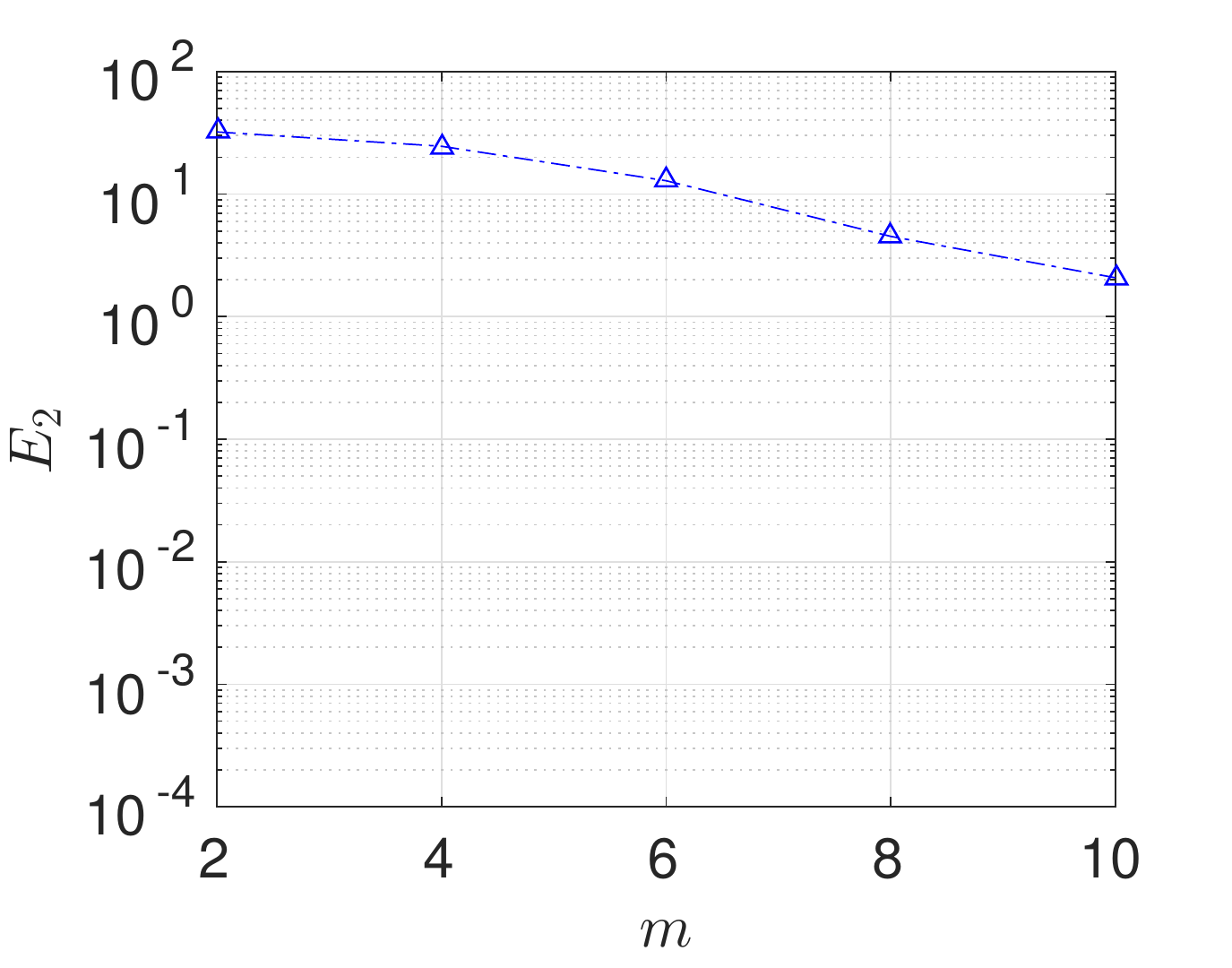}
 \includegraphics[height=6.cm]{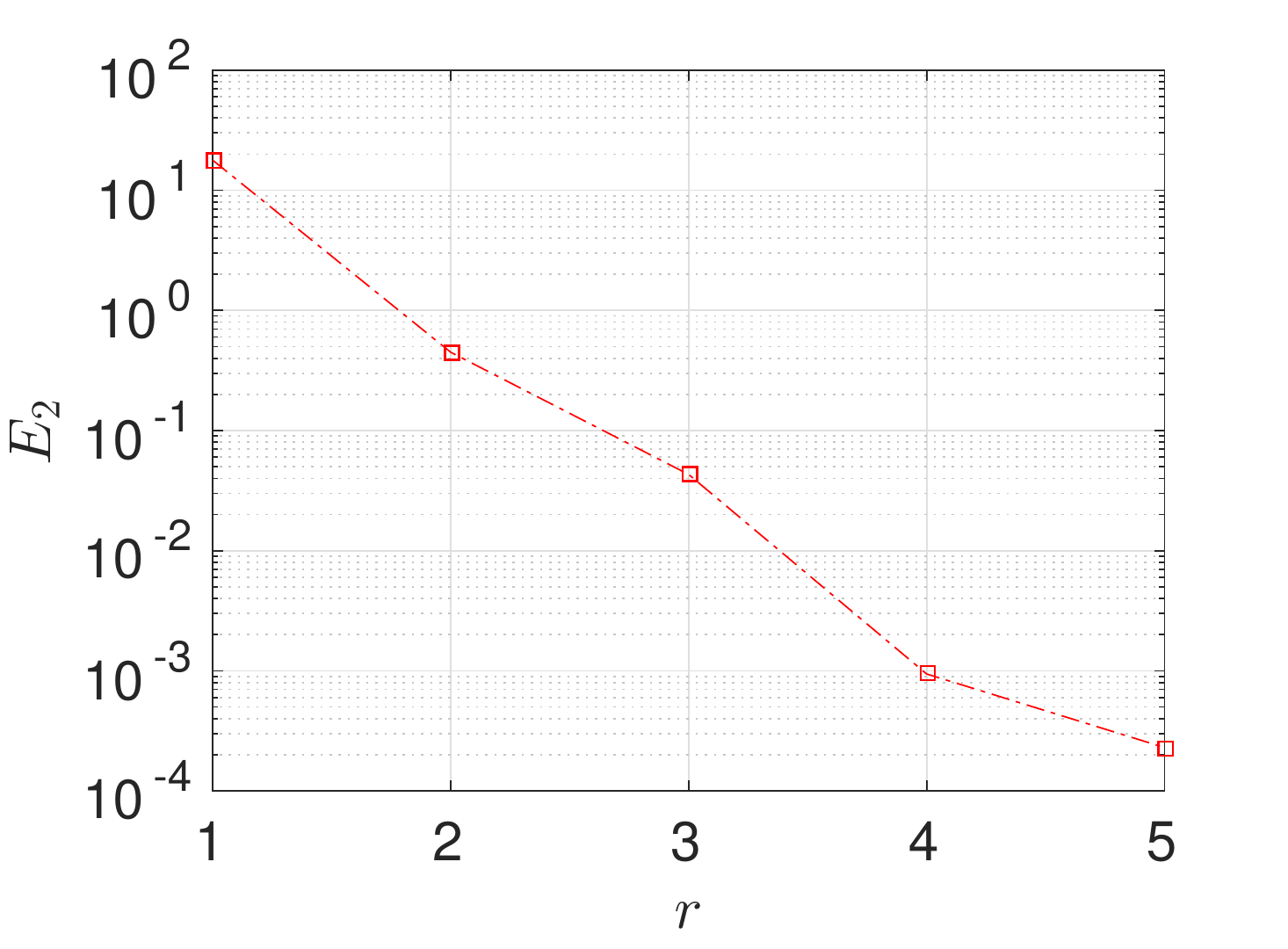}}
 
\caption{Approximation of the sine functional \eqref{SINEF}. 
Shown are the $L_2([-1,1]^m)$ errors obtained by using Porter's 
polynomial functionals \eqref{Porter_interpolant_noncardinal} 
of second-order (left) with non-cardinal basis,
and canonical tensor decomposition (right). Specifically, 
we study convergence as a function of the number of 
interpolation nodes in Porter's method: $m=10$ and sparse grids 
level 5 yields $\#\widetilde{S}^{(m+1)}_2=1673$ nodes (see Table \ref{tab:thenumber}). In the canonical tensor decomposition 
method we show convergence of \eqref{functional-SSE} as 
a function of the separation rank ($r$). The dimension of the 
test function space is chosen to be $m=10$. 
}
\label{fig:porterSSEsine}
\end{figure}

\begin{figure}[!ht]
\centerline{
\includegraphics[height=3.5cm]{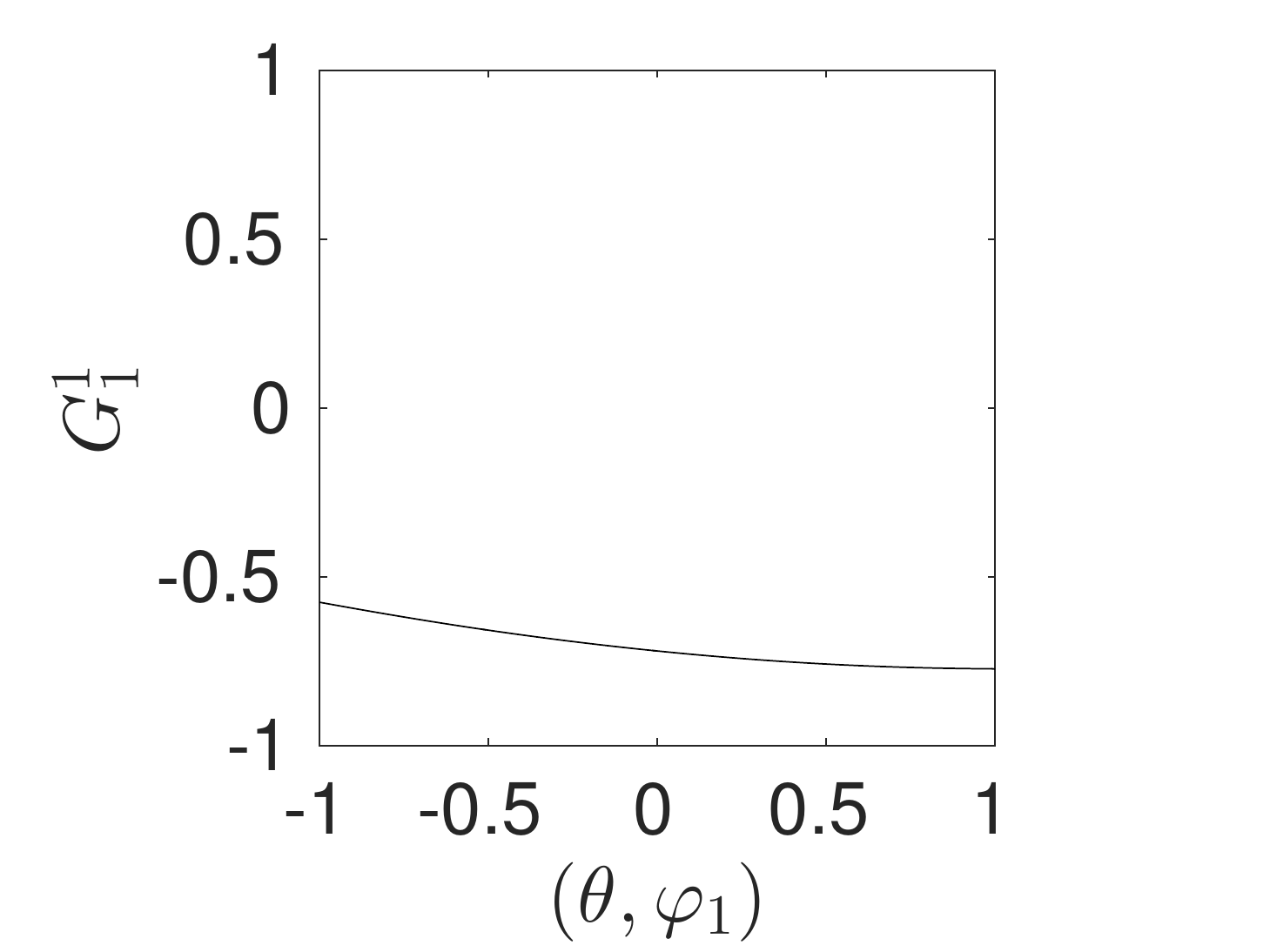}
\includegraphics[height=3.5cm]{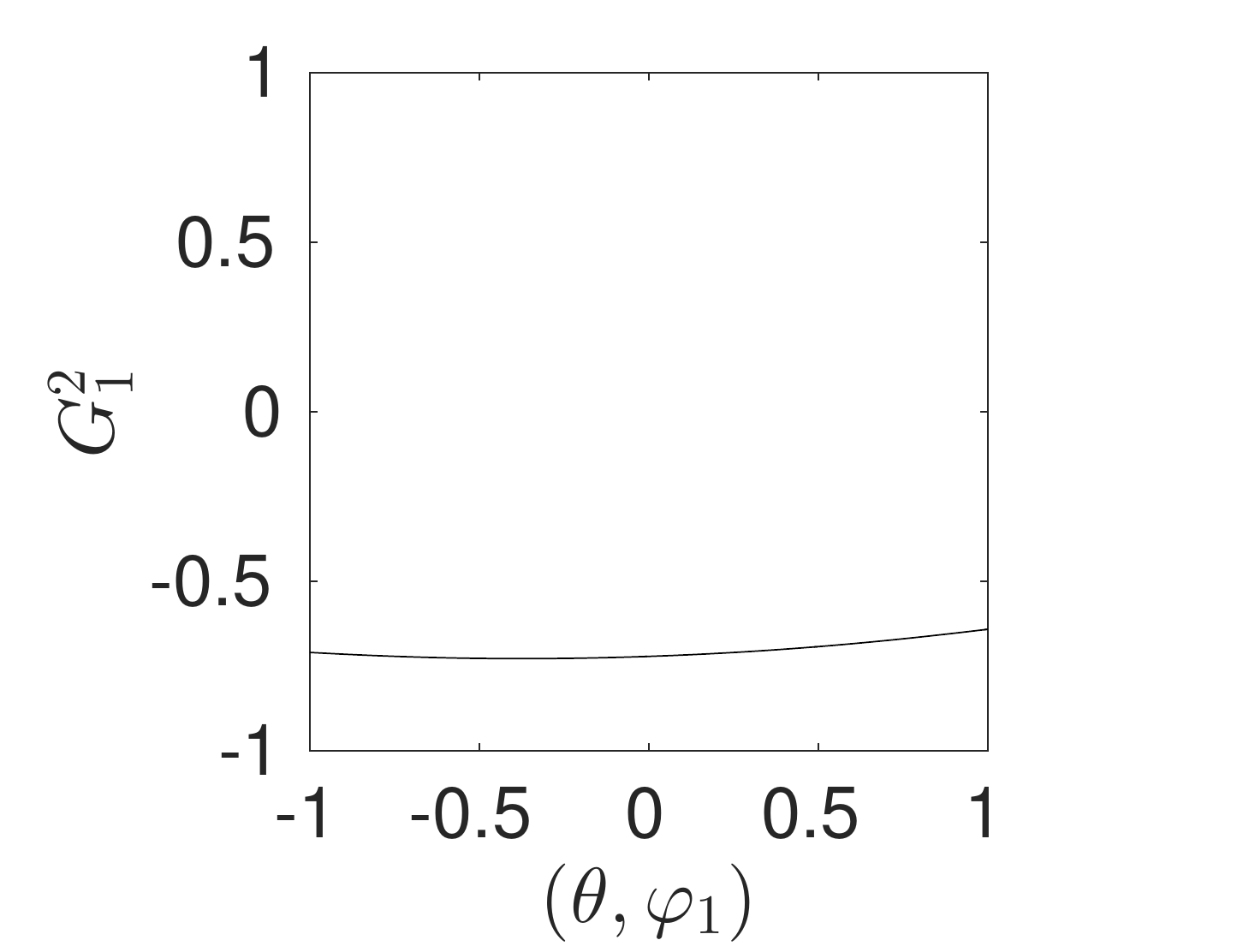}
\includegraphics[height=3.5cm]{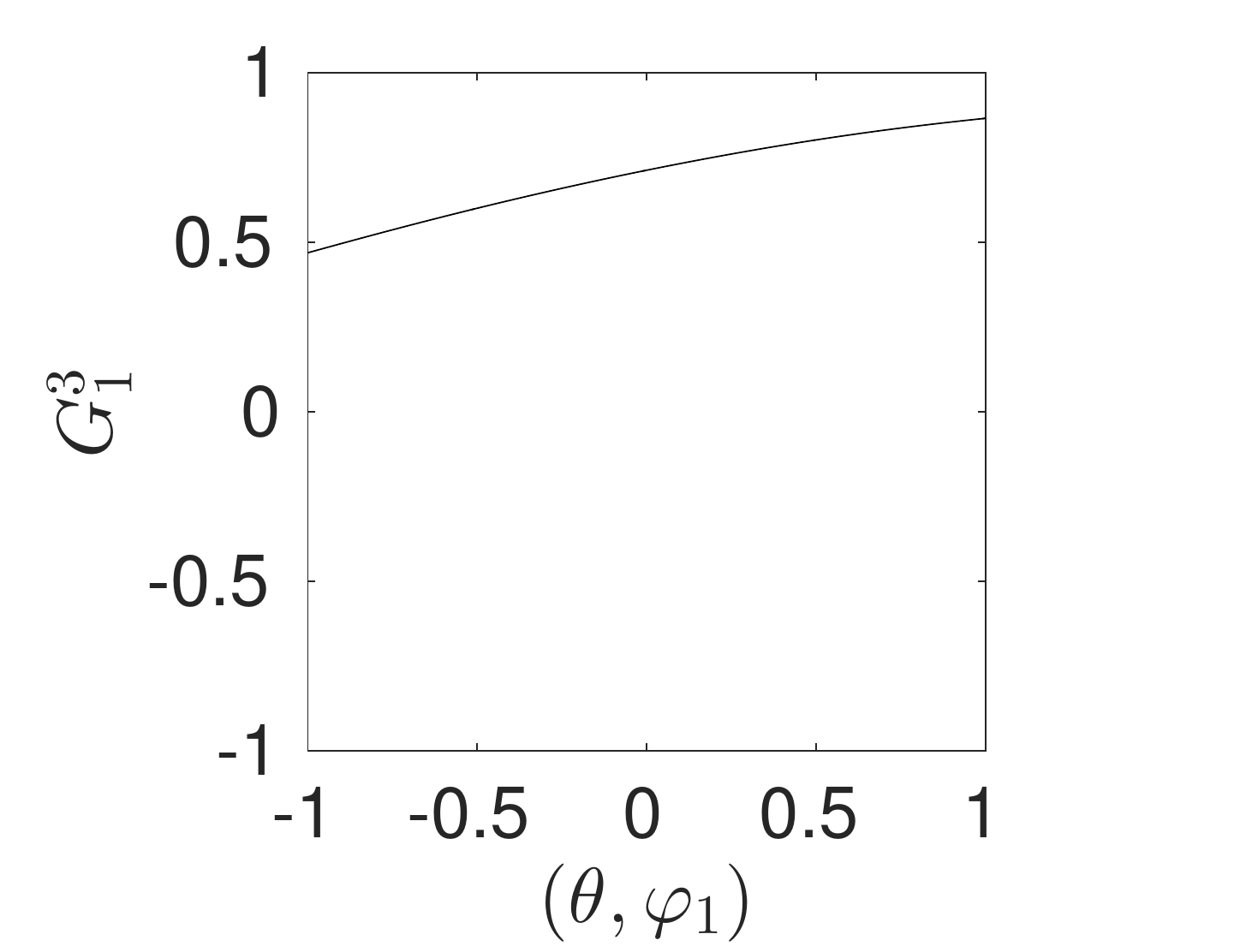}
}
\centerline{
\includegraphics[height=3.5cm]{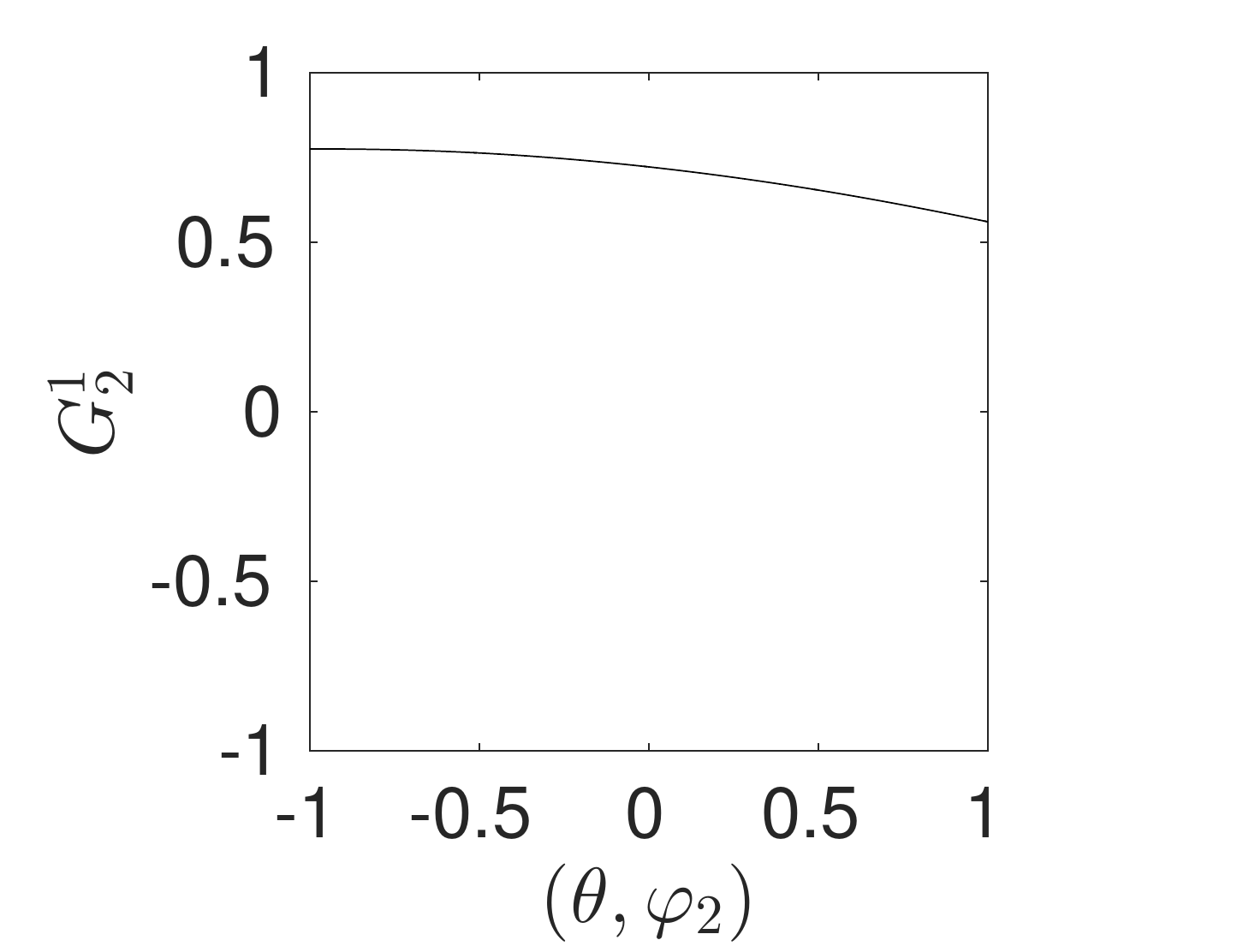}
\includegraphics[height=3.5cm]{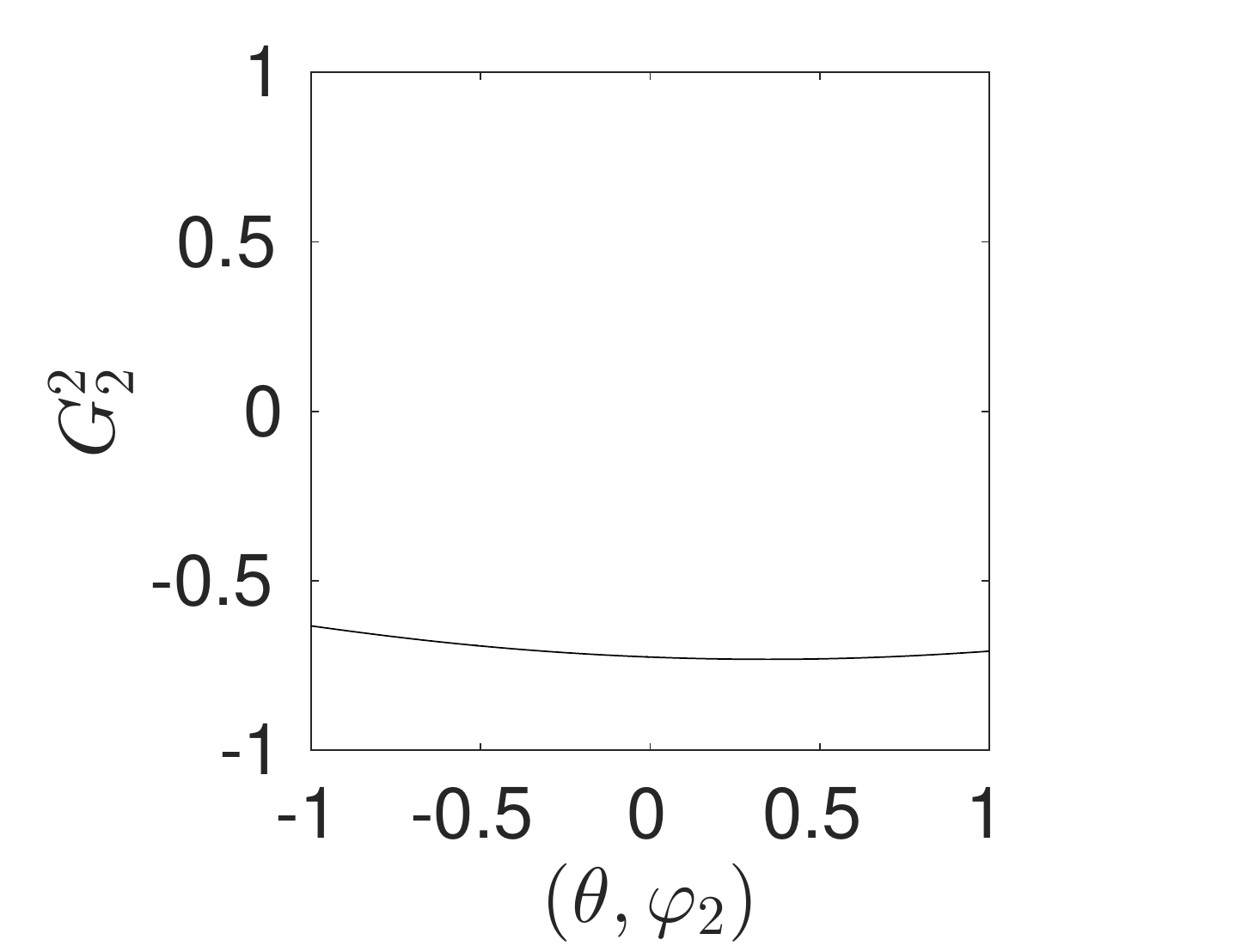}
\includegraphics[height=3.5cm]{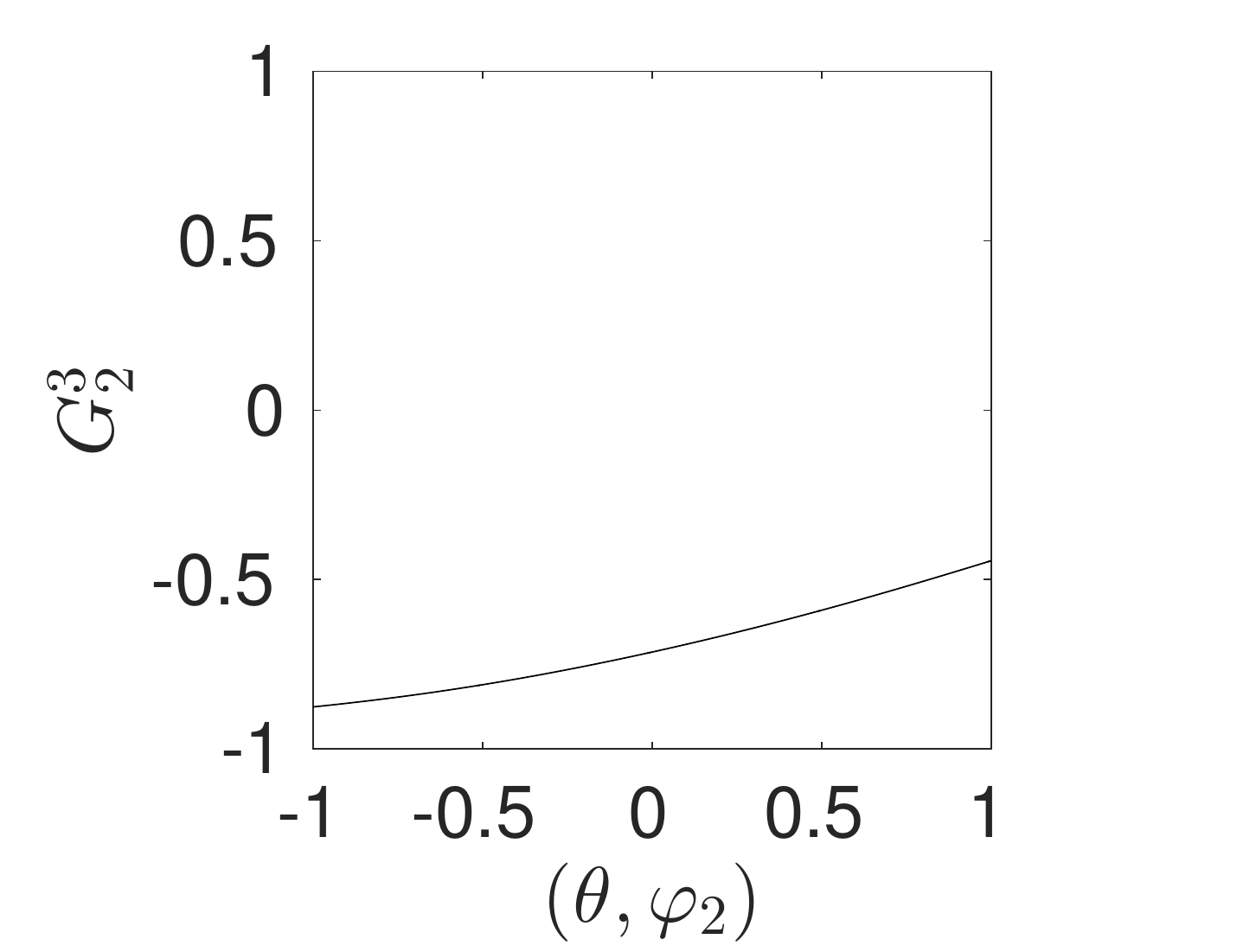}
}
\centerline{
\includegraphics[height=3.5cm]{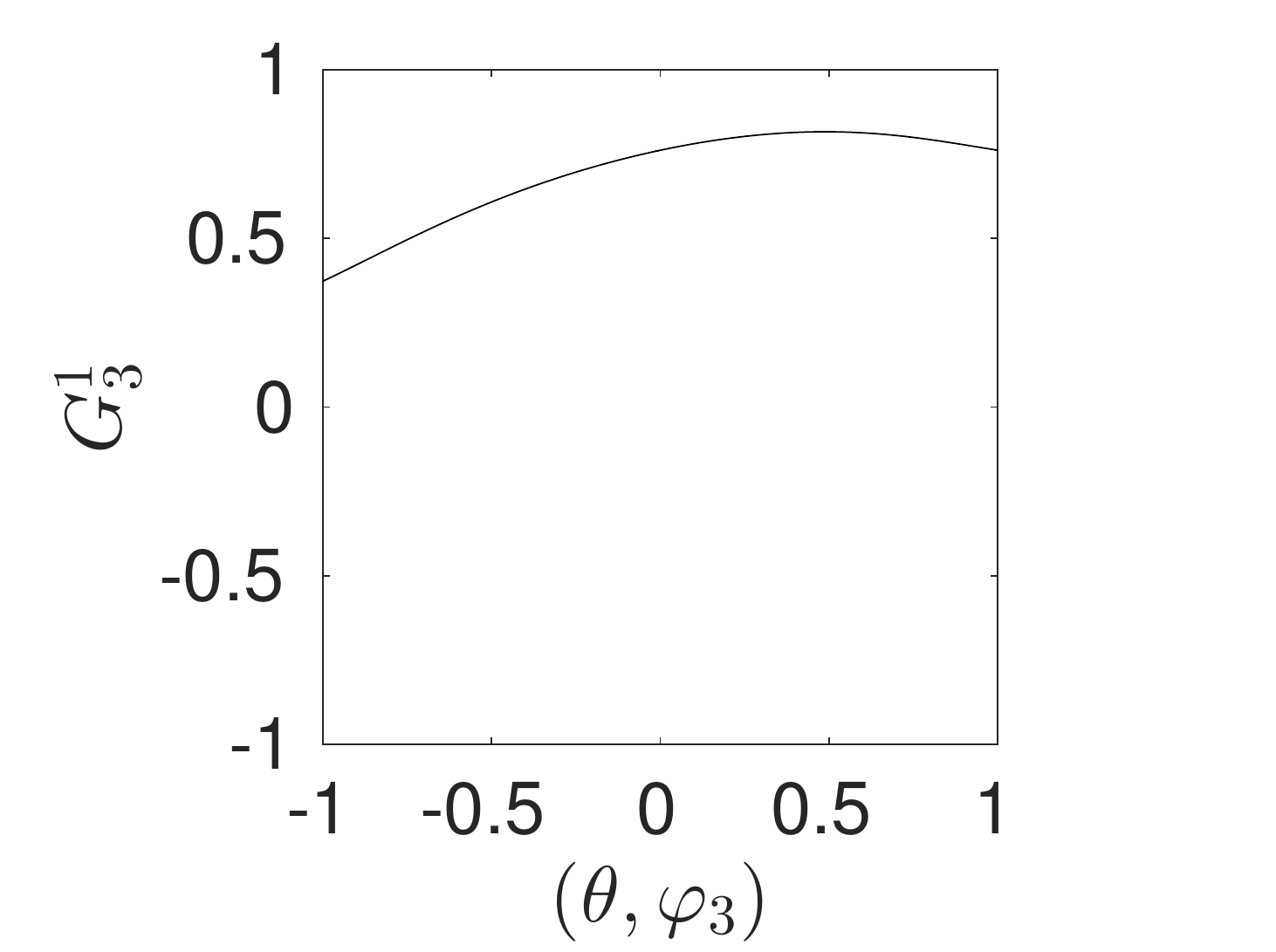}
\includegraphics[height=3.5cm]{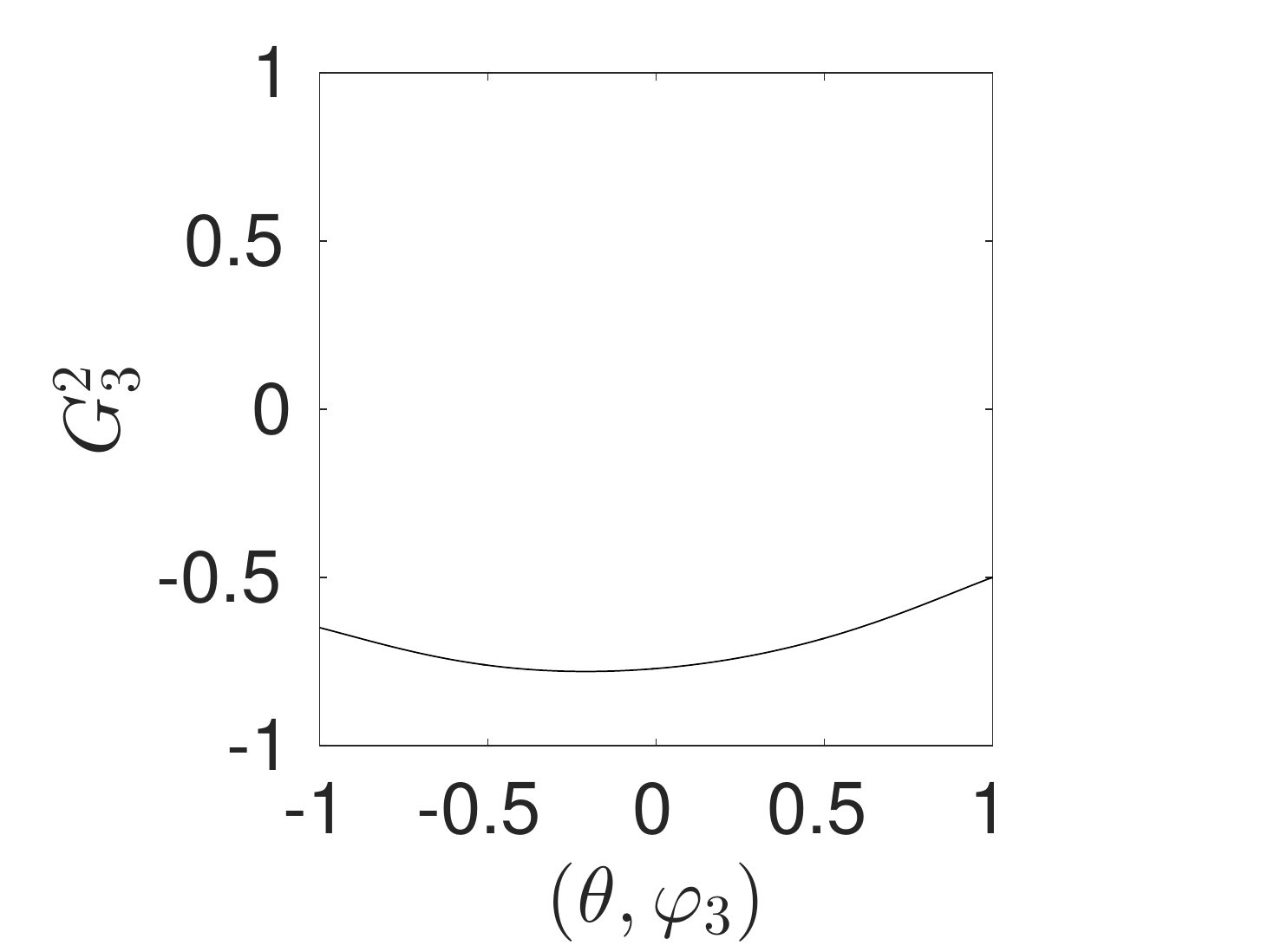}
\includegraphics[height=3.5cm]{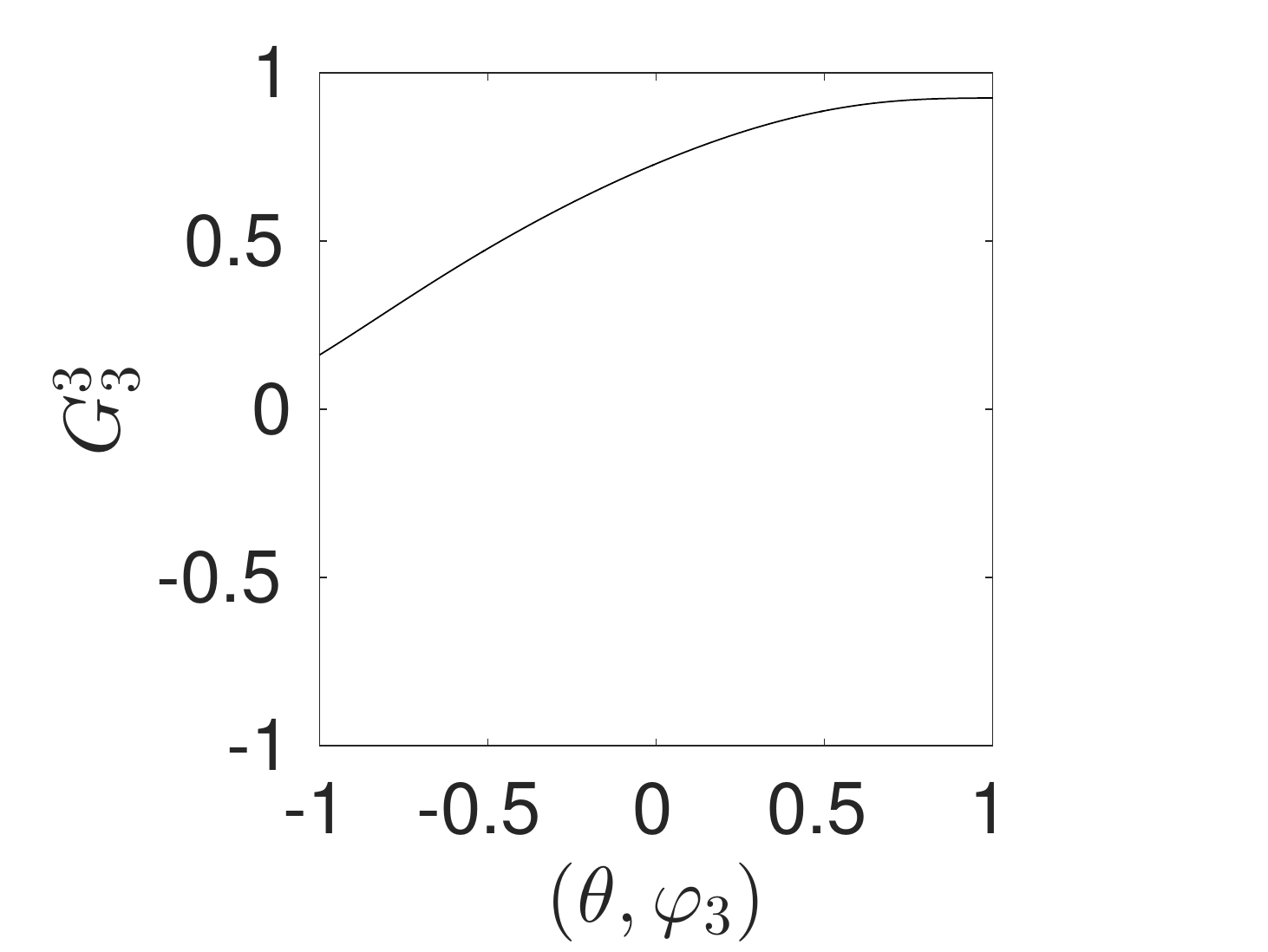}
}

\caption{Tensor components of the sine functional \eqref{SINEF} 
versus $(\theta,\varphi_k)$.}
\label{fig:G}
\end{figure}

Consider the nonlinear functional 
\begin{equation}
F([\theta])=\sin\left((K_1,\theta)\right),\qquad 
(K_1,\theta)=\int_0^{2\pi}K_1(x)\theta(x)dx,
\label{DINE}
\end{equation}
where $K_1(x)$ is as in equation \eqref{K1}. 
We represent $F([\theta])$ in the function space
\begin{equation}
D_m=\textrm{span}\{\varphi_1,...,\varphi_m\}
\end{equation}
spanned by a finite-dimensional orthonormal basis. 
This yields the cylindrical representation (see Section \ref{sec:tensor})
\begin{equation}
f(a_1,...,a_m) = \sin\left(\sum_{k=1}^m s_k(\varphi_k,\theta)\right),
 \qquad s_k=\int_0^{2\pi}K_1(x)\varphi_k(x)dx.
\label{SINEF}
\end{equation}
In Figure \ref{fig:porterSSEsine} we study the accuracy of 
second-order Porter's polynomial functionals and 
canonical tensor decomposition \eqref{functional-SSE} in 
approximating the functional \eqref{SINEF}. 
Specifically, the polynomial functionals are constructed by using 
the set of nodes $\widetilde{S}^{(m+1)}_2$ defined in equation 
\eqref{SNq2}, where $a_{i_j}$ are sampled 
at Gauss-Hermite sparse grids of level 5. 
We recall that this set yields a 
rank-deficient matrix \eqref{matH}, which requires 
More-Penrose pseudoinversion 
(see Section \ref{sec:Hopf polynomial}). 
Correspondingly, the polynomial functionals do not 
interpolate \eqref{SINEF}. 
The canonical tensor decomposition, on the other hand, 
is based on Legendre polynomials of order $Q=6$. 
The functionals $G^{l}_k((\theta,\varphi_k))$ 
are shown in Figure \ref{fig:G}.

\subsubsection{Functional Derivatives}
\begin{figure}[!ht]
\centerline{\hspace{1.3cm}Analytical \hspace{3.8cm} Canonical Tensor Decomposition}
\vspace{0.cm}
\centerline{\hspace{-.7cm}
 \includegraphics[height=4.3cm]{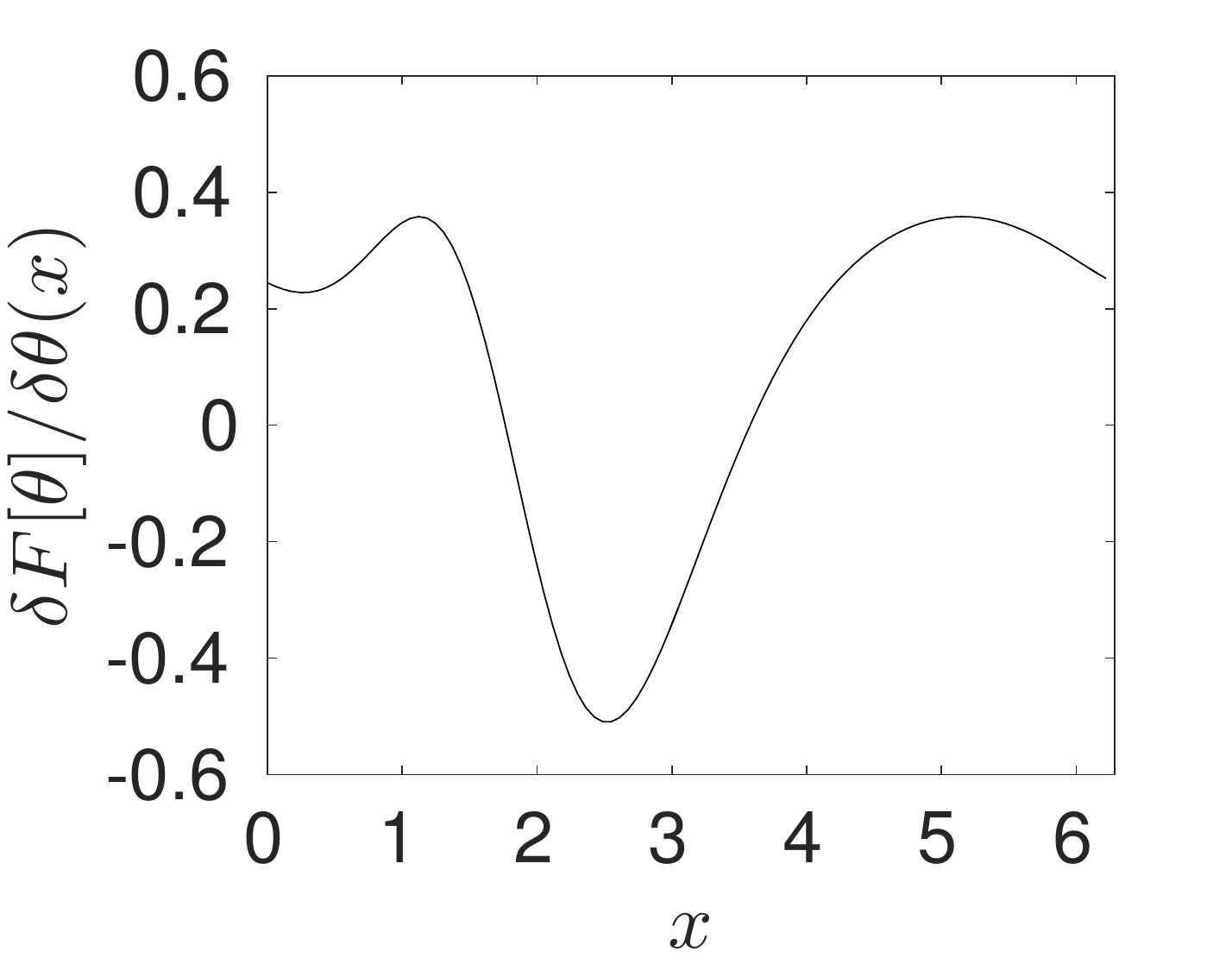}\hspace{1.8cm}
\includegraphics[height=4.3cm]{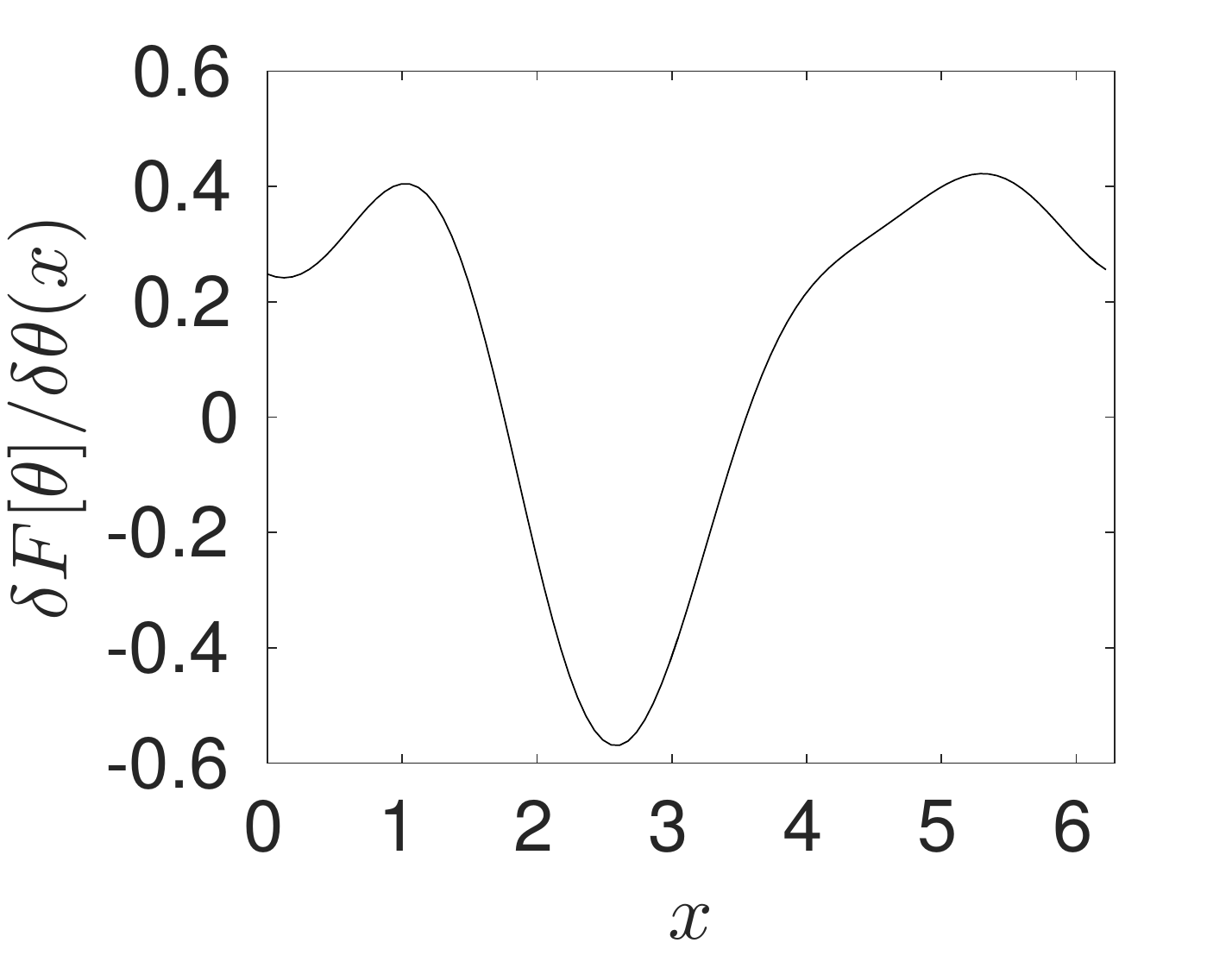}
 }

\centerline{
 \includegraphics[height=5.cm]{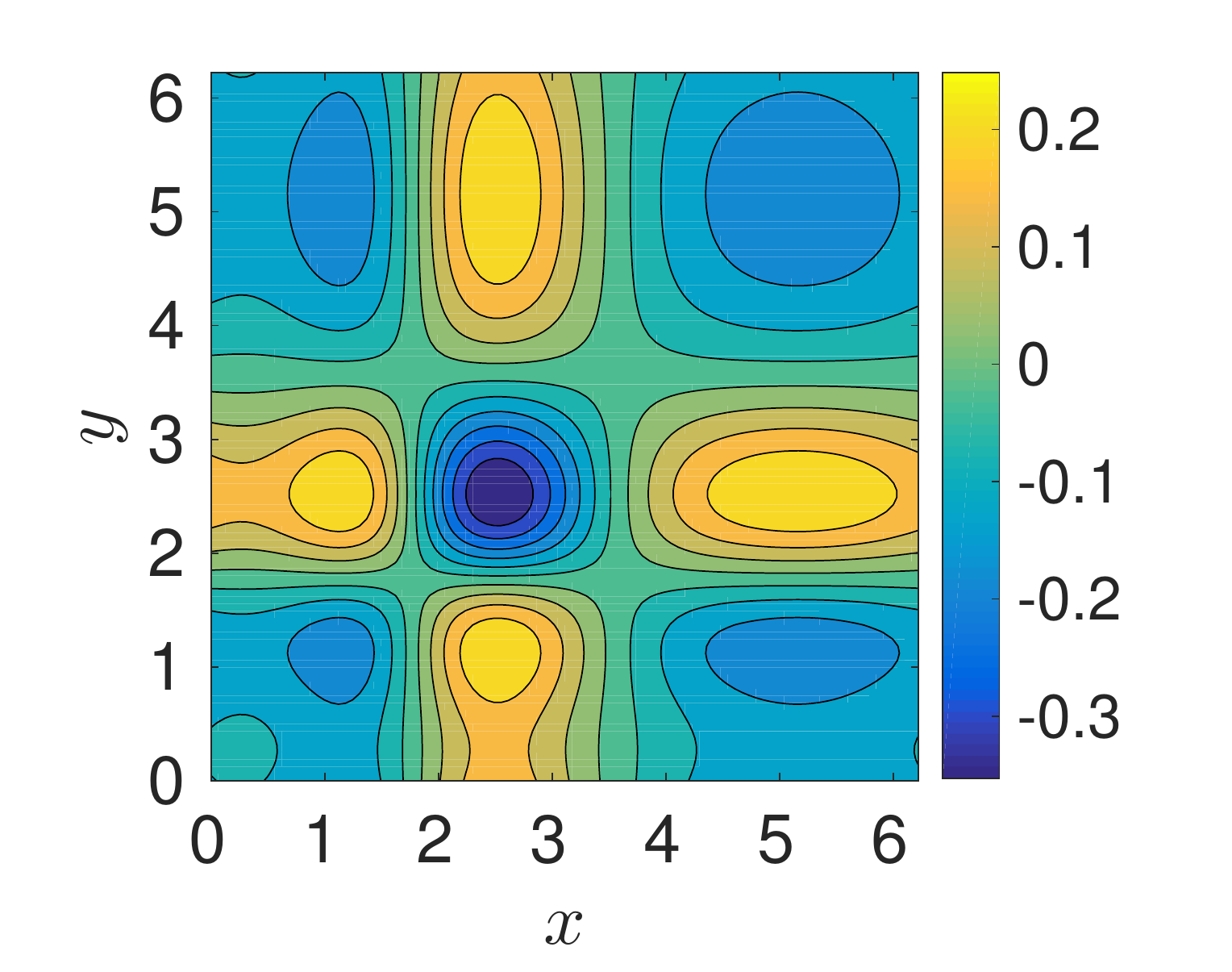}\hspace{1.cm}
 \includegraphics[height=5.cm]{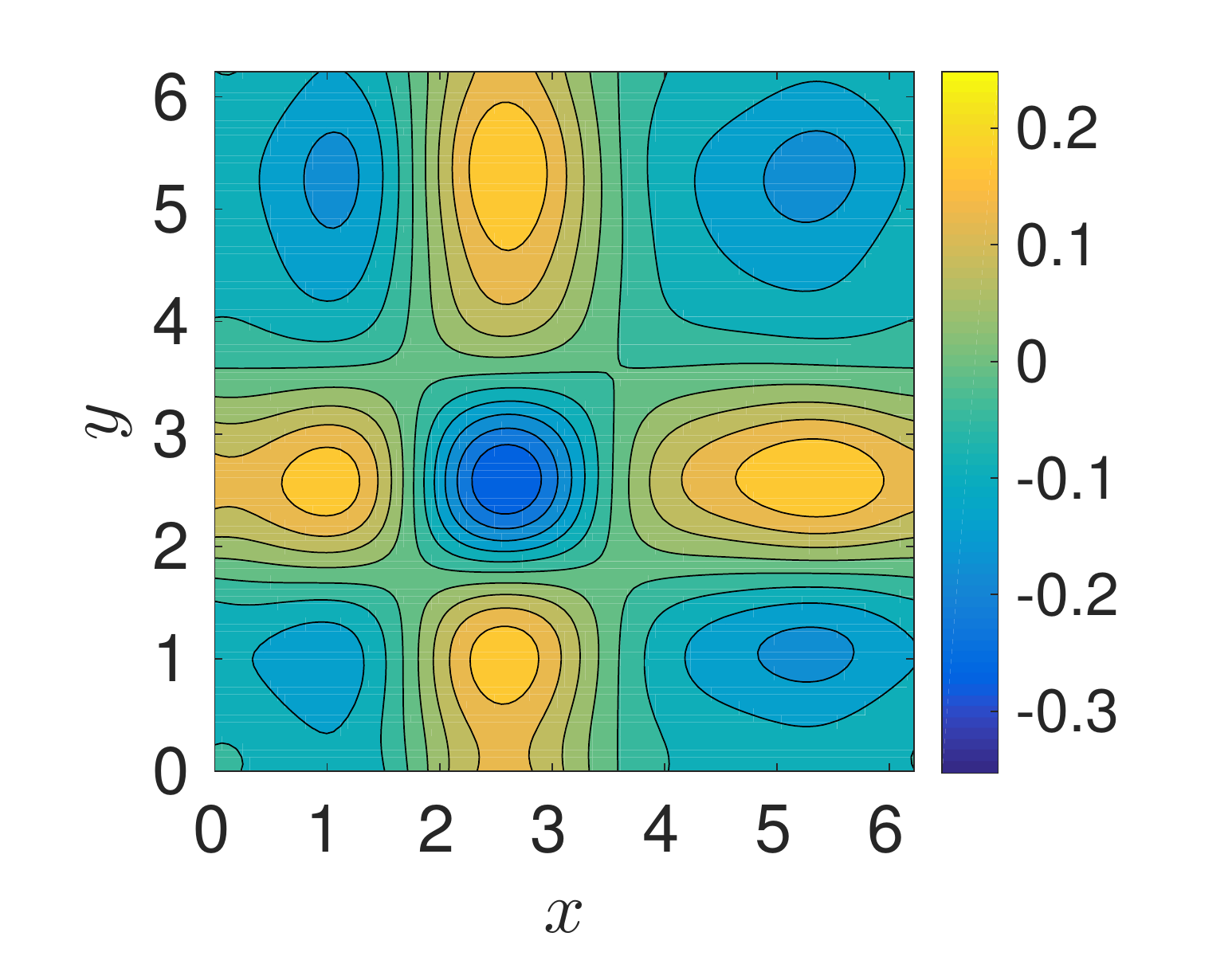}
 }
 
\caption{Sine functional \eqref{DINE}. First- and second-order 
functional derivatives evaluated at 
$\theta(x)=9(1+\sin(x)+\sin(2x))/10$. 
Specifically, we compare the analytical results 
\eqref{d1sin} and \eqref{d2sin} versus 
results obtained by canonical tensor decomposition 
(equations \eqref{ff1} and \eqref{ff5}).}
\label{fig:sineFD1}
\end{figure}

The first- and second-order functional derivatives of \eqref{DINE}
are easily obtained as 
\begin{equation}
\frac{\delta F([\theta])}{\delta \theta(x)}=K_1(x)\cos((K_1,\theta)),
\label{d1sin}
\end{equation}
\begin{equation}
\frac{\delta^2 F([\theta])}{\delta \theta(x)\delta \theta(y)}=
-K_1(x)K_1(y)\sin((K_1,\theta)).
\label{d2sin}
\end{equation}
Note that for each $\theta(x)$, such functional derivatives 
are basically a rescaled version of 
the functions $K_1(x)$ and $K_1(x)K_1(y)$.  
In Figure \ref{fig:sineFD1} 
we compare the exact functional derivatives versus those obtained 
by the canonical tensor decomposition with separation rank $r=4$. 
Recall that the functional derivatives can be approximated in the 
space of cylindrical functionals as (see Section \ref{sec:tensor}) 
\begin{equation}
\frac{\delta {F}([\theta])}{\delta\theta(x)}=\sum_{k=1}^m
\frac{\partial f}{\partial a_k}\varphi_k(x),
\label{SSEsinefd1}
\end{equation}
\begin{equation}
\frac{\delta^2 {F}([\theta])}{\delta\theta(x)\delta\theta(y)}=\sum_{k,j=1}^m
\frac{\partial^2 f}{\partial a_k\partial a_j}\varphi_k(x)\varphi_j(y).
\label{SSEsinefd2}
\end{equation}
In the case of canonical tensor expansions, 
$\partial f/\partial a_k$ and 
$\partial^2 f/\partial a_j\partial a_k$ are 
defined in \eqref{pd1} and \eqref{pd2}, respectively.

\section{Numerical Results: Functional Differential Equations}
\label{sec:numerical results functional equations}

Computing the numerical solution to a functional differential 
equation is a long standing open problem in mathematical 
physics. In this Section we address this problem with reference 
to linear functional equations in the 
form \eqref{linfde00}. In particular, we study the initial value 
problem for a prototype functional advection-reaction equation. 
We also develop the numerical discretization of the 
Navier-Stokes-Hopf functional equation 
(see Section \ref{sec:motivation}), and discuss 
its computational complexity.
{
\color{r}
\subsection{Advection-Reaction Functional Differential Equation}
\label{sec:ADVR}
Consider the following advection-reaction functional differential 
equation
\begin{equation}
\frac{\partial F([\theta],t)}{\partial t}+\int_{0}^{2\pi} \theta(x)\frac{\partial }{\partial x}\left( \frac{\delta F([\theta],t)}{\delta \theta(x)}\right)dx=H([\theta],t)
\label{FDE3}
\end{equation}
evolving from the initial condition 
\begin{equation}
F([\theta],0)=F_0([\theta]).
\label{fffer}
\end{equation}
Here $H([\theta],t)$ is a given functional reaction term. 
We assume that $D(F)$ (the domain of the functional
$F$) is a suitable space of functions\footnote{\color{r}We have seen in 
Section \ref{sec:FDEs approximation} that the solution 
to the initial value problem \eqref{FDE3}-\eqref{fffer} 
is strongly dependent on the choice of 
function space $D(F)$.}, e.g., the space 
of periodic functions in $[0,2\pi]$ or the space of 
infinitely differentiable functions in $[0,2\pi]$ such that 
$\theta(0)=0$. Evaluation of 
$F([\theta],t)$ and $H([\theta],t)$ in 
the finite-dimensional subspace 
\begin{equation}
D_m = \left\{\theta(x)\in D(F) \,\left|\, \theta(x)=\sum_{k=1}^m a_k \varphi_k(x)  \right.  
\right\}\subseteq D(F),
\label{Dm}
\end{equation}
where $\{\varphi_1,...,\varphi_m\}$ is an orthonormal basis,
yields the following multivariate functions
\begin{align}
f(a_1,...,a_m,t)=F([\theta],t) \qquad h(a_1,...,a_m,t)=H([\theta]), \qquad \theta\in D_m.
\end{align}
We also recall that the functional derivative $\delta F([\theta],t)/\delta \theta(x)$  can be expressed in $D_m$ as (see Eq. \eqref{gg2}) 
\begin{equation}
 \frac{\delta F([\theta],t)}{\delta \theta(x)}=\sum_{j=1}^m 
\frac{\partial f}{\partial a_j}\varphi_j(x).
\end{equation}
A substitution of these expression back into \eqref{FDE3} 
yields the following initial value problem for a multivariate first-order PDE
\begin{equation}
\frac{\partial f}{\partial t}+\sum_{j=1}^m\left(\sum_{k=1}^m 
C_{jk}a_k\right)\frac{\partial f}{\partial a_j}=h,
\qquad f(a_1,...,a_m,0)=f_0(a_1,...,a_m),
\label{PDE-advR}
\end{equation}
where 
\begin{equation}
C_{jk}=\int_{0}^{2\pi} \varphi_k(x)\frac{d\varphi_j(x)}{dx}dx 
\label{matC}
\end{equation}
The entries $C_{ij}$ depend on the 
choice of $D(F)$, and correspondingly $D_m$. For instance, 
if we assume that $D(F)$ is the space of infinitely differentiable 
periodic functions in $[0,2\pi]$ then the matrix \eqref{matC} 
is centro-skew-symmetric. 

\vs
\noindent
{\em Example 1:} 
Another example of advection FDE is
\begin{equation}
\frac{\partial F}{\partial t}+\int_{0}^{2\pi} \theta(x)\frac{\partial }{\partial x} \left[\frac{\delta F}{\delta \theta(x)}\right]dx=
\int_{0}^{2\pi} \theta(x)\frac{\partial^2 }{\partial x^2}\left[ \frac{\delta F}{\delta \theta(x)}\right]dx.
\label{FDE3_bis}
\end{equation}
The discrete form such equation is 
\begin{equation}
\frac{\partial f}{\partial t}+\sum_{k,j=1}^m 
a_k\left(C_{jk}-F_{jk}\right)\frac{\partial f}{\partial a_j}=0,
\end{equation} 
where $C_{jk}$ is defined in \eqref{matC} while  
\begin{equation}
F_{kj}= \int_{0}^{2\pi} \varphi_k(x)\frac{d^2\varphi_j(x)}{dx^2}dx.
\end{equation} 

\subsubsection{Analytical Solution}
\label{sec:ADVanalyticalSolution}
The analytical solution to the initial value 
problem \eqref{PDE-advR} can be 
computed by using the method of 
characteristics \cite{Rhee}. To this end, let $\bm C$ 
be the matrix with entries \eqref{matC}, and $\bm a$ 
be the vector of coordinates $(a_1, ..., a_m)$. 
Then the solution to \eqref{PDE-advR} is  
\begin{align}
f(\bm a,t)=f_0\left(e^{-t\bm C }\bm a\right) + 
\int_0^t h\left(e^{s\bm C }\bm a\right)ds.
\label{solutionF}
\end{align}
In the particular case where the reaction term  $h$
is zero we obtain
\begin{equation}
f(a_1,...,a_m,t) = f_0\left(a_1(t), ..., 
a_m(t)\right),
\label{solutionFF}
\end{equation}  
where the (inverse) flow map $a_i(t;a_1,...,a_m)$ is 
given by\footnote{We recall that the matrix exponential 
appearing in \eqref{gg1i} can be represented as 
\begin{equation}
e^{-t\bm C}= \bm U^T e^{-t\bm \Lambda} \bm U,
\end{equation}
where $\bm U$ is the matrix of eigenvectors 
of $\bm C$ (columnwise) and $\bm \Lambda$ is the 
diagonal matrix of eigenvalues.} 
\begin{equation}
a_i (t;a_1,...,a_m)=\sum_{j=1}^m Z_{ij}(t)a_j,\qquad 
\bm Z(t)=e^{-t\bm C}.
\label{gg1i}
\end{equation}
The solution to the functional equation \eqref{FDE3} 
can be obtained by taking the continuum limit of 
\eqref{solutionFF}, i.e., by sending $m$ to infinity. 
This yields 
\begin{equation}
F([\theta],t)=F_0([\theta(x,t)]),
\label{FDEsolution}
\end{equation}
where
\begin{equation}
\theta(x)=\sum_{k=1}^\infty a_k\varphi_k(x),\qquad \theta(x,t)=
\sum_{k=1}^\infty a_k(t)\varphi_k(x).
\label{theta_xt}
\end{equation}
The coefficients $a_k(t)$ are obtained by applying the 
semigroup $Z_{ij}(t)$ to $a_k$ (see equation \eqref{gg1i}).
The analytical expression \eqref{FDEsolution} can be 
written more rigorously in terms of the action of 
a semigroup $U(t)$ \cite{Engel} to $\theta(x)$, i.e., 
\begin{equation}
F([\theta],t)=F_0([U(t)\theta(x)]).
\label{semiG}
\end{equation}
As we will see, such semigroup defines a translation 
in the space of functions $D(F)$. Such translation can be 
generated by rotations or contractions, depending on 
the space of functions $D(F)$ we consider. 
Hereafter we discuss this matter in more detail.

\paragraph{Periodic Function Spaces}
Assume that the function 
space $D(F)$ (domain of the solution functional $F$), 
is the space of infinitely differentiable periodic 
functions in $[0,2\pi]$.  In this case, the matrix $C_{ij}$ 
defined in \eqref{matC} is  skew-symmetric, 
thanks to the periodicity of $\varphi_k$
(just integrate \eqref{matC} by parts ). 
Therefore, by the spectral theorem, $C_{ij}$ it has purely 
imaginary eigenvalues $\lambda_k=k i$, 
$k\in \mathbb{Z}$.
Since $\bm C$ is skew-symmetric we have 
that  $\exp[-t\bm C]$ is orthogonal, i.e., it defines 
an isometry in $\mathbb{R}^m$. Such isometry 
generates a translation in the space of periodic 
functions with group velocity equal 
to one. In other words, we have 
\begin{align} 
\theta(x,t) =\theta(x-t). 
\end{align}
Therefore, the analytical solution to the functional differential 
equation \eqref{FDE3} (with $H=0$) in the space of periodic functions $D(F)$ is 
\begin{equation}
F([\theta],t)=F_0([\theta(x-t)]).
\label{generalsolution}
\end{equation}
From this equation, we see that if $F_0$ is invariant under 
translation, i.e., $F_0([\theta(x-t)])=F_0([\theta(x)])$ then 
the solution functional is constantly equal to 
the initial condition $F_0([\theta])$. This is discussed in more detail 
the following two examples. 

\vs
\noindent
{\em Example 1:}
Consider the initial condition 
\begin{equation}
F_0([\theta])=\sin\left(\int_{0}^{2\pi}\theta(x)dx\right).
\label{FDE0i}
\end{equation}
We have seen in Section \ref{sec:approximability_of_functionals} that this 
nonlinear functional is approximable by a one-dimensional  
function relative to the standard Fourier (modal) 
basis in $[0,2\pi]$. Specifically, we obtained 
\begin{equation}
f_0(a_0,...,a_m)= \sin\left(\sqrt{2\pi} a_0\right),
\end{equation}
independently on $m$. This implies that the analytical solution 
to the initial value problem \eqref{PDE-advR} is 
\begin{equation}
f(a_0,...,a_m,t)=\sin\left(\sqrt{2\pi} \sum_{k=0}^m 
Z_{0k}(t) a_k\right), \quad \textrm{where}\quad 
Z_{0k}(t)=\left[e^{-t \bm C}\right]_{0k}.
\label{as1}
\end{equation}
It can be shown that  if we sort the basis 
elements as in \eqref{D2m} then $Z_{0k}(t)=\delta_{k0}$. 
Therefore, the solution to \eqref{FDE3}-\eqref{FDE0i} (with $H=0$) 
in the space of periodic functions is 
\begin{equation}
F([\theta],t)=F_0([\theta]);
\label{thesolutionfunctional}
\end{equation}
i.e., the constant functional. 
In general, the solution to a FDE is independent on the 
way represent the test function space $D(F)$. 
Thus, it shouldn't be surprising that we obtain exactly 
the same result if we consider a finite-dimensional 
expansion in terms of nodal trigonometric polynomials. 
In this case, the initial condition functional can be 
discretized as 
\begin{equation}
f_0(a_0,...,a_m)= \sin\left(\eta \sum_{k=0}^m a_k\right),\quad 
\textrm{where}\quad \eta = \int_{0}^{2\pi} \varphi_k(x) dx = \left(\frac{2\pi}{m+1}\right)
\end{equation}
and the solution is
\begin{equation}
f(a_0,...,a_m,t)= \sin\left(\eta \sum_{k=0}^m a_k(t)\right), 
\label{as111}
\end{equation}
where $a_k(t)$ is defined in \eqref{as1}. 
Since $\exp(-t\bm C)$ is an orthogonal matrix we have 
that 
\begin{equation}
\sum_{k=0}^m a_k(t) = \sum_{k=0}^m a_k.
\label{cdfg}
\end{equation}
In the limit $m\rightarrow \infty$ \eqref{as111} and \eqref{cdfg} 
imply \eqref{thesolutionfunctional}. Such result can also 
be obtained by directly noting that if $\theta$ is 
periodic in $[0,2\pi]$ then 
\begin{equation}
\int_0^{2\pi} \theta(x)dx = \int_0^{2\pi} \theta(x-t)dx\qquad 
\textrm{for all $t\in \mathbb{R}$.}
\label{symF}
\end{equation}
Substituting this into \eqref{FDE0i} and  \eqref{generalsolution} yields \eqref{thesolutionfunctional}.

\vs
\noindent
{\em Example 2:}
Consider the initial condition 
\begin{equation}
F_0([\theta])=\exp\left[-\int_{0}^{2\pi}\theta(x)^2dx\right].
\label{expFun}
\end{equation}
where $\theta(x)$ is in the space of infinitely differentiable 
periodic functions in $[0,2\pi]$. Represent $\theta(x)$ 
in a finite-dimensional space spanned by 
any orthonormal periodic basis in $[0,2\pi]$. 
This yields the following multivariate function corresponding to $F_0([\theta])$ (see also Section \ref{sec:approximability_of_functionals})
\begin{equation}
f_0(a_1,...,a_m)= \prod_{k=1}^m e^{-a_k^2}.
\label{iic1}
\end{equation}
The analytical solution to the multivariate 
PDE \eqref{PDE-advR} with initial condition \eqref{iic1} and 
$h=0$ is (see equation \eqref{solutionFF})  
\begin{equation}
f(a_1,...,a_m,t)=\prod_{k=1}^m 
e^{-a_k(t)^2},
\label{abs1}
\end{equation}
where
\begin{equation}
a_k(t)=\sum_{j=1}^m Z_{kj}(t) a_j\qquad \textrm{and}
\qquad Z_{kj}(t)=\left[e^{-t \bm C}\right]_{kj}.
\label{abs2}
\end{equation}  
In a collocation setting where the the 
basis functions are normalized nodal trigonometric 
polynomials it is easy to see that 
\begin{equation}
f(a_1,...,a_m,t) = f(a_1,...,a_m) \qquad \textrm{for all $t\geq 0$}
\end{equation}
In fact, $a_j$ are the rescaled values of 
$theta(x)$ at node $x_j\in[0,2\pi]$,
the rescaling coefficient being the norm of the 
trigonometric polynomial. 
Since $\theta(x)$ is periodic, when $t$ 
increases we have that values  $a_j$ are 
just shifted to another location in 
$[0,2\pi]$, leaving the product in \eqref{abs1} constant.    
In other words, the solution to the FDE \eqref{FDE3} (with $H=0$)
corresponding to the initial condition \eqref{expFun} is 
\begin{equation}
F([\theta],t) = F_0([\theta]), 
\label{thesolutionfunctional1}
\end{equation}
i.e., the constant functional.
As before,  this  result can be obtained by 
noting that if $\theta(x)$ is periodic in $[0,2\pi]$ then 
\begin{equation}
\int_0^{2\pi} \theta(x)^2dx = \int_0^{2\pi} \theta(x-t)^2dx\qquad 
\textrm{for all $t\in \mathbb{R}$.}
\label{symF1}
\end{equation}
A substitution of \eqref{symF1}  into \eqref{expFun} yields \eqref{thesolutionfunctional1}.

\paragraph{Polynomial Function Spaces}
To obtain non-trivial solutions to equation \eqref{FDE3}, 
we consider the following space of functions 
\begin{equation}
D(F) = \left\{\theta \in {C}^{\infty}([0,2\pi])\, 
|\,\theta(0)=0\right\}
\label{Fspace9}
\end{equation}
and the initial condition \eqref{expFun}. 
We generate an orthonormal polynomial basis spanning $D(F)$ 
by orthonormalizing the modified Chebyshev 
basis $x T_k(x/\pi-1)$ in $[0,2\pi]$ 
though the Gram-Schmidt procedure. 
The basis functions we obtain are shown 
in Figure \ref{fig:basis_eigen}.
\begin{figure}[t]
\centerline{(a) \hspace{5cm} (b)\hspace{5cm}(c)}
\centerline{
\includegraphics[height=5cm]{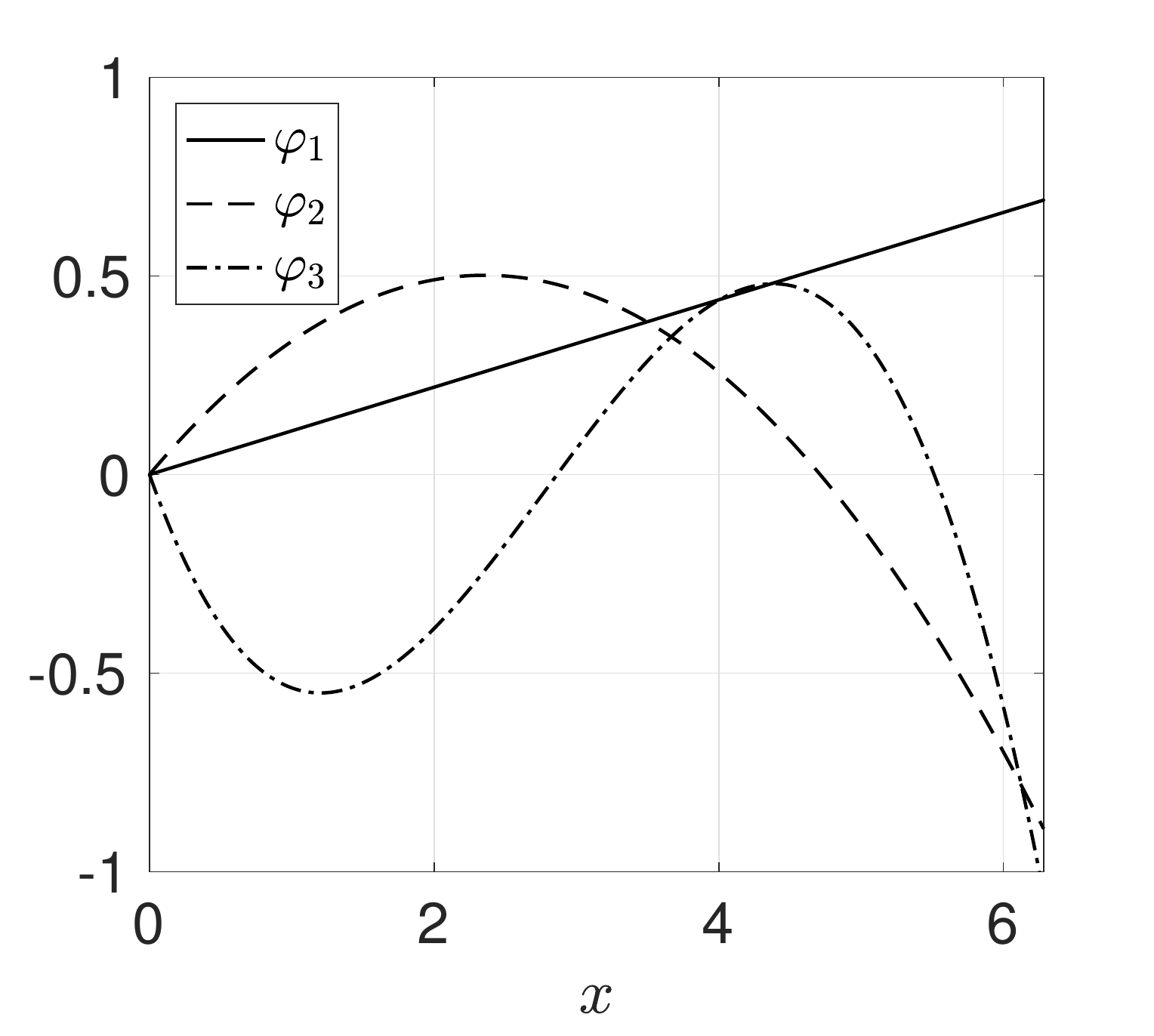}
\includegraphics[height=5cm]{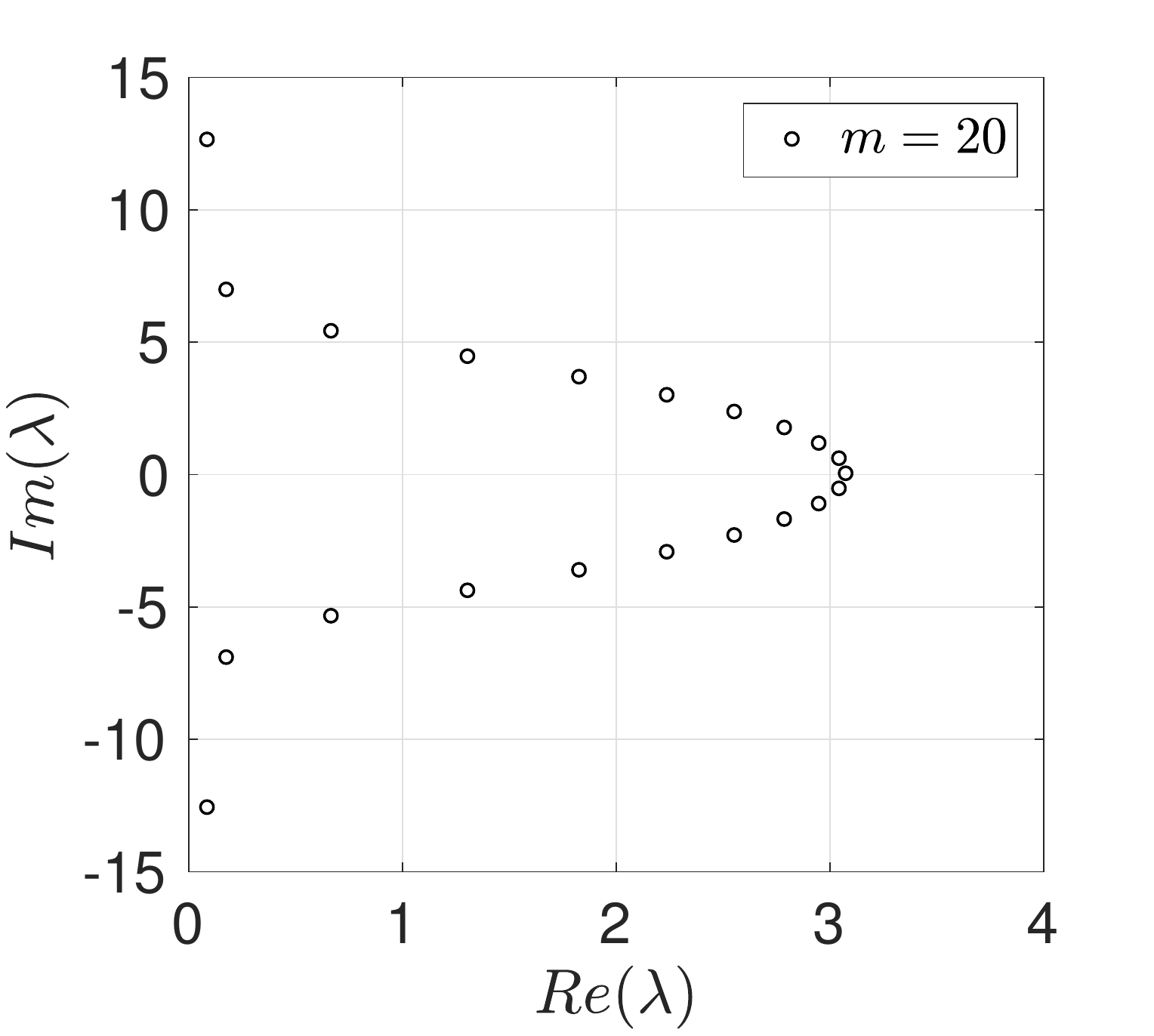}
\includegraphics[height=5cm]{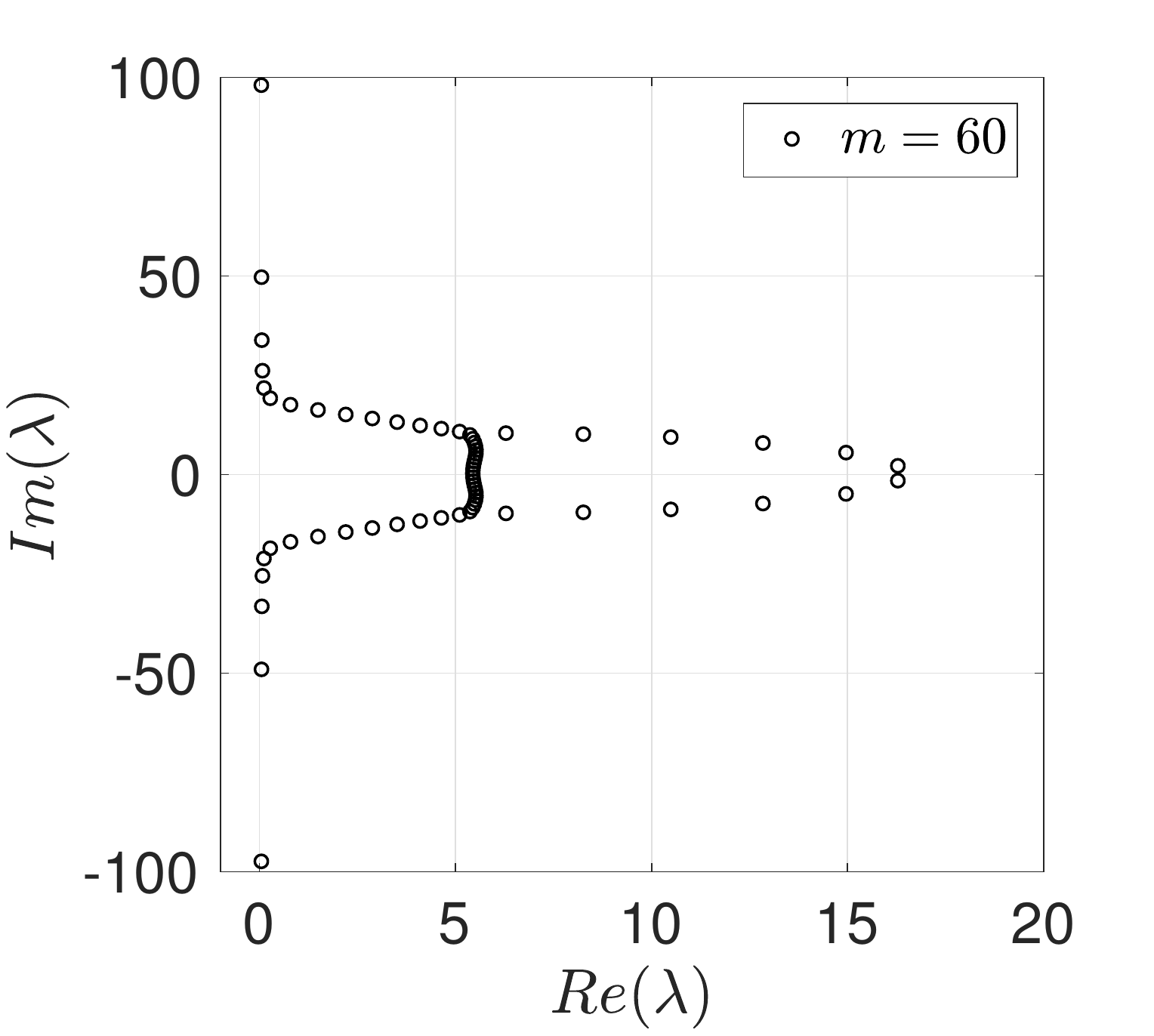}
}
\caption{\color{r}(a) First three orthonormal polynomial basis functions 
spanning the function space  \eqref{Fspace9}; (b)-(c) spectrum of 
the matrix \eqref{matC} for different number of variables $m$.}
\label{fig:basis_eigen}
\end{figure}
This allows us to define the finite-dimensional function space
\begin{equation}
D_m =\left\{\theta\in C^{\infty}([0,2\pi])\,|\, \theta(x) =\sum_{k=1}^m a_k\varphi_k(x) \right\}\subseteq D(F), 
\label{finite_dim_FS}
\end{equation}
where $\varphi_k$ are the orthonormal basis functions
shown in Figure \ref{fig:basis_eigen}.
The matrix $C_{ij}$ in this case is not skew-symmetric, and it has 
eigenvalues with positive real part (see Figure \ref{fig:basis_eigen}).
This implies that $\bm Z(t)=\exp[-t\bm C]$ is a contraction map 
that takes any function $\theta(x)\in D(F)$ and continuously deforms 
it to $\theta(x)=0$ (see Figure \ref{fig:function_df}).
\begin{figure}[t]
\centerline{(a)\hspace{6cm} (b)}
\centerline{
\includegraphics[height=5.5cm]{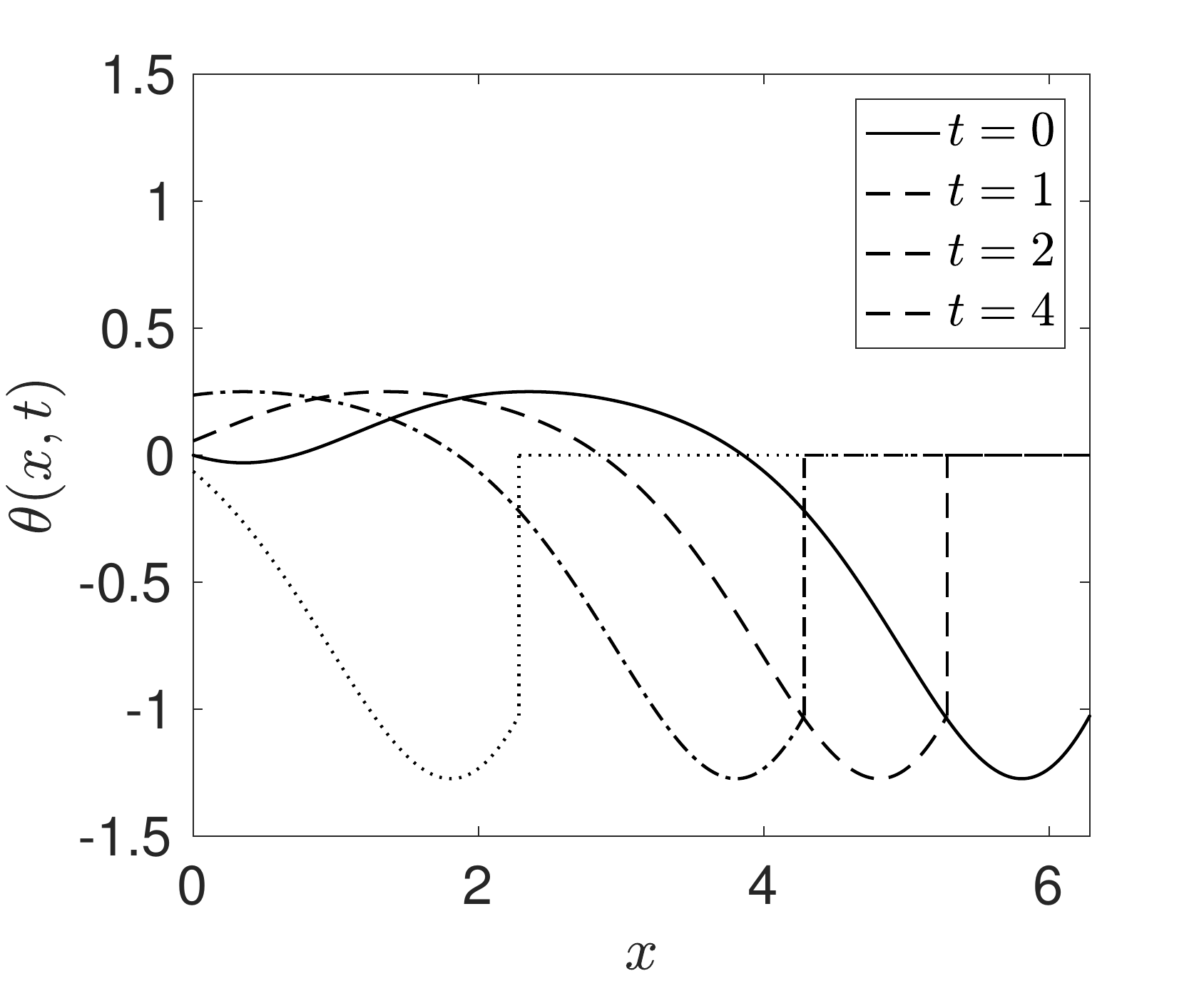}
\includegraphics[height=5.5cm]{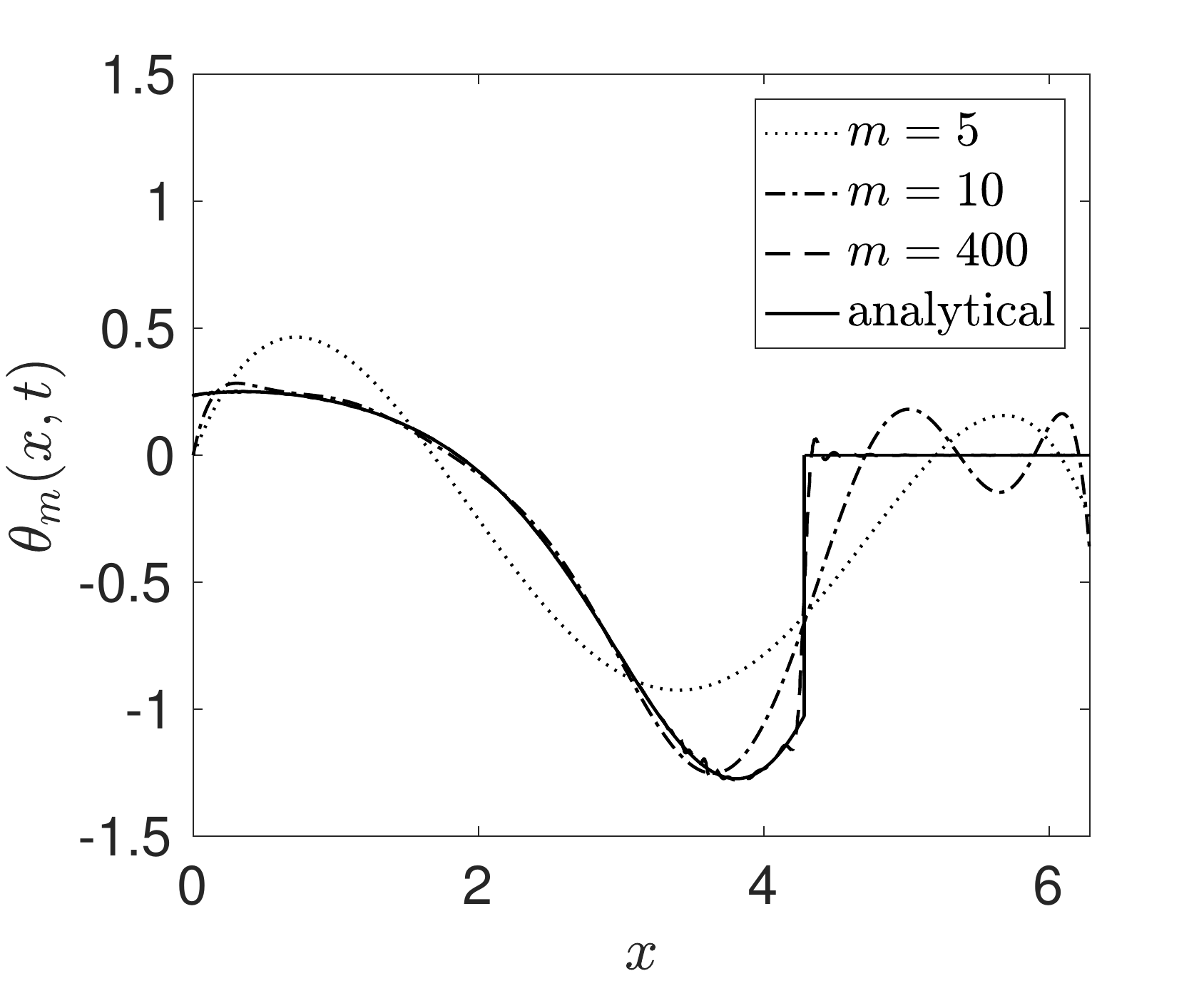}
}
\caption{\color{r}
(a) Evolution of  $\theta(x)= x\exp[-\sin(x/2)^2]\cos(x
+4)/4$ under the (infinite-dimensional) 
semigroup generated by the FDE \eqref{FDE3}, 
i.e, $U(t)\theta(x)$  (see Eq. \eqref{semiG});   
(b) Convergence of $\theta_m(x,t)$ to  $\theta(x,t)$ 
as we increase the number of dimensions $m$ at $t=2$.
}

\label{fig:function_df}
\end{figure}
From a dynamical system viewpoint, 
$\exp[-t\bm C]\bm a$ is in fact stable spiral.  
This means that if we set the initial condition as in 
\eqref{abs1}, then we obtain $f(a_1,...,a_m,t)\rightarrow 1$ 
everywhere as $t\rightarrow \infty$. The speed at which 
$f$ goes to one can be bounded by the spectral radius of $\bm C$.
The solution to the functional equation \eqref{FDE3}, with initial 
condition \eqref{expFun} is 
\begin{equation}
F([\theta],t)= \exp\left[-\int_{0}^{2\pi} \theta(x,t)^2dx\right],
\label{eq:SF}
\end{equation}
where $\theta(x,t)$ is defined in \eqref{theta_xt}.
More explicitly, 
\begin{equation}
F([\theta],t)=
\begin{cases}
\displaystyle \exp\left[-\int_t^{2\pi}\theta(x)^2dx\right] & t\in[0,2\pi]\\
\displaystyle 1 & t> 2\pi
\end{cases}.
\label{THESOLUTION}
\end{equation}

\paragraph{Regularity of the Solution Functional}
The characteristic system associated with the 
first-order PDE \eqref{PDE-advR} is
\begin{equation}
\left\{
\begin{array}{l}
\displaystyle \frac{da_k}{dt}=\sum_{j=1}^m a_j 
C_{kj}\qquad 
C_{kj}= \int_{0}^{2\pi}\varphi_j(x)\frac{d \varphi_k(x)}{dx}dx\\
\displaystyle a_k(0)=\int_{0}^{2\pi}\theta(x)\varphi_k(x)dx
\end{array}\right.
\end{equation}
Given any test function $\theta(x)$ in the function 
space \eqref{Fspace9}, e.g., the function \eqref{thef}, 
the semigroup $\exp(-t\bm C)$ pushes forward in time its 
Fourier coefficients, yielding the function 
$\theta(x,t)$ defined in equation \eqref{theta_xt}. In particular, 
if we consider the initial condition \eqref{thef} then $\theta(x,t)$ 
is shown in Figure \ref{fig:function_df}(a). As easily seen, such 
function has a shock discontinuity moving leftwards with 
velocity equal to one towards the origin as time increases. 
Specifically we have $\theta(x,t)=\theta(x-t)$. 
Remarkably, $\theta(x,t)$ is not in $D_m$ if $t>0$. In fact, such 
function does not satisfy the boundary condition $\theta(0,t)=0$ 
($t>0$) and it has a shock discontinuity. 
In other words, the semigroup generated by the FDE \eqref{FDE3}
immediately pushes $\theta(x)\in D_m$ 
out of $D_m$. This has important consequences 
when we aim at approximating $\theta(x,t)$ with elements 
of $D_m$. In particular, we need to use a high resolution 
to resolve the jump at $x=0$ and the shock in $[0,2\pi]$ 
(see Figure \ref{fig:function_df}(b)).
However, we emphasize that the singularities we 
just mentioned do not have any serious effect on 
the regularity of the solution functional \eqref{eq:SF}. 
In fact, such functional involves integration in 
$x$, which is very-well defined for bounded functions 
with a finite number of discontinuities 
(see Figure \ref{fig:analytical_sol2}). Also, $F([\theta],t)$ is 
continuous in $\theta$ and smooth in time, 
thanks to the properties of the exponential semigroup. 
In particular, from equation \eqref{THESOLUTION} we 
see that 
\begin{equation}
\theta_m\rightarrow \theta^* \quad 
\Rightarrow 
\quad 
\left|F([\theta_m],t)- F([\theta^*],t)\right|\rightarrow 0
\end{equation}
i.e., the the solution functional is continuous. Moreover,
if we restrict the set of admissible test functions to a 
function space $D(F)$ that is ``close enoug h'' to $D_m$, e.g. 
in the sense of \eqref{functionaldeviation} or \eqref{knw}, then 
we are allowed to say that the functional is approximable $D_m$.
In Figure \ref{fig:analytical_sol2} we show convergence of 
$F([\theta_m],t)$ to $F([\theta^*],t)$ as we increase $m$, 
where 
\begin{equation}
\theta^*(x)= \frac{x}{4}\exp\left[-\sin\left(\frac{x}{2}\right)^2\right]
\cos(x+4),
\label{thef}
\end{equation}
and $\theta_m(x)$ is the projection of $\theta^*(x)$ in the finite 
dimensional function space $D_m$ spanned by the orthonormal 
polynomial basis shown in Figure \ref{fig:basis_eigen}. The function 
\eqref{thef} is shown Figure \ref{fig:function_df} (case $t=0$). 
\begin{figure}[t]
\centerline{(a)\hspace{6.5cm} (b)}
\centerline{
\includegraphics[height=5.5cm]{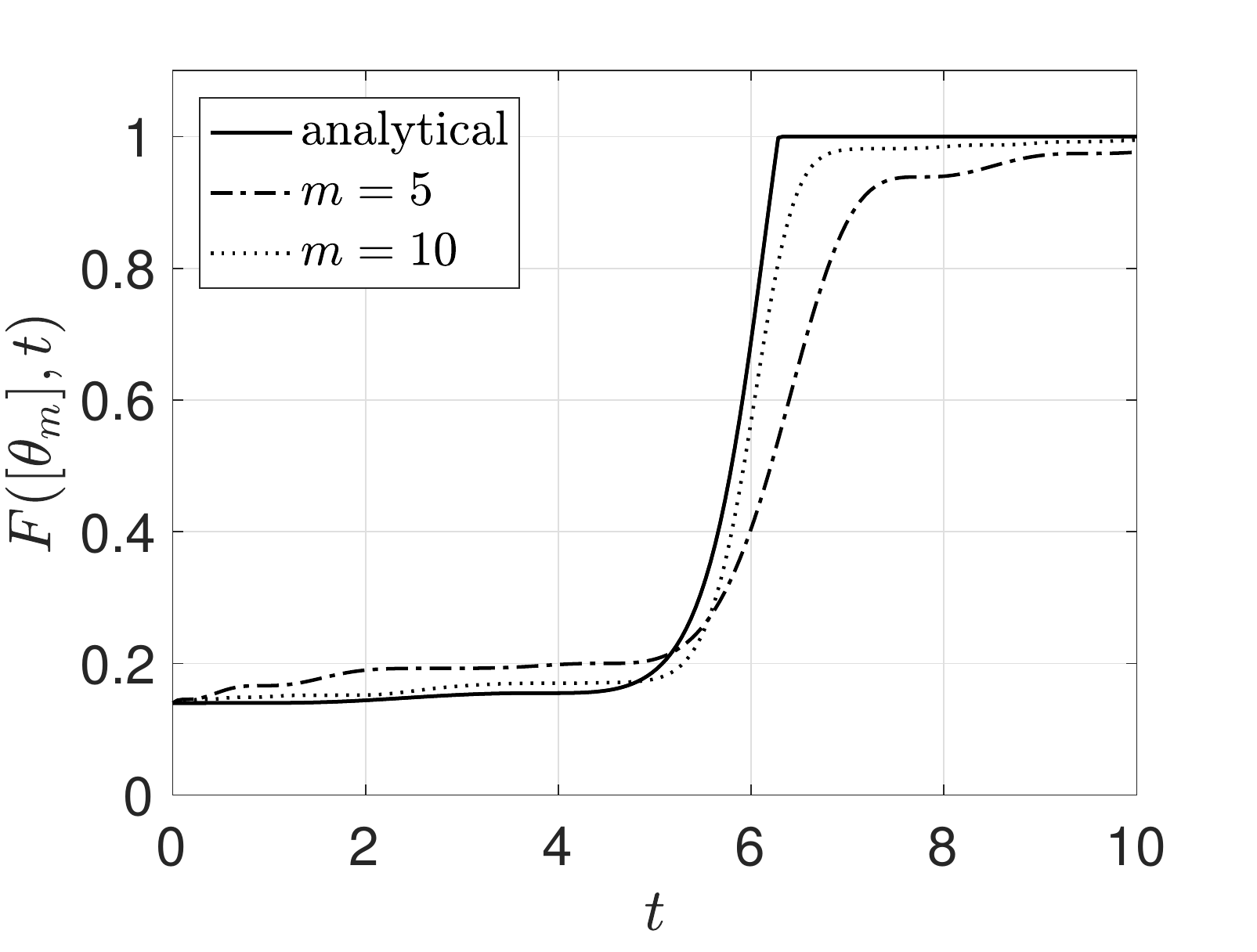}
\includegraphics[height=5.5cm]{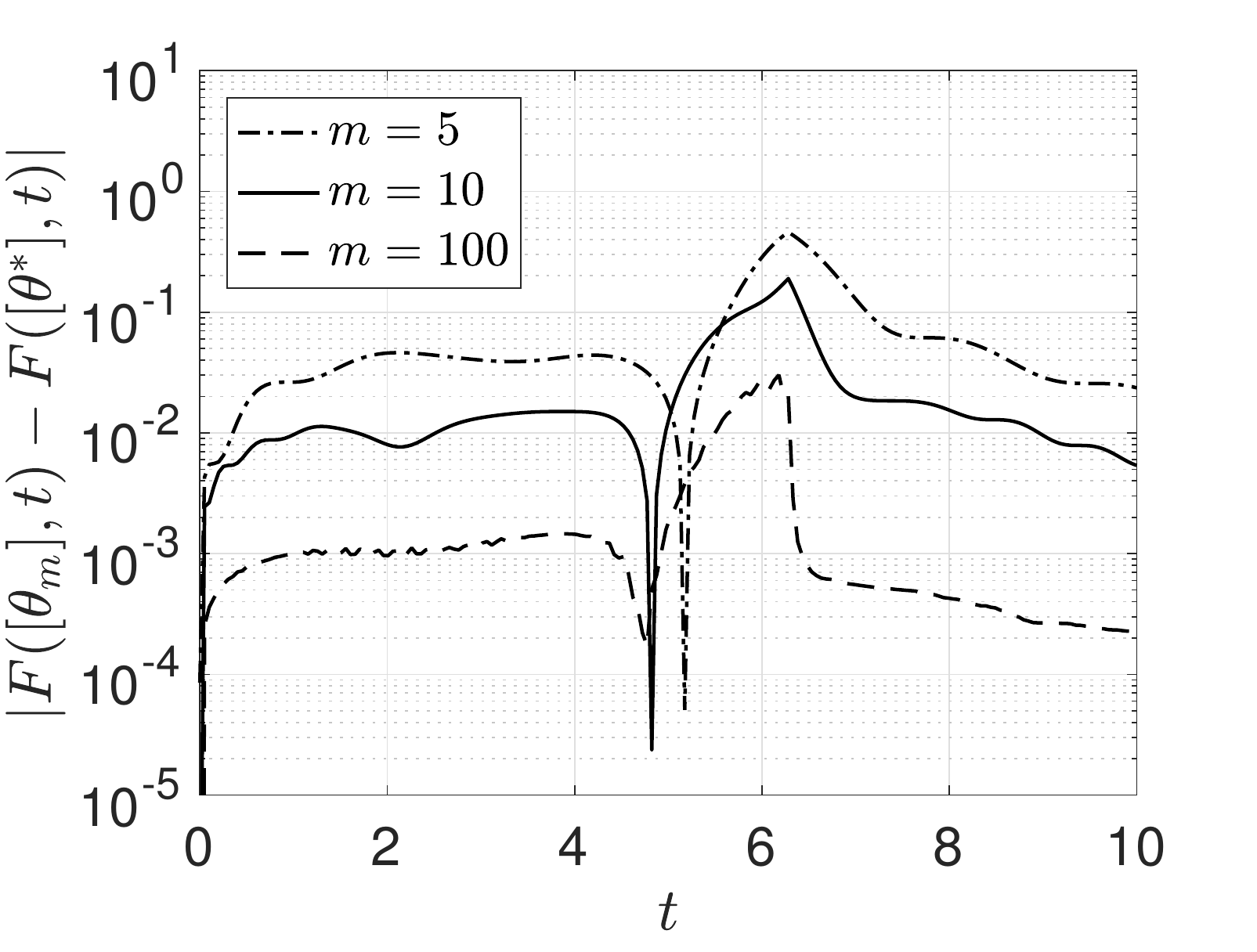}
}
\caption{\color{r}
Functional convergence with the number of 
dimensions $m$. We set $\theta^*$ as in \eqref{thef} and 
evaluate $F([\theta_m],t)$ and the $L_{\infty}$ error 
$| F([\theta^*],t)-F([\theta_m],t)|$ versus time for different $m$, 
where $\theta_m(x)$ is the projection of $\theta^*(x)$ in the finite 
dimensional function space $D_m$ spanned by the orthonormal 
polynomial basis shown in Figure \ref{fig:basis_eigen}. It is seen that $m=10$ yields reasonably 
accurate results.}
\label{fig:analytical_sol2}
\end{figure}

\paragraph{Functional Derivatives}
The first-order functional derivative of the solution 
functional \eqref{THESOLUTION} is 
\begin{equation}
\frac{\delta F([\theta],t) }{\delta \theta(x)}= 
\begin{cases}
\displaystyle -2 \theta(x) \exp\left[-\int_t^{2\pi}\theta(x)^2 dx\right] & t\in[0,2\pi],\\
\displaystyle -2 \theta(x) & t>2\pi.
\end{cases}
\label{FunctionalDerivative}
\end{equation}
Note that the functional derivative at time $t$ evaluated at 
$\theta(x)$ is simply a rescaled version of $\theta(x)$, 
where the scaling factor grows from 
$2\exp\left[\int_0^{2\pi}\theta(x)dx\right]$ (at $t=0$) 
to $2$ (at $t=2\pi$).
Such derivative can be expressed in $D_m$ as 
\begin{equation}
\left.\frac{\delta F([\theta],t)}{\delta \theta(x)}\right|_{D_m}=
\sum_{j=1}^m\frac{\partial f}{\partial a_j}\varphi_j(x),
\label{ffd3}
\end{equation}
where 
\begin{equation}
f(a_1,...,a_m,t)=\prod_{j=1}^m \exp\left[-\left(\sum_{k=1}^m 
Z_{jk}(t)a_k\right)^2\right].
\end{equation}
More explicitly,
\begin{equation}
\left.\frac{\delta F([\theta],t)}{\delta \theta(x)}\right|_{D_m}=
-2f(a_1,...,a_m,t) \sum_{p=1}^m \varphi_p(x) \sum_{k,j=1}^m 
Z_{jp}(t) Z_{jk}(t)a_k.
\label{FDeRT}
\end{equation}
In Figure \ref{fig:sol2D} we plot the analytical solution \eqref{eq:SF} in 
two dimensions. 
\begin{figure}
\centerline{\footnotesize$t=0$\hspace{3.2cm}$t=1$\hspace{3.2cm}$t=3$\hspace{3.2cm}$t=6$}
\centerline{
\includegraphics[height=3.5cm]{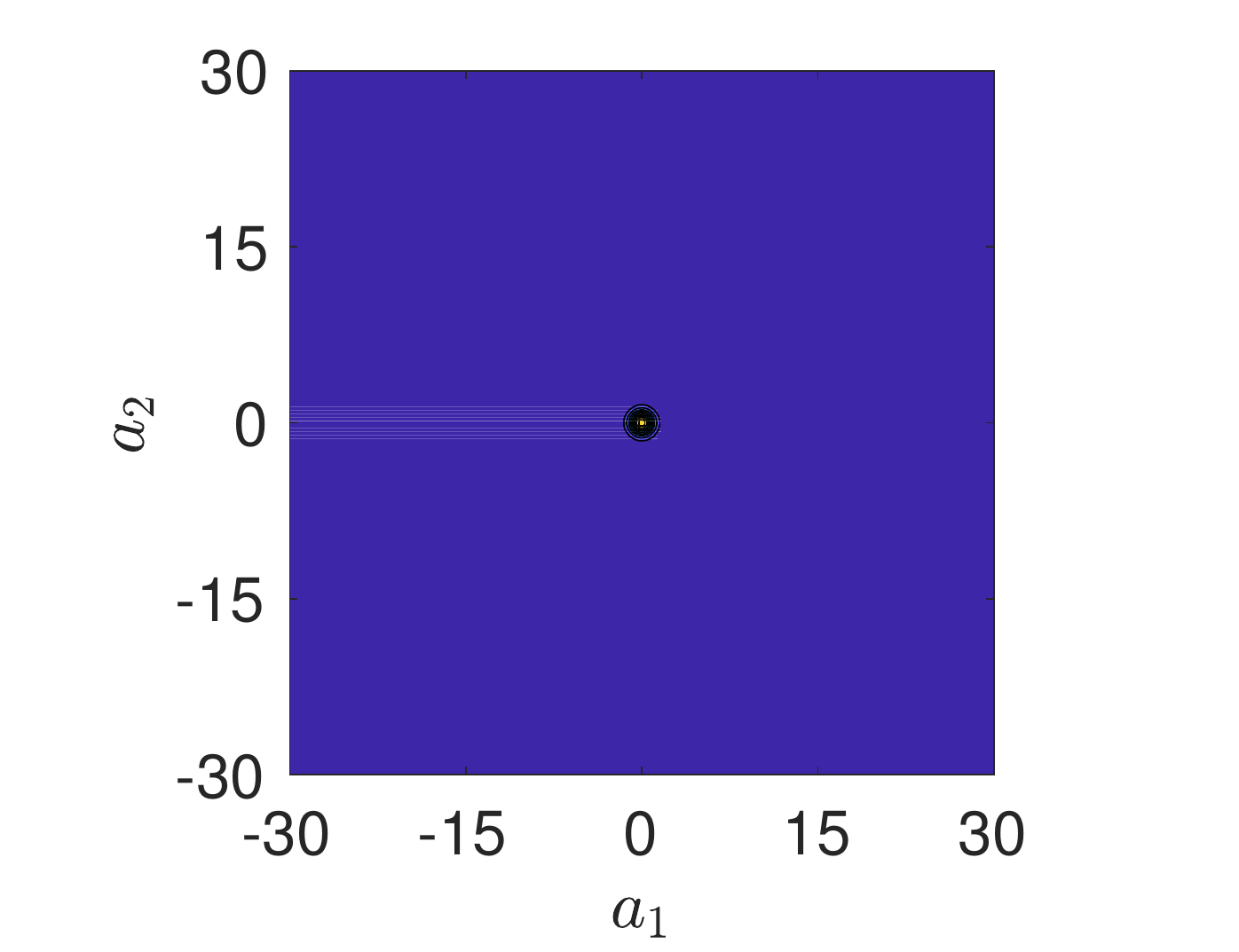}\hspace{-0.7cm}
\includegraphics[height=3.5cm]{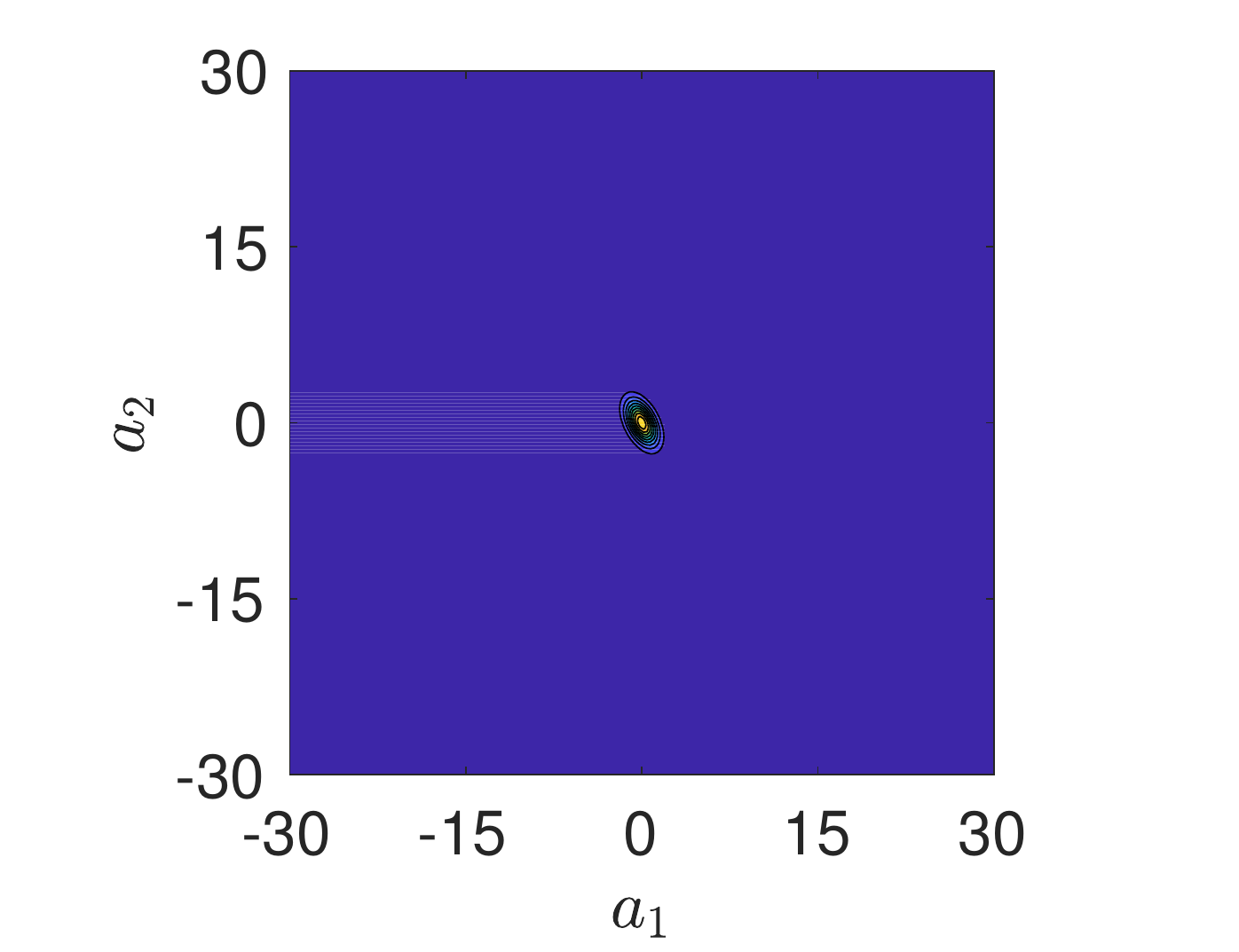}\hspace{-0.7cm}
\includegraphics[height=3.5cm]{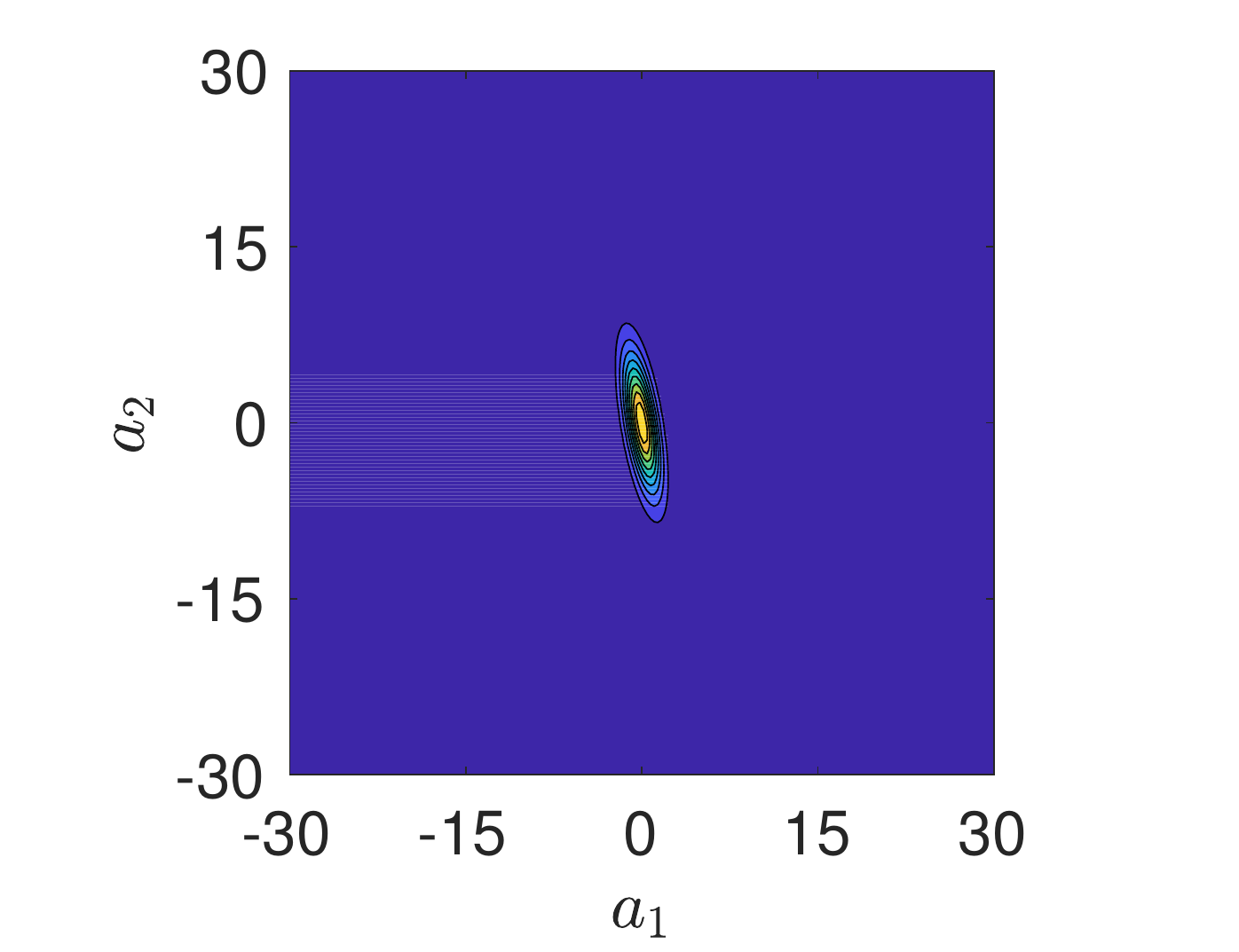}\hspace{-0.7cm}
\includegraphics[height=3.5cm]{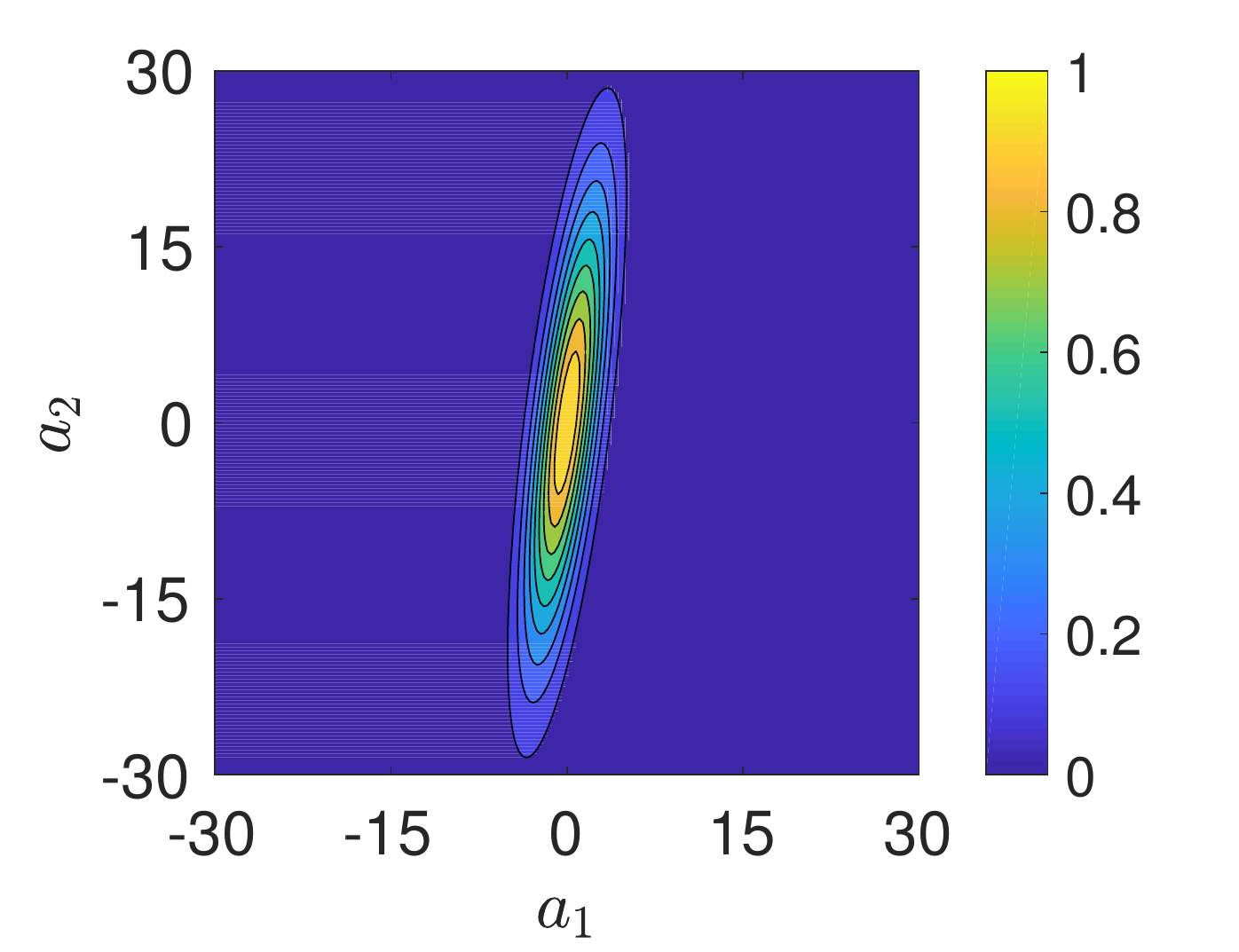}
}

\caption{\color{r}
Solution to the multivariate PDE \eqref{PDE-advR} 
(with $h=0$) in $m=2$ dimensions. The stable spiral at the origin of the characteristic system attracts every curve in the space of modes and ultimately yields $f=1$ everywhere after a transient (see also Figure \ref{fig:analytical_sol2}).}
\label{fig:sol2D}
\end{figure}
The side length of the hypercube cube that encloses 
any level set of the solution at time $t$ depends on 
the number of dimensions $m$. In particular, if we are 
interested in estimating the size of the hypercube that encloses 
the set $\{\bm a\in \mathbb{R}^m \,|\, f(\bm a,t)\geq \epsilon\}$, then 
we can use the formula \footnote{\color{r} 
The inequality \eqref{est:b} can 
be obtained by noticing that the solution $f(\bm a,t)$ is in the form 
\begin{equation}
f(\bm a,t) = e^{-\bm a^T \bm S(t) \bm a}\qquad \bm S(t)= \bm Z(t)^T \bm Z(t), 
\end{equation}
i.e., 
\begin{equation}
-\log(f(\bm a,t)) = \bm a^T \bm S(t)\bm a.
\end{equation}
By diagonalizing $\bm S(t)$ we obtain 
\begin{equation}
(\bm V(t) \bm a)^T \bm \Lambda(t) (\bm V(t) \bm a) = - \log(f(\bm a,t)).
\end{equation}
Upon definition of the rotated coordinate system 
$\bm y(t)=\bm V(t) \bm a$, we find that the largest 
semi-axis of the ellipse representing the $\epsilon$ level 
set of $f(\bm a,t)$ is 
\begin{equation}
\frac{1}{2}\sqrt{\frac{-\log(\epsilon)}{\lambda_{min}}}, 
\end{equation} 
where $\lambda_{min}$ is the smallest eigenvalue of $\bm S(t)$.
This formula assumes that there is no rotation in the Gaussian 
function $f(\bm a,t)$ during the dynamics and therefore it provides a 
conservative upper bound $b$ that coincides with the largest semi-axis 
of the ellipsoidal level set. On the other hand, there exist a 
rotation of the ellipsoid that minimizes the size of the aforementioned 
hypercube. Such rotation  aligns the largest semi-axis with the diagonal 
of the hypercube. 
We recall that the diagonal of a hypercube in dimension $m$ has length 
$\sqrt{m}\ell$, where $\ell$ is the side length of the hypercube.
Therefore, the upper and the lower bound estimates for $b$, i.e, the 
half side length of the hypercube that encloses the $-\log(\epsilon)$ 
level set of the solution are 
\begin{equation}
\frac{1}{2}\sqrt{\frac{-\log(\epsilon)}{ m \lambda_{min}}}\leq b 
\leq \frac{1}{2}\sqrt{\frac{-\log(\epsilon)}{\lambda_{min}}}.
\end{equation}
} 
\begin{equation}
\frac{1}{2}\sqrt{\frac{-\log(\epsilon)}{ m \lambda_{min}}}\leq b 
\leq \frac{1}{2}\sqrt{\frac{-\log(\epsilon)}{\lambda_{min}}}.
\label{est:b}
\end{equation}
where $\lambda_{min}$ is the smallest eigenvalue of the matrix $\bm Z(t)^T \bm Z(t)$ (see Equation \eqref{gg1i}).
In Figure \ref{fig:enclosing_cube} we plot 
the upper and the lower bound estimates of 
half of the side length of the hypercube that encloses the $10^{-10}$ 
level set of the solution at $t=1$ 
versus the number of dimensions $m$. It is seen that 
the size of the hypercube increases exponentially fast with $m$. 
This has important consequences when it comes to numerical 
simulations. In particular, if we perform simulations with 
far field boundary conditions, then the size of the computational 
domain should  be chosen large enough to accommodate the 
support of the solution throughout the simulation time interval 
of interest. 

\begin{figure}[t]
\centerline{\hspace{0.2cm}$t=1$ \hspace{7.cm} $t=3$}
\centerline{
\includegraphics[height=6cm]{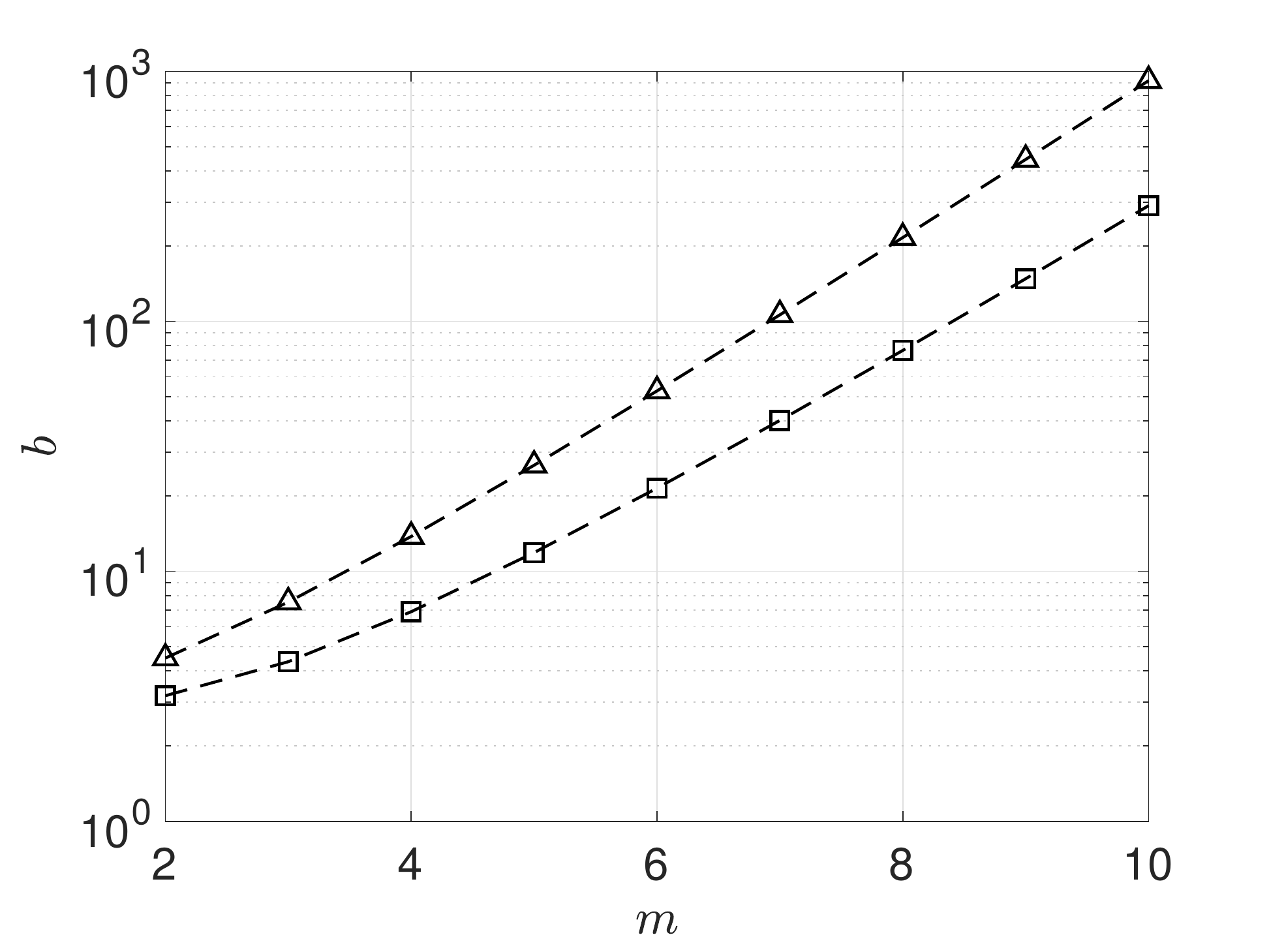}
\includegraphics[height=6cm]{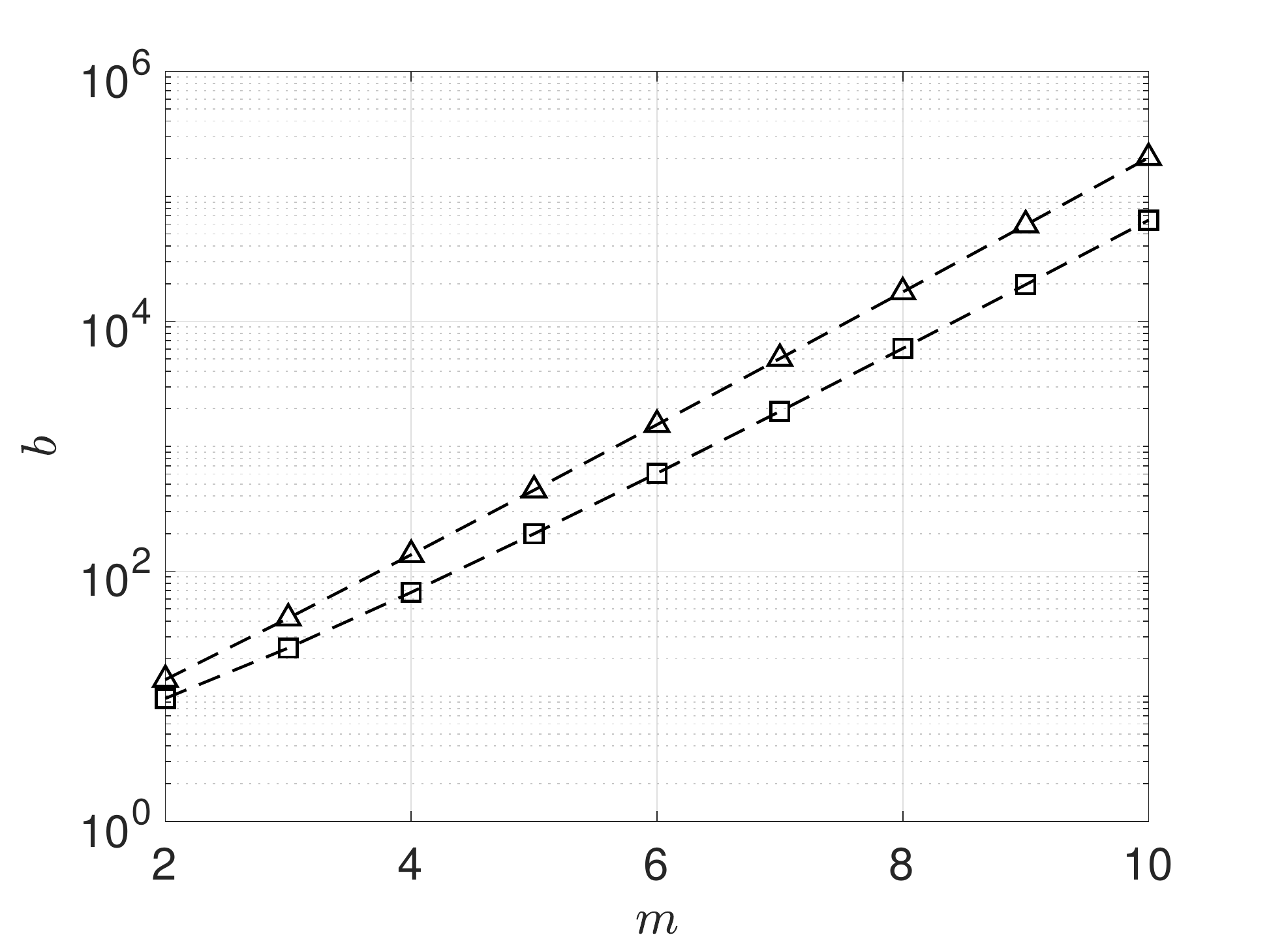}
}
\caption{\color{r}Upper and lower bound estimates \eqref{est:b} of 
the half side length of the hypercube that encloses the $10^{-10}$ 
level set of the solution at $t=1$ and $t=3$ versus the number of 
dimensions $m$. It is seen that $b$ increases 
exponentially fast with $m$. 
The computational domain should be at least as 
big as $[-b,b]^m$, where $b$ is the number given by the 
upper curve. For example, if we are aiming at resolving the 
numerical solution in the case $m=6$ within the time interval 
$[0,3]$, then we should consider the computational 
domain $[-b,b]^6$ with $b\simeq 1600$.}
\label{fig:enclosing_cube}
\end{figure}

\subsection{Numerical Discretization}
\label{sec:numericalDiscretizationADVR}
Consider the multivariate PDE \eqref{PDE-advR} with $h=0$, 
hereafter written in the operator form
\begin{equation}
\frac{\partial f}{\partial t}=L f,
\label{PDE11}
\end{equation}
where
\begin{equation}
L= - \sum_{j,k=1}^ma_kC_{jk}\frac{\partial }{\partial a_j}.
\label{LadvR1}
\end{equation}
Note that $L$ is a {\em separable} linear 
operator with separation rank $r_L=m^2$.
In fact,  $L$ can be written in the form 
\begin{equation}
L= \sum_{q=1}^{m^2} \alpha_q L^q_1(a_1)\cdots L_m^q(a_m),
\label{eq:operators}
\end{equation}
for suitable one-dimensional linear operators $L^q_j(a_j)$ defined in Table \ref{tab:ordering}.
\begin{table}\color{r}
\centering
\begin{tabular}{c|ccccc}
$q$ & $\alpha_q$ & $L_1^q$ & $L_2^q$ & $\cdots$ & $L_m^q$\\
\hline\\
$1$ &  $-C_{11}$ &  $a_1\partial_{a_1}$ & $1$ & $\cdots$& $1$\\
$2$ &$-C_{12}$ &  $      \partial_{a_1}$ & $a_2$ & $\cdots$& $1$\\
$\vdots $&  $\vdots $ & $\vdots$ & $\vdots$ &$\vdots$ &$\vdots$ \\
$m$ &$-C_{1m}$ & $\partial_{a_1}$       & $1$ & $\cdots$& $a_m$\\
$m+1$ &$-C_{21}$ & $a_1$       & $\partial_{a_2}$ & $\cdots$& $1$\\
$m+2$ &$-C_{22}$ & $1$       & $a_2\partial_{a_2}$ & $\cdots$& $1$\\
$\vdots $&  $\vdots $ & $\vdots$ & $\vdots$ &$\vdots$ &$\vdots$ \\
$2m$ &$-C_{2m}$ & $1$       & $\partial_{a_2}$ & $\cdots$& $a_m$\\
$\vdots $&  $\vdots $ & $\vdots$ & $\vdots$ &$\vdots$ &$\vdots$ \\
$m^2-m+1$ &$-C_{m1}$ & $a_1$ & $1$ & $\cdots$& $\partial_{a_m}$\\
$m^2-m+2$ &$-C_{m2}$ & $1$ & $a_2$ & $\cdots$& $\partial_{a_m}$\\
$\vdots $&  $\vdots $ & $\vdots$ & $\vdots$ &$\vdots$ &$\vdots$ \\
$m^2$ &$-C_{mm}$ & $1$ & $1$ & $\cdots$& $a_m\partial_{a_m}$\\
\end{tabular}
\caption{\color{r}Ordering of the linear operators defined 
in  equation \eqref{eq:operators}.}
\label{tab:ordering}
\end{table}
The solution to the multivariate PDE \eqref{PDE11} can be 
represented by using any tensor series expansion. In particular, 
hereafter we consider the hierarchical Tucker format (Section \ref{sec:CP}) and the canonical polyadic tensor 
decomposition (Section \ref{sec:CTD-ALS}).

\subsubsection{HT and CP Algorithms  with  Explicit Time Stepping}
\label{sec:HTD}
Let us discretize the PDE \eqref{PDE11} in time by using any explicit 
time stepping scheme, for example the second-order Adams-Bashforth scheme (see Section \ref{sec:ADI_SSE})
\begin{equation}
f_{n+2}=f_{n+1}+\frac{\Delta t}{2} L\left(3f_{n+1}-f_n\right).
\label{AB2_}
\end{equation}
In this setting, we see that the only operations needed to 
compute $f_{n+2}$ with tensor methods are: i) addition, 
ii) application of a (separable) linear operator, and iii) rank 
reduction\footnote{On the other hand, if we employ 
implicit time-discretization schemes such as \eqref{CNFDE} we 
end up solving linear systems with tensor methods. This was 
addressed, e.g., in \cite{Ortega,Etter,Reynolds}.}, the last 
operation being the most important among all three. 
From a computational viewpoint, it would be also very useful if we  
could split the tensor operations yielding $f_{n+2}$ into sequences 
of tensor operations followed by rank reduction. In this way 
we could minimize the storage requirements and the overall 
computational cost. For example, we could split \eqref{AB2_} 
as 
\begin{enumerate}
\item Compute a low rank representation of $w_{n+1}=(3f_{n+1} - f_{n})$,
\item Compute a low rank representation of 
$q_{n+1}= L w_{n+1}$,
\item Compute a low rank representation of 
$f_{n+2}=f_{n+1} +\Delta t q_{n+1}/2$.
\end{enumerate}
Are we allowed to do so? Unfortunately no. 
Splitting tensor operations into 
sequences of tensor operations followed by 
rank reduction usually yields severe cancellation 
errors.
In some cases, this problem can be overcome. For example, 
an efficient and robust algorithm that allows us to split 
sums and rank reduction operations was recently proposed 
in \cite{Kressner2014} in the context of 
hierarchical Tucker formats. The algorithm leverages on 
the block diagonal structure that arises when adding 
hierarchical Tucker format. 
Also, the vector resulting from application of the separable 
linear operator $L$ to the hierarchical Tucker 
format $w_{n+1}$ at point ii) above 
can be computed very efficiently if we have available 
a hierarchical Tucker representation of the operator $L$. 
Such representation can be easily constructed in a Fourier 
collocation setting by vectorizing all operators in table 
\ref{tab:ordering}. Each one-dimensional 
operator is represented relative to a basis trigonometric (nodal) 
polynomials (Fourier spectral collocation method) \cite{Hesthaven} 
In practice, we simply need to convert the 
vectorized canonical polyadic series of $L$ into a hierarchical
Tucker expansion, which is a relative standard operation.

\subsubsection{CP-ALS Algorithm with Implicit Time Stepping}
\label{sec:CTD-ALS}
Let us discretize the PDE \eqref{PDE11} in time by using the 
Crank-Nicolson method. This yields 
\begin{equation}
\left[I-\frac{\Delta t}{2} L \right]f_{n+1}=
\left[I+\frac{\Delta t}{2} L \right]f_n.
\label{CNFDE}
\end{equation}
This equation can be written in a compact notation as 
\begin{equation}
Af_{n+1}=Bf_n,
\end{equation}
where  
\begin{equation}
A = I-\frac{\Delta t}{2} L, \quad \textrm{and}\quad B = I+\frac{\Delta t}{2} L.
\end{equation}
By using the definition of $L$ given in \eqref{eq:operators} 
(see also Table \ref{tab:ordering}), it is clear that both $A$ and $B$ are 
separable operators in the form 
\begin{align}
A= 1\cdots 1 +\frac{\Delta t}{2} \sum_{i,j=1}^m a_iC_{ji}
\frac{\partial }{\partial a_j},\label{AL}\\
B= 1\cdots 1 - \frac{\Delta t}{2} \sum_{i,j=1}^m  a_iC_{ji}
\frac{\partial }{\partial a_j}.\label{BL}
\end{align}
These operators have separation rank $r_A=r_B=m^2+1$ and can 
be conveniently written as
\begin{align}
A= \sum_{q=0}^{m^2} \eta_q E^q_1(a_1)\cdots E^q_m(a_m), \qquad 
B= \sum_{q=0}^{m^2}  \zeta_q E^q_1(a_1)\cdots E^q_m(a_m),\label{BLs}
\end{align}
where all quantities are defined in Table \ref{tab:AL}.
The difference with Table \ref{tab:ordering} is that 
we added one row (the zeroth one) to represent the 
identity operator $1\cdots 1$, and we rescaled 
all coefficients $C_{ij}$.
\begin{table}\color{r}
\centering
\begin{tabular}{c|cccccc}
$q$ &  $\eta_q$ & $\zeta_q$ & $E_{1}^q$ & $E_{2}^q$ & $\cdots$ & $E_{m}^q$\\
\hline\\
$0$ &  $1$ &$1$ &  $1$ & $1$ & $\cdots$& $1$\\
$1$ &  $\Delta t C_{11}/2$ &  $-\Delta t C_{11}/2$ &  $a_1\partial_{a_1}$ & $1$ & $\cdots$& $1$\\
$2$ &$\Delta t C_{12}/2$ &$-\Delta t C_{12}/2$&  $\partial_{a_1}$ & $a_2$ & $\cdots$& $1$\\
$\vdots $&  $\vdots $ & $\vdots$ & $\vdots$ &$\vdots$ &$\vdots$ \\
$m$ &$\Delta tC_{1m}/2$ & $-\Delta tC_{1m}/2$ & $\partial_{a_1}$       & $1$ & $\cdots$& $a_m$\\
$m+1$ &$\Delta tC_{21}/2$ & $-\Delta tC_{21}/2$& $a_1$       & $\partial_{a_2}$ & $\cdots$& $1$\\
$m+2$ &$\Delta tC_{22}/2$ & $-\Delta tC_{22}/2$ & $1$       & $a_2\partial_{a_2}$ & $\cdots$& $1$\\
$\vdots $&  $\vdots $ & $\vdots$ & $\vdots$ &$\vdots$ &$\vdots$ \\
$2m$ &$\Delta tC_{2m}/2$ &$-\Delta tC_{2m}/2$& $1$       & $\partial_{a_2}$ & $\cdots$& $a_m$\\
$\vdots $&  $\vdots $ & $\vdots$ & $\vdots$ &$\vdots$ &$\vdots$ \\
$m^2-m+1$ &$\Delta tC_{m1}/2$&$-\Delta tC_{m1}/2$ & $a_1$ & $1$ & $\cdots$& $\partial_{a_m}$\\
$m^2-m+2$ &$\Delta tC_{m2}/2$ &$-\Delta tC_{m2}/2$& $1$ & $a_2$ & $\cdots$& $\partial_{a_m}$\\
$\vdots $&  $\vdots $ & $\vdots$ & $\vdots$ &$\vdots$ &$\vdots$ \\
$m^2$ &$\Delta tC_{mm}/2$&$-\Delta tC_{mm}/2$ & $1$ & $1$ & $\cdots$& $a_m\partial_{a_m}$\\
\end{tabular}
\caption{\color{r}Ordering of the linear operators $A$ and $B$ defined 
in  \eqref{BLs}.}
\label{tab:AL}
\end{table}
A substitution of the CP decomposition 
\begin{equation}
\hat{f}_{n+1} = \sum_{l=1}^r \prod_{k=1}^m{G^l_k(a_k,t_{n+1})}
\label{CP-numerical}
\end{equation}
into equation \eqref{CNFDE} yields the residual 
\begin{equation}
R = A\hat{f}_{n+1}-B\hat{f}_n.
\label{rRr}
\end{equation}
Minimization of the $L_2$ norm of \eqref{rRr} with respect to 
${\bm\beta}_q(t_{n+1})$ yields the linear systems 
of equations
\begin{equation}
 \bm M^L_q \bm \beta_q(t_{n+1}) = \bm M^R_q \bm \beta_q(t_{n}),\qquad q=1,...,m
\label{ALS2}
\end{equation}
where 
\begin{equation}
\bm M^L_q = \sum_{e,z=0}^{m^2} K^L_{ez}
\left[
{\displaystyle \substack{\vspace{0.2cm}\\\vspace{0.04cm} m\\\bigcirc\\k=1\\k\neq q}}
%\prod_{\substack{k=1\\k\neq q}}^m 
\hat{\bm \beta}_{k}(t_{n+1})^T \bm 
E^{ez}_k \hat {\bm \beta}_{k}(t_{n+1})
\right]^T\otimes \left[{\bm E}^{ez}_q\right]^T, 
\label{ML}
\end{equation}
\begin{equation}
\bm M^R_q = \sum_{e,z=0}^{m^2} {K}^R_{ez}
\left[ 
{\substack{\vspace{0.2cm}\\\vspace{0.04cm} m\\\bigcirc\\k=1\\k\neq q}}
%\prod_{\substack{k=1\\k\neq q}}^m 
\hat{\bm \beta}_{k}(t_{n})^T 
{\bm E}^{ez}_k \hat{\bm \beta}_{k}(t_{n+1})
\right]^T\otimes \left[{\bm E}^{ez}_q\right]^T, 
\label{MR}
\end{equation}
and
\begin{equation}
\left[{\bm E}^{ez}_q\right]_{sh}=\int_{-b}^b E_{q}^e(a)\phi_s(a) E_{q}^z(a)\phi_h(a)da.
\label{eq:77}
\end{equation}
In equations \eqref{ML}-\eqref{MR}, 
$\bigcirc$ denotes the Hadamard matrix 
product, $\otimes$ is the Kroneker matrix product, 
$\hat {\bm \beta}_{k}$ is the matrix version of $\bm \beta_k$, i.e., 
\begin{equation}
\hat {\bm \beta}_{k}(t_n)=
\left[\begin{array}{ccc}
\beta^1_{k1}(t_n)&\cdots& \beta^r_{k1}(t_n)\\
\vdots  &\ddots & \vdots \\
\beta^1_{kQ}(t_n) &\cdots & \beta^r_{kQ} (t_n)
\end{array}
\right],
\end{equation}
$\bm E^{ez}_q$ is the  $Q\times Q$ matrix \eqref{eq:77}, 
and  $K^L_{ez}$ and ${K}^R_{ez}$ are entries of the 
matrices
\begin{equation}
\bm K^L = \bm \eta \bm \eta^T,\quad \textrm{and}\quad 
\bm K^R = \bm \eta \bm \zeta^T,
\end{equation}
where $\bm \eta$ and $\bm \zeta$ are column vectors with entries $\eta_q$ and $\zeta_q$ defined in Table \ref{tab:AL}. 

\paragraph{Computing the Matrix System}
There are many of symmetries we can exploit when 
constructing the separated series expansion of the 
operators $A$ and $B$ in \eqref{BLs}. 
Indeed, a closer look at Table \ref{tab:AL} suggests that 
if we employ the same series expansion in each 
variable $a_j$ (e.g.,  a trigonometric series)   
then the number of operators in \eqref{BLs} that 
we effectively need to compute reduces to the 
following four 
\begin{equation}
1, \qquad a_j, \qquad \frac{\partial}{\partial a_j}, \qquad a_j\frac{\partial}{\partial a_j}.
\end{equation}
This means that the number of terms that are effectively 
different in the fundamental matrix \eqref{eq:77} are {\em only 12} 
(9 if we are willing to employ matrix transposes). 
Specifically,   
\begin{equation}
\begin{array}{l l l l}
\displaystyle 1.\, \int_{-b}^b  \phi_s \phi_h da \quad & 
\displaystyle 2.\,\int_{-b}^b a \phi_s \phi_h da, \quad & 
\displaystyle 3.\,\int_{-b}^b    \phi_s\frac{d\phi_h}{da} da,\quad &
\displaystyle 4.\,\int_{-b}^b a \phi_s \frac{d\phi_h}{da}da,\\
\displaystyle 5.\, \int_{-b}^b a^2 \phi_s \phi_h da, \quad  &
\displaystyle 6\,\int_{-b}^b a^2  \phi_s\frac{d\phi_h}{da} da, \quad &
\displaystyle 7.\,\int_{-b}^b \frac{d\phi_s}{da} \phi_h da, \quad & 
\displaystyle 8.\, \int_{-b}^b  a \frac{d\phi_s}{da} \phi_h da, \\  
\displaystyle 9.\,\int_{-b}^b\frac{d\phi_s}{da}\frac{d\phi_h}{da} da, \quad &
\displaystyle 10.\,\int_{-b}^b a \frac{d\phi_s}{da} \frac{d\phi_h}{da} da, \quad & 
\displaystyle 11.\, \int_{-b}^b  a^2 \frac{d\phi_s}{da} \phi_h da, \quad & 
\displaystyle 12.\,\int_{-b}^b a^2 \frac{d\phi_s}{da}\frac{d\phi_h}{da} da.
\end{array}
\label{eq:all integrals}
\end{equation}
All these integrals can be pre-computed and stored as $Q\times Q$ 
matrices (see Eq. \eqref{eq:50}). For each $e$, $z$ 
and $q$, the tensor \eqref{eq:77} corresponds to one of the $12$ 
integrals above. Such map, denoted as $g^{ez}_q$, 
takes in the triple $(e,z,q)$, 
where $q\in \{1,...,m\}$ and $e,z\in\{0,...,m^2\}$, 
and it returns a number between $1$ and $12$ identifying 
which integral in the set \eqref{eq:operators} corresponds 
to the tensor entry in \eqref{eq:77}. For example, if we sort 
the operators  as in Table \ref{tab:AL} then in $m=2$ 
dimensions we have
\begin{equation}
g^{ez}_{1}=\left[
\begin{array}{ccccc}
     1    & 8 &    7&     2 &    1\\   
     4   & 12 &   10 &    6 &    4\\
     3   & 10   &  9   &  4    & 3\\
     2   & 11   &  8   &  5   &  2\\
     1   &  8    & 7    & 2    & 1
\end{array}
\right]\qquad \qquad 
g^{ez}_{2}=\left[
\begin{array}{c c c c c}
    1    & 1    & 2     &7     &8\\
     1     &1   &  2    & 7    & 8\\
     2     &2   &  5    & 8    &11\\
     3     &3   &  4    & 9    &10\\
     4     &4   &  6    &10   & 12
\end{array}
\right],
\label{dim2}
\end{equation}
while in dimension $m=3$ we have, e.g., 
\begin{equation}
g^{ez}_{1}=\left[
\begin{array}{cccccccccc}
     1 &    8  &   7  &   7    &  2  &   1  &   1  &   2 &    1&     1\\
     4 &   12  &  10&    10&   6  &   4  &   4  &   6 &    4&     4\\
     3 &   10  &   9 &    9  &   4  &   3  &   3  &   4 &    3&     3\\
     3 &   10  &   9 &    9  &   4  &   3  &   3  &   4 &    3&     3\\
     2 &   11  &   8 &    8  &   5  &   2  &   2  &   5 &    2&     2\\
     1 &    8   &  7  &   7   &   2  &   1  &   1 &    2 &    1&     1\\
     1 &    8   &  7  &   7   &   2  &   1  &   1 &    2 &    1&     1\\
     2 &   11  &   8  &   8  &   5  &   2  &   2 &    5 &    2&     2\\
     1 &    8   &  7  &   7   &   2  &   1  &   1 &    2 &    1&     1\\
     1 &    8   &  7  &   7   &   2  &   1  &   1 &    2 &    1&     1
     \end{array}
     \right].
     \label{dim3}
\end{equation}
A combinatorial argument shows that the number of 
entries equal to $1$, $2$, $3$, etc.,  in each matrix 
$g_{1}^{ez}$, $g_{2}^{ez}$, ..., $g_{m}^{ez}$ is the same 
(for fixed $m$). For instance, in \eqref{dim2} we 
have $4$ ones, $2$ threes, $3$ fours, $1$ 
five, etc. This is very useful when we break the sum 
in $e$ and $z$ in \eqref{ML} and \eqref{MR} into 
multiple sums, and use the associative property of the 
tensor product to reduce the number of operations. 
By using the map $g_q^{ez}$ we can immediately 
identify each matrix $\bm E_q^{ez}$. For example, 
in the case $m=2$ we have (see Eq. \eqref{dim2})
\begin{equation}
\left[{\bm E}^{11}_1\right]_{sh}=\int_{-b}^b \phi_s(a)\phi_h(a) da, \qquad 
\left[{\bm E}^{12}_1\right]_{sh}=\int_{-b}^b a\frac{d\phi_s(a)}{da}\phi_h(a) da,\quad ...\quad .
\end{equation}

\paragraph{Summary of the Algorithm}
We first compute all integrals in \eqref{eq:all integrals} 
and store them in $12$ matrices $Q\times Q$,
$Q$ being the number of degrees of freedom in each 
variable (e.g., collocation points of Fourier modes). 
We also set up the map between such set 
of matrices and any element of the tensor \eqref{eq:77}. 
Such map basically takes in the triple $(q,e,z)$, 
where $q\in \{1,...,m\}$ and $e,z\in\{0,...,m^2\}$, 
and it returns a number between $1$ and $12$ identifying 
which integral in the set \eqref{eq:operators} 
corresponds to the tensor entry in \eqref{eq:77}. In 
a matrix setting, this basically allows us to efficiently 
compute each matrix $\bm E_q^{ez}$ appearing 
in \eqref{ML} and \eqref{MR}. 
Next, we compute the compute the canonical tensor 
decomposition of the initial condition $f_0(a_1,...,a_m)$,
by applying the methods we described in Section \ref{sec:CP}. 
This gives us the set of vectors 
$\{\bm \beta_1(t_0),...,\bm \beta_m(t_0)\}$. 
With such vectors available, we can build the matrices 
$\bm M^L_1$ and $\bm M^R_1$ defined in \eqref{ML} and 
\eqref{MR}. To this end, we need an initial guess 
for $\{\bm \beta_1(t_1),...,\bm \beta_m(t_1)\}$ which we can 
take to be equal to $\{\bm \beta_1(t_0),...,\bm \beta_m(t_0)\}$, 
or a small random perturbation of it. 
With $\bm M^L_1$ and $\bm M^R_1$ in place and set, we 
can solve the linear system \eqref{ALS2} and 
update $\bm \beta_1(t_1)$. At this point 
we recompute  $\bm M^L_2$ and $\bm M^R_2$ (with the updated 
$\bm \beta_1(t_1)$) and solve for $\bm \beta_2(t_1)$. 
We repeat this process for $q=3,...,m$ and iterate over 
and over among all variables until convergence. 
Parallel versions of the ALS algorithm 
were recently proposed by Karlsson {\em et al.} 
in \cite{Karlsson}.

\subsubsection{Long-Term Integration} 
\label{sec:long-term-integration} 
A rigorous error analysis of the HT and CP-ALS 
algorithms to solve the multivariate 
PDE \eqref{PDE11} goes beyond the scope of this report
(see \cite{Bachmayr} for a recent account). 
It it useful, however, to point out a few things on 
the nature of the discretization error, in particular 
on how the temporal local truncation error depends on 
the dimension $m$. To this end, let us first recall that 
the local truncation error at time $t_{n+1}$ of 
the second-order Adams Bashforth 
(AB2) scheme applied to the linear PDE \eqref{PDE11} is 
$5\Delta t^3 L^3 f(\eta,\bm a)/12$, where $\eta$ is 
some time instant between $t_n$ and $t_{n+1}$, $f(\eta,\bm a)$ 
is the exact solution \eqref{solutionFF}, and $L$ is defined in \eqref{LadvR1}. The operator $L^2$ can be explicitly written as  
\begin{align}
L^2 =& \sum_{i,j,l,p=1}^m x_lC_{pl}\frac{\partial }{\partial a_p}\left(
x_jC_{ij}\frac{\partial }{\partial a_j}\right)\nonumber\\
=& \sum_{i,j,l,p=1}^m x_lC_{pl}C_{ij}\left(\delta_{jp}\frac{\partial }{\partial a_j}+x_j\frac{\partial^2}{\partial a_j\partial a_p}\right),
\label{L2}
\end{align}
while $L^3$ has the form 
\begin{align}
L^3 =& \sum_{i,j,l,p,q,z=1}^m x_q C_{zq}C_{pl}C_{ij}\delta_{lq}\left(\delta_{jp}\frac{\partial }{\partial a_j}+x_j\frac{\partial^2}{\partial a_j\partial a_p}\right)+\nonumber\\
& \sum_{i,j,l,p,q,z=1}^m x_q x_l C_{zq}C_{pl}C_{ij}\left(\delta_{jp}\frac{\partial^2 }{\partial a_j\partial a_z}+\delta_{jz}\frac{\partial^2}{\partial a_j\partial a_p}+x_j\frac{\partial^3}{\partial a_j\partial a_p\partial a_z}\right).
\label{L3}
\end{align}
From the last expression it is clear that any 
inaccuracy in the computation of the derivatives  
adds up to the temporal truncation error with 
at least with factor $m^6$. Indeed, as shown 
in Figure \ref{fig:truncation}, the norm of $L^3 f$ 
grows faster than $m^6$. This has important consequences 
on the accuracy attainable with the AB2 time-integration 
scheme applied to high-dimensional linear PDEs.
\begin{figure}
\centerline{
\includegraphics[height=6.5cm]{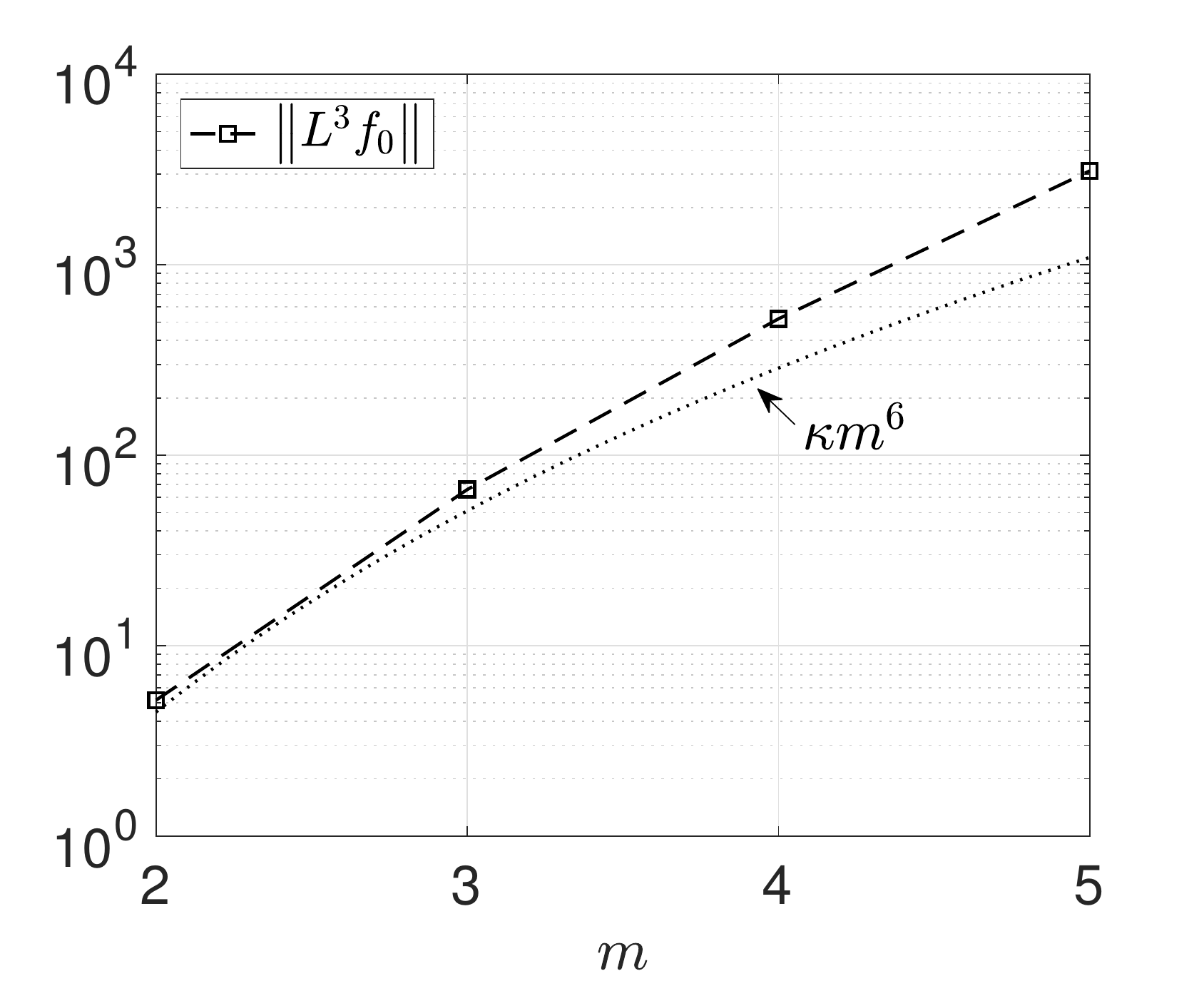}
}
\caption{\color{r} Norm of  $L^3 f$ at initial time versus 
$m$. The operator  $L^3$, defined \eqref{L3}, has been 
discretized  in space by using an accurate 
Fourier spectral method with $N=600$ points 
in each dimension. It is seen that the norm of $L^3 f_0$ grows 
faster than $m^6$. This has important consequences 
on the local truncation error generated by AB2 time-integrator 
applied to the PDF \eqref{PDE11}.}
\label{fig:truncation}
\end{figure}
In fact, suppose that $\left\| L^3 f\right\|$ grows 
like $m^6$ not just at initial time but at each time 
step (it actually grows faster than that), i.e.,  
\begin{equation}
\left\| L^3 f\right\| \sim \kappa m^6,
\label{te}
\end{equation}
where $\kappa$ is a suitable constant. If 
we want the AB2 scheme to operate 
at constant truncation error for different number of 
dimensions $m$ then we need to guarantee that 
\begin{equation}
m_1^6 \Delta t_1^3 = m_2^6 \Delta t_2^3.
\label{re}
\end{equation}
For example, if $m_1=2$ and $m_2=10$ then we have 
\begin{equation}
\Delta t_2 =  \frac{\Delta t_1}{25}.
\end{equation}
In other words, if we want our time integrator to operate at constant 
truncation error, then we need to run the simulation in $10$ dimensions
with a time step that is roughly $25$ times smaller 
than the one we emply in the simulation in $2$ dimensions.
This can tax the computational resources quite substantially. 
In fact, suppose we are interested in integrating our PDE up to $T=1$, 
and we set $\Delta t_1=10^{-4}$ in $2$ dimensions, i.e.,  
$10^4$ time steps. Assuming that the local truncation error 
is the same at each time step (see Eq. \eqref{te}), at the end 
of the integration period we accumulated an error of approximately
$10^{-6} \kappa$. The same error is roughly 
attained at  $T=0.04$ if we integrate the PDE in 
$10$ dimensions. 
In fact, the time step $\Delta t_2=\Delta t_1/25$ guarantees a 
constant local truncation error which adds up to $10^{-6} \kappa$
after just $0.04$ time units. The local truncation error 
manifest itself as {\em numerical diffusion} 
which eventually dissipates the numerical solution 
to zero. 
Note that in this simple calculation we did not take into account 
the accuracy of the rank reduction process in the CP-ALS 
and HT algorithms, which takes place at each time step.

\subsection{Numerical Results}
We solve the multivariate PDE \eqref{PDE11} in 
the hypercube $[-b,b]^m$, with $b=60$ and 
variable $m$. Such domain is chosen large 
enough to accommodate periodic (zero) boundary 
conditions in the integration period of interest. 
We study both the  HT (hierarchical Tucker) and the 
CP-ALS schemes we discussed in Section \ref{sec:HTD} 
and Section \ref{sec:CTD-ALS}. Specifically, we implemented
a Fourier collocation method with $600$ nodes in each variable 
and explicit AB2 time stepping.
To study the accuracy of the numerical solution, 
we consider the time-dependent relative error
\begin{equation}
\epsilon_m(t) = \left|\frac{F([\theta_m^*],t)-\hat{F}([\theta_m^*],t)}{F([\theta_m^*],t)}\right|,
\label{relative_error}
\end{equation}  
where $F$ is the analytical solution \eqref{eq:SF}, $\hat{F}$ 
is the numerical solution we obtained by using the 
CP-ALS or the HT  algorithms in the test function 
space $D_m$ (i.e., $m$ with independent variables) 
and with  separation rank $r$. 
The test function $\theta^*$ in \eqref{relative_error} 
is defined as 
 \begin{equation}
\theta_m^{*}(x)= h \sum_{j=1}^m \varphi_j(x),\qquad 
h= 0.698835274542439,
\label{point}
\end{equation}
where $\varphi_j(x)$ are the 
orthonormal polynomials shown 
in Figure \ref{fig:basis_eigen}.  
The accuracy of the CP-ALS and the HT algorithms  
is studied in Figure \ref{fig:ALS-CP-HT}, where 
we plot the relative pointwise error \eqref{relative_error}
for different separation ranks $r$ and for 
different number of dimensions. 
It is seen that, as expected,
the accuracy of the numerical solution 
increases as we increase the separation rank. 
Also, as we increase the number of dimensions from 
$2$ to $6$ the relative error increases, in agreement 
with the results of Section \ref{sec:long-term-integration}
(we emply a constant $\Delta t=2.55\times 10^{-4}$ 
in all our simulations).
It is worthwhile emphasizing that the CP-ALS 
is a randomized algorithm which requires 
initialization at each time-step. This means that 
results of simulations with the same nominal 
parameters may be different. On the other hand, tensor
methods based on multivariate/distributed singular 
value decomposition, such as the hierarchical Tucker 
decomposition \cite{Grasedyck2017}, do 
not suffer from this issue. 
The CP-ALS algorithm is faster than HT but, as we just said,
accuracy control may be an issue.  
The separation rank of both the CP-ALS 
and HT algorithms are computed adaptively up 
to the maximum value $r_{max}$ specified 
in the legend of Figure \ref{fig:ALS-CP-HT}. 
\begin{figure}
\centerline{\hspace{0.5cm}CP-ALS\hspace{6.5cm} HT}
\noindent
\centerline{
\rotatebox{90}{\hspace{2cm}two dimensions}\hspace{0.2cm}
\includegraphics[height=5.5cm]{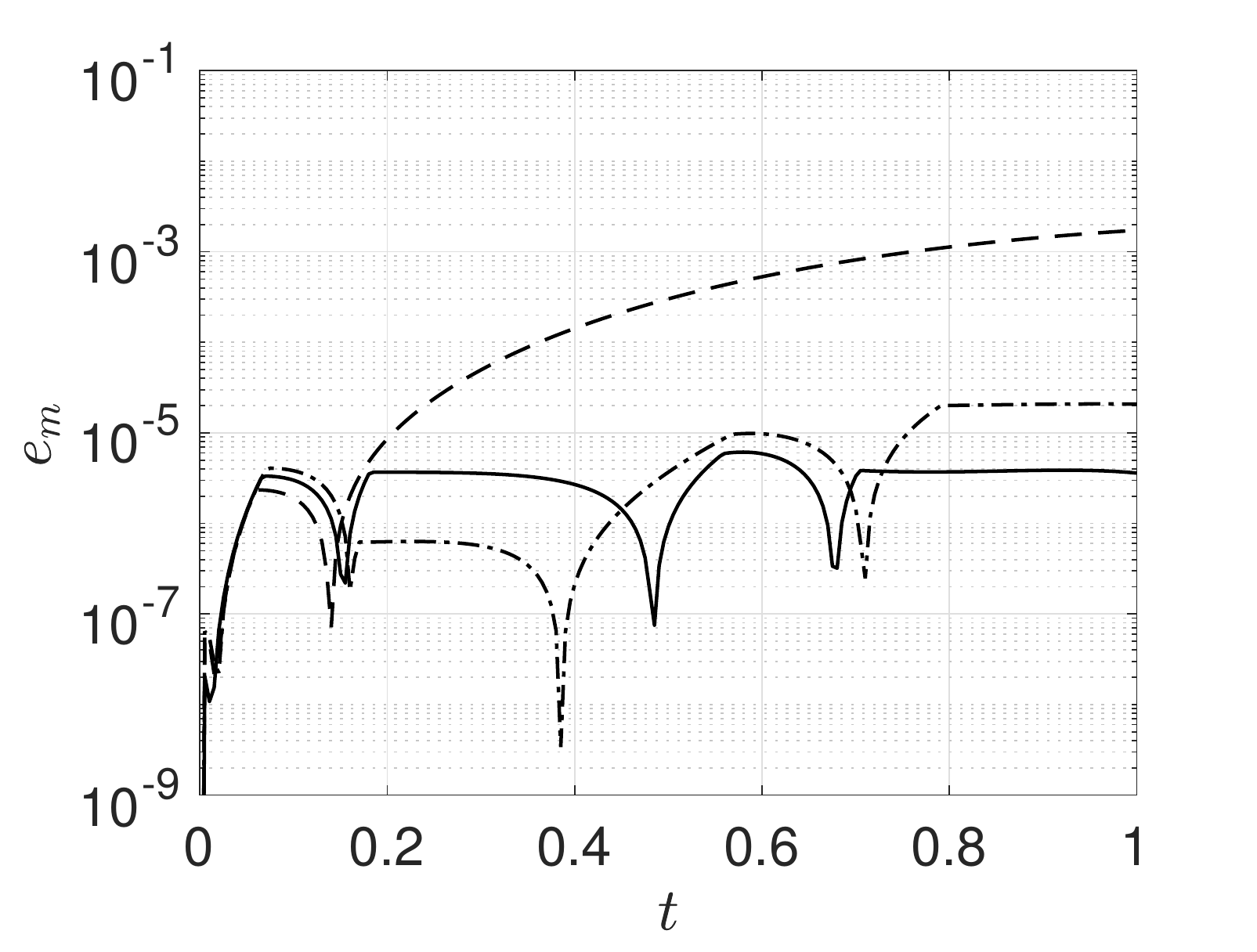}
\includegraphics[height=5.5cm]{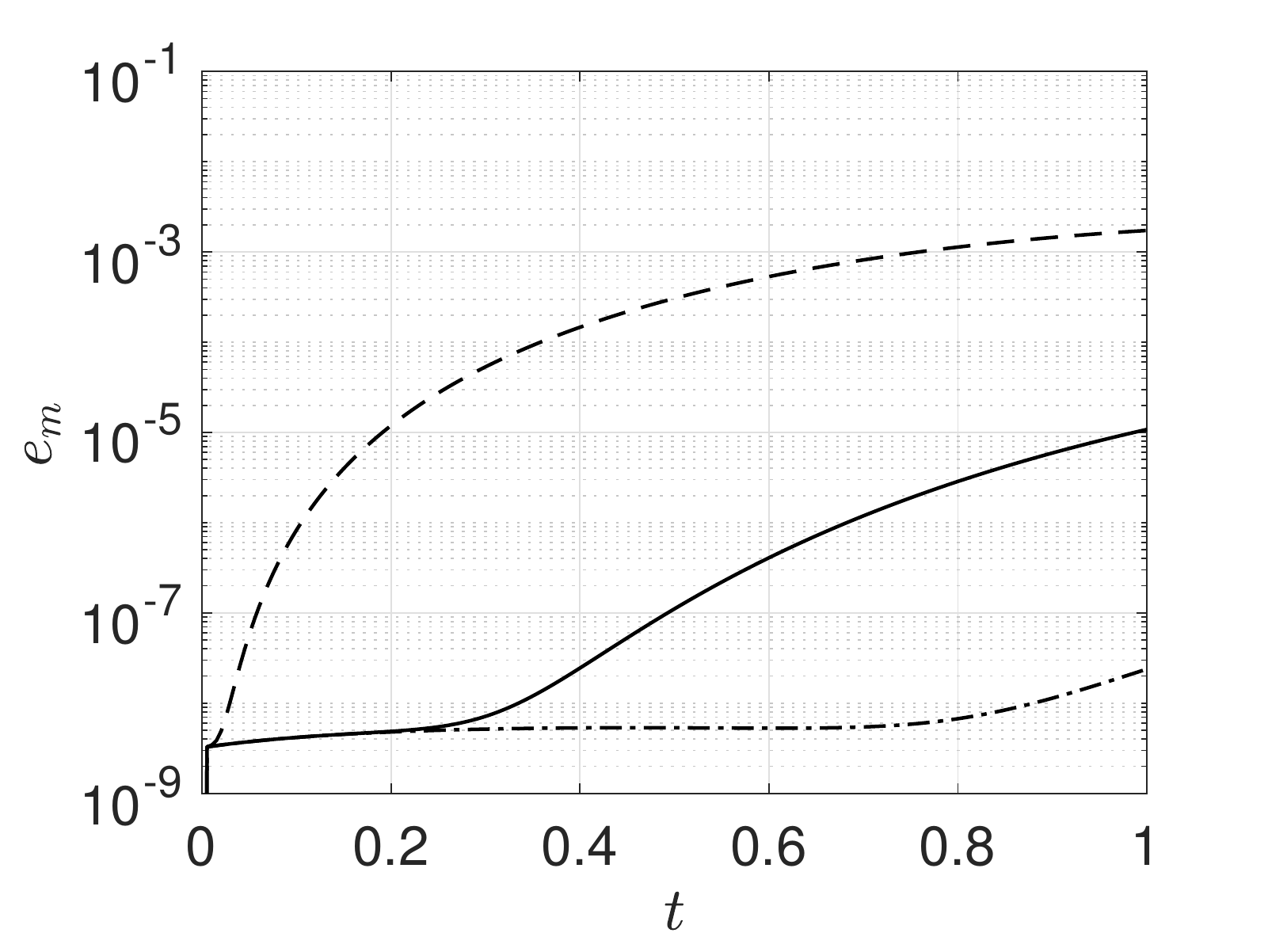}
}
\noindent
\centerline{
\rotatebox{90}{\hspace{2cm}three dimensions}\hspace{0.2cm}
\includegraphics[height=5.5cm]{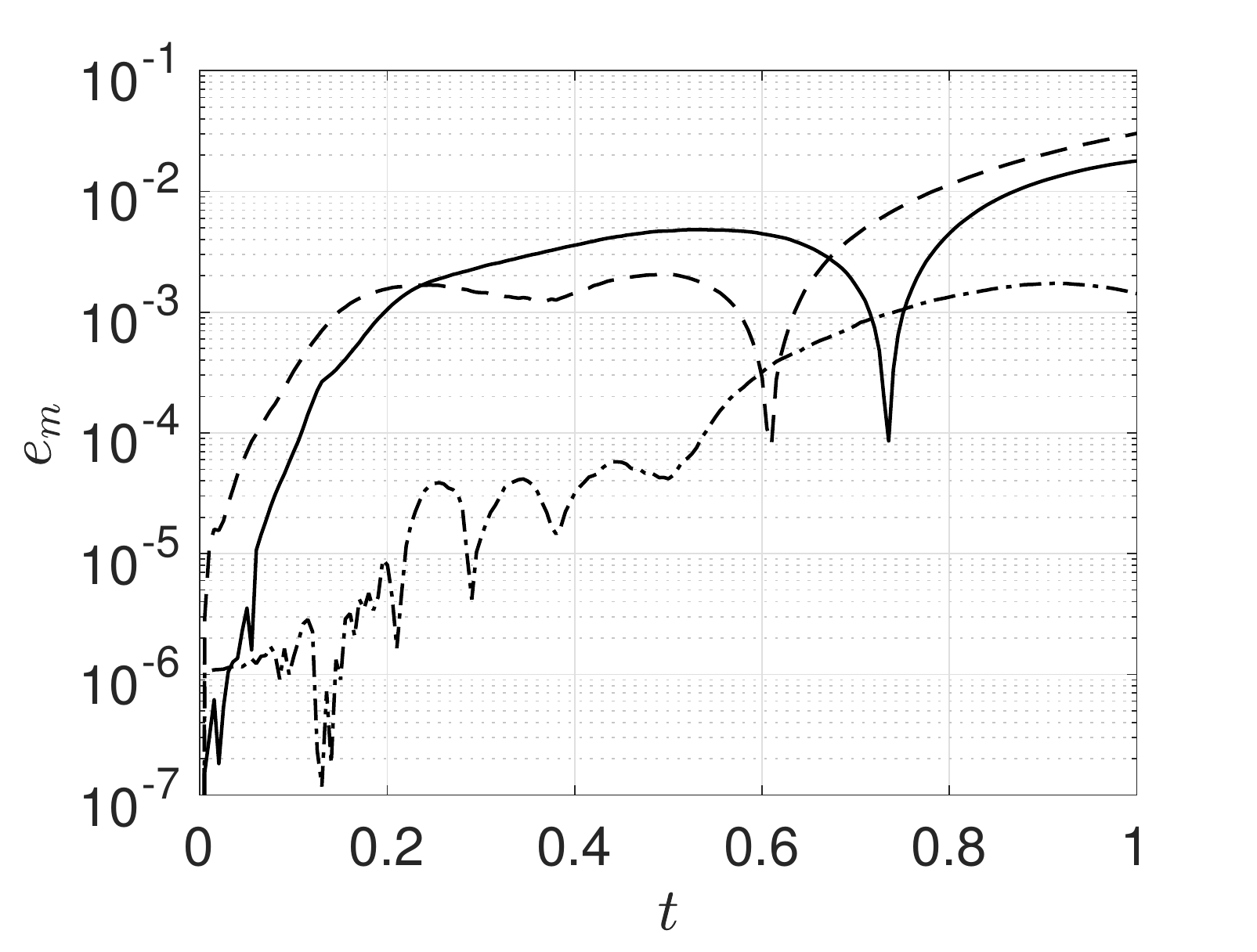}
\includegraphics[height=5.5cm]{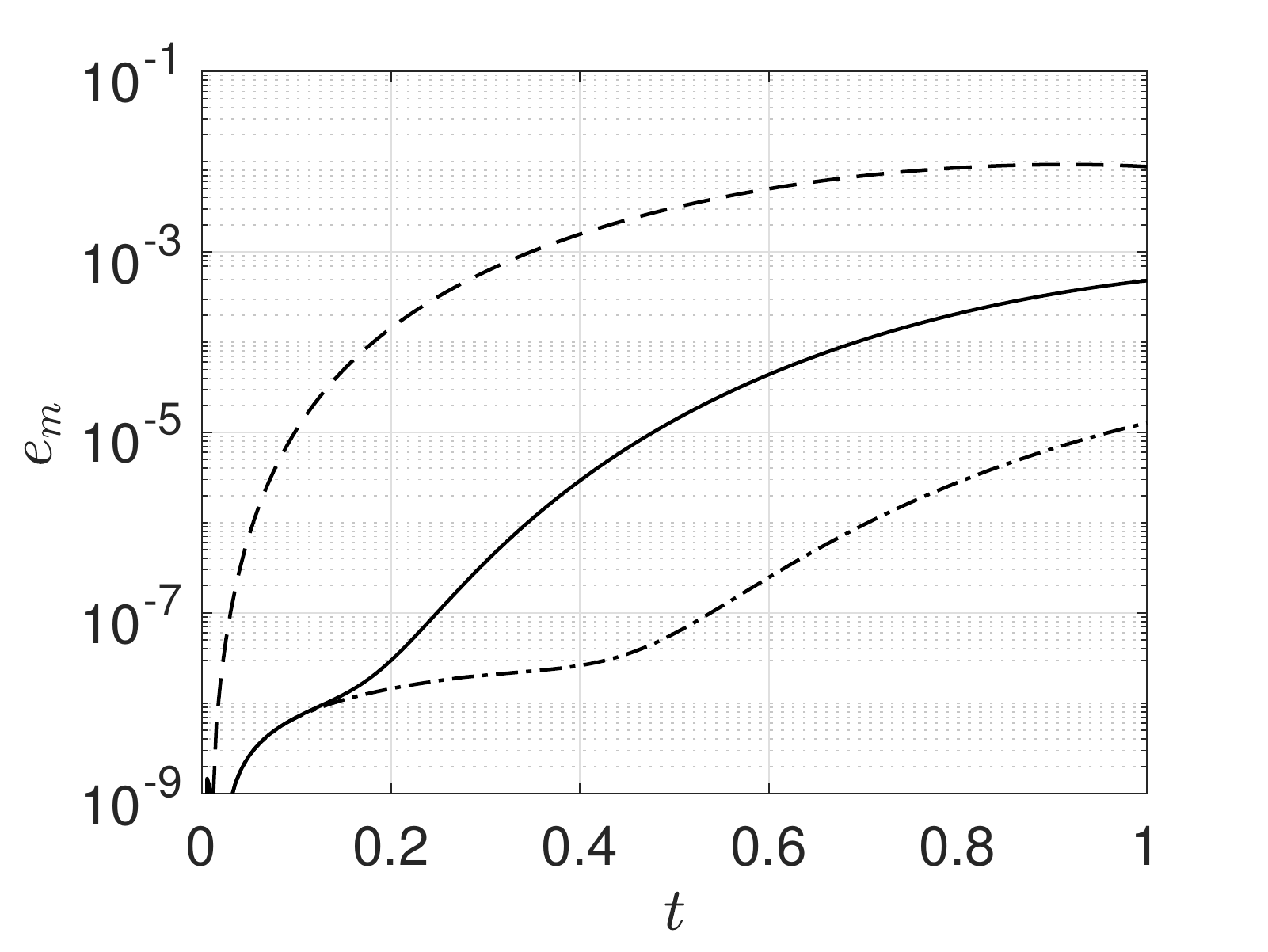}
}

\centerline{
\rotatebox{90}{\hspace{2cm}six dimensions}\hspace{0.2cm}
\includegraphics[height=5.5cm]{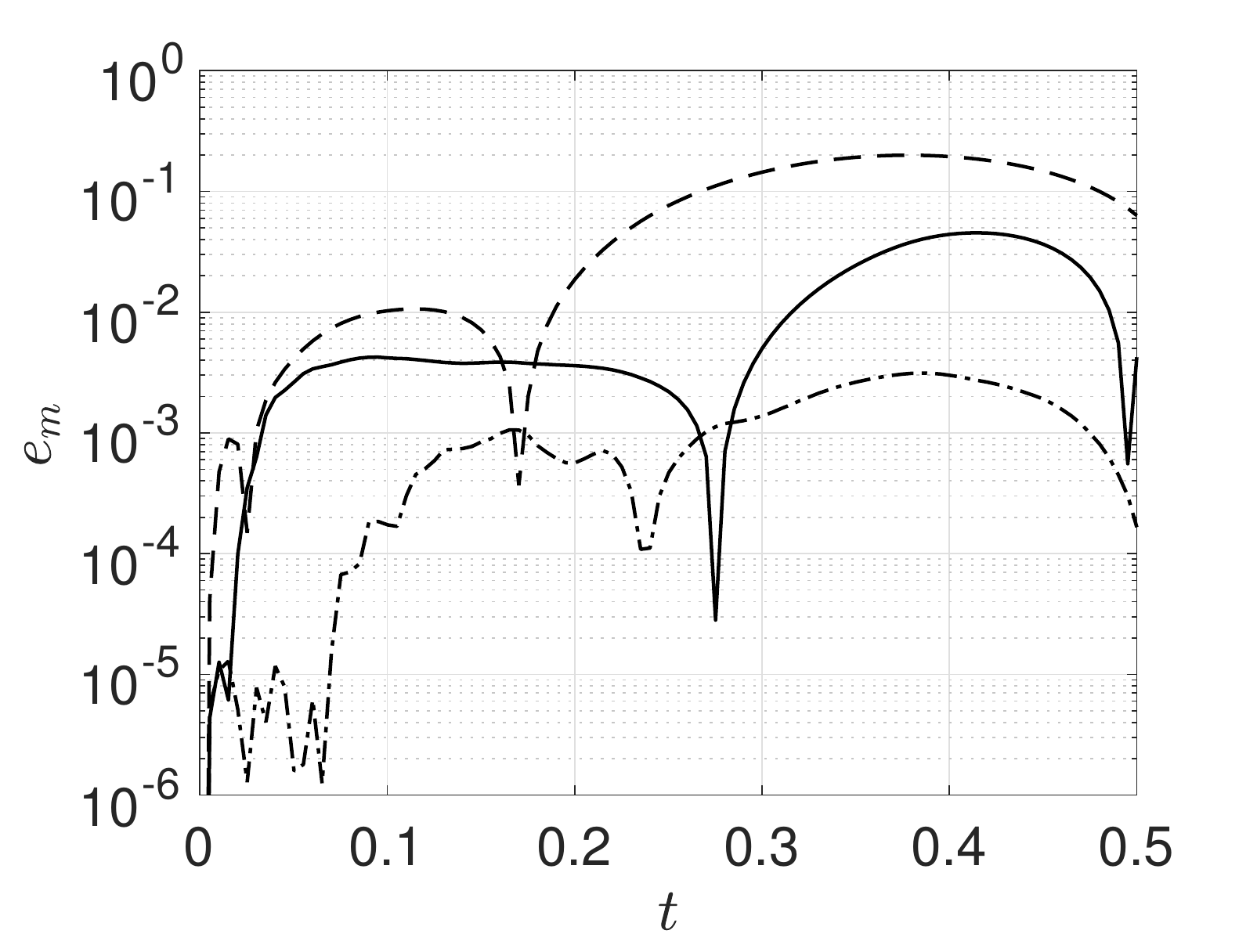}
\includegraphics[height=5.5cm]{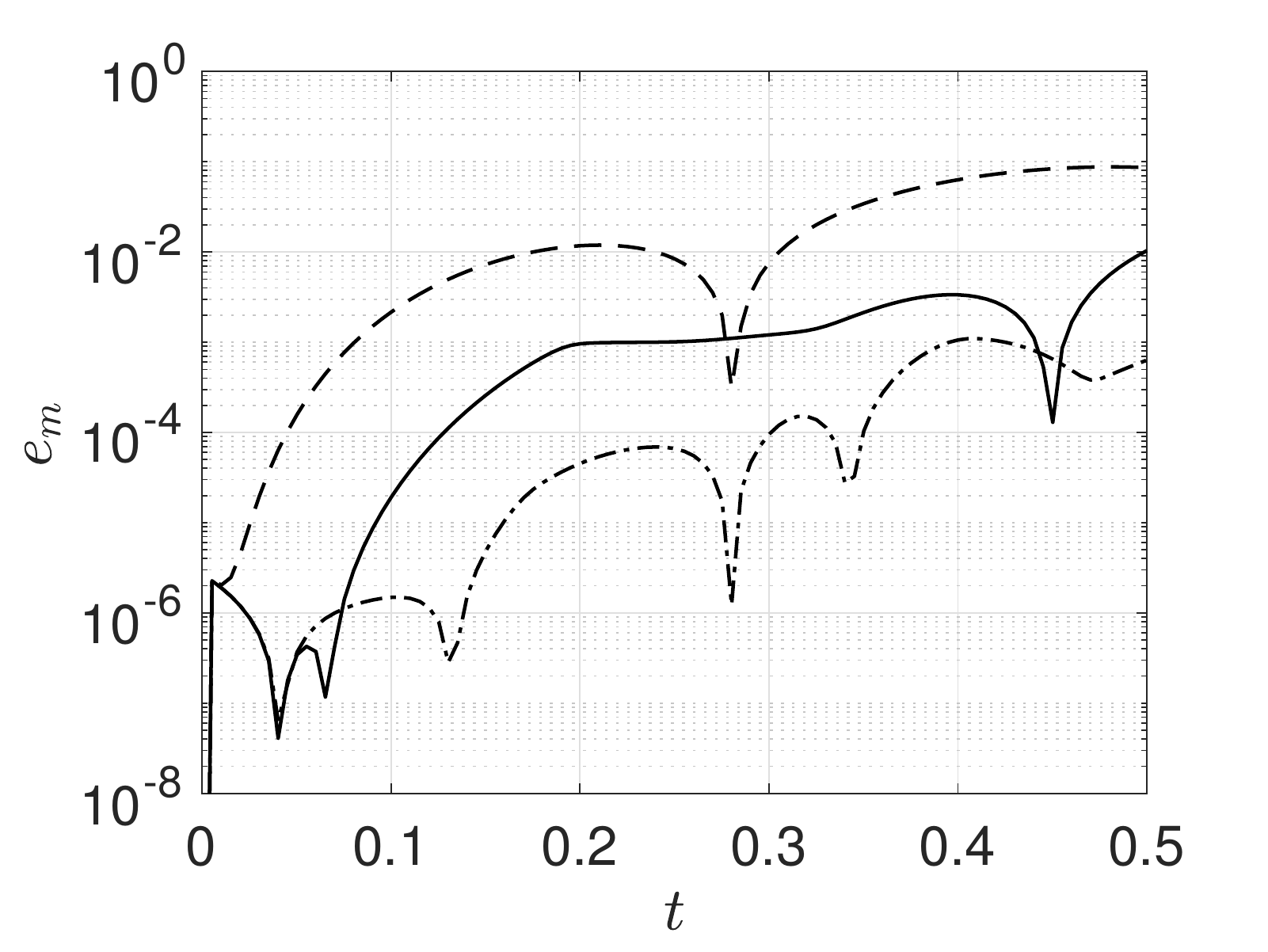}
}
\centerline{$--$ $r_{max}=4$\hspace{1cm} --- $r_{max}=8$
\hspace{1cm} $-\cdot -$ $r_{max}=12$
}
\caption{\color{r}Accuracy of the numerical solution 
to the FDE \eqref{FDE3} in the finite dimensional function space 
\eqref{finite_dim_FS} for different number of dimensions. 
Specifically, we plot the relative $L_\infty$ 
error \eqref{relative_error} versus time 
we obtained by using the CP-ALS and HT algorithms 
we discussed in Section \ref{sec:numericalDiscretizationADVR}. 
The separation rank of both CP-ALS and HT is computed adaptively 
at each time step up to the maximum value $r_{max}$.}
\label{fig:ALS-CP-HT}
\end{figure}
Note that in two dimensions the CP-ALS algorithm 
results in error plots that looks very similar when 
$r_{max}=8$ and $r_{max}=12$. This is because the 
separation rank is less than $8$ in both cases 
throughout the simulation up to $t=1$. 
The variability of the results is related to the 
random initialization required by the ALS algorithm 
at each time step. On the other hand, 
in three and six dimensions the error plots we obtain 
for  $r_max=4$ and $r_max=8$ are of the same 
order of magnitude because such separation ranks 
are achieved after just few time steps and they 
are not sufficient to accurately represent the 
multivariate solution. 
A similar phenomenon attributable to the separation rank 
is observed in the HT simulations. In particular, 
the time instant at which the HT tensor series 
requires a separation rank higher than $8$ can be clearly 
identified, i.e., at $t=0.2$ in two dimensions, at $t=0.1$ 
in three dimensions, and at $t=0.05$ in six dimensions.  
In Figure \ref{fig:CP_MODES} we plot 
the time evolution of the first three CP 
modes $G_6^1$, $G_6^2$ and $G_6^3$ 
(see equation \eqref{CP-numerical}) 
representing the dynamics of the solution functional in 
the variable $a_6=(\theta,\varphi_6)$. 
\begin{figure}
\centerline{
\rotatebox{90}{\hspace{1.8cm}$t=0.01$ }\hspace{0.7cm}
\includegraphics[height=4.3cm]{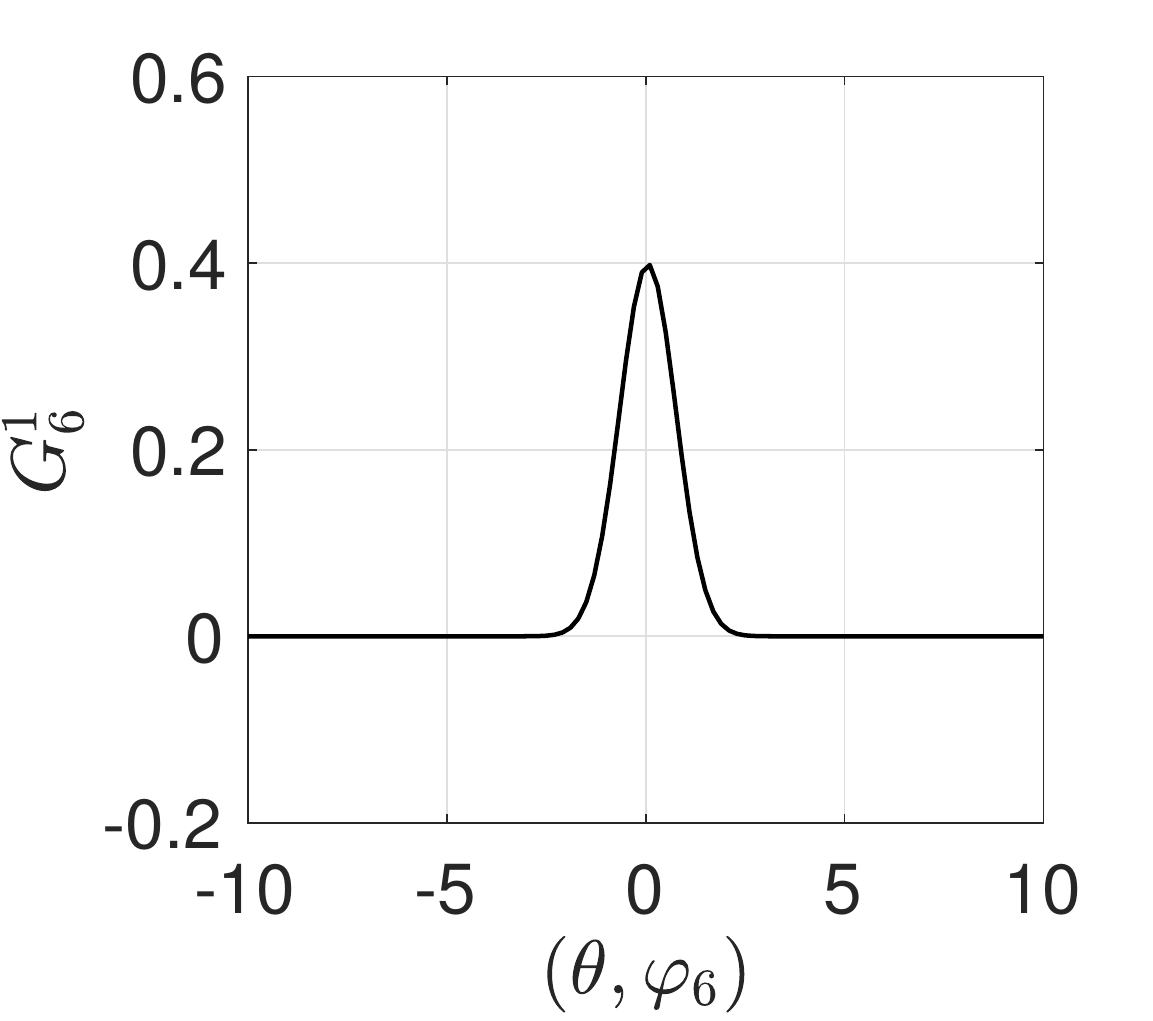}
\includegraphics[height=4.3cm]{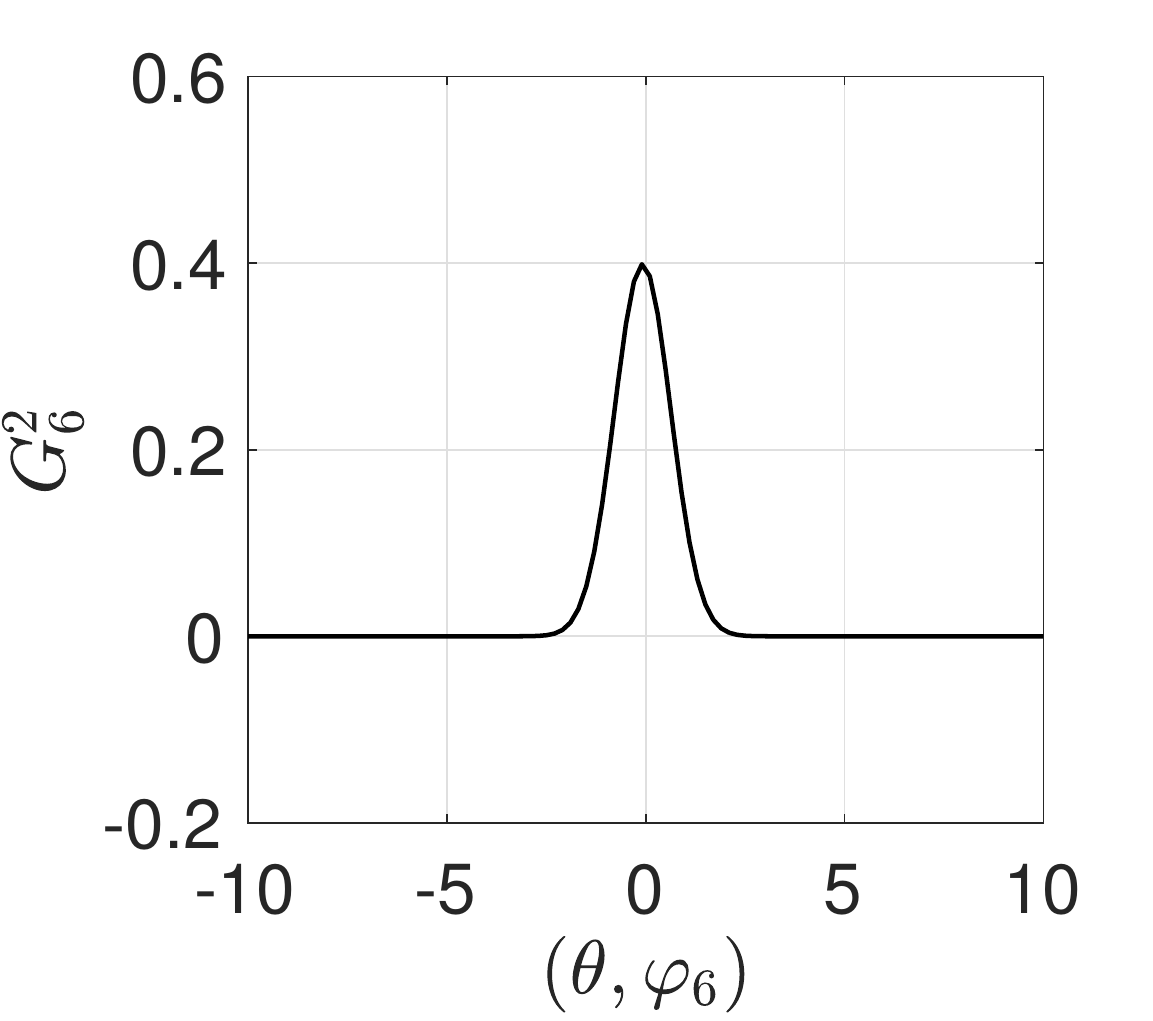}
\includegraphics[height=4.3cm]{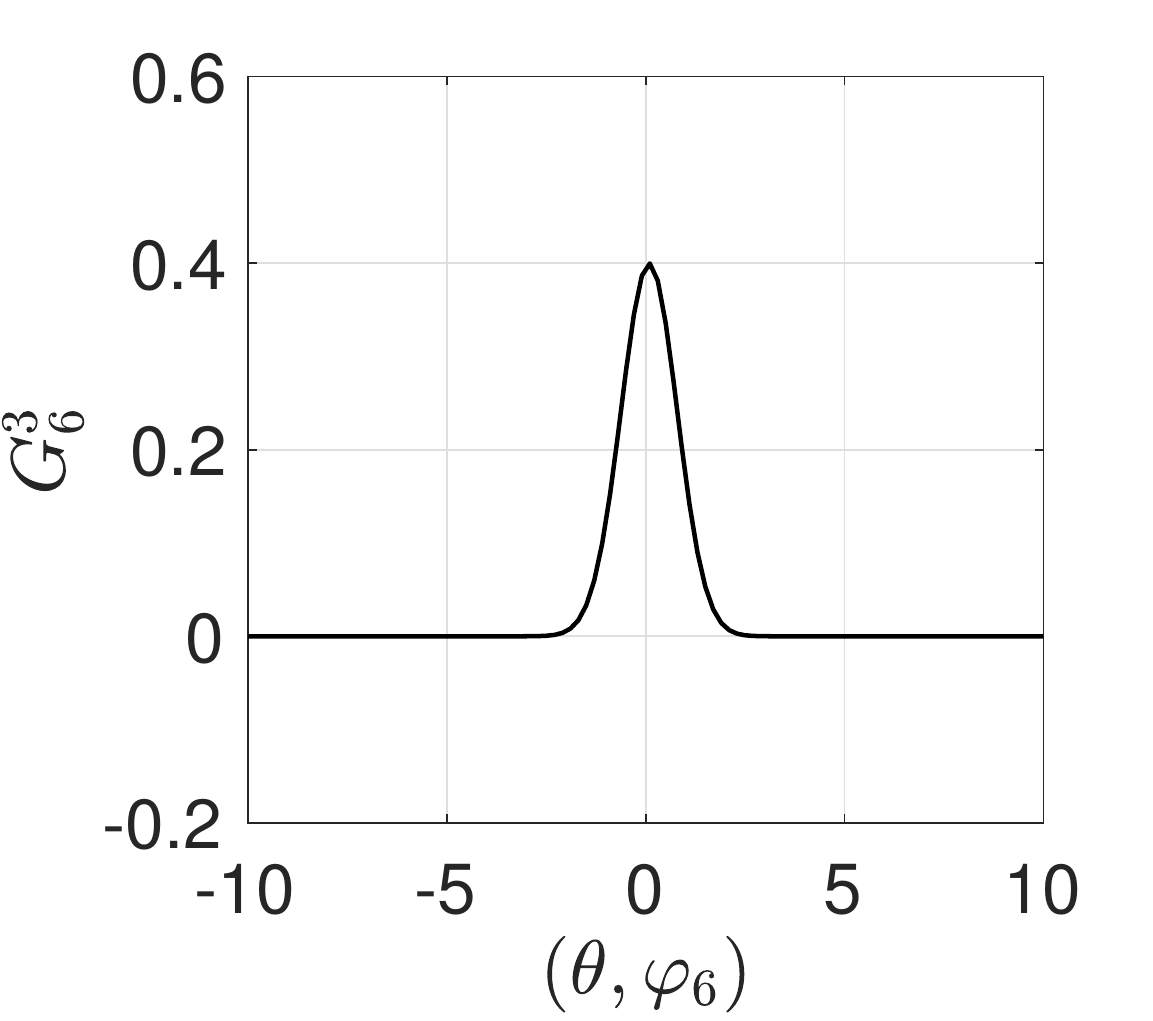}
}
\centerline{
\rotatebox{90}{\hspace{1.8cm}$t=0.5$ }\hspace{0.7cm}
\includegraphics[height=4.3cm]{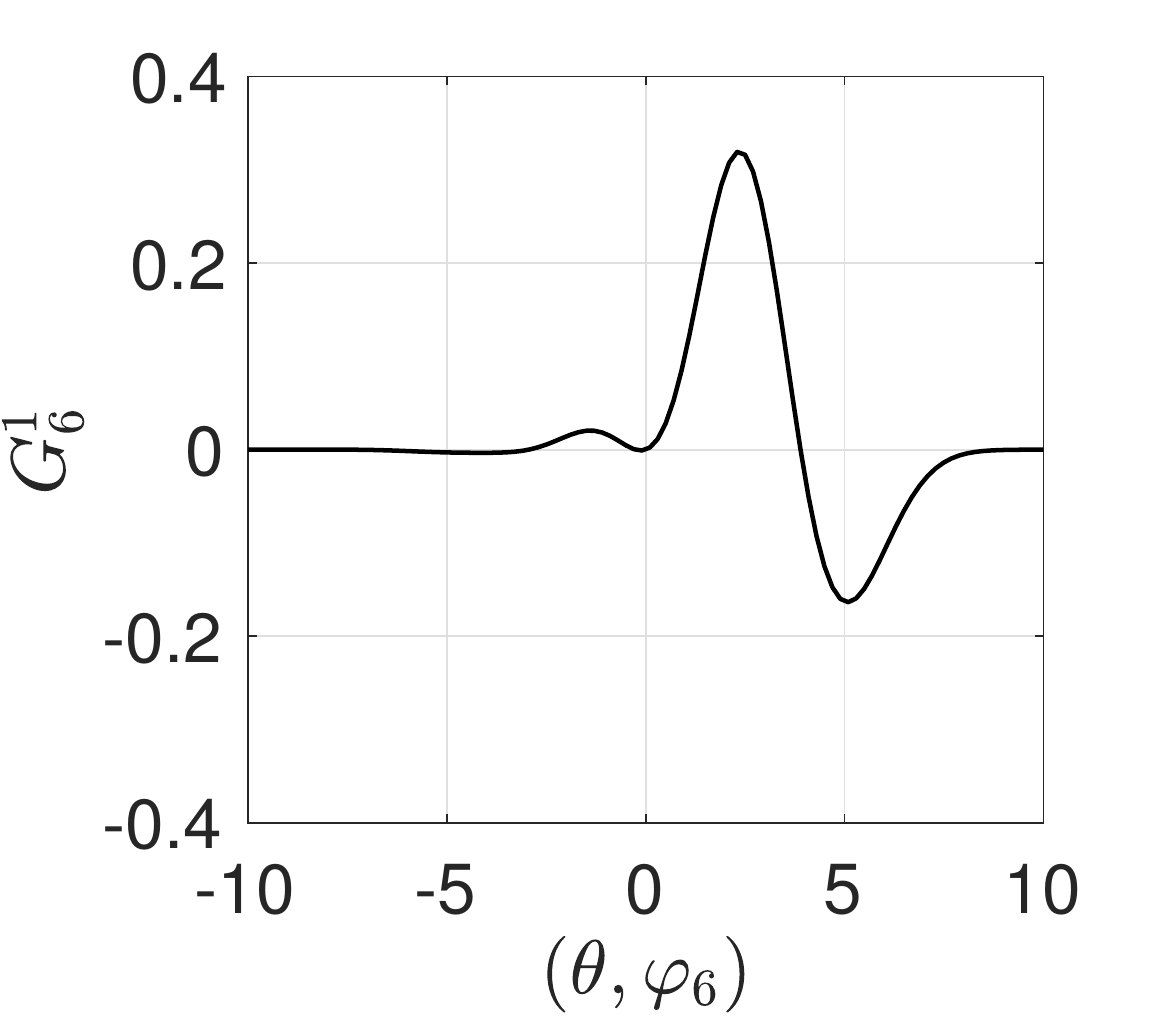}
\includegraphics[height=4.3cm]{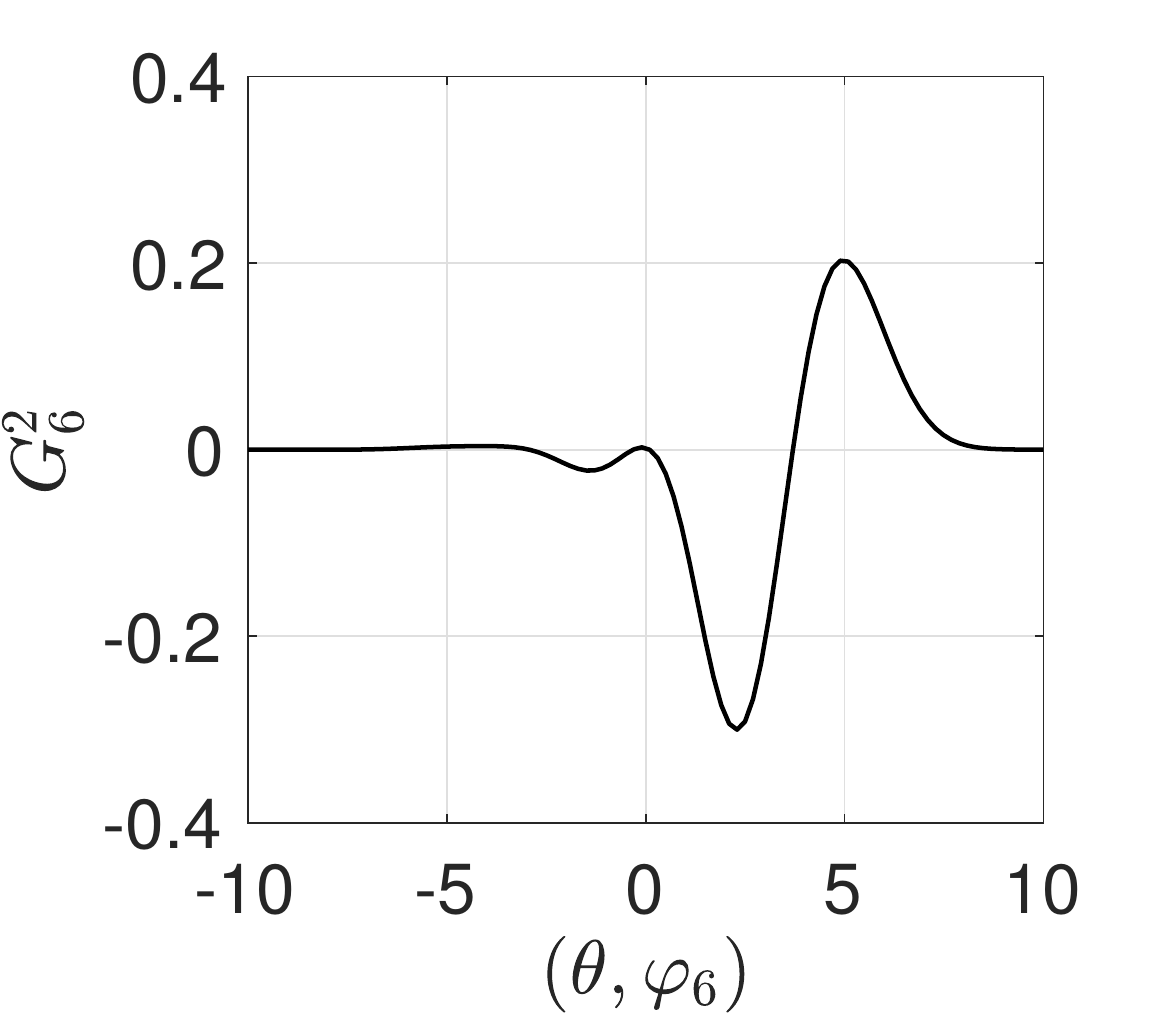}
\includegraphics[height=4.3cm]{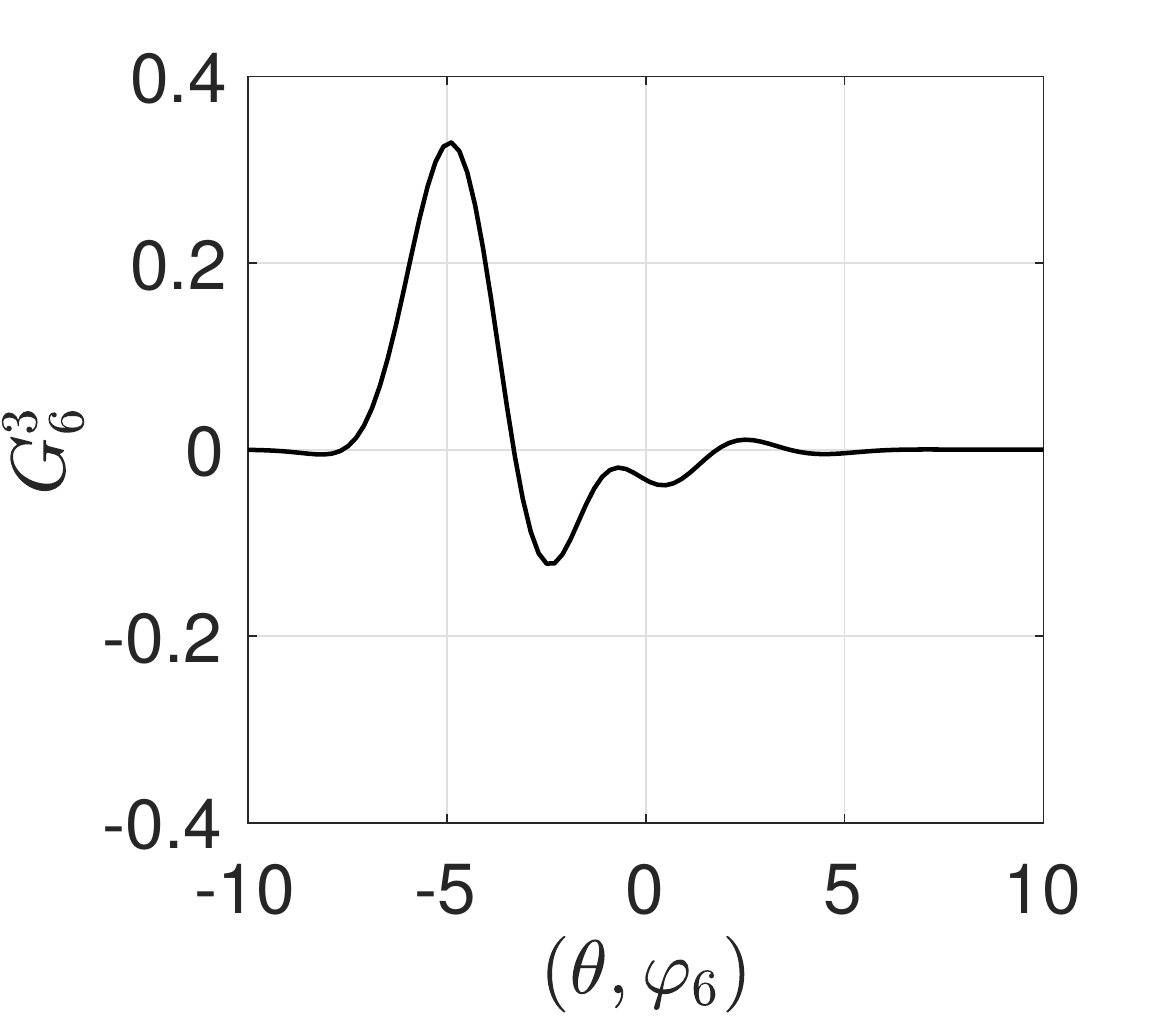}
}
\caption{\color{r}Time evolution of the first three CP 
modes representing the dynamics along $a_6=(\theta,\varphi_6)$ 
in the CP expansion \eqref{CP-numerical} of the 
solution functional. These modes 
are obtained by solving numerically the FDE \eqref{FDE3} 
in the function space $D_6$ (see Eq. \eqref{finite_dim_FS}) 
with the CP-ALS algorithm.}
\label{fig:CP_MODES}
\end{figure}

\vs
\noindent
{\em Remark:} If $\theta(x)$ is not in the function space 
$D_m$  (see Eq. \eqref{finite_dim_FS}), but can be represented in $D_m$ with accuracy then 
the solution to the multivariate PDE \eqref{PDE11} 
provides an approximation of the solution 
to the full FDE \eqref{FDE3} at such $\theta(x)$. 
For instance, consider the following test function 
\begin{equation}
\theta(x)= \sin(x).  
\end{equation}
The Fourier coefficients of $\theta$ relative to the 
orthonormal basis shown in Figure \ref{fig:basis_eigen}
are  plotted in Figure \eqref{fig:spectrum_and_accuracy}. 
Clearly $\theta(x)$ is not in $D_6$, 
but it can be approximated well in $D_6$. This means 
that the solution to the FDE at $\sin(x)$ can be approximated 
by the solution of the six-dimensional PDE arising when we
 evaluate the FDE \eqref{FDE3} $D_6$.
\begin{figure}
\centerline{(a) \hspace{7cm} (b)}
\centerline{
\includegraphics[height=6cm]{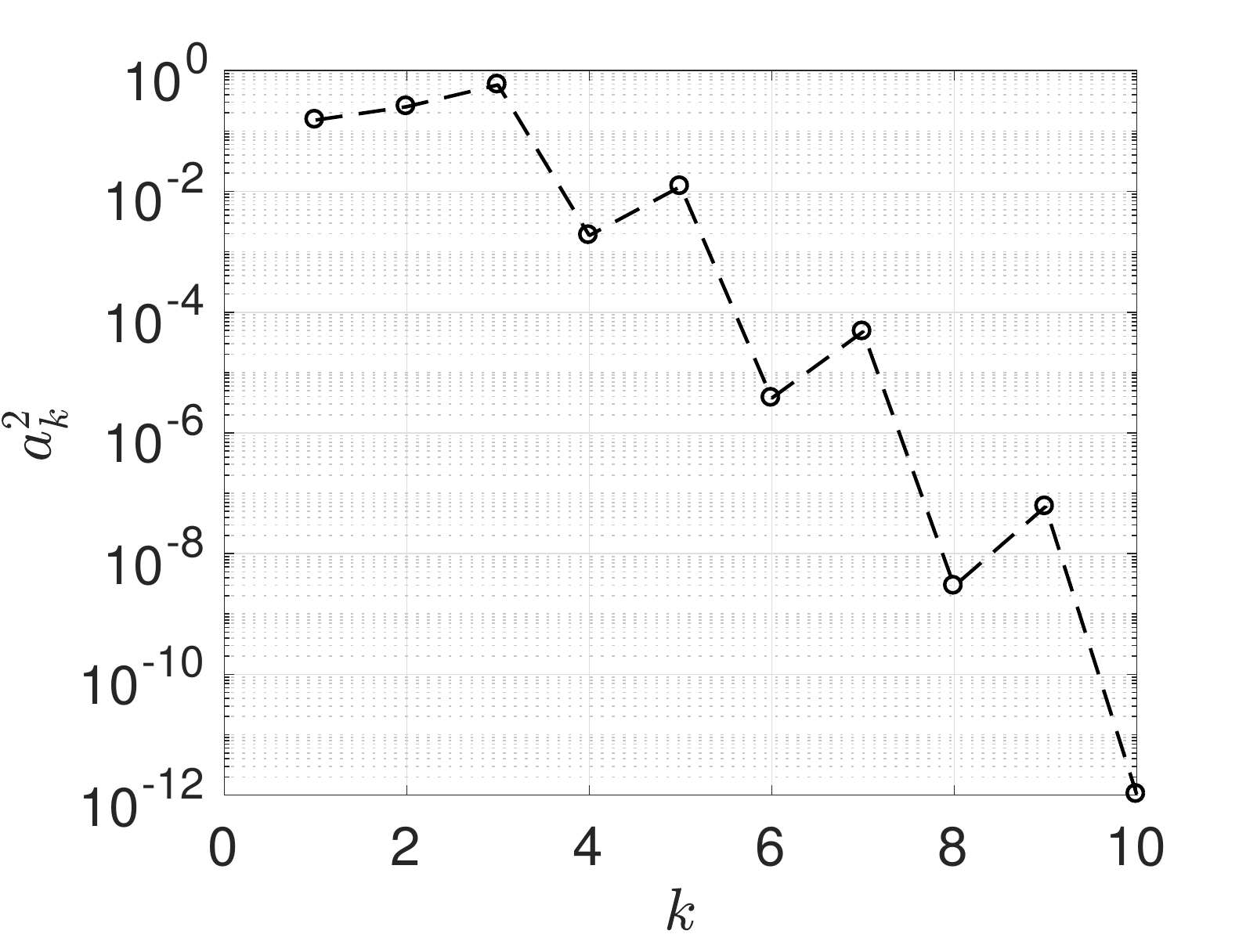}
\includegraphics[height=6cm]{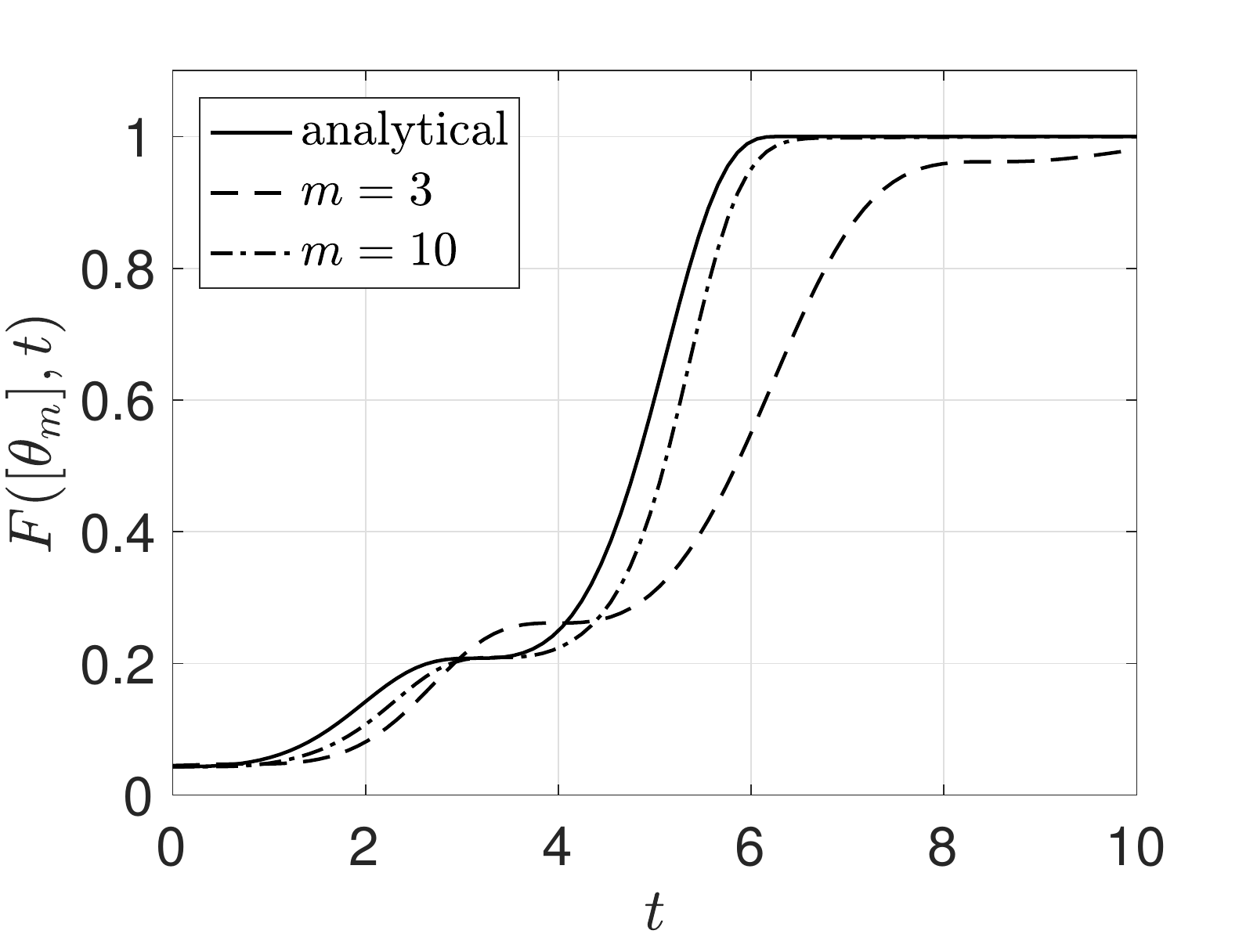}
}
\caption{\color{r}(a) Fourier  coefficients of $\sin(x)$ relative to the 
orthonormal polynomial basis shown in Figure \ref{fig:basis_eigen}. (b) 
Comparison between the solution to the FDE 
in $D_3$ and $D_{10}$ and the exact solution 
\eqref{THESOLUTION} evaluated at $\theta(x)=\sin(x)$.}
\label{fig:spectrum_and_accuracy}
\end{figure}

\paragraph{Computing Functional Derivatives} 
In Figure \ref{fig:functional derivative} we compare the 
exact first-order functional derivative \eqref{FunctionalDerivative} 
at $t=0.4$ with the numerical approximation we obtained 
with the HT algorithm in $m=3$ and $m=6$ dimensions. 
The separation rank of the tensor series expansions is set 
to $r=12$. As expected, the numerical approximation 
converges to the exact solution as we increase 
the number of dimensions. 
\begin{figure}[t]
\centerline{\hspace{0.7cm}$\theta(x)=\sin(x)$ \hspace{5.2cm} 
$\theta(x)=xe^{-x^2/4}$}
\noindent
\centerline{
\includegraphics[height=5.5cm]{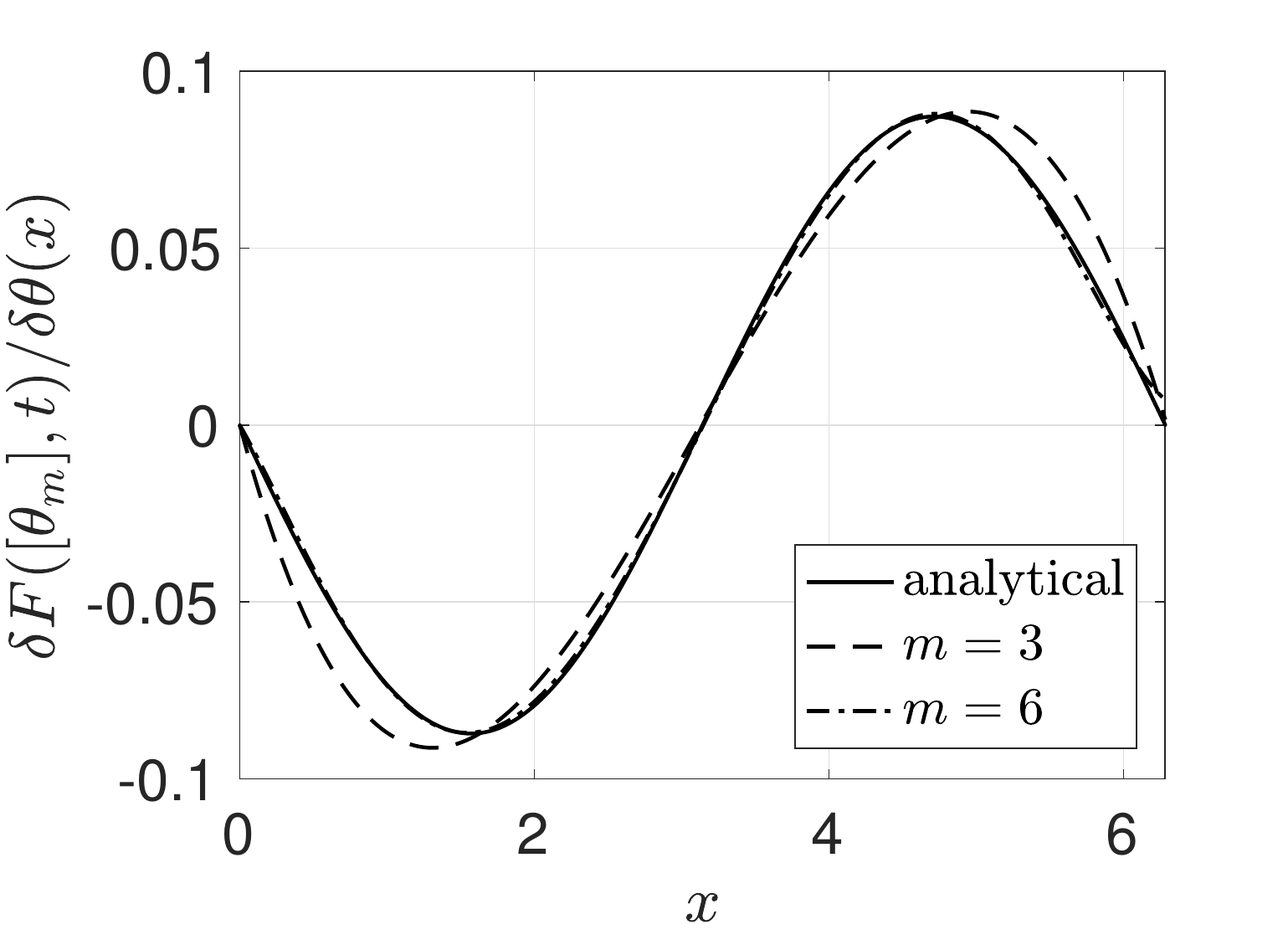}
\includegraphics[height=5.5cm]{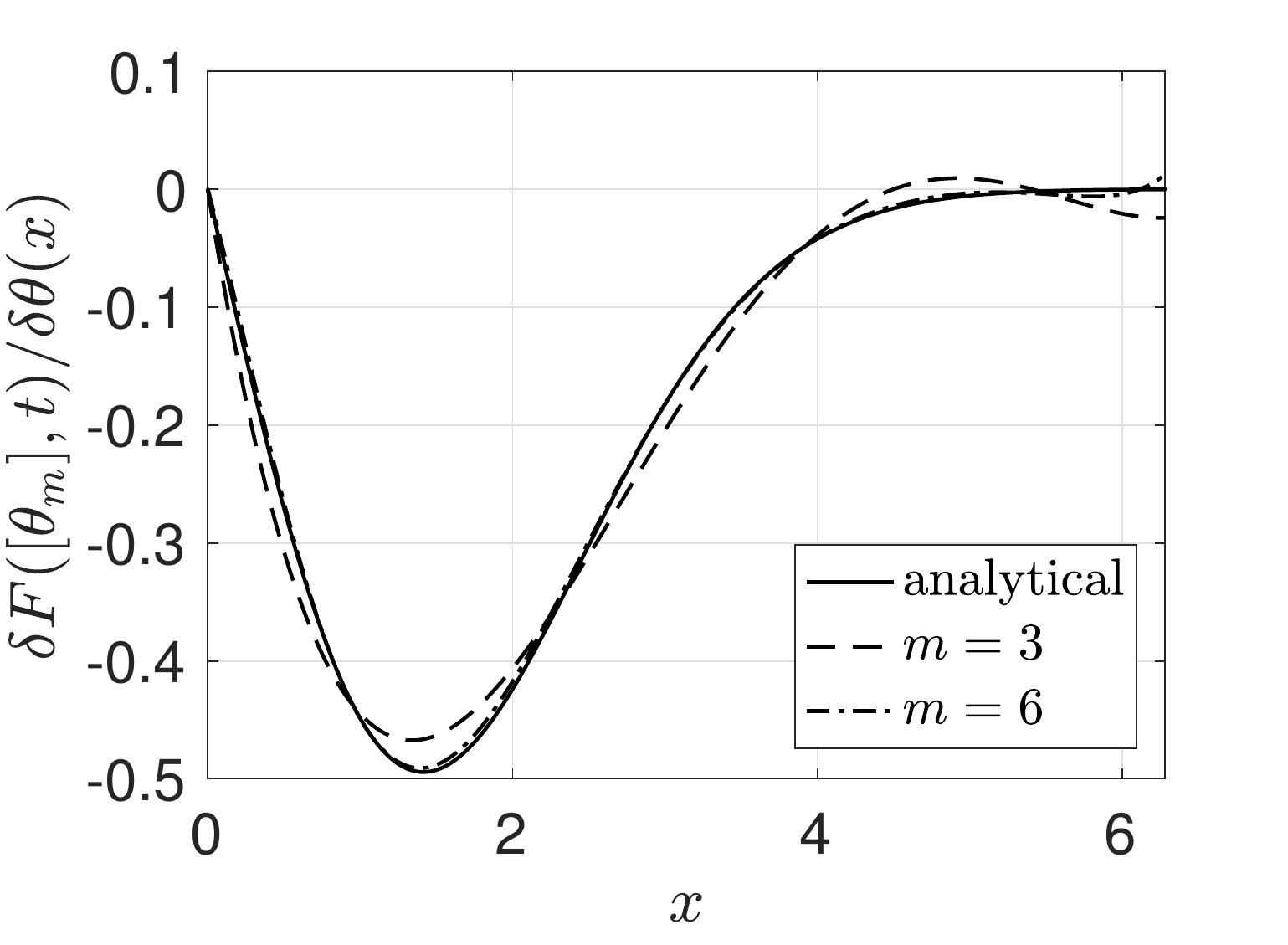}
}
\caption{\color{r}First-order functional derivative of the solution 
functional \eqref{THESOLUTION} 
at $t=0.4$, evaluated at the test functions $\theta(x)=\sin(x)$ and 
$\theta(x)=xe^{-x^2/4}$.
We plot the exact result \eqref{FunctionalDerivative} and 
the numerical approximation \eqref{FDeRT} we obtained 
in the function spaces $D_3$ and $D_6$ 
(see Eq. \eqref{finite_dim_FS}). The numerical approximation is 
computed by solving the mutivariate PDE \eqref{PDE-advR} in 
three and six dimensions.}
\label{fig:functional derivative}
\end{figure}

\subsection{The Navier-Stokes-Hopf Functional Equation}
\label{sec:Navier-Stokes-Hopf}
In this Section we discuss approximation of the 
Navier-Stokes-Hopf functional equation  
\begin{equation}
\frac{\partial \Phi([\bm \theta],t)}{\partial t}=\sum_{k=1}^3\int_V
\theta_k(\bm x)\left[i\sum_{j=1}^3\frac{\partial }{\partial x_j}
\left(\frac{\delta^2 \Phi([\bm \theta],t)}{\delta \theta_k(\bm x)
\delta\theta_j(\bm x)}\right)+\nu\nabla^2 \left(\frac{\delta 
\Phi([\bm \theta],t)}{\delta \theta_k(\bm x)}\right) \right]d \bm 
x.
\label{h1}
\end{equation}
In this formulation, $V$ is a periodic box 
and $\bm \theta(\bm x)$ is chosen in a 
divergence-free space of test functions. 
The main advantage of using such 
divergence-free space is that the pressure term
drops out, just as in the classical Navier-Stokes 
equations\footnote{\color{r}
We recall that the pressure functional in 
the Navier-Stokes-Hopf equation is defined as 
(see \cite{Monin2}, p. 749)
\begin{equation}
\Pi([\bm \theta],\bm x,t) = \left<p(\bm x,t;\omega)\exp
\left[i\int_V \bm \theta(\bm x)\cdot \bm u(\bm x,t;\omega) d\bm 
x\right]\right>.
\end{equation}
Therefore, 
\begin{align}
\int_V \bm \theta\cdot\left<\nabla p \exp
\left[i\int_V \bm \theta\cdot \bm u d\bm x\right]\right>d\bm x=&
\left<\int_V \bm \theta\cdot \nabla p d\bm x \exp
\left[i\int_V \bm \theta\cdot \bm u d\bm x\right]\right>.
\end{align}
If $\bm \theta$ is divergence-free, i.e., 
$\nabla\cdot \bm \theta = 0$ in $V$, and satisfies the 
tangency condition $\bm \theta\cdot \hat{\bm n}=0$ 
($\hat{\bm n}$ is the outward unit vector orthogonal 
to the boundary of $V$) we have
\begin{align}
\int_V \bm \theta(\bm x)\cdot \nabla p(\bm x,t;\omega) d\bm x = 
\int_{\partial V} p(\bm x,t;\omega) \bm \theta(\bm x)\cdot 
\hat{\bm n} d\bm x = 0.
\label{integg}
\end{align}
The integral \eqref{integg} is zero also 
if $\bm \theta$ is divergence-free and the 
domain $V$  is a periodic box. 
In fact, in this case $\bm \theta(\bm x)$ is periodic, 
$p$ is periodic, and therefore the integral 
along the boundary $\partial V$  vanishes.}. 
To develop a discretization of the Hopf equation \eqref{h1}
it is convenient to first address the question of 
how to represent divergence-free spaces of 
periodic functions in 2D and 3D.

\subsubsection{Symmetries of the Solution Functional} 
The divergence-free constraint in the velocity field 
induces a certain number of symmetries in the 
Hopf functional. 
Let first assume that the flow develops in a bounded 
region $V$, i.e., $\bm u\cdot \widehat {\bm n}=0$, 
where $\widehat {\bm n}$ is the outward unit vector 
normal to the boundary of $ V$. 
Such region could be any volume co-moving with the fluid. 
In this assumption, 
\begin{equation}
\Phi([\bm \theta+\nabla \varphi]) = \Phi([\bm \theta])\label{div}
\end{equation}
for any $\bm \theta(\bm x)$  and any $\varphi (\bm x)$. In fact, 
\begin{equation}
\int_V \bm u \cdot \nabla \varphi d\bm x= \int_{\partial V}\varphi \bm u\cdot \hat {\bm n} d\bm x=0.
\label{29}
\end{equation}
If the boundary conditions assigned 
to $\bm u$ on $\partial V$ (boundary of $V$) are different from 
$\bm u\cdot \hat{\bm n}=0$, then \eqref{div} is still 
valid, but not for any $
\varphi(\bm x)$. For example, if $\bm u$ is not 
orthogonal to $\hat{\bm n}$ 
along some part of the boundary $\partial V$, then it is sufficient to set $ \varphi=0$ along such boundary and use the fact that $\nabla \cdot \bm u=0$ to conclude that \eqref{div} is still valid.
Now let us consider the Helmholtz decomposition
\begin{equation}
\bm \theta =\bm \eta + \nabla \varphi, 
\label{dec}
\end{equation}
where $\nabla \cdot \bm \eta=0$ in the domain $V$, and $\varphi = 0$ at the boundary $\partial V$. By taking the divergence of \eqref{dec} we conclude that  $\varphi$ satisfies the Dirichlet boundary value problem 
\begin{equation}
\nabla^2 \varphi = \nabla\cdot \bm \theta\quad \textrm{(in $V$)},\qquad 
\varphi=0  \quad \textrm{(in $\partial V$)},
\label{dvp}
\end{equation}
which has a solution. With the decomposition \eqref{dec} available, we have
\begin{equation}
\Phi([\bm \eta+\nabla \varphi ]) = \Phi([\bm \eta]).
\label{symm}
\end{equation}
In fact,  
\begin{equation}
\int_V (\bm \eta+\nabla \varphi) \cdot \bm u d\bm x= 
\int_V \bm \eta\cdot \bm u d\bm x+ 
\int_V \nabla \varphi \cdot \bm u d\bm x.
\end{equation} 
However, by applying the Gauss theorem 
(recall that $\nabla \cdot \bm u=0$)
\begin{equation}
\int_V \nabla \varphi \cdot \bm u dV=\int_{\partial V} \varphi
\bm u\cdot \widehat {\bm n}  d\bm x =0
\label{w1}
\end{equation}
Therefore \eqref{symm} holds. The boundary conditions for $\bm \eta=\bm 
\theta-\nabla \varphi$ depend on the gradient of the solution 
to \eqref{dvp} at the boundary $\partial V$.
Next, consider the unique Helmholtz-Hodge 
decomposition (see \cite{Bhatia} or \cite{ChorinMarsden}, p. 36) of the 
field $\bm \theta$ in the form \eqref{dec},
where  $\nabla \cdot \bm \eta=0$ and $\bm \eta$ is tangent to $\partial V$, 
i.e., $\bm \eta\cdot \hat{\bm n}=0$ ($\hat{\bm n}$ outward unit vector normal 
to $\partial V$). 
If there is no flow across the boundary $\partial V$ ($\bm u\cdot 
\hat{\bm n}=0$), or if the boundary conditions of $\varphi$ are chosen such that 
the integral at the right hand side of \eqref{w1} is zero then  we have the 
symmetry \eqref{symm}. For example, if $\varphi$ is constant along the 
boundary then the divergence-free requirement on $\bm u$ implies that the 
integral at the right hand side of \eqref{w1} is zero. Also, if $\bm \theta$ and $
\bm u$ are periodic on a box then $\varphi$ is periodic and the boundary integral 
in \eqref{w1} is zero.
In fact, the field $\varphi$ arising from the Helmoltz-Hodge 
decomposition of $\bm \theta$ is the unique solution to the  Neuman problem (see \cite{ChorinMarsden}, p. 36)
\begin{equation}
\nabla^2 \varphi = \nabla\cdot \bm \theta\quad \textrm{(in $V$)},\qquad 
\frac{\partial \varphi}{\partial n}=\bm \theta \cdot \widehat {\bm n} 
\quad \textrm{(in $\partial V$)}.
\end{equation}
Clearly, if $\bm \theta$ is periodic then $\varphi$ is periodic and therefore $\bm \eta$ is periodic too. In other words, If we evaluate $\Phi$ on a space of periodic functions then we have the invariance 
\begin{equation}
\Phi([\bm \theta])=\Phi([\bm \eta]),
\end{equation}
 where $\bm \eta$ is periodic and divergence free ($\bm \eta$ is the 
 divergence-free part arising from the unique Helmholtz-Hodge 
 decomposition of $\bm \theta$.

\subsubsection{Divergence-Free Function Spaces} 
\label{sec:DivFreeSF}
There has been a significant research activity in identifying 
bases for divergence-free spaces of functions. 
For example, Deriaz and Perrier \cite{Deriaz1,Deriaz} have 
developed an effective algorithm to construct divergence-free 
and curl-free wavelets in 2D and 3D with various types of 
boundary conditions. 
Other divergence-free bases can be constructed in terms 
of radial basis functions \cite{Fuselier}, trigonometric 
polynomials \cite{Landriani}, or eigenvalue problems with 
appropriate boundary conditions conditions \cite{Venturi_IJHMT,Ozisik}.  
Hereafter we discuss how to construct a divergence-free 
basis for two-dimensional periodic flows. 
To this end, we consider the tensor product basis
\begin{equation}
\psi_n(x,y)=l_{j(n)}(x)l_{i(n)}(y),
\end{equation}
where $j(n)$ and $i(n)$ are suitable sequences of integer number while $l_k(x)$ are trigonometric polynomials. Any scalar-valued periodic function on the square (such as the streamfunction) can be represented as 
\begin{equation}
\Psi(x,y)=\sum_{k=1}^M \alpha_k \psi_j(x,y).
\end{equation}
Next define the divergence-free basis 
\begin{equation}
\bm \Gamma_k (x,y)=\left(\frac{\partial \psi_k(x,y) }{\partial y},\, 
-\frac{\partial \psi_k(x,y) }{\partial x} \right).
\label{basis}
\end{equation}
 It is clear that each basis element $\bm \Gamma_k$ is 
 divergence-free by construction, i.e., $\nabla \cdot \bm \Gamma_k=0$. However, the basis \eqref{basis} is not orthogonal nor normalized relative to the $L_2([-b,b]^2)$ inner product
 \begin{align}
\left(\bm \Gamma_i,\bm \Gamma_j\right)=&\int_{-b}^b\int_{-b}^b 
\bm \Gamma_i\cdot \bm \Gamma_j dxdy,\nonumber\\
=&\int_{-b}^b\int_{-b}^b 
\left(
\frac{\partial \psi_i}{\partial y}\frac{\partial \psi_j}{\partial y}+
\frac{\partial \psi_i}{\partial x}\frac{\partial \psi_j}{\partial x}
\right) dxdy.\nonumber
 \end{align}
To generate a divergence-free and orthonormal basis one 
could use the Gram-Schmidt orthogonalization. 
As we shall see hereafter, such procedure preserves 
the divergence-free character of each basis element. 
We first normalize $\bm \Gamma_1$ as 
\begin{equation}
\widehat{\bm \Gamma}_1(x,y) = \frac{\bm \Gamma_1(x,y)}{\left\|\bm \Gamma_1(x,y)\right\|}, \qquad \left\|\bm \Gamma_1\right\|=
\sqrt{(\bm \Gamma_1,\bm \Gamma_1)}.
\end{equation}  
Clearly, $\nabla \cdot \widehat {\bm \Gamma}_1=0$ since the $L_2$ norm of $\bm \Gamma_1$ is just a real number.
Next, define
\begin{equation}
\bm \Theta_2= \bm \Gamma_2- (\bm \Gamma_2,\widehat{\bm 
\Gamma}_1)\widehat{\bm \Gamma}_1,\qquad 
\widehat{\bm \Gamma}_2(x,y)=\frac{\bm \Theta_2(x,y)}{\left\|\bm
\Theta_2(x,y)\right\|}.
\end{equation}
As before, $\nabla \cdot \widehat{\bm \Gamma}_2=0$. Moreover 
$\widehat{\bm \Gamma}_2$ is orthogonal to $\widehat{\bm \Gamma}_1$, 
i.e., $(\widehat{\bm \Gamma}_1,\widehat{\bm \Gamma}_2)=0$. 
The algorithm proceeds with the computation of 
\begin{equation}
\bm \Theta_3= \bm \Gamma_3- (\bm \Gamma_3,\widehat{\bm 
\Gamma}_2)\widehat{\bm \Gamma}_2-(\bm \Gamma_3,\widehat{\bm 
\Gamma}_1)\widehat{\bm \Gamma}_1 ,\qquad \widehat{\bm \Gamma}_3(x,y)=\frac{\bm \Theta_3(x,y)}{\left\|\bm \Theta_3(x,y)\right\|}.
\end{equation}
Thanks to the fact that $\bm \Gamma_3$, $\widehat{\bm \Gamma}_2$ 
and $\widehat{\bm \Gamma}_1$ are divergence free, we have 
that $\nabla \cdot \widehat{\bm \Gamma}_3=0$. 
Moreover $\widehat{\bm \Gamma}_3$ is orthogonal to both 
$\widehat{\bm \Gamma}_2$ and $\widehat{\bm \Gamma}_1$. 
In other words, $\{\widehat{\bm \Gamma}_1,\widehat{\bm \Gamma}_2,
\widehat{\bm \Gamma}_3\}$ is a divergence-free 
orthonormal system.  Proceeding in a similar way, we can 
construct the divergence-free orthonormal basis we 
were looking for, and define the following finite-dimensional 
divergence-free space of functions
\begin{equation}
D_M=\textrm{span}\{\bm \Gamma_1, \cdots, \bm \Gamma_M\}.
\label{DM}
\end{equation}
An element of $D_M$ is in the form 
\begin{equation}
\bm \theta(x,y)=\sum_{k=1}^Ma_k\hat{\bm \Gamma}_k(x,y).
\label{divfree}
\end{equation}

\subsubsection{Analytical Solution to the Characteristic Function Equation}
A substitution of \eqref{divfree} into the 
Navier-Stokes-Hopf equation \eqref{h1} yields the 
multivariate (complex-valued) PDE\footnote{\color{r}
We emphasize that equation \eqref{NSCHF} has 
exactly the same structure as the characteristic function 
equation we obtain for the one-dimensional 
Burgers equation.  To show this, it is 
sufficient to discretize the Burgers-Hopf equation
\begin{equation}
\frac{\partial \Phi([\theta],t)}{\partial t}=\int_{-b}^b\theta(x)\left[i\frac{\partial }{\partial x}\frac{\delta^2 \Phi([\theta],t)}{\delta \theta(x)^2}+\nu\frac{\partial^2 }{\partial x^2} \frac{\delta \Phi([\theta],t)}{\delta \theta(x)} \right]d x, 
\label{h2}
\end{equation}
in the finite-dimensional space $D_m$ 
of periodic functions in $[-b,b]$. To 
this end, consider the series expansion
\begin{equation}
\theta(x)=\sum_{k=1}^m a_k\varphi_k(x).
\label{sexP}
\end{equation}
where $\varphi_k(x)$ are orthonormal trigonometric 
polynomials. Substituting \eqref{sexP} into 
\eqref{h2} and evaluating all functional 
derivatives in $D_m$ yields an equation in the form \eqref{NSCHF}.}
\begin{equation}
\frac{\partial \phi}{\partial t}= 
i \sum_{p,j,k=1}^M A_{pjk}a_p 
\frac{\partial^2\phi}{\partial a_k\partial a_j}+
\nu \sum_{k,p=1}^M a_pB_{pk}\frac{\partial \phi}{\partial a_k},
\label{NSCHF}
\end{equation}
where 
\begin{equation}
B_{pk}=\int_V \widehat{\bm \Gamma}_p\cdot 
\nabla^2\widehat{\bm \Gamma}_k d \bm x, \qquad 
A_{pjk}=\int_V \widehat{\bm \Gamma}_p\cdot 
\left[\left(\widehat{\bm \Gamma}_k\cdot \nabla\right)\widehat{\bm \Gamma}_j \right]d \bm x.
\label{coefficients}
\end{equation}
By using integration by parts we can simplify $B_{pk}$ to
\begin{align}
B_{pk}&=\int _{V}\left[\widehat \Gamma^{(x)}_{p} \left(\frac{\partial ^2 \widehat \Gamma^{(x)}_k}{\partial x^2}+\frac{\partial ^2 \widehat \Gamma^{(x)}_k}{\partial y^2}\right)+ 
\widehat \Gamma^{(y)}_{p} \left(\frac{\partial ^2 \widehat \Gamma^{(y)}_k}{\partial x^2}+\frac{\partial ^2 \widehat \Gamma^{(y)}_k}{\partial y^2}\right)\right]
d\bm x, \nonumber \\
&=-\int _{V}\left[\frac{\partial\widehat  \Gamma^{(x)}_{p}}{\partial x} \frac{\partial \widehat \Gamma^{(x)}_k}{\partial x}+\frac{\partial \widehat \Gamma^{(x)}_{p}}{\partial y} \frac{\partial \widehat \Gamma^{(x)}_k}{\partial y}
+ \frac{\partial \widehat \Gamma^{(y)}_{p}}{\partial x} \frac{\partial \widehat \Gamma^{(y)}_k}{\partial x}+\frac{\partial \widehat \Gamma^{(y)}_{p}}{\partial y} \frac{\partial \widehat \Gamma^{(y)}_k}{\partial y}\right]
d\bm x.
\end{align}
Similarly, the coefficients $A_{pjk}$ can be written as 
\begin{align}
A_{pkj}&=\int _{V}\left[
\Gamma^{(x)}_{p} \left(\Gamma^{(x)}_k
\frac{\partial \Gamma_j^{(x)}}{\partial x}+\Gamma^{(y)}_k
\frac{\partial \Gamma_j^{(x)}}{\partial y}\right)+ 
\Gamma^{(y)}_{p} \left(\Gamma^{(x)}_k
\frac{\partial \Gamma_j^{(y)}}{\partial x}+\Gamma^{(y)}_k
\frac{\partial \Gamma_j^{(y)}}{\partial y}
\right)
\right]
d\bm x.
\end{align}
%
%
% FROM HERE
By taking the inverse Fourier  transform of the characteristic 
function equation \eqref{NSCHF} we obtain 
\begin{equation}
\frac{\partial p(\bm u,t)}{\partial t}+\sum_{p=1}^M
\frac{\partial }{\partial u_p}\left[\left(\nu\sum_{k=1}^M u_k B_{pk} -
\sum_{k,j=1}^M u_ku_jA_{pkj} \right) p(\bm u,t)\right]=0.
\label{jointPDFCh}
\end{equation}
This is the multivariate first-order PDE that governs the 
evolution of the finite-dimensional approximation of the probability density functional\footnote{Equation \eqref{jointPDFCh}
can be obtained by evaluating the Navier-Stokes 
probability density functional equation \cite{Dopazo,Obrien} 
in the finite-dimensional test function space \eqref{DM}.}.
The formal solution to \eqref{jointPDFCh} is
\begin{equation}
p(\bm u,t)=p_0\left(\bm U(t,\bm u)\right)\exp\left(-\int_0^t \nabla\cdot 
\bm G(\bm u(\tau,\bm U)) d\tau \right),
\label{solutionPDf}
\end{equation}
where 
\begin{equation}
G_p(\bm u) = \nu\sum_{k=1}^M u_k B_{pk} - \sum_{k,j=1}^M u_ku_jA_{pkj}, \qquad p=1,...,M.
\end{equation}
The flow map $\bm u(t,\bm U)$ and its inverse $\bm U(t,\bm u)$ 
are defined by the solution to the ODE system
\begin{equation}
\frac{d\bm u}{dt }=\bm G(\bm u), \qquad \bm u(0)=\bm U.
\label{ODE1}
\end{equation}
By taking the Fourier transform of \eqref{solutionPDf}, we obtain 
the following analytical solution to the discretized 
Navier-Stokes-Hopf equation \eqref{NSCHF}
\begin{equation}
\phi(\bm a,t)=\int_{-\infty}^\infty\cdots \int_{-\infty}^\infty
e^{i\bm u\cdot \bm a}p_0\left(\bm U(t,\bm u)\right)\exp\left(-\int_0^t \nabla\cdot 
\bm G(\bm u(\tau,\bm U)) d\tau \right) d\bm u.
\label{solutionChF}
\end{equation}
Numerical methods to compute the flow map $\bm u(t, \bm U)$ and 
its inverse $\bm U(t,\bm u)$ in equations \eqref{solutionPDf} and 
\eqref{solutionChF} often rely on ray tracing, i.e., numerical 
solutions to the \eqref{ODE1} on a Cartesian mesh 
of initial conditions. Such methods are
not efficient in high-dimensions \cite{You2017}, and may be
improved significantly by representing flow maps 
using, e.g., tensor formats.

\paragraph{Functional Derivatives} 
The the first-order functional derivatives of the Hopf characterisitic functional (evaluated in the space of divergence-free functions) can be approximated  as 
\begin{align}
\frac{\delta \Phi([\bm \theta],t)}{\delta \theta^{(x)}(\bm x)}=\sum_{k=1}^M
\frac{\partial \phi}{\partial a_k}\widehat\Gamma_k^{(x)}(\bm x),\qquad 
\frac{\delta \Phi([\bm \theta],t)}{\delta \theta^{(y)}(\bm x)}=&\sum_{k=1}^M\frac{\partial \phi}{\partial a_k}\widehat\Gamma_k^{(y)}(\bm x).
\end{align}
Similarly, the second-order derivatives can be written as
\begin{align}
\frac{\delta^2 \Phi([\bm \theta],t)}{\delta \theta^{(x)}(\bm x) \delta \theta^{(x)}(\bm y)}=&\sum_{k,j=1}^M
\frac{\partial^2 \phi}{\partial a_k a_j}\widehat\Gamma_k^{(x)}(\bm x)\widehat\Gamma_j^{(x)}(\bm y),\\
\frac{\delta^2 \Phi([\bm \theta],t)}{\delta \theta^{(y)}(\bm x) \delta \theta^{(y)}(\bm y)}=&\sum_{k,j=1}^M
\frac{\partial^2 \phi}{\partial a_k a_j}\widehat\Gamma_k^{(y)}(\bm x)\widehat\Gamma_j^{(y)}(\bm y),\\
\frac{\delta^2 \Phi([\bm \theta],t)}{\delta \theta^{(x)}(\bm x) \delta \theta^{(y)}(\bm y)}=&\sum_{k,j=1}^M
\frac{\partial^2 \phi}{\partial a_k a_j}\widehat\Gamma_k^{(x)}(\bm x)\widehat\Gamma_j^{(y)}(\bm y).
\end{align}
Evaluating  these derivatives at $\bm \theta =0$ is equivalent to evaluate the derivatives of the characteristic function at $(a_1,...,a_M)=(0,...,0)$. This yields, the following representation of the mean and cross correlation of the velocity field 
\begin{equation}
\left<\bm u(\bm x,t)\right>=\sum_{k=1}^M \left.
\frac{\partial \phi(\bm a,t)}{\partial a_k}\right|_{(0,...,0)} \widehat{\bm \Gamma}_k(\bm x),
\end{equation}
\begin{equation}
\left<u^{(x)}(\bm x,t)u^{(y)}(\bm y,t)\right>=\sum_{k,j=1}^M \left.
\frac{\partial^2 \phi(\bm a,t)}{\partial a_k \partial a_j}\right|_{(0,...,0)}  \widehat\Gamma^{(x)}_k(\bm x)\widehat\Gamma^{(y)}_j(\bm y).
\end{equation}

\paragraph{Computational Complexity} 
The number of dimensions $M$ appearing 
in the characteristic function equation \eqref{NSCHF} and 
the joint PDF equation \eqref{jointPDFCh} coincides with the  
number of number of degrees of freedom 
we employ in the discretization of the velocity field. 
For instance, if we consider the classical 
two-dimensional Kolmogorov flow \cite{Lucas} 
represented on a $128 \times 128$ Fourier basis, 
then $M=16384$ ($128^2$). Computing the solution 
to such high-dimensional linear PDEs (i.e.,  \eqref{NSCHF} 
or \eqref{jointPDFCh}) obviously requires 
parallel algorithms and a highly-efficient tensor methods. 
If we employ operator splitting in time, then we can easily take 
care of the linear part in  \eqref{jointPDFCh} -- i.e., the 
one depending linearly on $u_k$ -- by integrating out 
the corresponding dynamics with the method 
of characteristics. 
To this end, let us first 
consider an orthonormal divergence-free basis that 
diagonalizes the matrix $B_{pk}$ in 
\eqref{coefficients}. Such basis exists and it can be 
computed by standard linear algebra 
techniques  \cite{ZhengDong2011} 
(simultaneous diagonalization of two quadratic forms).  
Relative to the new basis, the joint PDF equation \eqref{jointPDFCh} can be written as 
\begin{equation}
\frac{\partial p(\bm u,t)}{\partial t} = -\sum_{p=1}^M
\frac{\partial }{\partial u_p}\left[\left(\nu u_pB_{pp}- \sum_{j,k=1}^M u_ku_jA_{pkj}\right)p(\bm u,t)\right],
\end{equation}
where, with some abuse of notation, we denoted 
by $B_{pp}$ and $A_{pkj}$ the entries of  $B_{pk}$ 
and $A_{pkj}$ in \eqref{coefficients} 
relative to the new basis. The last equation can be 
written in the operator form 
\begin{equation}
\frac{\partial p(\bm u,t)}{\partial t} = L(\bm u) p(\bm u,t)+ Q(\bm u)p(\bm u,t),
\label{opjPDF}
\end{equation}
where 
\begin{equation}
L(\bm u) = - \sum_{p=1}^M
\nu u_p B_{pp}\frac{\partial }{\partial u_p} - B_{pp}\nu + \sum_{j,p=1}^M u_j
\left(A_{ppj}+A_{pjp}\right), \qquad 
Q(\bm u)=
\sum_{k,j,p=1}^M u_ku_jA_{pkj}\frac{\partial }{\partial u_p}.
\end{equation}
Note that $L(\bm u)$ depends linearly on $\bm u$
 while $Q(\bm u)$ is quadratic in $\bm u$. Both $L$ and $Q$ are separable linear operators. 
The formal solution to \eqref{opjPDF} is\footnote{We recall 
that $\exp[t(L+Q)]$ is the Frobenious-Perron 
operator associated with the quadratic 
dynamical system \eqref{ODE1} (see \cite{Dominy2017}).}
\begin{equation}
p(\bm u,t)=e^{t\left(L(\bm u)+Q(\bm u)\right)}p(\bm u,0).
\label{prop}
\end{equation}
By using operator splitting in time, e.g., the classical 
second-order Strang splitting, we can write
\begin{equation}
p(\bm u,t_n)=e^{\Delta t L(\bm u)/2 }e^{\Delta t Q(\bm u) }e^{\Delta t L(\bm u)/2 } p(\bm u,t_{n-1}).
\end{equation}
Clearly, the action of the semigroup $\exp [t L(\bm u)]$ can be 
computed exactly since the flow map corresponding  
characteristic system associated with the first order linear 
PDE $\dot{g} = L g$ is trivial. Specifically,  for 
any $g_0$ and any $t$ we have 
\begin{equation}  
e^{t L(\bm u)}g_0(\bm u)= 
g_0\left(e^{-t\nu B_{11}}u_1,...,e^{-t\nu B_{MM}}u_M\right)
\exp\left[- t \nu \sum_{j=1}^M B_{jj} + \sum_{j=1}^M \int_{0}^t 
e^{t\nu B_{jj}}u_{j0}\sum_{p=1}^M \left(A_{ppj}+A_{pjp}\right)\right].
\end{equation}
The action of the semigroup $\exp[t Q(\bm u)]$ is much more 
complicated to compute. Indeed, the flow map of  quadratic 
dynamical systems is, in general, not known explicitly. For small 
$\Delta t$ one may introduce the approximation 
$\exp[\Delta Q(\bm u)]\simeq I + \Delta t Q(\bm u)$ 
(or any higher-order one arising, e.g., from multi-step methods) 
and leverage on the separability of the operator $Q(\bm u)$ 
(rank $M^3$).
%

%\subsubsection{Probability Density Functional Equation }
%By taking the continuum limit of the joint PDF equation 
%\eqref{jointPDFCh} we obtain 
%the following equation governing the probability density functional 
%of any stochastic solution to the Navier Stokes equations
%\begin{align}
%\frac{\partial P([\bm a],t)}{\partial t} + \int_V \bm a(\bm x) \frac{\delta}{\delta \bm a(\bm x)}\left(\left[-(\bm a(x)\cdot \nabla) 
%a_j(\bm x) +\nu\nabla^2 \bm a_j(\bm x)\right] P([\bm a],t)\right)d\bm x=0.
%\label{PrDFeq2NS}
%\end{align}
%This equation can be derived obtained by taking the functional Fourier 
%transfrom of the Hopf equation \eqref{h1}, 
%and then performing functional integration by parts. 

}

\subsection*{Acknowledgments}
This research was supported by the Air Force Office 
of Scientific Research grant FA9550-16-1-0092.

%\newpage
\appendix
\section{Functional Fourier Transform}
\label{app:functional fourier transform}
Functional Fourier transforms can be defined as a 
continuum limit of multi-variate 
Fourier transforms. To introduce this concept 
in a simple and intuitive way, let us consider the 
expression of the joint probability density function
of $n$ random variables in terms of the joint characteristic function
\begin{equation}
p(a_1,...,a_n)=\frac{1}{(2\pi)^{n/2}}
\int_{-\infty}^\infty\cdots\int_{-\infty}^\infty e^{-i(a_1b_1+\cdots a_nb_n)}
\phi(b_1,...,b_n)db_1\cdots db_n. \label{fint}
\end{equation}
If we think of $a_i$ and $b_i$ as values of two 
continuous functions $a(x)$ and $b(x)$ at 
locations $x_i\in\mathbb{R}$, then it makes sense to 
consider what happens to $p$ and $\phi$ as we send the 
number of collocation points to infinity. 
In this limit, the joint probability 
density function $p(a_1,...,a_n)$ becomes a probability density 
functional $P([a(x)])$, while the joint characteristic 
function $\phi(b_1,...,b_n)$ becomes a (Hopf) characteristic 
functional. This allows us to define the functional Fourier transform as
\begin{equation}
 P([a(x)])=
\int e^{-i
\int a(x)b(x)dx}
\Phi([b(x)])\mathcal{D}[b(x)],\label{Fint}
\end{equation}
where 
\begin{equation}
\mathcal{D}[b(x)]=\lim_{n\rightarrow\infty}\frac{1}{(2\pi)^{n/2}}\prod_{i=1}^n db(x_i).
\end{equation}
Equation \eqref{Fint} establishes the 
connection between the Hopf functional and the probability 
density functional of a random field. It also allows us 
to define special functionals, such as the Dirac 
delta functional 
\begin{equation}
 \delta[a(x)]=\int e^{-i
\int a(x)b(x)dx}\mathcal{D}[b(x)].\label{functional_Dirac}
\end{equation}
By using this definition, it is easy to show that the 
probability density functional can be expressed as an average 
of Dirac delta functional 
\begin{align}
 P([a(x)])=&\left<\delta[a(x)-u(x;\omega)]\right>\nonumber\\
=&\left<\int e^{-i
\int a(x)b(x)dx+i\int u(x;\omega)b(x)dx}\mathcal{D}[b(x)]\right>\nonumber\\
=&\int e^{-i
\int a(x)b(x)dx}\left<e^{i\int u(x;\omega)b(x)dx}\right>\mathcal{D}[b(x)]\nonumber\\
=&\int e^{-i
\int a(x)b(x)dx}\Phi([b(x)])\mathcal{D}[b(x)].\nonumber
\end{align}
%
% \subsection*{Transformation of Functional Derivatives}
% 
% Another interesting question is how functional derivatives 
% transform under the functional Fourier transform. To answer 
% this question, let us first recall the following continuum 
% limits 
% \begin{equation}
%  \frac{\partial \phi}{\partial t}=
%  \sum_{i,j=1^n}b_iD^{(2)}_{ij}\frac{\partial \phi}{\partial b_j}
% \end{equation}
% 

\section{Evaluation of Functional Integrals}
\label{sec:functional integrals}
Functional integrals in the form\footnote{In \eqref{fi} $F$ 
is a measurable functional on a complete metric 
space $D(F)$ and $W>0$ is a functional integral 
measure.}
\begin{equation}
\int_{D(F)} F([\theta])W([\theta])\mathcal{D}[\theta]
\label{fi}
\end{equation}
arise naturally in many branches of physics 
and mathematics, e.g., in quantum and statistical mechanics, 
field theory, quantum optics, solid state physics, and financial 
mathematics \cite{Kleinert,Justin,Beran,Amit}.
Analytical results for functional integrals 
are available only in few exceptional cases, 
and therefore one usually has to resort to numerical approximation 
techniques \cite{Egorov,Popov}. 
Perhaps, the most classical one is Monte Carlo (MC), 
which gives results as an ensemble average over a 
large number of realizations of the test 
function $\theta(x)$, i.e., 
\begin{equation} 
\int F([\theta])W[\theta]\mathcal{D}[\theta]\simeq 
\frac{1}{N} \sum_{k=1}^N F([\theta_k(x)]).
\end{equation}
Each term in the series corresponds to a specific 
{\em path} $\theta_k(x)$, which is drawn from 
the probability measure $W([\theta])$. 
Monte Carlo and quasi-Monte Carlo methods \cite{Dick}
have slow convergence rate and therefore 
they require significant computational resources to achieve 
good accuracy.
In a renormalized perturbation theory setting \cite{McComb}, 
the functional integral can be expanded relative to a certain parameter
and the terms in the series expansion have a structure determined 
by the order of the perturbation parameter. This yields well-know diagrammatic 
representations of functional (path) integrals, as a sum over 
all possible paths consistent with the order of the perturbation 
parameter \cite{Itzykson,Kleinert}.
Alongside statistical methods and perturbation series 
expansions, deterministic techniques were also proposed 
to evaluate functional integrals. These include, 
path summation based on short time propagators \cite{Rosa-Clot}, fast 
Fourier transforms \cite{Eydeland}, HDMR expansions 
\cite{Wasilkowski} and functional approximation 
techniques \cite{Lobanov}. 
An intuitive way to calculate functional integrals 
is to consider a finite-dimensional subset of $D_m\subset D(F)$, e.g., 
the linear span of the basis $\{\varphi_1,...,\varphi_m\}$. 
In this setting, $\theta(x)$ admits the representation 
\begin{equation}
\theta(x)=\sum_{k=1}^ma_k\varphi_k(x),\qquad a_k=\left(\theta,\varphi_k\right).
 \label{g11}
\end{equation}
If we substitute \eqref{g11} into \eqref{fi} we obtain 
the multivariate integral
\begin{equation}
\int \cdots \int f\left(a_1,...,a_m\right)w(a_1,...,a_m)
da_1\cdots da_m,
\label{ffi1}
\end{equation}
where 
\begin{equation}
f(a_1,...,a_m)=F\left(\left[\sum_{k=1}^ma_k\varphi_k(x)\right]\right), \qquad w(a_1,...,a_m)=W\left(\left[\sum_{k=1}^ma_k\varphi_k(x)\right]\right).
\end{equation}
This form includes both ``modal'' and ``nodal'' 
discretizations of the functional integral \eqref{fi}. 
For example, if \eqref{g11} is 
an interpolant through $m$ nodes along the $x$ axis, 
then $a_j=\theta(x_j)$ and \eqref{ffi1} becomes 
\begin{equation}
\int \cdots \int  f\left(\theta(x_1),...,\theta(x_m)\right)
w(\theta(x_1),...,\theta(x_m))d\theta(x_1)\cdots d\theta(x_m).
\label{fi2}
\end{equation}
In this formulation, we obtain the functional integral \eqref{fi} 
as the limit of an infinite number of nodes $m\rightarrow \infty$. 
This is the viewpoint taken in quantum mechanics and field theory  
to define path integrals \cite{Kleinert,Justin}, i.e., integrals 
over suitable trajectories of functions.
{\color{r}
From a mathematical viewpoint, the limiting procedure 
defining the functional integral 
measure in terms of an infinite products of elementary 
measures should be handled with care. In fact,  the classical Lebesgue measure does not exist in spaces of infinite 
dimension \cite{Marzucchi}. On the other hand, Gaussian measures 
are still well defined in such setting. This is why we included
$W([\theta])$ in \eqref{continuous_ip_W0}. The argument leading 
to the result on non-existence of an analogue to the Lebesgue 
measure in infinite dimension is related to the argument showing 
that the Heine-Borel theorem does not hold in infinite-dimensional 
normed linear spaces

\paragraph{Remark} Roughly speaking, the meaning 
of $\mathcal{D}[\theta]$ in \eqref{fi} is: 
integrate over all possible degrees of freedom defining 
$\theta$ in the class of functions $D(F)$. If $\theta$ is 
in the span of a finite-dimensional basis (e.g., \eqref{g11}), 
then there is no problem in replacing \eqref{Fint2} with 
an integral over the range of the expansion coefficients , 
as we did in \eqref{ffi1}. A more tricky situation arises 
when $\theta$ is truly infinite-dimensional. In this case, 
the functional integral \eqref{fi} is an integral over 
an infinite number of variables. 
To ensure convergence the functional integral 
measure $W([\theta])$ must be carefully chosen. 
}
{\color{r}
\subsection{Functional Integrals of Cylindrical Functionals}
\label{app:change of variables in SSE}
Functional integrals involving cylindrical functionals 
(see Section \ref{sec:tensor}) can be written in the 
general form
\begin{equation}
\int_{D(f)} f\left((\theta,\varphi_1),...,(\theta,\varphi_m)\right)W([\theta])\mathcal{D}[\theta].
\label{Fint2}
\end{equation}

\noindent
{\em Example 1:} Consider the nonlinear functional 
\begin{equation}
f((\theta,\sin(x)))=\frac{1}{1+(\theta,\sin(x))^2}\qquad 
(\theta,\sin(x))=\int_{0}^{2\pi} \theta(x)\sin(x)dx, 
\end{equation}
in the space of infinitely-differentiable periodic 
functions in $[0,2\pi]$, i.e., 
\begin{equation}
\theta\in D(f)=\{\theta\in {C}^{\infty}([0,1]),\, 
\theta(0)=\theta(2\pi)\}.
\end{equation}
Also, consider the functional integral measure
\begin{equation}
W([\theta]) = \kappa e^{-(\theta,\theta)^2},
\end{equation}
where $\kappa$ is a (possibly infinite) normalization 
constant. We expand $\theta$ in a Fourier series, 
\begin{equation}
\theta(x) = a_0 + \sum_{k=1}^\infty a_k \sin(kx)+
\sum_{k=1}^\infty b_k \cos(kx).
\end{equation}
In this representation, the functional integral measure 
becomes
\begin{equation}
w(a_0,a_1,b_1, ... ) = 2\sqrt{\pi} e^{-4\pi^2 a_0^2}
\lim_{j\rightarrow \infty}\prod_{k=1}^j \pi^k e^{-\pi^2 (a^2_k+b^2_k)}
\end{equation}
and the integral \eqref{Fint2} can be computed analytically. 
The result is 
\begin{align}
\int_{D(f)} \frac{W([\theta])\mathcal{D}[\theta] }{1+(\theta,\sin(x))^2} 
&= 
\int_{-\infty}^\infty\cdots \int_{-\infty}^\infty
\frac{w(a_0,a_1,b_1, ... )}{1+\pi^2a_1^2} da_0\prod_{k=1}^
\infty da_kdb_k\nonumber\\
&=\sqrt{\pi}\int_{-\infty}^{\infty}\frac{e^{-\pi^2a_1^2}}{1+\pi^2a_1^2} da_1\nonumber\\
&=- \sqrt{\pi}e(\textrm{erf}(1) - 1),
\end{align}
where $\textrm{erf}$ is the error 
function and $e$ is the Napier number. 
\begin{flushright}
$\Box$
\end{flushright}

\noindent
An interesting question is whether integrals in the form 
\eqref{Fint2} can be simplified and reduced to 
integrals over the range of $\{(\theta,\varphi_1), ..., (\theta,\varphi_m)\}$, i.e., a subset of $\mathbb{R}^m$. 
The answer to this question is, in general, negative, as 
have just shown in the previous example. However, 
if we pick $\theta$ in a finite-dimensional function 
space, e.g., spanned by 
the basis $\{\xi_1(x),...,\xi_N(x)\}$, i.e., 
\begin{equation}
\theta(x)=\sum_{k=1}^N a_k\xi_k(x),
\label{p1}
\end{equation}
then the functional integral effectively reduces to 
a finite-dimensional integral over the range of the 
expansion coefficients $a_1,...,a_N$.
A substitution of \eqref{p1} 
into \eqref{Fint2} yields 
\begin{equation}
\int_{D_N(f)} f\left((\theta,\varphi_1),...,(\theta,\varphi_m)\right)
W([\theta])\mathcal{D}[\theta]\simeq \int\cdots \int 
g(a_1,...,a_N)w(a_1,..,a_N)da_1\cdots da_N,
\label{uyt}
\end{equation}
where  
\begin{equation}
g(a_1,...,a_N)=f\left(\sum_{k=1}^N a_k \alpha_{1k},...,
\sum_{k=1}^N a_k \alpha_{mk}\right), \qquad 
\alpha_{ij}=(\xi_k,\varphi_j).
\end{equation}
In particular, if the functions $\xi_k(x)$ coincide with 
$\varphi_k(x)$, and $\{\varphi_1,...,\varphi_N\}$ is 
orthonormal with respect to the inner product $(,)$, then 
we obtain
\begin{equation}
\int_{D_N(f)} f\left((\theta,\varphi_1),...,(\theta,\varphi_m)\right)
W([\theta])\mathcal{D}[\theta]= \int\cdots \int f\left(a_1,...,a_m\right)
w(a_1,..,a_N) da_1\cdots a_N.
\label{measure}
\end{equation}
If $w(a_1,..,a_N)$ is a probability density, then \eqref{measure}
implies that 

\begin{equation}
\int_{D_N(f)} f\left((\theta,\varphi_1),...,(\theta,\varphi_m)\right)
W([\theta])\mathcal{D}[\theta] =\int\cdots \int f\left(a_1,...,a_m\right)
w(a_1,..,a_m) da_1\cdots a_m.
\end{equation}
The the last integral can be computed if the measure $w$ is 
separable, and $f$ is represented in terms of a tensor decomposition, e.g., a canonical tensor series or an 
HDMR expansion (see Section \ref{sec:tensor} and Section \ref{sec:HDMR}). 
}

\section{Derivation of Functional Differential Equations}
\label{sec:derivation}
In this Appendix we briefly discuss how to derive functional 
differential equations, and their equivalence to systems of 
PDEs, or PDEs in an infinite number of variables.

\subsection{The Generating Functional Approach}
The seminal work of Martin, Siggia and Rose \cite{Martin}  
opened the possibility to apply quantum field theoretic methods, 
such as Feynman diagrams and Schwinger-Dyson equations, 
to classical and statistical physics. The key idea 
is to construct a generating functional (action functional) 
based on the equations of motion 
\cite{Jensen,Phythian,Daniele_JMathPhys} and 
then apply a Heisenberg operator theory which parallels 
the Schwinger formalism of quantum field theory. This 
allows to obtain closed equations for quantities of 
interest such as the correlation functions and  
response functions (averaged Green's functions). 
For a through description of the method see the 
excellent review paper of Jensen \cite{Jensen}.

\subsection{From PDEs to Functional Equations: The Method of Continuum Limits}
\label{sec:continuum_limits}
In Section \ref{sec:FDEs approximation} we have discussed how  
derive probability density functional equations and Hopf 
functional equations as continuum limits of joint PDF equations 
and joint characteristic function equations, 
respectively. Such procedure is rather formal, but 
easy to follow and effective in many cases. More importantly,
it can be applied to many different linear and nonlinear PDEs. 

{
\color{r}
\vs
\noindent
{\em Example 1:}
Let us consider the Kuramoto-Sivashinsky equation 
\begin{equation}
\frac{\partial u}{\partial t}+\frac{1}{2}\frac{\partial u^2}{\partial x}+
\frac{\partial^2 u }{\partial x^2}+\frac{\partial^4 u }{\partial x^4}=0
\label{eq3}
\end{equation}
with random initial data $u(x,0;\omega)=u_0(x;\omega)$.
In a collocation setting, the solution to \eqref{eq3} can be written as 
\begin{equation}
 u(x,t;\omega)=\sum_{k=1}^n u_k(t;\omega)l_k(x),
 \label{seq3}
\end{equation}
where $u_k(t;\omega)=u(x_k,t;\omega)$ while 
$l_j(x)$ are suitable basis functions, e.g., Lagrange characteristic 
polynomials through the points $\{x_1,...,x_n\}$. 
Substituting \eqref{seq3} into \eqref{eq3} and setting 
the residual at collocation points equal to zero yields
the system 
\begin{equation}
 \frac{d u_k}{dt}+\frac{1}{2}\sum_{j=1}^n D^{(1)}_{kj}u^2_j+\sum_{j=1}^n\left(D^{(2)}_{kj}+D^{(4)}_{kj}\right)u_j=0 \qquad k=1,...,n
\end{equation}
where $D^{(1)}_{kj}$, $D^{(2)}_{kj}$  and are first-, 
second-, and fourth-order differentiation matrices, 
respectively. The joint characteristic function 
of $\{u_1,...,u_n\}$ is 
\begin{equation}
\phi(\theta_1,...,\theta_n,t)=\left<\exp\left[i\sum_{j=1}^n \theta_j u_j(t;\omega)\right]\right> 
\end{equation}
and it satisfies the evolution equation
\begin{align}
\frac{\partial \phi}{\partial t}=&-\frac{i}{2}\sum_{k,j=1}^n \theta_k D^{(1)}_{kj}\left<u^2_j
\exp\left[i\sum_{j=1}^n \theta_j u_j)\right]\right>- i\sum_{k,j=1}^n \theta_k (D^{(2)}_{kj}+D^{(4)}_{kj})\left<u_j
\exp\left[i\sum_{j=1}^n \theta_j u_j)\right]\right>\nonumber\\
= &\frac{i}{2}\sum_{k,j=1}^n \theta_k D^{(1)}_{kj}
\frac{\partial^2 \phi}{ \partial \theta^2_j}-
\sum_{k,j=1}^n \theta_k (D^{(2)}_{kj}+D^{(4)}_{kj})\frac{\partial \phi}{\partial \theta_j}.
\end{align}
By taking the continuum limit of this equation, i.e., by sending $n$ 
to infinity, we  formally obtain
\begin{equation}
\frac{\partial \Phi}{\partial t}=\int_a^b\theta(x)
\left[\frac{i}{2}\frac{\partial }{\partial x}\frac{\delta^2\Phi([\theta])}{\delta\theta(x)^2} 
-\frac{\partial^2 }{\partial x^2}\frac{\delta \Phi([\theta])}{\delta\theta(x)}-\frac{\partial^4 }{\partial x^4}\frac{\delta \Phi([\theta])}{\delta\theta(x)}\right]dx.
\end{equation}
The continuum-limit argument just invoked is
formal and can be used for other finite-dimensional 
equations (see \cite{Waclawczyk}). 
In the next Section we discuss the inverse operation, i.e., how to 
transform functional differential equations into PDEs 
with a finite number of variables. 
}

\subsubsection{From Functional Equations to PDEs}
\label{app:equivalenceHopf}
The restriction of functionals and functional differential equations 
to finite-dimensional function spaces yields multivariate fields 
and multivariate PDEs, respectively. Such PDEs can be obtained directly 
by using the finite dimensional theory, i.e., functional equations on 
finite-dimensional function spaces are {\em in toto} equivalent to 
multivariate PDEs.
To show this, let us consider the finite-dimensional Hilbert space
\begin{equation}
D_m=\textrm{span}\{\varphi_1(x),...,\varphi_m(x)\},
 \label{SN}
\end{equation}
where ${\varphi_1, ..., \varphi_m}$ are orthonormal basis 
functions. Let us represent the test function $\theta(x)$ 
relative to such orthonormal basis
\begin{equation}
\theta_m(x)=\sum_{k=1}^m a_k\varphi_k(x), \qquad a_k=(\theta,\varphi_k)
\label{exP}
\end{equation}
Restricting the domain of the Hopf functional \eqref{hopf_functional} 
to the finite-dimensional space \eqref{SN} (see Figure \ref{fig:1}) 
yields the joint characteristic function
\begin{equation}
\phi(a_1,...,a_m)=\left<e^{i(a_1U_1(\omega)+\cdots+ a_m U_m(\omega))}\right>,
\qquad U_k(\omega)=\int_a^b u(x;\omega)\varphi_k(x)dx.
\label{Uk}
\end{equation}
Note that the linear combination of $U_k(\omega)$ yields a 
finite-dimensional approximation of $u(x,\omega)$, i.e., 
\begin{equation}
\widehat{u}(x;\omega)=\sum_{k=1}^m U_k(\omega)\varphi_k(x).
\label{findim}
\end{equation}
Based on this expansion, the mean field and the correlation function 
of $\widehat{u}(x,\omega)$ can be easily expressed as 
\begin{align}
\left<\widehat{u}(x,\omega)\right>=&
\sum_{k=1}^m \left<U_k(\omega)\right>\varphi_k(x)=
\sum_{k=1}^m \frac{1}{i}\left.\frac{\partial \phi}{\partial a_k}\right|_{0}
\varphi_k(x), \label{meanh}\\
\left<\widehat{u}(x,\omega)\widehat{u}(y;\omega)\right>=&
\sum_{k,j=1}^m\left<U_k(\omega)U_j(\omega)\right>\varphi_k(x)\varphi_j(x)=-
\sum_{k,j=1}^m\left.\frac{\partial^2 \phi}{\partial a_k\partial a_j}\right|_{0}
\varphi_k(x)\varphi_j(x)\label{varh}.
\end{align}
Similar equations hold for higher-order moments. Next we 
consider the restriction of functional equations to 
finite-dimensional function spaces. 

\vs
\noindent
{\em Example 1:} 
Evaluating the Hopf functional  
equation in the finite-dimensional space $D_m$ yields a 
characteristic function equation that can be obtained 
directly by using simpler methods, e.g., by differentiating the 
characteristic function. To show this, consider  
the diffusion equation $\partial u/\partial t=\partial^2 u/\partial x^2$, 
in a bounded domain $[a,b]$, with Dirichlet boundary conditions 
and a random initial condition. The corresponding Hopf equation is
\begin{equation}
\frac{\partial \Phi([\theta],t)}{\partial t}=\int_{a}^b\theta(x)
\frac{\partial^2}{\partial x^2}\left[\frac{\delta \Phi([\theta],t)}
{\delta \theta(x)}\right]dx.
\label{o1}
\end{equation}
Now, we restrict $\theta$ to the finite-dimensional function space
$D_m$. Evaluating the Hopf functional for $\theta\in D_m$ yields 
the multivariate characteristic function
\begin{equation}
\phi(a_1,....,a_m,t)=\Phi([\theta_m],t)\qquad \theta_m\in D_m.
\label{o2}
\end{equation}
A substitution of \eqref{ha1} and \eqref{o2} into \eqref{o1} yields
\begin{equation}
 \frac{\partial \phi}{\partial t}=
\sum_{k,j=1}^m a_k\frac{\partial \phi}{\partial a_j}
H_{kj},\qquad H_{kj}=\int_{a}^b\varphi_k(x)
\frac{\partial^2 \varphi_j(x)}{\partial x^2}dx.
\label{coincide}
\end{equation}
This is an equation for the joint characteristic function of the 
coefficients  $\{U_k(\omega)\}_{k=1,...,m}$ defined in \eqref{Uk}. 
Such equation can be obtained directly by using 
the semidiscrete (Galerkin) form of the diffusion equation
\begin{equation}
\frac{du_j(t;\omega)}{dt}=\sum_{k=1}^m u_k(t;\omega)\int_{a}^b\varphi_k(x)
\frac{\partial^2 \varphi_j(x)}{\partial x^2}dx\qquad j=1, ...,m. 
\end{equation}
Indeed, by differentiating the joint characteristic function 
of $\{u_1,...,u_m\}$ we obtain
\begin{align}
\frac{\partial \phi}{\partial t} = &
\sum_{k=1}^m i\left<\frac{du_k}{dt} a_k\exp
\left[i\sum_{j=1}^m u_k(t;\omega) a_k\right]\right>, \nonumber\\ 
=&\sum_{k,j=1}^m a_kH_{kj}\frac{\partial \phi}{\partial a_j},
\end{align}
which coincides with \eqref{coincide}.

\subsubsection{From Functional Equations to Systems of PDEs}

In this Section we discuss the connection between functional 
differential equations and systems of PDEs. The main idea is 
to evaluate the FDE on function spaces of increasing 
dimension (see, e.g., Section \ref{sec:interpolationnodes})
To illustrate the method, we restrict our 
attention to the simple Hopf equation \eqref{o1}. 
We have seen already that evaluating the Hopf 
functional at $\theta(x)=\alpha\delta(y-x)$ 
yields the one-point characteristic function 
of $u(y,t;\omega)$ 
\begin{equation}
\Phi([\theta],t)=\left<e^{i \alpha \int_a^b u(x,t;\omega)
\delta(y-x)dx}\right>=\left<e^{i \alpha u(y,t;\omega)}\right>.
\end{equation}
Similarly, if we consider
\begin{equation}
 \theta(x) \in \textrm{span}\{\delta(y_1-x),...,\delta(y_m-x)\}
\end{equation}
then we have that $\Phi([\theta],t)$ coincides with the 
joint characteristic function of $\{u(y_k,t;\omega)\}_{k=1,..,m}$.
Now, let us see what happens to the Hopf equation 
when we evaluate it at 
$\theta(x)=\alpha_1 \delta(y_1-x)$. We have
\begin{align}
\frac{\partial \phi(\alpha_1,y_1,t)}{\partial t}=&\alpha_1 \int_{a}^b\delta(y_1-x)
\frac{\partial^2}{\partial x^2}\left[\frac{\delta \Phi([\theta],t)}
{\delta \theta(x)}\right]_{\theta=\alpha_1\delta(y_1-x)}dx,\nonumber\\
=&i\alpha_1\left<\frac{\partial ^2 u(y_1,t;\omega)}{\partial y_1^2}
e^{i\alpha u(y_1,t;\omega)}\right>.
\label{oo1}
\end{align}
This equation can also  be obtained directly by differentiating 
the one-point characteristic function $\phi(\alpha,y,t)$ with 
respect to time and then using the diffusion equation.
Note that the term at the right hand side 
of \eqref{oo1} cannot be easily expressed in terms 
of the one-point characteristic function. However, if we use 
the two-point function we have 
\begin{equation}
\left[\frac{\partial^3\phi(\alpha_1,\alpha_2,y_1,y_2,t)}
{\partial \alpha_2 \partial y_2^2}\right]_{\substack{\alpha_2=0\\y_2=y_1}}=
\left<\frac{\partial ^2 u(y,t;\omega)}{\partial y^2}
e^{i\alpha u(y,t;\omega)}\right>.
\end{equation}
Therefore, we obtain 
\begin{align}
\frac{\partial \phi(\alpha_1,y_1,t)}{\partial t}=&
i\alpha_1\left[\frac{\partial^3
\phi(\alpha_1,\alpha_2,y_1,y_2,t)}
{\partial \alpha_2 \partial y_2^2}\right]_{\substack{\alpha_2=0
\\y_2=y_1}}
\end{align}
At this point we need an evolution equation 
for the two point characteristic function 
$\phi(\alpha_1,\alpha_2,y_1,y_2,t)$. 
Such equation may be obtained from \eqref{o1} by using the 
test function
$\theta(x)=\alpha_1\delta(y_1-x)+\alpha_2\delta(y_2-x)$.
If we perform the calculation we find that such evolution depends on 
the three-point characteristic function. 
If we keep get on going  we obtain 
{\em an infinite hierarchy} of finite-dimensional characteristic 
function PDEs. Taking the Fourier transform of such hierarchy 
yields the well-known Bogoliubov-Born-Green-Kirkwood-Yvon (BBGKY)
hierarchy of probability density function equations 
\cite{Montgomery,HeyrimJCP_2014}.
Other functional equations yield other types of 
hierarchies of PDEs which should be studied on a case-by-case 
basis. 
Another approach to convert an FDE to an infinite hierarchy of 
PDEs was developed my Thomas Lundgren in \cite{Lundgren}. 
The hierarchy of PDEs is known as Lundgren-Monin-Novikov 
hierarchy, and it can be constructed by a direct approach 
\cite{Wacawczyk1,Hosokawa,Friedrich} or by 
functional integration of the probability density 
functional equation (see Section \ref{sec:PDENFUNEQ}). 

%%%%%%

%\newpage
\bibliographystyle{plain}
\bibliography{report}

\end{document}